\documentclass[10pt]{report}
\usepackage[top=1.2in, bottom=1in, left=1in, right=1in,includefoot,heightrounded]{geometry}

\usepackage[utf8]{inputenc}
\usepackage{graphicx}
\graphicspath{ {images/} }
\usepackage{caption}
\usepackage{subcaption}
\usepackage{float}
\usepackage{fancyhdr}
\renewcommand{\headrulewidth}{0pt}
\usepackage{amsmath}
\usepackage{amsthm}
\usepackage{amssymb, mathrsfs, xcolor}
\usepackage{enumerate}
\usepackage{verbatim}
\usepackage{enumitem}
\makeatletter
\newcommand{\mylabel}[2]{#2\def\@currentlabel{#2}\label{#1}}
\makeatother

\usepackage{indentfirst}
\usepackage{hyperref}
\hypersetup{
    colorlinks=true,
    linkcolor=blue,
    citecolor=brown,
    filecolor=magenta,
    urlcolor=cyan,
}
\usepackage{titlesec}
\usepackage{fancyhdr}
\usepackage{mathtools}
\usepackage{CJKutf8}
\usepackage{enumitem}
\usepackage{mleftright}

\numberwithin{equation}{section}
\newtheorem{mainthm}{Theorem}
\newtheorem{keythm}{Theorem}[chapter]
\newtheorem{thm}{Theorem}[section]

\newtheorem{cor}[thm]{Corollary}
\newtheorem{lem}[thm]{Lemma}
\newtheorem{prop}[thm]{Proposition}

\newtheorem{ques}[thm]{Question}
\newtheorem{prob}[thm]{Problem}

\newtheorem{cor-}[keythm]{Corollary}
\newtheorem{lem-}[keythm]{Lemma}
\newtheorem{prop-}[keythm]{Proposition}
\newtheorem{ques-}[keythm]{Question}
\newtheorem{prob-}[keythm]{Problem}

\theoremstyle{definition}
\newtheorem{maindefn}[mainthm]{Definition}
\newtheorem{defn}[thm]{Definition}
\newtheorem{conv}[thm]{Convention}
\newtheorem{note}[thm]{Notation}
\newtheorem*{defn*}{Definition}

\newtheorem{remark}[thm]{Remark}
\newtheorem{example}[thm]{Example}

\newtheorem{defn-}[keythm]{Definition}
\newtheorem{conv-}[keythm]{Convention}
\newtheorem{note-}[keythm]{Notation}
\newtheorem{remark-}[keythm]{Remark}
\newtheorem{example-}[keythm]{Example}

\newtheorem*{remark*}{Remark}

\newcommand{\1}{\mathbf{1}}

\newcommand{\Af}{\mathfrak{A}}
\newcommand{\B}{\mathbb{B}}
\newcommand{\Bs}{\mathscr{B}}
\newcommand{\Bf}{\mathfrak{B}}
\newcommand{\C}{\mathbb{C}}
\newcommand{\Cc}{\mathcal{C}}
\newcommand{\Co}{\mathscr{C}}
\newcommand{\D}{\mathscr{D}}
\newcommand{\Do}{\mathcal{D}}

\newcommand{\E}{\mathcal{E}}
\newcommand{\Ex}{\mathbf{E}}
\newcommand{\Eb}{\mathbb{E}}
\newcommand{\F}{\mathcal{F}}

\newcommand{\Fs}{\mathscr{F}}
\newcommand{\fb}{\mathbf{f}}
\newcommand{\gb}{\mathbf{g}}

\newcommand{\Ga}{\mathbf{\Gamma}}
\newcommand{\Gc}{\mathcal{G}}

\newcommand{\Hb}{\mathbb{H}}

\newcommand{\Ic}{\mathcal{I}}
\newcommand{\Mf}{\mathfrak{M}}
\newcommand{\Mc}{\mathcal{M}}
\newcommand{\Nf}{\mathfrak{N}}

\newcommand{\Oh}{\mathcal{O}}
\newcommand{\Pv}{\mathbf{P}}
\newcommand{\Pc}{\mathcal{P}}
\newcommand{\Pf}{\mathfrak{P}}

\newcommand{\R}{\mathbb{R}}
\newcommand{\Rf}{\mathfrak{R}}
\newcommand{\Rc}{\mathcal{R}}
\newcommand{\Se}{\mathcal{S}}
\newcommand{\Sp}{\mathbb{S}}
\newcommand{\Sc}{\mathscr{S}}
\newcommand{\Sf}{\mathfrak{S}}

\newcommand{\Tc}{\mathcal{T}}
\newcommand{\Tf}{\mathfrak{T}}
\newcommand{\U}{\mathcal{U}}

\newcommand{\Us}{\mathscr{U}}
\newcommand{\V}{\mathcal{V}}
\newcommand{\Vs}{\mathscr{V}}

\newcommand{\X}{\mathfrak{X}}
\newcommand{\Xs}{\mathscr{X}}
\newcommand{\Ys}{\mathscr{Y}}
\newcommand{\Z}{\mathbb{Z}}
\newcommand{\Zf}{\mathfrak{Z}}
\newcommand{\Zs}{\mathscr{Z}}
\newcommand{\Gr}{\mathbf{Gr}}

\newcommand{\LogL}{\mathrm{LogL}}
\newcommand{\logl}{\mathrm{logl}}
\newcommand{\Lip}{\mathrm{Lip}}
\newcommand{\codiff}{\vartheta}

\newcommand{\re}{\operatorname{Re}}
\newcommand{\im}{\operatorname{Im}}
\newcommand{\Su}{\mathsf{S}}
\newcommand{\De}{\mathsf{P}}
\newcommand{\id}{\operatorname{id}}
\newcommand{\Inv}{\mathrm{Inv}}

\newcommand{\sgn}{\operatorname{sgn}}
\newcommand{\essup}{\mathop{\operatorname{essup}}}
\newcommand{\eps}{\varepsilon}

\newcommand{\Span}{\operatorname{Span}}

\newcommand{\Vol}{\operatorname{Vol}}
\newcommand{\loc}{\mathrm{loc}}
\newcommand{\pv}{\operatorname{p.\!v.}}
\newcommand{\supp}{\operatorname{supp}}
\newcommand{\dist}{\operatorname{dist}}

\newcommand{\divg}{\operatorname{div}}
\newcommand{\rank}{\operatorname{rank}}

\newcommand{\Coorvec}[1]{\frac\partial{\partial#1}}
\newcommand{\Lie}[1]{\mathrm{Lie}_{#1}}

\newcommand{\und}[1]{{\textcolor{red}{#1}}}

\title{On Rough Frobenius-type Theorems}
\author{Liding Yao}
\date{1 Aug 2022}

\begin{document}
\pagenumbering{gobble}

\fancyhf{} 
\fancyhead[RO,R]{\thepage} 
\renewcommand{\headrulewidth}{0pt}

\begin{center}
    \Large
    \textbf{On Rough Frobenius-type Theorems and Their H\"older Estimates}
    
    \let\thefootnote\relax\footnotetext{The work in this paper represents the author's PhD thesis at the University of Wisconsin-Madison. It includes the results in \textsf{\href{https://arxiv.org/abs/2002.07973}{arXiv:2002.07973}}, \textsf{\href{https://arxiv.org/abs/2004.07288}{arXiv:2004.07288}}, \textsf{\href{https://arxiv.org/abs/2202.07729}{arXiv:2202.07729}} and a part of \textsf{\href{https://arxiv.org/abs/2105.10120}{arXiv:2105.10120}}.}
    
    
    \vspace{0.4cm}
    {Liding Yao}
    
    \vspace{0.9cm}
    \textbf{Abstract}

\end{center}

    The thesis is dedicated to giving several results on Frobenius-type theorems in non-smooth settings, and giving H\"older regularity estimates for the respective coordinate systems. 
    
    For the real Frobenius theorem, we extend the definition of involutivity to non-Lipschitz tangent subbundles using generalized functions. We prove the Frobenius Theorem with sharp regularity when the subbundle is log-Lipschitz: if $\mathcal V$ is a log-Lipschitz involutive subbundle of rank $r$, then for any $\varepsilon>0$, locally there is a homeomorphism $\Phi(u,v)$ such that $\Phi,\frac{\partial\Phi}{\partial u^1},\dots,\frac{\partial\Phi}{\partial u^r}\in\mathscr C^{1-\varepsilon}$, and $\mathcal V$ is spanned by the continuous vector fields $\Phi_*\frac\partial{\partial u^1},\dots,\Phi_*\frac\partial{\partial u^r}$.
    
    We also develop a singular version of the Frobenius theorem on log-Lipschitz vector fields. If $X_1,\dots,X_m$ are log-Lipschitz vector fields such that $[X_i,X_j]=\sum_{k=1}^mc_{ij}^kX_k$ for some generalized functions $c_{ij}^k$ that can be written as the derivatives of log-Lipschitz functions, then for any point $p$ there is a $C^1$-manifold $\mathfrak S$ containing $p$ such that $X_1,\dots,X_m$ span the tangent space at every point in $\mathfrak S$.
    
    Nirenberg's famous complex Frobenius theorem gives necessary and sufficient conditions on a locally integrable structure for when the manifold is locally diffeomorphic to $\mathbb R^r\times\mathbb C^m\times \mathbb R^{N-r-2m}$ through a coordinate chart $F$ in such a way that the structure is locally spanned by $F^*\frac\partial{\partial t^1},\dots,F^*\frac\partial{\partial t^r},F^*\frac\partial{\partial z^1},\dots,F^*\frac\partial{\partial z^m}$, where we have given $\mathbb R^r\times\mathbb C^m \times\mathbb R^{N-r-2m}$ coordinates $(t,z,s)$.  When the structures are differentiable, we give the optimal H\"older-Zygmund regularity for the coordinate charts which achieve this realization.  Namely, if the structure has H\"older-Zygmund regularity of order $\alpha>1$, then the coordinate chart $F$ that maps to $\mathbb R^r\times\mathbb C^m \times\mathbb R^{N-r-2m}$ may be taken to have H\"older-Zygmund regularity of order $\alpha$, and this is sharp.  Furthermore, we can choose this $F$ in such a way that the vector fields $F^*\frac\partial{\partial t^1},\dots,F^*\frac\partial{\partial t^r},F^*\frac\partial{\partial z^1},\dots,F^*\frac\partial{\partial z^m}$ on the original manifold have H\"older-Zygmund regularity of order $\alpha-\varepsilon$ for every $\varepsilon>0$, and we give an example to show that the regularity for $F^*\frac\partial{\partial z}$ is optimal.
    
    Similarly we give a counterexample for the $C^k$-version of the Newlander-Nirenberg theorem: we give an example of $C^k$-integrable almost complex structure that does not admit a corresponding $C^{k+1}$-complex coordinate system.
    
    Combining the technics from the log-Lipschitz real Frobenius theorem and the sharp complex Frobenius theorem, we show that if a complex Frobenius structure $\mathcal S$ is $\mathscr C^\alpha$ ($\frac12<\alpha\le1$) such that $\mathcal S+\bar{\mathcal S}$ is log-Lipschitz, then for every $0<\eps<2\alpha-1$ there is a homeomorphism $\Phi(t,z,s)$ such that $\Phi,\frac{\partial\Phi}{\partial t},\frac{\partial\Phi}{\partial z}\in\mathscr C^{2\alpha-1-\eps}$, and $\mathcal S$ is spanned by the $\mathscr C^{2\alpha-1-\eps}$ vector fields $\Phi_*\frac\partial{\partial t^1},\dots,\Phi_*\frac\partial{\partial t^r},\Phi_*\frac{\partial}{\partial z^1},\dots,\Phi_*\frac\partial{\partial z^m}$. 
    
    On the quantitative side, we show that if $X_1,\dots,X_m$ are log-Lipschitz real vector fields such that $[X_i,X_j]=\sum_{k=1}^mc_{ij}^kX_k$ holds for some $c_{ij}^k\in\mathscr C^{\alpha-1}$ where $1<\alpha<2$, then on each leaf where $X_1,\dots,X_m$ span the tangent spaces we can find a regular parameterization $\Phi$ such that $\Phi^*X_1,\dots,\Phi^*X_m$ are $\mathscr C^\alpha$, and their $\mathscr C^\alpha$ norm depend only on the diffeomorphic invariant quantities of $X_1,\dots,X_m$, like the norms $\|c_{ij}^k\|_{C^{0,\alpha-1}_X}$.

\tableofcontents



\chapter{Introduction}
\pagenumbering{arabic} 


Let $\Mf$ be a smooth manifold, and let $\Se\le T\Mf$ (or $\le \C T\Mf$) be a real or complex tangent subbundle. We call $\Se$ \textbf{involutive}, if for any (real or complex) vector fields $X,Y$ that are sections of $\Se$, their Lie bracket $[X,Y]$ is also a section of $\Se$.

The Frobenius-type theorems describe how certain types of involutive subbundles are related to the span of coordinate vector fields. In the thesis, we are going to talk about the following theorems in nonsmooth settings: Let $\Mf$ be a $n$-dimensional smooth manifold,  $\Se\le \C T\Mf$ be a smooth involutive subbundle, and $p\in\Mf$ be a fixed point.
\begin{itemize}[parsep=-0.3ex]
    \item The real Frobenius theorem: If\footnote{We use $\overline{\Se}$ as the complex conjugate subbundle, see Definition \ref{Defn::ODE::ConjBundle}.} $\Se=\overline{\Se}$, then there is a smooth parameterization $\Phi(t,s):\Omega\subseteq\R^r\times\R^{n-r}\to\Mf$ near $p$ such that $\frac{\partial\Phi}{\partial t^1}(t,s),\dots,\frac{\partial \Phi}{\partial t^r}(t,s)$ span $\Se_{\Phi(t,s)}$ for all $(t,s)\in\Omega$. Here $r=\rank\Se$.
    \item The Newlander-Nirenberg theorem: If $\Se\oplus\overline{\Se}=\C T\Mf$, then there is a smooth parameterization $\Phi(z):\Omega\subseteq\C^m\to\Mf$ near $p$ such that $\frac{\partial\Phi}{\partial z^1}(z),\dots,\frac{\partial \Phi}{\partial z^m}(z)$ span $\Se_{\Phi(z)}$ for all $z\in\Omega$. Here $m=n/2$.
    \item Nirenberg's complex Frobenius theorem: If  $\Se+\overline{\Se}$ is also a smooth subbundle and is involutive, then there is a smooth parameterization $\Phi(t,z,s):\Omega\subseteq\R^r\times\C^m\times\R^{n-r-2m}\to\Mf$ near $p$ such that $\frac{\partial\Phi}{\partial t^1}(t,z,s),\dots,\frac{\partial \Phi}{\partial t^r}(t,z,s),\frac{\partial\Phi}{\partial z^1}(t,z,s),\dots,\frac{\partial \Phi}{\partial z^m}(t,z,s)$ span $\Se_{\Phi(t,z,s)}$ for all $(t,z,s)\in\Omega$. Here $r=\rank(\Se\cap\overline{\Se})$ and $m=\rank\Se-r$. This result generalizes the above two special cases.
\end{itemize}

In the thesis we discuss several estimates of these type of theorems in the non-smooth setting, as well as some examples that indicate some of the H\"older exponents are optimal.

We remark that for the results in Section \ref{Section::MainResults}, Theorem \ref{MainThm::LogFro} follows from \cite{YaoLLFro}; Theorems \ref{MainThm::CpxFro} and \ref{MainThm::Sharpddz} follow from \cite{YaoCpxFro}; Theorem \ref{MainThm::CounterNN} follows from \cite{LidingCounterNN}; and Theorem \ref{MainThm::QuantFro} is the part of \cite{StreetYaoVectorFields}. The Theorems \ref{MainThm::SingFro}, \ref{MainThm::RoughFro1} and \ref{MainThm::RoughFro2} are new results.


\section{The Main Results}\label{Section::MainResults}

For $\alpha\in(0,\infty)$, let $\Co^\alpha$ be the class of H\"older-Zygmund functions\footnote{For non-integer exponents, the H\"older-Zygmund space agrees with the classical H\"older space. When $\alpha\in\Z_+$ we have $C^{\alpha-1,1}\subsetneq \Co^\alpha$. See Remark \ref{Rmk::Hold::ZygVsLip}.} of order $\alpha$ (See Definitions \ref{Defn::Intro::DefofHold}). Let $\Co^{\alpha-}=\bigcap_{\eps>0}\Co^{\alpha-\eps}$ and $\Co^{\alpha+}=\bigcup_{\eps>0}\Co^{\alpha+\eps}$. We use $\Co^\infty=C^\infty$ as the class of smooth functions.

Let $U\subseteq\R^n$ be an open subset, we denote by $\Co^\Lip$ the class of Lipschitz functions.
We say a continuous function $f:U\to\R^m$ is bounded \textbf{log-Lipschitz}, denoted by $f\in\Co^\LogL(U;\R^m)$, if
\begin{equation*}\label{Eqn::Intro::ClogNorm}
    \|f\|_{\Co^\LogL(U;\R^m)}:=\sup_{x\in U}|f(x)|+\sup_{x,y\in U;0<|x-y|<1}|f(x)-f(y)|\Big(|x-y|\log\frac e{|x-y|}\Big)^{-1}<\infty.
\end{equation*}

We say $f:U\to\R^m$ is bounded \textbf{little log-Lipschitz}, denoted by $f\in\Co^\logl(U;\R^m)$, if in addition $\lim_{r\to0}\sup_{x,y\in U;0<|x-y|<r}|f(x)-f(y)|(|x-y|\log\frac e{|x-y|})^{-1}=0$. See Definition \ref{Defn::Hold::CLogLandLip}.

Note that we have $\Co^{1+\eps}\subsetneq\Co^\Lip\subsetneq\Co^1\subsetneq\Co^\LogL\subsetneq\Co^{1-}$ and $\Co^\Lip\subsetneq\Co^\logl\subsetneq\Co^\LogL$.

For $\Co^{\frac12+}$-subbundles (in particular for log-Lipschitz subbundles) we introduce the notion of distributional\footnote{Throughout the thesis, we use \textbf{distributions} for \textbf{generalized functions}, which are defined to be linear functionals on some test function spaces. For the corresponding terminology in differential geometry, we always use \textbf{tangent subbundles}, or \textbf{singular tangent subbundles}.} involutivity as follows. 

\setcounter{mainthm}{-1}
\begin{maindefn}\label{Defn::Intro::DisInv}
Let $\Mf$ be a $\Co^{\frac32+}$ manifold and let $\Se\le \C T\Mf$ be a $\Co^{\frac12+}$ complex tangent subbundle. We say $\Se$ is \textbf{involutive in the sense of distributions}, if for every $\Co^{\frac12+}$ vector fields $X=\sum_{i=1}^nX^i\Coorvec{x^i}$, $Y=\sum_{j=1}^nY^j\Coorvec{x^j}$ on $\Se$ and every $\Co^{\frac12+}$ 1-form $\theta=\sum_{k=1}^n\theta_kdx^k$ on $\Se^\bot$, the following generalized function
\begin{equation}\label{Eqn::Intro::DefInvEqn}
    \langle\theta,[X,Y]\rangle:=\sum_{i,j=1}^n\Big(\theta_iX^j\frac{\partial Y^i}{\partial x^j}-\theta_jY^i\frac{\partial X^j}{\partial x^i}\Big),
\end{equation}
is identically zero in the sense of distributions on $\Mf$.
\end{maindefn}

Here $\Se^\bot:=\{(p,l)\in \C T^*\Mf:p\in \Mf,\ l(v)=0,\forall v\in\Se_p\}\le \C T^*\Mf$ is the dual bundle of $\Se$. 
Note that the right hand side of \eqref{Eqn::Intro::DefInvEqn} is not necessarily a $L^1_\loc$-function. It makes sense as a distribution, see Lemma \ref{Lem::Hold::MultLoc} and Corollary \ref{Cor::Hold::[X,Y]WellDef}. Also see Definition \ref{Defn::DisInv::DefFunVF}. Replacing real subbundles by complex subbbundles, the same definition also work for complex subbundles. 

We see that distributional involutivity is the natural generalization to the classical involutivity. See \ref{Section::DisInv::CharInv} for a detailed discussion.

The result on real Frobenius theorem is the following:

\begin{mainthm}[The Frobenius theorem for log-Lipschitz subbundles]\label{MainThm::LogFro}
Let $\Mf$ be $n$-dimensional $C^{1,1}$ manifold, and let $\V\le T\Mf$ be a rank-$r$ log-Lipschitz real tangent subbundle, which is involutive in the sense of distributions.

Then for any $p\in \Mf$, there is a neighborhood $\Omega\subseteq\R^r_t\times\R^{n-r}_s$ of $0$ and a continuous map $\Phi:\Omega\to \Mf$ such that
\begin{enumerate}[parsep=-0.3ex,label=(\arabic*)]
    \item\label{Item::MainThm::LogFro::0}$\Phi(0)=p$.
    \item\label{Item::MainThm::LogFro::Top}$\Phi:\Omega\to\Mf$ is a topological embedding, and $\frac{\partial\Phi}{\partial t^1},\dots,\frac{\partial\Phi}{\partial t^r}:\Omega\to T\Mf$ are continuous. 
    
    {\normalfont In order words the pushforward vector fields $\Phi_*\Coorvec{t^j}=\frac{\partial\Phi}{\partial t^j}\circ\Phi^\Inv$ ($1\le j\le r$) are defined pointwise.} 
    \item\label{Item::MainThm::LogFro::Span}$\V|_{\Phi(\Omega)}$ is spanned by $\Phi_*\Coorvec{t^1},\dots,\Phi_*\Coorvec{t^r}$. Equivalently, for any $(t,s)\in \Omega$, $\frac{\partial\Phi}{\partial t^1}(t,s),\dots,\frac{\partial\Phi}{\partial t^r}(t,s)\in T_{\Phi(t,s)}\Mf$ form a basis of the (real) linear subspace $\V_{\Phi(t,s)}$.
\end{enumerate}
Moreover,
\begin{enumerate}[parsep=-0.3ex,label=(\arabic*)]\setcounter{enumi}{3}
    \item\label{Item::MainThm::LogFro::BigLReg} For any $\eps>0$ there is an open subset $\Omega_\eps\subseteq\Omega$ of $(0,0)$ such that $\Phi\in\Co^{1-\eps}(\Omega_\eps;\Mf)$, its inverse $\Phi^\Inv\in\Co^{1-\eps}(\Phi(\Omega_\eps);\Mf)$ and the partial derivatives $\frac{\partial\Phi}{\partial t^1},\dots,\frac{\partial\Phi}{\partial t^r}\in\Co^{1-\eps}(\Omega_\eps;T\Mf)$. 
    \item[\mylabel{Item::MainThm::LogFro::LittleLReg}{(4')}] Suppose $\V$ is  little log-Lipschitz, then automatically $\Phi\in\Co^{1-}(\Omega;\Mf)$, $\Phi^\Inv\in\Co^{1-}(\Phi(\Omega);\R^n)$ and  $\frac{\partial\Phi}{\partial t^1},\dots,\frac{\partial\Phi}{\partial t^r}\in\Co^{1-}(\Omega;T\Mf)$. 
    \item[\mylabel{Item::MainThm::LogFro::C1Reg}{(4'')}] Let $\alpha\in\{\Lip\}\cup(1,\infty]$. Suppose $\Mf\in\Co^\infty$ and $\V\in\Co^\alpha$ for some $\alpha>1$, then automatically $\Phi\in\Co^\alpha(\Omega;\Mf)$, $\Phi^\Inv\in\Co^\alpha(\Phi(\Omega);\R^n)$ and  $\frac{\partial\Phi}{\partial t^1},\dots,\frac{\partial\Phi}{\partial t^r}\in\Co^\alpha(\Omega;T\Mf)$. 
\end{enumerate}
\end{mainthm}
Here $\Co^{k+\Lip}=C^{k,1}=W^{k+1,\infty}$, $k=0,1,2,\dots$ are the classical Lipschitz spaces. See Definition \ref{Defn::Hold::CLogLandLip}.

\begin{remark*}For a general version of Theorem \ref{MainThm::LogFro}, see Theorem \ref{KeyThm::RealFro::ImprovedLogFro}. Here Theorem \ref{MainThm::LogFro} is slightly stronger than \cite[Theorems 1.4 and 1.7]{YaoLLFro}, where we treat the cases of log-Lipschitz and little log-Lipschitz separately.
\end{remark*}

The result Theorem \ref{MainThm::LogFro} \ref{Item::MainThm::LogFro::BigLReg} shows that both $\Phi$ and $\Phi^\Inv$ are $\Co^{1-\eps}$. In general $\Phi\in \Co^{1-\eps}$ does not imply $\Phi^\Inv\in \Co^{1-\eps}$ and vice versa, since we do not have the inverse function theorem for $\Co^{1-\eps}$-maps. This $\Co^{1-\eps}$-regularity is sharp, in the sense that there exists a log-Lipschitz subbundle such that either $\Phi\notin \Co^{1-}$ or $\Phi^\Inv\notin \Co^{1-}$, see Corollary \ref{Cor::Exmp::SharpRealFro} \ref{Item::Exmp::SharpRealFro::LogL}.

By considering the setting of singular tangent subbundles, as an analog to \cite{MontanariMorbidelliSingularFrobenius} we have the following
\begin{mainthm}[A singular Frobenius theorem on log-Lipschitz vector fields]\label{MainThm::SingFro}
Let $\Mf$ be a $C^{1,1}$-manifold, let $X_1,\dots,X_m$ be some locally log-Lipschitz vector fields on $\Mf$. Suppose there are some functions $c_{ij}^k\in \Co^{\LogL-1}_\loc(\Mf)$ (see Definitions \ref{Defn::Hold::LogL-1} and \ref{Defn::DisInv::DefFunVF}) such that in the sense of distributions in $\Mf$,
    \begin{equation}\label{Eqn::MainThm::SingFro::Inv}
        [X_i,X_j]=\sum_{k=1}^mc_{ij}^kX_k,\quad 1\le i,j\le m.
    \end{equation}
    
    Then the span of $X_1,\dots,X_q$ is pointwisely integrable. That is, for any $p\in\Mf$ there is a $r$-dimensional $C^1$ embedded submanifold $\Sf$ containing $p$, where $r$ is the rank of $\Span(X_1|_p,\dots,X_q|_p)\le T_p\Mf$, such that $X_1|_\Sf,\dots,X_q|_\Sf$ are vector fields on $\Sf$ and span the tangent space of $\Sf$ at every point.
    
    Moreover, the Carnot-Carath\'eodory ball $B_{(X_1,\dots,X_m)}(p,\infty)$ is an immerse submanifold, where $\Sf\subset B_X(p,\infty)$ is an open subset with respect to the (immersed) differential structure.
\end{mainthm}
Here for a collection $X=(X_1,\dots,X_m)$ of $C^0$-vector fields on $\Mf$ and for $R\in(0,\infty]$, we use $B_X(p,R)=\{q\in\Mf:\dist_X(p,q)<R\}$ where
{\small\begin{equation}\label{Eqn::Intro::CCDist}
        \dist_X(x,y):=\inf\Big\{d>0:\exists \gamma\in C^{0,1}([0,d];\Mf),\gamma(0)=x,\gamma(d)=y\text{ s.t. }\dot\gamma(t)=\sum_{j=1}^ma_j(t)X_j(\gamma(t))\text{ where } \sum_{j=1}^m|a_j(t)|^2\le1\Big\}.
\end{equation}}

In fact we can take $\Sf=B_X(p,R)$ for every $R<\infty$. The set $B_X(p,\infty)$ may not be an embedded submanifold as we can consider an irrational torus in $\Sp^1\times\Sp^1$.

\begin{remark*}
See Section \ref{Section::DisInv::Sing} for a further discussion of singular tangent subbundles and singular Frobenius theorems. I would like to thank Sylvain Lavau and  Andrew Lewis for many informative discussion on the related concepts and theorems.
\end{remark*}
\medskip

For complex Frobenius theorems we have the following:

\begin{mainthm}[Sharp H\"older regularity for Nirenberg's complex Frobenius theorem]\label{MainThm::CpxFro}
Let $\alpha\in(1,\infty)$, let $\Mf$ be an $n$-dimensional smooth manifold, and let $\Se$ be a $\Co^\alpha$-complex Frobenius structure such that $\rank(\Se\cap \bar\Se)=r$ and $\rank\Se=r+m$. Then for any $p\in \Mf$, there is a neighborhood $\Omega\subseteq\R^r_t\times\C^m_z\times\R^{n-r-2m}_s$ of $0$ and  a $\Co^\alpha$-regular parameterization $\Phi:\Omega\to\Mf$ (i.e. $\Phi$ and $\Phi^\Inv:\Phi(\Omega)\to\R^n$ are both $\Co^\alpha$) near $p$, such that 
\begin{enumerate}[parsep=-0.3ex,label=(\arabic*)]
    \item\label{Item::MainThm::CpxFro::0}$\Phi(0)=p$.
    \item\label{Item::MainThm::CpxFro::Span}For any $(t,z,s)\in \Omega$, $\frac{\partial\Phi}{\partial t^1}(t,z,s),\dots,\frac{\partial\Phi}{\partial t^r}(t,z,s),\frac{\partial\Phi}{\partial z^1}(t,z,s),\dots,\frac{\partial\Phi}{\partial z^m}(t,z,s)\in\C T_{\Phi(t,z,s)}\Mf$ form a basis of the complex linear subspace $\Se_{\Phi(t,z,s)}$.
    \item\label{Item::MainThm::CpxFro::Reg1} $\Phi_*\Coorvec{t^1},\dots,\Phi_*\Coorvec{t^r},\Phi_*\Coorvec{z^1},\dots,\Phi_*\Coorvec{z^m}$ are $\Co^{\alpha-}$-complex vector fields on $\Phi(\Omega)$.
\end{enumerate}
Suppose there is a $\beta\in[\alpha+1,\infty)$ such that $\Se\cap\bar\Se$ is a $\Co^\beta$-complex tangential subbundle. Then we can find a $\Phi$ such that in addition to the results above, we have
\begin{enumerate}[parsep=-0.3ex,label=(\arabic*)]\setcounter{enumi}{3}
    \item\label{Item::MainThm::CpxFro::Reg2} $\Phi_*\Coorvec{t^1},\dots,\Phi_*\Coorvec{t^r}$ are $\Co^{\beta}$-vector fields on $\Phi(\Omega)$.
\end{enumerate}

\end{mainthm}We remark that the constructions of $\Phi$ for \ref{Item::MainThm::CpxFro::0} - \ref{Item::MainThm::CpxFro::Reg1} and for \ref{Item::MainThm::CpxFro::0} - \ref{Item::MainThm::CpxFro::Reg2} are different.

Combining the real and the complex Frobenius Theorem (Theorems \ref{MainThm::LogFro} and \ref{MainThm::CpxFro}), along with the definition of distributional involutivity (Definition \ref{Defn::Intro::DisInv}), we have the following:
\begin{mainthm}[Estimates of non-differentiable complex Frobenius structures I]\label{MainThm::RoughFro1}
Let $\alpha\in(\frac12,1]$, let $\Mf$ be an $n$-dimensional smooth manifold, and let $\Se$ be a $\Co^\alpha$-complex Frobenius structure such that $\rank(\Se\cap \bar\Se)=r$ and $\rank\Se=r+m$. Assume the subbundle $\Se+\bar\Se$ is log-Lipschitz. Then for any $\eps>0$ and any $p\in\Mf$  there is a neighborhood $\Omega=\Omega_\eps\subseteq\R^r_t\times\C^m_z\times\R^{n-r-2m}_s$ of $0$ and  a map $\Phi:\Omega \to\Mf$ near $p$, such that
\begin{enumerate}[parsep=-0.3ex,label=(\arabic*)]
    \item\label{Item::MainThm::RoughFro1::0}$\Phi(0)=p$.
    \item\label{Item::MainThm::RoughFro1::Top}$\Phi:\Omega\to\Mf$ is a topological embedding, and $\frac{\partial\Phi}{\partial t^1},\dots,\frac{\partial\Phi}{\partial t^r},\frac{\partial\Phi}{\partial z^1},\dots,\frac{\partial\Phi}{\partial z^m}:\Omega\to\C T\Mf$ are continuous. 
    \item\label{Item::MainThm::RoughFro1::Span}$\Se|_{\Phi(\Omega)}$ is spanned by $\Phi_*\Coorvec{t^1},\dots,\Phi_*\Coorvec{t^r},\Phi_*\Coorvec{z^1},\dots,\Phi_*\Coorvec{z^m}$.
    \item\label{Item::MainThm::RoughFro1::Reg1}$\Phi\in\Co^{2\alpha-1-\eps}(\Omega;\Mf)$, $\Phi^\Inv\in\Co^{2\alpha-1-\eps}(\Phi(\Omega);\Omega)$ and $\Phi_*\Coorvec{t^1},\dots,\Phi_*\Coorvec{t^r},\Phi_*\Coorvec{z^1},\dots,\Phi_*\Coorvec{z^m}$ are $\Co^{2\alpha-1-\eps}$ vector fields on $\Phi(\Omega)$.
\end{enumerate}
Suppose there is a $\beta\in\{1+\Lip\}\cup(2,\infty)$ such that $\Se\cap\bar\Se\in\Co^\beta$. Then  we can find a (possibly different) $\Phi$ such that in addition to the results above, we have
\begin{enumerate}[parsep=-0.3ex,label=(\arabic*)]\setcounter{enumi}{4}
    \item\label{Item::MainThm::RoughFro1::Reg2} $\Phi_*\Coorvec{t^1},\dots,\Phi_*\Coorvec{t^r},\Phi_*\Coorvec{z^1},\dots,\Phi_*\Coorvec{z^m}$ are $\Co^\beta$ vector fields on $\Phi(\Omega)$.
\end{enumerate}
\end{mainthm}

If $\Se+\bar\Se$ is little-Log Lipschitz, then all the $\Co^{2\alpha-1-\eps}$ regularity estimates in Theorem \ref{MainThm::RoughFro1} can be improved by $\Co^{(2\alpha-1)-}$. If $\Se+\bar\Se$ is Lipschitz, then they can all improved to be $\Co^{2\alpha-1}$:
\begin{mainthm}[Estimates of non-differentiable complex Frobenius structures II]\label{MainThm::RoughFro2}
Let $\alpha\in(\frac12,1]$, and let $\Mf,\Se$ be as in the assumption of Theorem \ref{MainThm::RoughFro1}. 

Suppose the subbundle $\Se+\bar\Se$ is little log-Lipschitz. Then for any $p\in\Mf$  there is a neighborhood $\Omega\subseteq\R^r_t\times\C^m_z\times\R^{n-r-2m}_s$ of $0$ and a map $\Phi:\Omega \to\Mf$ near $p$, such that the results \ref{Item::MainThm::RoughFro1::0}, \ref{Item::MainThm::RoughFro1::Top} and \ref{Item::MainThm::RoughFro1::Span} in Theorem \ref{MainThm::RoughFro1} hold, with regularity estimate:
\begin{enumerate}[parsep=-0.3ex,label=(\arabic*')]\setcounter{enumi}{3}
    \item\label{Item::MainThm::RoughFro2::Reg1}$\Phi\in\Co^{(2\alpha-1)-}(\Omega;\Mf)$, $\Phi^\Inv\in\Co^{(2\alpha-1)-}(\Phi(\Omega);\Omega)$ and $\Phi_*\Coorvec{t^1},\dots,\Phi_*\Coorvec{t^r},\Phi_*\Coorvec{z^1},\dots,\Phi_*\Coorvec{z^m}$ are $\Co^{(2\alpha-1)-}$ vector fields on $\Phi(\Omega)$.
\end{enumerate}

Suppose $\alpha<1$ and the subbundle $\Se+\bar\Se$ is Lipschitz. Then for any $p\in\Mf$  there is a neighborhood $\Omega\subseteq\R^r_t\times\C^m_z\times\R^{n-r-2m}_s$ of $0$ and a map $\Phi:\Omega \to\Mf$ near $p$, such that the results \ref{Item::MainThm::RoughFro1::0}, \ref{Item::MainThm::RoughFro1::Top} and \ref{Item::MainThm::RoughFro1::Span} in Theorem \ref{MainThm::RoughFro1} hold, with regularity estimate:
\begin{enumerate}[parsep=-0.3ex,label=(\arabic*'')]\setcounter{enumi}{3}
    \item\label{Item::MainThm::RoughFro2::Reg2}$\Phi\in\Co^{2\alpha-1}(\Omega;\Mf)$, $\Phi^\Inv\in\Co^{2\alpha-1}(\Phi(\Omega);\Omega)$ and $\Phi_*\Coorvec{t^1},\dots,\Phi_*\Coorvec{t^r},\Phi_*\Coorvec{z^1},\dots,\Phi_*\Coorvec{z^m}$ are $\Co^{2\alpha-1}$ vector fields on $\Phi(\Omega)$.
\end{enumerate}

Suppose there is a $\beta\in\{1+\Lip\}\cup(2,\infty)$ such that $\Se\cap\bar\Se\in\Co^\beta$. Then  we can find a (possibly different) $\Phi$ such that \ref{Item::MainThm::RoughFro1::Reg2} in Theorem \ref{MainThm::RoughFro1} also holds.
\end{mainthm}


For a general version of Theorems \ref{MainThm::CpxFro}, \ref{MainThm::RoughFro1} and \ref{MainThm::RoughFro2}, see Theorems \ref{KeyThm::CpxFro1} and \ref{KeyThm::CpxFro2} in Section \ref{Section::PfCpxFro}. Also see Section \ref{Section::CpxFro::Overview} for an overview of the proofs and their historical remarks.

\bigskip
In Theorem \ref{MainThm::CpxFro} (that $\alpha>1$), the results $\Phi\in\Co^\alpha$ and $\Phi_*\Coorvec z\in\Co^{\alpha-}$ are all optimal, in the sense that given $\alpha>1$, we have examples of complex Frobenius structures $\Se$, such that if $\Phi$ is a $C^1$-chart representing $\Se$, then $\Phi\notin\Co^{\alpha+\eps}$ for any $\eps>0$, and $\Phi_*\Coorvec z\notin\Co^\alpha$, respectively. See Section \ref{Section::ExampleOverview}. In Section \ref{Section::Sharpddz} we construct an example that meets the following:
\begin{mainthm}[Sharpness of $F^*\Coorvec z\in\Co^{\alpha-}$ in Theorem \ref{MainThm::CpxFro}]\label{MainThm::Sharpddz}
Let $\alpha>1$. There is a $\Co^\alpha$-complex Frobenius structure $\Se$ on $\R^3$, such that $\rank\Se =1$, $\rank(\Se+\bar\Se)=2$, and there does not exist a $C^1$-coordinate chart $F:U\subseteq\R^3\to\C^1_z\times\R^1_s$ near $0\in\R^3$ such that $F^*\Coorvec z$ spans $\Se|_U$ and $F^*\Coorvec z\in\Co^\alpha(U;\C^3)$.
\end{mainthm}
See also Corollary \ref{Cor::Exmp::ExdSharpddz} for a generalization.
The construction for Theorem \ref{MainThm::Sharpddz} is modified from a counterexample to $C^k$-regularity for the Newlander-Nirenberg Theorem:
	\begin{mainthm}[A counterexample to $C^k$-regularity for the Newlander-Nirenberg theorem]\label{MainThm::CounterNN}
	For every positive integer $k$ there is a $C^k$ integrable almost complex structure $J$ on $\R^{2n}$ such that we cannot find a neighborhood $U\subset\R^{2n}$ of $0$ and a complex coordinate chart $\phi\in C^{k+1}(U;\C^n)$ which maps $J$ to the standard complex structure.
	\end{mainthm}

So far all results above are qualitative, since we do not have additional objects to label the scale. When we are handed in a specified collection of generators $X_1,\dots,X_m$ for our subbundle $\V\le T\Mf$ (or $\Se\le\C T\Mf$), we can further ask what is the behaviour of the vector fields under the pullback of $\Phi$, and whether we are allowed to choose a good $\Phi$ such that the new vector fields $\Phi^*X_1,\dots,\Phi^*X_m$ are ``normalized''. 

In a collaboration with my advisor Brian Street, we manage to give the following:
\begin{mainthm}[A quantitative Frobenius theorem in low regularity]\label{MainThm::QuantFro}
Let $n,m\ge1$, $\alpha\in(1,2)$, $\mu_0>0$ and $M_0>0$. There are constants $\widehat K=\widehat K(n,\alpha,\mu_0,M_0)>0$ and $K_0=K_0(n,m,\alpha,\mu_0,M_0)>0$ that satisfy the following:

Let $\Mf$ be a $C^{1,1}$-manifold and let $p_0\in\Mf$. Let $X_1,\dots,X_m$ be log-Lipschitz real vector fields on $\Mf$. Suppose $n:=\rank\Span(X_1(p_0),\dots,X_m(p_0))$ and by re-ordering the indices we have  the following:
\begin{itemize}[parsep=-0.3ex]
    \item $|X_{i_1}\wedge\dots\wedge X_{i_n}(p_0)|\le M_0|X_1\wedge\dots\wedge X_n(p_0)|$ for all $1\le i_1,\dots,i_n\le m$.
    \item The canonical map $t\mapsto e^{t^1X_1+\dots+t^nX_n}(p_0)$ is defined for $t\in B^n(0,\mu_0)$. 
    \item For every $0<r\le\mu_0$ and every $p$ such that $e^{t\cdot X}(p)$ is defined for $t\in B^n(0,r)$, we have $e^{t\cdot X}(p)\neq p$ for all $t\in B^n(0,r)\backslash\{0\}$ and $p\in B_X(p_0,r)$.
    \item We have $(c_{ij}^k)_{\substack{1\le i,j\le m\\1\le k\le n}}$ such that
    \begin{equation}\label{Eqn::Quant::AssumptionM0}
        [X_i,X_j]=\sum_{k=1}^nc_{ij}^kX_k\quad\text{with }\sup_{i,j,k}\|c_{ij}^k\|_{C^{0,\alpha-1}_X(\Mf)}<M_0.
    \end{equation}
\end{itemize}

Then for any $p_0\in \Mf$ there is a map $\Phi:B^n(0,1)\to\Mf$ such that
\begin{itemize}[parsep=-0.3ex]
    \item $\Phi(0)=p_0$ and $\Phi$ is $\Co^{2-}$-diffeomorphism onto its image.
    \item $\Phi^*X'=[\Phi^*X_1,\dots,\Phi^*X_n]^\top$ and $\Phi^*X''=[\Phi^*X_{n+1},\dots,\Phi^*X_m]^\top$ are collections of $\Co^\alpha$-vector fields on $B^n(0,1)$ that can be written as
    \begin{equation}\label{Eqn::Quant::Khat}
        \Phi^*X'=\widehat K (I+A')\nabla,\quad \Phi^*X''=\widehat KA''\nabla.
    \end{equation}
    Where $A=\begin{pmatrix}A'\\A''\end{pmatrix}\in\Co^\alpha(\B^n;\R^{m\times n})$ satisfy $A'(0)=0$, $\|A'\|_{\Co^\alpha(\B^n;\R^{n\times n})}\le\tfrac12$ and $\|A''\|_{\Co^\alpha(\B^n;\R^{(m-n)\times n})}\le K_0$.
\end{itemize}
\end{mainthm}

Here $C^{0,\alpha-1}_X(\Mf)\subset C^0_\loc(\Mf)$ is the H\"older space along with vector fields $X=(X_1,\dots,X_m)$, whose norm is given by the following: recall \eqref{Eqn::Intro::CCDist} for the Carnot-Carath\'eodory distance,
\begin{equation}\label{Eqn::Intro::CsXNorm}
   \|f\|_{C^{0,\alpha-1}_X(\Mf)}:=\|f\|_{C^0}+\sup\limits_{x,y\in \Mf;x\neq y}\frac{|f(x)-f(y)|}{\dist_X(x,y)^{\alpha-1}}.
\end{equation}

The results for $\alpha\ge2$ is also true, with some modifications on the spaces $C^{0,\alpha-1}_X$. The case $\alpha>2$ is already known in \cite[Theorem 2.14]{CoordAdaptedII}. There is also an analogy involving complex vector fields, see \cite{StreetSubHermitian}.

\section{Notations and Organizations}\label{Section::Convention}

We use \textbf{distributions} for generalized functions, rather than tangential subbundles on manifolds.

Given a bijection $f:U\to V$, we use $f^\Inv:V\to U$ as the inverse function of $f$. We do not use the notion $f^{-1}$ in order to reduce confusion when $(\cdot)^{-1}$ is used for inverting a matrix.

We use $A\Subset B$ to denote that $A$ is a relatively compact subset of $B$.

By $A \lesssim B$ we mean that $A \leq CB$ where $C$ is a constant independent of $A,B$. We use $A \approx B$ for ``$A \lesssim B$ and $B \lesssim A$''. And we use $A\lesssim_xB$ to emphasize the constant dependence on quantity $x$.

We use $\B^n=B^n(0,1)$ as the unit ball in $\R^n$ and $r\B^n=B^n(0,r)$ for $r>0$. We use a complex cone $\Hb^n:=\{x+iy\in\C^n:x\in\B^n,4|y|<1-|x|\}\subseteq\C^n$.

For functions $f:U\to\R$, $g:V\to\R$, we use the tensor notation $f\otimes g:U\times V\to\R$ as $(f\otimes g)(x,y):=f(x)g(y)$.

\medskip
We denote by $\Co^\alpha$ the class of H\"older-Zygmund functions:
\begin{defn}\label{Defn::Intro::DefofHold}
    Let $\alpha\in(0,\infty]$, let $\Omega\subset\R^n$ be a convex open subset and let $\Xs$ be a Banach space. The space $\Co^\alpha(\Omega;\Xs)$ of bounded ($\Xs$ vector-valued) H\"older-Zygmund functions on $\Omega$, is given recursively as follows:
    \begin{enumerate}[parsep=-0.3ex,label=(\roman*)]
    \item\label{Item::Intro::DefofHold::<1} For $0<\alpha<1$: $\Co^\alpha(\Omega;\Xs)$ consists of those continuous functions $f:\Omega\to \Xs$ such that:
    \begin{equation*}
         \textstyle\|f\|_{\Co^\alpha(\Omega;\Xs)}:=\sup_{x\in \Omega} |f(x)|_\Xs + \sup_{\substack{x,y\in\Omega;x\neq y}} |x-y|^{-\alpha}|f(x)-f(y)|_\Xs<\infty.
    \end{equation*}
    \item\label{Item::Intro::DefofHold::=1} For $\alpha=1$: $\Co^1(\Omega;\Xs)$ consists of those continuous functions $f:\Omega\to \Xs$ such that:
    \begin{equation*}
         \textstyle\|f\|_{\Co^1(\Omega;\Xs)}:=\sup_{x\in \Omega} |f(x)|_\Xs + \sup_{\substack{x,y\in\Omega;x\neq y}} |x-y|^{-1}|\frac{f(x)+f(y)}2-f(\frac{x+y}2)|_\Xs<\infty.
    \end{equation*}
    \item\label{Item::Hold::DefofHold::>1} For $\alpha>1$: $\Co^\alpha(\Omega;\Xs)$ consists of $f\in\Co^{\alpha-1}(\Omega;\Xs)$ such that the following norm is finite: $$\textstyle\|f\|_{\Co^\alpha(\Omega;\Xs)}:=\|f\|_{\Co^{\alpha-1}(\Omega;\Xs)}+\sum_{j=1}^n\|\partial_jf\|_{\Co^{\alpha-1}(\Omega;\Xs)}.$$
    \item\label{Item::Intro::DefofHold::-} We use $\Co^\infty(\Omega;\Xs):=\bigcap_{\beta>0}\Co^\beta(\Omega;\Xs)$ and $\Co^{\alpha-}(\Omega;\Xs):=\bigcap_{\beta<\alpha}\Co^\beta(\Omega;\Xs)$, both endowed with the standard limit Fr\'echet topologies. That is, we say $f_k\to f_0$ in $\Co^{\alpha-}$ if $f_k\to f_0$ in $\Co^\beta$ for all $\beta<\alpha$. See also Remark \ref{Rmk::Hold::Topology}.
    \item\label{Item::Intro::DefofHold::+} We use $\Co^{\alpha+}(\Omega;\Xs):=\bigcup_{\beta>\alpha}\Co^\beta(\Omega;\Xs)$, whose topologies are not used in the thesis.
    \item\label{Item::Intro::DefofHold::loc} We use $\Co^\alpha_\loc(\Omega;\Xs)$, $\Co^{\alpha-}_\loc(\Omega;\Xs)$ and $\Co^{\alpha+}_\loc(\Omega;\Xs)$ to be the spaces of all $f\in C^0_\loc(\Omega;\Xs)$ such that for every precompact open $\Omega'\Subset\Omega$, we have $f\in\Co^\alpha(\Omega';\Xs)$, $f\in\Co^{\alpha-}(\Omega';\Xs)$ and $f\in\Co^{\alpha+}(\Omega';\Xs)$ respectively. (See Remark \ref{Rmk::Hold::Topology}.)
    \end{enumerate}

    We use $\Co^\beta(\Omega)=\Co^\beta(\Omega;\R)$ and $\Co^\beta_\loc(\Omega)=\Co^\beta_\loc(\Omega;\R)$ for $\beta\in\{\alpha,\alpha-,\alpha+:\alpha>0\}\cup\{\infty\}$.
    
    We use $\|f\|_{\Co^\alpha(\Omega)}$ either when $f:\Omega\to\R$ is a function or when the codomain of $f$ is clear. And we use $\|f\|_{\Co^\alpha}$ when both the domain and the codomain are clear.
\end{defn}
We define $\Co^\alpha(\Omega;\Xs)$ for $\alpha\le0$ in Definition \ref{Defn::Hold::NegHold} using the Besov space $\Bs_{\infty\infty}^\alpha$.

Let $U\subseteq\R^m_x$ and $V\subseteq\R^n_s$ be two open subsets. For $\alpha,\beta\in(0,\infty]$, we denote by $\Co^{\alpha,\beta}(U,V)$ the space of functions on $U\times V$ that is bounded $\Co^\alpha$ in $x$ and bounded $\Co^\beta$ in $s$ separately, see Definition \ref{Defn::Hold::MixHold}. For $\alpha,\beta\in(-\infty,\infty]$ we denote by $\Co^\alpha\Co^\beta(U,V)$ the space of functions on $U\times V$ that is bounded $\Co^\alpha$ in $x$ and bounded $\Co^\beta$ in $s$ simultaneously, see Definition \ref{Defn::Hold::BiHold}. In most of the cases we work on the spaces like $\Co^\alpha L^\infty\cap\Co^\beta\Co^\gamma(U,V)$.

We use $\Co^{\LogL}(\Omega)$ and $\Co^\Lip(\Omega)$ as the spaces of functions on $\Omega$ with are bounded log-Lipschitz and bounded Lipschitz respectively. For $k\ge1$, we use $\Co^{k+\LogL}(\Omega)$ and $\Co^{k+\Lip}(\Omega)$ as the spaces of function whose first $k$-derivatives are log-Lipschitz/ Lipschitz. See Definition \ref{Defn::Hold::CLogLandLip}.

\medskip

For the index of regularity, we use the following extended index set from Definition \ref{Defn::Hold::ExtIndex}:
\begin{equation*}
    \R_\Eb:=\{\alpha,\alpha-:\alpha\in\R\}\cup\{k+\LogL,k+\Lip:k=-1,1,2,\dots\}\cup\{\infty\}.
\end{equation*}
For $\sigma_1,\sigma_2\in\R_\Eb$, we define order $\sigma_1\le\sigma_2$ if $\Co^{\sigma_1}\supseteq\Co^{\sigma_2}$.

We denote by $\R_\Eb^+:=\{\sigma\in\R_\Eb:\sigma>\Lip-1\}$ the set of positive generalized indices. Equivalently $\alpha\in\R_\Eb^+$ if and only if $\alpha\in\R_\Eb$ and $\Co^\alpha\subset C^0$.

We denote by $\R^-:=\{\alpha-:\alpha\in\R\}$ and $\R_+^-:=\{\alpha-:\alpha\in\R_+\}$.

We use the convention $\infty=\infty+1=\infty-1=\infty+=\infty-=\frac\infty2$.

For $\alpha,\beta\in\R\cup\{\infty\}$, we use the formal operations $(\alpha-)+\beta=\alpha+(\beta-)=(\alpha-)+(\beta-):=(\alpha+\beta)-$; and if $\alpha,\beta\in\R_+\cup\{\infty\}$, we use $(\alpha-)\cdot\beta=\alpha\cdot(\beta-)=(\alpha-)\cdot(\beta-):=(\alpha\beta)-$. See Convention \ref{Conv::Hold::ExtIndOp}.

For $\alpha,\beta\in\R_\Eb$ such that $\alpha,\beta>0$, we use $\alpha\circ\beta$ to be the unique positive index in $\R_\Eb$ such that $f\circ g\in\Co^{\alpha\circ\beta}_\loc(\R^n)$ for every $m,n\ge1$, $f\in\Co^\alpha_\loc(\R^m)$ and $g\in\Co^\beta_\loc(\R^n;\R^m)$. See Definition \ref{Defn::Hold::CompIndex} and Corollary \ref{Cor::Hold::CompOp}.

In Section \ref{Section::EllipticPara} we use notations $\beta^{\sim\alpha}$, $\beta^{\wedge\alpha}$ for $\alpha\in(\frac12,\infty)$ and $\beta\in\R_+\cup\R_+^-$, that labels the H\"older exponents for the coordinate charts and the coordinate vector fields (see  \eqref{Eqn::EllipticPara::Betas}): for $\gamma\in\R_+$,
\begin{align*}
    \gamma^{\sim\alpha}:=\min\left(\gamma,\alpha+1\right),\quad \gamma^{\wedge\alpha}:=\min\left(\gamma-,\alpha+1\right),&\quad\text{for }\alpha\in(1,\infty);\\
    \gamma^{\sim\alpha}=\gamma^{\wedge\alpha}:=\min\left(\gamma-,(2-\tfrac1\alpha)\gamma,\alpha+1\right),&\quad\text{for }\alpha\in(\tfrac12,1];
    \\
    (\gamma-)^{\sim\alpha}=(\gamma-)^{\wedge\alpha}:=\min\left(\gamma-,(2-\tfrac1\alpha)\gamma-,\alpha+1\right),&\quad\text{for }\alpha\in(\tfrac12,\infty).
\end{align*}

\bigskip

For the matrix space $\R^{m\times n}$, we use the standard $\ell^2$-matrix norm, that is the operator norm of the linear map $(\R^n,\ell^2)\to(\R^m,\ell^2)$:
\begin{equation}\label{Eqn::Intro::MatrixNorms}
    |A|_{\R^{m\times n}}=|A|_{\ell^2}:=\sup\limits_{v\in\R^n\backslash\{0\}}|A\cdot v|_{\R^m}/|v|_{\R^n}.
\end{equation}

Note that for $A\in\R^{m\times n}$ and $B\in\R^{n\times p}$ we have $|AB|\le |A||B|$. In particular $\R^{n\times n}$ is a Banach algebra.

\medskip
Given a coordinate chart $x=(x^1,\dots,x^n)$, we write  $dx:=[dx^1,\dots,dx^n]$ as a row vector and $\Coorvec{x}:=\begin{bmatrix}\partial_{x^1}\\\vdots\\\partial_{x^n}\end{bmatrix}$ as a column vector. We may use $X$ sometimes as a single vector field, and sometimes as a (column) collection of vector fields $X=[X_1,\dots,X_r]^\top$.  We use these conventions when there is no ambiguity.

In a mixed real and complex domain $\Omega\subseteq\R^r\times\C^m$ with standard (real and complex) coordinate system $(t,z)=(t^1,\dots,t^r,z^1,\dots,z^m)$, 
for a function $f$ on $\Omega$, we use $\nabla f,\nabla_zf,\partial_zf$ as functions taking values in column vectors, which have rows $r+2m$, $2m$ and $m$ respectively as follows,
\begin{equation}\label{Eqn::Intro::ColumnNote}
    \begin{gathered}\textstyle
    \nabla f=\nabla_{t,z}f=\partial_{t,z,\bar z}f=\big[\frac{\partial f}{\partial t^1},\dots,\frac{\partial f}{\partial t^r},\frac{\partial f}{\partial z^1},\dots,\frac{\partial f}{\partial z^m},\frac{\partial f}{\partial \bar z^1},\dots,\frac{\partial f}{\partial \bar z^m}\big]^\top,
\\\textstyle
\nabla_zf=\partial_{z,\bar z}f=\big[\frac{\partial f}{\partial z^1},\dots,\frac{\partial f}{\partial z^m},\frac{\partial f}{\partial \bar z^1},\dots,\frac{\partial f}{\partial \bar z^m}\big]^\top,
\qquad \partial_zf=\frac{\partial f}{\partial z}=\big[\frac{\partial f}{\partial z^1},\dots,\frac{\partial f}{\partial z^m}\big]^\top.
\end{gathered}
\end{equation}

Notice that $df$ is a row vector while $\nabla f$ is a column vector.

For a map $F:\Mf\to\R^r\times\C^m$, we usually use the notation $F=(F',F'')$ where $F':\Mf\to\R^r$ and $F'':\Mf\to\C^m$ are the corresponding components of $F$.

On a manifold, we refer a \textbf{topological coordinate chart} to a map $F:U\subseteq\Mf\to\R^n$ which is homeomorphism onto its image. We refer \textbf{topological parameterization} to a map $\Phi:\Omega\subseteq\R^n\to\Mf$ which is homeomorphism onto its image. For $\alpha>1$, we say that $F$ is a $\Co^\alpha$ coordinate chart, equivalently that $\Phi$ is a $\Co^\alpha$ regular parameterization, if they and their inverses are both $\Co^\alpha$.

\bigskip\bigskip
The paper is organized as follows: the proof of theorems are mostly separated into four parts:
\begin{itemize}
    \item The real Frobenius-type theorems: Theorems \ref{MainThm::LogFro} and \ref{MainThm::SingFro}. The proofs are done in Section \ref{Section::RealFro}. They require Sections \ref{Section::AppThm} and \ref{Section::FlowComm}.
    \item The ``pure complex'' part of the complex Frobenius theorems: Theorem \ref{KeyThm::EllipticPara} in Section \ref{Section::EllipticPara}, which is the key part to the proof of Theorems \ref{MainThm::CpxFro}, \ref{MainThm::RoughFro1} and \ref{MainThm::RoughFro2}. The proof requires Sections \ref{Section::BiHoldSec}, \ref{Section::HoloFro} and \ref{Section::SecHolLap}. 
    \item The quantitative Frobenius theorem: Theorem \ref{MainThm::QuantFro}. The proof is done in Section \ref{Section::PfQuantFro}. It requires Sections \ref{Section::CanonicalCoordinates}. The key is to prove Theorem \ref{KeyThm::Rough1Form}, which is the  Section \ref{Section::Rough1Form}.
    \item The examples: Theorems \ref{MainThm::Sharpddz} and \ref{MainThm::CounterNN}. The proofs are mostly self-contained.
\end{itemize}

Once the real part and the pure complex part are done, we get the proof of Theorems \ref{MainThm::CpxFro}, \ref{MainThm::RoughFro1} and \ref{MainThm::RoughFro2} in Section \ref{Section::PfCpxFro}. The proof is essentially the combinations of Theorem \ref{MainThm::LogFro} (see also Theorems \ref{KeyThm::RealFro::ImprovedLogFro} and \ref{KeyThm::RealFro::BLFro}) and Theorem \ref{KeyThm::EllipticPara}, though tracing the indices of composition maps, which are done in Section \ref{Section::SingHoldComp}. See also Theorems \ref{KeyThm::CpxFro1} and \ref{KeyThm::CpxFro2} for more general versions of Theorems \ref{MainThm::CpxFro}, \ref{MainThm::RoughFro1} and \ref{MainThm::RoughFro2}.

The Sections \ref{Section::HoldSec} and \ref{Section::BiHoldSec} give discussions of H\"older-Zygmund spaces in the bi-parameter setting and set up necessary technics, which are the based stones of the thesis. A major part of the results come from \cite[Section 2]{YaoCpxFro}. Here in Section \ref{Section::BiHoldMult} the Proposition \ref{Prop::Hold::MultLow} is not appeared in arXiv before.


In Section \ref{Section::AppThm} we prove Theorem \ref{Thm::Hold::ApproxThm}, a result from \cite[Section 2.2]{YaoLLFro} that show good approximation of the product functions if their product is identically zero.
We then use Theorem \ref{Thm::Hold::ApproxThm} to Section \ref{Section::FlowComm} and show the flow commutativity between two commutative log-Lipschitz vector fields, see Proposition \ref{Prop::ODE::StrFlowComm}. The flow commuting is the key to the proof of Theorems \ref{MainThm::LogFro} and \ref{MainThm::SingFro}. We prove Theorem \ref{MainThm::LogFro} in Section \ref{Section::RealFro::PfLogFro} and Theorem \ref{MainThm::SingFro} in Section \ref{Section::RealFro::PfSingFro}.
In order to prove Theorems \ref{MainThm::CpxFro}, \ref{MainThm::RoughFro1} and \ref{MainThm::RoughFro2}, we generalize Theorem \ref{MainThm::LogFro} into the bi-layer case: the Theorem \ref{KeyThm::RealFro::BLFro} in Section \ref{Section::RealFro::BLFro}.

Most results in Section \ref{Section::FlowComm} (the part after Definition \ref{Defn::ODE::SubUnitCurve}) are not appeared in arXiv before.

In Section \ref{Section::CanonicalCoordinates} we give some basic properties of the canonical coordinates of the second kind. The proof of these results are not appeared in arXiv before but they are mostly inspired from \cite{CoordAdapted}. This section is only used in the proof of  Theorem \ref{MainThm::QuantFro} in Section \ref{Section::PfQuantFro}.

Section \ref{Section::SecHolLap} follows from \cite[Appendix B]{YaoCpxFro}, where we construct an inverse Laplacian that can be holomophically extended to some complex domain. The proof is inspired by \cite{Analyticity}. The main result in this section, Proposition \ref{Prop::HolLap}, is only used in Proposition \ref{Prop::EllipticPara::AnalyticPDE} in Section \ref{Section::EllipticPara::AnalPDE}.

Chapter \ref{Chapter::DisInv} discusses some distributional characterizations for sections and involutivity, along with some discussion on rough integrability and singular subbundles. The results in Sections \ref{Section::DisInv::CharSec} and \ref{Section::DisInv::CharInv} come from \cite[Section 5]{YaoLLFro}. This part is independent of the later sections and the proof of the main theorems. Also the Sections \ref{Section::DisInv::Integrable} and \ref{Section::DisInv::Sing} are not appeared in arXiv before.

Chapter \ref{Chapter::Malgrange} gives two theorems using Malgrange's methods of factorizing coordinate changes: Theorem \ref{KeyThm::EllipticPara} in Section \ref{Section::EllipticPara} and Theorem \ref{KeyThm::Rough1Form} in Section \ref{Section::Rough1Form}. See overview of their proofs in Section \ref{Section::EllipticPara::Overview} and Section \ref{Section::Rough1FormOV}, respectively. 
Both proofs consist of finding an intermediate coordinate change using some existence theorem of elliptic PDEs, and showing the objects in the new coordinate chart have improved regularity using some interior regularity theorem of nonlinear elliptic PDEs. In Theorem \ref{KeyThm::EllipticPara}, to deal with the case $\alpha=1$ or when the subbundle is $\Co^{\alpha,\beta-}$, we need an extra result of Schauder's estimate in Section \ref{Section::EllipticPara::Scaling}, which is not appeared in arXiv before.

In Section \ref{Section::RealFro} we prove the real Frobenius type theorems. We prove Theorem \ref{MainThm::LogFro} in Section \ref{Section::RealFro::PfLogFro}, Theorem \ref{MainThm::SingFro} in Section \ref{Section::RealFro::PfSingFro} and a generalization of Theorem \ref{MainThm::LogFro} in Section \ref{Section::RealFro::BLFro}. Section \ref{Section::RealFro::PDEPartLogFro} comes from \cite[Section 6.1]{YaoLLFro}, which gives a PDE interpretation of Theorem \ref{MainThm::LogFro} and is independent of the other sections.

In Section \ref{Section::CpxFro} we give a general version of complex Frobenius theorem, the Theorems \ref{KeyThm::CpxFro1} and \ref{KeyThm::CpxFro2}. We give the proof using by combining Theorem \ref{KeyThm::EllipticPara} in Section \ref{Section::EllipticPara} and Theorem \ref{KeyThm::RealFro::BLFro} in Section \ref{Section::RealFro::BLFro}.

Section \ref{Section::Examples} lists several examples which indicate the sharpness of some H\"older exponents in Theorems \ref{MainThm::LogFro} and \ref{MainThm::CpxFro}. We prove Theorem \ref{MainThm::CounterNN} in Section \ref{Section::CountNN} and Theorem \ref{MainThm::Sharpddz} in Section \ref{Section::Sharpddz}. Section \ref{Section::CanCoor-1Reg} discusses the regularity loss of the canonical coordinates, which is independent of the other subsections in Section \ref{Section::Examples}.

In Section \ref{Section::QuantFro} we prove Theorem \ref{MainThm::QuantFro}. The proof follows from \cite[Section 9]{StreetYaoVectorFields}.

\chapter{The Tools And The Necessities}
\section{On Single parameter H\"older-Zygmund Spaces}\label{Section::HoldSec}

In this section we recall some properties of H\"older-Zygmund spaces. The results in this section are mostly for technical purposes.

\subsection{Basic characterizations and paraproducts}\label{Section::BasicHold}
Recall the classical (single-parameter) H\"older-Zygmund spaces in Definition \ref{Defn::Intro::DefofHold}. We first recall some results from functions and distributions that take value in Banach spaces.

\begin{defn}
Let $\Omega\subseteq\R^n$ be an open subset and let $\Xs$ be a complete locally convex topological vector space. The space $\D'(U;\Xs)$ of $\Xs$-valued distributions on $U$, is the set of all continuous linear map $F:C_c^\infty(U)\to\Xs$. The space $\Sc'(\R^n;\Xs)$ of $\Xs$-valued tempered distributions,  is the set of all continuous linear map $F:\Sc(\R^n)\to\Xs$.
\end{defn}

Let $U\subseteq\R^n$ and $V\subseteq\R^q$ be two open subsets. By the Schwartz Kernel Theorem (see for example \cite[Corollary 51.6 and Theorem 51.7]{TrevesTopologicalSpaces}), we have canonical correspondence $\D'(U;\D'(V))\cong\D'(U\times V)$ and $\Sc'(\R^n;\Sc'(\R^q))\cong\Sc'(\R^n\times\R^q)$. Thus in applications we can limit our focus on distributions defined on the product domain $U\times V$.


On the space $\Sc'(\R^n;\Xs)$ we can talk about the Littlewood-Paley decomposition of $\Xs$-valued functions
\begin{defn}\label{Defn::Hold::DyadicResolution}
A \textit{(Fourier) dyadic resolution of unity} for  $\R^n$ is a sequence of Schwartz functions $\phi=(\phi_j)_{j=0}^\infty$ on $\R^n$ that are all real-valued even functions, such that
\begin{itemize}[nolistsep]
    \item The Fourier transform $\hat\phi_0(\xi)=\int_{\R^n}\phi(x)e^{-2\pi ix\xi}dx$ satisfies $\supp\hat\phi_0\subset\{|\xi|<2\}$ and $\hat\phi_0\big|_{\{|\xi|<1\}}\equiv1$.
    \item For $j\ge1$, $\phi_j(x):=2^{nj}\phi_0(2^jx)-2^{n(j-1)}\phi_0(2^{j-1}x)$ for $x\in\R^n$.
\end{itemize}
We define the Littlewood-Paley block operators $\{\Su_j,\De_j\}_{j=0}^\infty=\{\Su_j^\phi,\De_j^\phi\}_{j=0}^\infty$ associated to $\phi$, as
\begin{equation}\label{Eqn::Hold::DefBlockOP}
    \De_jf=\De_j^\phi:=\phi_j\ast f,\qquad \Su_jf=\Su_j^\phi:=\sum_{k=0}^j\De_kf=2^{jn}\phi_0(2^j\cdot)\ast f,\quad j\ge0,\ f\in\Sc'(\R^n).
\end{equation}
\end{defn}
Immediately we see that
\begin{equation}\label{Eqn::Hold::RmkDyaSupp}
    \supp\phi_0\subset\{|\xi|<2\},\quad\supp\phi_j\subset\{2^{j-1}<|\xi|<2^{j+1}\},\quad j\ge1.
\end{equation}
In other words,
\begin{equation}\label{Eqn::Hold::LPCharComp}
    \De_j\De_k=0,\quad\text{if }|j-k|\ge2;\qquad\De_j\Su_k=\begin{cases}\Su_k,&\text{if }j\le k-2,\\0,&\text{if }j\ge k+2.\end{cases}
\end{equation}

Let $\lambda\in\Sc(\R^n)$ be a Schwartz function and let $F:\Sc(\R^n)\to\Xs$ be a $\Xs$-valued tempered distribution. The convolution $\lambda\ast F:\Sc(\R^n)\to\Xs$ is defined to be the same as the one for numerical  tempered distribution: for $u\in\Sc(\R^n)$,
\begin{equation*}
    \textstyle(\lambda\ast F)(u):=F(\widetilde\lambda\ast u)\in\Xs,\quad\text{where }\widetilde\lambda(x)=\lambda(-x).
\end{equation*}

One can see that when $\Xs$ is a Banach and $F\in L^1(\R^n;\Xs)$ is a strongly measurable and Bochner integrable function, we have $\lambda\ast F\in L^1(\R^n;\Xs)$ and $\lambda\ast F(x)=\int_{\R^n}\lambda(x-y)F(y)dy$ with the right hand side be a convergence Bochner integral. See \cite[Section 2.4.d]{AnalBana}.

\begin{lem}\label{Lem::Hold::OperatorExtension}
Let $\alpha\in(0,\infty)$, let $\Omega\subseteq\R^n$ be an open set.
Let $T:\Xs\to\Ys$ be a bounded linear operator between two Banach spaces. Then $(\tilde Tf)(x):=T(f(x))$ define a bounded linear map $\tilde T:\Co^\alpha(\Omega;\Xs)\to\Co^\alpha(\Omega;\Ys)$. Moreover $\|\tilde T\|_{\Co^\alpha(\Xs)\to\Co^\alpha(\Ys)}=\|T\|_{\Xs\to\Ys}$ has the same operator norm.
\end{lem}
\begin{proof}
We only prove the case $\alpha\in(0,1)$ since the proof of other cases follow from the same argument. Clearly $\|\tilde T\|\ge\|T\|$. By Definition \ref{Defn::Intro::DefofHold} for $f\in\Co^\alpha(\Omega;\Xs)$ we have
\begin{align*}
    \|\tilde Tf\|_{\Co^\alpha(\Omega;\Ys)}&=\sup_{x\in \Omega}|T(f(x))|_\Ys+\sup_{x,y\in \Omega;x\neq y}|T(f(x))-T(f(y))|_\Ys\cdot|x-y|^\alpha
    \\
    &\le \|T\|_{\Xs\to\Ys}\sup_{x\in \Omega}|f(x)|_\Xs+\|T\|_{\Xs\to\Ys}\sup_{x,y\in \Omega;x\neq y}|f(x)-f(y)|_\Xs\cdot|x-y|^\alpha=\|T\|_{\Xs\to\Ys}\|f\|_{\Co^\alpha(\Omega;\Xs)}.
\end{align*}
Thus $\tilde T$ is bounded with the same operator norm.
\end{proof}

\begin{lem}\label{Lem::Hold::HoldChar}
    Let $\alpha>0$ and let $\phi=(\phi_j)_{j=0}^\infty$ be a dyadic resolution. Let $\Xs$ be an arbitrary Banach space.
    \begin{enumerate}[parsep=-0.3ex,label=(\roman*)]
        \item\label{Item::Hold::HoldChar::LPHoldChar} $f\in\Co^\alpha(\R^n;\Xs)$ if and only if $\phi_j\ast f\in C^0(\R^n;\Xs)$ for every $j\ge0$ and $\sup_{j\ge0;x\in \Xs}2^{j\alpha}|\phi_j\ast f(x)|_\Xs<\infty$. Moreover there is a $C=C(\phi,\alpha)>0$ that does not depend on $f$ such that
        \begin{equation}\label{Eqn::Hold::HoldChar::LPHoldChar}
            \textstyle C^{-1}\|f\|_{\Co^\alpha(\R^n;\Xs)}\le\sup_{j\ge0;x\in \Xs}2^{j\alpha}|\phi_j\ast f(x)|_\Xs\le C\|f\|_{\Co^\alpha(\R^n;\Xs)},\quad\forall f\in\Co^\alpha(\R^n).
        \end{equation}
        \item\label{Item::Hold::HoldChar::DomainChar} Let $\Omega\subset\R^n$ be a bounded Lipschitz domain. Then $\Co^\alpha(\Omega;\Xs)$ has an equivalent norm
        \begin{equation*}
            f\mapsto\inf\{\|\tilde f\|_{\Co^\alpha(\R^n;\Xs)}:\tilde f\in\Co^\alpha(\R^n;\Xs),\tilde f\big|_\Omega=f\}.
        \end{equation*}
        \item\label{Item::Hold::HoldChar::02} Let $\Omega\subseteq\R^n$ be either the whole space or a bounded smooth domain. Suppose $0<\alpha<2$, then $\Co^\alpha(\Omega;\Xs)$ has an equivalent norm
        \begin{equation}\label{Eqn::Hold::HoldChar::02Norm}
            \textstyle f\mapsto\sup_{x\in\Omega}|f(x)|_\Xs+\sup_{x,y\in\Omega;x\neq y}\big|\frac{f(x)+f(y)}2-f\big(\frac{x+y}2\big)\big|_\Xs\cdot|x-y|^{-\alpha}.
        \end{equation}
    \end{enumerate}
\end{lem}
See \cite[Chapters 2.2.2, 2.5.7 and 3.4.2]{Triebel1} for the proof of $\Xs=\R$. For general Banach space $\Xs$, the convolution $\phi_j\ast f(x)=\int\phi_j(x-y)f(y)dy$ is defined via Bochner integral (see for example \cite[Chapter 2]{AnalBana}).
\begin{proof}[Proof of Lemma \ref{Lem::Hold::HoldChar} \ref{Item::Hold::HoldChar::LPHoldChar} and \ref{Item::Hold::HoldChar::02} with $\Omega=\R^n$]
For convenience we let $\|f\|_{\Co^\alpha}'$ to be \eqref{Eqn::Hold::HoldChar::02Norm} and $\|f\|_{\alpha;\phi}:=\sup_{j\ge0;x\in \Xs}2^{j\alpha}|\phi_j\ast f(x)|_\Xs$.

    Since $\phi_j$ are even functions and satisfies $\int\phi_j=0$ for $j\ge0$, for $x\in\R^n$,
    \begin{align}
    \label{Eqn::Hold::HoldChar:Pf::Proof01}&\phi_j\ast f(x)=\int f(x-y)\phi_j(y)dy=\int (f(x-y)-f(x))\phi_j(y)dy\\
    \label{Eqn::Hold::HoldChar:Pf::Proof02}=&\frac12\big(\int (f(x-y)-f(x))\phi_j(y)dy+\int (f(x+y)-f(x))\phi_j(-y)dy\big)=\int\big(\frac{f(x+y)+f(x-y)}2-f(x)\big)\phi_j(y)dy.
\end{align}

    Recall for $\alpha\in(0,1)$, $\|f\|_{\Co^\alpha}=\|f\|_{C^0}+\sup_{x,h\in\R^n}|h|^{-\alpha}|f(x)-f(x-h)|_\Xs$.   
    Clearly $\|\phi_0\ast f\|_{L^\infty}\le\|\phi_0\|_{L^1}\|f\|_{C^0}\le\|\phi_0\|_{L^1}\|f\|_{\Co^\alpha}$. By \eqref{Eqn::Hold::HoldChar:Pf::Proof01}, for $j\ge 1$,
\begin{equation}\label{Eqn::Hold::HoldChar:Pf::Bs<C^s01}
\begin{aligned}
    \|\phi_j\ast f\|_{L^\infty}&\textstyle\le \int|f(x-y)-f(x)|_\Xs|\phi_j(y)|dy\le \|f\|_{\Co^\alpha}\int|y|^\alpha|\phi_j(y)|dy
    \\
    &\textstyle=2^{(1-j)\alpha}\|f\|_{\Co^\alpha}\int|y|^\alpha|\phi_1(y)|dy\approx 2^{-j\alpha}\|f\|_{\Co^\alpha}.
\end{aligned}
\end{equation}
Thus $\|f\|_{\alpha;\phi}\lesssim\|f\|_{\Co^\alpha}$.

To see $\|f\|_{\Co^\alpha}\lesssim\|f\|_{\alpha;\phi}$, we denote $\mu_j=\sum_{k=j-1}^{j+1}\phi_j$ with the convention $\phi_{-1}=0$. So $\phi_j\ast f=\mu_j\ast\phi_j\ast f$ for all $j\ge0$.
Clearly $\|f\|_{C^0}\le\sum_{j=0}^\infty\|\phi_j\ast f\|_{L^\infty}\le\sum_{j=0}^\infty 2^{-j\alpha}\|f\|_{\alpha;\phi}\lesssim\|f\|_{\alpha;\phi}$. Let $k\in\Z_+$ and let $h\in\R^n$ satisfies $2^{-k-\frac12}\le |h|\le 2^{-k+\frac12}$, what we need is to show $|f(x+h)-f(x)|\le 2^{-k\alpha}$.

Write $f=f_k^0+f_k^\infty$ where 
\begin{equation}\label{Eqn::Hold::HoldChar:Pf::Tmp}
    f_k^0:=\sum_{j=0}^k\phi_j\ast f,\quad f_k^\infty:=\sum_{j=k+1}^\infty\phi_j\ast f.
\end{equation}

For $f_k^\infty$, clearly $|f_k^\infty(x+h)-f_k^\infty(x)|\le 2\|f_k^\infty\|_{L^\infty}\le\sum_{j=k+1}^\infty\|\phi_j\ast f\|_{L^\infty}\lesssim 2^{-k\alpha}\|f\|_{\alpha;\phi}$. For $f_k^0$, we use $f_k^0(x+h)-f_k^0(x)=\int_0^1h\cdot(\nabla f_k^0)(x+th)dt$, so
\begin{align*}
    &|f_k^0(x+h)-f_k^0(x)|\le\int_0^1|h||(\nabla f_k^0)(x+th)|dt\le 2^{\frac12-k}\|\nabla f_k^0\|_{L^\infty}=2^{\frac12-k}\sum_{j=0}^k\|\nabla\mu_j\ast\phi_j\ast f\|_{L^\infty}
    \\
    \lesssim&2^{-k}\sum_{j=0}^k2^j\|\phi_j\ast f\|_{L^\infty}\lesssim 2^{-k}2^{(1-s)k}\|f\|_{\alpha;\phi}=2^{-k\alpha}\|f\|_{\alpha;\phi}.
\end{align*}
This completes the proof of \ref{Item::Hold::HoldChar::LPHoldChar} for $\alpha\in(0,1)$.

\noindent The proof of \ref{Item::Hold::HoldChar::02} for $\Omega=\R^n$ is similar.  Similar to \eqref{Eqn::Hold::HoldChar:Pf::Bs<C^s01}, by \eqref{Eqn::Hold::HoldChar:Pf::Proof02} we get $\|f\|_{\alpha;\phi}\lesssim\|f\|_{\Co^\alpha}'$:
\begin{equation*}
    \|\phi_j\ast f\|_{L^\infty}\le\int\big|\tfrac{f(x+y)+f(x-y)}2-f(x)\big||\phi_j(y)|dy\le\|f\|_{\Co^\alpha}'\int|y|^s|\phi_j(y)|dy\approx2^{-j\alpha}\|f\|_{\Co^\alpha}'.
\end{equation*}

To see $\|f\|_{\Co^\alpha}'\lesssim\|f\|_{\alpha;\phi}$, it suffices to show  $\big|\tfrac{f(x+h)+f(x-h)}2-f(x)\big|\lesssim 2^{-k\alpha}\|f\|_{\alpha;\phi}$ when $2^{-\frac12-k}\le|h|\le 2^{\frac12-k}$. We write $f=f_k^0+f_k^\infty$ as in \eqref{Eqn::Hold::HoldChar:Pf::Tmp}. Clearly $\big|\tfrac{f_k^\infty(x+h)+f_k^\infty(x-h)}2-f_k^\infty(x)\big|\le 2\|f_k^\infty\|_{L^\infty}\le\sum_{j=k+1}^\infty\|\phi_j\ast f\|_{L^\infty}\lesssim 2^{-k\alpha}\|f\|_{\alpha;\phi}$. For $f_k^0$, we use $$f_k^0(x+h)-2f_k^0(x)+f_k^0(x-h)=\int_{[0,1]^2}\sum_{i,j=1}^nh_ih_j(\partial_{ij}^2f_k^0)(x+(t+s-1)h)dtds.$$ Therefore
\begin{align*}
    &|f_k^0(x+h)-2f_k^0(x)+f_k^0(x-h)|_\Xs\le\int_{[0,1]^2}|h|^2|(\nabla^2 f_k^0)(x+(t+s-1)h)|_\Xs dtds\le 2^{1-2k}\|\nabla^2f_k^0\|_{L^\infty(\Xs\otimes\R^{n^2})}
    \\
    \le&2^{1-2k}\sum_{j=0}^k\|\nabla^2\mu_j\ast\phi_j\ast f\|_{L^\infty(\Xs\otimes\R^{n^2})}\lesssim2^{-2k}\sum_{j=0}^k2^{2j}\|\phi_j\ast f\|_{L^\infty(\Xs\otimes\R^{n^2})}\lesssim 2^{-2k}2^{(2-s)k}\|f\|_{\alpha;\phi}=2^{-k\alpha}\|f\|_{\alpha;\phi}.
\end{align*}
This completes the proof of \ref{Item::Hold::HoldChar::02} for $\Omega=\R^n$, in particular we get \ref{Item::Hold::HoldChar::LPHoldChar} for $\alpha=1$.

To prove the case $\alpha>1$ in \ref{Item::Hold::HoldChar::LPHoldChar}, it suffices to show $\|f\|_{\alpha;\phi}\approx\|f\|_{\alpha-1;\phi}+\|\nabla f\|_{\alpha-1;\phi}$. We recall the notation $\mu_j$ from above, which gives for every $j\ge0$,
\begin{equation}\label{Eqn::Hold::HoldChar:Pf::GradBdd}
    \|\nabla\phi_j\ast f\|_{L^\infty(\Xs\otimes\R^n)}=\|\nabla\mu_j\ast\phi_j\ast f\|_{L^\infty(\Xs\otimes\R^n)}\le\|\nabla\mu_j\|_{L^1}\|\phi_j\ast f\|_{L^\infty(\Xs)}\approx2^j\|\phi_j\ast f\|_{L^\infty(\Xs)}.
\end{equation}
Thus $\|\nabla f\|_{\alpha-1;\phi}\lesssim\|f\|_{\alpha;\phi}$ and we get $\|f\|_{\alpha-1;\phi}+\|\nabla f\|_{\alpha-1;\phi}\lesssim\|f\|_{\alpha;\phi}$.

Conversely, we can use $\phi_j\ast f=(I-\Delta)(I-\Delta)^{-1}(\phi_j\ast f)$ where $\Delta=\sum_{k=1}^n\partial_{x_k}^2$:
\begin{align*}
    &\|\phi_j\ast f\|_{L^\infty(\Xs)}=\|(I-\Delta)^{-1}\mu_j\ast (I-\Delta)(\phi_j\ast f)\|_{L^\infty(\Xs)}
    \\
    =&\|(I-\Delta)^{-1}\mu_j\|_{L^1}\|\phi_j\ast f\|_{L^\infty(\Xs)}+\sum_{k=1}^n\|(I-\Delta)^{-1}\partial_k\mu_j\|_{L^1}\|\partial_k\phi_j\ast f\|_{L^\infty(\Xs)}
    \\
    \approx&2^{-2j}\|\phi_j\ast f\|_{L^\infty(\Xs)}+\sum_{k=1}^n2^{-j}\|\partial_k\phi_j\ast f\|_{L^\infty(\Xs)}\lesssim2^{-j}(\|\phi_j\ast f\|_{L^\infty(\Xs)}+\|\nabla\phi_j\ast f\|_{L^\infty(\Xs)}).
\end{align*}
This implies $\|f\|_{\alpha;\phi}\lesssim\|f\|_{\alpha-1;\phi}+\|\nabla f\|_{\alpha-1;\phi}$ and finish the proof of \ref{Item::Hold::HoldChar::LPHoldChar}.  
\end{proof}

To prove \ref{Item::Hold::HoldChar::DomainChar} and \ref{Item::Hold::HoldChar::02} for bounded smooth $\Omega$, we use the following extension operator from \cite{Triebel1}:

Let $M$ be a fixed integer and let $(a_j,b_j)_{j=1}^M$ be numbers such that $b_j>0$ and $\sum_{j=1}^Ma_j(-b_j)^k=1$ for all integer $-M/2<k<M/2$. We define the ``half-plane extension'' $E_{\mathbf H}=E_{\mathbf H;M,\{a_j,b_j\}}$ as
\begin{equation}\label{Eqn::Hold::HalfPlaneExt}
    E_{\mathbf H}f(x',x_n):=\begin{cases}f(x',x_n)&x_n>0,\\\sum_{j=1}^M a_j\cdot f(x',-b_jx_n)&x_n<0.\end{cases}
\end{equation}
Following Definition \ref{Defn::Intro::DefofHold} we see that $E_{\mathbf H}:\Co^\alpha(\{x_n>0\};\Xs)\to\Co^\alpha(\R^n;\Xs)$ is bounded for all $0<\alpha<M/2$. See also \cite[Chapter 2.9]{Triebel1}, we leave the checking and computations to reader.

We pick an open cover $U_1,\dots,U_N\subseteq\R^n$ of $\partial \Omega$, and invertible affine linear transformations  $\Phi_\nu:\B^n\xrightarrow{\sim}U_\nu$ such that $\Phi_\nu(\B^n\cap\{x_n>0\})=U_\nu\cap\Omega$. We pick $\chi_0\in C_c^\infty(\Omega)$ and $\chi_\nu\in C_c^\infty(U_\nu;[0,1])$ ($1\le \nu\le N$) to be a partition of unity for $\Omega$. Define
\begin{equation}\label{Eqn::Hold::SmoothExt}
       E_\Omega f:=\chi_0f+\sum_{\nu=1}^N\sqrt{\chi_\nu} (E_{\mathbf H}[(\sqrt{\chi_\nu} f)\circ\Phi_\nu]\circ\Phi_\nu^\Inv).
\end{equation}

From direct computations we see that
\begin{lem}[{\cite[Sections 2.9.2 and 3.3.4]{Triebel1}}]\label{Lem::Hold::SingExtLem}
For $E_\Omega$ defined in \eqref{Eqn::Hold::SmoothExt}, $E_\Omega:\Co^\alpha(\Omega;\Xs)\to\Co^\alpha(\R^n;\Xs)$ is bounded linear for $0<\alpha<\frac M2$ as well.
\end{lem}

\begin{proof}[Proof of Lemma \ref{Lem::Hold::HoldChar} \ref{Item::Hold::HoldChar::DomainChar} and \ref{Item::Hold::HoldChar::02} for $\Omega$ bounded smooth]
The direction $\|f\|_{\Co^\alpha(\R^n;\Xs)}\le\inf_{\tilde f|_\Omega=f}\|\tilde f\|_{\Co^\alpha(\R^n;\Xs)}$ is trivial. The boundedness $E_\Omega:\Co^\alpha(\Omega;\Xs)\to\Co^\alpha(\R^n;\Xs)$ gives $\inf_{\tilde f|_\Omega=f}\|\tilde f\|_{\Co^\alpha(\R^n;\Xs)}\le\|E_\Omega f\|_{\Co^\alpha(\R^n;\Xs)}\lesssim\|f\|_{\Co^\alpha(\R^n;\Xs)}$. Thus \ref{Item::Hold::HoldChar::DomainChar} holds. The same argument shows \ref{Item::Hold::HoldChar::02} for $\Omega$ bounded smooth.
\end{proof}


Following the notations in \cite[Definition 2.3.1/2(i)]{Triebel1} for $\alpha>0$ the H\"older-Zygmund spaces $\Co^\alpha(\Omega)$ are indeed the Besov spaces $\Bs_{\infty\infty}^\alpha(\Omega)$. Thus it is natural to define $\Co^\alpha$-spaces for $\alpha\le0$.

\begin{defn}\label{Defn::Hold::NegHold}Let $\alpha\le0$, let $\Xs$ be a Banach space, and let $(\phi_j)_{j=0}^\infty$ be a dyadic resolution of unity in $\R^n$. We define $\Co^\alpha:=\Bs_{\infty\infty}^\alpha$. More precisely, $\Co^\alpha(\R^n;\Xs)$ is the space of all temper distributions $f\in\Sc'(\R^n;\Xs)$ such that
\begin{equation*}
\textstyle    \|f\|_{\Co^\alpha(\R^n;\Xs)}:=\sup_{j\ge0}2^{j\alpha}\|\phi_j\ast f\|_{L^\infty(\R^n;\Xs)}<\infty.
\end{equation*}
Let $\Omega\subset\R^n$ be a bounded domain, we define $\Co^\alpha(\Omega;\Xs):=\{\tilde f|_\Omega:\tilde f\in\Co^\alpha(\R^n;\Xs)\}$ with norm 
$$\|f\|_{\Co^\alpha(\Omega;\Xs)}:=\inf\{\|\tilde f\|_{\Co^\alpha(\R^n;\Xs)}:\tilde f|_\Omega=f\}.$$

We use $\Co^{\alpha-}(\Omega;\Xs)=\bigcap_{\beta<\alpha}\Co^\beta(\Omega;\Xs)$ with standard (projective) limit topology.

We also  define $\Co^\alpha_\loc(\Omega;\Xs)$ and $\Co^{\alpha-}_\loc(\Omega;\Xs)$ to be the spaces of functions/distributions $f\in\D'(\Omega;\Xs)$ such that for every precompact open $\Omega'\Subset\Omega$, we have $f\in\Co^\alpha(\Omega';\Xs)$ and $f\in\Co^{\alpha-}(\Omega';\Xs)$ respectively.
\end{defn}
The $\Co^\alpha$-norms depend on  $\phi=(\phi_j)_{j=0}^\infty$, while the spaces themselves are not. See \cite[Chapter 2.3.2]{Triebel1}.
\begin{remark}\label{Rmk::Hold::GradBdd}
\begin{enumerate}[parsep=-0.3ex,label=(\roman*)]
    \item\label{Item::Hold::GradBdd} By \eqref{Eqn::Hold::HoldChar:Pf::GradBdd} we see that $\nabla:\Co^\alpha(\R^n)\to\Co^{\alpha-1}(\R^n;\R^n)$ is bounded for all $\alpha\in\R$. Thus by Lemma \ref{Lem::Hold::HoldChar} \ref{Item::Hold::HoldChar::DomainChar} we see that $\nabla:\Co^\alpha(\Omega)\to\Co^{\alpha-1}(\Omega;\R^n)$ is bounded for all $\alpha\in\R$ and all bounded smooth domain $\Omega\subseteq\R^n$.
    \item Lemma \ref{Lem::Hold::OperatorExtension} is still true for $\alpha\le0$ when $\Omega=\R^n$ since we have $(\phi_j\ast\tilde Tf)(x)=T(\phi_j\ast f(x))$. For bounded smooth domains, by the existence of extension operator we know $\tilde T:\Co^\alpha(\Omega;\Xs)\to\Co^\alpha(\Omega;\Ys)$ is bounded, but we do not know whether $\|\tilde T\|=\|T\|$ is still true. 
    \item In Section \ref{Section::BiHoldSec} we consider the spaces $\Co^\beta(\Omega_1;\Co^\alpha(\Omega_2))$. In our discussions only the case $\beta>0$ is needed.
\end{enumerate}
\end{remark}

\begin{lem}\label{Lem::Hold::ZygExample}
$x\log|x|\in\Co^1_\loc(\R)$, $\log|x|\in\Co^0(\R)$ and $\pv\frac1x\in\Co^{-1}(\R)$.
\end{lem}
\begin{proof}
To show $x\log|x|\in\Co^1_\loc(\R)$ one can use the direct computation via Definition \ref{Defn::Intro::DefofHold} \ref{Item::Intro::DefofHold::=1}. We use scaling via Lemma \ref{Lem::Hold::HoldChar} \ref{Item::Hold::HoldChar::02}.

Let $(\phi_j)_{j=0}^\infty$ be a dyadic resolution. We have $\phi_j(x)=2^{j-1}\phi_1(2^{j-1}x)$ and $\int x^k\phi_j(x)dx=0$ for $j\ge1$ and $k\ge0$. Let $f(x):=\log|x|$. Since $f\in L^1+L^\infty$ and $\phi_j\in L^\infty\cap L^1$, we know $f\ast\phi_j\in L^\infty(\R)$ for each $j\ge0$. By scaling we have $f(cx)=f(x)+\log|c|$ for $c>0$. Therefore
\begin{equation*}
    \textstyle f\ast\phi_j(x)=f(2^{j-1}\cdot)\ast\phi_1(x)=(f\ast\phi_1)(x)+(j-1)\log2\int_\R \phi_1(y)dy=(f\ast\phi_1)(x).
\end{equation*}
Therefore $\|f\ast\phi_j\|_{L^\infty(\R)}=\|f\ast\phi_1\|_{L^\infty(\R)}<\infty$ for all $j\ge1$. We conclude that $\log|x|\in\Co^0(\R)$.

We see that $\pv\frac1x=\frac d{dx}\log|x|$ because for every $g\in C_c^\infty(\R)$, $$\textstyle-\int(\log|x|)g'(x)dx=\lim_{\eps\to0}\int_{|x|>\eps}\frac1xg(x)dx-(g(\eps)\log\eps-g(-\eps)\log\eps)=\lim_{\eps\to0}\int_{|x|>\eps}\frac1xg(x)dx=\langle\pv\frac1x,g\rangle.$$
By Remark \ref{Rmk::Hold::GradBdd} \ref{Item::Hold::GradBdd} $\pv\frac1x\in\Co^{-1}(\R)$.

Since $x\log|x|-x=\int_0^x\log|t|dt$, and clearly $x\in C^1_\loc(\R)\subset\Co^1_\loc(\R)$, we get $x\log|x|\in\Co^1_\loc(\R)$ and finish the proof.
\end{proof}
\begin{remark}\label{Rmk::Hold::ZygVsLip}
    $x\log|x|$ is an example in $\Co^1$ that does not belong to $C^{0,1}=\Co^\Lip$. Taking definite integrals $k$ times we see that $\Co^{k+1}\subsetneq C^{k,1}(=\Co^{k+\Lip})$ for all integers $k$.
\end{remark}

\begin{remark}[The $\Co^{\alpha-}$ and $\Co^{\alpha+}$ topologies]\label{Rmk::Hold::Topology}
In a topological vector space $\Xs$, recall that a set $\Bs\subset\Xs$ is (von Neumann) bounded, if for any open neighborhood $\Us\subseteq\Xs$ of $0$, there is a $M>0$ such that $\Bs\subseteq M\cdot\Us(=\{M\cdot x:x\in\Us\})$. 

\begin{enumerate}[parsep=-0.3ex,label=(\roman*)]
    \item\label{Item::Hold::Topology::-} On $\Co^{\alpha-}(\Omega)$ we use the standard limit topology a.k.a. the initial topology: $\Us\subseteq\Co^{\alpha-}(\Omega)$ is open if and only if for every $\eps>0$ there is an open set $\Us_\eps\subseteq\Co^{\alpha-\eps}(\Omega)$ such that $\Us=\Us_\eps\cap \Co^{\alpha-}(\Omega)$. It is a separable Fr\'echet space and the topology is metrizable. Moreover, 
    \begin{itemize}[nolistsep]
        \item $\big\{\{f\in\Co^{\alpha-}(\Omega):\|f\|_{\Co^{\alpha-\eps}}<\delta\}:\eps,\delta>0\big\}$ form a neighborhood basis of $0\in\Co^{\alpha-}(\Omega)$.
        \item $f_j\xrightarrow{\Co^{\alpha-}(\Omega)}f_0$ if and only if $f_j\xrightarrow{\Co^{\alpha-\eps}(\Omega)}f_0$ for every $\eps>0$.
        \item $\Bs\subset \Co^{\alpha-}(\Omega)$ is a bounded set if and only if $\Bs\subset \Co^{\alpha-\eps}(\Omega)$ is bounded for all $\eps>0$, i.e. for every $\eps>0$ there is a $M_\eps>0$ such that $f\in\Bs$ implies $\|f\|_{\Co^{\alpha-\eps}}<M_\eps$.
    \end{itemize}
    
    One can view $\Co^{\alpha-}(\Omega)$ as the inverse limit $\Co^{\alpha-}(\Omega)=\varprojlim\{\Co^{\alpha-\eps}(\Omega):\eps>0\}$, whose arrows are natural inclusions.
    
    Similarly we have $\Co^{\alpha}_\loc(\Omega)=\varprojlim\{\Co^{\alpha}(\Omega'):\Omega'\Subset\Omega\}$ and $\Co^{\alpha-}_\loc(\Omega)=\varprojlim\{\Co^{\alpha-}(\Omega'):\Omega'\Subset\Omega\}$, both of whose arrows are natural restriction maps.
    \item We do not use the topologies of $\Co^{\alpha+}(\Omega)$ in the thesis. But we point out that in general $\Co^{\alpha+}$ is endowed with the standard colimit topology a.k.a. the final topology: $\Us\subseteq\Co^{\alpha+}(\Omega)$ is open if and only if for every $\eps>0$, the subset $\Us\cap\Co^{\alpha+\eps}(\Omega)$ is open in $\Co^{\alpha+\eps}(\Omega)$. This is not a metrizable space.
    
    Alternatively, using Besov spaces, we have duality $\Bs_{11}^{-\alpha}(\R^n)'=\Co^\alpha(\R^n)$ for all $\alpha\in\R$, see \cite[Theorem 2.11.2(i)]{Triebel1}. We can define $\Bs_{11}^{(-\alpha)-}(\R^n)$ and its topology analogous to \ref{Item::Hold::Topology::-}. One can show that the colimit topology of $\Co^{\alpha+}(\R^n)$ coincides with the weak-* topology as the dual space of $\Bs_{11}^{(-\alpha)-}(\R^n)$.

\end{enumerate}

We leave the proof to the above properties to readers.
\end{remark}


For H\"older functions with negative indices it is sometimes more convenient to write it as derivatives of classical functions.
\begin{lem}\label{Lem::Hold::NegHold=SumGood}
Let $\alpha\in\R$, and let $\Omega\subseteq\R^n$ be either the total space or a bounded smooth domain. Then $\Co^{\alpha-1}(\Omega)=\{g_0+\sum_{j=1}^n\partial_j g_j:g_0,\dots,g_n\in\Co^\alpha(\Omega)\}$.
\end{lem}
\begin{proof}
For a $f\in\Co^{\alpha-1}(\Omega)$, by Definition \ref{Defn::Hold::NegHold} there is an extension $\tilde f\in\Co^{\alpha-1}(\R^n)$. Take $\tilde g_0:=(I-\Delta)^{-1}f$ and $\tilde g_j:=-\partial_j\tilde g_0$ for $1\le j\le n$, we see that $\tilde g_0,\dots,\tilde g_n\in\Co^\alpha(\R^n)$. Hence for $g_j:=\tilde g_j|_\Omega\in\Co^\alpha(\Omega)$, $0\le j\le n$, we have $f=g_0+\sum_{j=1}^n\partial_jg_j$.

Conversely, let $g_0,\dots,g_n\in\Co^\alpha(\Omega)$, by Lemma \ref{Lem::Hold::HoldChar} \ref{Item::Hold::HoldChar::DomainChar} and Definition \ref{Defn::Hold::NegHold} they have extensions $\tilde g_0,\dots,\tilde g_n\in\Co^\alpha(\R^n)$. Thus by Remark \ref{Rmk::Hold::GradBdd} \ref{Item::Hold::GradBdd} $\tilde g_0+\sum_{j=1}^n\partial_j\tilde g_j\in\Co^{\alpha-1}(\R^n)$. Taking restrictions to $\Omega$ we get $g_0+\sum_{j=1}^n\partial_j g_j\in\Co^{\alpha-1}(\Omega)$, finishing the proof.
\end{proof}

In Sections \ref{Section::EllipticPara::Scaling} and \ref{Section::Rough1FormSum} we need to do scaling on functions.
\begin{lem}\label{Lem::Hold::ScalingLem}
Let $\alpha>0$ and $\mu_0>0$, there is a $C=C_{n,\alpha,\mu_0}>0$ such that for any Banach space $\Xs$ and any $f\in\Co^\alpha(\B^n;\Xs)$,
\begin{equation}\label{Eqn::Hold::ScalingLemma}
    \|f(\mu\cdot)\|_{\Co^\alpha(\B^n;\Xs)}\le |f(0)|_{\Xs}+C\mu^{\min(\alpha,\frac12)}\|f\|_{\Co^\alpha(\B^n;\Xs)},\quad\forall0<\mu<\mu_0.
\end{equation}
\end{lem}
\begin{proof}
By Definition \ref{Defn::Intro::DefofHold}, the constant function $x\mapsto f(0)$ always has norm $\|f(0)\|_{\Co^\alpha(\B^n;\Xs)}=|f(0)|_\Xs$. Thus replacing $f(x)$ by $f(x)-f(0)$ we can assume $f(0)=0$.

By taking a scaling $x\mapsto\mu_0 x$, we can assume $\mu_0=1$ without loss of generality. Thus, $f$ is defined on the unit ball.

For $\mu\in(0,1]$ set $f_\mu(x):=f(\mu x)$. 
For $x\in\B^n$ and $\mu\in(0,1]$, by Definition \ref{Defn::Intro::DefofHold} \ref{Item::Intro::DefofHold::<1},
\begin{equation}\label{Eqn::Hold::ScalingLemma::ProofSup}
    |f_\mu(x)|_\Xs=|f(\mu x)-f(0)|_\Xs\lesssim\|f\|_{\Co^{\min(\alpha,\frac12)}(\B^n;\Xs)}|\mu x-0|^{\min(\alpha,\frac12)}\lesssim\|f\|_{\Co^\alpha(\B^n)}\mu^{\min(\alpha,\frac12)}. 
\end{equation}
When $\alpha\in(0,2)$, using Lemma \ref{Lem::Hold::HoldChar} \ref{Item::Hold::HoldChar::02}, for $x_1,x_2\in\B^n$,
\begin{equation}\label{Eqn::Hold::ScalingLemma::Proof02}
    \textstyle|\frac {f_\mu(x_1)+f_\mu(x_2)}2-f_\mu(\frac{x_1+x_2}2)|_\Xs=\left|\frac {f(\mu x_1)+f(\mu x_2)}2-f(\mu\frac{x_1+x_2}2)\right|_\Xs\lesssim_\alpha\|f\|_{\Co^{\min(\alpha,\frac12)}}|\mu(x_1-x_2)|^\alpha\le\mu^\alpha\|f\|_{\Co^\alpha}|x_1-x_2|^\alpha.
\end{equation}
    Combining \eqref{Eqn::Hold::ScalingLemma::ProofSup} and \eqref{Eqn::Hold::ScalingLemma::Proof02}, we get \eqref{Eqn::Hold::ScalingLemma} for the case $0<\alpha<2$, since
    $$\textstyle\|f_\mu\|_{\Co^\alpha(\B^n)}\approx\sup\limits_{x\in\B^n}|f_\mu(x)|_\Xs+\sup\limits_{x_1,x_2\in\B^n}|x_1-x_2|^{-\alpha}\left|\frac {f_\mu(x_1)+f_\mu(x_2)}2-f_\mu(\frac{x_1+x_2}2)\right|_\Xs\lesssim_\alpha\mu^{\min(\alpha,\frac12)}\|f\|_{\Co^\alpha(\B^n)}.$$

    For $\alpha\ge2$, we proceed by induction. We prove the result for $\alpha\in [l,l+1)$, for $l\in \{1,2,\ldots\}$.  The base case, $l=1$ was shown above. We assume the result for $l-1$ and prove it for $l$.
    
    Assume $\alpha\in[l,l+1)$ where $l\ge2$. Note that $\nabla f_\mu(x)=\mu(\nabla f)(\mu x)$, so $\|\partial_{x^j}(f_\mu)\|_{\Co^{\alpha-1}(\B^n)}=\mu\|(\partial_{x^j}f)_\mu\|_{\Co^{\alpha-1}(\B^n)}\le \|(\partial_{x^j}f)_\mu\|_{\Co^{\alpha-1}(\B^n)}$ for $j=1,\dots,n$. Here $(\partial_{x^j}f)_\mu(x)=(\partial_{x^j}f)(\mu x)$. 
    
    By the inductive hypothesis $\|f_\mu\|_{\Co^{\alpha-1}(\B^n)}\le C_{\alpha-1}\mu^{\frac12}\|f\|_{\Co^{\alpha-1}(\B^n)}$ and for $j=1,\dots,n$, $\|(\partial_{x^j}f)_\mu\|_{\Co^{\alpha-1}(\B^n)}\le C_{\alpha-1}\mu^{\frac12}\|\partial_{x^j}f\|_{\Co^{\alpha-1}(\B^n)}$. So by Definition \ref{Defn::Intro::DefofHold} \ref{Item::Hold::DefofHold::>1} we get 
    \begin{align*}
        &\|f_\mu\|_{\Co^\alpha(\B^n;\Xs)}\approx\|f_\mu\|_{\Co^{\alpha-1}(\B^n;\Xs)}+\sum_{j=1}^n\|\partial_{x^j}(f_\mu)\|_{\Co^{\alpha-1}(\B^n)}
        \\&\lesssim\mu^\frac12\Big(\|f\|_{\Co^{\alpha-1}(\B^n)}+\sum_{j=1}^n\|\partial_{x^j}f\|_{\Co^{\alpha-1}(\B^n)}\Big)\approx \mu^\frac12\|f\|_{\Co^{\alpha}(\B^n)}=\mu^{\min(\alpha,\frac12)}\|f\|_{\Co^{\alpha}(\B^n)},
    \end{align*}
    completing the proof.
\end{proof}

For products of functions and distributions the paraproduct is of great importance:
\begin{lem}[Paraproduct decomposition]\label{Lem::Hold::LemParaProd}
Let $(\phi_j)_{j=0}^\infty\subset\Sc(\R^n)$ be a dyadic resolution. For bounded continuous functions $u,v$ defined on $\R^n$, consider the following decomposition
\begin{equation}\label{Eqn::Hold::LemParaProd::ParaDecomp1}
    uv=\sum_{j=0}^\infty\sum_{j'=0}^{j-3}(\phi_j\ast u)(\phi_{j'}\ast v)+\sum_{j'=0}^\infty\sum_{j=0}^{j'-3}(\phi_j\ast u)(\phi_{j'}\ast v)+\sum_{|j-j'|\le 2}(\phi_j\ast u)(\phi_{j'}\ast v).
\end{equation}

Then for every $l\ge0$, 
    \begin{equation}\label{Eqn::Hold::LemParaProd::ParaDecomp2}
    \phi_l\ast(uv)=\phi_l\ast\bigg(\sum_{j=l-2}^{l+2}\sum_{j'=0}^{j-3}(\phi_j\ast u)(\phi_{j'}\ast v)+\sum_{j'=l-2}^{l+2}\sum_{j=0}^{j'-3}(\phi_j\ast u)(\phi_{j'}\ast v)+\sum_{\substack{j,j'\ge l-3\\|j-j'|\le 2}}(\phi_j\ast u)(\phi_{j'}\ast v)\bigg).
\end{equation}
\end{lem}

\begin{proof}
\eqref{Eqn::Hold::LemParaProd::ParaDecomp2} can be obtained by investigating the Fourier support condition in Definition \ref{Defn::Hold::DyadicResolution}. Note that $\phi_l\ast\left((\phi_j\ast u)(\phi_{j'}\ast v)\right)\neq0$ only when $\supp\hat\phi_l\cap(\supp\hat\phi_j+\supp\hat\phi_{j'})\neq\varnothing$, and by \eqref{Eqn::Hold::RmkDyaSupp} we have on the following:
$$
\supp\hat\phi_l\subset\begin{cases}\{2^{l-1}<|\xi|<2^{l+1}\},&l\ge1,\\\{|\xi|<2\},&l=0;\end{cases}\ \supp\hat\phi_j+\supp\hat\phi_{j'}\subset\begin{cases}\{2^{\max(j,j')-2}<|\xi|<2^{\max(j,j')+2}\},&|j-j'|\ge3,\\\{|\xi|<2^{\max(j,j')+2}\},&|j-j'|\le2.\end{cases}$$

For details, see  \cite[Proof of Theorem 2.8.2]{Triebel1} or \cite[Proof of Theorem 2.52]{BahouriCheminDanchin}. Note that their choice of decompositions \eqref{Eqn::Hold::LemParaProd::ParaDecomp1} may have slight differences.
\end{proof}
\begin{defn}\label{Defn::Hold::ParaOp}
We define $\Pf(u,v):=\sum_{j=0}^\infty\sum_{j'=0}^{j-3}(\phi_j\ast u)(\phi_{j'}\ast v)$ and $\Rf(u,v):=\sum_{|j-j'|\le 2}(\phi_j\ast u)(\phi_{j'}\ast v)$ to be the paraproduct operators associated with $\phi$.
\end{defn}
\begin{remark}\label{Rmk::Hold::LemParaProd}
    From \eqref{Eqn::Hold::LemParaProd::ParaDecomp1} we see that $uv=\Pf(u,v)+\Pf(v,u)+\Rf(u,v)$. In \cite[Chapter 2.8.1, (2.41)]{BahouriCheminDanchin} this is denoted by $uv=T(v,u)+T(u,v)+R(u,v)$.
\end{remark}

\begin{lem}\label{Lem::Hold::ParaBdd} The paraproduct operators have the following boundedness:
\begin{align}
\label{Eqn::Hold::ParaBdd::PfBdd1}
    \Pf:\Co^\alpha(\R^n)\times L^\infty(\R^n)\to\Co^\alpha(\R^n),&\qquad\alpha\in\R;
    \\
\label{Eqn::Hold::ParaBdd::PfBdd2}
    \Pf:\Co^\alpha(\R^n)\times \Co^\beta(\R^n)\to\Co^{\alpha+\beta}(\R^n),&\qquad\alpha\in\R,\quad\beta<0;
    \\
\label{Eqn::Hold::ParaBdd::RfBdd}
    \Rf:\Co^\alpha(\R^n)\times \Co^\beta(\R^n)\to\Co^{\alpha+\beta}(\R^n),&\qquad\alpha,\beta\in\R,\quad\alpha+\beta>0.
\end{align}

In particular $\Pf:\Co^\alpha(\R^n)\times \Co^\beta(\R^n)\to\Co^\alpha(\R^n)$ holds for $\alpha\in\R$, $\beta>0$ and $\Rf :\Co^\alpha(\R^n)\times L^\infty(\R^n)\to\Co^\alpha(\R^n)$ holds for $\alpha>0$.

\end{lem}
Also see Lemma \ref{Lem::Hold::ImprovedParaBdd} for a precise control of their norms.
\begin{proof}
See \cite[Theorem 3.3.2(ii)]{Triebel1} for \eqref{Eqn::Hold::ParaBdd::PfBdd1}, \eqref{Eqn::Hold::ParaBdd::PfBdd2} and \cite[Corollary 2.86]{BahouriCheminDanchin} for \eqref{Eqn::Hold::ParaBdd::RfBdd}. Note that by Lemma \ref{Lem::Hold::HoldChar} and Definition \ref{Defn::Hold::NegHold} H\"older-Zygmund spaces are a special case to Besov spaces $\Co^\alpha(U)=\Bs_{\infty\infty}^\alpha(U)$ for  all $\alpha\in\R$.
\end{proof}

For the product itself we have
\begin{lem}\label{Lem::Hold::Product}
Let $\alpha>0$ and $\beta\in(-\alpha,\alpha]$. Let $\Omega\subseteq\R^n$ be either the total space of a bounded smooth domain.
\begin{enumerate}[nolistsep,label=(\roman*)]
    \item\label{Item::Hold::Product::Hold1} There is a $C=C(\Omega,\alpha,\beta)\ge0$ such that $\|fg\|_{\Co^\beta(\Omega)}\le C\|f\|_{\Co^\alpha(\Omega)}\|g\|_{\Co^\beta(\Omega)}$ for all $f\in\Co^\alpha(\Omega)$, $g\in\Co^\beta(\Omega)$.
    \item\label{Item::Hold::Product::Hold2}Let $\Xs$ be a Banach algebra (that is $\Xs$ is a Banach space with product structure such that $\|uv\|_\Xs\le\|u\|_\Xs\|v\|_\Xs$). Then $\|fg\|_{\Co^\alpha(\Omega;\Xs)}\le \|f\|_{\Co^\alpha(\Omega;\Xs)}\|g\|_{L^\infty(\Omega;\Xs)}+\|f\|_{L^\infty(\Omega;\Xs)}\|g\|_{\Co^\alpha(\Omega;\Xs)}$ holds for all $f,g\in\Co^\alpha(\Omega;\Xs)$.
    
    In particular $\|fg\|_{\Co^\alpha(\Omega;\Xs)}\le 2\|f\|_{\Co^\alpha(\Omega;\Xs)}\|g\|_{\Co^\alpha(\Omega;\Xs)}$.
\end{enumerate}
\end{lem}
\begin{proof}
    To prove \ref{Item::Hold::Product::Hold1}, by Lemma \ref{Lem::Hold::ParaBdd} we see that $(u,v)\mapsto \Pf(u,v)$, $(u,v)\mapsto \Pf(v,u)$ and $(u,v)\mapsto\Rf(u,v)$ are all bounded as $\Co^\alpha(\R^n)\times \Co^\beta(\R^n)\to\Co^\beta(\R^n)$, by Lemma \ref{Lem::Hold::LemParaProd} we get $(u,v)\mapsto uv:\Co^\alpha(\R^n)\times \Co^\beta(\R^n)\to\Co^\beta(\R^n)$. Using Lemma \ref{Lem::Hold::HoldChar} \ref{Item::Hold::HoldChar::DomainChar} and taking restriction on $\Omega$, we get the boundedness $\Co^\alpha(\Omega)\times \Co^\beta(\Omega)\to\Co^\beta(\Omega)$.
    
    The result \ref{Item::Hold::Product::Hold2} follows from direct computation using Definition \ref{Defn::Intro::DefofHold}. We remark that there is no non-trivial in the result. We leave the detail to reader.
\end{proof}

We also remark that the product operator has the local property.
\begin{lem}\label{Lem::Hold::MultLoc}
Let $\Omega\subseteq\R^n$ be an arbitrary open set, and let $\alpha,\beta\in\R$ satisfy $\alpha+\beta>0$. 
\begin{enumerate}[parsep=-0.3ex,label=(\roman*)]
    \item\label{Item::Hold::MultLoc::Loc} Suppose $\tilde f_1,\tilde f_2\in \Co^\alpha(\R^n)$ and $\tilde g_1,\tilde g_2\in \Co^\beta(\R^n)$ satisfy $\tilde f_1|_\Omega=\tilde f_2|_\Omega$ and $\tilde g_1|_\Omega=\tilde g_2|_\Omega$. Then $(\tilde f_1\tilde g_1)|_\Omega=(\tilde f_2\tilde g_2)|_\Omega$ as distributions on $\Omega$.
    \item\label{Item::Hold::MultLoc::WellDef} In particular for $f\in \Co^\alpha_\loc(\Omega)$ and $g\in \Co^\beta_\loc(\Omega)$, the product $fg$ is well-defined distribution in $\Co^{\min(\alpha,\beta)}_\loc(\Omega)$.
\end{enumerate}
\end{lem}
Also see \cite[Section 4]{BookProduct} for more details.
\begin{proof}

\ref{Item::Hold::MultLoc::Loc}: Let $i_1,i_2\in\{1,2\}$. By Lemma \ref{Lem::Hold::HoldChar} \ref{Item::Hold::HoldChar::LPHoldChar} and Definition \ref{Defn::Hold::NegHold} we have
\begin{equation*}
    \|\Su_\sigma f_{i_1}-f_{i_1}\|_{\Co^{\alpha-\delta}}\le\sum_{j=\sigma+1}^\infty\|\phi_j\ast f_{i_1}\|_{\Co^{\alpha-\delta}}\lesssim_{\phi,\alpha,\delta} \sum_{j=\sigma+1}^\infty2^{-j\alpha+j(\alpha-\delta)}\|f_{i_1}\|_{\Co^\alpha}\approx_{\delta} 2^{-\sigma\delta}\|f_{i_1}\|_{\Co^\alpha}\xrightarrow{\sigma\to\infty}0.
\end{equation*} Therefore $\lim_{\sigma\to\infty}\Su_\sigma f_{i_1}\xrightarrow{\Co^{\alpha-\delta}(\R^n)}f_{i_1}$ and similarly $\lim_{\sigma\to\infty}\Su_\sigma g_{i_2}\xrightarrow{\Co^{\beta-\delta}(\R^n)}g_{i_2}$ for all $\delta>0$. Thus by Lemma \ref{Lem::Hold::Product} \ref{Item::Hold::Product::Hold1} we have convergence $\lim_{\sigma\to\infty}\Su_\sigma\tilde f_{i_1}\Su_\sigma\tilde g_{i_2}=\tilde f_{i_1}\tilde g_{i_2}$ in $\Co^{\min(\alpha,\beta)-\delta}(\R^n)$ for every $\delta>0$. By \cite[Lemma 4.2.2]{BookProduct} we see that 
$$(\tilde f_1\tilde g_1-\tilde f_2\tilde g_2)|_\Omega=((\tilde f_1-\tilde f_2)\tilde g_1+\tilde f_2(\tilde g_1-\tilde g_2))|_\Omega=\lim_{\sigma\to\infty}\big(\Su_\sigma(\tilde f_1-\tilde f_2)\cdot\Su_\sigma\tilde g_1+\Su_j\tilde f_2\cdot\Su_\sigma(\tilde g_1-\tilde g_2)\big)|_\Omega=0,\quad \text{in }\D'(\Omega).$$

\noindent\ref{Item::Hold::MultLoc::WellDef}: Let $f\in \Co^\alpha_\loc(\Omega)$ and $g\in \Co^\beta_\loc(\Omega)$. By \ref{Item::Hold::MultLoc::Loc} it suffices to prove the well-definedness of $fg|_{\Omega'}\in \Co^{\min(\alpha,\beta)}(\Omega')$ on any precompact open $\Omega'\Subset\Omega$.

Let $\chi\in C_c^\infty(\Omega)$ be such that $\chi|_{\Omega'}\equiv1$. By taking zero extension outside $\Omega$ we have $\chi f\in \Co^\alpha(\R^n)$ and $\chi g\in \Co^\beta(\R^n)$. By Lemma \ref{Lem::Hold::Product} \ref{Item::Hold::Product::Hold1} we have $\chi f\cdot \chi g\in \Co^\beta(\R^n)$. Since $\chi f|_{\Omega'}=f|_{\Omega'}$ and $\chi g|_{\Omega'}=g|_{\Omega'}$, by \ref{Item::Hold::MultLoc::Loc} we see that $fg|_{\Omega'}=\chi f\cdot \chi g|_{\Omega'}\in \Co^{\min(\alpha,\beta)}(\Omega')$ is defined, finishing the proof.
\end{proof}
\begin{cor}\label{Cor::Hold::[X,Y]WellDef}
Let $\frac12<\alpha<1$, and let $\Omega\subseteq\R^n$ be an open set. If $X$ and $Y$ are $\Co^\alpha_\loc$ vector fields on $\Omega$, then $[X,Y]\in \Co^{\alpha-1}_\loc(\Omega;\R^n)$ is a vector field with distributional coefficients. In particular $[X,Y]$ is defined when $X$ and $Y$ are log-Lipschitz.
\end{cor}

Recall that a vector field $X$ on $\Omega\subseteq\R^n$ is a vector-valued map $X=(X^1,\dots,X^n):\Omega\to\R^n$ as well as a differential operator $Xf=\sum_{j=1}^nX^j\partial_jf$. For $X=\sum_{i=1}^nX^i\Coorvec{x^i}$ and $Y=\sum_{j=1}^nY^j\Coorvec{x^j}$ the Lie bracket $[X,Y]$ is given by
\begin{equation}\label{Eqn::Hold::LieBracketEqn}
    [X,Y]=\sum_{i,j=1}^n\Big(X^j\frac{\partial Y^i}{\partial x^j}-Y^j\frac{\partial X^i}{\partial x^j}\Big)\Coorvec{x^i}.
\end{equation}
\begin{proof} Write $X=:\sum_{j=1}^nX^j\Coorvec{x^j}$ and $Y=:\sum_{j=1}^nY^j\Coorvec{y^j}$, we have $X^j,Y^j\in \Co^\alpha_\loc(\Omega)$ and $\frac{\partial Y^i}{\partial x^j},\frac{\partial X^i}{\partial x^j}\in \Co^{\alpha-1}_\loc(\Omega)$ for $1\le i,j\le n$. Applying Lemma \ref{Lem::Hold::MultLoc} \ref{Item::Hold::MultLoc::WellDef} since $\alpha>1-\alpha$, we get $X^j\frac{\partial Y^i}{\partial x^j}-Y^j\frac{\partial X^i}{\partial x^j}\in \Co^{\alpha-1}_\loc(\Omega)$ for all $1\le i,j\le n$. By \eqref{Eqn::Hold::LieBracketEqn}, we conclude that $[X,Y]\in \Co^{\alpha-1}_\loc(\Omega;\R^n)$.
\end{proof}

In the thesis we often discuss Lipschitz spaces and log-Lipschitz spaces.

\begin{defn}\label{Defn::Hold::CLogLandLip}
Let $\Omega\subseteq\R^n$ be an open subset and $\Xs$ be a Banach space. We denote by $\Co^\Lip(\Omega;\Xs)$ the class of $\Xs$-valued Lipschitz functions with norm
\begin{equation*}
    \|f\|_{\Co^\Lip(\Omega;\Xs)}:=\sup_{x\in \Omega}|f(x)|_\Xs+\sup_{x,y\in \Omega;x\neq y}|f(x)-f(y)|_\Xs\cdot|x-y|^{-1}.
\end{equation*}
We say a continuous function $f:\Omega\to\Xs$ is bounded \textbf{log-Lipschitz}, denoted by $f\in\Co^\LogL(\Omega;\Xs)$, if
\begin{equation*}
    \|f\|_{\Co^\LogL(\Omega;\Xs)}:=\sup_{x\in \Omega}|f(x)|_\Xs+\sup_{x,y\in \Omega;0<|x-y|<1}|f(x)-f(y)|_\Xs\Big(|x-y|\log\frac e{|x-y|}\Big)^{-1}<\infty.
\end{equation*}

We say $f:\Omega\to\Xs$ is bounded \textbf{little log-Lipschitz}, denoted by $f\in\Co^\logl(\Omega;\Xs)$, if in addition $\lim_{r\to0}\sup_{x,y\in U;0<|x-y|<r}|f(x)-f(y)|(|x-y|\log\frac e{|x-y|})^{-1}=0$.

Let $\sigma\in\{\Lip,\LogL\}$, we denoted by $\Co^{k+\sigma}(\Omega;\Xs)$ the space of all $f\in\Co^{\sigma}(\Omega;\Xs)$ whose first $k$-order (Fr\'echet) derivatives are $\Co^\sigma$, with norm $\|f\|_{\Co^{k+\sigma}(\Omega;\Xs)}:=\sum_{|\nu|\le k}\|\partial^\nu f\|_{\Co^\sigma(\Omega;\Xs)}$.
\end{defn}

One can check that $\Co^\logl(\Omega;\Xs)\subseteq\Co^\LogL(\Omega;\Xs)$ is a closed subspace. We have correspondence of the classical notation $\Co^{k+\Lip}=C^{k,1}=W^{k+1,\infty}$.

\begin{defn}\label{Defn::Hold::LogL-1}
Let $\phi=(\phi_j)_{j=0}^\infty$ be a fixed dyadic resolution and let $\Xs$ be a Banach space. We define $\Co^{\LogL-1}(\R^n;\Xs)$ to be the space of all tempered distributions $f\in\Sc'(\R^n;\Xs)$ such that 
\begin{equation*}
    \|f\|_{\Co^{\LogL-1}(\R^n;\Xs)}:=\sup_{j\ge1}j^{-1}\|\Su_jf\|_{L^\infty(\R^n;\Xs)}<\infty.
\end{equation*}
For an arbitrary open set $\Omega\subseteq\R^n$, we  define $\Co^{\LogL-1}(\Omega;\Xs):=\{\tilde f|_\Omega:\tilde f\in\Co^{\LogL-1}(\R^n;\Xs)\}$ with $\|f\|_{\Co^{\LogL-1}(\Omega;\Xs)}=\inf_{\tilde f|_\Omega=f}\|_{\Co^{\LogL-1}(\R^n;\Xs)}$.

We denote $\Co^{\Lip-1}(\Omega;\Xs):=L^\infty(\Omega;\Xs)$.
\end{defn}The $\Co^{\LogL-1}$-norm depends on the choice of the dyadic resolution $\phi$, but the norm topology is not. See Lemma \ref{Lem::Hold::CharLogL-1} \ref{Item::Hold::CharLogL-1::Norm} below.
\begin{lem}\label{Lem::Hold::CharLogL-1}
\begin{enumerate}[parsep=-0.3ex,label=(\roman*)]
    \item\label{Item::Hold::CharLogL-1::Norm} Let $(\phi_j)_{j=0}^\infty$ and $(\psi_j)_{j=0}^\infty$ be two dyadic resolutions, then there is a $C=C(n,\phi,\psi)>0$ such that $\sup_{j\ge1}j^{-1}\|\Su^\phi_jf\|_{L^\infty(\R^n;\Xs)}\le C\sup_{j\ge1}j^{-1}\|\Su^\psi_jf\|_{L^\infty(\R^n;\Xs)} $. In other words different choices of dyadic resolutions result in equivalent norms for $\Co^{\LogL-1}(\R^n;\Xs)$ and $\Co^{\LogL-1}(\Omega;\Xs)$.
    \item\label{Item::Hold::CharLogL-1::Grad} Let $\Omega\subseteq\R^n$ be either the total space or a bounded smooth domain. Then the product map $\Co^\alpha(\Omega)\times\Co^{\LogL-1}(\Omega)\to\Co^{\LogL-1}(\Omega)$ is bounded bilinear for every $\alpha>0$. Moreover taking derivative $\nabla:\Co^{\LogL}(\Omega)\to \Co^{\LogL-1}(\Omega;\R^n)$ is also bounded linear.
    \item\label{Item::Hold::CharLogL-1::CompLin}Let $L$ be an invertible linear map on $\R^n$. Then $[f\mapsto f\circ L]:\Co^{\LogL-1}(\R^n)\to\Co^{\LogL-1}(\R^n)$. Moreover $\|f\circ L\|_{\Co^{\LogL-1}}\le C\max(1,\log|L|_{\ell^2})\cdot\|f\|_{\Co^{\LogL-1}}$ where $C$ depends only on $\phi$.
\end{enumerate}
\end{lem}
Here for the linear map $L$ we use matrix norm in \eqref{Eqn::Intro::MatrixNorms}.
We are going to extend \ref{Item::Hold::CharLogL-1::CompLin} to the composition of general diffeomorphism in Lemma \ref{Lem::Hold::LogL-1Comp}.
\begin{proof}\ref{Item::Hold::CharLogL-1::Norm}:
Using the property \eqref{Eqn::Hold::RmkDyaSupp} we have $\Su_{j+2}^\psi\Su_j^\phi=\Su_j^\phi$, thus 
$$j^{-1}\|\Su_j^\phi f\|_{L^\infty}\le j^{-1}\|\Su_j^\phi\|_{C^0\to C^0}\|\Su_{j+2}^\psi f\|_{L^\infty}\le 3\|2^{jn}\phi_0(2^j\cdot)\|_{L^1}(j+2)^{-1}\|\Su_{j+2}^\psi f\|_{L^\infty},\quad\forall j\ge1.$$
Taking sup over $j$ we get $\sup_{j\ge1}j^{-1}\|\Su^\phi_jf\|_{L^\infty(\R^n;\Xs)}\lesssim_{\phi}\sup_{j\ge1}j^{-1}\|\Su^\psi_jf\|_{L^\infty(\R^n;\Xs)}$. 

Thus the norms (using $\phi$ or $\psi$) for $\Co^{\LogL-1}(\R^n;\Xs)$ are equivalent. Taking restrictions to $\Omega$ we get the equivalent norms result for  $\Co^{\LogL-1}(\Omega;\Xs)$.

\medskip\noindent\ref{Item::Hold::CharLogL-1::Grad}: By passing to extensions it suffices to consider $\Omega=\R^n$. The boundedness of $\nabla:\Co^\LogL\to\Co^{\LogL-1}$ follows from \cite[Page 119 Example]{BahouriCheminDanchin}.

For $f\in\Co^\alpha(\R^n)$ and $g\in\Co^{\LogL-1}(\R^n)$, we use Lemma \ref{Lem::Hold::LemParaProd} to estimate $fg=\Pf(f,g)+\Pf(f,g)+\Rf(f,g)$.

Since $\Co^{\LogL-1}(\R^n)\hookrightarrow\Co^{-\alpha/2}(\R^n)$, by Lemma \ref{Lem::Hold::ParaBdd} we have $\Pf,\Rf:\Co^\alpha\times\Co^{-\alpha/2}\to\Co^{\alpha/2}\hookrightarrow\Co^{\LogL-1}$. It remains to prove the boundedness for $(f,g)\mapsto\Pf(g,f)$. Indeed for $\sigma\ge1$,
\begin{align*}
    &\Su_\sigma\Pf(g, f)=\Su_\sigma\sum_{j=0}^{\sigma+2}\De_jg\cdot\Su_{j-3}f=\Su_\sigma(\Su_{\sigma+2}g\cdot\Su_{\sigma-1}f)-\Su_\sigma\sum_{j=0}^{\sigma+1}\Su_jg\cdot\De_{j-2}f;
    \\
    \Rightarrow&\|\Su_\sigma(g,f)\|_{L^\infty}\lesssim\|g\|_{\Co^{\LogL-1}}\|f\|_{\Co^\alpha}\Big((\sigma+2)-\sum_{j=0}^{\sigma+1}j2^{-\alpha(j-2)}\Big)\lesssim\sigma\|g\|_{\Co^{\LogL-1}}\|f\|_{\Co^\alpha}.
\end{align*}
Thus we get \ref{Item::Hold::CharLogL-1::Grad}.



\medskip\noindent\ref{Item::Hold::CharLogL-1::CompLin}: Using integration by substitution we have $\phi_j\ast(f\circ L)(x)=((\phi_j\circ L^\Inv)\ast f(Lx))\cdot|\det L|^{-1}$. In other words $\Su^\phi_j(f\circ L)(x)=|\det L|^{-1}\cdot\Su^{\phi\circ L^\Inv}_jf(Lx)$.

We have the Fourier transform $(\phi_0\ast L^\Inv)^\wedge(\xi)=|\det L|\hat\phi_0(L\xi)$. Thus by \eqref{Eqn::Hold::RmkDyaSupp} $\supp(\phi_0\ast L^\Inv)^\wedge\subset\{|\xi|\le 2^{j_0-1}\}$ where $j_0\ge1$ is the smallest number such that $2^{j_0-2}\ge|L|_{\ell^2}$. Therefore $\Su_j^{\phi\circ L^\Inv}\Su_{j+j_0}^\phi=\Su_j^{\phi\circ L^\Inv}$ for all $j\ge0$, and we have
\begin{align*}
    \|\Su^\phi_j(f\circ L)\|_{L^\infty}&=|\det L|^{-1}\|\Su^{\phi\circ L^\Inv}_jf\|_{L^\infty}=|\det L|^{-1}\|\Su^{\phi\circ L^\Inv}_j\circ\Su_{j+j_0}^\phi f\|_{L^\infty}
    \\
    &\le|\det L|^{-1}\|\phi_0\circ L^\Inv\|_{L^1}(j+j_0)\|f\|_{\Co^{\LogL-1}_\phi}=(j+j_0)\|\phi_0\|_{L^1}\|f\|_{\Co^{\LogL-1}_\phi}.
    \end{align*}
Since $j+j_0\le j_0j\le\max(1,2+\log_2|L|_{\ell^2})\cdot j$, taking the supremum over $j\ge1$ we complete the proof.
\end{proof}

\begin{defn}\label{Defn::Hold::ExtIndex}
We use the extended index set
\begin{equation*}
    \R_\Eb:=\{\alpha,\alpha-:\alpha\in\R\}\cup\{k+\LogL,k+\Lip:k=-1,0,1,2,\dots\}\cup\{\infty\}.
\end{equation*}
For two generalized indices $\sigma_1,\sigma_2\in\R_\Eb$, we define order $\sigma_1\le\sigma_2$ if $\Co^{\sigma_1}\supseteq\Co^{\sigma_2}$.

We denote by $\R_\Eb^+:=\{\sigma\in\R_\Eb:\sigma>\Lip-1\}$ the set of positive generalized indices. Equivalently $\alpha\in\R_\Eb^+$ if and only if $\alpha\in\R_\Eb$ and $\Co^\alpha\subset C^0$.

We denote by $\R^-:=\{\alpha-:\alpha\in\R\}$ and $\R_+^-:=\{\alpha-:\alpha\in\R_+\}$.
\end{defn}
Clearly $\Co^{1+}\subsetneq\Co^\Lip\subsetneq\Co^1$ and $\Co^\LogL\subsetneq\Co^{1-}$. The Zygmund space $\Co^1$ is contained in the log-Lipschitz space $\Co^\LogL$ (see \cite[Proposition 2.107]{BahouriCheminDanchin}). Therefore $\R_\Eb$ is a total order set and we have order $(k+1)-<k+\LogL<k+1<k+\Lip<k+1+\eps$ for all $\eps>0$ and all $k=-1,0,1,2,\dots$.

\begin{remark}
We do not include the $\logl$, the index for little log-Lipschitz in the index set $\R_\Eb$. On one hand $\Co^\logl$ is just the closure of smooth functions in $\Co^\LogL$, which reduces many discussion to the log-Lipschitz index $\LogL$. On the other hand, the little log-Lipschitz space $\Co^\logl$ and the Zygmund space $\Co^1$ do not have inclusions each other, which will make the index set not totally ordered.
\end{remark}
\begin{conv}\label{Conv::Hold::ExtIndOp}
We use the convention $\infty-=\infty+1=\infty-1=\frac12\infty=\infty$.

For $\alpha,\beta\in\R\cup\{\infty\}$, we use the convention $(\alpha-)-=\alpha-$. We use the formal operations $(\alpha-)+\beta=\alpha+(\beta-)=(\alpha-)+(\beta-):=(\alpha+\beta)-$; and if $\alpha,\beta\in\R_+\cup\{\infty\}$, we use $(\alpha-)\cdot\beta=\alpha\cdot(\beta-)=(\alpha-)\cdot(\beta-):=(\alpha\beta)-$.
\end{conv}

Later in Sections \ref{Section::RealFro::BLFro} and \ref{Section::PfCpxFro}, we consider the products and compositions of function spaces with the extended indices.

We can extend Lemma \ref{Lem::Hold::Product} to spaces of extended indices.
\begin{lem}\label{Lem::Hold::ProdforExt}
Let $\Omega\subseteq\R^n$ be either the total space or a bounded smooth domain. Let $\alpha,\beta\in\R_\Eb^+$. Then the product map $\Co^\alpha(\Omega)\times\Co^\beta(\Omega)\to\Co^{\min(\alpha,\beta)}(\Omega)$ is always bilinearly continuous.
\end{lem}

\begin{proof}
First we prove special case $\alpha=\beta$. The boundedness for $\alpha\in\R_+$ is the result of Lemma \ref{Lem::Hold::Product}. For $\alpha\in\{\LogL,\Lip\}$ the standard argument shows $\|fg\|_{\Co^\alpha}\le\|f\|_{L^\infty}\|g\|_{\Co^\alpha}+\|f\|_{\Co^\alpha}\|g\|_{L^\infty}$ which implies the boundedness. For $\alpha=k+\LogL$ or $k+\Lip$, the boundedness follows from the product rule $\nabla(fg)=f\nabla g+g\nabla f$ since by definition $\|fg\|_{\Co^{\alpha}}\approx\|\nabla(fg)\|_{\Co^{\alpha-1}}+\|fg\|_{\Co^{\alpha-1}}$.

By Definition \ref{Defn::Intro::DefofHold} \ref{Item::Intro::DefofHold::-}, we have $\Co^{\alpha-}(\Omega)\times\Co^{\alpha-}(\Omega)\to\Co^{\alpha-}(\Omega)$ for all $\alpha\in\R_+$. Similarly $\Co^\infty(\Omega)\times\Co^\infty(\Omega)\to\Co^\infty(\Omega)$ holds as well.

Now if $\alpha<\beta$ then the inclusion map $\Co^\beta(\Omega)\hookrightarrow\Co^\alpha(\Omega)$ is bounded, thus $\Co^\alpha(\Omega)\times\Co^\beta(\Omega)\hookrightarrow\Co^{\min(\alpha,\beta)}(\Omega)\times\Co^{\min(\alpha,\beta)}(\Omega)$ is also bounded. This completes the proof.
\end{proof}
\begin{remark}\label{Rmk::Hold::Product::Hold2}
Using Definition \ref{Defn::Hold::CLogLandLip}, Lemma \ref{Lem::Hold::Product} \ref{Item::Hold::Product::Hold2} is still valid if $\alpha\in\{k+\LogL,k+\Lip:k=0,1,2,\dots\}$, namely $\|fg\|_{\Co^{k+\sigma}(\Omega;\Xs)}\le\|f\|_{L^\infty(\Omega;\Xs)}\|g\|_{\Co^{k+\sigma}(\Omega;\Xs)}+\|f\|_{\Co^{k+\sigma}(\Omega;\Xs)}\|g\|_{L^\infty(\Omega;\Xs)}$ holds for every $\sigma\in\{\LogL,\Lip\}$, $k\ge0$, $\Omega\subset\R^n$ open, and every Banach algebra $\Xs$. We leave the details to reader.
\end{remark}

In application we need the result on matrix-valued functions. Recall the matrix norm from \eqref{Eqn::Intro::MatrixNorms}.
\begin{lem}\label{Lem::Hold::CramerMixed}
Let $m,n\ge1$, let $\Omega\subseteq\R^n$ be an open subset, and let $\Xs(\Omega;\C^{m\times m})\subseteq L^\infty(\Omega;\C^{m\times m})$ be a Banach space of matrix functions that is closed under product and contains the identity matrix. 
\begin{enumerate}[label=(\roman*)]
    \item\label{Item::Hold::CramerMixed::Prod1} Let $C_1>0$ be the constant such that $\|fg\|_\Xs\le C_1\|f\|_\Xs\|g\|_\Xs$. 

Then the map of inverting the matrices
\begin{equation}\label{Eqn::Hold::CramerMixed::InvertingMatrix1}
    [A\mapsto A^{-1}]:\{A\in \Xs(\Omega;\C^{m\times m}):\|A-I_m\|_{\Xs(\Omega;\C^{m\times m})}<\tfrac1{2C_1}\}\to\Xs (\Omega;\C^{m\times m}),
\end{equation} 
is real-analytic and in particular Fr\'echet differentiable. Moreover $\|A^{-1}-I_m\|_{\Xs}\le2\|A-I_m\|_\Xs$ for all $A\in\Xs(\Omega;\C^{m\times m})$ such that $\|A-I_m\|_{\Xs(\Omega;\C^{m\times m})}<\tfrac1{2C_1}$.

In particular for any $\alpha\in(0,\infty)\cup\{k+\LogL,k+\Lip:k=0,1,2,\dots\}$ there exists a $c_0=c_0(\Omega,\alpha,m)>0$ such that $\|A^{-1}-I_m\|_{\Co^\alpha}\le 2\|A-I_m\|_{\Co^\alpha}$, for all $A\in \Co^\alpha(\Omega;\C^{m\times m})$ such that $ \|A-I_m\|_{\Co^\alpha}\le c_0$.

\item\label{Item::Hold::CramerMixed::Prod2}Suppose there is a $C_2>0$ such that $\|fg\|_\Xs\le C_2(\|f\|_\Xs\|g\|_{L^\infty}+\|f\|_{L^\infty}\|g\|_\Xs)$. Then \eqref{Eqn::Hold::CramerMixed::InvertingMatrix1} can be upgraded to the following map, which is also real-analytic (in particular Fr\'echet differentiable).
\begin{equation}\label{Eqn::Hold::CramerMixed::InvertingMatrix2}
    [A\mapsto A^{-1}]:\{A\in \Xs(\Omega;\C^{m\times m}):\|A-I_m\|_{L^\infty(\Omega;\C^{m\times m})}<\tfrac1{2C_2}\}\to\Xs (\Omega;\C^{m\times m}).
\end{equation} 
Moreover $\|A^{-1}-I_m\|_{\Xs}\le4\|A-I_m\|_\Xs$ for all $A\in\Xs(\Omega;\C^{m\times m})$ such that $\|A-I_m\|_{L^\infty(\Omega;\C^{m\times m})}<\tfrac1{2C_2}$.

\end{enumerate}
\end{lem}
\begin{proof}
\ref{Item::Hold::CramerMixed::Prod1}: By assumption $\|(A-I_m)^k\|_\Xs=\|(A-I_m)^{k-1}(A-I_m)\|_\Xs\le C_1\|(A-I_m)^{k-1}\|_\Xs\|A-I_m\|_\Xs$ for each $k\ge1$. Thus $\|(A-I_m)^k\|_\Xs\le C_1^{k-1}\|A-I_m\|_\Xs$ for all $k\ge1$. Using power series $A^{-1}=(I_m-(I_m-A))^{-1}=\sum_{k=0}^\infty(I_m-A)^k$ and applying the assumption $\|A-I\|_\Xs\le1/(2C_1)$, we have
\begin{equation}\label{Eqn::Hold::CramerMixed::Tmp}
    \|A^{-1}-I_m\|_\Xs\le\sum_{k=1}^\infty\|(I_m-A)^k\|_\Xs\le\sum_{k=1}^\infty C_1^{k-1}\|I_m-A\|_\Xs^k\le\|A-I_m\|_\Xs\sum_{k=1}^\infty\frac1{2^{k-1}}=2\|A-I_m\|_\Xs.
\end{equation}
In particular and the power series converges in $\Xs$, and $A^{-1}\in\Xs$ holds. Therefore the map $[A\mapsto A^{-1}]$ is real-analytic.

By Lemma \ref{Lem::Hold::ProdforExt} $\Co^\alpha(\Omega)$ is closed under product, taking $\Xs=\Co^\alpha$ the estimate applies and we get the constant $c_0$.

\medskip\noindent\ref{Item::Hold::CramerMixed::Prod2}: Now $\|(A-I_m)^k\|_\Xs\le C_2(\|(A-I_m)^{k-1}\|_{L^\infty}\|A-I_m\|_\Xs+\|(A-I_m)^{k-1}\|_\Xs\|A-I_m\|_{L^\infty})$, we get, $$\|(A-I_m)^k\|_\Xs\le C_2\|(A-I_m)\|_{L^\infty}^{k-1}\|A-I_m\|_\Xs+C_2\|(A-I_m)^{k-1}\|_\Xs\|A-I_m\|_{L^\infty},\quad k=1,2,3,\dots.$$
Using induction on $k$ we get 
\begin{equation*}
    \|(A-I_m)^k\|_\Xs\le kC_2^{k-1}\|A-I_m\|_{L^\infty}^{k-1}\|A-I_m\|_\Xs,\quad k=1,2,3,\dots.
\end{equation*}
Therefore when $\|A-I\|_{L^\infty}<1/(2C_2)$, similar to \eqref{Eqn::Hold::CramerMixed::Tmp}, we have
\begin{equation*}
    \|A^{-1}-I_m\|_\Xs\le\sum_{k=1}^\infty kC_2^{k-1}\|I_m-A\|_{L^\infty}^{k-1}\|I_m-A\|_\Xs\le\|A-I_m\|_\Xs\sum_{k=1}^\infty\frac k{2^{k-1}}=4\|A-I_m\|_\Xs.
\end{equation*}
In particular $A^{-1}$ exists and the inverting matrix map \eqref{Eqn::Hold::CramerMixed::InvertingMatrix2} is real-analytic.
\end{proof}

\subsection{The compositions}\label{Section::SingHoldComp}

The results for compositions is a bit complicated especially for indices smaller than 1.
\begin{defn}\label{Defn::Hold::CompIndex}
Let $\alpha,\beta\in\R_\Eb^+$. We define $\alpha\circ\beta$ to be largest index $\gamma\in\R_\Eb^+$ such that $[(f,g)\mapsto f\circ g]:\Co^\alpha_\loc(\R^m)\times\Co^\beta_\loc(\R^n;\R^m)\to\Co^\gamma_\loc(\R^n)$ holds for all integer $m,n\ge1$.
\end{defn}

In order to let $\alpha\circ\beta$ make sense, we need to show that the set $\{\gamma\in\R_\Eb^+:[(f,g)\mapsto f\circ g]:\Co^\alpha_\loc(\R^m)\times\Co^\beta_\loc(\R^n;\R^m)\to\Co^\gamma_\loc(\R^n),\ \forall m,n\ge1\}\subset\R_\Eb^+$ always has a maximum in $\R_\Eb^+$. This is guaranteed by Corollary \ref{Cor::Hold::CompOp}.

When $\alpha,\beta\in\R_\Eb$ is such that $\Co^{\alpha\circ\beta}$ is a Banach space (rather than a general Fr\'echet space or locally convex space, namely $\alpha\circ\beta\notin\{\gamma-:\gamma>0\}\cup\{\infty\}$), we have the following quantitative estimates.
\begin{prop}\label{Prop::Hold::QComp}
Let $\alpha,\beta\in\R_+\cup\{k+\LogL,k+\Lip:k=0,1,2,\dots\}$ and let $m,n\ge1$ be two integers. For any $M>1$ there is a $K_1=K_1(m,n,\alpha,\beta,M)>0$ such that for every $g\in\Co^{\alpha}(\B^n;\R^m)$ satisfies $g(\B^n)\subseteq \B^m$ and $\|g\|_{\Co^{\alpha}(\B^n;\R^m)}\le M$, the following holds
\begin{enumerate}[parsep=-0.3ex,label=(\roman*)]
    \item\label{Item::Hold::QComp::>1} For $\max(\alpha,\beta)\ge\Lip$ and $\{\alpha,\beta\}\neq\{1,\Lip\} $, for every $f\in\Co^\beta(\B^m)$, we have $f\circ g\in \Co^{\min(\alpha,\beta)}(\B^n)$ with $\|f\circ g\|_{\Co^{\min(\alpha,\beta)}(\B^n)}\le K_1\|f\|_{\Co^\beta(\B^m)}$.
    
    \item\label{Item::Hold::QComp::<1} For $\alpha,\beta\in(0,1)$, for every $f\in\Co^\beta(\B^m)$, we have $f\circ g\in \Co^{\alpha\beta}(\B^n)$ with $\|f\circ g\|_{\Co^{\alpha\beta}(\B^n)}\le K_1\|f\|_{\Co^\beta(\B^m)}$.
\end{enumerate}
\end{prop}

\begin{proof}See \cite[Lemma 5.8]{CoordAdaptedII} for the case $\alpha,\beta\in\R$, $\alpha>1$. For the others we start with the case $0<\alpha\le1<\beta$.

By Lemma \ref{Lem::Hold::HoldChar} \ref{Item::Hold::HoldChar::02} we have $\| f\circ g\|_{\Co^\alpha(\B^n)}\approx_{n,\alpha}\| f\circ g\|_{C^0(\B^n)}+\sup\limits_{x_0\neq x_1}|\tfrac{ f( g(x_0))+ f( g(x_1))}2- f( g(\tfrac{x_0+x_1}2))|\cdot|x_0-x_1|^{-\alpha}$ and $|\tfrac{ f( g(x_0))+ f( g(x_1))}2- f(\tfrac{ g(x_0)+ g(x_1)}2)|\lesssim_{m,\beta}\| f\|_{\Co^{\min(\beta,\frac32)}}|g(x_0)-g(x_1)|^{\min(\beta,\frac32)}$ for $x_0,x_1\in\B^n$. Therefore
\begin{align*}
    &|\tfrac{ f( g(x_0))+ f( g(x_1))}2- f( g(\tfrac{x_0+x_1}2))|\le|\tfrac{ f( g(x_0))+ f( g(x_1))}2- f(\tfrac{ g(x_0)+ g(x_1)}2)|+| f(\tfrac{ g(x_0)+ g(x_1)}2)- f( g(\tfrac{x_0+x_1}2))|
    \\
    \lesssim&_{m,\beta}\| f\|_{\Co^{\min(\beta,\frac32)}}| g(x_0)- g(x_1)|^{\min(\beta,\frac32)}+\| f\|_{C^{0,1}}|\tfrac{ g(x_0)+ g(x_1)}2- g(\tfrac{x_0+x_1}2)|
    \\
    \lesssim&_{n,\alpha}\| f\|_{\Co^{\min(\beta,\frac32)}}\| g\|_{\Co^{\alpha/\min(\beta,\frac32)}}^{\min(\beta,\frac32)}|x_0-x_1|^{\alpha}+\| f\|_{C^{0,1}}\| g\|_{\Co^\alpha}|x_0-x_1|^\beta\lesssim_M\| f\|_{\Co^{\beta}}|x_0-x_1|^{\alpha}.
\end{align*}

In particular we can choose $K_1(m,n,\alpha,\beta,M)\approx_{m,n,\alpha,\beta} (1+M)^{\min(\beta,\frac32)}$.

Now we have $\Co^\alpha\circ\Co^\beta\subseteq\Co^{\min(\alpha,\beta)}$ for all $\alpha,\beta\in\R_+$ such that $\max(\alpha,\beta)>1$. When $\max(\alpha,\beta)\in\{k+\LogL,k+\Lip:k\ge1\}$, in this case  $\min(\alpha,\beta)\in(0,1)$ is a real number and we have the embedding $\Co^{\max(\alpha,\beta)}\hookrightarrow\Co^\frac32$, thus by enlarging $K_1$, the estimate $\|f\circ g\|_{\Co^{\min(\alpha,\beta)}}\le K_1(\frac32)\|f\|_{\Co^{3/2}}\le K_1(\beta)\|f\|_{\Co^\beta}$. Therefore we get \ref{Item::Hold::QComp::>1} for all $\max(\alpha,\beta)>1$ and $\min(\alpha,\beta)\in\R_+$.

We postpone the proof of $\min(\alpha,\beta)\in\{k+\LogL,k+\Lip:k=1,2,\dots\}$ to the end. 

When $\min(\alpha,\beta)\le \Lip$ and $\max(\alpha,\beta)\ge\Lip$, using the embedding $\Co^{\max(\alpha,\beta)}\hookrightarrow\Co^\Lip$, it suffices to consider the case $\max(\alpha,\beta)=\Lip$. 

We are about the prove the case where $\max(\alpha,\beta)=\Lip$, $\min(\alpha,\beta)\in(0,1)\cup\{\LogL,\Lip\}$. To unify the argument we consider general modulus of continuity, which is a strictly increasing map $\mu:[0,1]\to[0,1]$ such that $\mu(0)=0$. For an open set $\Omega$ we denote by $\Co^\mu(\Omega)$ the space of continuous function $f$ such that 
\begin{equation}\label{Eqn::Hold::CMuSpace}
    \|f\|_{\Co^\mu(\Omega)}:=\|f\|_{C^0(\Omega)}+\sup_{0<|x-y|<1}\frac{|f(x)-f(y)|}{\mu(|x-y|)}.
\end{equation}
Clearly for modulus
\begin{equation*}
    \mu_\alpha(t):=t^\alpha,\quad 0<\alpha\le1;\qquad\mu_{log}(t):=t\log(e/t),
\end{equation*}
we have $\Co^{\mu_\alpha}=\Co^\alpha$ when $0<\alpha<1$, $\Co^{\mu_1}=\Co^\Lip$ and $\Co^{\mu_{\log}}=\Co^{\LogL}$, all of whose corresponding norms are equivalent. Note that $\mu_\alpha,\mu_{log}$ are all satisfy doubling condition, in particular $\sup_{0<t<1/M}\mu(Mt)/\mu(t)<\infty$ for $\mu=\mu_\alpha,\mu_{log}$ and for all $M>1$.

Let $\mu,\nu$ be two modulus of continuity, the composition $\mu\circ\nu:[0,1]\to[0,1]$ is also a modulus of continuity. Suppose $\sup_{0<t<1/M}\mu(Mt)/\mu(t)<\infty$. Then for every $f\in\Co^\mu(\B^m)$ and $g\in\Co^\nu(\B^n;\B^m)$ such that $\|g\|_{\Co^\nu}<M$, we have
\begin{align*}
    &\|f\circ g\|_{\Co^{\mu\circ\nu}(\B^n)}=\|f\circ g\|_{C^0(\B^n)}+\sup_{0<|x_0-x_1|<1}\frac{|f(g(x_0))-f(g(x_1))|}{\mu\circ\nu(|x_0-x_1|)}
    \\
    \le&\|f\|_{C^0(\B^m)}+\|f\|_{\Co^\mu(\B^m)}\sup_{0<|x_0-x_1|,|g(x_0)-g(x_1)|<1}\frac{\mu(|g(x_0)-g(x_1)|)}{\mu\circ\nu(|x_0-x_1|)}+\sup_{0<|x_0-x_1|<1\le |g(x_0)-g(x_1)|}\frac{2\|f\|_{C^0(\B^m)}}{\mu\circ\nu(|x_0-x_1|)}
    \\
    \le& \|f\|_{\Co^\mu(\B^m)}\big(1+\sup_{0<t<1/M}\mu(Mt)/\mu(t)+2\cdot \mu(1/M)^{-1}\big)\lesssim_{\mu,M} \|f\|_{\Co^\mu(\B^m)}.
\end{align*}
In particular we can choose $K_1=1+3M^\alpha$ if $\mu(t)=\mu_\alpha(t)=t^\alpha$ and $K_1=1+3M\log M$ if $\mu(t)=\mu_{log}(t)=t\log(e/t)$.

Take $\mu,\nu\in\{\mu_\alpha,\mu_{log}:0<\alpha\le1\}$, we get \ref{Item::Hold::QComp::>1} for $\min(\alpha,\beta)\in(0,1)\cup\{\LogL,\Lip\}$ and $\max(\alpha,\beta)=\Lip$, and \ref{Item::Hold::QComp::<1}.

Finally we prove the case $\min(\alpha,\beta)\in\{k+\LogL,k+\Lip:k=1,2,\dots\}$, using the embedding $\Co^{\max(\alpha,\beta)}\hookrightarrow\Co^{\min(\alpha,\beta)}$, it suffices to consider $\alpha=\beta\in\{k+\LogL,k+\Lip:k=1,2,\dots\}$.

Now $\Co^\alpha\subset C^1$, by chain rule $\nabla(f\circ g)=((\nabla f)\circ g)\cdot\nabla g$. When $k=1$ and $\sigma\in\{\LogL,\Lip\}$, for the proven case $\Co^\sigma\circ\Co^\Lip\subseteq\Co^\sigma$ we have
\begin{equation*}
    \|f\circ g\|_{\Co^{1+\sigma}(\B^n)}=\|f\circ g\|_{\Co^\sigma}+\|\nabla (f\circ g)\|_{\Co^\sigma}\le K_1(\|g\|_{\Co^\Lip})(\|f\|_{\Co^\sigma}+\|\nabla f\|_{\Co^\sigma}\|\nabla g\|_{\Co^\sigma}).
\end{equation*}
Since $\|g\|_{\Co^\Lip}+\|\nabla g\|_{\Co^\sigma}\lesssim\|g\|_{\Co^{1+\sigma}}$, by tracking the constant we get the constant $K_1$.

The case $k\ge2$ is done by the induction along with the fact from Lemma \ref{Lem::Hold::ProdforExt} that $\|f\circ g\|_{\Co^{k+\sigma}}\le\|(\nabla f)\circ g\|_{\Co^{k-1-\sigma}}\|\nabla g\|_{\Co^{k-1-\sigma}}$.
\end{proof}

The following shows the complete classification of $\alpha\circ\beta$.
\begin{cor}\label{Cor::Hold::CompOp}
Let $\alpha,\beta\in\R_\Eb^+$ be two generalized indices. The operator $\alpha\circ\beta$ is given by the following:
\begin{enumerate}[nolistsep,label=(\roman*)]
    \item\label{Item::Hold::CompOp::Lip} When $\max(\alpha,\beta)\ge\Lip$ and $\{\alpha,\beta\}\neq\{1,\Lip\}$, we have $\alpha\circ\beta=\max(\alpha,\beta)$.
    \item\label{Item::Hold::CompOp::<1} When $\max(\alpha,\beta)<1$, $\alpha\circ\beta=\alpha\beta$ in the sense of Convention \ref{Conv::Hold::ExtIndOp}. 
    \item\label{Item::Hold::CompOp::Log1} When $\max(\alpha,\beta)\in\{1-,\LogL,1\}$, $\alpha\circ\beta=\min(\alpha-,\beta-)$.
    \item\label{Item::Hold::CompOp::Log2} $1\circ\Lip=\Lip\circ 1=\LogL$.
\end{enumerate}
Moreover, let $U_1\subseteq\R^n$ and $V_2,U_2\subseteq\R^m$ be open subsets such that $V_2\Subset U_2$ is precompact, let $\Us\subset\Co^\alpha_\loc(U_2)$ and $\Vs\subset\Co^\beta_\loc(U_1;V_2)$ be two bounded subsets with respect to their topologies, then $\{f\circ g:f\in\Us,g\in\Vs\}\subset\Co^{\alpha\circ\beta}_\loc(U_1)$ is also a bounded subset.
\end{cor}
Here $\Co^\beta_\loc(U_1;U_2)$ is the subset of $g\in\Co^\beta_\loc(U_1;\R^m)$ such that $g(U_1)\subseteq U_2$.
\begin{remark}\label{Rmk::Hold::CompOp::Prod}
In Corollary \ref{Cor::Hold::CompOp} \ref{Item::Hold::CompOp::<1},  we formally say $\alpha\circ\beta=\alpha\beta$ by setting $(\gamma-)\cdot\delta=\gamma\cdot(\delta-)=(\gamma-)\cdot(\delta-):=(\gamma\delta)-$  for real numbers $\gamma,\delta$, following from Convention \ref{Conv::Hold::ExtIndOp}. 
\end{remark}
\begin{remark}
    In our setting for a bounded set $\Vs\subset\Co^\beta_\loc(U_1;V_2)$, the closure $\overline{\Vs}\subset\Co^\beta_\loc(U_1;\R^m)$ may not be contained in $\Co^\beta_\loc(U_1;V_2)$. 
    
    Let $\chi\in C_c^\infty(-2,2)$ be such that $\chi|_{[-1,1]}\equiv1$. Then for any $\beta>0$, $\{c\cdot\chi:-1<c<1\}\subset\Co^\beta_\loc(\R;(-1,1))$ is a bounded set, but its closure $\{c\cdot\chi:-1\le c\le 1\}\not\subset\Co^\beta_\loc(\R;(-1,1))$. 
    
    Note that under the usual subset topology, $\Co^\beta_\loc(U_1;\overline{V_2})$ is a completed topological vector space, but $\Co^\beta_\loc(U_1;V_2)$ is completed unless $V_2=\R^m$.
\end{remark}
\begin{proof}[Proof of Corollary \ref{Cor::Hold::CompOp}]
Clearly $\alpha\circ\beta\le\min(\alpha,\beta)$ whenever the left hand side is defined, because for $f\in\Co^\alpha(\R^n)$ and $g\in\Co^\beta(\R^n;\R^m)$ we have $f\circ\id=f$ and $\id\circ g=g$ where $\id\in\Co^\infty_\loc(\R)$. Taking $f$ such that $f\notin\Co^\sigma_\loc(\R)$ for all $\sigma<\alpha$ we see that $\alpha\circ\beta\le\alpha$, similar argument for $g$ shows $\alpha\circ\beta\le\beta$.

By passing to local, Proposition \ref{Prop::Hold::QComp} shows that when $\alpha,\beta\in(0,\infty)\cup\{k+\LogL,k+\Lip:k=0,1,2,\dots\}$, \begin{enumerate}[parsep=-0.3ex,label=(\alph*)]
    \item\label{Item::Hold::CompOp:Pf::PfLip} $\Co^\alpha_\loc(\R^m)\circ\Co^\beta_\loc(\R^n;\R^m)\subseteq\Co^{\min(\alpha,\beta)}_\loc(\R^n)$ for all $\max(\alpha,\beta)\ge\Lip$ and $\{\alpha,\beta\}\neq\{1,\Lip\}$;
    \item\label{Item::Hold::CompOp:Pf::Pf<1} $\Co^\alpha_\loc(\R^m)\circ\Co^\beta_\loc(\R^n;\R^m)\subseteq\Co^{\alpha\beta}_\loc(\R^n)$ for all $\alpha,\beta\in(0,1)$.
\end{enumerate}
Moreover the above two composition maps map every bounded subset of $\Co^\alpha_\loc(\R^m)\times\Co^\beta_\loc(\R^n;\R^m)$ to bounded subset of $\Co^{\min(\alpha,\beta)}_\loc(\R^n)$ or $\Co^{\alpha\beta}_\loc(\R^n)$, respectively.

When $\alpha,\beta$ satisfy the assumptions of \ref{Item::Hold::CompOp::Lip}, and if $\min(\alpha,\beta)\in\{\gamma-:\gamma\in\R\}$, then taking $\alpha_1,\beta_1\in(0,\infty)\cup\{k+\LogL,k+\Lip:k=0,1,2,\dots\}$ such that $\Lip\le\max(\alpha_1,\beta_1)\le\max(\frac32,\max(\alpha,\beta)) $ and $\min(\alpha_1,\beta_1)\in\{\gamma+\eps,\gamma-\eps\}$ for number $\eps>0$, we see that $\Co^{\alpha_1}_\loc\circ\Co^{\beta_1}_\loc\subseteq\Co^{\min(\alpha_1,\beta_1)}_\loc$ from \ref{Item::Hold::CompOp:Pf::PfLip}. Letting $\eps\to0$ we see that $\Co^\alpha_\loc(\R^m)\circ\Co^\beta_\loc(\R^n;\R^m)\subseteq\Co^{\min(\alpha,\beta)}_\loc(\R^n)$ for all extended indices such that $\max(\alpha,\beta)\ge\Lip$ and $\{\alpha,\beta\}\neq\{1,\Lip\}$. Moreover the composition map bounded subset to bounded subset as well. 

Therefore $\min(\alpha,\beta)\le\alpha\circ\beta\le\min(\alpha,\beta)$, which give the proof of \ref{Item::Hold::CompOp::Lip}.

When $\alpha,\beta\in(0,1)$, since $|t|^\alpha\circ|t|^\beta=|t|^{\alpha\beta}\notin\Co^{\alpha\beta-}$ we see that $\alpha\circ\beta\le\alpha\beta$ if the left hand side is defined. By \ref{Item::Hold::CompOp::<1} we see that $\alpha\beta\le\alpha\circ\beta$, thus $\alpha\circ\beta=\alpha\beta$. The standard limit argument shows that $\alpha\circ\beta=\alpha\beta$ holds for $\alpha$ or $\beta\in\{\gamma-:\gamma\in(0,1)\}$, in the sense of Remark \ref{Rmk::Hold::CompOp::Prod}. This proves \ref{Item::Hold::CompOp::<1}. Note that in this case the composition map $\Co^\alpha_\loc(\R^m)\times\Co^\beta_\loc(\R^n;\R^m)\to\Co^{\alpha\beta}_\loc(\R^n)$ still preserve boundedness of the subsets.

By Lemma \ref{Lem::Hold::ZygExample} $f_1(t):=t\log|t|$ satisfies $f_1\in\Co^1\backslash\Co^\Lip$, but $f_1\circ f_1(t)=t\log^2|t|+t(\log|t|+\log|\log|t||)$ is not log-Lipschitz. Thus $1\circ 1\le1-$. And for $\gamma\in(0,1)$, $|t|^\gamma\circ f_1(t)=|t|^\gamma|\log|t||^\gamma$, $f_1(|t|^\gamma)=|t|^\gamma\gamma\log|t|$, neither of which belongs to $\Co^\gamma$. Thus $\gamma\circ1,1\circ\gamma\le\gamma-$ for all $\gamma\in(0,1)$.

On the other hand by \ref{Item::Hold::CompOp:Pf::Pf<1} we see that $\Co^{\gamma-}\circ\Co^{1-}\subseteq\Co^{\gamma-}$ and $\Co^{1-}\circ\Co^{\gamma}\subseteq\Co^{\gamma-}$ for all $\gamma\in(0,1)$, thus for $\max(\alpha,\beta)\in\{1-,\LogL,1\}$ we have $\min(\alpha-,\beta-)\le\alpha\circ\beta\le\min(\alpha-,\beta-)$, which gives \ref{Item::Hold::CompOp::Log1}. The boundedness of the map follows from \ref{Item::Hold::CompOp:Pf::Pf<1} along with the limit argument.

Finally the function $g_1(t):=\max(t,0)$ is locally Lipschitz, but $g_1\circ f_1(t)=f_1\circ g_1(t)=\max(t,0)\cdot\log|t|$ is log-Lipschitz but not Zygmund-1 ($\Co^1$). Thus $1\circ\Lip,\Lip\circ 1\le\LogL$. By \ref{Item::Hold::CompOp:Pf::PfLip} we see that $\LogL\circ\Lip=\Lip\circ\LogL=\LogL$, therefore  $\LogL\le 1\circ\Lip\le\LogL$ and $\LogL\le\Lip\circ 1\le\LogL$, which gives \ref{Item::Hold::CompOp::Log2}. The boundedness of the map follows from \ref{Item::Hold::CompOp:Pf::PfLip}.

Now we have completely classified the operators $\alpha\circ\beta$ for all positive $\alpha,\beta\in\R_\Eb^+$. Moreover we know the composition map $\Co^\alpha_\loc(\R^m)\times\Co^\beta_\loc(\R^n;\R^m)\subseteq\Co^{\alpha\circ\beta}_\loc(\R^n)$ preserves boundedness of their subsets.

Finally for every open sets $U_1\subseteq\R^n$ and $U_2\subseteq\R^m$, and for any precompact $V_1\Subset U_2$, take $\chi_1\in C_c^\infty(U_1)$ and $\chi_2\in C_c^\infty(U_2)$ be such that $\chi_1|_{V_1}\equiv1$, $\chi_2|_{V_2}\equiv1$, we see that $f\circ g|_{V_1}=(\chi_2f)\circ(\chi_1g)|_{V_1}$ for all $f:U_2\to\R$ and $g:U_1\to V_2$. On the other hand if $[f\mapsto \chi_2 f]:\Co^\alpha_\loc(U_2)\to\Co^\alpha(\R^m)$ and $[g\mapsto\chi_1g]:\Co^\beta_\loc(U_1;U_2)\to\Co^\beta(\R^n;\R^m)$ are continuous linear maps which preserve boundedness. Therefore by the boundedness of the composition map $\Co^\alpha_\loc(\R^m)\circ\Co^\beta_\loc(\R^n;\R^m)\to\Co^{\alpha\circ\beta}_\loc(\R^n)$ we see that $\{(\chi_2f)\circ(\chi_1g):f\in\Us,g\in\Vs\}\subset\Co^{\alpha\circ\beta}_\loc(\R^n)$ is bounded. Since $V_1$ is arbitrary, we conclude that $\{f\circ g:f\in\Us,g\in\Vs\}\subset\Co^{\alpha\circ\beta}_\loc(U_1)$ is a bounded subset.
\end{proof}

If the second element of the composition is a diffeomorphisms, we are allow to take the first element as distributions. This is useful in defining distributions on manifolds.
\begin{lem}\label{Lem::Hold::PushForwardFuncSpaces}
Let $\kappa\in(1,\infty]$, $\alpha\in(0,\kappa]$ and $\beta\in(\min(1-\kappa,-\alpha),\alpha]$. Let $U,V\subseteq\R^n$ be two open sets and let $\varphi:U\to V$ be a $\Co^\kappa$-diffeomorphism (i.e. $\varphi,\varphi^\Inv$ are both $\Co^\kappa_\loc$-maps). Let $\rho\in \Co^\alpha_\loc(U)$.

If $f\in \Co^\beta_\loc(V)$ then $\rho\cdot(f\circ\varphi)\in \Co^\beta_\loc(U)$ as well.
\end{lem}
\begin{remark}
    Here $f\circ\varphi$ is a well-defined distribution in the sense that if $\lim_{\sigma\to\infty}f_\sigma=f$ in $\Co^{\beta-}_\loc(V)$ then $\lim_{\sigma\to\infty}f_\sigma\circ\varphi=f\circ\varphi$ in $\Co^{\beta-}_\loc(U)$ holds. See also \cite[Chapter 2.10]{Triebel1}. 
\end{remark}
\begin{proof}
We only need to prove the case $\alpha=\kappa$. Indeed in this case gives the property $f\circ\varphi\in\Co^\beta_\loc(U)$. Since $\alpha+\beta>0$, by Lemma \ref{Lem::Hold::MultLoc} we get $\rho(f\circ\varphi)\in\Co^\beta_\loc(U)$ holds.

We process by induction. We prove the result for $\beta>\max(1-\kappa,-l)$ with $l=0,1,2,\dots$. The base $l=0$ is automatically true since in this case $\beta>0$ and we can apply Corollary \ref{Cor::Hold::CompOp} \ref{Item::Hold::CompOp::Lip}. We assume the result for $l-1$ and prove it for $l$.

Now let $\beta>\max(1-\kappa,-l)$ and $f\in\Co^\beta_\loc(V)$.
By passing to a precompact subset and applying Lemma \ref{Lem::Hold::NegHold=SumGood}, we can write $f=g_0+\sum_{j=1}^n\partial_jg_j$ where $g_0,\dots,g_n\in \Co^{\beta+1}_\loc(V)$. By the inductive hypothesis, since $\beta+1>\max(1-\kappa,1-l)$, we have $g_j\circ\varphi\in \Co^{\beta+1}_\loc(U)$ for $j=0,\dots,n$.

Denote $\psi=(\psi^1,\dots,\psi^n):=\varphi^\Inv:V\to U$, which is a $\Co^\kappa$-map by assumption. Using chain rule and product rule,
\begin{equation*}
    \rho\cdot((\partial_jg_j)\circ\varphi)=\sum_{k=1}^n(\partial_k(g_j\circ\varphi))\cdot\rho\cdot((\partial_j\psi^k)\circ\varphi).
\end{equation*}
We have $g_j\circ\varphi\in \Co^{\beta+1}_\loc(U)$, by Remark \ref{Rmk::Hold::GradBdd} \ref{Item::Hold::GradBdd} we get $\partial_k(g_j\circ\varphi)\in \Co^{\beta}_\loc(U)$. Since $\rho((\partial_j\psi^k)\circ\varphi)\in \Co^{\kappa-1}_\loc(U)$.
Thus by Lemma \ref{Lem::Hold::MultLoc} \ref{Item::Hold::MultLoc::WellDef}, since $\beta\in (1-\kappa,0]$, we have the product $(\partial_k(g_j\circ\varphi))\cdot\rho((\partial_j\psi^k)\circ\varphi)\in \Co^{\beta-1}_\loc(U)$. Taking sum over $k=1,\dots,n$ we get $\rho\cdot((\partial_jg_j)\circ\varphi)\in \Co^\beta_\loc(U)$ for $j=1,\dots,n$.


Clearly $\rho\cdot(g_0\circ\varphi)\in \Co^{\beta+1}_\loc(U)\subset\Co^{\beta}_\loc(U)$. Taking sum over $j=0,\dots,n$ we get the result.
\end{proof}

To define $\Co^{\LogL-1}$-distribution on manifold (for Theorem \ref{MainThm::SingFro}), we need a similar result to Lemma \ref{Lem::Hold::PushForwardFuncSpaces} as well. But the proof is more complicated.
\begin{lem}\label{Lem::Hold::LogL-1Comp}
Let $\alpha>0$ and let $U,V\subseteq\R^n$ be two open sets. Let $f\in\Co^{\LogL-1}_\loc(V)$ and let $\Phi:U\xrightarrow{\sim}V$ be a $\Co^{\alpha+1}$-diffeomorphism. Then $f\circ\Phi\in\Co^{\LogL-1}_\loc(U)$ is defined.
\end{lem}
\begin{proof}Without loss of generality we assume $0<\alpha<1$. By passing to local we can assume $f\in\Co^{\LogL-1}_c(V)$. By shrinking $U$ and $V$ we can assume $\Phi$ and $\Phi^\Inv$ both have bounded $\Co^{\alpha+1}$-norms.

Let $(\phi_j)_{j=0}^\infty$ be a fixed dyadic resolution that defines the $\Co^{\LogL-1}$ norm. The integration by substitution formula yields
\begin{equation*}
    \phi_j\ast(f\circ\Phi)(x)=\int_U\phi_j(x-t)f(\Phi(t))dt=\int_V\phi_j(x-\Phi^\Inv(s))\cdot f(s)\cdot|\det \nabla(\Phi^\Inv)(s)|ds,\quad x\in\R^n,\quad j\ge0.
\end{equation*}
Since $\nabla\Phi^\Inv\in\Co^\alpha(V;\R^{n\times n})$ is pointwise invertible, by Lemma \ref{Lem::Hold::CharLogL-1} \ref{Item::Hold::CharLogL-1::Grad} we see that $g(s):=f(s)|\det\nabla\Phi^\Inv(s)|$ is a distribution in $\Co^{\LogL-1}_c(V)$. For $j\ge0$, let $h_j:=\Delta^{-1}(g-\Su^\phi_jg)=\sum_{k=j+1}^\infty(\Delta^{-1}\phi_k)\ast g$. Using $\phi_k\ast g=\phi_k\ast\Su^\phi_{k+1}g$ from \eqref{Eqn::Hold::LPCharComp} we have for $j\ge1$,
\begin{equation*}
    \|\nabla h_j\|_{L^\infty}\le \sum_{k=j+1}^\infty\|(\nabla\Delta^{-1}\phi_k)\ast\Su^\phi_{k+2}g\|_{L^\infty}\le\sum_{k=j+1}^\infty\|\nabla\Delta^{-1}\phi_k\|_{L^1}\|\Su^\phi_{k+2}g\|_{L^\infty}\lesssim_\phi\sum_{k=j+1}^\infty 2^{-k}k\|g\|_{\Co^{\LogL-1}}\lesssim_g j2^{-j}.
\end{equation*}

Let $\eta_j:=\sum_{k=0}^j\phi_k$, we have $\Su^\phi_j(f\circ\Phi)=\int_V\eta_j(x-\Phi^\Inv(s))\cdot g(s)ds$ for all $j\ge0$. Let $\chi\in C_c^\infty(V)$ be such that $\chi|_{\supp g}\equiv1$, then we have for $j\ge1$,
\begin{align*}
    &\|\Su^\phi_j(f\circ\Phi)\|_{L^\infty(\R^n)}\le\sup_{x\in\R^n}\Big|\int_V\chi(s)\eta_j(x-\Phi^\Inv(s))\cdot(\Su^\phi_jg(s)+\Delta h_j(s))ds\Big|
    \\
    \le&\sup_{x\in\R^n}\int_V|\chi(s)\eta_j(x-\Phi^\Inv(s))\cdot\Su^\phi_jg(s)|ds+\sup_{x\in\R^n}\Big|\int_V\nabla_s(\chi(s)\eta_j(x-\Phi^\Inv(s)))\cdot\nabla h_j(s))ds\Big|
    \\
    \le&\sup_{x\in\R^n}\|\eta_j(x-\Phi^\Inv(\cdot))\|_{L^1(V)}\|\Su^\phi_jg\|_{L^\infty}+\sup_{x\in\R^n}\|\eta_j(x-\Phi^\Inv(\cdot))\cdot\chi\|_{W^{1,1}(V)}\|\nabla h_j\|_{L^\infty}
    \\
    \le&\|\eta_j\|_{L^1}\|\det\nabla\Phi\|_{L^\infty(U)}\cdot j\|g\|_{\Co^{\LogL-1}}+\|\eta_j\|_{W^{1,1}}\|\det\nabla\Phi\|_{L^\infty(U)}\|\chi\|_{C^1}\|\nabla h_j\|_{L^\infty}
    \\
    \lesssim&_{\phi,\Phi,g}j+2^j\cdot(j2^{-j})\approx j.
\end{align*}
Taking supremum over $j\ge1$ we get $f\circ\Phi\in\Co^{\LogL-1}(\R^n)$ hence $f\circ\Phi\in\Co^{\LogL-1}_c(U)$. 
\end{proof}
\begin{remark}
    The same proof also work for the case $f\in\Co^\beta$ for $-1<\beta<0$.
\end{remark}

\subsection{Inverse Laplacian and inverse functions}
\begin{note}\label{Note::Hold::Codiff}
We use the co-differential, $\codiff=\codiff_{\R^n}$, which is a linear operator
taking $k$ forms to $k-1$ forms, satisfying for $1\leq i_1<i_2<\cdots<i_k\leq n$,
\begin{equation*}
    \codiff (f dx^{i_1}\wedge \cdots \wedge dx^{i_k})
    =\sum_{l=1}^k \frac{\partial f}{\partial x^{i_l}} (-1)^{l} dx^{i_1}\wedge \cdots \wedge dx^{i_{l-1}}\wedge dx^{i_{l+1}}\wedge \cdots \wedge dx^{i_k}.
\end{equation*}

\end{note}
In particular on 1-forms, $-\codiff$ is the divergence operator, namely
\begin{equation*}
    \text{For }\lambda=\sum_{i=1}^n\lambda_idx^i,\quad\codiff\lambda=-\sum_{i=1}^n \frac{\partial\lambda_i}{\partial x^i}.
\end{equation*}
For any form $\omega$, we have $d(\codiff\omega)+\codiff(d\omega)=-\Delta\omega$, where
$\Delta=\sum_{i=1}^n\partial_{x^i}^2$ is the classical Laplacian  acting on the components of $\omega$; in this
setting $\Delta$ is called the (negative) Hodge Laplacian.

\begin{note}\label{Note::Hold::Newtonian}
We use $\Ga$ as the Newtonian potential in $\R^n$, namely
\begin{equation}\label{Eqn::Hold::Newtonian}
    \Ga(x):=
    \begin{cases}\frac{|x|}2,&n=1,\\\frac1{2\pi}\log|x|,&n=2,\\-|\mathbb S^{n-1}|^{-1}|x|^{2-n},&n\ge3.\end{cases}
\end{equation}
\end{note}
We will often convolve functions with $k$-forms.
Formally, if $\omega=\sum_{1\leq i_1<\cdots<i_k\leq n} \omega_{i_1,\ldots, i_k}dx^{i_1}\wedge \cdots \wedge dx^{i_k}$ is a $k$-form, and $\phi$ is a function, we set $\phi*\omega=\sum_{1\leq i_1<\cdots<i_k\leq n} (\phi*\omega_{i_1,\ldots, i_k})dx^{i_1}\wedge \cdots \wedge dx^{i_k}$.

\begin{lem}\label{Lem::Hold::GreensOp}
Let $U\subset\R^n$ be a bounded open subset, then for every $\alpha\in\R$ the map $[f\mapsto(\Ga\ast f)|_U]:\Co^\alpha_c(U;\wedge^kT^*U)\to\Co^{\alpha+2}(U;\wedge^kT^*U)$ is continuous for all $0\le k\le n$. Moreover there is a $C=C(U,\alpha)>0$ such that $\|\Ga\ast f\|_{\Co^{\alpha+2}(U)}\le C\|f\|_{\Co^\alpha(U)}$ for all $f\in \Co^\alpha_c(U;\wedge^kT^*U)$.
\end{lem}
\begin{proof}By passing to the coordinate components it suffices to prove the case $k=0$.

Let $\chi\in \Sc(\R^n)$ whose Fourier transform is $C_c^\infty$. We write $\Ga_0:=\chi\ast\Ga$ and $\Ga_\infty:=(\delta_0-\chi)\ast\Ga$. Since $\Ga$  is well-known to be a fundamental solution for
    the Laplacian, we have the Fourier transform $((I-\Delta)\Ga_\infty)^\wedge(\xi)=(1-\hat\chi(\xi))(4\pi^2|\xi|^{-2}+1)$, which is a bounded smooth function. Therefore by H\"ormander-Mikhlin multiplier theorem (see \cite[Theorem 2.3.7]{Triebel1}), $[f\mapsto (I-\Delta)\Ga_\infty\ast f]:\Bs_{\infty\infty}^\alpha(\R^n)\to\Bs_{\infty\infty}^\alpha(\R^n)$. By \cite[Theorem 2.3.8]{Triebel1} we get $(I-\Delta)^{-2}:\Bs_{\infty\infty}^\alpha(\R^n)\to\Bs_{\infty\infty}^{\alpha+2}(\R^n)$.

Recall from Lemma \ref{Lem::Hold::HoldChar} and Definition \ref{Defn::Hold::NegHold} that $\Co^\alpha(\R^n)=\Bs_{\infty\infty}^\alpha(\R^n)$, so $[f\mapsto\Ga_\infty\ast f]:\Co^\alpha_c(U)\subset\Co^\alpha(\R^n)\to\Co^{\alpha+2}(\R^n)$ is bounded linear.

On the other hand $\supp\hat\Ga_0\subseteq\supp\hat\chi$ is compact, thus $\Ga_0\in \Co^\infty_\loc(\R^n)$. Clearly $[f\mapsto\Ga_0\ast f]:\mathscr E(\R^n)\to \Co^\infty_\loc(\R^n)$ is continuous and thus $[f\mapsto\Ga_0\ast f]:\Co^\alpha_c(U)\to \Co^{\alpha+2}_\loc(\R^n)$ is also bounded linear.

Now $[f\mapsto\Ga\ast f]:\Co^\alpha_c(U)\to\Co^{\alpha+2}_\loc(\R^n)$ is bounded linear, taking restriction to $U$ we get the continuity.

Let $\tilde U\subset\R^n$ be a larger bounded open subset such that $U\Subset\tilde U$ is precompact. We know the Banach subspace $\{f\in\Co^\alpha(\R^n):f|_{U^c}\equiv0\}\subset\Co^\alpha(\R^n)$ is a bounded subset of $\Co^\alpha_c(\tilde U)$. Since $[f\mapsto(\Ga\ast f)|_{\tilde U}]:\Co^\alpha_c(\tilde U)\to\Co^{\alpha+2}(\tilde U)$ is continuous, we know the following map is bounded:
\begin{equation*}
    [f\mapsto(\Ga\ast f)|_{U}]:\{f\in\Co^\alpha(\R^n):f|_{U^c}\equiv0\}\hookrightarrow\Co^\alpha_c(\tilde U)\to\Co^{\alpha+2}(\tilde U)\xrightarrow{(-)|_{U}}\Co^{\alpha+2}(U).
\end{equation*}
Thus such constant $C(U,\alpha)>0$ exists.
\end{proof}
We also consider  the Dirichlet problem of the Laplacian equation:
\begin{lem}[The Dirichlet Problem]\label{Lem::Hold::DiriSol}
    Let $\alpha>0$, and let $U$ be a bounded domain with smooth boundary. Then for $f\in\Co^{\alpha-2}(U)$ and $g\in\Co^\alpha(\partial U)$ there is a unique $u\in\Co^{\alpha}(U)$ such that $\Delta u=f$ and $u\big|_{\partial U}=g$.
    
    Moreover, the solution map $(f,g)\mapsto u$ is continuous linear map $\Co^{\alpha-2}(U)\oplus\Co^\alpha(\partial U)\rightarrow \Co^{\alpha}(U)$.
\end{lem}
See \cite[Theorem 15]{DirichletBoundedness} for a proof of Lemma \ref{Lem::Hold::DiriSol}. 

\begin{remark}\label{Rmk::Hold::ZygmundonSphere}
In Lemma \ref{Lem::Hold::DiriSol}, for $\alpha>0$, $\Co^\alpha(\partial U)$ is the Zymgund-H\"older space on manifold. One can use Definition \ref{Defn::DisInv::DefFunVF} to define $\Co^\alpha$-functions on it.
Note that the sphere is a compact manifold, and therefore $\Co^\alpha_\loc(\partial U)=\Co^\alpha(\partial U)$.
The norm can be defined using any finite atlas,
and the equivalence class of the norm does not depend on 
the choice of atlas. Moreover we have
\begin{equation}\label{Eqn::Hold::HZNormforSphere}
    \|f\|_{\Co^\alpha(\partial U)}\approx_\alpha\inf\{\|\tilde f\|_{\Co^\alpha(U)}:\tilde f\big|_{\partial U}=f\},
\end{equation}
since the trace operator $(\cdot)|_{\partial U}:\Co^\alpha(U)\to\Co^\alpha(\partial U)$ is continuous and surjective; see \cite[Theorem 2.7.2]{Triebel1}.

In the case of unit ball $\B^n$ and its boundary $\Sp^{n-1}=\partial\B^n$, given a function $f\in\Co^\alpha(\Sp^{n-1})$, we can take $\tilde f(x)=f(\frac{x}{|x|})\chi(|x|)$ where $\chi\in C_c^\infty(-\frac12,2)$ such that $\chi(1)=1$, then we have $\tilde f\big|_{\Sp^{n-1}}=f$ and $\|\tilde f\|_{\Co^\alpha(\B^n)}\approx_\alpha\|f\|_{\Co^\alpha(\Sp^{n-1})}$.
\end{remark}

\begin{note}\label{Note::Hold::DiriSol}
Let $U\subset\R^n$ be a bounded smooth domain. We denote $\Pc_0=\Pc_0^U$ to be the Laplacian solution operator with zero Dirichlet boundary. Namely $\Pc_0$ is the map $(f,0)\mapsto u$ in Lemma \ref{Lem::Hold::DiriSol}.
\end{note}

In Propositions \ref{Prop::Rough1Form::ExistPDE} and \ref{Prop::Rough1Form::ExistPDE} we need the following quantitative version inverse function theorem and the pushforwards of the exact forms.
\begin{prop}\label{Prop::Hold::QIFT}
    Let $\alpha\in(0,\infty)$ and $\beta\in[\alpha,\alpha+1]$ be two real numbers, then there is a constant $K_2=K_2(n,\alpha,\beta)>1$ satisfying the following:
    
    Suppose $R\in\Co^{\alpha+1}(\B^n;\R^n)$ satisfies $R\big|_{\partial\B^n}=0$ and $\|R\|_{\Co^{\alpha+1}}\le K_2^{-1}$, then the map $F:=\id+R:\B^n\to\R^n$ is a $\Co^{\alpha+1}$-diffeomorphism of $\B^n$. Moreover,
    \begin{enumerate}[parsep=-0.3ex,label=(\roman*)]
        \item\label{Item::Hold::QIFT::PhiEst} Let $\Phi=(\phi^1,\dots,\phi^n):\B^n\xrightarrow{\sim}\B^n$ be the inverse map of $F$. Then $\|\Phi\|_{\Co^{\alpha+1}(\B^n;\R^n)}\le K_2$ and  $\|\nabla\Phi-I_n\|_{\Co^{\alpha}(\B^n;\R^{n\times n})}\le K_2\|R\|_{\Co^{\alpha+1}(\B^n;\R^n)}$.
        \item\label{Item::Hold::QIFT::dFormEst} If $\lambda\in\Co^\alpha_c(\B^n;T^*\B^n)$ satisfies $\supp\lambda\subsetneq\frac12\B^n$ and $d\lambda\in\Co^{\beta-1}(\B^n;\wedge^2T^*\B^n)$, then $\supp F_*\lambda\subsetneq\frac34\B^n$ and 
        \begin{equation}\label{Eqn::Hold::QIFT::dFormEst}
            \|d(F_*\lambda)\|_{\Co^{\beta-1}(\B^n;\wedge^2T^*\B^n)}\le K_2\|d\lambda\|_{\Co^{\beta-1}(\B^n;\wedge^2T^*\B^n)}.
        \end{equation}
    \end{enumerate}
\end{prop}
\begin{proof}
We let $K_2$ be a large constant which may change from line to line.
In particular, we will choose $K_2$ large enough such that $\|R\|_{\Co^{\alpha+1}}\le K_2^{-1}$ implies $\| R\|_{C^0}+\|\nabla R\|_{C^0}\le\frac14$ and $\|\nabla R\|_{\Co^\alpha}\le c_0(\B^n,\alpha,n)$, where $c_0$ is in Lemma \ref{Lem::Hold::CramerMixed} \ref{Item::Hold::CramerMixed::Prod1}.

By Lemma \ref{Lem::Hold::CramerMixed} \ref{Item::Hold::CramerMixed::Prod1}, $\nabla F(x)=I+\nabla R(x)$ is an invertible matrix for every $x\in\B^n$, and we have 
\begin{equation}\label{Eqn::Hold::QIFT::Proof1}
\|(\nabla F)^{-1}-I\|_{\Co^\alpha}=\|(I+\nabla R)^{-1}-I\|_{\Co^\alpha}\le 2\|\nabla R\|_{\Co^\alpha}.
 \end{equation}

Since $\|\nabla R\|_{C^0}\le\frac14$, we have $|R(x_1)-R(x_2)|\le\|\nabla R\|_{C^0}|x_1-x_2|\le\frac14|x_1-x_2|$, which implies 
\begin{equation}\label{Eqn::Hold::QIFT::Proof1.5}
    |F(x_1)-F(x_2)|\ge|x_1-x_2|-|R(x_1)-R(x_2)|\ge\textstyle\frac34|x_1-x_2|.
\end{equation}

This implies $F$ is injective. By the Inverse Function Theorem, we know $F:\B^n\xrightarrow{\sim} F(\B^n)$ is a $\Co^{\alpha+1}$-diffeomorphism.

The assumption $R\big|_{\partial\B^n}=0$ gives $F(\partial\B^n)=\partial\B^n$. Since $F(\B^n)$ is contractible and $\overline{F(\B^n)}\supset\partial\B^n$, we get that $F(\B^n)=\B^n$. We conclude $F$ is a $\Co^{\alpha+1}$-diffeomorphism on $\B^n$.

\medskip
\noindent Proof of \ref{Item::Hold::QIFT::PhiEst}: 
First, we claim that there is a $C_1(n,\alpha)>0$, which does not 
depend on $R$, such that whenever $R$ satisfies the assumptions
of the proposition, we have 
\begin{equation}\label{Eqn::Hold::QIFT::1::Tmp1}
\|\Phi \|_{\Co^{\alpha+1}(\B^n;\R^n)}\leq C_1.
\end{equation}
In particular $\|\Phi \|_{\Co^{\alpha+1}}\le K_2$ holds by taking $K_2\ge C_1$.

Since $\Phi$ is the inverse map of $F$, by \cite[Lemma 5.9]{CoordAdaptedII} or by induction on $\alpha$ with chain rule, we know $\|\Phi\|_{\Co^{\alpha+1}(\B^n;\R^n)}$ only depends on $n$, $\alpha$, $\|F\|_{\Co^{\alpha+1}(\B^n;\R^n)}$, and $\|(\nabla F)^{-1}\|_{C^0(\B^n;\R^{n\times n})}$. 
We will show that $\|F\|_{\Co^{\alpha+1}(\B^n;\R^n)}$ and $\|(\nabla F)^{-1}\|_{C^0(\B^n;\R^{n\times n})}$ have bounds that do not depend on $R$.

We have $\|F\|_{\Co^{\alpha+1}(\B^n;\R^n)}\le\|\id\|_{\Co^{\alpha+1}(\B^n)}+\|R\|_{\Co^{\alpha+1}(\B^n;\R^n)}\le \|\id\|_{\Co^{\alpha+1}(\B^n)}+\frac14$. The right hand side of this inequality does not depend on $R$.

By \eqref{Eqn::Hold::QIFT::Proof1.5}, $|F(x_1)-F(x_2)|\ge\frac34|x_1-x_2|$ implies $\sup\limits_{x\in\B^n}|(\nabla F(x))^{-1}|\le\frac43$, so $\|(\nabla F)^{-1}\|_{C^0(\B^n;\R^{n\times n})}\le\frac43$, which does not depend on $R$ as well.  This establishes \eqref{Eqn::Hold::QIFT::1::Tmp1}.

Note that the identity matrix $I$ can be viewed as a constant function defined on the unit ball. Since $\Phi$ is a $\Co^{\alpha+1}$-diffeomorphism on $\B^n$, we get the equality $I=I\circ\Phi$ as a matrix function on $\B^n$.

By the chain rule $I=\nabla(F\circ\Phi)=((\nabla F)\circ\Phi)\cdot\nabla\Phi$, so $\nabla\Phi-I=(\nabla F)^{-1}\circ\Phi-I=((\nabla F)^{-1}-I)\circ\Phi$.

By Proposition \ref{Prop::Hold::QComp} \ref{Item::Hold::QComp::>1}, we know  $\|G\circ\Phi\|_{\Co^\alpha(\B^n;\R^{n\times n})}\le K_1(n^2,n,\alpha,\alpha,C_1)\|G\|_{\Co^\alpha(\B^n;\R^{n\times n})}$ holds when \eqref{Eqn::Hold::QIFT::1::Tmp1} is satisfied, where $K_1$ is the constant in Proposition \ref{Prop::Hold::QComp}.
Therefore by taking $G=(\nabla F)^{-1}-I$, we see that
\begin{equation*}
    \| \nabla \Phi-I \|_{\Co^\alpha(\B^n;\R^{n\times n})} = \|((\nabla F)^{-1}-I)\circ \Phi\|_{\Co^\alpha} \leq K_1\|(\nabla F)^{-1}-I\|_{\Co^\alpha}
    \leq 2K_1 \| \nabla R\|_{\Co^{\alpha}(\B^n;\R^{n\times n})}\le  2K_1 \|R\|_{\Co^{\alpha+1}(\B^n;\R^n)},
\end{equation*}
and we obtain $\| \nabla \Phi-I \|_{\Co^\alpha}\le K_2\|R\|_{\Co^{\alpha+1}}$ by replacing $K_2$ with $\max(K_2,2K_1)$.


\medskip
\noindent Proof of \ref{Item::Hold::QIFT::dFormEst}: Let $\lambda\in\Co_c^\alpha(\B^n;T^*\B^n)$ be as in the assumption of \ref{Item::Hold::QIFT::dFormEst}. In particular, $\supp\lambda\subsetneq\frac12\B^n$.

By the assumption $\|R\|_{C^0}\le\frac14$, we have $F(\frac12\B^n)\subseteq\frac12\B^n+\frac14\B^n=\frac34\B^n$, so $\supp F_*\lambda=F(\supp\lambda)\subseteq\frac34\B^n$.

We define a 1-form $\rho$ and a function $\xi$ by
$$\rho=\sum_{i=1}^n\rho_idx^i:=-\Ga\ast \codiff d\lambda,\quad \xi:=-\Ga\ast \codiff\lambda.$$
$\rho$ and $\xi$ are globally defined in $\R^n$ because $\lambda$ is compactly supported in $\frac12\B^n$. 

By Lemma \ref{Lem::Hold::GreensOp}, we have $\rho\in\Co^\beta_\loc(\R^n;T^*\R^n)$, $\xi\in\Co^{\alpha+1}_\loc(\R^n)$ and  
\begin{equation}\label{Eqn::Hold::QIFT::RhoisBoundedbyDTheta}
    \|\rho\|_{\Co^\beta(\B^n;T^*\B^n)}\lesssim_{\beta}\|\codiff d\lambda\|_{\Co^{\beta-2}}\lesssim_\beta\|d\lambda\|_{\Co^{\beta-1}},\quad\|\xi\|_{\Co^\alpha(\B^n)}\lesssim_\alpha\|\codiff\lambda\|_{\Co^{\alpha-1}}\lesssim_\alpha\|\lambda\|_{\Co^\alpha}.
\end{equation}

Since convolution commutes with $d$ and $\codiff$, by construction $\lambda-\rho=-\Ga\ast(d\codiff+\codiff d)\lambda+\Ga\ast\codiff\lambda=-d(\Ga\ast\codiff f)=d\xi$. By direct computation along with Lemma \ref{Lem::Hold::Product} $F_*df=d(f\circ F^\Inv)=dF_*f$ for all $f\in\Co^{\alpha+1}_\loc(\R^n)\subset C^1_\loc(\R^n)$, therefore
\begin{equation}\label{Eqn::Hold::QIFT::RestrictionOn1Ball}
    d(F_*\lambda)=d(F_*\rho)+d(F_*d\xi)=d(F_*\rho)+d^2(\xi\circ F^\Inv)=d(F_*\rho)=d\Big(\sum_{i=1}^n(\rho_i\circ\Phi)d\phi^i\Big),\quad\text{on }\B^n.
\end{equation}
Thus, to prove \eqref{Eqn::Hold::QIFT::dFormEst}, by \eqref{Eqn::Hold::QIFT::RhoisBoundedbyDTheta} and the fact that $\supp dF_*\lambda\subseteq\frac34\B^n$ it suffices to show 
\begin{equation}\label{Eqn::Hold::QIFT::Proof2}
    \|d(F_*\rho)\|_{\Co^{\beta-1}(\frac34\B^n;\wedge^2T^*\R^n)}\lesssim\|\rho\|_{\Co^\beta(\B^n;T^*\B^n)}.
\end{equation}

Fix $\chi\in C_c^\infty(\B^n)$ such that $\chi|_{\frac34\B^n}\equiv1$. For each $1\le i\le n$, set $$\tilde\rho_i:=\chi(\rho_i\circ\Phi),\quad\tilde\phi^i:=\chi\phi^i.$$
So $\tilde\rho_i\in\Co^\beta_c(\B^n)$ and $d\tilde\phi^i\in\Co^\alpha_c(\B^n;T^*\B^n)$ are globally defined 1-forms for each $i$, such that 
\begin{equation}\label{Eqn::Hold::QIFT::RestrictionOn3/4Ball}
    \sum_{i=1}^n(\tilde\rho_id\tilde\phi^i)\big|_{\frac34\B^n}=\sum_{i=1}^n\left((\rho_i\circ\Phi)d\phi^i\right)\big|_{\frac34\B^n}=(F_*\rho)\big|_{\frac34\B^n}.
\end{equation}

By \eqref{Eqn::Hold::QIFT::RestrictionOn3/4Ball} and \eqref{Eqn::Hold::QIFT::RestrictionOn1Ball} we have $$\|d(F_*\lambda)\|_{\Co^{\beta-1}(\frac34\B^n;\wedge^2T^*\R^n)}=\Big\|\sum_{i=1}^nd(\tilde\rho_id\tilde\phi^i)\Big\|_{\Co^{\beta-1}(\frac34\B^n;\wedge^2T^*\R^n)}\le\sum_{i=1}^n\|d(\tilde\rho_id\tilde\phi^i)\|_{\Co^{\beta-1}(\R^n;\wedge^2T^*\R^n)}.$$
By Lemma \ref{Lem::Hold::Product} we have $\|\tilde\rho_i\|_{\Co^{\beta}(\R^n)}\lesssim_\beta\|\chi\|_{\Co^\beta(\B^n)}\|\rho_i\circ\Phi\|_{\Co^\beta(\B^n)}$ and $\|\tilde\phi^i\|_{\Co^{\alpha+1}(\R^n)}\lesssim_\alpha\|\chi\|_{\Co^{\alpha+1}(\B^n)}\|\Phi\|_{\Co^{\alpha+1}(\B^n)}$, by the result \ref{Item::Hold::QIFT::PhiEst} and Proposition \ref{Prop::Hold::QComp} \ref{Item::Hold::QComp::>1} we have $\|\Phi\|_{\Co^{\alpha+1}(\B^n;\R^n)}\lesssim1$ and $\|\rho_i\circ\Phi\|_{\Co^\beta(\B^n)}\lesssim\|\rho_i\|_{\Co^\beta(\B^n)}$. Combining them we get
\begin{equation}\label{Eqn::Hold::QIFT::BddTrhoTphi}
    \|\tilde\rho_i\|_{\Co^\beta(\R^n)}\lesssim_{\alpha,\beta}\|\rho_i\|_{\Co^\beta(\B^n)},\quad\|\tilde\phi^i\|_{\Co^{\alpha+1}(\R^n)}\lesssim_\alpha1,\quad1\le i\le n.
\end{equation}

Thus, to obtain \eqref{Eqn::Hold::QIFT::Proof2} and complete the proof, it 
suffices to show
\begin{equation}\label{Eqn::Hold::QIFT::Proof2::Tmp1}
\|d(\tilde\rho_id\tilde\phi^i)\|_{\Co^{\beta-1}(\R^n;\wedge^2T^*\R^n)}\lesssim_{\alpha,\beta}\|\tilde\rho_i\|_{\Co^\beta(\B^n)},\quad 1\le i\le n. 
\end{equation}

To prove \eqref{Eqn::Hold::QIFT::Proof2::Tmp1} we define 1-form $\tau_i$ on $\R^n$ by $$\tau_i:=\Pf(\tilde\rho_i,d\tilde\phi^i)-\Pf(\tilde\phi^i,d\tilde\rho_i)+\Rf(\tilde\rho_i,d\tilde\phi^i),\quad1\le i\le n.$$

By Lemma \ref{Lem::Hold::LemParaProd}, we have $\tilde\rho_id\tilde\phi^i=\Pf(\tilde\rho_i,d\tilde\phi^i)+\Pf(d\tilde\phi^i,\tilde\rho_i)+\Rf(\tilde\rho_i,d\tilde\phi^i)$. Since as vector-valued functions $d\Pf(\tilde\phi^i,\tilde\rho_i)=\Pf(d\tilde\phi^i,\tilde\rho_i)+\Pf(\tilde\phi^i,d\tilde\rho_i)$ holds for $1\le i\le n$ on $\R^n$, we have
\begin{equation}\label{Eqn::Hold::QIFT::Proof2::Tmp2}
d\tau_i=d\Pf(\tilde\rho_i,d\tilde\phi^i)-d\Pf(\tilde\phi^i,d\tilde\rho_i)+d\Rf(\tilde\rho_i,d\tilde\phi^i)=d\Pf(\tilde\rho_i,d\tilde\phi^i)+d\Pf(d\tilde\phi^i,\tilde\rho_i)+d\Rf(\tilde\rho_i,d\tilde\phi^i)=d(\tilde\rho_id\tilde\phi^i).
\end{equation}

We claim that
\begin{equation}\label{Eqn::Hold::QIFT::Proof2::Tmp3}
    \| \tau_i\|_{\Co^{\beta}(\R^n;T^*\R^n)}\lesssim
    \|\tilde\rho_i\|_{\Co^\beta(\R^n;T^*\R^n)},\quad 1\le i\le n.
\end{equation}
By Lemma \ref{Lem::Hold::LemParaProd}, along with the fact that $\|\tilde\phi^i\|_{\Co^{\alpha+1}}\lesssim1$, we get $\|\Pf(\tilde\rho_i,d\tilde\phi^i)\|_{\Co^\beta}\lesssim\|\tilde\rho_i\|_{\Co^\beta}$ and $\|\Rf(\tilde\rho_i,d\tilde\phi^i)\|_{\Co^\beta}\lesssim\|\Rf(\tilde\rho_i,d\tilde\phi^i)\|_{\Co^{\alpha+\beta}}\lesssim\|\tilde\rho_i\|_{\Co^\beta}$. 

\medskip

To complete the proof of \eqref{Eqn::Hold::QIFT::Proof2::Tmp3}, we
need to show 
\begin{equation}\label{Eqn::Hold::QIFT::Proof2::Tmp4}
    \|\Pf(\tilde\phi^i,d\tilde\rho_i)\|_{\Co^\beta(\R^n;T^*\R^n)}\lesssim_{\alpha,\beta}\|\tilde\rho_i\|_{\Co^\beta(\R^n;T^*\R^n)},\quad1\le i\le n.
\end{equation}
We separate the proof of \eqref{Eqn::Hold::QIFT::Proof2::Tmp4} into two cases:  $\beta>1$ and $0<\beta\leq 1$. We apply Lemma \ref{Lem::Hold::ParaBdd}.

For the case $\beta>1$, \eqref{Eqn::Hold::QIFT::Proof2::Tmp4} follows from \eqref{Eqn::Hold::ParaBdd::PfBdd1} with \eqref{Eqn::Hold::QIFT::BddTrhoTphi} that $\tilde\phi^i\in\Co^{\alpha+1}\subseteq\Co^\beta$ and $d\tilde\rho_i\in\Co^{\beta-1}\subsetneq L^\infty$. 
For the case $\beta\le1$, our assumption $0<\alpha\le\beta$ implies that $\alpha\in(0,1]$. So $\beta-\alpha-1<\beta-1\leq 0$ and we have $d\tilde\rho_i\in\Co^{\beta-1}\subsetneq \Co^{\beta-\alpha-1}$. By \eqref{Eqn::Hold::ParaBdd::PfBdd2} along with $\tilde\phi^i\in\Co^{\alpha+1}$ from \eqref{Eqn::Hold::QIFT::BddTrhoTphi}, we get $\Pf(\tilde\phi^i,d\tilde\rho_i)\in\Co^{(\alpha+1)+(\beta-\alpha-1)}=\Co^\beta$, which gives \eqref{Eqn::Hold::QIFT::Proof2::Tmp4} and establishes \eqref{Eqn::Hold::QIFT::Proof2::Tmp3}.

Using \eqref{Eqn::Hold::QIFT::Proof2::Tmp2} and \eqref{Eqn::Hold::QIFT::Proof2::Tmp3},
we see 
$\|d(\tilde\rho_id\tilde\phi^i)\|_{\Co^{\beta-1}(\R^n;\wedge^2T^*\R^n)}=\|d\tau_i\|_{\Co^{\beta-1}}\lesssim\|\tau_i\|_{\Co^\beta}\lesssim\|\tilde\rho_i\|_{\Co^\beta}$, establishing \eqref{Eqn::Hold::QIFT::Proof2::Tmp1}
and completing the proof (by possibly increasing $K_2$).
\end{proof}

In Section \ref{Section::EllipticPara::Scaling}, we need the following existence of second order elliptic operators with zero Dirichlet boundary:
\begin{prop}\label{Prop::Hold::StrongLapInv}
Let $m,n\ge1$ and $\alpha\in(\frac12,1]$. Let $U\subseteq\R^n$ be a bounded smooth domain. Then there is a $K_3=K_3(U,m,\alpha)>0$ that satisfies the following:

Let $\displaystyle L_B:=\sum_{j,k=1}^{n}B^1_{jk}\partial_j(B^2_{jk}\partial_k)+\sum_{l=1}^{n}B^3_l\partial_l$ be a second order differential operator on $\C^m$-valued functions, with $B^1_{jk},B^2_{jk}\in\Co^\alpha(U;\C^{m\times m})$ and $B^3_l\in\Co^{\alpha-1}(U;\C^{m\times m})$. Suppose 
\begin{equation}\label{Eqn::Hold::StrongLapInv::Assumption}
    \|B\|_{\alpha}=\|B\|_{\alpha;\phi,E}:=\sum_{j,k=1}^{n}\sum_{\nu=1}^2\|B^\nu_{jk}\|_{\Co^\alpha(U;\C^{m\times m})}+\sum_{l=1}^{n}\|B^3_l\|_{\Co^{\alpha-1}(U;\C^{m\times m})}<K_3^{-1}.
\end{equation}

Then the following map is (both side) invertible,
\begin{equation*}
    \Delta+L_B:\{F\in\Co^{\beta+1}(U;\C^m):F|_{\partial U}\equiv0\}\to\Co^{\beta-1}(U;\C^m),\quad\forall\beta\in(1-\alpha,\alpha].
\end{equation*}

\end{prop}

The difficult part is to show the invertibility when $\beta$ closed to $1-\alpha$, whereas $K_3$ cannot depend on $\beta$.

In this subsection only, we use Littlewood-Paley characterizations to define $\Co^\beta$-norms. More precisely:
\begin{defn}\label{Defn::Hold::HoldLPNorm}
Let $\phi=(\phi_j)_{j=0}^\infty$ be a dyadic resolution (see Definition \ref{Defn::Hold::DyadicResolution}) of $\R^{n}$ and let $E$ be a fixed extension operator for $U$ such that $E:\Co^\beta(U)\to\Co^\beta(\R^{n})$ is bounded linear for all $-1\le\beta\le2$. Let $\Xs$ be a finite dimensional vector space, for $\tilde f\in\Co^\beta(\R^{n};\Xs)$ and $f\in\Co^\beta(U;\Xs)$ we use the following norms
\begin{equation}\label{Eqn::Hold::HoldLPNorm}
    \|\tilde f\|_{\Co^\beta_\phi(\Xs)}:=\sup_{j\ge0}2^{j\beta}\|\phi_j\ast\tilde f\|_{L^\infty(\R^{n};\Xs)},\quad \|f\|_{\Co^\beta_{\phi,E}(\Xs)}:=\|Ef\|_{\Co^\beta_\phi(\Xs)}.
\end{equation}
\end{defn}

By Lemma \ref{Lem::Hold::SingExtLem} such $E$ exists. By Lemma \ref{Lem::Hold::HoldChar} and Definition \ref{Defn::Hold::NegHold}, $\|\cdot\|_{\Co^\beta_\phi(\Xs)}$ is an equivalent norm of $\Co^\beta(\R^{n};\Xs)$ for all $-1\le\beta\le2$.

Clearly for every  $-1\le\beta_0\le\beta_1\le2$, 
\begin{equation*}
    \|\tilde f\|_{\Co^{\beta_0}_\phi}\le\|\tilde f\|_{\Co^{\beta_1}_\phi},\quad\| f\|_{\Co^{\beta_0}_{\phi,E}}\le\| f\|_{\Co^{\beta_1}_{\phi,E}},\quad\forall\tilde f\in\Co^{\beta_1}(\R^{n}),\quad f\in\Co^{\beta_1}(U).
\end{equation*}
Note that both inequalities have no extra constant.

For $\beta\in(0,2]$, we see that for every $\tilde f\in\Co^\beta(\R^{n})$,
\begin{equation}\label{Eqn::Hold::LPNormLInftyBdd}
    \textstyle\|\tilde f\|_{L^\infty}\le\sum_{j=0}^\infty\|\phi_j\ast\tilde f\|_{L^\infty}\le\|\tilde f\|_{\Co^\beta_\phi}\sum_{j=0}^\infty2^{-j\beta}=\|\tilde f\|_{\Co^\beta_\phi}\tfrac1{1-2^{-\beta}}\le\frac3{\beta}\|\tilde f\|_{\Co^\beta_\phi}.
\end{equation}
Thus $\|f\|_{L^\infty}\le\frac3\beta\|f\|_{\Co^\beta_{\phi,E}}$ as well.

Moreover by \cite[Theorem 2.4.7(i)]{Triebel1}, $\{\Co^\beta(\R^{n}),\|\cdot\|_{\Co^\beta_\phi}:-1\le\beta\le2\}$ forms an exact complex interpolation family, see \cite[Definition 1.2.2/2]{TriebelInterpolation}. In particular, for any $\beta_0,\beta_1,\gamma_0,\gamma_1\in[-1,2]$ and for any linear operator $\Tc$ such that $\Tc:\Co^{\beta_j}(\R^{n})\to\Co^{\gamma_j}(\R^{n})$ is bounded for $j=0,1$, then $\Tc:\Co^{(1-\theta)\beta_0+\theta\beta_1}(\R^{n})\to\Co^{(1-\theta)\gamma_0+\theta\gamma_1}(\R^{n})$ is also bounded for all $0<\theta<1$, with operator norm
\begin{equation}\label{Eqn::Hold::Interpo}
    \|\Tc\|_{\Co^{(1-\theta)\beta_0+\theta\beta_1}_\phi\to\Co^{(1-\theta)\gamma_0+\theta\gamma_1}_\phi}\le\|\Tc\|_{\Co^{\beta_0}_\phi\to\Co^{\gamma_0}_\phi}^{1-\theta}\|\Tc\|_{\Co^{\beta_1}_\phi\to\Co^{\gamma_1}_\phi}^{\theta}.
\end{equation}
There is no extra constant in \eqref{Eqn::Hold::Interpo} as well.

Applying Lemma \ref{Lem::Hold::LapInvBdd} via scaling, $\Delta:\{f\in\Co^{\beta+1}(U):f|_{\partial U}\}\to\Co^{\beta-1}(U)$ is invertible for all $0\le\beta\le1$. Recall that we denote its inverse by $\Pc_0:\Co^{\beta-1}(U)\to \Co^{\beta+1}(U)$ in Notation \ref{Note::Hold::DiriSol}.

By \eqref{Eqn::Hold::HoldLPNorm} $\|\Pc_0\|_{\Co^{\beta-1}_{\phi,E}\to\Co^{\beta+1}_{\phi,E}}=\|\tilde f\mapsto E\Pc_0(\tilde f|_{U})\|_{\Co^{\beta-1}_\phi\to\Co^{\beta+1}_\phi}$. Therefore by \eqref{Eqn::Hold::Interpo}, there is a $C_0=C_0(U,m,\phi,E)>0$ such that
\begin{equation}\label{Eqn::Hold::TildePcBdd}
    \|\Pc_0\|_{\Co^{\beta-1}_{\phi,E}\to\Co^{\beta+1}_{\phi,E}}\le C_0,\quad\forall 0\le\beta\le1.
\end{equation}

 We need an improvement of Lemma \ref{Lem::Hold::ParaBdd}.
\begin{lem}[{\cite[Theorems 2.82 and 2.85]{BahouriCheminDanchin}}]\label{Lem::Hold::ImprovedParaBdd}
Let $\Pf$ and $\Rf$ be the paraproduct operators in Definition \ref{Defn::Hold::ParaOp} associated with the same $(\phi_j)_{j=0}^\infty$ in Definition \ref{Defn::Hold::HoldLPNorm}. There is $C_2=C_2(n,\phi)>0$ such that
\begin{align}
    \|\Pf(f,g)\|_{\Co^\beta_\phi}&\le C_2\|f\|_{\Co^\beta_\phi}\|g\|_{L^\infty},
    &&\forall \beta\in[-1,2],\quad f\in\Co^\beta(\R^{n}),\ g\in L^\infty(\R^{n});
    \\
    \|\Pf(f,g)\|_{\Co^{\beta+\gamma}_\phi}&\le \frac{C_2}{-\gamma}\|f\|_{\Co^\beta_\phi}\|g\|_{\Co^\gamma_\phi},
    &&\forall \beta\in[-1,2],\quad \gamma\in[-1,0),\ f\in\Co^\beta(\R^{n}),\ g\in \Co^\gamma(\R^{n});
    \\
    \|\Rf(f,g)\|_{\Co^{\beta+\gamma}_\phi}&\le \frac{C_2\|f\|_{\Co^\beta_\phi}\|g\|_{\Co^\gamma_\phi}}{\beta+\gamma},
    &&\forall \beta,\gamma\in[-1,2],\ \beta+\gamma>0,\ f\in\Co^\beta(\R^{n}),\ g\in \Co^\gamma(\R^{n}).
\end{align}
\end{lem}

\begin{proof}[Proof of Proposition \ref{Prop::Hold::StrongLapInv}]In the proof we use $\Co^\gamma_\phi$-norm from Definition \ref{Defn::Hold::HoldLPNorm} instead of the classical $\Co^\gamma$-norm in Definition \ref{Defn::Intro::DefofHold}. 

Recall in Notation \ref{Note::Hold::DiriSol} $\Pc_0:\Co^{\beta-1}(U;\C^m)\to \Co^{\beta+1}(U;\C^m)$ is the inverse operator to $\Delta$ with zero Dirichlet boundary, for all $\beta\in(1-\alpha,\alpha]$. Thus if $\Delta+L_B$ is invertible, then $(\Delta+L_B)^{-1}=(\id+\Pc_0 L_B)^{-1}\Pc_0$, and $\id+\Pc_0 L_B$ is an isomorphism on $\{F\in\Co^{\beta+1}(U;\C^m):F|_{\partial U}\equiv0\}$.

By \eqref{Eqn::Hold::TildePcBdd} $\Pc_0$ is uniformly bounded. Therefore to prove \eqref{Eqn::Hold::StrongLapInv::Bdd} it suffices to show $\id+\Pc_0 L_B$ is an isomorphism on $\Co^{\beta+1}(U;\C^m)$ when \eqref{Eqn::Hold::StrongLapInv::Assumption} is satisfied. 

We also assume that $K_3$ can depend on $\phi$ and $E$ in the proof. The extra dependence can be eliminated when we switch the $\|\cdot\|_{\Co^\alpha_{\phi,E}}$ norm back to the $\|\cdot\|_{\Co^\alpha(U)}$ norm.

We claim that, by taking $K_3$ suitably large, we have
\begin{equation}\label{Eqn::Hold::StrongLapInv::Bdd}
    \|(\id+\Pc_0 L_B)^{-1}\|_{\Co^{\beta-1}_{\phi,E}(\C^m)\to \Co^{\beta+1}_{\phi,E}(\C^m)}\le K_3\cdot (\beta+\alpha-1)^{-1},\quad\forall \beta\in(1-\alpha,\alpha].
\end{equation}

For $F=(F^1,\dots,F^m)$, we write 
$$(L_BF)^i=\sum_{j,k=1}^{n}\sum_{u,v=1}^mB^{1,iu}_{jk}\partial_j(B^{2,uv}_{jk}\partial_kF^v)+\sum_{l=1}^{n}\sum_{u=1}^mB^{3,iu}_l\partial_lF^u,\quad i=1,\dots,m.$$
We take decomposition $L_B=L_{B,1}+L_{B,2}$ where for $i=1,\dots,m$,
\begin{align*}
    (L_{B,1}F)^i:=&\sum_{j,k=1}^{n}\sum_{u,v=1}^m\left(\Pf \big(EB^{1,iu}_{jk},E\partial_j(B^{2,uv}_{jk}\partial_kF^v)\big)+\Pf \big(E\partial_j(B^{2,uv}_{jk}\partial_kF^v),EB^{1,iu}_{jk}\big)\right)\Big|_{U}
    \\
    &+\sum_{l=1}^{n}\sum_{u=1}^m\left(\Pf(EB^{3,iu}_l,E\partial_lF^u)+\Pf(E\partial_lF^u,EB^{3,iu}_l)\right)\Big|_{U};
    \\
    (L_{B,2}F)^i:=&\sum_{j,k=1}^{n}\sum_{u,v=1}^m\Rf\big(EB^{1,iu}_{jk},E\partial_j(B^{2,uv}_{jk}\partial_kF^v)\big)\Big|_{U}+\sum_{l=1}^{n}\sum_{u=1}^m\Rf(EB^{3,iu}_l,E\partial_lF^u)\Big|_{U}.
\end{align*}

By Lemma \ref{Lem::Hold::ImprovedParaBdd} $\|B^{2,uv}_{jk}\partial_kF^v\|_{\Co^\beta_{\phi,E}}$ is uniformly bounded in $\alpha\in(\frac12,1]$ and $\beta\in(1-\alpha,\alpha]$:
\begin{align*}
    &\|B^{2,uv}_{jk}\partial_kF^v\|_{\Co^\beta_{\phi,E}}=\|\Pf(EB^{2,uv}_{jk},E\partial_kF^v)+\Pf(E\partial_kF^v,EB^{2,uv}_{jk})+\Rf(EB^{2,uv}_{jk},E\partial_kF^v)\|_{\Co^\beta_{\phi}}
    \\
    \le& C_2\min\big(\|EB^{2,uv}_{jk}\|_{\Co^\beta_\phi}\|E\partial_kF^v\|_{L^\infty},\tfrac1{\alpha-\beta}\|EB^{2,uv}_{jk}\|_{\Co^\alpha_\phi}\|E\partial_kF^v\|_{\Co^{\beta-\alpha}_\phi}\big)
    \\&+C_2\|EB^{2,uv}_{jk}\|_{L^\infty}\|E\partial_kF^v\|_{\Co^\beta_\phi}+\tfrac{C_2}\alpha\|EB^{2,uv}_{jk}\|_{\Co^\alpha_\phi}\|E\partial_kF^v\|_{\Co^{\beta-\alpha}_\phi}&\hspace{-0.3in}\text{(by Lemma \ref{Lem::Hold::ImprovedParaBdd})}
    \\
    \le&C_2(\min(\tfrac3\beta,\tfrac1{\alpha-\beta})+1+\tfrac1\alpha)\|EB^{2,uv}_{jk}\|_{\Co^\alpha_\phi}\|E\partial_kF^v\|_{\Co^\beta_\phi}\le11C_2\|B^{2,uv}_{jk}\|_{\Co^\alpha_{\phi,E}}\|\partial_kF^v\|_{\Co^\beta_{\phi,E}}.                
    &\text{(by \eqref{Eqn::Hold::LPNormLInftyBdd})}
\end{align*}

Similarly $L_{B,1}:\Co^{\beta+1}(U;\C^m)\to \Co^{\beta-1}(U;\C^m)$ is uniformly bounded in $\alpha\in(\frac12,1],\beta\in(1-\alpha,\alpha]$: 
\begin{align*}
    \|L_{B,1}F\|_{\Co^{\beta-1}_\phi}\lesssim&_E\tfrac1{\alpha+1-\beta}\|EB^1\|_{\Co^\alpha_\phi}\|E\nabla(B^2\otimes\nabla F)\|_{\Co^{\beta-1-\alpha}_\phi}+\|EB^1\|_{L^\infty}\|E\nabla(B^2\otimes\nabla F)\|_{\Co^{\beta-1}_\phi}
    \\
    &+\min\big(\|EB^3\|_{\Co^{\beta-1}_\phi}\|E\nabla F\|_{L^\infty},\tfrac1{\alpha-\beta}\|EB^3\|_{\Co^{\alpha-1}_\phi}\|E\nabla F\|_{\Co^{\beta-\alpha}}\big)+\|EB^3\|_{\Co^{-1}_\phi}\|E\nabla F\|_{\Co^\beta_\phi}
    \\
    \le&\big(\tfrac1{\alpha+1-\beta}+\tfrac1\alpha\big)\|B_1\|_{\Co^\alpha_{\phi,E}}\|\nabla(B^2\otimes\nabla F)\|_{\Co^{\beta-1}_{\phi,E}}+\big(\min(\tfrac1\beta,\tfrac1{\alpha-\beta})+1\big)\|B^3\|_{\Co^{\alpha-1}_{\phi,E}}\|\nabla F\|_{\Co^\beta_{\phi,E}}
    \\
    \lesssim&_{U,m}\|B\|_\alpha(1+\|B\|_\alpha)\|F\|_{\Co^{\beta+1}_{\phi,E}}.
\end{align*}
Here we use $\frac1{\alpha+1-\beta}<\frac1{3/2-1}=2$ and $\min(\frac1\beta,\frac1{\alpha-\beta})\le\frac2\alpha<4$.

$L_{B,2}:\Co^{\beta+1}(U;\C^m)\to \Co^{\alpha-1}(U;\C^m)$ is not uniformly bounded but gains a few derivatives: we have $\alpha+\beta>1$,
\begin{align*}
    \|L_{B,2}F\|_{\Co^{\alpha-1}_{\phi,E}}&\le\|L_{B,2}F\|_{\Co^{\alpha+\beta-1}_{\phi,E}}\lesssim_{U,m,E}\tfrac1{\alpha+\beta-1}\big(\|B^1\|_{\Co^\alpha_{\phi,E}}\|\nabla(B^2\otimes\nabla F)\|_{\Co^{\beta-1}_{\phi,E}}+\|B^3\|_{\Co^{\alpha-1}_{\phi,E}}\|\nabla F\|_{\Co^\beta_{\phi,E}}\big)
    \\
    &\lesssim_{U,m}(\alpha+\beta-1)^{-1}\|B\|_\alpha(1+\|B\|_\alpha)\|F\|_{\Co^{\beta+1}_{\phi,E}}.
\end{align*}

Combining \eqref{Eqn::Hold::TildePcBdd}, there is a $C_3=C_3(U,m,\phi,E)>1$ such that for every $\alpha\in(\frac12,1]$, $\beta\in(1-\alpha,\alpha]$ and $F\in\Co^{\beta+1}(U;\C^m)$,
\begin{equation}\label{Eqn::Hold::ImprovedParaBdd::PfTmp}
    \|\Pc_0 L_{B,1}F\|_{\Co^{\beta+1}_{\phi,E}}\le C_3\|B\|_\alpha(1+\|B\|_\alpha)\|F\|_{\Co^{\beta+1}_{\phi,E}},\quad\|\Pc_0 L_{B,2}F\|_{\Co^{\alpha+1}_{\phi,E}}\le\frac{ C_3\|B\|_\alpha(1+\|B\|_\alpha)}{\alpha+\beta-1}\|F\|_{\Co^{\beta+1}_{\phi,E}}.
\end{equation}

Take $K_3:=\frac{8C_3}{2\alpha-1}$, we see that $K_3=K_3(U,m,\phi,E,\alpha)\in(0,1)$. When \eqref{Eqn::Hold::StrongLapInv::Assumption} is satisfied, \eqref{Eqn::Hold::ImprovedParaBdd::PfTmp} implies that $\|\Pc_0 L_{B,1}\|_{\Co^{\beta+1}_{\phi,E}\to \Co^{\beta+1}_{\phi,E}}\le\frac14$ and $\|\Pc_0 L_{B,2}\|_{\Co^{\beta+1}_{\phi,E}\to \Co^{\beta+1}_{\phi,E}}\le\frac14\cdot\frac{2\alpha-1}{\alpha+\beta-1}$. In particular $\|\Pc_0 L_B\|_{\Co^{\alpha+1}_{\phi,E}\to \Co^{\alpha+1}_{\phi,E}}\le\frac14+\frac14<\frac12$ is a contraction mapping.

Therefore for every $\alpha\in(\frac12,1]$ and $\beta\in(1-\alpha,\alpha]$, using power series $(\id+T)^{-1}=\sum_{k=0}^\infty(-1)^kT^k$, provided that the right hand side converges, taking $T\in\{\Pc_0 L_{B,1},(\id+\Pc_0 L_{B,1})^{-1}\Pc_0 L_{B,2}\}$ we get
{\small\begin{align*}
    \|(\id+\Pc_0 L_{B,1})^{-1}\|_{\Co^{\beta+1}_{\phi,E}\to \Co^{\beta+1}_{\phi,E}}&\textstyle\le\sum_{k=0}^\infty\|\Pc_0 L_B\|_{\Co^{\beta+1}_{\phi,E}\to \Co^{\beta+1}_{\phi,E}}^k\le\sum_{k=0}^\infty4^{-k}\le\frac43.
    \\
    \|(\id+(\id+\Pc_0 L_{B,1})^{-1}\Pc_0 L_{B,2})^{-1}\|_{\Co^{\alpha+1}_{\phi,E}\to \Co^{\alpha+1}_{\phi,E}}&\textstyle\le\sum\limits_{k=0}^\infty\|(\id+\Pc_0 L_{B,1})^{-1}\|_{\Co^{\alpha+1}_{\phi,E}\to \Co^{\alpha+1}_{\phi,E}}^k\|\Pc_0 L_{B,2}\|_{\Co^{\alpha+1}_{\phi,E}\to \Co^{\alpha+1}_{\phi,E}}^k\le\sum\limits_{k=0}^\infty3^{-k}\le\frac32.
\end{align*}
}

Finally for the $\Co^{\beta+1}$-boundedness of $(\id+\Pc_0 L_B)^{-1}$, we use  
\begin{align*}
    (\id+\Pc_0 L_B)^{-1}=&(\id+\Pc_0 L_{B,1}+\Pc_0 L_{B,2})^{-1}=(\id+(\id+\Pc_0 L_{B,1})^{-1}\Pc_0 L_{B,2})^{-1}(\id+\Pc_0 L_{B,1})^{-1}
    \\
    =&(\id+\Pc_0 L_{B,1})^{-1}-((\id+\Pc_0 L_{B,1})^{-1}\Pc_0 L_{B,2})(\id+(\id+\Pc_0 L_{B,1})^{-1}\Pc_0 L_{B,2})^{-1}(\id+\Pc_0 L_{B,1})^{-1}.
\end{align*}
Therefore
\begin{align*}
    \|(\id+\Pc_0 L_B)^{-1}\|_{\Co^{\beta+1}_{\phi,E}\to \Co^{\beta+1}_{\phi,E}}\le&\|(\id+\Pc_0 L_{B,1})^{-1}\|_{\Co^{\beta+1}_{\phi,E}\to \Co^{\beta+1}_{\phi,E}}\big(1+\|\Pc_0 L_{B,2}\|_{\Co^{\beta+1}_{\phi,E}\to \Co^{\alpha+1}_{\phi,E}}
\\
&\times\|(\id+(\id+\Pc_0 L_{B,1})^{-1}\Pc_0 L_{B,2})^{-1}\|_{\Co^{\alpha+1}_{\phi,E}\to \Co^{\alpha+1}_{\phi,E}}\|(\id+\Pc_0 L_{B,1})^{-1}\|_{\Co^{\alpha+1}_{\phi,E}\to \Co^{\alpha+1}_{\phi,E}}\big)
    \\
\le&\tfrac43\times\big(1+\tfrac{2\alpha-1}{4(\alpha+\beta-1)}+\tfrac32\cdot\tfrac43\big)\le\tfrac{13}3(\alpha+\beta-1)^{-1}.
\end{align*}
Since $K_3=\frac{8C_3}{2\alpha-1}>\frac3{13}$, we get \eqref{Eqn::Hold::TildePcBdd} and finish the proof.
\end{proof}

\subsection{An approximation theorem on vanishing products}\label{Section::AppThm}

Recall the Littlewood-Paley operators $\De_\sigma=\phi_\sigma\ast(-)$ and $\Su_\sigma=\sum_{j=0}^\sigma\phi_\sigma\ast(-)=2^{n\sigma}\phi_0(2^\sigma\cdot)\ast(-)$ in Definition \ref{Defn::Hold::DyadicResolution}, where $\phi=(\phi_j)_{j=0}^\infty $ is a dyadic resolution. 

The log-Lipschitz Frobenius theorems (Theorems \ref{MainThm::LogFro} and \ref{MainThm::SingFro}) rely heavily on the following:

\begin{thm}\label{Thm::Hold::ApproxThm}
Let $\alpha\in(0,1)$ and $\beta\in(0,\alpha)$. Let $m,n\ge1$. There is a $C=C_{n,m,\phi,\alpha,\beta}>0$ such that, if $\fb\in \Co^\alpha(\R^n;\R^m)$ and $\gb\in \Co^{-\beta}(\R^n;\R^m)$ satisfy $\fb\cdot \gb=0$ in the sense of distributions, then
\begin{equation}\label{Eqn::Hold::ApproxThmEqn}
    \|(\Su_\sigma\fb)\cdot(\Su_\sigma\gb)\|_{C^0(\R^n)}\le C\cdot2^{-(\alpha-\beta)\sigma}\|\fb\|_{\Co^\alpha(\R^n;\R^m)}\|\gb\|_{\Co^{-\beta}(\R^n;\R^m)},\qquad \forall\sigma\ge1.
\end{equation}
\end{thm}
By Lemma \ref{Lem::Hold::Product} \ref{Item::Hold::Product::Hold1}, we know the equality $\fb\cdot\gb=0$ holds in $\Co^{-\beta}(\R^n)$. Note that we already have the convergence $\lim\limits_{\sigma\to\infty}(\Su_\sigma\fb)\cdot(\Su_\sigma\gb)=0$ in $ \Co^{-\beta-\delta}(\R^n)$, for every $\delta>0$. 
Theorem \ref{Thm::Hold::ApproxThm} says that this sequence  converges uniformly with quantitative speed, which is better than the convergence in some spaces of distributions.

In the proof we use the following modified paraproduct formula, which can be deduced using \eqref{Eqn::Hold::LPCharComp} and \eqref{Eqn::Hold::LemParaProd::ParaDecomp2}:
\begin{equation}\label{Eqn::Hold::PfApThm::BonyDecomp}
\De_l(\fb\cdot\gb)=\De_l\sum_{j=l-2}^{l+2}(\De_j\fb)\cdot(\Su_{j-6}\gb)+\De_l\sum_{k=l-2}^{l+2}(\Su_{k-6}\fb)\cdot(\De_k\gb)+\De_l\sum_{j,k\ge l-2;|j-k|\le5}(\De_j\fb)\cdot(\De_k\gb).
\end{equation}
We leave the deduction to readers.

\begin{proof}[Proof of Theorem \ref{Thm::Hold::ApproxThm}]We set $\|\fb\|_{\Co^\alpha}=\|\gb\|_{\Co^{-\beta}}=1$ without loss of generality.

Clearly $\|2^{n\sigma}\phi_0(2^\sigma\cdot)\|_{L^1}=\|\phi_0\|_{L^1}$ for all $\sigma\ge0$, and $\|\phi_\sigma\|_{L^1}=\|2^{n(\sigma-1)}\phi_1(2^{\sigma-1}\cdot)\|_{L^1}=\|\phi_1\|_{L^1}$ for all $\sigma\ge1$. Therefore we have the following uniform bounds,
\begin{equation}\label{Eqn::Hold::PfApThm::PSC0Norm}
    \|\Su_\sigma\|_{C^0\to C^0}\le\|2^{n\sigma}\phi_0(2^\sigma\cdot)\|_{L^1}=\|\phi_0\|_{L^1},\quad\|\De_\sigma\|_{C^0\to C^0}\le\|\phi_\sigma\|_{L^1}\le\|\phi_0\|_{L^1}+\|\phi_1\|_{L^1},\quad\forall\sigma\ge0.
\end{equation}

By \eqref{Eqn::Hold::RmkDyaSupp}, $\Su_\sigma f$ and $\Su_\sigma g$ both have Fourier supports in $\{|\xi|<2^{\sigma+1}\}$. Thus when $l\ge\sigma+3$, $\supp(\Su_\sigma f\cdot\Su_\sigma g)^\wedge\subseteq\{|\xi|<2^{\sigma+2}\}$  is disjoint from $\supp\hat\phi_l\subseteq\{|\xi|>2^{l-1}\}$, which means $\De_l(\Su_\sigma f\cdot\Su_\sigma g)=0$. Therefore we have
\begin{equation}\label{Eqn::Hold::PfApThm::Decomp}
    \Su_\sigma \fb\cdot\Su_\sigma \gb=\sum_{l=0}^\infty\De_l(\Su_\sigma \fb\cdot\Su_\sigma \gb)=\sum_{l=0}^{\sigma+2}\De_l(\Su_\sigma \fb\cdot\Su_\sigma \gb)=\Su_{\sigma-4}(\Su_\sigma \fb\cdot\Su_\sigma \gb)+\sum_{l=\sigma-3}^{\sigma+2}\De_l(\Su_\sigma \fb\cdot\Su_\sigma \gb).
\end{equation}
Our goal is to prove $\|\Su_{\sigma-4}(\Su_\sigma \fb\cdot\Su_\sigma \gb)\|_{C^0(\R^n)}\lesssim2^{-(\alpha-\beta)\sigma}$ and $\sum_{l=\sigma-3}^{\sigma+2}\|\De_l(\Su_\sigma \fb\cdot\Su_\sigma \gb)\|_{C^0(\R^n)}\lesssim2^{-(\alpha-\beta)\sigma}$, where the implicit constants depend only on $m,n,\alpha,\beta$ and $\phi$.

\medskip
By \eqref{Eqn::Hold::PfApThm::BonyDecomp} and the assumption $\fb\cdot\gb=0$, we know,
\begin{equation}\label{Eqn::Hold::PfApThm::Para0}
\De_l\sum_{j=l-2}^{l+2}(\De_j\fb)\cdot(\Su_{j-6}\gb)+\De_l\sum_{k=l-2}^{l+2}(\Su_{k-6}\fb)\cdot(\De_k\gb)+\De_l\sum_{j,k\ge l-2;|j-k|\le5}(\De_j\fb)\cdot(\De_k\gb)=0,\quad l=0,1,2,\dots.
\end{equation}

For $\sigma,l\ge0$, we have decomposition
\begin{align}
    \label{Eqn::Hold::PfApThm::Para-1}
    \De_l(\Su_\sigma\fb\cdot\Su_\sigma\gb)=\De_l(\Su_\sigma\fb\cdot\Su_\sigma\gb)-\De_l(\fb\cdot\gb)=&\De_l\sum_{j=l-2}^{l+2}\big((\De_j\Su_\sigma\fb)\cdot(\Su_{j-6}\Su_\sigma\gb)-\De_j\fb\cdot\Su_{j-6}\gb\big)
    \\
    \label{Eqn::Hold::PfApThm::Para-2}
    &+\De_l\sum_{k=l-2}^{l+2}\big((\Su_{k-6}\Su_\sigma\fb)\cdot(\De_k\Su_\sigma\gb)-\Su_{k-6}\fb\cdot\De_k\gb\big)
    \\
    \label{Eqn::Hold::PfApThm::Para-3}
    &+\De_l\sum_{j,k\ge l-2;|j-k|\le5}\big((\De_j\Su_\sigma\fb)\cdot(\De_k\Su_\sigma\gb)-\De_j\fb\cdot\De_k\gb\big).
\end{align}

When $l\le \sigma-4$ and $|j-l|\le 2$, we have $j\le\sigma-2$. By \eqref{Eqn::Hold::LPCharComp} we have $\De_j\Su_\sigma=\De_j$ and $\Su_{j-6}\Su_\sigma=\Su_{j-6}$, so in this case \eqref{Eqn::Hold::PfApThm::Para-1} is zero. Similarly \eqref{Eqn::Hold::PfApThm::Para-2} vanishes as well. Thus for $\Su_{\sigma-4}(\Su_\sigma \fb\cdot\Su_\sigma \gb)$ in \eqref{Eqn::Hold::PfApThm::Decomp} we have
\begin{align*}
    \Su_{\sigma-4}(\Su_\sigma \fb\cdot\Su_\sigma \gb)
    =&\sum_{l=0}^{\sigma-4}\De_l\sum_{j,k\ge l-2;|j-k|\le5}\big((\De_j\Su_\sigma\fb\cdot \De_k\Su_\sigma\gb)-\De_j\fb\cdot\De_k\gb\big)
    \\
    =&\sum_{l=0}^{\sigma-4}\De_l\sum_{j,k\ge \sigma-2;|j-k|\le5}\big((\De_j\Su_\sigma\fb)\cdot(\De_k\Su_\sigma\gb)-\De_j\fb\cdot\De_k\gb\big)\\
    =&\Su_{\sigma-4}\sum_{j,k\ge \sigma-2;|j-k|\le5}\big((\De_j\Su_\sigma\fb)\cdot (\De_k\Su_\sigma\gb)-\De_j\fb\cdot\De_k\gb\big).
\end{align*}

Using \eqref{Eqn::Hold::PfApThm::PSC0Norm} and Lemma \ref{Lem::Hold::HoldChar} \ref{Item::Hold::HoldChar::LPHoldChar} we have an estimate for \eqref{Eqn::Hold::PfApThm::Para-3}: for every $l\ge0$,
\begin{equation}\label{Eqn::Hold::PfApThm::Term3C0}
    \begin{aligned}
    &\sum_{\substack{j,k\ge l-2\\|j-k|\le5}}\big\|(\De_j\Su_\sigma\fb)\cdot (\De_k\Su_\sigma\gb)-\De_j\fb\cdot\De_k\gb\big\|_{C^0}\le \sum_{\substack{j,k\ge l-2\\|j-k|\le5}}\big(\|\Su_\sigma\|_{C^0\to C^0}^2+1\big)\|\De_j\fb\|_{C^0}\|\De_k\gb\|_{C^0}
    \\
    \lesssim&_\phi\sum_{\substack{j,k\ge l-2\\|j-k|\le5}}\|\De_j\fb\|_{C^0}\|\De_k\gb\|_{C^0}\lesssim_{\phi,\alpha,\beta}\sum_{\substack{j,k\ge l-2\\|j-k|\le5}}2^{-j\alpha}2^{k\beta}\lesssim\sum_{q\ge l}2^{-q\alpha}2^{q\beta}\lesssim2^{-l(\alpha-\beta)}.
    \end{aligned}
\end{equation}

Taking $l=\sigma$ in \eqref{Eqn::Hold::PfApThm::Term3C0} and using \eqref{Eqn::Hold::PfApThm::PSC0Norm} again we get
\begin{align}\label{Eqn::Hold::PfApThm::BddS}
    \|\Su_{\sigma-4}(\Su_\sigma \fb\cdot\Su_\sigma \gb)\|_{C^0}\le\|\Su_{\sigma-4}\|_{C^0\to C^0}\sum_{\substack{j,k\ge l-2\\|j-k|\le5}}\big\|(\De_j\Su_\sigma\fb)\cdot (\De_k\Su_\sigma\gb)-\De_j\fb\cdot\De_k\gb\big\|_{C^0}\lesssim2^{-\sigma(\alpha-\beta)}.
\end{align}
This gives the first estimate $\|\Su_{\sigma-4}(\Su_\sigma \fb\cdot\Su_\sigma \gb)\|_{C^0}\lesssim 2^{-\sigma(\alpha-\beta)}$ for \eqref{Eqn::Hold::PfApThm::Decomp}.

\medskip Next we prove $\|\De_l(\Su_\sigma \fb\cdot\Su_\sigma \gb)\|_{C^0(\R^n)}\lesssim2^{-\sigma(\alpha-\beta)}$ when $\sigma-3\le l\le \sigma+2$. In this case \eqref{Eqn::Hold::PfApThm::Para-1} and \eqref{Eqn::Hold::PfApThm::Para-2} do not vanish.

The term \eqref{Eqn::Hold::PfApThm::Para-1} is still good: by \eqref{Eqn::Hold::PfApThm::PSC0Norm} and Lemma \ref{Lem::Hold::HoldChar} \ref{Item::Hold::HoldChar::LPHoldChar} we have
\begin{equation}\label{Eqn::Hold::PfApThm::Term1C0}
    \begin{aligned}
    \Big\|\sum_{j=l-2}^{l+2}\big((\De_j\Su_\sigma\fb)(\Su_{j-6}\Su_\sigma\gb)-\De_j\fb\cdot\Su_{j-6}\gb\big)\Big\|_{C^0}&\lesssim\sum_{j=l-2}^{l+2}\big(\|\Su_\sigma\|_{C^0\to C^0}^2+1\big)\|\De_j\fb\|_{C^0}\|\Su_{j-6}\gb\|_{C^0}
    \\
    &\lesssim\sum_{j=l-2}^{l+2}2^{-j\alpha}2^{j\beta}\approx2^{-l(\alpha-\beta)}.
\end{aligned}
\end{equation}

Combining \eqref{Eqn::Hold::PfApThm::Term1C0} and \eqref{Eqn::Hold::PfApThm::Term3C0} we have
\begin{equation}\label{Eqn::Hold::PfApThm::Term0C0}
\begin{aligned}
&\Big\|\De_l\sum_{j=l-2}^{l+2}(\De_j\fb)\cdot(\Su_{j-6}\gb)+\De_l\sum_{j,k\ge l-2;|j-k|\le5}(\De_j\fb)\cdot(\De_k\gb)\Big\|_{C^0}
\\
\lesssim&\|\De_l\|_{C^0\to C^0}\Big(\sum_{j=l-2}^{l+2}2^{-j\alpha}2^{j\beta}+\sum_{j,k\ge l-2;|j-k|\le5}2^{-j\alpha}2^{k\beta}\Big)\approx2^{-l(\alpha-\beta)}.
\end{aligned}
\end{equation}

By assumption $\sigma-3\le l\le\sigma+2$, so the right hand side above is bounded by $2^{-\sigma(\alpha-\beta)}$.

It remains to estimate \eqref{Eqn::Hold::PfApThm::Para-2}. Note that for $|k-l|\le 2$ and $l\le\sigma+2$ we have $k\le\sigma+4$. By \eqref{Eqn::Hold::LPCharComp} we have $\Su_{k-6}\Su_\sigma=\Su_{j-6}$. Therefore we can rewrite the summand in \eqref{Eqn::Hold::PfApThm::Para-2} as, 
\begin{equation}
    \textstyle(\Su_{k-6}\Su_\sigma\fb)\cdot(\De_k\Su_\sigma\gb)-\Su_{k-6}\fb\cdot\De_k\gb
    =\Su_{k-6}\fb\cdot(\De_k\Su_\sigma\gb-\De_k\gb)=\Su_{k-6}\fb\cdot\big(\De_k\sum_{\mu=\sigma+1}^{\sigma+5}\De_\mu\gb\big).
\end{equation}
Here for the last  equality above, we use the fact $\De_k(\Su_\sigma-\id)=\De_k\sum_{\mu=\sigma+1}^\infty\De_\mu=\De_k\sum_{\mu=\sigma+1}^{\sigma+5}\De_\mu$.
We define
$$\textstyle\psi_\sigma:=\sum_{\mu=\sigma+1}^{\sigma+5}\phi_\mu,\quad\sigma\ge0.$$Thus $\De_k\sum_{\mu=\sigma+1}^{\sigma+5}\De_\mu\gb=\De_k(\psi_\sigma\ast\gb)=\psi_\sigma\ast\De_k\gb$. 
Therefore for $x\in\R^n$,
\begin{equation}\label{Eqn::Hold::PfApThm::Para2Est}
    \begin{aligned}
    &\Su_{k-6}\fb\cdot(\psi_\sigma\ast\De_k\gb )(x)=\int_{\R^n}(\Su_{k-6}\fb)(x)\cdot(\De_k\gb )(x-y)\psi_\sigma(y)dy
    \\
    =&\big(\psi_\sigma\ast(\Su_{k-6}\fb\cdot\De_k\gb)\big)(x)+\int_{\R^n}\psi_\sigma(y)\big((\Su_{k-6}\fb)(x)-(\Su_{k-6}\fb)(x-y)\big)\cdot(\De_k\gb )(x-y)dy
    \\
    =&\big(\psi_\sigma\ast(\Su_{k-6}\fb\cdot\De_k\gb)\big)(x)+\int_{\R^n}\Big(\int_0^1y\psi_\sigma(y)\cdot(\nabla\Su_{k-6}\fb)(x-ty)dt\Big)\cdot(\De_k\gb )(x-y)dy.
\end{aligned}
\end{equation}

By scaling we have  $\|x\phi_\mu(x)\|_{L^1_x}=2^{1-\mu}\|x\phi_1(x)\|_{L^1_x}$ and $\|\phi_\mu\|_{L^1}=\|\phi_1\|_{L^1}$ for all $\mu\ge1$, thus
\begin{equation}\label{Eqn::Hold::PfApThm::PsiNorm}
    \|x\psi_\sigma(x)\|_{L^1_x}=2^{-\sigma}\|x\psi_0(x)\|_{L^1_x}\approx_\phi2^{-\sigma},\quad\|\psi_\sigma\|_{L^1}\approx_\phi1,\quad\forall \sigma\ge0.
\end{equation}
Combining \eqref{Eqn::Hold::PfApThm::PsiNorm} with \eqref{Eqn::Hold::PfApThm::Para0} and \eqref{Eqn::Hold::PfApThm::Para2Est}, for $\sigma-3\le l\le\sigma+2$ we have
\begin{align*}
    &\Big\|\De_l\sum_{k=l-2}^{l+2}\big((\Su_{k-6}\Su_\sigma\fb)\cdot(\De_k\Su_\sigma\gb)-\Su_{k-6}\fb\cdot\De_k\gb\big)\Big\|_{C^0}
    \\
    \le&\Big\|\De_l\sum_{k=l-2}^{l+2}\psi_\sigma\ast(\Su_{k-6}\fb\cdot\De_k\gb)\Big\|_{C^0}+\|\De_l\|_{C^0\to C^0}\sum_{k=l-2}^{l+2}\int_{\R^n}|y\psi_\sigma(y)|\|\nabla\Su_{k-6}\fb\|_{C^0}\|\De_k\gb\|_{C^0}dy
    \\
    \le&\Big\|\psi_\sigma\ast\De_l\Big(\sum_{j=l-2}^{l+2}\De_j\fb\cdot\Su_{j-6}\gb+\sum_{\substack{|j-k|\le5\\j,k\ge l-2}}\De_j\fb\cdot\De_k\gb\Big)\Big\|_{C^0}+\|\phi_l\|_{L^1}\sum_{k=l-2}^{l+2}\|y\psi_\sigma\|_{L^1_y}\|\Su_{k-6}\nabla \fb\|_{C^0}\|\De_k\gb\|_{C^0}&\text{by }\eqref{Eqn::Hold::PfApThm::Para0}
    \\
    \lesssim&\|\psi_\sigma\|_{L^1}2^{-\sigma(\alpha-\beta)}+\|x\mapsto x\psi_\sigma(x)\|_{L^1}\sum_{k=l-2}^{l+2}\sum_{\mu=0}^{k-6}\|\De_\mu\nabla\fb\|_{C^0}\|\De_k\gb\|_{C^0}&\text{by }\eqref{Eqn::Hold::PfApThm::Term0C0}
    \\
    \lesssim&\|\psi_\sigma\|_{L^1}2^{-\sigma(\alpha-\beta)}+2^{1-\sigma}\|x\mapsto x\psi_1(x)\|_{L^1}\sum_{\mu=0}^\sigma\|\nabla\fb\|_{\Co^{\alpha-1}}2^{\mu(1-\alpha)}\|\gb\|_{\Co^{-\beta}}2^{l\beta}&\text{by }\eqref{Eqn::Hold::HoldChar::LPHoldChar}
    \\\lesssim&2^{-\sigma(\alpha-\beta)}+2^{-\sigma}2^{\sigma(1-\alpha)}2^{\sigma\beta}\lesssim2^{-\sigma(\alpha-\beta)}.&\text{by }\eqref{Eqn::Hold::PfApThm::PsiNorm}
\end{align*}
Here for the last inequality we use the fact $\|\nabla\fb\|_{\Co^{\alpha-1}}\lesssim\|\fb\|_{\Co^\alpha}=1$ from Remark \ref{Rmk::Hold::GradBdd} \ref{Item::Hold::GradBdd}.

Combining the above estimate with \eqref{Eqn::Hold::PfApThm::Term3C0} and \eqref{Eqn::Hold::PfApThm::Term1C0} we get $\|\De_l(\Su_\sigma \fb\cdot\Su_\sigma \gb)\|_{C^0(\R^n)}\lesssim2^{-\sigma(\alpha-\beta)}$ when $\sigma-3\le l\le \sigma+2$, which is the second estimate for \eqref{Eqn::Hold::PfApThm::Decomp}.

We therefore have $\|\Su_\sigma\fb\cdot\Su_\sigma\gb\|_{C^0}\lesssim2^{-\sigma(\alpha-\beta)}$ and finish the proof.
\end{proof}

 In application we take $\fb$ and $\gb$ in Theorem \ref{Thm::Hold::ApproxThm} to be the coefficients of vector fields and their derivatives. The theorem give a systematical way of approximating vector fields.

\begin{prop}\label{Prop::Hold::InvVFApt}
Let $m\ge1$ and $0<\alpha,\beta<1$ be satisfy $\alpha>\max(\frac12,\beta)$. Let $\Omega\subseteq\R^n$ be an open set, Let $\Omega'\Subset\Omega$ be a precompact open subset and let $\chi\in C_c^\infty(\Omega)$ satisfies $\chi|_{\Omega'}\equiv1$. 

Let $X_1,\dots,X_m\in \Co^\alpha_\loc(\Omega;\R^n)$ be vector fields on $\Omega$. Assuming there are distributions $(c_{ij}^k)_{i,j,k=1}^m\subset \Co^{-\beta}_\loc(\Omega)$ such that
\begin{equation}\label{Eqn::Hold::InvVFApt::InvEqn1}
    [X_i,X_j]=\sum_{k=1}^mc_{ij}^kX_k,\quad 1\le i,j\le m.
\end{equation}

For $\sigma=0,1,2,\dots$, we define 
\begin{equation}\label{Eqn::Hold::InvVFApt::AppVFDef}
    X^\sigma_i:=\Su_\sigma(\chi X_i)=(\Su_\sigma(\chi a_i^1),\dots,\Su_\sigma(\chi a_i^n)),\quad c_{ij}^{k\sigma}:=\Su_\sigma(\chi c_{ij}^k),\quad 1\le i,j,k\le m.
\end{equation}

Then there is a constant constant $C$ that does not depend on $\sigma$ (but can depend on $X_i,c_{ij}^k,\chi,\phi$), such that
\begin{equation}\label{Eqn::Hold::InvVFApt::Main}
    \Big\|[X^\sigma_i,X^\sigma_j]-\sum_{k=1}^mc_{ij}^{k\sigma}X_k^\sigma\Big\|_{C^0(\Omega';\R^n)}\le C2^{-\sigma\min(\alpha-\beta,2\alpha-1)},\quad\forall \sigma\ge0.
\end{equation}
\end{prop}

Recall from Lemmas \ref{Lem::Hold::MultLoc} \ref{Item::Hold::MultLoc::WellDef} and Corollary \ref{Cor::Hold::[X,Y]WellDef}, $c_{ij}^kX_k\in \Co^{-\beta}_\loc(\Omega;\R^n)$ and $[X_i,X_j]\in \Co^{\alpha-1}_\loc(\Omega;\R^n)$ are defined as distributions.

\begin{proof}We can assume $\beta\ge1-\alpha$, hence the right hand side of \eqref{Eqn::Hold::InvVFApt::Main} becomes $C2^{-\sigma(\alpha-\beta)}$. By computation,
\begin{equation}\label{Eqn::Hold::InvVFApt::InvEqn2}
    [\chi X_i,\chi X_j]-(X_i\chi)\cdot \chi X_j+(X_j\chi)\cdot\chi X_i-\sum_{k=1}^m\chi c_{ij}^k\cdot\chi X_k=0\quad\text{on }\R^n,\quad 1\le i,j\le m.
\end{equation}
By writing $X_i=:\sum_{l=1}^na^l\Coorvec{x^l}$ where $a^l\in \Co^\alpha_\loc(\Omega)$, \eqref{Eqn::Hold::InvVFApt::InvEqn2} becomes
\begin{equation}\label{Eqn::Hold::InvVFApt::InvEqn3}
    \sum_{q=1}^n\Big(\chi a_i^q\frac{\partial(\chi a_j^l)}{\partial x^q}-\chi a_j^q\frac{\partial(\chi a_i^l)}{\partial x^q}\Big)+\sum_{q=1}^n\Big(\chi a_i^l\cdot a_j^q\frac{\partial\chi}{\partial x^q}-\chi a_j^l\cdot a_i^q\frac{\partial\chi}{\partial x^q}\Big)-\sum_{k=1}^m\chi c_{ij}^k\cdot\chi a_k^l=0,\ 1\le i,j\le m,\ 1\le l\le n.
\end{equation}

We see that the above equality holds in $\Co^{\min(\alpha-1,-\beta)}(\R^n)=\Co^{-\beta}(\R^n)$.

In order to apply Theorem \ref{Thm::Hold::ApproxThm}, we define a function $\fb_{ij}^l\in \Co^\alpha(\R^n;\R^{2n+m+1})$ and a distribution  $\gb_{ij}^l\in \Co^{-\beta}(\R^n;\R^{2n+m+1})$ for $1\le i,j\le m$ and $1\le l\le n$, as
\begin{equation*}\arraycolsep=0.1pt
    \begin{array}{ccccccccccccccc}
     \fb_{ij}^l:=\big(&\chi a_i^1&,\dots,&\chi a_i^n&,&-\chi a_j^1&,\dots,&-\chi a_j^n&,&-\chi a_1^l&,\dots,&-\chi a_m^l&,&\textstyle\sum_{q=1}^n\big(\chi a_i^la_j^q\partial_q\chi-\chi a_j^l a_i^q\partial_q\chi\big)&\big),
     \\
     \gb_{ij}^l:=\big(&\partial_1(\chi a_j^l)&,\dots,&\partial_n(\chi a_j^l)&,&\partial_1(\chi a_i^l)&,\dots,&\partial_n(\chi a_i^l)&,&\chi c_{ij}^1&,\dots,&\chi c_{ij}^m&,&\centering1&\big).
\end{array}
\end{equation*}

The condition \eqref{Eqn::Hold::InvVFApt::InvEqn3} is now the same as saying $\fb_{ij}^l\cdot\gb_{ij}^l=0\in \Co^{-\beta}(\R^n)$ for all $1\le i,j\le m$, $1\le l\le n$. 

By Theorem \ref{Thm::Hold::ApproxThm} we have $\|\Su_\sigma\fb_{ij}^l\cdot\Su_\sigma\gb_{ij}^l\|_{C^0}\lesssim 2^{-\sigma(\alpha-\beta)}$ for all $i,j,l$ and $\sigma$, in other words,
\begin{equation*}
    \bigg\|[X_i^\sigma,X_j^\sigma]+\Su_\sigma(X_j\chi)\cdot X_i^\sigma-\Su_\sigma(X_i\chi)\cdot X_j^\sigma-\sum_{k=1}^mc_{ij}^{k\sigma}X_k^\sigma\bigg\|_{C^0(\R^n)}\lesssim2^{-\sigma(\alpha-\beta)},\quad 1\le i,j\le m,\quad\sigma\ge0.
\end{equation*}

To prove \eqref{Eqn::Hold::InvVFApt::Main} it remains to show that $\|\Su_\sigma(X_j\chi)\cdot X_i^\sigma-\Su_\sigma(X_i\chi)\cdot X_j^\sigma\|_{C^0(\Omega')}\lesssim 2^{-\sigma(\alpha-\beta)}$ for all $1\le i,j\le m$. Indeed, by assumption $\chi|_{\Omega'}\equiv1$, so $X_i\chi|_{\Omega'}\equiv0$. Hence for $1\le i\le m$,
\begin{equation}\label{Eqn::Hold::InvVFApt::Tmp}
    \begin{aligned}
    \|\Su_\sigma(X_i\chi)\|_{C^0(\Omega')}=&\|\Su_\sigma(X_i\chi)-X_i\chi\|_{C^0(\Omega')}\le \|\Su_\sigma(X_i\chi)-X_i\chi\|_{C^0(\R^n)}\le\sum_{k=\sigma+1}^\infty\|\De_\sigma(X_i\chi)\|_{C^0(\R^n)}
    \\\lesssim&\|X_i\chi\|_{\Co^\alpha(\R^n)}2^{-\sigma\alpha}\lesssim_{\alpha,X,\chi}2^{-\sigma(\alpha-\beta)}.
\end{aligned}
\end{equation}

Therefore
\begin{align*}
    \sum_{i,j=1}^m\|\Su_\sigma(X_j\chi)\cdot X_i^\sigma-\Su_\sigma(X_i\chi)\cdot X_j^\sigma\|_{C^0(\Omega';\R^n)}&\le\sum_{i,j=1}^m\|\Su_\sigma(X_j\chi)\|_{C^0(\Omega')}\|X_i^\sigma\|_{C^0(\Omega;\R^n)}\\
    &\lesssim 2^{-\sigma(\alpha-\beta)}\sup_{1\le i\le m}\|\chi X_i\|_{C^0(\R^n)}\lesssim2^{-\sigma(\alpha-\beta)}.
\end{align*}
This completes the proof.
\end{proof}

As a direct application to Proposition \ref{Prop::Hold::InvVFApt}, we see that for log-Lipschitz subbundles, the distributional involutivity (see Definition \ref{Defn::Intro::DisInv}) implies the strongly asymptotically involutivity (see Definition \ref{Defn::RealFro::Intro::AsyInv}).

\begin{cor}\label{Cor::Hold::CorAsyInv}
Following the assumptions of Proposition \ref{Prop::Hold::InvVFApt}, and assuming $X_1,\dots,X_m\in\Co^\LogL_\loc(\Omega;\R^n)$, then $\lim_{\nu\to\infty}\|X_i^\nu-X_i\|_{C^0(\Omega';\R^n)}=0$ for all $1\le i\le m$, and there is a $t_0>0$ such that
\begin{equation}\label{Eqn::Hold::CorAsyInv::Eqn1}
    \lim\limits_{\nu\to\infty}\max_{1\le i,j\le m}\Big\|[X^\nu_i,X^\nu_j]-\sum_{k=1}^mc_{ij}^{k\nu}X_k^\nu\Big\|_{C^0(\Omega';\R^n)}\exp\Big(t_0\cdot\max_{1\le l\le m}\|\nabla X_l^\nu\|_{C^0(\Omega';\R^{n\times n})}\Big)=0.
\end{equation}

In particular, for the case of two commutative vector fields ($m=2$ and $c_{ij}^k\equiv0$), namely $[X_1,X_2]=0$ in $\Omega$, we have
\begin{equation}\label{Eqn::Hold::CorAsyInv::Eqn2}
        \exists t_0>0,\ \lim\limits_{\nu\to\infty}\|[X_1^\nu,X_2^\nu]\|_{C^0(\Omega';\R^n)}e^{t_0\mleft(\|\nabla X_1^\nu\|_{C^0(\Omega';\R^{n\times n})}+\|\nabla X_2^\nu\|_{C^0(\Omega';\R^{n\times n})}\mright)}=0.
\end{equation}
\end{cor}
\begin{remark}\label{Rmk::Hold::CorAsyInv}For a distributional involutive log-Lipschitz subbundle $\V\le T\Mf$, by Lemma \ref{Lem::ODE::GoodGen} \ref{Item::ODE::GoodGen::InvComm} locally we can find a log-Lipschitz basis $(X_1,\dots,X_r)$ for $\V$ ($r=\rank\V$) such that $[X_i,X_j]=0$ for all $1\le i,j\le r$. So by \eqref{Eqn::Hold::CorAsyInv::Eqn1} with $c_{ij}^k=c_{ij}^{k\nu}\equiv0$, we see that the condition \eqref{Eqn::RealFro::Intro::AsyInv1} is satisfied (see Section \ref{Section::RealFro::HisRmk}). We conclude that a log-Lipschitz involutive subbundle is always strongly asymptotic involutive (see Definition \ref{Defn::RealFro::Intro::AsyInv}).
\end{remark}

\begin{proof}[Proof of Corollary \ref{Cor::Hold::CorAsyInv}]
By assumption $(X_i^\nu -X_i)|_{\Omega'}=(X_i^\nu -\chi X_i)|_{\Omega'}=(\Su_\nu(\chi X_i)-\chi X_i)|_{\Omega'}$ for all $1\le i\le m$. Thus similar to \eqref{Eqn::Hold::InvVFApt::Tmp},
\begin{equation*}
    \|X_i^\nu-X_i\|_{C^0(\Omega';\R^n)}\le\|\Su_\nu(\chi X_i)-\chi X_i\|_{C^0(\Omega';\R^n)}\le\sum_{\nu+1}^\infty\|\De_\nu(\chi X_i)\|_{C^0(\Omega';\R^n)}\lesssim_\alpha\|\chi X_i\|_{C^{0,\alpha}(\R^n;\R^n)}2^{-\nu\alpha}.
\end{equation*}
Let $\nu\to\infty$ we see that $\lim_{\nu\to\infty}\|X_i^\nu-X_i\|_{C^0(\Omega';\R^n)}=0$.

By Proposition \ref{Prop::Hold::InvVFApt} there is a $C_1>0$ and a $\delta>0$ (in fact $\delta=\min(\alpha-\beta,2\alpha-1)$) such that
\begin{equation*}
    \Big\|[X^\nu_i,X^\nu_j]-\sum_{k=1}^mc_{ij}^{k\nu}X_k^\nu\Big\|_{C^0(\Omega';\R^n)}\le C_12^{-\delta \nu},\quad\forall\nu=1,2,3,\dots.
\end{equation*}
By Lemma \ref{Lem::Hold::CharLogL-1} \ref{Item::Hold::CharLogL-1::Grad} we can find a $C_2>0$ such that
\begin{equation*}
    \|\nabla X_i^\nu\|_{C^0(\Omega';\R^n)}\le\|\nabla X_i^\nu\|_{C^0(\R^m;\R^n)}\le C_2\nu,\quad\forall \nu=1,2,3,\dots.
\end{equation*}
Therefore for $t_0>0$ we have
\begin{equation}\label{Eqn::Hold::CorAsyInv::Tmp}
    \Big\|[X^\nu_i,X^\nu_j]-\sum_{k=1}^mc_{ij}^{k\nu}X_k^\nu\Big\|_{C^0(\Omega';\R^n)}e^{t_0\cdot\max_{1\le l\le m}\|\nabla X_l^\nu\|_{C^0(\Omega';\R^{n\times n})}}\le C_12^{-\delta\nu}e^{t_0C_2\nu},\quad \nu=1,2,3,\dots.
\end{equation}
    
Taking $t_0=(2C_2)^{-1}\delta\log2$, the right hand side of \eqref{Eqn::Hold::CorAsyInv::Tmp} is $C_12^{-\frac\delta2\nu}$, which goes to $0$ as $\nu\to\infty$. Therefore we get \eqref{Eqn::Hold::CorAsyInv::Eqn1} for such $t_0$.

Replacing $t_0$ by $\frac {t_0}2$ we get \eqref{Eqn::Hold::CorAsyInv::Eqn2}.
\end{proof}

\section{The Bi-parameter H\"older-Zygmund Spaces}\label{Section::BiHoldSec}
In this section we set up some bi-parameter H\"older-Zygmund spaces and prove some results of products and compositions. 
\subsection{Basic definitions}

We define bi-parameter H\"older-Zygmund structures using vector valued formulation. 

\begin{defn}\label{Defn::Hold::BiHold}
Let $\alpha\in\R$ and $\beta>0$. Let $x=(x^1,\dots,x^n)$ and $s=(s^1,\dots,s^q)$ be the coordinate system of $\R^n$ and $\R^q$ respectively. Let $U\subseteq\R^n$ and $V\subseteq\R^q$ be two open subsets. Let $\Xs$ be a Banach space.
\begin{enumerate}[label=(\roman*),parsep=-0.3ex]
    \item\label{Item::Hold::BiHold::CalphaLinfty} We define $\Co^\alpha_x L^\infty_s(U,V;\Xs)=L^\infty_s\Co^\alpha_x(V,U;\Xs)$ to be the set of all $f\in\D'(U\times V)$ such that $x\mapsto\langle f(x,\cdot),u\rangle_{\D',C_c^\infty(V)}\in \Co^\alpha(U)$ for all $u\in C_c^\infty(V)$ and
    \begin{equation}\label{Eqn::Hold::BiHold::CalphaLinftyNorm}
        \|f\|_{\Co^\alpha L^\infty(U,V;\Xs)}:=\sup_{u\in C_c^\infty(V);\|u\|_{L^1(V)}=1}\|x\mapsto\langle f(x,\cdot),u\rangle_{L^\infty(V),L^1(V)}\|_{\Co^\alpha_x(U;\Xs)}<\infty.
    \end{equation}
    \item\label{Item::Hold::BiHold::CalphaCbeta} We define $\Co^\alpha_x\Co^\beta_s(U,V;\Xs):=\Co^\beta(V;\Co^\alpha(U;\Xs))$ with the same norm given in Definition \ref{Defn::Intro::DefofHold}.
    \item\label{Item::Hold::BiHold::Domain} For $\Ys,\Zs\in\{\Co^\alpha L^\infty,\Co^\alpha\Co^\beta:\alpha\in\R,\beta>0\}$, we define $\Ys\cap\Zs(U,V;\Xs):=\Ys(U,V;\Xs)\cap\Zs(U,V;\Xs)$ and $\Ys+\Zs(U,V;\Xs):=\Ys(U,V;\Xs)+\Zs(U,V;\Xs)$ with norms \begin{align*}
        \|f\|_{\Ys\cap\Zs(U,V;\Xs)}&:=\|f\|_{\Ys(U,V;\Xs)}+\|f\|_{\Zs(U,V;\Xs)},\\ \|f\|_{\Ys+\Zs(U,V;\Xs)}&:=\inf\{\|f_1\|_{\Ys(U,V;\Xs)}+\|f_2\|_{\Zs(U,V;\Xs)}:f=f_1+f_2\}.
    \end{align*}
\end{enumerate}
\end{defn}
\begin{remark}\label{Rmk::Hold::RmkforBiHold}
\begin{enumerate}[parsep=-0.3ex,label=(\roman*)]
    \item\label{Item::Hold::RmkforBiHold::CxCs=CsCx} The $\Co^\alpha\Co^\beta$ spaces are symmetric in the sense that $\Co^\alpha\Co^\beta(U,V)=\Co^\beta\Co^\alpha(V,U)$ when $\alpha>0$. This follows from Definition \ref{Defn::Intro::DefofHold} since we have $\Co^\beta(V;\Co^\alpha(U))=\Co^\alpha(U;\Co^\beta(V))$.
    \item When $\alpha>0$, the pair $\langle f(x,\cdot),u\rangle_{\D',C_c^\infty(V)}$ can be defined pointwise for $x\in U$. When $\alpha\le 0$, the map $[x\mapsto\langle f(x,\cdot),u\rangle_{\D',C_c^\infty(V)}]$ can be realized as a distribution on $U$ in the way that $v\in C_c^\infty(U)\mapsto\langle f,v\otimes u\rangle_{\D',C_c^\infty(U\times V)}$, where $v\otimes u(x,s)=v(x)u(s)$.
    \item \eqref{Eqn::Hold::BiHold::CalphaLinftyNorm} is a rigorous way of saying $\essup_{s\in V}\|f(\cdot,s)\|_{\Co^\alpha(U)}<\infty$, see Lemma \ref{Lem::Hold::BiHoldChar}.
    \item\label{Item::Hold::RmkforBiHold::Extinthesametime}
    For $\alpha,\gamma\le0$, the definition $f\in\Co^\alpha L^\infty\cap \Co^\gamma\Co^\beta( U,V)$ means there are extensions $\tilde f_1\in \Co^\alpha L^\infty(U,\R^q)$ and $\tilde f_2\in \Co^\gamma\Co^\beta(U,\R^q)$ such that $\tilde f_1|_{U\times V}=\tilde f_2|_{U\times V}=f$. When $U$ is a smooth domains, by the existence of common extension operator we can choose $\tilde f_1=\tilde f_2:=E_xf$ where $E_x:\Co^\delta(U)\to\Co^\delta(\R^n)$ is extension operator which is bounded for $\delta\in\{\alpha,\gamma\}$.
\end{enumerate}
\end{remark}


Let $U,V$ be to arbitrary open sets, we use $\mathcal L(L^1(V),\Co^{\alpha}(U))$ as the space of all bounded linear operators $L^1(V)\to\Co^{\alpha}(U)$ endowed with the standard operator norm. This is in fact the same as $\Co^\alpha L^\infty(U,V)$, see Lemma \ref{Lem::Hold::CharCalphaInfty1}.

Let $T:\Co^{\alpha_1}(U_1)\to\Co^{\alpha_2}(U_1)$ be a bounded linear operator, by Lemma \ref{Lem::Hold::OperatorExtension} and Definition \ref{Defn::Hold::BiHold} \ref{Item::Hold::BiHold::CalphaCbeta}, the operator $\tilde Tf(x,s):=T(f(\cdot,s))(x)$ defines a bounded linear operator $\tilde T:\Co^{\alpha_1}\Co^\beta(U_1,V)\to\Co^{\alpha_2}\Co^\beta(U_1,V)$ for $\beta>0$ with the same operator norms. This is also true if we replace $\Co^\beta$ by $L^\infty$.

To summary, for $\Co^\alpha L^\infty$-spaces we have the following:
\begin{lem}\label{Lem::Hold::CharCalphaInfty1}
Let $\alpha_1,\alpha_2\in\R$ and let $U_1,U_2,V$ be arbitrary open sets.
\begin{enumerate}[parsep=-0.3ex,label=(\roman*)]
    \item\label{Item::Hold::CharCalphaInfty1::VectMap} There is a canonical correspondence $\Co^{\alpha_1}L^\infty(U_1,V)\cong\mathcal L(L^1(V),\Co^{\alpha_1}(U_1))$ with the respective norms are the same.
    \item\label{Item::Hold::CharCalphaInfty1::OptExt} Let $T:\Co^{\alpha_1}(U_1)\to\Co^{\alpha_2}(U_2)$ be a bounded linear operator, then there is a unique operator $\tilde T:\Co^{\alpha_1}(U_1,V)\to\Co^{\alpha_2}(U_2,V)$ such that 
\begin{equation}\label{Eqn::Hold::OperatorExtensionInfty}
    [y\in U_2\mapsto\langle(\tilde Tf)(y,\cdot),u\rangle_{\D',C_c^\infty(V)}]=T[x\in U_1\mapsto\langle f(x,\cdot),u\rangle_{\D',C_c^\infty(V)}],\quad\forall f\in\Co^{\alpha_1}(U_1),\quad u\in C_c^\infty(V).
\end{equation}

Moreover $\tilde T$ and $T$ have the same operator norms.
\end{enumerate}
\end{lem}
\begin{proof}
Let $F:L^1(V)\to\Co^{\alpha_1}(U_1)$ be a bounded linear operator, clearly $F:C_c^\infty(V)\to \Co^{\alpha_1}(U_1)\subset\D'(U_1)$ is bounded. Thus by the Schwartz Kernel Theorem $F$ corresponds to a distribution $f\in\D'(U_1\times V)$ such that $\langle F(u),v\rangle_{\D',C_c^\infty(U)}=\langle f,u\otimes v\rangle_{\D',C_c^\infty(U\times V)}$. Since $\|x\mapsto\langle f(x,\cdot),u\rangle_{\D',C_c^\infty(V)}\|_{\Co^{\alpha_1}(U_1)}=\|F(u)\|_{\Co^{\alpha_1}(U_1)}$, taking supremum over $u\in C_c^\infty$ we get $\|f\|_{\Co^{\alpha_1} L^\infty}\le\|F\|_{L^1\to\Co^{{\alpha_1}}}$. In particular such $F$ gives a $f\in\Co^{\alpha_1}L^\infty(U_1,V)$

Conversely, for a $f\in\Co^{\alpha_1}L^\infty(U_1,V)$, we have a map $F:C_c^\infty(V)\to\Co^{\alpha_1}(U_1)$ given by $F(u)=\langle f,u\rangle_{\D',C_c^\infty(V)}$. By \eqref{Eqn::Hold::BiHold::CalphaLinftyNorm}, $\|F(u)\|_{\Co^{\alpha_1}(U_1)}\le\|f\|_{\Co^{\alpha_1}L^\infty(U_1,V)}\|u\|_{L^1(V)}$. Since $C_c^\infty(V)\subset L^1(V)$ is dense, $F$ uniquely extends to a bounded linear map $F:L^1(V)\to\Co^{\alpha_1}(U_1)$ with $\|F\|_{L^1\to\Co^{\alpha_1}}\le\|f\|_{\Co^{\alpha_1}L^\infty}$. This finishes the proof of \ref{Item::Hold::CharCalphaInfty1::VectMap}.

Now the right hand side of \eqref{Eqn::Hold::OperatorExtensionInfty} defines a bounded linear map $\mathcal L(L^1(V),\Co^{\alpha_1}(U_1))\to \mathcal L(L^1(V),\Co^{\alpha_2}(U_2))$, which has the same operator norm to $T$. By \ref{Item::Hold::CharCalphaInfty1::VectMap} this descends to a unique bounded linear map $\Co^{\alpha_1}(U_1,V)\to\Co^{\alpha_2}(U_2,V)$, which is our $\tilde T$.  This finishes the proof of \ref{Item::Hold::CharCalphaInfty1::OptExt}.
\end{proof}

\begin{note}\label{Note::Hold::ConvVar}
For functions $f(x,s)$ and $\rho(x)$, we use $\ast_x$ as the convolution acting on $x$-variable: $$\rho\ast_xf(x,s)=(\rho\ast f(\cdot,s))(x).$$ 
\end{note}
Note that if $f\in\Sc'(\R^n_x\times \R^q_s)$ and $\rho\in\Sc(\R^n)$, then $\rho\ast_xf$ is a well-defined tempered distribution in $\R^n_x\times\R^q_s$, since $\rho\ast_xf=(\rho\otimes\delta_{0^q})\ast f$. Here $\delta_{0^q}$ is the Direc delta measure of the origin in $\R^q$.

\begin{remark}
    In the thesis we use a different formulation of bi-parameter spaces from \cite[Section 2.1]{YaoCpxFro}. These two definitions coincide by the Lemma \ref{Lem::Hold::BiHoldChar} below.
    
    Definition \ref{Defn::Hold::BiHold} \ref{Item::Hold::BiHold::CalphaCbeta} is better than the characterization in Lemma \ref{Lem::Hold::BiHoldChar} \ref{Item::Hold::BiHoldChar::Domain} because of Lemma \ref{Lem::Hold::OperatorExtension}: for an operator $T:\Co^{\alpha_1}(U_1)\to\Co^{\alpha_2}(U_2)$ and its lift $\tilde T:\Co^{\alpha_1}\Co^\beta(U_1,V)\to\Co^{\alpha_2}\Co^\beta(U_1,V)$, we have $\|\tilde T\|=\|T\|$ by Lemma \ref{Lem::Hold::OperatorExtension}. However we do not know the relation between $\|\tilde T\|$ and $\|T\|$ if we use \eqref{Eqn::Hold::BiHoldChar::Domain} to define $\|\tilde T\|$.
\end{remark}
\begin{lem}\label{Lem::Hold::BiHoldChar}
Let $\phi=(\phi_j)_{j=0}^\infty\subset\Sc(\R^n)$ and $\psi=(\psi_k)_{k=0}^\infty\subset\Sc(\R^q)$ be two dyadic resolutions. Let $\alpha\in\R$ and $\beta>0$. Let $U\subseteq\R^n$ and $V\subseteq\R^q$ be two bounded smooth domains.

\begin{enumerate}[label=(\roman*),parsep=-0.3ex]
    \item\label{Item::Hold::BiHoldChar::CalphaLinfty} $f\in \Co^\alpha L^\infty(\R^n,V)$ if and only if $\phi_j\ast_xf\in L^\infty(\R^n\times V)$ for all $j\ge0$ and $\sup_{j\ge0}2^{j\alpha}\|\phi_j\ast_xf\|_{L^\infty}<\infty$. Moreover $f\mapsto\sup_{j\ge0}2^{j\alpha}\|\phi_j\ast_xf\|_{L^\infty}$ is an equivalent norm for $\Co^\alpha L^\infty(\R^n,V)$.
    \item\label{Item::Hold::BiHoldChar::CalphaCbeta} $f\in\Co^\alpha\Co^\beta(\R^n,\R^q)$ has the following equivalent norm:
    \begin{equation*}
        f\mapsto\sup_{j,k\ge0}2^{j\alpha+k\beta}\|(\phi_j\otimes\psi_k)\ast f\|_{L^\infty(\R^n\times\R^q)},\quad\text{where }\phi_j\otimes\psi_k(x,s)=\phi_j(x)\psi_k(s).
    \end{equation*}
    \item\label{Item::Hold::BiHoldChar::Domain} Let $\Xs\in\{L^\infty,\Co^\beta\}$. $\Co^\alpha\Xs(U,V)$ has equivalent norm
    \begin{equation}\label{Eqn::Hold::BiHoldChar::Domain}
        f\mapsto\inf\{\|\tilde f\|_{\Co^\alpha\Xs(\R^n,\R^q)}:\tilde f|_{U\times V}=f\}.
    \end{equation}
    \item\label{Item::Hold::BiHoldChar::CommCAlphaLInfty}For $\alpha>0$, $\Co^\alpha L^\infty(U,V)=\Co^\alpha(U;L^\infty(V))$ and the corresponding norms are equivalent. And the common space equals to the following
    \begin{equation}\label{Eqn::Hold::BiHoldChar::CommCAlphaLInfty::Tmp}
        \{f\in L^\infty(U\times V):f(\cdot,s)\in \Co^\alpha(U)\text{ a.e. }s\in V\text{ with }\essup_{s\in V}\|f(\cdot,s)\|_{\Co^\alpha(U)}<\infty\}.
    \end{equation}
\end{enumerate}

\end{lem}
\begin{remark}
    We postpone the proof of Lemma \ref{Lem::Hold::BiHoldChar} after the proof of Lemma \ref{Lem::Hold::CommuteExt}. Thus from then on we can view $f\in\Co^\alpha L^\infty\cap\Co^\gamma\Co^\beta(U,V)$ as a concrete function/distribution on product space $U\times V$ satisfying the corresponding Littlewood-Paley characterization, rather than a vector-valued map. Thus by taking coordinate components components, we see that Lemma \ref{Lem::Hold::BiHoldChar} is also true if we replace $\Co^\alpha\Xs(U,V)$ ($\Xs\in\{L^\infty,\Co^\beta$) by $\Co^\alpha\Xs(U,V;\R^m)$ or $\Co^\alpha\Xs(U,V;\C^m)$.
\end{remark}

On the domain case we need the following result of extension operators.

\begin{lem}[Common extensions]\label{Lem::Hold::CommuteExt}
Let $R>0$, and let $U\subset\R^n$, $V\subset\R^q$ be two bounded domain with smooth boundaries.
Then there are extension operators $E_x:L^\infty(U)\to L^\infty(\R^n)$ and $E_s:L^\infty(V)\to L^\infty(\R^q)$ such that \begin{enumerate}[nolistsep,label=(\alph*)]
    \item\label{Item::Hold::CommuteExt::BddofExt} $E_x:\Co^\alpha(U)\to\Co^\alpha(\R^n)$ and $E_s:\Co^\alpha(V)\to\Co^\alpha(\R^q)$ are both bounded linear for all $-R<\alpha<R$.
    \item\label{Item::Hold::CommuteExt::Comm} $E_xE_s,E_sE_x:\Co^{(-R)+}L^\infty(U,V)\to \Co^{(-R)+}L^\infty(\R^n,\R^q)$ are well-defined and equal.
\end{enumerate}
Moreover $\tilde E:=E_xE_s$ is an extension operator for $U\times V$ that has boundedness 
\begin{equation}\label{Eqn::Hold::CommuteExt::Bdd}
    \tilde E: \Co^\alpha\Xs(U,V)\to  \Co^\alpha\Xs(\R^n,\R^q),\quad-R<\alpha<R,\quad \Xs\in\{L^\infty,\Co^\beta:0<\beta<R\}.
\end{equation}
\end{lem}
\begin{proof}
We define $E_x$ and $E_s$ as in \eqref{Eqn::Hold::HalfPlaneExt} and \eqref{Eqn::Hold::SmoothExt} by taking $\Omega=U$ and $\Omega=V$ respectively and taking $M=2\lceil R\rceil$. By \cite[Theorems 2.9.2, 2.9.4 (i) and 3.3.4 (i)]{Triebel1} $E_x$ and $E_s$ are both $L^\infty$ and $\Co^\alpha$-bounded for $-R<\alpha<R$. The property \ref{Item::Hold::CommuteExt::Comm} follows immediately from the expression themselves.

One can see that for the adjoint of $E_s$ we have $E_s^*:L^1(\R^q)\to L^1(V)$ with $\|E_s^*\|_{L^1(\R^q)\to L^1(V)}=\|E_s\|_{L^\infty(V)\to L^\infty(\R^q)}\le\sum_{j=1}^M|a_j|$ where $(a_j)$ are in \eqref{Eqn::Hold::HalfPlaneExt}. So for $U'\in\{U,\R^n\}$,
\begin{equation*}
    \|x\mapsto\langle E_sf(x,\cdot),u\rangle_{L^\infty,L^1(\R^q)}\|_{\Co^\alpha(U')}=\|x\mapsto\langle f(x,\cdot),E_s^*u\rangle_{L^\infty,L^1(V)}\|_{\Co^\alpha(U')}\le\sum_{j=1}^M|a_j|\|f\|_{\Co^\alpha L^\infty(U',V)}\|u\|_{L^1(\R^q)}.
\end{equation*}
Therefore we have the boundedness $E_s:\Co^\alpha L^\infty(U',V)\to\Co^\alpha L^\infty(U',\R^n)$.

For $V'\in\{V,\R^q\}$, by Lemma \ref{Lem::Hold::OperatorExtension} $E_x:\Co^\alpha(U)\to\Co^\alpha(\R^n)$ descends to $E_x:\Co^\beta(V';\Co^\alpha(U))\to\Co^\beta(V';\Co^\alpha(\R^n))$, i.e. $E_x:\Co^\alpha\Co^\beta(U,V')\to\Co^\alpha\Co^\beta(\R^n,V')$. By Definition \ref{Defn::Hold::BiHold} \ref{Item::Hold::BiHold::CalphaLinfty} we have $E_x:\Co^\alpha L^\infty(U,V')\to\Co^\alpha L^\infty(\R^n,V')$ immediately. Therefore \eqref{Eqn::Hold::CommuteExt::Bdd} is done by the following:
\begin{equation*}
    E_sE_x:\Co^\alpha\Xs(U,V)\xrightarrow{E_x}\Co^\alpha\Xs(\R^n,V)\xrightarrow{E_s}\Co^\alpha\Xs(\R^n,\R^q).\qedhere
\end{equation*}
\end{proof}

\begin{proof}[Proof of Lemma \ref{Lem::Hold::BiHoldChar}]By the Schwartz Kernel Theorem (see \cite[Theorem 51.7]{TrevesTopologicalSpaces}) $\Co^\alpha L^\infty(\R^n,V)$ can be identified as a subspace of $\D'(\R^n\times V)$. Thus we can view the element $f:L^1(V)\to\Co^\alpha(\R^n)$ as a distribution on $\R^n\times V$.

\smallskip
\noindent\ref{Item::Hold::BiHoldChar::CalphaLinfty}: By Lemma \ref{Lem::Hold::HoldChar} \ref{Item::Hold::HoldChar::LPHoldChar} (if $\alpha>0$) and Definition \ref{Defn::Hold::NegHold} (if $\alpha\le0$) we have $\|f(\cdot,s)\|_{\Co^\alpha(\R^n)}\approx2^{j\alpha}\|\phi_j\ast_xf(\cdot,s)\|_{L^\infty}$ for every $s\in V$. Since $(L^1)^*=L^\infty$, for a $g\in L^\infty(\R^n\times V)$ we have $\|g(x,\cdot)\|_{L^\infty(V)}=\sup_{\|u\|_{L^1}=1}\langle g(x,\cdot),u\rangle_{L^\infty, L^1}$. Therefore
\begin{equation*}
    \|f\|_{\Co^\alpha L^\infty(\R^n,V)}\approx\sup_{\|u\|_{L^1}=1}\sup_{j\ge0}2^{j\alpha}\sup_{x\in\R^n}\langle\phi_j\ast_xf(x,\cdot),u\rangle_{L^\infty, L^1(V)}=\sup_{j\ge0}2^{j\alpha}\|\phi_j\ast_xf\|_{L^\infty(\R^n\times V)}.
\end{equation*}
Here we use that fact that $\phi_j\ast_xf$ is continuous along $x$-variable, which means the sup and the essential sup of $x$ are the same.

\smallskip
\noindent\ref{Item::Hold::BiHoldChar::CalphaCbeta}: Since $\beta>0$, by Lemma \ref{Lem::Hold::HoldChar} \ref{Item::Hold::HoldChar::LPHoldChar} again we have
\begin{equation*}
    \|f\|_{\Co^\alpha\Co^\beta(\R^n,\R^q)}\approx\sup_{k\ge0}\sup_{s\in\R^q}2^{k\beta}\|\psi_k\ast_sf(\cdot,s)\|_{\Co^\alpha(\R^n)}\approx\sup_{j,k\ge0}2^{j\alpha+k\beta}\|\psi_k\ast_s(\phi_j\ast_x f)\|_{L^\infty(\R^n\times\R^q)}.
\end{equation*}
This completes the proof since
\begin{equation}\label{Eqn::Hold::ConvolutionCommute}
    \psi_k\ast_s(\phi_j\ast_x f)=\phi_j\ast_x(\psi_k\ast_s f)=(\phi_j\otimes\psi_k)\ast f.
\end{equation}

The proof of \ref{Item::Hold::BiHoldChar::Domain} follows immediately from Lemma \ref{Lem::Hold::CommuteExt}.

Finally to prove \ref{Item::Hold::BiHoldChar::CommCAlphaLInfty}, clearly by definition, the space $\Co^\alpha(U;L^\infty(V))$ and \eqref{Eqn::Hold::BiHoldChar::CommCAlphaLInfty::Tmp} are equal.
By \ref{Item::Hold::BiHoldChar::CalphaLinfty} and Lemma \ref{Lem::Hold::HoldChar} \ref{Eqn::Hold::HoldChar::LPHoldChar} we see that $\Co^\alpha L^\infty(\R^n,V)=\Co^\alpha(\R^n;L^\infty(V))=\{f\in L^\infty(\R^n\times V):\essup_{s\in V}\|f(\cdot,s)\|<\infty\}$ with equivalent norms. Using \ref{Item::Hold::BiHoldChar::Domain} we can replace $\R^n$ by bounded smooth domain $U$.
\end{proof}
\begin{remark}\label{Rmk::Hold::RmkforBiHold::notVectMeas} Now $f\in L^\infty\Co^\alpha(U,V)$ is saying that $f(\cdot,s)\in\Co^\alpha(U)$ almost every $s\in V$ and is essentially bounded. However $L^\infty\Co^\alpha(U,V)$ does not coincide with the vector-valued space $L^\infty(V;\Co^\alpha( U))$. In classical definitions $L^\infty(U;\Xs)$ is the set of all $\Xs$-valued strong measurable function (quotiented by almost everywhere) which is essentially bounded. Such functions in $L^\infty(U;\Co^\alpha(V))$ must be \textit{almost everywhere separable valued}. But $\Co^\alpha(U)$ is \textit{not} a separable space. So $L^\infty(U;\Co^\alpha(V))\subsetneq L^\infty\Co^\alpha(U,V)$ and the inclusion is strict, and moreover $L^\infty\Co^\alpha(U,V)\not\subset L^1(V;\Co^\alpha(U))$. See for example \cite[Remark 1.2.16]{AnalBana}.
\end{remark}
\begin{remark}\label{Rmk::Hold::BiHoldInterpo} Let $\alpha_0,\alpha_1\in\R$ and $\beta_0,\beta_1>0$. 
    Clearly $\min(2^{j\alpha_0+k\beta_0},2^{j\alpha_1+k\beta_1})\le2^{j((1-\theta)\alpha_0+\theta\alpha_1)+k((1-\theta)\beta_0+\theta\beta_1)}$ for $\theta\in[0,1]$. Thus we have inclusion: for $U,V$ either the total space or bounded smooth domains, 
    \begin{equation}\label{Eqn::Hold::RmkforBiHold::Interpo}
        \Co^{\alpha_0}\Co^{\beta_0}\cap\Co^{\alpha_1}\Co^{\beta_1}(U,V)\subset \Co^{\alpha_\theta}\Co^{\beta_\theta}(U,V),\quad\alpha_\theta:=(1-\theta)\alpha_0+\theta\alpha_1,\ \beta_\theta:=(1-\theta)\beta_0+\theta\beta_1.
    \end{equation}
    In particular $\Co^{\alpha}L^\infty\cap \Co^0\Co^\beta(U,V)\subset\Co^{(1-\theta)\alpha}\Co^{\theta\beta}(U,V)$ for all $\alpha,\beta>0$.
\end{remark}

As a corollary to Lemma \ref{Lem::Hold::BiHoldChar}, the paraproduct operators can be adapted to bi-parameter H\"older spaces.
\begin{lem}\label{Lem::Hold::VectPara}
Let $(s^1,\dots,s^q)$ be the standard coordinate system for $\R^q$ and let $\phi=(\phi_j)_{j=0}^\infty$ be a dyadic resolution for $\R^q$. 

$U,V,W$ be three open sets and let $\Pi:\Co^{\alpha_1}(U)\times\Co^{\alpha_2}(V)\to\Co^{\alpha_3}(W)$ be a bounded bilinear operator where $\alpha_1,\alpha_2,\alpha_3\in\R$. Then the $s$-variable paraproduct operators 
\begin{equation*}
    \Pf^\Pi(f,g):=\sum_{k=0}^\infty\sum_{k'=0}^{k-3}\Pi(\phi_k\ast_sf,\phi_{k'}\ast_s g),\quad \Rf^\Pi(f,g):=\sum_{|k-k'|\le2}\Pi(\phi_k\ast_sf,\phi_{k'}\ast_s g),
\end{equation*}
define bilinearly bounded maps $\Pf^\Pi,\Rf^\Pi:\Co^{\alpha_1} \Co^\beta(U,\R^q)\times \Co^{\alpha_2} L^\infty(V,\R^q)\to\Co^{\alpha_3}\Co^\beta(W,\R^q)$, for every $\beta>0$.
\end{lem}
\begin{proof}
Direct computation using Lemmas \ref{Lem::Hold::LemParaProd} and \ref{Lem::Hold::ParaBdd}.
\end{proof}

Next we define the mixed H\"older-Zygmund space $\Co^{\alpha,\beta}_{x,s}$ as the following
\begin{defn}\label{Defn::Hold::MixHold}
Let $\alpha,\beta\in\R_\Eb^+$ be two generalized indices. Let $U\subseteq\R^n$ and $V\subseteq\R^q$ be two open sets, and let $\Xs$ be a Banach space. We define $\Co^{\alpha,\beta}_{x,s}(U,V;\Xs)$ to be the set of all continuous functions $f:U\times V\to\Xs$ such that $\{f(x,\cdot):x\in U\}\subset \Co^\beta(V;\Xs)$ and $\{f(\cdot,s):s\in V\}\subset\Co^\alpha(U;\Xs)$ are two bounded sets with respect to the topologies of $\Co^\beta(V;\Xs)$ and $\Co^\alpha(U;\Xs)$ respectively.

When $\alpha,\beta\in\R_+$, we use the norm,
$$\textstyle\|f\|_{\Co^{\alpha,\beta}_{x,s}(U,V;\Xs)}:=\sup_{x\in U}\|f(x,\cdot)\|_{\Co^\beta(V;\Xs)}+\sup_{s\in V}\|f(\cdot,s)\|_{\Co^\alpha(U;\Xs)}.$$

We use abbreviation $f\in \Co^{\alpha,\beta}_{x,s}$ when the domain is clear. We use $f\in\Co^{\alpha,\beta}$ when the variables are also clear.
\end{defn}
\begin{remark}
    \begin{enumerate}[parsep=-0.3ex,label=(\roman*)]
        \item The norm definition is symmetric and we have $\Co^{\alpha,\beta}_{x,s}=\Co^{\beta,\alpha}_{s,x}$. But in applications we mostly refer the left variable as space variable, and the right variable as parameter. For the notation $\Co^{\alpha,\beta}$, we usually use $\alpha\ge\beta$.
        \item When $\alpha=\beta$, we have $\Co^{\alpha,\alpha}=\Co^\alpha$. See Lemma \ref{Lem::Hold::CharMixHold} \ref{Item::Hold::CharMixHold::HoldbyComp} below. So our $\Co^{\alpha,\beta}$-spaces is a generalization to the single variable cases. 
        \item We can define three variables mixed H\"older-Zygmund spaces and the version on manifold analogously. See Definition \ref{Defn::Hold::MoreMixHold} \ref{Item::Hold::MoreMixHold::3Mix} and Definition \ref{Defn::ODE::MixHoldMaps}. 
    \end{enumerate}
\end{remark}

For $\gamma\in\R_+$, the spaces $\Co^{\gamma-,\beta}$ and $\Co^{\gamma+,\beta}$ are just for convenient purpose:
\begin{lem}\label{Lem::Hold::LimitArgument}
Let $\gamma\in\R_+$ and $\beta\in\R_\Eb^+$. Let $U\subseteq\R^n$ and $V\subseteq\R^q$ be two open subsets. Then $f\in\Co^{\gamma-,\beta}(U,V)$ if and only if $f\in\Co^{\gamma-\eps,\beta}(U,V)$ for all $0<\eps<\gamma$.
\end{lem}
Thus the discussion of $\Co^{\alpha,\beta}$-regularities for positive $\alpha,\beta\in\R_\Eb^+$ can be reduced to the case where $\alpha,\beta\in(0,\infty)\cup\{k+\LogL,k+\Lip:k=0,1,2,\dots\}$.
\begin{proof}

$f\in\Co^{\gamma-,\beta}(U,V)$ if and only if $\{f(\cdot,s):s\in V\}\subset\Co^{\gamma-}(U)$ is bounded, if and only if $\{f(\cdot,s):s\in V\}\subset\Co^{\gamma-\eps}(U)$ is bounded for all $\eps>0$, if and only if $f\in\Co^{\gamma-\eps,\beta}(U,V)$ for all $\eps>0$.
\end{proof}

\begin{lem}\label{Lem::Hold::CharMixHold}
Let $\alpha,\beta\in(0,\infty)$, and let $U\subset\R^m$, $V\subset\R^n$ be two bounded smooth domains. Then
\begin{enumerate}[nolistsep,label=(\roman*)]
    \item\label{Item::Hold::CharMixHold::Mix=Bi} $\Co^{\alpha,\beta}(U,V)=\Co^\alpha L^\infty(U,V)\cap \Co^\beta L^\infty(V,U)$ with equivalent norms.
    \item\label{Item::Hold::CharMixHold::HoldbyComp} $\Co^{\alpha,\alpha}(U,V)=\Co^\alpha(U\times V)$ with  equivalent norms.
\end{enumerate}
\end{lem}
\begin{proof}
\ref{Item::Hold::CharMixHold::Mix=Bi}: By Lemma \ref{Lem::Hold::CommuteExt} we have an extension operator $\tilde E:=E_xE_s=E_sE_x:L^\infty\Co^\beta(U,V)\to L^\infty\Co^\beta(\R^n,\R^q)$, thus by Lemma \ref{Lem::Hold::BiHoldChar} \ref{Item::Hold::BiHoldChar::Domain},
\begin{equation}\label{Eqn::Hold::CharMixHold::Tmp}
    \textstyle\sup_{x\in U}\|f(x,\cdot)\|_{\Co^\beta(V)}\ge\essup_{x\in U}\|f(x,\cdot)\|_{\Co^\beta(V)}=\|E_sf\|_{\Co^\beta L^\infty(\R^q,U)}\approx\|\tilde Ef\|_{\Co^\beta L^\infty(\R^n,\R^q)}\approx\|f\|_{L^\infty\Co^\beta(U,V)}.
\end{equation}
Similarly $\sup_{s\in V}\|f(\cdot,s)\|_{\Co^\alpha(U)}\gtrsim\|f\|_{\Co^\alpha L^\infty(U,V)}$. Therefore $\|f\|_{\Co^{\alpha,\beta}(U,V)}\gtrsim\|f\|_{\Co^\alpha L^\infty\cap L^\infty\Co^\beta(U,V)}$.

For the converse direction, it suffices to show the first inequality in \eqref{Eqn::Hold::CharMixHold::Tmp} is indeed an equality.
By Remark \ref{Rmk::Hold::BiHoldInterpo} we have $\Co^\alpha L^\infty\cap L^\infty\Co^\beta(U,V)\subset\Co^{\frac\alpha2}\Co^\frac\beta2(U,V)\subset C^0(U,V)$. Thus $f\in \Co^\alpha L^\infty\cap L^\infty\Co^\beta(U,V)$ can be pointwise defined, and the essential supremum is indeed the supremum. This completes the proof.

\smallskip
\noindent\ref{Item::Hold::CharMixHold::HoldbyComp}: Clearly $\|f(x,\cdot)\|_{\Co^\alpha(V)}\le\|f\|_{\Co^\alpha(U\times V)}$ and $\|f(\cdot,s)\|_{\Co^\alpha(U)}\le\|f\|_{\Co^\alpha(U\times V)}$ for each $x\in U$ and $s\in V$. Therefore $\|f\|_{\Co^{\alpha,\alpha}(U,V)}\le 2\|f\|_{\Co^\alpha(U\times V)}$.

For the converse, we can choose $(\phi_j)_{j=0}^\infty$ and $(\psi_k)_{k=0}^\infty$ as dyadic resolutions for $\R^n_x$ and $\R^q_s$ respectively, such that $\supp\hat\phi_0\subset B^n(0,\sqrt2)$ and $\supp\hat\psi_0\subset B^q(0,\sqrt2)$ additionally. Define $(\pi_l)_{l=0}^\infty$ as $\pi_0=\phi_0\otimes\psi_0$ and $\pi_l=2^{l(n+q)}\pi_0(2^l\cdot)-2^{(l-1)(n+q)}\pi_0(2^{l-1}\cdot)$ for $l\ge1$. We see that $(\pi_l)_l$ is a dyadic resolutions for $\R^{n+q}$ and
\begin{equation*}
    \textstyle\pi_l=\sum_{\max(j,k)=l}\phi_j\otimes\psi_k=\phi_l\otimes\delta_{0\in\R^q}+\delta_{0\in\R^n}\otimes\psi_l-\sum_{\min(j,k)=l}\phi_j\otimes\psi_k,\quad\forall\ l\ge0.
\end{equation*}
Therefore for $\tilde f\in\Co^\alpha(\R^{n+q})$, by Lemma \ref{Lem::Hold::HoldChar} \ref{Item::Hold::HoldChar::LPHoldChar} and \eqref{Eqn::Hold::ConvolutionCommute},
\begin{align*}
    &\textstyle\|\pi_l\ast \tilde f\|_{L^\infty}\le\|\phi_l\ast_x\tilde f\|_{L^\infty}+\|\psi_l\ast_s\tilde f\|_{L^\infty}+\sum_{j=l}^\infty\|\psi_l\ast_s(\phi_j\ast_x\tilde f)\|_{L^\infty}+\sum_{k=l+1}^\infty\|\phi_l\ast_x(\psi_k\ast_s\tilde f)\|_{L^\infty}
    \\
    \lesssim&_\alpha\|\tilde f\|_{\Co^\alpha_xL^\infty_s}2^{-l\alpha}+\|\tilde f\|_{L^\infty_x\Co^\alpha_s}2^{-l\alpha}+\|\psi_l\|_{L^1}\|\tilde f\|_{\Co^\alpha_xL^\infty_s}\sum_{j=l}^\infty2^{-j\alpha}+\|\phi_l\|_{L^1}\|\tilde f\|_{L^\infty_x\Co^\alpha_s}\sum_{k=l+1}^\infty2^{-k\alpha}\lesssim\|\tilde f\|_{\Co^\alpha_xL^\infty_s\cap L^\infty_x\Co^\alpha_s}2^{-l\alpha}.
\end{align*}
Thus $\|\tilde f\|_{\Co^\alpha(\R^n\times\R^q)}\lesssim_\alpha\|\tilde f\|_{\Co^\alpha L^\infty\cap L^\infty\Co^\alpha(\R^n,\R^q)}$. Restricting to the domain $U\times V$ we finish the proof.
\end{proof}




\begin{remark}\label{Rmk::Hold::SimpleHoldbyCompFact}Let $\alpha,\beta>0$ and let $U,V$ be two bounded smooth domains. Since $\Co^0\subset L^\infty\subset\Co^1$, using Lemma \ref{Lem::Hold::BiHoldChar} we have $\Co^\alpha L^\infty\cap\Co^1\Co^\beta(U,V)\subset \Co^{\alpha,\beta}(U,V)\subset \Co^\alpha L^\infty\cap\Co^0\Co^\beta(U,V)$.
\end{remark}

\begin{remark}[Comparison of different working spaces]\label{Rmk::Hold::ComparisonSpaces}
    The bi-parameter structures used in \cite[Section 2]{Gong} and \cite[Section 4.4]{NijenhuisWoolf} are the $\Cc^{\alpha,\beta}$-spaces, which are similar to our $\Co^{\alpha,\beta}$-spaces. We recall them here:
    
    \textit{Let $k\ge l\ge0$ be integers, let $r,s\in[0,1)$ and let $U\subseteq\R^n_x$, $V\subseteq\R^q_s$ be two open (bounded smooth) subsets. The space $\Cc^{k+r,l+s}(U,V)$ is the set of all continuous functions $f:U\times V\to\R$ whose norm below is finite:}
    \begin{equation*}
        \textstyle\|f\|_{\Cc^{k+r,l+s}(U,V)}:=\sum_{j=0}^l\sum_{i=0}^{k-j}\big(\sup_{s\in V}\|\nabla_x^i\nabla_s^jf(\cdot,s)\|_{C^{0,r}(U)}+\sup_{x\in U}\|\nabla_x^i\nabla_s^jf(x,\cdot)\|_{C^{0,s}(V)}\big).
    \end{equation*}
    \textit{Here $C^{0,\gamma}$ is the standard H\"older space $(0\le\gamma<1)$.}

    One can see that $\Cc^{k+r,l+s}(U,V)=\Co^{k+r,l+s}(U,V)$ when $r,s\neq0$ and $\Cc^{k+r,l+s}(U,V)\subsetneq\Co^{k+r,l+s}(U,V)$ when $r$ or $s=0$. We leave the proof to readers. 
    
    The $\Co^{\alpha,\beta}$-spaces is better than the $\Cc^{\alpha,\beta}$-spaces in some sense. On one hand $\Cc^{\alpha,\beta}$ are only defined for $\lfloor\alpha\rfloor\ge\lfloor\beta\rfloor$, while $\Co^{\alpha,\beta}$ are defined for all positive $\alpha,\beta$. On the other hand, the discussions of $\Cc^{\alpha,\beta}$ are always separated between the integer cases and non-integer cases, while for $\Co^{\alpha,\beta}$-spaces the discussion can be unified as we have complex interpolations,
    \begin{equation}\label{Eqn::Hold::CpxInt}
        [\Co^{\alpha_0,\beta_0}(U,V),\Co^{\alpha_1,\beta_1}(U,V)]_\theta=\Co^{(1-\theta)\alpha_0+\theta\alpha_1,(1-\theta)\beta_0+\theta\beta_1}(U,V),\quad\text{for all }\alpha_0,\alpha_1,\beta_0,\beta_1>0,\quad\theta\in[0,1].
    \end{equation}
    We also leave the proof of \eqref{Eqn::Hold::CpxInt} to readers.
\end{remark}

Here are more definitions.
\begin{defn}\label{Defn::Hold::MoreMixHold}
Let $\alpha,\beta\in\R_\Eb^+$. Let $U\subseteq\R^n$, $V\subseteq\R^q$ be two open sets. We define the following spaces:
\begin{enumerate}[nolistsep,label=(\roman*)]
    \item $\Co^{\alpha,\beta}_{\loc}(U,V):=\{f\in C^0_\loc(U\times V):f\in\Co^{\alpha,\beta}(U',V'),\ \forall U'\Subset U,\ V'\Subset V\}$.
    \item\label{Item::Hold::MoreMixHold::3Mix}Let $\gamma\in\R_\Eb^+$. Let $W\subseteq\R^m$ be an open convex set. We define $\Co^{\alpha,\beta,\gamma}(U,V,W)$ to be the set of all continuous functions $f:U\times V\times W\to\R$ such that $\{f(\cdot,s,t):s\in V,t\in W\}\subset\Co^\alpha(U)$, $\{f(x,\cdot,t):x\in U,t\in W\}\subset\Co^\beta(V)$ and $\{f(x,s,\cdot):x\in U,s\in V\}\subset\Co^\gamma(W)$ are all bounded subsets with respect to their topologies.

    When $\alpha,\beta,\gamma\in(0,\infty)$ are numbers, we use the norm
\begin{align*}
    \|f\|_{\Co^{\alpha,\beta,\gamma}_{x,s,t}(U,V,W)}:=&\sup_{x\in U,s\in V}\|f(x,s,\cdot)\|_{\Co^\gamma(W)}+\sup_{x\in U,t\in W}\|f(x,\cdot,t)\|_{\Co^\beta(V)}+\sup\limits_{s\in V,t\in W}\|f(\cdot,s,t)\|_{\Co^\alpha(U)}.
\end{align*}
\end{enumerate}
\end{defn}

\begin{remark}\label{Rmk::Hold::RmkCAlpBetGam}By Lemma \ref{Lem::Hold::CharMixHold} \ref{Item::Hold::CharMixHold::Mix=Bi} we can say that $\Co^{\alpha,\beta,\gamma}_{x,s,t}=\Co^\alpha_xL^\infty_{s,t}\cap\Co^\beta_sL^\infty_{x,t}\cap\Co^\gamma_tL^\infty_{x,s}=\Co^\alpha_xL^\infty_{s,t}\cap L^\infty_x\Co^{\beta,\gamma}_{s,t}$, where $\Co^\alpha_xL^\infty_{s,t}\cap\Co^\beta_sL^\infty_{x,t}\cap\Co^\gamma_tL^\infty_{x,s} $ and $\Co^\alpha_xL^\infty_{s,t}\cap L^\infty_x\Co^{\beta,\gamma}_{s,t}$ are defined analogously by Definitions \ref{Defn::Hold::BiHold} and \ref{Defn::Hold::MixHold}.
\end{remark}

For $0<\alpha<\infty$, $0<\beta<\infty$, that $f\in\Co^{\alpha+1,\beta}_{x,s}$ does not imply $\nabla_x f\in\Co^{\alpha,\beta}_{x,s}$.
\begin{example}\label{Exam::Hold::GradMixHold}
Let $\alpha,\beta>0$. We consider the following $f:\R_x\times\R_s\to\C$:
$$f(x,s):=\sum_{l=1}^\infty 2^{-l(\alpha+1)}\exp(i2^{l}x+i2^{\frac{\alpha+1}\beta l}s).\quad\text{ We have }\frac{\partial f}{\partial x}(x,s)=\sum_{l=1}^\infty i2^{-l\alpha}\exp(i2^{l}x+i2^{\frac{\alpha+1}\beta l}s).$$

One can check that $f\in \Co^{\alpha+1}_xL^\infty_s$ and $f\in L^\infty_x\Co^\beta_s$. So $f\in\Co^{\alpha+1,\beta}_{x,s}$.

Let $\psi=(\psi_k)_{k=0}^\infty$ be a dyadic resolution. By \eqref{Eqn::Hold::RmkDyaSupp} $\supp\hat\psi_k=\supp2^{k-1}\hat\psi_1\subsetneq(-2^{k+1},-2^{k-1})\cup(2^{k-1},2^{k+1})$ for $k\ge1$, so there is an $\eps_0>0$ such that $\supp\hat\psi_k\subseteq 2^k\cdot\big([-2+\eps_0,-\frac12-\eps_0]\cup[\frac12+\eps_0,2-\eps_0]\big)$. Note that there are infinitely many integer pairs $(k,l)$ such that $|k-\frac{\alpha+1}\beta l|<\frac{\eps_0}2$, and for such a pair $(\phi_k\ast_sf)(x,s)=i2^{-l\alpha}\exp(i2^{l}x+i2^{\frac{\alpha+1}\beta l}s)$. 

So  $\|\psi_k\ast_sf\|_{L^\infty(\R^2)}\approx2^{-k\frac{\alpha\beta}{\alpha+1}}$ holds for infinite many $k\ge0$, which means $f\notin L^\infty_x\Co^{\frac{\alpha\beta}{\alpha+1}+\eps}_s$ for any $\eps>0$.
\end{example}
In general, for $0<\alpha,\beta<\infty$ and a function $f\in\Co^{\alpha+1,\beta}_{x,s}$, the best we can guarantee is $\nabla_xf\in\Co^{\alpha,\frac{\alpha\beta}{\alpha+1}}_{x,s}$, where the index $\frac{\alpha\beta}{\alpha+1}$ is optimal.
For applications we need the following endpoint optimal result.
\begin{lem}\label{Lem::Hold::GradMixHold}
Let $\alpha,\beta\in(0,\infty)$ and let $U\subseteq\R^n_x$, $V\subseteq\R^q_s$ be two bounded convex smooth domains. Then $\nabla_x:\Co^{\alpha+1,\beta}(U,V)\to\Co^{\alpha,\frac{\alpha\beta}{\alpha+1}-}(U,V;\R^n)$ is bounded.
In particular, if $f\in\Co^{\infty,\beta}(U,V)$, then $\nabla_xf\in\Co^{\infty,\beta-}(U,V;\R^n)$.
\end{lem}
\begin{proof}
Clearly $\nabla_x:\Co^{\alpha+1,\beta}_{x,s}\subset \Co^{\alpha+1}_xL^\infty_s\cap \Co^0_x\Co^\beta_s\to \Co^\alpha_xL^\infty_s\cap \Co^{-1}_x\Co^\beta_s$ is bounded. By Remark \ref{Rmk::Hold::BiHoldInterpo} we have embedding $\Co^\alpha_xL^\infty_s\cap \Co^{-1}_x\Co^\beta_s\subset\Co^{(1-\theta)\alpha-\theta}_x\Co^{\theta\beta}_s$ for all $\theta\in(0,1)$. Thus $\Co^\alpha_xL^\infty_s\cap \Co^{-1}_x\Co^\beta_s\subset\bigcap_{0<\theta<\frac\alpha{\alpha+1}}\Co^{(1-\theta)\alpha-\theta}_x\Co^{\theta\beta}_s\subset L^\infty_x\Co^{\frac{\alpha\beta}{\alpha+1}-}_s$. This proves the boundedness $\nabla_x:\Co^{\alpha+1,\beta}_{x,s}\to\Co^{\alpha,\frac{\alpha\beta}{\alpha+1}-}_{x,s}$.

Let $\alpha\to+\infty$ we get $\nabla_x:\Co^{\infty,\beta}_{x,s}\to\Co^{\infty,\beta-}_{x,s}$.
\end{proof}

On the other hand, if $f(x,s)$ is, say a harmonic function in $x$, then the assumption $f\in\Co^{\alpha+1,\beta}_{x,s}$ can imply $\nabla_xf\in\Co^{\alpha,\beta}_{x,s}$. In fact we have $f\in\Co^\infty_x\Co^\beta_s$ by the following.
\begin{lem}\label{Lem::Hold::NablaHarm}
Let $\beta\in\R_\Eb^+$, and let $U\subset\R^n_x$ and $V\subset\R^q_s$ be two bounded convex open sets with smooth boundaries. Suppose $f\in \Co^{-1}\Co^\beta(U,V)$ satisfies $\Delta_xf=0$. Then $f\in\Co^\infty_\loc\Co^\beta(U,V)$. In particular $f\in\Co^{\beta+1,\beta}_\loc(U,V)$ and $\nabla_xf\in\Co^\beta_\loc(U\times V;\R^n)$.

As a corollary, let $\tilde U\subseteq\C^m_z$ be an open set. If $g\in \Co^\beta_\loc(\tilde U\times V;\C)$ is holomorphic in $z$, then $\partial_zg\in\Co^\beta_\loc(\tilde U\times V;\C^m)$ holds.
\end{lem}

Here for $\beta\in\R_\Eb^+\backslash\R_+$, we use the space $\Co^\alpha\Co^\beta(U,V):=\Co^\beta(V;\Co^\alpha(U))$ and $\Co^\alpha_\loc\Co^\beta(U,V):=\Co^\beta(V;\Co^\alpha_\loc(U))$ adapted from Definition \ref{Defn::Hold::BiHold} \ref{Item::Hold::BiHold::CalphaCbeta}.

\begin{proof}Since the result is local in $x\in U$, it is enough to consider the case where $U=\B^n$ is the unit ball. By Definition \ref{Defn::Intro::DefofHold} \ref{Item::Intro::DefofHold::-} it suffices to consider the case $\beta\in\R_+\cup\{k+\LogL,k+\Lip:k=0,1,2,\dots\}$.

Applying the spherical average formula, (or by Lemma \ref{Lem::SecHolLap::HLLem} \ref{Item::SecHolLap::HLLem::Harm}) we have
\begin{equation*}
    |\nabla^kh(x)|\lesssim_{n,k}(1-|x|)^{-k-1}\|\nabla^kh\|_{\Co^{-k-1}(\B^n)}\lesssim_k(1-|x|)^{-k-1}\|h\|_{\Co^{-1}(\B^n)},\quad\forall h\in\Co^{-1}(\B^n),\ \Delta h=0.
\end{equation*}

When $\beta\in(0,1)$, for $s_1,s_2\in V$, by taking $h_{s_1,s_2}(x)=\frac{f(x,s_1)-f(x,s_2)}{|s_1-s_2|^\beta}$ we see that for every $\Delta_xf=0$, $$\|f\|_{\Co^k\Co^\beta(r\B^n,V)}=\sup_{s_1,s_2\in V}\|h_{s_1,s_2}\|_{\Co^k(r\B^n)}\lesssim_k(1-r)^{-k-1}\sup_{s_1,s_2\in V}\|h_{s_1,s_2}\|_{\Co^{-1}(\B^n)}=(1-r)^{-k-1}\|f\|_{\Co^{-1}\Co^\beta(\B^n,V)}.$$

Similarly arguments show that $\|f\|_{\Co^k\Co^\beta(r\B^n,V)}\lesssim_k(1-r)^{-k-1}\|f\|_{\Co^{-1}\Co^\beta(\B^n,V)}$ for $\beta\in[1,\infty)\cup\{k+\LogL,k+\Lip:k=0,1,2,\dots\}$, since all $\Co^\beta$-norms are defined by a supremum over some linear functionals.

Therefore $f\in\Co^k\Co^\beta(r\B^n,V)$ for all $r<1$ and all $k$, which means $f\in\Co^\infty_\loc\Co^\beta(\B^n,V)$. The results $f\in\Co^{\beta+1,\beta}_\loc$ and $\nabla_xf\in\Co^\beta_\loc$ follow immediately.

When $g(z,s)$ is holomorphic in $z$, we know $g$ is complex harmonic in $z$-variable. So we reduce it to the real case by considering $\tilde U\subseteq\R^{2m}$.
\end{proof}

\subsection{Propositions on bi-parameter products}\label{Section::BiHoldMult}
In this subsection we compute the mixed regularity of products using bi-parameter paraproducts. 

As a motivation, in Section \ref{Section::EllipticPara::ExistPDE} we need to solve a non-linear elliptic pde system to find a suitable coordinate change that has best possible component-wise regularity. We start with a coefficient function $A\in \Co^{\alpha,\beta}(U,V)$ and we want to find a solution $H\in\Xs^1(U,V)$ in  a suitable bi-parameter H\"older space $\Xs^1$. To show the existence of $H$ we need spaces $\Xs^1$, $\Xs^0$ and $\Xs^{-1}$ that are as small as possible, and satisfy the following (see Proposition \ref{Prop::EllipticPara::ExistPDE}):
\begin{itemize}[parsep=-0.3ex]
    \item $\nabla_x:\Xs^1(U,V)\to\Xs^0(U,V;\R^n)$, $\nabla_x:\Xs^0(U,V)\to\Xs^{-1}(U,V;\R^n)$ are bounded; and $\Delta_x:\Xs^1(U,V)\to\Xs^{-1}(U,V)$ is right invertible.
    \item We have product maps $\Xs^0(U,V)\times\Xs^{-1}(U,V)\to\Xs^{-1}(U,V)$ and $\Xs^0(U,V)\times\Xs^0(U,V)\to\Xs^0(U,V)$.
    \item $\Xs^0(U,V)\supseteq\Co^{\alpha,\beta}(U,V)$.
\end{itemize}

 In application $\Xs^1$ labels the mixed regularity for the coordinate change, and $\Xs^0$ labels the mixed regularity for the pushforward vector fields which are the generators for a given subbundle.

We shall see that for $k=1,0,-1$, we can pick $\Xs^k=\Co^{\alpha+k}_xL^\infty_s\cap\Co^k_x\Co^\beta_s$ when $\alpha>1$, and  $\Xs^k=\Co^{\alpha+k}_xL^\infty_s\cap\Co^k_x\Co^\beta_s+\Co^{\alpha+k}_x\Co^{(2-1/\alpha)\beta}_s$ when $\frac12<\alpha<1$. Note that $\Co^{\alpha+k}_x\Co^{(2-1/\alpha)\beta}_s\subseteq\Co^{k}_x\Co^{\beta}_s$ when $\alpha>1$.

Firstly we extend Lemma \ref{Lem::Hold::Product} \ref{Item::Hold::Product::Hold1} to the bi-parameter case.
\begin{lem}\label{Lem::Hold::LinftyMult}
Let $U\subseteq\R^n$ and $V\subseteq\R^q$ be either the total spaces or bounded smooth domains. The product map $[(f,g)\mapsto fg]:\Co^\alpha L^\infty(U,V)\times\Co^\beta L^\infty(U,V)\to\Co^{\min(\alpha,\beta)}L^\infty(U,V)$ is bounded for all $\alpha,\beta\in\R$ such that $\alpha+\beta>0$.
\end{lem}
\begin{proof}
By Lemma \ref{Lem::Hold::BiHoldChar} \ref{Item::Hold::BiHoldChar::Domain} it suffices to prove $U=\R^n$ and $V=\R^q$. By symmetry we can assume $\beta\le\alpha$.

Let $(\phi_j)_{j=0}^\infty$ be a dyadic resolution for $\R^n$. By Lemma \ref{Lem::Hold::BiHoldChar} \ref{Item::Hold::BiHold::CalphaLinfty} it suffices to prove $\|\phi_j\ast_x(fg)\|_{L^\infty(\R^n\times\R^q)}\lesssim 2^{-j\beta}\|f\|_{\Co^\alpha L^\infty}\|g\|_{\Co^\beta L^\infty}$. Indeed by Lemma \ref{Lem::Hold::LemParaProd},
\begin{align*}
    \|\phi_l\ast_x(fg)\|_{L^\infty}&\le\|\phi_l\|_{L^1}\bigg(\sum_{j=l-2}^{l+2}\sum_{j'=0}^{j-3}+\sum_{j'=l-2}^{l+2}\sum_{j=0}^{j'-3}+\sum_{\substack{j,j'\ge l-3\\|j-j'|\le 2}}\bigg)\|\phi_j\ast_x f\|_{L^\infty}\|\phi_{j'}\ast_x g\|_{L^\infty}
    \\
    &\lesssim\|f\|_{\Co^\alpha L^\infty}\|g\|_{\Co^\beta L^\infty}\bigg(\sum_{j=l-2}^{l+2}\sum_{j'=0}^{j-3}+\sum_{j'=l-2}^{l+2}\sum_{j=0}^{j'-3}+\sum_{\substack{j,j'\ge l-3\\|j-j'|\le 2}}\bigg)2^{-j\alpha-j\beta}\approx 2^{-j\beta}\|f\|_{\Co^\alpha L^\infty}\|g\|_{\Co^\beta L^\infty}.
\end{align*}
Here we use the assumption $\alpha+\beta >0$ and $\beta\le\alpha$.
\end{proof}

In this subsection, we fix $(\phi_j)_{j=0}^\infty\subset\Sc(\R^n)$ and $(\psi_k)_{k=0}^\infty\subset\Sc(\R^q)$ as two dyadic resolutions for their spaces. For bounded functions $f$ and $g$ defined on $\R^n \times\R^q $, we denote
\begin{equation*}
    f_{jk}:=(\phi_j\otimes \psi_k)\ast f,\quad g_{jk}:=(\phi_j\otimes\psi_k)\ast g,\quad j,k\ge0.
\end{equation*}

We write $\displaystyle fg=\bigg(\sum\limits_{j'+3\le j}+\sum\limits_{j+3\le j'}+\sum\limits_{|j-j'|\le2}\bigg)\bigg(\sum\limits_{k'+3\le k}+\sum\limits_{k+3\le k'}+\sum\limits_{|k-k'|\le2}\bigg)(f_{jk}g_{j'k'})$ as $3\times 3$ sums:

\begin{equation}\label{Eqn::Hold::BigSum}
    \begin{aligned}
    fg=\sum_{\mu,\nu=1}^3 S^{\mu\nu}=&\sum\limits_{j=3}^\infty\sum\limits_{j'=0}^{j-3}\sum\limits_{k=3}^\infty\sum\limits_{k'=0}^{k-3}f_{jk}g_{j'k'}+\sum\limits_{j=3}^\infty\sum\limits_{j'=0}^{j-3}\sum\limits_{k'=3}^\infty\sum\limits_{k=0}^{k'-3}f_{jk}g_{j'k'}+\sum\limits_{j=3}^\infty\sum\limits_{j'=0}^{j-3}\sum\limits_{|k-k'|\le2}f_{jk}g_{j'k'}\\
    &+\sum\limits_{j'=3}^\infty\sum\limits_{j=0}^{j'-3}\sum\limits_{k=3}^\infty\sum\limits_{k'=0}^{k-3}f_{jk}g_{j'k'}+\sum\limits_{j'=3}^\infty\sum\limits_{j=0}^{j'-3}\sum\limits_{k'=3}^\infty\sum\limits_{k=0}^{k'-3}f_{jk}g_{j'k'}+\sum\limits_{j'=3}^\infty\sum\limits_{j=0}^{j'-3}\sum\limits_{|k-k'|\le2}f_{jk}g_{j'k'}\\
    &+\sum\limits_{|j-j'|\le2}\sum\limits_{k=3}^\infty\sum\limits_{k'=0}^{k-3}f_{jk}g_{j'k'}+\sum\limits_{|j-j'|\le2}\sum\limits_{k=3}^\infty\sum\limits_{k=0}^{k'-3}f_{jk}g_{j'k'}+\sum\limits_{|j-j'|\le2}\sum\limits_{|k-k'|\le2}f_{jk}g_{j'k'}.
\end{aligned}
\end{equation}
Here $S^{\mu\nu}(f,g)$,  $1\le \mu,\nu\le 3$, is defined to be the $\mu$-th row $\nu$-th column term in \eqref{Eqn::Hold::BigSum}.


The discussion on $\alpha>1$ follows from the following

\begin{prop}\label{Prop::Hold::Mult}
Let $\alpha,\beta\in(0,\infty)$, $\eta\in(-\alpha,\alpha]$ and $\theta\in[0,\alpha-\eta)$. There is  $C=C_{\alpha,\beta,\eta,\theta}>0$ such that 
$$\|fg\|_{\Co^{\alpha-\theta}L^\infty\cap\Co^{-\eta-\theta}\Co^\beta(\R^n,\R^q)}\le C\|f\|_{\Co^{\alpha}L^\infty\cap\Co^{-\eta}\Co^\beta(\R^n,\R^q)}\|g\|_{\Co^{\alpha-\theta}L^\infty\cap\Co^{-\eta-\theta}\Co^\beta(\R^n,\R^q)},$$
for all $f\in \Co^{\alpha}L^\infty\cap\Co^{-\eta}\Co^\beta(\R^n,\R^q)$ and $g\in\Co^{\alpha-\theta}L^\infty\cap\Co^{-\eta-\theta}\Co^\beta(\R^n,\R^q)$.
\end{prop}

For applications, we use the cases $(\eta,\theta)=(1,0)$ and $(\eta,\theta)=(0,0)$ to prove Corollary \ref{Cor::EllipticPara::EstofH} that the left hand side of \eqref{Eqn::EllipticPara::ExistenceH} behave well. And we use $(\eta,\theta)=(0,1)$ to prove Proposition \ref{Prop::Hold::CompThm} which is required in Corollary \ref{Cor::EllipticPara::EstofH} \ref{Item::EllipticPara::EstofH::Phi}.
\begin{proof}

By Lemma \ref{Lem::Hold::LinftyMult}, since $\alpha-\theta>\eta\ge-\alpha$, we have boundedness $[(f,g)\mapsto fg]:\Co^{\alpha}L^\infty\cap\Co^{-\eta}\Co^\beta(\R^n,\R^q)\times \Co^{\alpha-\theta}L^\infty\cap\Co^{-\eta-\theta}\Co^\beta(\R^n,\R^q)\to \Co^{\alpha-\theta}L^\infty\cap\Co^{-\eta-\theta}\Co^\beta(\R^n,\R^q)$.
Our goal is to show the bilinear control $\|fg\|_{\Co^{-\eta-\theta} \Co^\beta}\lesssim\|f\|_{\Co^\alpha L^\infty \cap\Co^{-\eta} \Co^\beta}\|g\|_{\Co^{\alpha-\theta} L^\infty \cap\Co^{-\eta-\theta} \Co^\beta }$.

%


By assumption $\alpha>\max(0,\eta+\theta)$ and $\eta>-\alpha$. Using Lemma \ref{Lem::Hold::ParaBdd} we have the boundedness:
\begin{equation*}
    \Pf_x:\Co^{-\eta}(\R^n)\times\Co^{\alpha-\theta}(\R^n)\to\begin{cases}\Co^{-\eta}(\R^n)&0\le\theta<\alpha\\\Co^{(-\eta)-\eps}(\R^n)(\forall\eps>0)&\alpha=\theta\\\Co^{\alpha-\eta-\theta}(\R^n)&\alpha<\theta\le\alpha-\eta\end{cases}\quad \subseteq\Co^{-\eta-\theta}(\R^n).
\end{equation*}
Similarly argument shows that for all $\alpha,\beta,\eta,\theta$ that satisfy the assumption:
\begin{gather*}
    \Pf_x,\Rf_x\text{ are bounded }\Co^{\alpha-\theta}\times\Co^{-\eta}\to\Co^{-\eta-\theta},\ \Co^\alpha\times\Co^{-\eta-\theta}\to\Co^{-\eta-\theta},\ \Co^{-\eta-\theta}\times\Co^\alpha\to\Co^{-\eta-\theta}.
\end{gather*}
Thus applying Lemma \ref{Lem::Hold::VectPara} with $\Pi(u,v)\in\{\Pf_x(u,v),\Pf_x(v,u),\Rf_x(u,v)\}$ in the lemma, we get:
\begin{align*}
    S^{11},S^{13},S^{21},S^{23},S^{31},S^{33}&:\Co^{-\eta}\Co^\beta(\R^n,\R^q)\times\Co^{\alpha-\theta}L^\infty(\R^n,\R^q)\to\Co^{-\eta-\theta}\Co^\beta(\R^n,\R^q);
    \\
    S^{12},S^{13},S^{22},S^{23},S^{32},S^{33}&:\Co^\alpha L^\infty(\R^n,\R^q)\times\Co^{-\eta-\theta}\Co^\beta(\R^n,\R^q)\to\Co^{-\eta-\theta}\Co^\beta(\R^n,\R^q).
\end{align*}
 Therefore $\|S^{\mu\nu}(f,g)\|_{\Co^{-\eta-\theta}\Co^\beta}\lesssim\|f\|_{\Co^{\alpha-\theta} L^\infty \cap\Co^{-\eta-\theta} \Co^\beta}\|g\|_{\Co^{\alpha-\theta} L^\infty \cap\Co^{-\eta-\theta} \Co^\beta}$ for all $\mu,\nu\in\{1,2,3\}$, finishing the proof.
\end{proof}

\begin{remark}
    The proof above is shorter than the proof in \cite[Proposition 2.30]{YaoCpxFro}.
\end{remark}

\begin{cor}\label{Cor::Hold::CorMult}
Let $\alpha,\beta,\eta,\theta$ satisfy the assumption of Proposition \ref{Prop::Hold::Mult}. Let $U\subset\R^m$ and $V\subset\R^n$ be two bounded open convex subsets with smooth boundaries. Then the following product is bilinearly bounded:
\begin{equation}\label{Eqn::Hold::CorMult}
    [(f,g)\to fg]:\Co^\alpha L^\infty \cap\Co^{-\eta} \Co^\beta(U,V)\times\Co^{\alpha-\theta} L^\infty \cap\Co^{-\eta-\theta} \Co^\beta(U,V)\to \Co^{\alpha-\theta} L^\infty \cap\Co^{-\eta-\theta} \Co^\beta(U,V).
\end{equation}
In particular it is bilinearly bounded between the following spaces:
\begin{enumerate}[nolistsep,label=(\roman*)]
    \item\label{Item::Hold::CorMult::Prin} $\Co^\alpha L^\infty \cap\Co^0 \Co^\beta(U,V)\times\Co^{\alpha-1} L^\infty \cap\Co^{-1} \Co^\beta(U,V)\to \Co^{\alpha-1} L^\infty \cap\Co^{-1} \Co^\beta(U,V)$ for $\alpha>1$.
    \item\label{Item::Hold::CorMult::0} $\Co^\alpha L^\infty \cap\Co^0 \Co^\beta(U,V)\times \Co^\alpha L^\infty \cap\Co^0 \Co^\beta(U,V)\to \Co^\alpha L^\infty \cap\Co^0 \Co^\beta(U,V)$ for $\alpha>0$.
    \item\label{Item::Hold::CorMult::-1} $\Co^\alpha L^\infty \cap\Co^{-1} \Co^\beta(U,V)\times \Co^\alpha L^\infty \cap\Co^{-1} \Co^\beta(U,V)\to \Co^\alpha L^\infty \cap\Co^{-1} \Co^\beta(U,V)$ for $\alpha>1$.
\end{enumerate}

\end{cor}
\begin{proof}
Indeed \ref{Item::Hold::CorMult::Prin}, \ref{Item::Hold::CorMult::0} and \ref{Item::Hold::CorMult::-1} are the cases where $(\eta,\theta)=(0,1)$, $(\eta,\theta)=(0,0)$ and $(\eta,\theta)=(1,0)$ in \eqref{Eqn::Hold::CorMult} respectively.

By Lemma \ref{Lem::Hold::CommuteExt} there is an extension operator such that $\tilde E:\Xs(U,V)\to\Xs(\R^n,\R^q)$ and  $\tilde E:\Ys(U,V)\to\Ys(\R^n,\R^q)$. Here we use $\Xs=\Co^\alpha L^\infty \cap\Co^{-\eta} \Co^\beta$ and $\Ys=\Co^{\alpha-\theta} L^\infty \cap\Co^{-\eta-\theta} \Co^\beta$.
Applying Proposition \ref{Prop::Hold::Mult}, 
\begin{equation*}
    \|fg\|_{\Ys(U,V)}\le\|\tilde Ef\tilde Eg\|_{\Ys(\R^n,\R^q)}\lesssim\|\tilde Ef\|_{\Xs(\R^n,\R^q)}\|\tilde Eg\|_{\Ys(\R^n,\R^q)}\lesssim\|f\|_{\Xs(U,V)}\|g\|_{\Ys(U,V)}.
\end{equation*}
So we obtain \eqref{Eqn::Hold::CorMult}.
\end{proof}

For $\frac12<\alpha<1$, the bi-parameter spaces are unfortunately a bit convoluted.
\begin{prop}\label{Prop::Hold::MultLow}
Let $\frac12<\alpha<1$ and $\beta>0$. Let $U\subseteq\R^n$ and $V\subseteq\R^q$ be either the total spaces or bounded smooth domains. Then the product map
{\small\begin{equation*}
    \Co^{\alpha}L^\infty\cap\Co^0\Co^\beta+\Co^\alpha\Co^{(2-\frac1\alpha)\beta}(U,V)\times\Co^{\alpha-1}L^\infty\cap(\Co^{-1}\Co^\beta+\Co^{\alpha-1}\Co^{(2-\frac1\alpha)\beta}(U,V)\to\Co^{\alpha-1}L^\infty\cap\Co^{-1}\Co^\beta+\Co^{\alpha-1}\Co^{(2-\frac1\alpha)\beta}(U,V)
\end{equation*}}
is bounded bilinear.
\end{prop}
\begin{remark}
\begin{enumerate}[parsep=-0.3ex,label=(\roman*)]
\item There is no need to put parentheses on the spaces: we have $\Co^{\alpha}L^\infty\cap(\Co^0\Co^\beta+\Co^\alpha\Co^{(2-\frac1\alpha)\beta})=(\Co^{\alpha}L^\infty\cap\Co^0\Co^\beta)+\Co^\alpha\Co^{(2-\frac1\alpha)\beta}$ because $\Co^\alpha\Co^{(2-\frac1\alpha)\beta}\subset \Co^{\alpha}L^\infty$. Later it more convenient to consider $\Co^{\alpha}L^\infty\cap(\Co^0\Co^\beta+\Co^\alpha\Co^{(2-\frac1\alpha)\beta})$.

    \item For $\frac12<\alpha<1$, by Remarks \ref{Rmk::Hold::BiHoldInterpo} and \ref{Rmk::Hold::SimpleHoldbyCompFact} $\Co^{\alpha,\beta}\subset\Co^\alpha L^\infty\cap\Co^{\eps+1-\alpha}\Co^{(2-\frac{1+\eps}\alpha)\beta}$ holds for all $0<\eps<2\alpha-1$. Applying Proposition \ref{Prop::Hold::Mult} with $\eta=\alpha-1-\eps$ and $\theta=1$, we have the product map
    $$\Co^\alpha L^\infty\cap\Co^{\eps+1-\alpha}\Co^{(2-\frac{1+\eps}\alpha)\beta}\times \Co^{\alpha-1} L^\infty\cap\Co^{\eps-\alpha}\Co^{(2-\frac{1+\eps}\alpha)\beta}\to \Co^{\alpha-1} L^\infty\cap\Co^{\eps-\alpha}\Co^{(2-\frac{1+\eps}\alpha)\beta},\quad\forall0<\eps<2\alpha-1.$$ 
    \quad Thus for $f,g\in\Co^{\alpha,\beta}(\R^n,\R^q)$ we get $ (-\Delta)^{-\frac12}(f\cdot\nabla g)\in\Co^{\alpha,(2-\frac1\alpha)\beta-}$ which is endpoint optimal. By introducing the term $\Co^\alpha\Co^{(2-\frac1\alpha)\beta}$ in Proposition \ref{Prop::Hold::MultLow}, we get $(-\Delta)^{-\frac12}(f\cdot\nabla g)\in\Co^{\alpha,(2-\frac1\alpha)\beta}$, which is sharp. See Example \ref{Exam::Hold::MultUnBdd}.
\end{enumerate}
\end{remark}
\begin{proof}[Proof of Proposition \ref{Prop::Hold::MultLow}]
Firstly by Lemma \ref{Lem::Hold::VectPara}, we see that $\|fg\|_{\Co^{\alpha-1}L^\infty}\lesssim\|f\|_{\Co^\alpha L^\infty}\|g\|_{\Co^{\alpha-1}L^\infty}$.

It is more convenient to introduce $\Pf'(u,v):=\Pf(v,u)$, thus we have $uv=\Pf(u,v)+\Pf'(u,v)+\Rf(u,v)$.

By Lemma \ref{Lem::Hold::ParaBdd} $\Pf_x,\Pf_x'$ are both bounded as $\Co^{\alpha}\times\Co^{-1}\to\Co^{-1}$ and $\Co^{\alpha-1}\times\Co^0\to\Co^{-1}$. Therefore by Lemma \ref{Lem::Hold::VectPara}, again with  $\Pi(u,v)\in\{\Pf_x(u,v),\Pf_x(v,u),\Rf_x(u,v)\}$ in the lemma,
\begin{align*}
    S^{11},S^{13},S^{21},S^{23}&:\Co^0\Co^\beta(\R^n,\R^q)\times\Co^{\alpha-1}L^\infty(\R^n,\R^q)\to\Co^{-1}\Co^\beta(\R^n,\R^q);
    \\
    S^{12},S^{13},S^{22},S^{23}&:\Co^\alpha L^\infty(\R^n,\R^q)\times\Co^{-1}\Co^\beta(\R^n,\R^q)\to\Co^{-1}\Co^\beta(\R^n,\R^q).
\end{align*}

Similarly, the boundedness $\Pf_x,\Pf_x',\Rf_x:\Co^\alpha\times\Co^{\alpha-1}\to\Co^{\alpha-1}$ gives
\begin{align*}
    S^{11},S^{13},S^{21},S^{23},S^{31},S^{33}&:\Co^\alpha\Co^{(2-\frac1\alpha)\beta}(\R^n,\R^q)\times\Co^{\alpha-1}L^\infty(\R^n,\R^q)\to\Co^{\alpha-1}\Co^{(2-\frac1\alpha)\beta}(\R^n,\R^q);
    \\
    S^{12},S^{13},S^{22},S^{23},S^{32},S^{33}&:\Co^\alpha L^\infty(\R^n,\R^q)\times\Co^{\alpha-1}\Co^{(2-\frac1\alpha)\beta}(\R^n,\R^q)\to\Co^{\alpha-1}\Co^{(2-\frac1\alpha)\beta}(\R^n,\R^q).
\end{align*}

It remains prove the following:
\begin{equation}\label{Eqn::Hold::PfMult::Goal1}
    S^{31},S^{32},S^{33}:\Co^\alpha L^\infty\cap\Co^0\Co^\beta\times\Co^{\alpha-1}L^\infty\cap\Co^{-1}\Co^\beta\to\Co^0\Co^{(2-1/\alpha)\beta}(\subset\Co^{\alpha-1}\Co^{(2-1/\alpha)\beta}).
\end{equation}

The assumption $f\in\Co^\alpha L^\infty$ implies that $\|f_{jk}\|_{L^\infty}\lesssim_{\phi,\psi}\|f\|_{\Co^\alpha L^\infty}2^{-j\alpha}$, and the assumption $f\in\Co^0\Co^\beta$ implies that $\|f_{jk}\|_{L^\infty}\lesssim_{\phi,\psi}\|f\|_{\Co^0\Co^\beta}2^{j-k\beta}$. Similar control works for $g$. To conclude this, we have
\begin{gather}\label{Eqn::Hold::PfMult::AB}
    \|f_{jk}\|_{L^\infty}\lesssim\min(2^{-j\alpha},2^{-k\beta})\|f\|_{\Co^\alpha L^\infty\cap\Co^0 \Co^\beta},\quad\|g_{jk}\|_{L^\infty}\lesssim\min(2^{j(1-\alpha)},2^{j-k\beta})\|g\|_{\Co^{\alpha-1} L^\infty\cap\Co^{-1} \Co^\beta}.
\end{gather}



For convenience we normalize $f$ and $g$ so that $\|f\|_{\Co^\alpha L^\infty\cap\Co^0 \Co^\beta}=\|g\|_{\Co^{\alpha-\theta} L^\infty\cap\Co^{-\eta-\theta} \Co^\beta}=1$. Therefore what we need is to show for $\nu=1,2,3$,
\begin{equation}\label{Eqn::Hold::PfMult::Goal2}
    \textstyle S^{3\nu}_{lm}:=\|(\phi_l\otimes\psi_m)\ast S^{3\nu}(f,g)\|_{L^\infty}\lesssim_{\phi,\psi,\alpha,\beta}2^{-m(2-\frac1\alpha)\beta},\text{ for }l,m\ge0.
\end{equation}

Applying Lemma \ref{Lem::Hold::LemParaProd} on $x$-variable and $s$-variable separately we get
\begin{equation}\label{Eqn::Hold::BiParaProductDecomp}
    (\phi_l\otimes\psi_m)\ast (S^{31}+S^{32}+S^{33})(f,g)=\sum\limits_{\substack{|j-j'|\le2\\j,j'\ge l-3}}\Bigg(\sum\limits_{k=m-2}^{m+2}\sum\limits_{k'=0}^{k-3}+\sum\limits_{k'=m-2}^{m+2}\sum\limits_{k=0}^{k'-3}+\sum\limits_{\substack{|k-k'|\le2\\k,k'\ge m-3}}\Bigg)(f_{jk}g_{j'k'}).
\end{equation}

Notice that $f\in\Co^\alpha L^\infty$ and $g\in\Co^{\alpha-1}L^\infty$ give the following improvement of \eqref{Eqn::Hold::PfMult::AB}:
\begin{align*}
    \textstyle\big\|\sum_{k=1}^mf_{jk}\big\|_{L^\infty}&=\|2^{mq}\psi_0(2^m\cdot)\ast_s(\phi_j\ast_x f)\|_{L^\infty}\lesssim\|\phi_j\ast_xf\|_{L^\infty}\lesssim2^{-j\alpha}\|f\|_{\Co^\alpha L^\infty};
    \\
    \textstyle\big\|\sum_{k=1}^mg_{jk}\big\|_{L^\infty}&=\|2^{mq}\psi_0(2^m\cdot)\ast_s(\phi_j\ast_x g)\|_{L^\infty}\lesssim\|\phi_j\ast_xg\|_{L^\infty}\lesssim2^{j(1-\alpha)}\|g\|_{\Co^{\alpha-1} L^\infty}.
\end{align*}

Therefore for $l,m\ge0$,
\begin{align*}
    S^{31}_{lm}\le&\sum_{\substack{j,j'\ge l-3\\|j-j'|\le l-3}}\sum_{k=m-2}^{m+2}\|f_{jk}\|_{L^\infty}\bigg\|\sum_{k'=0}^{k-3}g_{j'k'}\bigg\|_{L^\infty}\lesssim\sum_{\substack{j,j'\ge l-3\\|j-j'|\le l-3}}\sum_{k=m-2}^{m+2}\min(2^{-j\alpha},2^{-k\beta})2^{j'(1-\alpha)}
    \\
    \approx&\sum_{l'\ge l}\min(2^{-l'\alpha},2^{-m\beta})2^{l'(1-\alpha)}=\sum_{l\le l'\le m\beta/\alpha}2^{l'(1-\alpha)-m\beta}+\sum_{l'> m\beta/\alpha}2^{-l'(2\alpha-1)}\lesssim2^{-m(2-\frac1\alpha)\beta};
    \\
    S^{32}_{lm}\le&\sum_{\substack{j,j'\ge l-3\\|j-j'|\le l-3}}\sum_{k'=m-2}^{m+2}\|g_{j'k'}\|_{L^\infty}\bigg\|\sum_{k=0}^{k'-3}f_{jk}\bigg\|_{L^\infty}\lesssim\sum_{\substack{j,j'\ge l-3\\|j-j'|\le l-3}}\sum_{k'=m-2}^{m+2}\min(2^{j'(1-\alpha)},2^{j'-k'\beta})2^{-j\alpha}
    \\
    \approx&\sum_{l'\ge l}\min(2^{l'(1-\alpha)},2^{l'-m\beta})2^{-l'\alpha}\lesssim2^{-m(2-\frac1\alpha)\beta};
    \\
    S^{33}_{lm}\le&\sum_{\substack{j,j'\ge l-3\\|j-j'|\le l-3}}\sum_{\substack{k,k'\ge m-3\\|k-k'|\le2}}\|f_{jk}\|_{L^\infty}\|g_{j'k'}\|_{L^\infty}\lesssim\sum_{\substack{j,j'\ge l-3\\|j-j'|\le l-3}}\sum_{\substack{k,k'\ge m-3\\|k-k'|\le2}}\min(2^{-j\alpha},2^{-k\beta})\min(2^{j'(1-\alpha)},2^{j'-k'\beta})
    \\
    \approx&\sum_{l'\ge l}\sum_{m'\ge m}\min(2^{-l'(2\alpha-1)},2^{l'-2m'\beta})=\sum_{m'\ge m}\big(\sum_{l\le l'\le m'\beta/\alpha}2^{l'-2m'\beta}+\sum_{l'>m'\beta/\alpha}2^{-l'(2\alpha-1)}\bigg)
    \\
    \lesssim&\sum_{m'\ge m}2^{-m'(2-1/\alpha)\beta}\lesssim2^{-m(2-\frac1\alpha)\beta}.
\end{align*}
We get \eqref{Eqn::Hold::PfMult::Goal2} hence \eqref{Eqn::Hold::PfMult::Goal1}, finishing the proof.
\end{proof}
We should point out that when $\frac12<\alpha<1$, the regularity loss in parameter is inevitable. 
\begin{example}\label{Exam::Hold::MultUnBdd}
For $\alpha,\beta>0$, $\eta>-\alpha$ and $0<\theta<2\alpha$. Let
$$f(x,s):=\sum_{j=1}^\infty 2^{-\frac{j\alpha\beta}{\alpha-\eta}} e^{2\pi i(2^{\frac{j\beta}{\alpha-\eta}}x+2^js)}\sum_{k=1}^{j-1}2^{\frac{\theta\beta}{\alpha-\eta}(k-j)}e^{2\pi i2^ks},\quad g(x,s):=\sum_{j=1}^\infty 2^{-\frac{j(\alpha-\theta)\beta}{\alpha-\eta}} e^{2\pi i(2^{\frac{j\beta}{\alpha-\eta}}x+2^js)}.$$
One can see that $f\in\Co^\alpha L^\infty\cap\bigcap_{\lambda\in\R}\Co^{(1-\lambda)\alpha}\Co^{\lambda\frac{\alpha\beta}{\alpha-\eta}}(\R,\R)$ and $g\in\Co^{\alpha-\theta}L^\infty\cap\bigcap_{\lambda\in\R}\Co^{(1-\lambda)(\alpha-\theta)}\Co^{\lambda\frac{(\alpha-\theta)\beta}{\alpha-\eta}}(\R,\R)$. Taking $\lambda=\frac{\alpha-\eta}\alpha$ and $\lambda=\frac{\alpha-\eta}{\alpha-\theta}$ respectively we get $f\in\Co^\alpha L^\infty\cap\Co^{-\eta}\Co^\beta$ and $g\in\Co^{\alpha-\theta} L^\infty\cap\Co^{-\eta-\theta}\Co^\beta$. 

Let $(\phi_j)_{j=0}^\infty\subset\Sc(\R)$ be a dyadic resolution, we have
\begin{equation*}
    (\phi_0\ast_s(fg))(x,s)=\sum_{j=0}^\infty\sum_{k=0}^{j-1}2^{-j(2\alpha-\theta)\frac\beta{\alpha-\eta}+\frac{\theta\beta(k-j)}{\alpha-\eta}}e^{2\pi i2^ks}=\sum_{k=0}^\infty 2^{-k(2\alpha-\theta)\frac\beta{\alpha-\eta}}e^{2\pi i2^ks}\frac{1}{2^{\frac{2\alpha\beta}{\alpha-\eta}}-1}.
\end{equation*}
Note that the assumption $\theta\le2\alpha$ is necessary otherwise the sum diverges.

Thus $\phi_0\ast_s(fg)\in L^\infty\Co^{(2\alpha-\theta)\frac\beta{\alpha-\eta}}(\R,\R)$ but $fg\notin\Co^{-M}\Co^\gamma(\R,\R)$ for all $M>0$ and $\gamma>(2\alpha-\theta)\frac\beta{\alpha-\eta}$.

Therefore when $\eta+\theta>\alpha$ we see that $\frac{2\alpha-\theta}{\alpha-\eta}<1$, thus $fg\notin \Co^{-M}\Co^\beta(\R,\R)$ for all $M>0$.  This example shows that in order for Proposition \ref{Prop::Hold::Mult} to be true we must assume $\eta+\theta\le\alpha$. 

One can find a more delicate example to show that the case $\eta+\theta=\alpha$ also fails. We leave the detail to readers.

As a special case, when $\eta=0$ and $\theta=1$, we see that $fg\notin\Co^{-M}\Co^{(2-1/\alpha)\beta+}(\R,\R)$ for all $M>0$. Thus the index $(2-1/\alpha)\beta$ in Proposition \ref{Prop::Hold::MultLow} cannot be improved.
\end{example}

In application we need Corollary \ref{Cor::Hold::CorMult} and Proposition \ref{Prop::Hold::MultLow} on matrix-valued functions,
\begin{lem}\label{Lem::Hold::CramerMixedPara}
For any $\alpha,\beta>0$, $0\le\eta<\alpha$ and any bounded smooth domains $U\subset\R^n$, $V\subset\R^q$, there exists a $\tilde c=\tilde c(U,V,\alpha,\beta,\eta,m)>0$ such that
\begin{enumerate}[parsep=-0.3ex,label=(\roman*)]
    \item $\|A^{-1}-I_m\|_{\Co^\alpha L^\infty}\le 2\|A-I_m\|_{\Co^\alpha L^\infty}$, for all $A\in \Co^\alpha L^\infty(U,V;\C^{m\times m})$ such that $ \|A-I_m\|_{\Co^\alpha L^\infty}\le \tilde c$.
    \item $\|A^{-1}-I_m\|_{\Co^\alpha\Co^\beta}\le 2\|A-I_m\|_{\Co^\alpha\Co^\beta}$, for all $A\in \Co^\alpha\Co^\beta(U,V;\C^{m\times m})$ such that $ \|A-I_m\|_{\Co^\alpha\Co^\beta}\le \tilde c$.
    \item $\|A^{-1}-I_m\|_{\Co^{\alpha,\beta}}\le 2\|A-I_m\|_{\Co^{\alpha,\beta}}$, for all $A\in \Co^{\alpha,\beta}(U,V;\C^{m\times m})$ such that $ \|A-I_m\|_{\Co^{\alpha,\beta}}\le \tilde c$.
    \item\label{Item::Hold::CramerMixedPara::Mix} $\|A^{-1}-I_m\|_{\Co^\alpha L^\infty\cap\Co^{-\eta}\Co^\beta}\le 2\|A-I_m\|_{\Co^\alpha L^\infty\cap\Co^{-\eta}\Co^\beta}$, for all $A\in \Co^\alpha L^\infty\cap\Co^{-\eta}\Co^\beta(U,V;\C^{m\times m})$ such that $ \|A-I_m\|_{\Co^\alpha L^\infty\cap\Co^{-\eta}\Co^\beta}\le \tilde c$.
    \item\label{Item::Hold::CramerMixedPara::SingHold} $\|A^{-1}-I_m\|_{\Co^\alpha}\le 2\|A-I_m\|_{\Co^\alpha}$, for all $A\in \Co^\alpha(U;\C^{m\times m})$ with $ \|A-I_m\|_{\Co^\alpha}\le \tilde c$. Here $\tilde c=\tilde c(U,\alpha,m)>0$.
\end{enumerate}
\end{lem}
\begin{proof}These are immediate from Lemmas \ref{Lem::Hold::CramerMixed} and \ref{Lem::Hold::Product}, where  \ref{Item::Hold::CramerMixedPara::Mix}  comes from Corollary \ref{Cor::Hold::CorMult}.
\end{proof}

\subsection{Compositions and inverse mapping with parameters}\label{Section::BiHoldComp}
In this subsection we work on the concept of compositions between parameter dependent maps. We are going to prove an auxiliary result, Proposition \ref{Prop::Hold::CompThm}, that is used in Corollary \ref{Cor::EllipticPara::EstofH} \ref{Item::EllipticPara::EstofH::Phi}. 

We recall that for a bijection $f:U_1\to U_2$, we denote its inverse map by $f^\Inv:U_2\to U_1$.

In this subsection, for maps $F(x,s):U\times V\subseteq\R^n\times\R^q\to\R^n$, usually we would denote by $\tilde F$  the map $\tilde F(x,s):=(F(x,s),s)$. Note that if $F(\cdot,s)$ has inverse for every $s$, then in this notation we have $$\tilde F^\Inv(y,s)=(F(\cdot,s)^\Inv(y),s).$$

Let $x=(x^1,\dots,x^n)$ and $y=(y^1,\dots,y^n)$ be two standard coordinate systems of $\R^n$, and let $s=(s^1,\dots,s^q)$ be the standard coordinate system of $\R^q$. We let $\Ic(x,s):=x$ for $x\in\R^n$ and $s\in\R^q$.

The quantitative composition estimate Proposition \ref{Prop::Hold::QComp} can be generalized to the case with parameters.
\begin{lem}\label{Lem::Hold::QPComp}
Let $m,n,q\ge1$ be integers. Let $\alpha,\beta,\gamma\in\R_+$ satisfy $\max(\alpha,\beta)>1$ and $\max(\beta,\gamma)>1$. For any $M>0$ there is a $K_4(m,n,q,\alpha,\beta,\gamma,M)>0$ satisfying the following:
    
    \smallskip\quad
    Suppose $g\in\Co^{\alpha,\gamma}(\B^n,\B^q;\R^m)$ satisfies $g(\B^n,s)\subseteq \B^m$ for all $s\in\B^q$ and $\|g\|_{\Co^{\alpha,\gamma}(\B^n;\R^m)}\le M$. Then for every $f\in\Co^{\beta,\gamma}(\B^m,\B^q)$, the function $h(x,s):=f(g(x,s),s)$ satisfies $h\in \Co^{\min(\alpha,\beta),\min(\beta,\gamma)}(\B^n,\B^q)$ with $\|h\|_{\Co^{\min(\alpha,\beta),\min(\beta,\gamma)}(\B^n,\B^q)}\le K_4\|f\|_{\Co^{\beta,\gamma}(\B^m,\B^q)}$.
\end{lem}
\begin{proof}
 Applying Proposition \ref{Prop::Hold::QComp} \ref{Item::Hold::QComp::>1}, when $\|g\|_{\Co^{\alpha,\gamma}}=\sup_s\|g(\cdot,s)\|_{\Co^\alpha}+\sup_x\|g(x,\cdot)\|_{\Co^\gamma}<M$ we have,
\begin{equation*}
\begin{aligned}
\|f(g(\cdot,s),s)\|_{\Co^{\min(\alpha,\beta)}(\B^n)}
\le&K_1(m,n,\alpha,\beta,M)\|f(\cdot,s)\|_{\Co^\beta(\B^m)}&\forall s\in\B^q;
\\
\|f(g(x,\cdot),s)\|_{\Co^{\min(\beta,\gamma)}(\B^q)}\le&K_1(m,q,\gamma,\beta,M)\|f(\cdot,s)\|_{\Co^{\beta}(\B^m)}&\forall x\in\B^n,s\in\B^q;
\\
\|f(g(x,s),\cdot)\|_{\Co^{\min(\beta,\gamma)}(\B^q)}\le&\textstyle C(q,\beta,\gamma)\sup_{y\in\B^m}\|f(y,\cdot)\|_{\Co^\gamma(\B^q)}&\forall x\in\B^n,s\in\B^q.
\end{aligned}
\end{equation*}Here $C(q,\beta,\gamma)>0$ is the constant such that $\|u\|_{\Co^{\min(\beta,\gamma)}(\B^q)}\le C\|u\|_{\Co^\gamma(\B^q)}$. Therefore by Lemma \ref{Lem::Hold::CharMixHold} \ref{Item::Hold::CharMixHold::HoldbyComp},
\begin{align*}
    &\|h\|_{\Co^{\min(\alpha,\beta),\min(\alpha,\gamma)}}\lesssim_{q,\beta}\sup_{s\in\B^q}\|h(\cdot,s)\|_{\Co^{\min(\alpha,\beta)}}+\sup_{x\in\B^n,s\in\B^q}(\|f(g(x,\cdot),s)\|_{\Co^{\min(\beta,\gamma)}}+\|f(g(x,s),\cdot)\|_{\Co^{\min(\beta,\gamma)}})
    \\
    \lesssim&_{m,n,q,\alpha,\beta,\gamma,M}\sup_{s\in\B^q}\|f(\cdot,s)\|_{\Co^\beta(\B^m)}+\sup_{s\in\B^q}\|f(\cdot,s)\|_{\Co^\beta(\B^m)}+\sup_{y\in\B^m}(\|f(y,\cdot)\|_{\Co^{\gamma}(\B^q)}\lesssim \|f\|_{\Co^{\beta,\gamma}(\B^m,\B^q)}.
\end{align*}
By keeping track of the constant above we obtain $K_4$.
\end{proof}

The quantitative inverse function theorem, Proposition \ref{Prop::Hold::QIFT} \ref{Item::Hold::QIFT::PhiEst}, can also be generalized to the case with parameters.

\begin{prop}\label{Prop::Hold::QPIFT}
Let $n,q\ge1$, $\alpha>0$ and $\beta\in(0,\alpha+1]$. Write $\Ic:\R^n\times\R^q\to\R^n$ as $\Ic(x,s):=x$. There is a $K_5=K_5(n,q,\alpha,\beta)>1$ that satisfies the following:

\medskip\quad Let $F\in\Co^{\alpha+1,\beta}(\B^n,\B^q;\R^n)$ satisfies $\|F-\Ic\|_{\Co^{\alpha+1,\beta}(\B^n,\B^q;\R^n)}<K_5^{-1}$ and $(F-\Ic)|_{(\partial\B^n)\times\B^q}\equiv0$. Then \begin{enumerate}[nolistsep,label=(\roman*)]
    \item\label{Item::Hold::QPIFT::Fbij} For each $s\in\B^q$, $F(\cdot,s):\B^n\to\B^n$ is bijective.
    \item\label{Item::Hold::QPIFT::Control} Define $\Phi(y,s):=F(\cdot,s)^\Inv(y)$ for $y\in\B^n$, $s\in\B^q$. Then $\Phi\in\Co^{\alpha+1,\beta}(\B^n,\B^q;\R^n)$ and satisfies 
    \begin{equation}\label{Eqn::Hold::QPIFT::Control}
        \|\Phi\|_{\Co^{\alpha+1,\beta}(\B^n,\B^q;\R^n)}\le K_5.
    \end{equation}
\end{enumerate}
\end{prop}

To prove Proposition \ref{Prop::Hold::QPIFT} we need a lemma on inverse function theorem with non-$C^1$ parameters.

\begin{lem}\label{Lem::Hold::QPIFTLow}Let $\alpha>0$ and $0<\beta<\min(\frac{2\alpha+1}{\alpha+1},\frac32)$. Let $F\in\Co^{\alpha+1,\beta}(\B^n_x,\B^q_s;\R^n_y)$ be a map such that $F(\cdot,s)$ is a $\Co^{\alpha+1}$-diffeomorphism onto its image and $F(\B^n,s)\supseteq\B^n$ for each $s\in\B^q$. Then $\Phi(y,s):=F(\cdot,s)^\Inv(y)$ satisfies $\Phi\in\Co^{\alpha+1,\beta}(\B^n_x,\B^q_s;\R^n_y)$.

Moreover, for any $M>1$ there is a $K_5=K_5(n,q,\alpha,\beta,M)>0$ that does not depend on $F$, such that if $\|F\|_{\Co^{\alpha+1,\beta}(\B^n,\B^q;\R^n)}<M$ and $\inf_{x\in \B^n,s\in\B^q}|\det(\nabla_x F(x,s))^{-1}|>M^{-1}$ then $\|\Phi\|_{\Co^{\alpha+1,\beta}(\B^n,\B^q;\R^n)}<K_5$.

\end{lem}
\begin{remark}
    The result is indeed true for all $0<\beta\le\alpha+1$. When $\beta>1$, the proof is simpler by applying the Inverse Function Theorem on the $\Co^\beta$-map $(x,s)\mapsto (F(x,s),s)$. In our case, we deal with the difference quotiences directly.
\end{remark}
\begin{proof}[Proof of Lemma \ref{Lem::Hold::QPIFTLow}]It suffices to prove the quantitative result: we assume $\|F(\cdot,s)\|_{\Co^{\alpha+1}}<M$ and $\inf_{x\in \B^n}|\det(\nabla_x F(x,s))^{-1}|>M^{-1}$ for each $s\in\B^q$. 

Applying \cite[Lemma 5.9]{CoordAdaptedII} we get $\|\Phi(\cdot,s)\|_{\Co^{\alpha+1}(\B^n;\R^n)}\lesssim_{n,\alpha,M}1$. 
Thus to prove the lemma, it suffices to prove $\|\Phi(y,\cdot)\|_{\Co^\beta(\B^q;\R^n)}\lesssim_{n,q,\alpha,\beta,M}1$ uniformly in $y\in\B^n$.
Without loss of generality, we assume $\alpha<1$ from now on, which means $\frac{2\alpha+1}{\alpha+1}<\frac32$.

Let $y\in\B^n$ and $s_0,s_1\in\B^q$. For $j=0,1$, applying mean value theorem to $F(\cdot,s_j)$,
\begin{equation*}\label{Eqn::PfQPIFT::MeanValue}
    \Phi(y,s_0)-\Phi(y,s_1)=\nabla_xF\big((1-\theta_j)\Phi(y,s_0)+\theta_j\Phi(y,s_1),s_j\big)^{-1}\cdot\big(F(\Phi(y,s_0),s_j)-F(\Phi(y,s_1),s_j))\big),\ \exists\theta_j\in(0,1).
\end{equation*}
Now $F\in\Co^{\alpha+1,\beta}(\B^n,\B^q;\R^n)$. By Lemma \ref{Lem::Hold::GradMixHold} (since $\frac{\alpha\beta}{2\alpha+1}<(\frac{\alpha\beta}{\alpha+1})-$), we have $\nabla_xF\in\Co^{\alpha,\frac{\alpha\beta}{2\alpha+1}}(\B^n,\B^q;\R^{n\times n})$ with $\|\nabla_xF\|_{\Co^{\alpha,\frac{\alpha\beta}{2\alpha+1}}_{x,s}}\lesssim_{n,q,\alpha,\beta}\|F\|_{\Co^{\alpha+1,\beta}_{x,s}}\le M$. By \cite[Lemma 5.7]{CoordAdaptedII} (or a modification of Lemma \ref{Lem::Hold::CramerMixed} \ref{Item::Hold::CramerMixed::Prod2}) and the cofactor representation of $(\nabla_xF)^{-1}$, we can find a $C_1=C_1(n,q,\alpha,\beta,M)>0$ such that
\begin{equation}\label{Eqn::PfQPIFT::C1}
    \sup_{s\in\B^q}\|\nabla_xF(\cdot,s)^{-1}\|_{\Co^\alpha(\B^n;\R^{n\times n})}+\sup_{x\in\B^n}\|\nabla_xF(x,\cdot)^{-1}\|_{\Co^\frac{\alpha\beta}{2\alpha+1}(\B^q;\R^{n\times n})}<C_1,
\end{equation}
when the assumptions $\|F\|_{\Co^{\alpha+1,\beta}(\B^n,\B^q;\R^n)}<M$ and $\inf_{x\in \B^n}|\det(\nabla_x F(x,s))^{-1}|>M^{-1}$ are satisfied. 

We separate the discussions between $\beta<1$ and $\beta\ge1$.

\medskip
\noindent\textit{Case $\beta\in(0,1)$}: It suffices to show $\sup_{y,s_0,s_1}|\Phi(y,s_0)-\Phi(y,s_1)|\lesssim_{n,q,\alpha,\beta,M}|s_0-s_1|^\beta$.

We have $y=F(\Phi(y,s_0),s_0)=F(\Phi(y,s_1),s_1)$ since $\Phi(\cdot,s_j)=F(\cdot,s_j)^\Inv$. Therefore
\begin{align*}
    &\textstyle|\Phi(y,s_0)-\Phi(y,s_1)|\le\big(\sup_{\B^n\times\B^q}|(\nabla_xF)^{-1}|\big)\cdot|F(\Phi(y,s_0),s_0)-F(\Phi(y,s_1),s_0))|
    \\
    \le&\textstyle C_1|F(\Phi(y,s_1),s_1)-F(\Phi(y,s_1),s_0))|\le C_1\sup_{x\in\B^n}\|F(x,\cdot)\|_{\Co^\beta(\B^q)}|s_0-s_1|^\beta\le C_1M|s_0-s_1|^\beta.
\end{align*}
This complete the proof of $0<\beta<1$.

\medskip
\noindent\textit{Case $\beta\in[1,\frac{2\alpha+1}{\alpha+1}]$}: It suffices to show $\sup_{y,s_0,s_1}|\Phi(y,s_0)+\Phi(y,s_1)-2\Phi(y,\frac{s_0+s_1}2)|\lesssim_{n,q,\alpha,\beta,M}|s_0-s_1|^\beta$.

Fix a $\tilde y\in\B^n$. For $s_0,s_1\in\B^q$ we denote $s_\frac12:=\frac{s_0+s_1}2$ and $x_i:=\Phi(\tilde y,s_i)$ for $i=0,\frac12,1$. 

Take $\eps:=1-\frac{\alpha+1}{2\alpha+1}\beta$, so $0<\eps<1-\frac\beta{\alpha+1}$. Since $\|F\|_{\Co^{\alpha+1,1-\eps}_{x,s}}\lesssim_\eps\|F\|_{\Co^{\alpha+1,\beta}_{x,s}}\le M$ we have $$\|F\|_{\Co^{\alpha+1,1-\eps}_{x,s}}\le\tilde M,\quad\text{for some }\tilde M=\tilde M(n,q,\alpha,M)>0.$$ 
Let $C_2:=K_5(n,q,\alpha,1-\eps,\tilde M)$. By the proven case $\beta=1-\eps$ above we have $\sup_{x\in\B^n}\|\Phi(x,\cdot)\|_{\Co^{1-\eps}(\B^q;\R^n)}<C_2$. In particular $|x_i-x_j|=|\Phi(\tilde y,s_i)-\Phi(\tilde y,s_j)|\le C_2\cdot|s_i-s_j|^{1-\eps}$.


By mean value theorem we can find $\xi_0$ and $\xi_1$ contained in the triangle with vertices $x_0,x_\frac12,x_1$ such that,
\begin{equation*}
    F(x_0,s_0)-F(x_\frac12,s_0)=\nabla_xF(\xi_0,s_0)\cdot(x_0-x_\frac12),\quad F(x_1,s_1)-F(x_\frac12,s_1)=\nabla_xF(\xi_1,s_1)\cdot(x_1-x_\frac12).
\end{equation*}
Clearly $|\xi_0-\xi_1|\le\max_{i,j=0,\frac12,1}|x_i-x_j|$. Therefore by \eqref{Eqn::PfQPIFT::C1} and the property $|x_i-x_j|\le C_2\cdot|s_i-s_j|^{\beta-1+\eps}$,
\begin{equation}\label{Eqn::PfQPIFT::CTmp}
    |\nabla_xF(\xi_0,s_0)^{-1}-\nabla_xF(\xi_1,s_1)^{-1}|\le C_1(|\xi_0-\xi_1|^\alpha+|s_0-s_1|^\frac{\alpha\beta}{2\alpha+1})\le C_1(C_2^\alpha+1)|s_0-s_1|^{\beta-1+\eps}.
\end{equation}
Here we use the fact that $|s_0-s_1|^\frac{\alpha\beta}{2\alpha+1}=|s_0-s_1|^{\beta-1+\eps}$ and $|\xi_0-\xi_1|^\alpha\le C_2^\alpha|s_0-s_1|^{\alpha(1-\eps)}\le C_2^\alpha|s_0-s_1|^{\beta-1+\eps}$.

Therefore using $\tilde y=F(x_j,s_j)$ for $j=0,\frac12,1$, we have
\begin{align*}
 &x_0+x_1-2x_\frac12=\nabla_xF(\xi_0,s_0)^{-1}\cdot(F(x_0,s_0)-F(x_\frac12,s_0))+\nabla_xF(\xi_1,s_1)^{-1}\cdot(F(x_1,s_1)-F(x_\frac12,s_1))
    \\
    =&\nabla_xF(\xi_0,s_0)^{-1}\cdot(2\tilde y-F(x_\frac12,s_0)-F(x_\frac12,s_1))+\big(\nabla_xF(\xi_0,s_0)^{-1}-\nabla_xF(\xi_1,s_1)^{-1}\big)\cdot(\tilde y-F(x_\frac12,s_1)).
\end{align*}

Applying \eqref{Eqn::PfQPIFT::C1} and \eqref{Eqn::PfQPIFT::CTmp} to the above equality:
\begin{align*}
    &|x_0+x_1-2x_\frac12|\le C_1|2F(x_\frac12,s_\frac12)-F(x_\frac12,s_0)-F(x_\frac12,s_1)|+C_1(C_2^\alpha+1)|s_0-s_1|^{\beta-1+\eps}|F(x_\frac12,s_\frac12)-F(x_\frac12,s_1)|
    \\
    \lesssim&_{q,\beta} C_1\|F\|_{L^\infty\Co^\beta}|s_0-s_1|^\beta+C_1(C_2^\alpha+1)\|F\|_{L^\infty\Co^{1-\eps}}|s_0-s_1|^{\beta-1+\eps+1-\eps}\le C_1(M+(C_2^\alpha+1)\tilde M)|s_0-s_1|^\beta.
\end{align*}
Thus $\sup_{y\in\B^n}\|\Phi(y,\cdot)\|_{\Co^\beta(\B^q)}\lesssim_{n,q,\alpha,\beta,M}1$ and we complete the proof.
\end{proof}


\begin{proof}[Proof of Proposition \ref{Prop::Hold::QPIFT}]We let $K_5$ be a large constant which may change from line to line.

Let $K_2=K_2(n,\alpha,\alpha)>0$ be the constant in Proposition \ref{Prop::Hold::QIFT}. So we can choose $C_1>0$ such that 
\begin{equation*}
    \|F-\Ic\|_{\Co^{\alpha+1,\beta}(\B^n,\B^q;\R^n)}<C_1^{-1}\quad\Rightarrow\quad\|F(\cdot,s)-\id_{\R^n}\|_{\Co^{\alpha+1}(\B^n;\R^n)}<K_2^{-1},\ \forall s\in\B^q.
\end{equation*}
By assumption $(F(\cdot,s)-\Ic(\cdot,s))|_{\partial\B^n}=(F(\cdot,s)-\id_{\R^n})|_{\partial\B^n}\equiv0$, so by Proposition \ref{Prop::Hold::QIFT} $F(\cdot,s):\B^n\to\B^n$ is a $\Co^{\alpha+1}$-diffeomorphism for each $s\in\B^n$ and we have 
\begin{equation}\label{Eqn::Hold::PfQPIFT::C1}
    \textstyle\|F-\Ic\|_{\Co^{\alpha+1,\beta}(\B^n,\B^q;\R^n)}<C_1^{-1}\quad\Rightarrow\quad\sup_{s\in\B^q}\|F(\cdot,s)^\Inv\|_{\Co^{\alpha+1}(\B^n;\R^n)}<K_2.
\end{equation} In particular $\Phi(y,s)$ is pointwise defined for $(y,s)\in\B^n\times\B^q$. This finishes the proof of \ref{Item::Hold::QPIFT::Fbij}.

\medskip
Let $c_0=c_0(\B^n,n,\alpha)\in(0,1)$ be the constant in Lemma \ref{Lem::Hold::CramerMixed} \ref{Item::Hold::CramerMixed::Prod1}, we have, for every $s\in\B^q$,
\begin{equation}\label{Eqn::Hold::PfQPIFT::InvMat}
    \|\nabla_xF(\cdot,s)-I_n\|_{\Co^\alpha(\B^n;\R^{n\times n})}<c_0\ \Rightarrow\ \|(\nabla_x F(\cdot,s))^{-1}-I_n\|_{\Co^\alpha(\B^n;\R^{n\times n})}<2c_0.
\end{equation}

Thus we can find a $C_2\ge\max(2c'_0{}^{-1},C_1,K_2,1)$, such that
\begin{equation}\label{Eqn::Hold::PfQPIFT::C2}
    \|F-\Ic\|_{\Co^{\alpha+1,\beta}(\B^n,\B^q;\R^n)}<C_2^{-1}\ \Rightarrow\ \sup_{x\in\B^n,s\in\B^q}|(\nabla_x F(x,s))^{-1}|\le\sup_{s\in\B^q}\|(\nabla_x F(\cdot,s))^{-1}\|_{\Co^\alpha(\B^n;\R^{n\times n})}<C_2.
\end{equation}

In particular $\inf_{x\in\B^n,s\in\B^q}|\det(\nabla_xF(x,s))^{-1}|>C_2^{-n}$. Thus when $0<\beta<\min(\frac{2\alpha+1}{\alpha+1},\frac32)$ (in particular when $0<\beta\le\frac{3\alpha+1}{2\alpha+1}$), we get $\|\Phi\|_{\Co^{\alpha+1,\beta}}\lesssim_{n,q,\alpha,\beta}1$ by applying Lemma \ref{Lem::Hold::QPIFTLow}.

It remains to prove the case $\beta\ge\frac{3\alpha+1}{2\alpha+1}$. By chain rule, we have
\begin{equation*}
    0_{q\times n}=\nabla_s(F(\Phi(y,s),s))=\big((\nabla_xF)(\Phi(y,s),s)\big)\cdot\nabla_s\Phi(y,s)+(\nabla_sF)(\Phi(y,s),s).
\end{equation*}
Therefore
\begin{equation}
    \label{Eqn::Hold::PfQPIFT::Chain}
    \nabla_s\Phi(y,s)=-((\nabla_xF)^{-1}\cdot\nabla_sF)\circ(\Phi(y,s),s).
\end{equation}
What we need is to show $\|\nabla_s\Phi(y,\cdot)\|_{\Co^{\beta-1}(\B^q;\R^{q\times n})}\lesssim_{n,q,\alpha,\beta}1$. We will proceed by induction on $k=\lceil\beta-\frac{3\alpha+1}{2\alpha+1}\rceil$. The base case $k=1$ has done since we have obtained $K_5(n,q,\alpha,\beta)>0$ for $\beta\le\frac{3\alpha+1}{2\alpha+1}$ from the above discussion.

Suppose the case $k\ge1$ is done. Let $\beta$ satisfies $\lceil\beta-\frac{3\alpha+1}{2\alpha+1}\rceil=k+1$. We denote $\gamma(\beta):=\max(\beta-1,\frac{3\alpha+1}{2\alpha+1})>1$. By induction hypothesis and $\|\Phi\|_{\Co^{\gamma(\beta)}_{(x,s)}}\lesssim \|\Phi\|_{\Co^{\alpha+1,\gamma(\beta)}_{x,s}}$  there is a $C_3=C_3(n,q,\alpha,\beta)>0$ such that
\begin{equation}\label{Eqn::Hold::PfQPIFT::Assume}
    \|F-\Ic\|_{\Co^{\alpha+1,\gamma(\beta)}}\le K_5(n,q,\alpha,\gamma(\beta))^{-1}\quad\Rightarrow\quad\|\Phi\|_{\Co^{\gamma(\beta)}(\B^n\times\B^q;\R^n)}\le C_3.
\end{equation}

 Applying Lemma \ref{Lem::Hold::QPComp} to \eqref{Eqn::Hold::PfQPIFT::Chain}, taking $f$ to be the components of $(\nabla_xF)^{-1}\cdot\nabla_sF\in\Co^{\beta-1}(\B^n\times\B^q;\R^{q\times n})$, we have, when the assumption of \eqref{Eqn::Hold::CompThm::Assum} is satisfied,
\begin{equation}\label{Eqn::Hold::PfQPIFT::Induction}
\begin{aligned}
\|\nabla_s\Phi\|_{\Co^{\beta-1}}
\lesssim&\sup_{x,s}\|((\nabla_xF)^{-1}\cdot\nabla_sF)(\Phi(y,\cdot),s)\|_{\Co^{\beta-1}}+\sup_{x,s}\|((\nabla_xF)^{-1}\cdot\nabla_sF)(\Phi(y,s),\cdot)\|_{\Co^{\beta-1}}
\\
\le&K_1(n,q,\gamma(\beta),\beta-1,C_3)\sup_{s\in\B^q}\|(\nabla_xF)^{-1}\nabla_sF(\cdot,s)\|_{\Co^{\beta-1}(\B^n)}
\\
&+K_1(q,q,\gamma(\beta),\beta-1,\|\id_{\B^q}\|_{\Co^{\gamma(\beta)}})\sup_{x\in\B^n}\|(\nabla_xF)^{-1}\nabla_sF(x,\cdot)\|_{\Co^{\beta-1}(\B^q)}.
\end{aligned}
\end{equation}

Similar to the argument of \eqref{Eqn::Hold::PfQPIFT::C2}, by considering the constants $c_0(\B^n,n,\beta-1)$ and $c_0(\B^q,n,\beta-1)$ in Lemma \ref{Lem::Hold::CramerMixed} \ref{Item::Hold::CramerMixed::Prod1}, we can find a $C_4=C_4(n,q,\alpha,\beta)>\max(C_2,C_3,K_5(n,q,\alpha,\gamma(\beta)))$, such that
\begin{equation*}
    \|F-\Ic\|_{\Co^{\alpha+1,\beta}(\B^n,\B^q;\R^n)}<C_4^{-1}\ \Rightarrow\ \sup_{s\in\B^q}\|(\nabla_xF(\cdot,s))^{-1}\|_{\Co^{\beta-1}(\B^n;\R^{n\times n})}+\sup_{x\in\B^n}\|(\nabla_xF(x,\cdot))^{-1}\|_{\Co^{\beta-1}(\B^q;\R^{n\times n})}<C_4.
\end{equation*}

Therefore $\|(\nabla_xF)^{-1}\nabla_sF\|_{\Co^{\beta-1}(\B^n\times\B^q;\R^{q\times n})}\lesssim_{n,q,\alpha,\beta}1$. Combining this fact with \eqref{Eqn::Hold::PfQPIFT::Induction} we see that $\sup_{y\in\B^n}\|\Phi(y,\cdot)\|_{\Co^\beta(\B^q;\R^n)}\lesssim_{n,q,\alpha,\beta}1$. Thus we get a $K_5(n,q,\alpha,\beta)$ for such $\beta$, finishing the induction and hence the whole proof.
\end{proof}

The result for inverse Laplacian can be adapted to parameter cases as well.

\begin{lem}\label{Lem::Hold::LapInvBdd}
Let $\alpha>-1$. Let $U\subset\R^n_x$ and $V\subset\R^q_s$ be two bounded convex open sets with smooth boundaries. Then the $x$-variable Laplacian $\Delta_x=\sum_{j=1}^n\partial_{x^j}^2$ has a bounded linear inverse $\Pc$ as the following map
$$\Delta_x:\{f\in\Co^{\alpha+1}_x\Xs_s(U,V):f|_{(\partial U)\times V}=0\}\to\Co^{\alpha-1}_x\Xs_s(U,V),\quad\Xs\in\{L^\infty,\Co^\beta:\beta>0\}.$$
In particular $\Pc:\Co^{\alpha-1}\Xs(U,V)\to\Co^{\alpha+1}\Xs(U,V)$ is bounded linear.
\end{lem}
\begin{proof}This is an immediate consequence of Lemma \ref{Lem::Hold::DiriSol} and Lemma \ref{Lem::Hold::OperatorExtension}.
\end{proof}


The following result will be used in Corollary \ref{Cor::EllipticPara::EstofH} \ref{Item::EllipticPara::EstofH::Phi} in Section \ref{Section::EllipticPara::ExistPDE}.

\begin{prop}\label{Prop::Hold::CompThm}
Let $\alpha>1$ and $\beta\in(0,\alpha+1]$. Then for any $\eps>0$ there exists a $\delta=\delta(n,q,\alpha,\beta,\eps)>0$ such that the following is true:

Suppose $F\in\Co^{\alpha+1}L^\infty\cap\Co^1\Co^\beta(\B^n,\B^q;\R^n)$ and $g\in\Co^\alpha L^\infty\cap\Co^0\Co^\beta(\B^n,\B^q)$ satisfy
\begin{equation}\label{Eqn::Hold::CompThm::Assum}
    (F-\Ic)|_{(\partial\B^n)\times\B^q}\equiv0\quad\text{and}\quad\|F-\Ic\|_{\Co^{\alpha+1}L^\infty\cap\Co^1\Co^\beta(\B^n,\B^q;\R^n)}+\|g\|_{\Co^\alpha L^\infty\cap\Co^0\Co^\beta(\B^n,\B^q)}<\delta.
\end{equation}
Then $\tilde F(x,s):=(F(x,s),s)$ is bijective on $\B^n\times\B^q$. 
Moreover, by setting $\tilde\Phi(y,s):=\tilde F^\Inv(y,s)$ for $(y,s)\in\B^n\times\B^q$, we have $g\circ\tilde\Phi\in \Co^\alpha L^\infty\cap\Co^{-1}\Co^\beta(\B^n,\B^q)$ and 
\begin{equation}\label{Eqn::Hold::CompThm::Concl}
    \|g\circ\tilde\Phi\|_{\Co^\alpha L^\infty\cap\Co^{-1}\Co^\beta(\B^n,\B^q)}<\eps.
\end{equation}
\end{prop}

\begin{remark}In Proposition \ref{Prop::Hold::CompThm} we may have $g\circ\tilde\Phi\in\Co^\alpha L^\infty \cap\Co^0\Co^\beta(\B^n,\B^q)$.
    But in applications, the Proposition \ref{Prop::EllipticPara::AnalyticPDE}, our result in Proposition \ref{Prop::Hold::CompThm} is enough.
\end{remark}

\begin{proof}[Proof of Proposition \ref{Prop::Hold::CompThm}]We have $\|F-\Ic\|_{\Co^{\alpha+1,\beta}(\B^n,\B^q;\R^n)}\lesssim\|F-\Ic\|_{\Co^{\alpha+1}L^\infty\cap\Co^1\Co^\beta(\B^n,\B^q;\R^n)}$ by Remark \ref{Rmk::Hold::SimpleHoldbyCompFact}. Therefore there is a $\delta_1>0$ such that if \eqref{Eqn::Hold::CompThm::Assum} is satisfied with $\delta=\delta_1$, then $\|F-\Ic\|_{\Co^{\alpha+1,\beta}(\B^n,\B^q;\R^n)}\le K_5^{-1}$, where $K_5=K_5(n,q,\alpha,\beta)>0$ is the constant in Proposition \ref{Prop::Hold::QPIFT}. Thus by Proposition \ref{Prop::Hold::QPIFT} \ref{Item::Hold::QPIFT::Fbij} $\Phi(y,s)=F(\cdot,s)^\Inv(y)$ is well-defined and $\|\Phi\|_{\Co^{\alpha+1,\beta}}\le K_5$ holds.

Applying Lemma \ref{Lem::Hold::QPComp} along with its constant $K_4$, we have
\begin{equation}\label{Eqn::Hold::PfCompThm::EstComp2}
    \|f\circ\tilde\Phi\|_{\Co^{\alpha+1,\beta}(\B^n,\B^q)}\le K_4(n,n,q,\alpha+1,\alpha+1,\beta)\|f\|_{\Co^{\alpha+1,\beta}(\B^n,\B^q)},\quad\forall  f\in\Co^{\alpha+1,\beta}(\B^n,\B^q).
\end{equation}
By taking $f$ to be the components of $F-\Ic$, we have $\|\Ic-\Phi\|_{\Co^{\alpha+1,\beta}(\B^n,\B^q;\R^n)}\lesssim\|F-\Ic\|_{\Co^{\alpha+1,\beta}(\B^n,\B^q;\R^n)}$. By $\Co^{\alpha+1,\beta}\subset\Co^{\alpha+1}L^\infty\cap\Co^0\Co^\beta$ in Remark \ref{Rmk::Hold::SimpleHoldbyCompFact}  and that $\nabla_y\Ic(y,s)=I_n$, we know $\|\nabla_y\Phi-I_n\|_{\Co^\alpha_y L^\infty_s\cap\Co^{-1}_y\Co^\beta_s}\lesssim\|F-\Ic\|_{\Co^{\alpha+1,\beta}_{x,s}}$. Thus, we can find a $C_2=C_2(n,q,\alpha,\beta)>0$ such that
\begin{equation*}
    \|\nabla_y\Phi-I_n\|_{\Co^\alpha L^\infty\cap\Co^{-1}\Co^\beta(\B^n,\B^q;\R^{n\times n})}\le C_2\|F-\Ic\|_{\Co^{\alpha+1,\beta}(\B^n,\B^q;\R^n)},\text{ when }\eqref{Eqn::Hold::CompThm::Assum}\text{ is satisfied with }\delta=\delta_1.
\end{equation*}

By Lemma \ref{Lem::Hold::CramerMixedPara} \ref{Item::Hold::CramerMixedPara::Mix} (with $\eta=1$), we see that if $C_2\|F-\Ic\|_{\Co^{\alpha+1,\beta}(\B^n,\B^q;\R^n)}$ is small enough, then $(\nabla_y\Phi)^{-1}\in\Co^\alpha L^\infty\cap\Co^{-1}\Co^\beta(\B^n,\B^q;\R^{n\times n})$. Therefore we can find a $\delta_2\in(0,\delta_1)$ such that
\begin{equation}\label{Eqn::Hold::PfCompThm::EstInvMat}
    \|(\nabla_y\Phi)^{-1}\|_{\Co^\alpha L^\infty\cap\Co^{-1}\Co^\beta(\B^n,\B^q;\R^{n\times n})}\le\delta_2^{-1},\text{ when }\eqref{Eqn::Hold::CompThm::Assum}\text{ is satisfied with }\delta=\delta_2.
\end{equation}

Now for $g\in\Co^\alpha L^\infty\cap\Co^0\Co^\beta(\B^n,\B^q)$, we take $G_j=\partial_{x^j}\Pc g$ where $\Pc:\Co^\alpha L^\infty\cap\Co^0\Co^\beta(\B^n,\B^q)\to\Co^{\alpha+2} L^\infty\cap\Co^2\Co^\beta(\B^n,\B^q)$ is the operator in Lemma \ref{Lem::Hold::LapInvBdd}. Thus $G=(G_1,\dots,G_n)\in \Co^{\alpha+1} L^\infty\cap\Co^1\Co^\beta(\B^n,\B^q;\R^n)$ and
\begin{equation*}
    g\circ\tilde\Phi=(\divg_xG)\circ\tilde\Phi,\qquad(\nabla_x G)\circ\tilde\Phi=\nabla_y(G\circ\tilde\Phi)\cdot(\nabla_y\Phi)^{-1}.
\end{equation*}

By Corollary \ref{Cor::Hold::CorMult} \ref{Item::Hold::CorMult::-1} we have
\begin{align*}
    \|g\circ\tilde\Phi\|_{\Co^\alpha L^\infty\cap\Co^{-1}\Co^\beta(\B^n,\B^q)}&\lesssim\|\nabla_y(G\circ\tilde\Phi)\|_{\Co^\alpha L^\infty\cap\Co^{-1}\Co^\beta(\B^n,\B^q;\R^{n\times n})}\|(\nabla_y\Phi)^{-1}\|_{\Co^\alpha L^\infty\cap\Co^{-1}\Co^\beta(\B^n,\B^q;\R^{n\times n})}.
\end{align*}

Therefore by \eqref{Eqn::Hold::PfCompThm::EstInvMat}, \eqref{Eqn::Hold::PfCompThm::EstComp2} and Remark \ref{Rmk::Hold::SimpleHoldbyCompFact}, we have: when \eqref{Eqn::Hold::CompThm::Assum} is satisfied with $\delta=\delta_2$,
\begin{align*}
    \|g\circ\tilde\Phi\|_{\Co^\alpha L^\infty\cap\Co^{-1}\Co^\beta(\B^n,\B^q)}&\lesssim\|G\circ\tilde\Phi\|_{\Co^\alpha L^\infty\cap\Co^0\Co^\beta(\B^n,\B^q;\R^n)}\lesssim\|G\circ\tilde\Phi\|_{\Co^{\alpha+1,\beta}(\B^n,\B^q;\R^n)}
    \\
    &\lesssim\|G\|_{\Co^{\alpha+1,\beta}(\B^n,\B^q;\R^n)}\lesssim\|G\|_{\Co^{\alpha+1} L^\infty\cap\Co^1\Co^\beta(\B^n,\B^q;\R^n)}\lesssim\|g\|_{\Co^{\alpha} L^\infty\cap\Co^0\Co^\beta(\B^n,\B^q)}.
\end{align*}

In other words, there is a $C_3=C_3(n,q,\alpha,\beta)>0$, such that
\begin{equation*}
    \|g\circ\tilde\Phi\|_{\Co^\alpha L^\infty\cap\Co^{-1}\Co^\beta(\B^n,\B^q)}\le C_3\|g\|_{\Co^{\alpha} L^\infty\cap\Co^0\Co^\beta(\B^n,\B^q)},\quad\forall\ g\in \Co^{\alpha} L^\infty\cap\Co^0\Co^\beta(\B^n,\B^q).
\end{equation*}
Now for every $\eps>0$, take $\delta=\min(C_3^{-1}\eps,\delta_2)$, we get \eqref{Eqn::Hold::CompThm::Concl} and complete the proof.
\end{proof}

In applications, we also need a local version of the composition estimate. We recall the notation $\alpha\circ\beta$ in Definition \ref{Defn::Hold::CompIndex} and Corollary \ref{Cor::Hold::CompOp}.
\begin{lem}\label{Lem::Hold::CompofMixHold}
Let $\alpha,\beta,\gamma,\delta\in\R_\Eb^+$ be positive generalized indices. Let $U_1\subseteq\R^n_x$, $V_1\subseteq\R^p_t$, $U_2\Subset\tilde U_2\subseteq\R^q_y$ and $V_2\subseteq\R^r_s$ be five open sets. 
\begin{enumerate}[parsep=-0.3ex,label=(\roman*)]
    \item\label{Item::Hold::CompofMixHold::Comp}  Let $g\in\Co^{\alpha,\beta}_\loc(U_1, V_1;U_2)$ and $f\in\Co^{\gamma,\delta}_\loc(\tilde U_2, V_2)$. We define $h:U_1\times V_1\times V_2\to\R $ as $h(x,t,s):=f(g(x,t),s)$. Then $h\in\Co^{\gamma\circ\alpha,\gamma\circ\beta,\delta}_{\loc}(U_1,V_1, V_2)$

    In particular $h\in\Co^{\min(\alpha,\gamma),\min(\beta,\gamma),\delta}_{\loc}(U_1,V_1, V_2)$ if $\max(\alpha,\gamma)>1$ and $\max(\beta,\gamma)>1$.
    
    \item\label{Item::Hold::CompofMixHold::InvMat} Let $A\in\Co^{\alpha,\beta}_{\loc}(U_1,V_1;\C^{m\times m})$ be a matrix function such that $A(x,s)\in\R^{n\times n}$ is invertible for every $(x,s)$. Then $A^{-1}\in\Co^{\alpha,\beta}_{\loc}(U_1,V_1;\C^{m\times m})$ as well.
    \item\label{Item::Hold::CompofMixHold::InvFun}Assume $\alpha,\beta\in\R_+$ are numbers. Let $F\in\Co^{\alpha+1,\beta}_{\loc}(U_1,V_1;\R^n)$ and $(x_0,t_0)\in U_1\times V_1$. Suppose $\nabla_xF(x_0,t_0)\in \R^{n\times n}$ is an invertible matrix. Then there are neighborhoods $U_1'\subseteq U_1$ of $x_0$, $V_1'\subseteq V_1$ of $t_0$ and $\Omega_1\subseteq\R^n$ of $F(x_0)$ such that
    \begin{itemize}[nolistsep]
        \item For each $t\in V'_1$, $F(\cdot,t):U'_1\to F(U'_1,t)\subseteq\R^n$ is bijective and $F(U'_1,t)\supseteq \Omega_1$.
        \item Set $\Phi(x,t):=F(\cdot,t)^\Inv(x)$, then $\Phi\in\Co^{\alpha+1,\min(\alpha+1,\beta)}_{\loc}(\Omega_1,V'_1;U'_1)$.
    \end{itemize}
    In particular for the map $\tilde F(x,t):=(F(x,t),t))$ we have that $\tilde F:U_1'\times V_1'\to \tilde F(U'_1\times V'_1)$ is homeomorphism.
\end{enumerate}
\end{lem}


\begin{proof}
\ref{Item::Hold::CompofMixHold::Comp}: Replacing $U_1,V_1,\tilde U_2,V_2$ by their precompact subsets, we can assume that $g\in\Co^{\alpha,\beta}(U_1, V_1;U_2)$ and $f\in\Co^{\gamma,\delta}(\tilde U_2, V_2)$. So $\{f(\cdot,s):s\in V_2\}\subset\Co^\gamma(\tilde U_2)$, $\{f(y,\cdot):y\in \tilde U_2\}\subset\Co^\delta(V_2)$, $\{g(\cdot,t):t\in V_1\}\subset\Co^\alpha(U_1;\R^r)$ and  $\{g(x,\cdot):x\in U_1\}\subset\Co^\alpha(V_1;\R^r)$ are all bounded subset.

Therefore by Corollary \ref{Cor::Hold::CompOp}, the following subsets are all bounded in their topologies:
\begin{gather*}
    \{h(x,t,\cdot):x\in U_1,t\in V_1\}\subseteq\{f(y,\cdot):y\in \tilde U_2\}\subset\Co^\delta(V_2);
    \\
    \{h(\cdot,t,s):t\in V_1,s\in V_2\}\subset\Co^{\gamma\circ\alpha}_\loc(U_1);\quad \{h(x,\cdot,s):x\in U_1,s\in V_2\}\subset\Co^{\gamma\circ\beta}_\loc(V_1).
\end{gather*}
This completes the proof of \ref{Item::Hold::CompofMixHold::Comp}.

\smallskip\noindent Proof of \ref{Item::Hold::CompofMixHold::InvMat}: Fix a $(x_0,s_0)\in U_1\times V_1$, it suffices to show $A^{-1}\in\Co^{\alpha,\beta}$ near $(x_0,s_0)$. Thus we can assume that $A\in\Co^{\alpha,\beta}(U_1,V_1;\C^{m\times m})$.

By Lemma \ref{Lem::Hold::Product} \ref{Item::Hold::Product::Hold2} and Remark \ref{Rmk::Hold::Product::Hold2}, for every open set $\Omega$ and every $\mu\in\R_+\cup\{k+\LogL,k+\Lip:k=0,1,2,\dots\}$, $\|B_1B_2\|_{\Co^\mu(\Omega;\C^{m\times m}}\le \|B_1\|_{L^\infty(\Omega;\C^{m\times m})}\|B_2\|_{\Co^\mu(\Omega;\C^{m\times m})}+\|B_1\|_{\Co^\mu(\Omega;\C^{m\times m})}\|B_2\|_{L^\infty(\Omega;\C^{m\times m})}$ for all $B_1,B_2\in\Co^\mu(\Omega;\C^{m\times m})$. Therefore by Lemma \ref{Lem::Hold::CramerMixed},
\ref{Item::Hold::CramerMixed::Prod2},
\begin{equation}\label{Eqn::Hold::CompofMixHold::Tmp}
    B\in \Co^\mu(\Omega;\C^{m\times m}),\ \|B-I_m\|_{L^\infty}<\tfrac12\Rightarrow B^{-1}\in \|B^{-1}-I\|_{\Co^\mu}\le4\|B-I\|_{\Co^\mu},\quad\forall\mu\in\R_+\cup\{k+\LogL,k+\Lip\}.
\end{equation}

By considering $\mu-\eps$ for $\mu\in\R_+$ and letting $\eps\to0^+$ in \eqref{Eqn::Hold::CompofMixHold::Tmp}, we see that for $\sigma\in\{\mu-:\mu\in\R_+\}$:
\begin{enumerate}[parsep=-0.3ex,label=(\alph*)]
    \item For every $B\in \Co^\sigma(\Omega;\C^{m\times m})$ such that $\|B-I_m\|_{L^\infty}<\tfrac12$, we have $B^{-1}\in \Co^\sigma(\Omega;\C^{m\times m})$.
    \item\label{Item::Hold::CompofMixHold::BddInvTmp} Moreover, if $\Us\subset \Co^\sigma(\Omega;\C^{m\times m})$ is a bounded subset (with respect to the $\Co^\sigma$-topology) such that $\Us\subset\{B\in \Co^\sigma(\Omega;\C^{m\times m}): \|B-I\|_{L^\infty}<\tfrac12\}$, then the set $\{B^{-1}:B\in\Us\}\subset \Co^\sigma(\Omega;\C^{m\times m})$ is also bounded.
\end{enumerate}

Now by passing to an invertible linear transform, we can assume that $A(x_0,s_0)=I_m$, thus $\|A-I_m\|_{L^\infty}<\frac12$ in a small neighborhood of $(x_0,s_0)$. By shrinking $U_1$ and $V_1$ (since we only need to prove the result near $(x_0,s_0)$), we assume that $|A(x,s)-I_m|_{\R^{n\times n}}<\frac12$ holds for all $(x,s)\in U_1\times V_1$.

By assumption $\{A(x,\cdot):x\in U_1\}\subset\Co^\beta(V_1;\C^{m\times m})$ and $\{A(\cdot,s):s\in V_1\}\subset\Co^\alpha(U_1;\C^{m\times m})$ are both bounded. Applying the property \ref{Item::Hold::CompofMixHold::BddInvTmp} with $\mu\in\{\alpha,\beta\}$ and $\Omega\in\{U_1,V_1\}$, we see that $\{A(x,\cdot)^{-1}:x\in U_1\}\subset\Co^\beta(V_1;\C^{m\times m})$ and $\{A(\cdot,s)^{-1}:s\in V_1\}\subset\Co^\alpha(U_1;\C^{m\times m})$ are also bounded. Therefore $A^{-1}\in\Co^{\alpha,\beta}(U_1,V_1;\C^{m\times m})$ holds, finishing the proof of \ref{Item::Hold::CompofMixHold::InvMat}.

\smallskip\noindent Proof of \ref{Item::Hold::CompofMixHold::InvFun}: Applying the Inverse Function Theorem on $F(\cdot,t_0)$, we can find a neighborhood $\tilde U_1\subseteq U_1$ of $x_0$ such that $F(\cdot,t_0):\tilde U_1\to F(\tilde U_1,t_0)$ is a $\Co^{\alpha+1}$-diffeomorphism. Shrinking $\tilde U_1$ if necessary, we can assume $F(\cdot,t_0),F(\cdot,t_0)^\Inv\in\Co^{\alpha+1}(\tilde U_1;\R^n)$ both have bounded $C^1$-norm.

By \eqref{Eqn::Hold::RmkforBiHold::Interpo} we have $F\in\Co^{\alpha+1}\Co^0\cap\Co^0\Co^\beta\subset\Co^{\frac\alpha2+1}_x\Co^{\frac{\alpha\beta}{2(\alpha+1)}}_s$, thus $[t\mapsto F(\cdot,t)]:V_1\to C^1(\tilde U_1;\R^n)$ is a continuous map (in fact H\"older) with respect to the $C^1$-norm topology. Since the set of bounded $C^1$-embedding map in $C^1(\tilde U_1;\R^n)$ is open, we can now find a neighborhood $\tilde V_1\subseteq V_1$ of $t_0$ such that $F(\cdot,t)$ is $C^1$-embedding for all $t\in\tilde V_1$. By Inverse Function Theorem again, the inverse map of $F(\cdot,t):\tilde U_1\to F(\tilde U_1,t)$ is also $\Co^{\alpha+1}$.

Take $U'_1:=\tilde U_1$ and take $\Omega_1$ such that $F(x_0,t_0)\in\Omega_1\Subset F(U'_1,t_0)$. By continuity there is a neighborhood $V'_1\subseteq\tilde V_1$ of $t_0$ such that $\Omega_1\Subset F(U'_1,t) $ for all $t\in V'_1$ as well. Now we get $U'_1,V'_1,\Omega_1$ as desired and $\Phi:\Omega_1\times V'_1\to U'_1$ is pointwise defined.

Since $F(\cdot,t)$ is $C^1$-embedding which is homeomorphic to its image, the result $\sup_{t\in V'_1}\|\Phi(\cdot,t)\|_{\Co^{\alpha+1}(\Omega_1;\R^n)}<\infty$ then follows from \cite[Lemma 5.9]{CoordAdaptedII}.

When $\beta\le1$, the result $\sup_{x\in \Omega_1}\|\Phi(x,\cdot)\|_{\Co^\beta(V_1;\R^n)}<\infty$ follows from Lemma \ref{Lem::Hold::QPIFTLow} via a scaling. When $\beta>1$, applying Inverse Function Theorem to the $\Co^{\min(\alpha+1,\beta)}$-map $[(x,t)\in U'_1\times V'_1\mapsto(F(x,t),t)]$, we see that its inverse $[(x,t)\in \Omega_1\times V'_1\mapsto(\Phi(x,t),t)]$ is also $\Co^{\min(\alpha+1,\beta)}$, thus $\sup_{x\in \Omega_1}\|\Phi(x,\cdot)\|_{\Co^{\min(\alpha+1,\beta)}(V_1;\R^n)}<\infty$.

Therefore either case we have $\Phi\in\Co^{\alpha+1}_xL^\infty_t\cap L^\infty_x\Co^{\min(\alpha+1,\beta)}_t=\Co^{\alpha+1,\min(\alpha+1,\beta)}_{x,t}$.

Now $\tilde F$ and $\tilde F^\Inv$ are both $\Co^{\alpha+1,\beta}\subset C^0$ so $\tilde F$ is homeomorphic to its image.
\end{proof}

\section{On Involutive Structures and ODE Flows}
In this section we review basic property of involutive tangent subbundle, give the regularity estimates of ODE flows, and prove the commuting property for commutative log-Lipschitz vector fields.

\subsection{Review of differentiable manifolds and subbundles}
\label{SectionConvDiffGeom}

For $\kappa>1$, we can define $\Co^\alpha$-maps between two $\Co^\kappa$-manifolds for $0<\alpha\le\kappa$, since by Corollary \ref{Cor::Hold::CompOp} we have compositions $\Co^\alpha_\loc\circ\Co^\kappa_\loc\subseteq\Co^\alpha_\loc$ and $\Co^\kappa_\loc\circ\Co^\alpha_\loc\subseteq\Co^\alpha_\loc$. (See also Definition \ref{Defn::DisInv::DefFunVF}.)
\begin{note}\label{Note::ODE::DiffMap}
    Let $\kappa>1$, $\alpha\in(0,\kappa]$, and let $\Mf,\Nf$ be two $\Co^\kappa$-differential manifolds.
    We use $\Co^\alpha_\loc(\Mf;\Nf)$ as the set of all (locally) $\Co^\alpha$-maps from $\Mf$ to $\Nf$. We do not use the topology of this space in the thesis.
\end{note}

In fact we can define $\Co^\alpha$-functions on $\Co^\kappa$-manifold whenever $\alpha>1-\kappa$, see Definition \ref{Defn::DisInv::DefFunVF} and Lemma \ref{Lem::Hold::PushForwardFuncSpaces}.

We can define the mixed regularity $\Co^{\alpha,\beta}$ on product manifolds.

\begin{defn}\label{Defn::ODE::MixHoldMaps}
Let $\kappa>1$ and $\alpha,\beta\in(0,\kappa]$. Let $\Mf^m,\Nf^n,\Pf^q$ be three $\Co^\kappa$-manifolds.

We say a continuous map $f:\Mf\times\Nf\to\Pf$ is $\Co^{\alpha,\beta}_{\loc}$, denoted by $f\in\Co^{\alpha,\beta}_{\loc}(\Mf,\Nf;\Pf)$, if for every $\Co^\kappa$-coordinate charts $\phi:U\subseteq\Mf\to\R^m$, $\psi:V\subseteq\Nf\to\R^n$ and $\rho:W\subseteq\Pf\to\R^q$ that satisfy $f(U\times V)\subseteq W $, we have 
\begin{equation}\label{Eqn::ODE::MixHoldMaps}
    \rho\circ f\circ(\phi^\Inv,\psi^\Inv)\in\Co^{\alpha,\beta}_{\loc}(\phi(U), \psi(V);\R^q).
\end{equation}

\end{defn}

\begin{remark}
\begin{enumerate}[parsep=-0.3ex,label=(\roman*)]
    \item The definition still has symmetry for indices: $\Co^{\alpha,\beta}_\loc(\Mf,\Nf;\Pf)=\Co^{\beta,\alpha}_\loc(\Nf,\Mf;\Pf)$.
    \item To check $f\in\Co^{\alpha,\beta}_\loc$ it suffices to show that \eqref{Eqn::ODE::MixHoldMaps} is valid on a coordinate cover $\{(\phi_i,\psi_i,\rho_i)\}_{i\in I}$ for $\Mf\times\Nf\times f(\Mf\times\Nf)\subseteq\Mf\times\Nf\times\Pf$ such that $f(U_i\times V_i)\subseteq W_i$ for each $i\in I$. The proof can be done using Lemma \ref{Lem::Hold::CompofMixHold} \ref{Item::Hold::CompofMixHold::Comp}. Such cover always exists because of the continuity of $f$. We leave the proof to readers.
    \item In applications we only consider differentiable maps localized near a fixed point. By passing to a precompact coordinate chart near a point, everything is reduced to the case of Euclidean space, and by shrinking the domains every functions or vector fields would have bounded H\"older-Zygmund regularities.
\end{enumerate}
\end{remark}
Definition \ref{Defn::ODE::MixHoldMaps} is applied for labeling component-wise regularity for vector fields  and parameterizations with parameters.
\begin{defn}
    Let $1<\alpha\le\kappa$ and let $\Mf$ be a $m$-dimensional $\Co^\kappa$-manifold. A \textbf{$\Co^\alpha$-regular parameterization} of $\Mf$ is a map $\Phi:\Omega\subseteq\R^m\to \Mf$ such that $\Phi$ is injective and $\Phi^\Inv:\Phi(\Omega)\subseteq \Mf\to\R^m$ is a $\Co^\alpha$-coordinate chart. A \textbf{topological parameterization} of $\Mf$ is a map $\Phi:\Omega\subseteq\R^m\to\Mf$ that is a topological embedding, i.e. $\Phi$ is continuous and injective, and $\Phi^\Inv:\Phi(\Omega)\to\R^m$ is also continuous.
\end{defn}

When $\kappa>2$, the (real) tangent bundle $T\Mf$ of a $\Co^\kappa$-manifold $\Mf$ has a natural $\Co^{\kappa-1}$-differential structure. Its complexification, the complex tangent bundle $\C T\Mf:=T\Mf\otimes\C= T\Mf\oplus iT\Mf$ is also a $\Co^{\kappa-1}$-differential manifold. So we can define (real or complex) $\Co^\alpha$-vector fields when $\alpha\in(0,\kappa-1]$.

Similarly, the cotangent bundle $T^*\Mf$ and the complex cotangent bundle $\C T^*\Mf$ are also $\Co^{\kappa-1}$-manifolds.

\begin{defn}\label{Defn::ODE::CpxSubbd}
    Let $\kappa>1$, $0<\alpha\le\kappa$ and let $\Mf$ be a $\Co^\kappa$-manifold. A \textbf{$\Co^\alpha$-complex tangent subbundle of rank $r$} is a subset $\Se\subseteq\C T\Mf$ satisfying that for any $p\in \Mf$,
    \begin{itemize}[parsep=-0.35ex]
        \item $\Se_p:=\Se\cap\C T_p\Mf$ is a $r$-dimensional complex vector subspace.
        \item There is an open neighborhood $U\subseteq \Mf$ of $p$ and $\Co^\alpha$-complex vector fields $X_1,\dots,X_r$ on $U$, such that $X_1|_q,\dots,X_r|_q$ form a complex linear basis for $\Se_q$ for all $q\in U$.
    \end{itemize}
    We denote it as $\Se\le\C T\Mf$. And we say that $(X_1,\dots,X_r)$ form a \textbf{$\Co^\alpha$-local basis} for  $\Se$ on $U$.
\end{defn}
We use ``subbundle'', or ``regular subbundle'' in Chapter \ref{Chapter::DisInv}, as the abbreviation of (real or complex) tangent subbundle.

The real tangent subbundles can be defined in the same way through replacing every $\C T\Mf$ by $T\Mf$. 
    But we focus mostly on complex subbundles since every real subbundle $\V\le T\Mf$ can be identified as a complex subbundle $\V\otimes\C\le\C T\Mf$.
    
Alternatively a rank $r$ $\Co^\alpha$-subbundle $\Se\le \C T\Mf$ can be equivalently defined as a $\Co^\alpha$-section of the complex Grassmannian bundle $\Gr_\C(r,\C T\Mf)=\coprod_{p\in\Mf}\{\text{rank }r\text{ complex spaces of }\C T_p\Mf\}$. For their equivalence we leave the proof to readers.

In Theorem \ref{KeyThm::EllipticPara} we consider tangent subbundles with mixed regularities.
\begin{defn}\label{Defn::ODE::CpxPaSubbd}
    Let $\kappa>1$, $0<\alpha,\beta\le\kappa-1$ and let $\Mf,\Nf$ be two $\Co^\kappa$-manifolds. A \textbf{$\Co^{\alpha,\beta}$-complex tangent subbundle of rank $r$} is a $\Co^{\min(\alpha,\beta)}$-subbundle $\Se\le \C T(\Mf\times\Nf)$ such that, for any $(u_0,v_0)\in\Mf\times\Nf$ there is a neighborhood $U\times V\subseteq \Mf\times\Nf$ of $(u_0,v_0)$ and complex vector fields $X_1,\dots,X_r\in\Co^{\alpha,\beta}_\loc(U,V;\C T(\Mf\times\Nf))$, such that $X_1|_{(u,v)},\dots,X_r|_{(u,v)}$ form a complex linear basis for $\Se_{(u,v)}$ for all $(u,v)\in U\times V$.
\end{defn}
Equivalently this is saying $\Se\in\Co^{\alpha,\beta}_\loc(\Mf,\Nf;\Gr_\C(r,\C T(\Mf\times \Nf)))$ as a section of the Grassmannian bundle. 

\medskip
Recall the notion of involutive structure in Definition \ref{Defn::Intro::DisInv}. Let $x=(x^1,\dots,x^n)$ be a local coordinate system, for $\Co^{\frac12+}$ vector fields $X=:\sum_{j=1}^nX^j\Coorvec{x^j}$ and $Y=:\sum_{k=1}^nY^k\Coorvec{x^k}$, we recall that
\begin{equation}\label{Eqn::ODE::LieBracketRecap}
    [X,Y]=\sum_{j,k=1}^n\Big(X^j\frac{\partial Y^k}{\partial x^j}-Y^j\frac{\partial X^k}{\partial x^j}\Big)\Coorvec{x^k}.
\end{equation}
When $X,Y$ are $C^1$, $[X,Y]\in C^0$ is pointwise defined. In general \eqref{Eqn::ODE::LieBracketRecap} is a vector field with distributional coefficients. In Proposition \ref{Prop::DisInv::CharSec} we will see that the $\Co^{(-\frac12)+}$-distribution $\langle\theta,[X,Y]\rangle=0$ if and only if $[X,Y]$ is a ``$\Co^{(-\frac12)+}$-vector field'' in section.

For completeness we give definition of the involutivity on subbundles with mixed regularities. 
\begin{defn}\label{Defn::ODE::InvMix}
Let $\kappa>2$, $\alpha\in(\frac12,\kappa-1]$ and $\beta\in(0,\kappa-1]$. Let $\Mf,\Nf$ be two $\Co^\kappa$-manifolds and $\Se$ be a $\Co^{\alpha,\beta}$-complex subbundle such that $\Se\le(\C T\Mf)\times\Nf$. We say $\Se$ is involutive, if for every $C^0$-sections $X,Y:\Mf\times\Nf\to\Se$ of $\Se$ that are $\Co^{\frac12+}_\loc$ in $\Mf$, and for every $C^0$-section $\theta:\Mf\times\Nf\to\Se^\bot$ that are $\Co^{\frac12+}_\loc$ in $\Mf$, the pairing $\langle\theta,[X,Y]\rangle:\Mf\times\Nf\to\C$ is identically zero.

Equivalently, we say $\Se$ is involutive, if $\Se|_{\Mf\times\{s\}}\le \C T\Mf$ is involutive in the sense of Definition \ref{Defn::Intro::DisInv}, for all $s\in\Nf$.
\end{defn}
\begin{remark}\label{Rmk::ODE::InvMix}Definition \ref{Defn::ODE::InvMix} does not require $X,Y$ or $\theta$ to be differentiable on $\Nf$, since when $\Se\le(\C T\Mf)\times\Nf$ the differentiations in \eqref{Eqn::Intro::DefInvEqn} are only taken on $\Mf$ but not on $\Nf$, and by Lemma \ref{Lem::Hold::MultLoc} one have the well-defined product map $C^0_\loc(\Nf;\Co^{\frac12+\eps}_\loc(\Mf))\times C^0_\loc(\Nf;\Co^{\eps-\frac12}_\loc(\Mf))\to C^0_\loc(\Nf;\Co^{\eps-\frac12}_\loc(\Mf))$ for all $\eps>0$. This is critical for the case $\beta\le\frac12$.
\end{remark}

An important fact for involutive structures is that we can pick a collection of good basis.

\begin{lem}[Canonical local basis]\label{Lem::ODE::GoodGen}
	Let $\alpha,\beta\in\R_\Eb^+\cup\{\logl\}$. Let $\Mf$ and $\Nf$ be two smooth-manifolds. Let $\Se$ be a $\Co^{\alpha,\beta}$-complex subbundle over $\Mf\times\Nf$ such that $\Se\le (\C T\Mf)\times\Nf$. Then
	\begin{enumerate}[parsep=-0.3ex,label=(\roman*)]
	\item\label{Item::ODE::GoodGen::InitialCoord} There exist numbers $r,m,q\ge0$. For every $p\in \Mf$ and $\lambda_0\in\Nf$, there is a (mixed real and complex) $\Co^\kappa$-coordinate system $(t,z,s)=(t^1,\dots,t^r,z^1,\dots,z^m,s^1,\dots,s^q):\tilde U\subseteq \Mf\to\R^r\times\C^m\times\R^q$ near $p$, such that 
	\begin{equation}\label{Eqn::ODE::GoodGen::DirectSumAssumption}
	\Se_{(p,\lambda_0)}\oplus\Span(\partial_{\bar z^1}|_p,\dots,\partial_{\bar z^m}|_p,\partial_{s^1}|_p,\dots,\partial_{s^q}|_p)=\C T_p\Mf.
	\end{equation}
	    \item\label{Item::ODE::GoodGen::Uniqueness}Let $(t,z,s):\tilde U\to\R^r\times\C^m\times\R^q$ be a coordinate system such that \eqref{Eqn::ODE::GoodGen::DirectSumAssumption} holds for some $(p,\lambda_0)\in\tilde U\times\Nf$. For any product set $U\times V\subseteq\{(p,\lambda_0)\in\tilde U\times\Nf:\eqref{Eqn::ODE::GoodGen::DirectSumAssumption}\text{ holds}\}$, there exist a unique $\Co^{\min(\alpha,\beta)}$-local basis $X=[X_1,\dots,X_{r+m}]^\top$ for $\Se$ on $U\times V$, such that the collection vectors $X'=[X_1,\dots,X_r]^\top$ and $X''=[X_{r+1},\dots,X_{r+m}]^\top$ are of the form 
	\begin{equation}\label{Eqn::ODE::GoodGen::GoodGenFormula}
	\begin{pmatrix}X'\\X''\end{pmatrix}=\begin{pmatrix}I_r&&A'&B'\\&I_m&A''&B''\end{pmatrix}\begin{pmatrix}\partial/\partial t\\\partial/\partial z\\\partial/\partial{\bar z}\\\partial/\partial s\end{pmatrix},\end{equation}
	where $A'\in\Co^{\alpha,\beta}_\loc(U,V;\C^{r\times m})$, $A''\in\Co^{\alpha,\beta}_\loc(U,V;\C^{m\times m})$, $B'\in\Co^{\alpha,\beta}_\loc(U,V;\C^{r\times q})$, $B''\in\Co^{\alpha,\beta}_\loc(U,V;\C^{m\times q})$.
	    \item\label{Item::ODE::GoodGen::InvComm} If $\alpha>\frac12$ and if $\Se$ is involutive, then $X_1,\dots,X_{r+m}$ in \ref{Item::ODE::GoodGen::Uniqueness} are pairwise commutative.
	    \item\label{Item::ODE::GoodGen::Real} Suppose $\Se=\bar\Se$, then we must have $m=0$, $A'=0$, $A''=0$, $B''=0$, and $B'$ is real-valued.
	\end{enumerate}

\end{lem}
\begin{proof}\ref{Item::ODE::GoodGen::InitialCoord}: Clearly $\Se_{(p,\lambda_0)}\cap\bar\Se_{(p,\lambda_0)}$ and $\Se_{(p,\lambda_0)}+\bar\Se_{(p,\lambda_0)}$ are both complex linear subspaces of $\C T_p\Mf$. We have $r:=\rank(\Se_{(p,\lambda_0)}\cap\bar\Se_{(p,\lambda_0)})$, $m:=\rank\Se_{(p,\lambda_0)}-r$ and $q:=\dim\Mf-\rank(\Se_{(p,\lambda_0)}+\bar\Se_{(p,\lambda_0)})$. Note that $\dim\Mf=r+2m+q$. See \cite[Lemma I.8.5]{Involutive} for details of \eqref{Eqn::ODE::GoodGen::DirectSumAssumption}. 

\medskip\noindent\ref{Item::ODE::GoodGen::Uniqueness}: The formula \eqref{Eqn::ODE::GoodGen::GoodGenFormula} follows from the standard linear algebra argument. See \cite[Lemma 1]{ShortFrobenius} or \cite[Page 18]{Involutive}. Note that by Lemma \ref{Lem::Hold::CompofMixHold} \ref{Item::Hold::CompofMixHold::InvMat}, if $M\in\Co^{\alpha,\beta}_\loc(U,V;\C^{(r+m)\times (r+m)})$ is a matrix function invertible at every point in the domain, then $M^{-1}\in\Co^{\alpha,\beta}_\loc(U,V;\C^{(r+m)\times (r+m)})$ holds as well. This is also true when $\alpha$ or $\beta=\logl$.

We leave the details to readers, where the argument can be found in the proof of Lemma \ref{Lem::SingFro::SingGoodGen} as well.

We also remark that in the statement of \ref{Item::ODE::GoodGen::Uniqueness}, $U\times V$ are arbitrary subset that does not depend on $\alpha,\beta$.

 If $\tilde X=[\tilde X_1,\dots,\tilde X_{r+m}]^\top$ is another local basis that has expression \eqref{Eqn::ODE::GoodGen::GoodGenFormula}, then  $\tilde X-X$ are linear combinations of $\partial_{\bar z}$ and $\partial_s$. But $\Se|_U\cap\Span(\partial_{\bar z},\partial_s)=\{0\}$, so $\tilde X=X$ is the unique collection.

\medskip\noindent\ref{Item::ODE::GoodGen::InvComm}: We use  $(\partial_{u^1},\dots,\partial_{u^{r+m}}):=(\partial_{t^1},\dots,\partial_{t^r},\partial_{z^1},\dots,\partial_{z^m})$, $(\partial_{v^1},\dots,\partial_{v^{m+q}}):=(\partial_{\bar z^1},\dots,\partial_{\bar z^m},\partial_{s^1},\dots,\partial_{s^q})$ and  $\begin{pmatrix}A'&B'\\A''&B''\end{pmatrix}=\left(f_j^k\right)_{\substack{1\le j\le r+m\\1\le k\le m+q}}$ for simplicity.

	In this notation, we see that $(\theta^k)_{k=1}^{m+q}:=(dv^k-\sum_{j=1}^{r+m}f_j^kdu^j)_{k=1}^{m+q}$ is a collection of $\Co^{\alpha,\beta}$ 1-forms that form a basis for $\Se^\bot|_{U\times V}$.
	
	To show that $X_1,\dots,X_{r+m}$ are commutative, by direct computation, for $1\le j,j'\le r+m$,
	\begin{equation}\label{Eqn::ODE::GoodGen::LieBraofX}
	    [X_j,X_{j'}]=\sum_{k=1}^{r+m}\bigg(\frac{\partial f_{j'}^k}{\partial u^j}-\frac{\partial f_j^k}{\partial u^{j'}}+\sum_{k'=1}^{m+q}\Big(f_j^{k'}\frac{\partial f_{j'}^k}{\partial v^{k'}}-f_{j'}^{k'}\frac{\partial f_j^k}{\partial v^{k'}}\Big)\bigg)\Coorvec{v^k}\in\Span\Big(\Coorvec{\bar z},\Coorvec{s}\Big).
	\end{equation}
	
	Therefore $\langle dv^k,[X_j,X_{j'}]\rangle=\langle\theta^k,[X_j,X_{j'}]\rangle=0$ for all $1\le k\le r+m$, which means the coefficients of $[X_j,X_{j'}]$ are all vanishing.
	
\medskip\noindent\ref{Item::ODE::GoodGen::Real}: When $\Se=\bar\Se$, the complex conjugate $\bar X_1,\dots,\bar X_{r+m}$ also form a local basis for $\Se$. By \ref{Item::ODE::GoodGen::Uniqueness} we see that $\bar X'=X'$, so $B'=\bar B'$ is a real-valued matrix.

By \eqref{Eqn::ODE::GoodGen::GoodGenFormula} we have $X''=\partial_z+A''\partial_{\bar z}+B''\partial_s$ and  $\bar X''=\partial_{\bar  z}+\bar A''\partial_{ z}+\bar B''\partial_s$. They are linearly independent since $\partial_z$ and $\partial_{\bar z}$ are linearly independent. So $X_1,\dots,X_{r+m},\bar X_{r+1},\dots,\bar X_{r+m}$ are all sections $\Se$ and are also linearly independent, which means $\rank\Se\ge r+2m$. Since $\rank\Se=r+m$ we must have $m=0$, which means $A'=0$, $A''=0$, $B''=0$ and $B'\in\Co^{\alpha,\beta}(U,V;\R^{r\times q})$.
\end{proof}

We recall the definition of the complex Frobenius structure.
\begin{defn}\label{Defn::ODE::CpxStr}
	Let $\alpha>\frac12$ and let $\Se\le\C T\Mf$ be a $\Co^\alpha$-complex subbundle. We call $\Se$ a \textbf{complex Frobenius structure} of $\Mf$, if $\Se+\bar \Se$ is a subbundle, and both $\Se$ and $\Se+\bar\Se$ are involutive.
\end{defn}
\begin{remark}
In \cite{Gong} $\Se$ is also called a \textit{Levi-flat CR vector bundle} on $\Mf$. When $\Se\cap\bar\Se=\{0\}$ and $\dim \Mf=1+2\rank\Se$, $\Se$ is indeed a classical CR-structure of $\Mf$.
\end{remark}
\begin{defn}\label{Defn::ODE::ConjBundle}
    Let $\Se\le\C T\Mf$ be a complex tangent subbundle, we use $\Se^\bot$ as the \textbf{dual bundle} of $\Se$, and $\bar\Se$ as the \textbf{complex conjugate} of $\Se$. They are given by 
    \begin{equation}\label{Eqn::ODE::DualConjBundle}
        \Se^\bot:=\{(p,\theta)\in\C T^*\Mf:\theta(v)=0,\forall v\in\C T_p\Mf\}\le\C T^*\Mf;\ \bar\Se:=\{(p,\bar v)\in\C T\Mf:(p,v)\in\Se\}\le\C T\Mf.
    \end{equation}
\end{defn}

\begin{remark}\label{Rmk::ODE::RmkInvStr}
\begin{enumerate}[label=(\roman*),parsep=-0.3ex]
    \item One can check that $\bar \Se$ is indeed a subbundle, since $\bar\Se_p\le\C T_p\Mf$ is  a complex subspace (which is closed under the scalar multiplication by $i$) for each $p\in \Mf$.
    \item\label{Item::ODE::RmkInvStr::NotBund} The fact that $\Se\le\C T\Mf$ is a subbundle does not necessarily implies $\Se\cap\bar \Se$ or $\Se+\bar\Se$ is a subbundle. For example, the \textit{Mizohata operator} $X=\partial_t+it\partial_x$ on $\R^2_{t,x}$ generates a rank 1 complex subbundle $\Se\le \C T\R^2$, but $\Se\cap\bar\Se$ and $\Se+\bar\Se$ do not have constant rank: $$(\Se\cap\bar\Se)_{(t,x)}=\begin{cases}0,&t\neq0,\\\C_t\times\{0\},&t=0,\end{cases}\quad (\Se+\bar\Se)_{(t,x)}=\begin{cases}\C_{t,x}^2,&t\neq0,\\\C_t\times\{0\},&t=0.\end{cases}$$
    \item\label{Item::ODE::RmkInvStr::MoreReg} Suppose $\Se+\bar\Se$ is a subbundle. Then it is possible that $\Se+\bar\Se$ is more regular than $\Se$. For example, in $\R^3$ with coordinates $(x,y,z)$, $\Se$ is spanned by $\partial_x+i(1+|y|)\partial_y$. $\Se$ is merely a Lipschitz subbundle, but $\Se+\bar\Se$ can be spanned by $\partial_x,\partial_y$. So $\Se+\bar\Se$ is a real-analytic subbundle. Similarly if $\Se\cap\bar\Se$ is a subbundle than it can be more regular than $\Se$ as well. In this case $\Se\cap\bar\Se=\{0\}$ is a rank 0 subbundle which formally speaking is real-analytic as well.
\end{enumerate}
\end{remark}

In the definition of complex Frobenius structure, we have no regularity assumption on $\Se+\bar\Se$. But the following lemma shows that once $\Se+\bar\Se$ has constant rank then $\Se\cap\bar\Se$ and $\Se+\bar \Se$ are both subbundles with at least $\Co^\alpha$-regularities.
\begin{lem}\label{Lem::ODE::S+IsBundle}
	Let $\alpha>0$ and let $\Mf$ be a $\Co^{\alpha+1}$-manifold. Let $\Se\le\C T\Mf$ be a $\Co^\alpha$-complex subbundle such that $\Se+\bar \Se$ has constant rank at every point, then $\Se+\bar\Se$ and $\Se\cap\bar\Se$ are both $\Co^\alpha$-complex tangent subbundles.
\end{lem}
\begin{proof}Fix a point $p\in\Mf$. By Lemma \ref{Lem::ODE::GoodGen} $\Se$ has a $\Co^\alpha$-local basis $(X_1,\dots,X_{r+m})$ that has the form \eqref{Eqn::ODE::GoodGen::GoodGenFormula}. 

At the point $p$, by \eqref{Eqn::ODE::GoodGen::GoodGenFormula} we see that $X_1|_p,\dots,X_{r+m}|_p,\bar X_{r+1}|_p,\dots,\bar X_{r+m}|_p$ are linearly independent and span $\Se_p+\bar\Se_p=(\Se+\bar\Se)_p$. By assumption $\Se+\bar\Se$ has constant rank, so $\rank(\Se+\bar\Se)=r+2m$. Clearly by \eqref{Eqn::ODE::GoodGen::GoodGenFormula} $X_1,\dots,X_{r+m},\bar X_{r+1},\dots,\bar X_{r+m}$ are still linear independent near $p$. Therefore $\Se+\bar\Se$ is a $\Co^\alpha$-subbundle since $X_1,\dots,X_{r+m},\bar X_{r+1},\dots,\bar X_{r+m}$ form a $\Co^\alpha$-local basis.

Since $\Se^\bot+\bar \Se^\bot=(\Se\cap\bar\Se)^\bot\le \C T^*\Mf$ and $\rank(\Se_q\cap\bar\Se_q)=2\rank\Se_q-\rank(\Se_q+\bar\Se_q)=r$ is constant along $q\in\Mf$, we see that $\Se^\bot+\bar \Se^\bot$ is complex cotangent subbundle of rank $2m+q$. Since $T\Mf$ and $T^*\Mf$ are locally $\Co^\alpha$ isomorphic and that $\Se^\bot$ is also $\Co^\alpha$, repeating the above argument for $(\Se^\bot,\Se^\bot+\bar\Se^\bot)$, we see that $(\Se\cap\bar\Se)^\bot$ is also a $\Co^\alpha$ subbundle. Taking duality back we get $\Se\cap\bar\Se\in\Co^\alpha$.
\end{proof}

\begin{defn}
Let $\Mf$ be a $m$-dimensional complex manifold. A holomorphic vector field on $\Mf$ is a complex vector field such that locally it has form $\sum_{k=1}^m f^k(z)\Coorvec{z^k}$ where $(z^1,\dots,z^m):U\subseteq\Mf\to\C^m$ is a complex (holomorphic) coordinate chart for $\Mf$, and $f^k$ are holomorphic functions on $U$.

A rank $r$ holomorphic tangent subbundle over $\Mf$ is a complex tangent subbundle $\Se\le \C T\Mf$ such that for any $p\in\Mf$ there is a neighborhood $U\subseteq\Mf$ of $p$ and holomorphic vector fields $Z_1,\dots,Z_r$ on $U$ that that form a local basis for $\Se|_U$.
\end{defn}

Equivalently $\Se$ is a holomorphic section of the Grassmannian bundle $\Gr_\C(r,T^{\Oh}\Mf)\subset \Gr_\C(r,\C T\Mf)$, where $T^\Oh \Mf=\Span(\Coorvec{z^1},\dots,\Coorvec{z^m})$ is the holomorphic tangent bundle. In literature we also use $T^{1,0}\Mf$ for $T^\Oh\Mf$.

We have a similar result to Lemma  \ref{Lem::ODE::GoodGen}.
\begin{lem}\label{Lem::ODE::HoloGoodGen}
Let $\Mf$ is a $(r+q)$-dimensional complex manifold and $\Se\le T^\Oh \Mf$ be a rank $r$ holomorphic subbundle. Then for any point $p\in\Mf$ there is a holomorphic chart  $(z,w)=(z^1,\dots,z^r,w^1,\dots,w^q):\tilde U\subseteq \Mf\to\C^r\times\C^n$ such that $\Se_p\oplus\Span(\partial_w|_p)=T^\Oh_p \Mf$. Here $T^\Oh \Mf$ is the holomorphic tangent bundle.

Moreover for such chart $(z,w)$ there is a neighborhood $U\subseteq\tilde U$ of $p$, and a unique holomorphic local basis $Z=[Z_1,\dots,Z_r]^\top$ for $\Se$ on $U$ that has the form $Z=\partial_z+A\partial_w$ where $A\in\Oh(U;\C^{r\times n})$.
\end{lem}
\begin{proof}The existence of $(z,w)$ follows the same but simpler argument from the existence of the chart $(t,z,s)$ in Lemma \ref{Lem::ODE::GoodGen}.

Fix $(z,w)$, we can find a holomorphic local basis $W=[W_1,\dots,W_r]^\top$ for $\Se$ near $p$ that has expression $W=M'\partial_z+M''\partial_w$ for some holomorphic matrix functions $M',M''$ defined near $p$. Since $\Se_p\oplus\Span(\partial_w|_p)=T^\Oh_p \Mf$, $M'(p)\in\C^{r\times r}$ is invertible. By continuity $M'$ is invertible in some neighborhood $U\subseteq\Mf$ of $p$. And by cofactor representation of a holomorphic matrix we know $M'^{-1}$ is also a holomorphic matrix function on $U$.

Take $Z=M'^{-1}W=\partial_z+M'^{-1}M''\partial_w$ in $U$, so $A:=M'^{-1}M''$ is the desired holomorphic matrix function.

Since $\Se$ is a complex tangent subbundle of $\Mf$ and we $\Se_p\oplus\Span(\partial_{\bar z}|_p,\partial_{\re w}|_p,\partial_{\im w}|_p)=\C T_p \Mf$. By Lemma \ref{Lem::ODE::GoodGen} \ref{Item::ODE::GoodGen::Uniqueness}, $A$ is uniquely determined by the mixed real and complex coordinate system $(z,\re w,\im w)$, thus determined by $(z,w)$.
\end{proof}

Finally we give an auxiliary result to illustrate the regularity of the pullback subbundle. Recall Definition \ref{Defn::ODE::CpxPaSubbd} for the subbundles with mixed regularity.
\begin{lem}\label{Lem::ODE::PullBackReg}
Let $\kappa>2$, $\alpha\in(0,\kappa-1]$, $\beta,\gamma\in\R_\Eb^+$ such that $\beta>\Lip$ and $\beta,\gamma\le\kappa$. Let $0\le r\le  m$, $q\ge0$, let $\Mf$ be a $(m+q)$-dimensional $\Co^\kappa$ manifold and $\Se\le \C T\Mf$ be a $\Co^\alpha$-complex tangent subbundle of rank $r$.

Let $\Omega=\Omega'\times\Omega''\subseteq\R^m_x\times\R^q_s$ be an open subset. Suppose that $\Phi\in\Co^{\beta,\gamma}_\loc(\Omega',\Omega'';\Mf)$ is a topological parameterization of $\Mf$ such that  $\Phi_*\Coorvec{x^1},\dots,\Phi_*\Coorvec{x^m}\in\Co^\alpha_\loc(\Phi(\Omega);T\Mf)$ and $\Se|_{\Phi(\Omega)}\le\Span(\Phi_*\Coorvec{x^1},\dots,\Phi_*\Coorvec{x^m})$ (here $\Span$ is for the complex span). Then $\Phi^*\Se\le (\C T\Omega')\times\Omega''$ is a subbundle that has regularity $\Co^{\min(\alpha,\beta),\alpha\circ\gamma}$. 

{\normalfont In particular by Corollary \ref{Cor::Hold::CompOp}, when $\beta\ge\alpha$},
\begin{itemize}[nolistsep]
    \item $\Phi^*\Se\in\Co^{\alpha}$ if $\gamma>\Lip$ and $\gamma\ge\alpha$.
    \item $\Phi^*\Se\in\Co^{\alpha,\alpha-}$ if $\alpha\le1$ and $\gamma=1-$.
    \item $\Phi^*\Se\in\Co^{\alpha,\alpha\gamma}$ if $\max(\alpha,\gamma)<1$.
\end{itemize}
\end{lem}
In general, say $\beta=\gamma<\alpha+1$, it is possible that $\Phi_*\Coorvec{x^{m+1}},\dots,\Phi_*\Coorvec{x^n}$ are only $\Co^{\beta-1}(\supset\Co^\alpha)$.
\begin{proof}
Let $(x_0,s_0)\in \Omega$, we need to show that $\Phi^*\Se$ has a $\Co^{\alpha,\alpha\circ\gamma}$-local basis near $(x_0,s_0)$.

By assumption $\Se$ has $\Co^\alpha$-local basis $(X_1,\dots,X_r)$ near $\Phi(u_0)\in\Mf$. Thus $(\Phi^*X_1,\dots,\Phi^*X_r)$ is a local basis for $\Phi^*\Se$ near $u_0$.

Since $\Se|_{\Phi(\Omega)}\le\Span(\Phi_*\Coorvec{x^1},\dots,\Phi_*\Coorvec{x^m})$, we can write $X_j=\sum_{k=1}^mf_j^k\cdot\Phi_*\Coorvec{x^k}$ for $1\le j\le r$, for some functions $f_j^k$ defined near $\Phi(u_0)$. Since $f_j^k$ are uniquely determined by $X_1,\dots,X_r$ and $\Phi$, we see that $f_j^k\in\Co^{\min(\alpha,\beta)}=\Co^\alpha$.

Taking pullback by $\Phi$, we have $\Phi^*X_j=\sum_{k=1}^m(f_j^k\circ\Phi)\Coorvec{x^k}$, since $\Phi\in\Co^{\beta,\gamma}_{x,s}$, by Lemma \ref{Lem::Hold::CompofMixHold} \ref{Item::Hold::CompofMixHold::Comp} we get $f_j^k\circ\Phi\in\Co^{\alpha\circ\beta,\alpha\circ\gamma}_{x,s}$, so $(\Phi^*X_1,\dots,\Phi^*X_r)$ is a $\Co^{\alpha\circ\beta,\alpha\circ\gamma}$-local basis for $\Se$. By Corollary \ref{Cor::Hold::CompOp} \ref{Item::Hold::CompOp::Lip} $\alpha\circ\beta=\min(\alpha,\beta)$, we finish the proof.
\end{proof}

\subsection{The flow commuting of log-Lipschitz vector fields}\label{Section::FlowComm}

In this subsection we estimate the ODE flow regularity for log-Lipschitz vector fields, and show that if two log-Lipschitz vector fields commute then their ODE flows also commute.

We define the log-Lipschitz modulus of continuity ${\mu_{\log}}:(0,1]\to(0,\infty)$ and a function $\Mc_{\mu_{\log}}:(0,1]\to[0,\infty)$ (also see \cite[(3.4)]{BahouriCheminDanchin}) as
\begin{equation}\label{Eqn::ODE::Mulog}
         {\mu_{\log}}(r):=r\log\tfrac er=r(1-\log r),\quad \Mc_{\mu_{\log}}(r):=\int_r^1\frac{ds}{{\mu_{\log}}(s)}=\log\log \frac er,\quad0< r\le 1.
\end{equation}
Thus \eqref{Eqn::Intro::ClogNorm} can be rewritten as
\begin{equation}\label{Eqn::ODE::ClogNorm}
    \|f\|_{\Co^{\LogL}( U;\R^m)}=\|f\|_{C^0( U;\R^m)}+\sup_{x,y\in  U;0<|x-y|<1}|f(x)-f(y)|\cdot{\mu_{\log}}(|x-y|)^{-1}.
\end{equation}

Clearly $\lim\limits_{r\to0}\Mc_{\mu_{\log}}(r)=+\infty$, so ${\mu_{\log}}$ satisfies the Osgood condition, namely
\begin{equation*}
    \int_0^1\frac{dr}{{\mu_{\log}}(r)}=+\infty.
\end{equation*}

For a log-Lipschitz vector field $X\in\Co^\LogL_\loc( U;\R^n)$ and a point $p\in U$, By Osgood's uniqueness result \cite{Osgood} (see \cite[Theorem 3.2]{BahouriCheminDanchin}),  the autonomous ODE $\dot\gamma(t)=X(\gamma(t))$, $\gamma(0)=p$ has a unique solution. 
\begin{defn}\label{Defn::ODE::FlowMap}
Let $X$ be a log-Lipschitz vector field on a $C^{1,1}$-manifold $\Mf$. The \textbf{flow map} of $X$, denoted by $\exp_X:\Do_X\subseteq\R\times \Mf\to \Mf$, is defined in the way that for any $p\in \Mf$, $t \mapsto\exp_X(t,p)$ is the unique maximal solution $\gamma=\gamma_p$ to the ODE 
$$\left\{\begin{array}{l}\dot\gamma(t)=X(\gamma(t))\in T_{\gamma(t)}\Mf,\\\gamma(0)=p.\end{array}\right.$$
Here $\Do_X\cap(\R\times\{p\})\hookrightarrow\R$ is the maximal existence interval of $\gamma$, which is an open interval containing $0$.

We denote $e^{tX}(p):=\exp_X(t,p)$.
\end{defn}
\begin{remark}There is no ambiguity to use $e^{tX}$ since $(t,p)\in\Do_X$ if and only if $(1,p)\in\Do_{tX}$ and $\exp_X(t,p)=\exp_{tX}(1,p)$ for all $(t,p)\in\Do_X$.
\end{remark}

Suppose in addition that $X$ has compact support, by taking the zero extension outside $ U$ we can say $X\in\Co^{\LogL}_c(\R^n;\R^n)$. In this case $e^{tX}(p)\equiv p$ holds for all $t\in\R$ and $p\in U^c$, so by the maximal existence theorem (see \cite[Proposition 3.11]{BahouriCheminDanchin}) we see that $e^{tX}:\R^n\to\R^n$ is defined for all $t\in\R$. Later on we often assume the vector fields to have compact supports in $\R^n$ so that we do not need to worry about the domains of their ODE flows.

In order to give the regularity estimate of $e^{tX}$, we use the following Gronwall-Osgood inequality:
\begin{lem}[{\cite[Lemma 3.4 and Corollary 3.5]{BahouriCheminDanchin}}]\label{Lem::ODE::GroOsgInq}
Let $\mu:(0,1]\to\R_+$ be a Osgood modulus of continuity, that is, $\mu$ is a continuous increasing function satisfying $\lim_{r\to0^+}\mu(r)=0$ and $\int_0^1\frac{dr}{\mu(r)}=+\infty$. Let $\tau\in\R_+$, $\gamma\in L^1([0,\tau];\R_+)$ and $\rho:[0,\tau]\to[0,1]$  satisfy
\begin{equation}\label{Eqn::ODE::GroOsgAss}
    \begin{gathered}
    \rho(0)<1;\qquad\rho(t)\le\rho(0)+\int_0^t\gamma(s)\mu(\rho(s))ds,\quad 0\le s\le\tau;
    \\
    \int_0^\tau\gamma(t)dt\le\Mc_\mu(\rho(0)),\quad\text{where }\Mc_\mu(t):=\int_t^1\frac{ds}{\mu(s)}.
\end{gathered}
\end{equation}

Then
\begin{equation}\label{Eqn::ODE::GroOsgCon}
    \rho(t)\le\Mc_\mu^{-1}\Big(\Mc_\mu(\rho(0))-\int_0^t\gamma(s)ds\Big),\quad\forall 0\le t\le\tau.
\end{equation}
\end{lem}

For the special case where $\mu={\mu_{\log}}$ we have the following:
\begin{lem}\label{Lem::ODE::Gronwall}
Let $0<a<1$, $b>0$ and $0<\tau\le\frac1b\log\log\frac ea$. Let $\rho:[0,\tau]\to\R$ be an integrable function such that $\rho(0)=a$ and
\begin{equation}\label{Eqn::ODE::Gronwall::Assumption}
    \rho(t)\le a+b\int_0^t\rho(s)(1-
    \log\rho(s))ds,\quad\forall 0\le t\le \tau.
\end{equation}
Then 
\begin{equation*}
    \rho(t)\le e\cdot a^{e^{-bt}},\quad\forall 0\le t\le\tau.
\end{equation*}
\end{lem}
Note that from the assumption of $\tau$ we have $e\cdot a^{e^{-bt}}\le1$.
\begin{proof}
Take $\gamma(t):\equiv b$ in Lemma \ref{Lem::ODE::GroOsgInq}. Note that by \eqref{Eqn::ODE::Mulog} $\int_0^\tau\gamma(s)ds=b\tau\le\log\log\frac ea=\Mc_{\mu_{\log}}(\rho(0))$, so the assumptions in \eqref{Eqn::ODE::GroOsgAss} are all satisfied. Thus we have the conclusion \eqref{Eqn::ODE::GroOsgCon}, since
\begin{equation*}
    \rho(t)\le \Mc_{\mu_{\log}}^{-1}(\Mc_{\mu_{\log}}(a)-bt)=e^{1-e^{\Mc_{\mu_{\log}}(a)-bt}}= e^{1-e^{\log\log\frac ea}\cdot e^{-bt}}=e^{1-e^{-bt}}\cdot a^{e^{-bt}}\le e\cdot a^{e^{-bt}},\quad\forall0\le t\le\tau.\qedhere
\end{equation*}
\end{proof}

Lemma \ref{Lem::ODE::Gronwall} establishes the H\"older estimate for log-Lipschitz ODE flows.

\begin{lem}\label{Lem::ODE::LogLipODEReg}
Let $X$ be a bounded log-Lipschitz vector field on $\R^n$, and let $ U\subseteq\R^n$ be a bounded open set. We consider the map $\exp_X(t,x):=e^{tX}(x)$ for $t\in\R$ and $p\in\R^n$.
\begin{enumerate}[parsep=-0.3ex,label=(\roman*)]
    \item\label{Item::ODE::LogLipODEReg::BigLL} For any $\eps>0$ there is a small $\tau_\eps>0$ such that $\exp_X\in \Co^{1+\LogL,1-\eps}((-\tau_\eps,\tau_\eps),U;\R^n)$. \textnormal{(also see \cite[Theorem 3.7]{BahouriCheminDanchin})}
    \item\label{Item::ODE::LogLipODEReg::BigZyg1}If $X\in\Co^{1-}$, then $\exp_X\in \Co^2L^\infty((-T,T),U;\R^n)$. In particular for any $\eps>0$ there is a small $\tau_\eps>0$ such that $\exp_X\in \Co^{2,1-\eps}((-\tau_\eps,\tau_\eps),U;\R^n)$.
    \item\label{Item::ODE::LogLipODEReg::LittleLL} If $X$ is little log-Lipschitz, then   $\exp_X\in \Co^{1+\LogL,1-}((-T,T), U;\R^n)$ for all $T>0$.
    \item\label{Item::ODE::LogLipODEReg::LittleZyg1} In particular if $X\in\Co^{\logl}\cap\Co^1$, then   $\exp_X\in \Co^{2,1-}((-T,T), U;\R^n)$ for all $T>0$.
\end{enumerate}
\end{lem}
\begin{proof}
Since $X$ has compact support, we know $\exp_X$ is globally defined. Since we have $\frac \partial{\partial t}\exp_X=X\circ\exp_X$, we have $|e^{t_1X}(x)-e^{t_2X}(x)|\le|t_1-t_2|\|X\|_{C^0}$ for all $x\in\R^n$ and $t_1,t_2\in\R$, thus $\{\exp_X(\cdot,x):x\in U\}\subset\Co^\Lip(-T,T)$ is bounded. By Proposition \ref{Prop::Hold::QComp} \ref{Item::Hold::QComp::>1}, $\{X\circ\exp_X(\cdot,x):x\in U\}\in\Co^\LogL(-T,T)$,
therefore $\exp_X\in\Co^{1+\LogL}L^\infty((-T,T),U;\R^n)$ for all $T>0$.

If in addition $X\in\Co^1$, then by  Proposition \ref{Prop::Hold::QComp} \ref{Item::Hold::QComp::>1} again, since $1\circ(1+\LogL)=1$, we have $\exp_X\in\Co^2L^\infty((-T,T),U;\R^n)$ for all $T>0$.

Thus to prove \ref{Item::ODE::LogLipODEReg::BigLL} and \ref{Item::ODE::LogLipODEReg::BigZyg1} it suffices to show the following:
\begin{equation}\label{Eqn::ODE::LogLipODEReg::BigLLClaim}
    \forall\eps>0,\ \exists\tau_\eps,C_1>0\quad\text{such that }|e^{tX}(x_1)-e^{tX}(x_2)|\le C_1|x_1-x_2|^{1-\eps},\quad\text{for } |t|<\tau_\eps\text{ and } x_1,x_2\in U.
\end{equation}
And if in addition $X$ is little log-Lipschitz, then to prove \ref{Item::ODE::LogLipODEReg::LittleLL} and \ref{Item::ODE::LogLipODEReg::LittleZyg1} it suffices to prove the following:
\begin{equation}\label{Eqn::ODE::LogLipODEReg::LittleLLClaim}
    \forall\eps, T>0,\ \exists C_2>0\quad\text{such that }|e^{tX}(x_1)-e^{tX}(x_2)|\le C_2|x_1-x_2|^{1-\eps},\quad\text{for } |t|<T\text{ and } x_1,x_2\in U.
\end{equation}

We now prove \eqref{Eqn::ODE::LogLipODEReg::BigLLClaim} and \eqref{Eqn::ODE::LogLipODEReg::LittleLLClaim}.

\medskip
\noindent Proof of \eqref{Eqn::ODE::LogLipODEReg::BigLLClaim}: Since $ U$ is a bounded set, by enlarging the constant $C_1$ it suffices to prove \eqref{Eqn::ODE::LogLipODEReg::BigLLClaim} with $|x_1-x_2|<\frac1e$. And by replacing $X$ by $-X$ we only need to prove the case $t\ge0$.

Note that by \eqref{Eqn::ODE::ClogNorm} and \eqref{Eqn::ODE::Mulog}, we have for every $x_1,x_2\in\R^n$ and $t\ge0$,
\begin{equation}
    |e^{tX}x_1-e^{tX}x_2|\le\Big|x_1-x_2+\int_0^t\big(X(e^{sX}x_1)-X(e^{sX}x_2)\big)ds\Big|\le|x_1-x_2|+\|X\|_{\Co^{\LogL}}\int_0^{t}{\mu_{\log}}\big(|e^{sX}x_1-e^{sX}x_2|\big)ds.
\end{equation}
For $x_1,x_2$ that satisfy $|x_1-x_2|<\frac1e$, we take 
\begin{equation*}
    a:=|x_1-x_2|,\quad b:=\|X\|_{\Co^{\LogL}},\quad \rho(t):=|e^{tX}(x_1)-e^{tX}(x_2)|,\quad \tau:=b^{-1}\cdot\log 2\le\tfrac1b\log\log\tfrac ea.
\end{equation*}

Now \eqref{Eqn::ODE::Gronwall::Assumption} is satisfied, and by Lemma \ref{Lem::ODE::Gronwall} we get
$|e^{tX}(x_1)-e^{tX}(x_2)|\le e|x_1-x_2|^{e^{-bt}}$ for all $0\le t\le\tau$ and $|x_1-x_2|<\frac1e$.

Therefore for $\eps>0$, taking $\tau_\eps=\min(\tau,\frac1b\ln\frac1{1-\eps})$ we get 
\begin{equation*}
    |e^{tX}(x_1)-e^{tX}(x_2)|\le e|x_1-x_2|^{e^{-bt}}\le e|x_1-x_2|^{e^{-b\tau_\eps}}\le e|x_1-x_2|^{1-\eps},\quad\forall0\le t\le\tau_\eps\text{ and }|x_1-x_2|<\tfrac1e.
\end{equation*}

Since $U$ is a bounded domain, let $L=\delta^{-1}\operatorname{diam}U$, we then have $ |e^{tX}(x_1)-e^{tX}(x_2)|\le eL|x_1-x_2|^{1-\eps}$ for all $x_1,x_2\in U$. We obtain \eqref{Eqn::ODE::LogLipODEReg::BigLLClaim} and complete the proof of \ref{Item::ODE::LogLipODEReg::BigLL} and \ref{Item::ODE::LogLipODEReg::BigZyg1}.

\medskip
\noindent Proof of \eqref{Eqn::ODE::LogLipODEReg::LittleLLClaim}:
Fix a small $\eps>0$ and a large $T>0$ from the assumption, the key is to find a $\delta=\delta(\eps, T)>0$ such that $|e^{tX}(x_1)-e^{tX}(x_2)|\le e|x_1-x_2|^{1-\eps}$ for all $|x_1-x_2|<\delta$ and $-T<t<T$. Replacing $X$ by $-X$ it suffices to estimate the case $0\le t<T$.

Let $b:=\frac 1T\log\frac1{1-\eps}>0$. By assumption $X$ is little log-Lipschitz, so there is a $\delta_1>0$ such that
\begin{equation*}
    |X(x_1)-X(x_2)|\le (b/e)^\frac1{1-\eps}|x_1-x_2|(1-\log|x_1-x_2|),\quad\forall x_1,x_2\in\R^n\text{ satisfying }|x_1-x_2|<\delta_1.
\end{equation*}

Take $\delta=\min(e^{1-e^{bT}},\delta_1)$, we have $T\le\frac1b\log\log\frac e\delta$. Now for any $x_1,x_2$ such that $|x_1-x_2|<\delta$, we take 
\begin{equation*}
    a:=|x_1-x_2|,\quad \rho(t):=|e^{tX}(x_1)-e^{tX}(x_2)|,\quad \tau:=T\le\tfrac1b\log\log\tfrac ea.
\end{equation*}
So \eqref{Eqn::ODE::Gronwall::Assumption} is satisfied, and by Lemma \ref{Lem::ODE::Gronwall} we get
$|e^{tX}(x_1)-e^{tX}(x_2)|\le e|x_1-x_2|^{e^{-bt}}\le e|x_1-x_2|^{1-\eps}$ for all $0\le t<T$ and $|x_1-x_2|<\delta$. 

Since $U$ is a bounded domain, let $L=\delta^{-1}\operatorname{diam}U$, we then have $|e^{tX}(x_1)-e^{tX}(x_2)|\le eL |x_1-x_2|^{1-\eps}$ for all $x_1,x_2\in U$. This gives \eqref{Eqn::ODE::LogLipODEReg::LittleLLClaim} and finishes the proof of \ref{Item::ODE::LogLipODEReg::LittleLL} and \ref{Item::ODE::LogLipODEReg::LittleZyg1}.
\end{proof}
\begin{remark}
In fact in \ref{Item::ODE::LogLipODEReg::LittleLL} we have $\exp_X\in\Co^{1+\logl,1-}$, if we define $\Co^{1+\logl}$ to be the space of function whose derivatives are all $\Co^\logl$.


In Lemma \ref{Lem::ODE::LogLipODEReg} \ref{Item::ODE::LogLipODEReg::BigLL}, even though $\exp_X\in \Co^{1-\eps}$ near $t=0$ holds for arbitrary $\eps\in(0,1)$, it is possible that $\exp_X$ cannot be $\Co^{1-}$ near $t=0$. See Corollary \ref{Cor::Exmp::SharpRealFro} \ref{Item::Exmp::SharpRealFro::LogL}.
\end{remark}

For a log-Lipschitz vector field $X\in \Co^{\LogL}_c(\R^n;\R^n)$ with compact support, we can consider the Schwartz approximation $X^\nu:=\Su_\nu X$ for $\nu\in\Z_{\ge0}$, where $\Su_\nu$ is the Littlewood-Paley summation operator given in Definition \ref{Defn::Hold::DyadicResolution}. The following result shows the uniform convergence of the log-Lipschitz integral curves.

\begin{defn}\label{Defn::ODE::SubUnitCurve}
    Let $\Mf$ be a $C^1$-manifold and let $X_1,\dots,X_m$ be continuous vector fields on $\Mf$. We say a curve $\gamma:[0,T]\to\Mf$ is \textit{a subunit curve along $X_1,\dots,X_m$}, if $\gamma$ is Lipschitz, and $\dot\gamma(t)\overset{a.e}{=}\sum_{j=1}^mu^j(t)X_j(\gamma(t))$ holds for some $u^1,\dots,u^m\in L^\infty[0,T]$ such that $\sum_{j=1}^m|u^j(t)|^2\le 1$ for almost every $0\le t\le T$.
\end{defn}
Thus the Carnot-Carath\'eodory distance for $X=(X_1,\dots,X_m)$ (see \eqref{Eqn::Intro::CCDist}) can be restated as
\begin{equation}\label{Eqn::ODE::CCDistNew}
    \dist_X(x,y)=\inf\{T>0:\exists\text{ subunit curve }\gamma:[0,T]\to\Mf\text{ along }X\text{ such that }\gamma(0)=p,\gamma(T)=q\}.
\end{equation}
\begin{lem}\label{Lem::ODE::IntCurveApprox}Let $m\ge1$. Let $(\phi_\sigma)_{\sigma=0}^\infty$ be a fixed dyadic resolution that defines $\{\Su_\sigma\}_{\sigma=0}^\infty$. For any $0<\delta<1$ and $M>1$ there exist a $0<\tau_1<1$ and a $C_1>1$ that satisfy the following:

Let $p\in\R^n$. Let $X_1,\dots,X_m\in\Co^{\LogL}(\R^n;\R^n)$ be log-Lipschitz vector fields such that $\sup_{1\le j\le m}\|X_j\|_{\Co^\LogL}<M$. Let $u=(u^1,\dots,u^m)\in L^\infty([0,1];\R^m)$ satisfies $\|u\|_{L^\infty}^2=\essup_{t\in[0,\tau_1]}\sum_{i=1}^m|u^j(t)|^2\le1$. Let $X_j^\sigma:=\Su_\nu X_j$ for $1\le j\le m$ and $\sigma=1,2,3\dots$. Consider the ODEs
\begin{equation}\label{Eqn::ODE::IntCurveApprox::GammaDef}
    \dot\gamma^\sigma(t)\overset{\text{a.e.}}{=}\sum_{i=1}^mu^i(t)X_i^\sigma(\gamma^\sigma(t)),\quad \dot\gamma(t)\overset{\text{a.e.}}{=}\sum_{i=1}^mu^i(t)X_i(\gamma(t)),\quad\gamma^\sigma(0)=\gamma(0)=p,\quad\text{for }\sigma\ge1.
\end{equation}
Then $\gamma^\sigma,\gamma:[0,\tau_1]\to \R^n$ are uniquely defined Lipschitz curves, and satisfy
\begin{equation}\label{Eqn::ODE::ODEforSing::GammaEst}
    \sup_{0\le t\le \tau_1}|\gamma^\sigma(t)-\gamma(t)|\le C_12^{-(1-\delta)\sigma},\quad\forall \sigma\ge1.
\end{equation}

In particular $\lim\limits_{\sigma\to\infty}\gamma^\sigma(t)=\gamma(t)$ uniformly for $t\in[0,\tau_1]$. And as a special cases, $\lim\limits_{\sigma\to\infty}e^{tX_j^\sigma}(p)=e^{tX_j}(p)$ holds uniformly for $1\le j\le m$, $t\in[-\tau_1,\tau_1]$ and $p\in\R^n$.
\end{lem}


\begin{proof}
The existence and uniqueness of $\gamma$ and $\gamma^\sigma$ can be founded in, for example \cite[Theorem 3.2]{BahouriCheminDanchin}. Our focus is to prove \eqref{Eqn::ODE::ODEforSing::GammaEst}.

Firstly for $1\le j\le m$,
\begin{align*}
    \|X^\sigma_j-X_j\|_{L^\infty}&=\sup_{x\in\R^n}|2^{n\sigma}\phi_0(2^\sigma\cdot)\ast X_j(x)-X_j(x)|\le \sup_{x\in\R^n}\int_{\R^n}|2^{n\sigma}\phi_0(2^\sigma y)||X_j(x-y)-X_j(x)|dy
    \\
    &\le\|X_j\|_{\Co^\LogL}\int_{\R^n}2^n\sigma|y|\max(\log\tfrac ey,1)|\phi_0(2^\sigma y)|dy\le C_{\phi}\|X_j\|_{\Co^\LogL}\cdot\sigma 2^{-\sigma}.
\end{align*}
Here $C_\phi>0$ is a constant depending only on $\phi_0$. Therefore
\begin{align*}
    |\dot\gamma^\sigma(t)-\dot\gamma(t)|&\le\max_{1\le j\le m}|X_j^\sigma(\gamma^\sigma(t))-X_j(\gamma(t))|\le\sum_{1\le j\le m}\Big(\|X_j^\sigma-X_j\|_{L^\infty}+|X_j(\gamma^\sigma(t))-X_j(\gamma(t))|\Big)
    \\
    &\le \max_{1\le j\le m}\|X_j\|_{\Co^\LogL}\Big(C_\phi\cdot\sigma 2^{-\sigma}+\mu_{\log}(|\gamma^\sigma(t)-\gamma(t)|)\Big).
\end{align*}
Applying Lemma \ref{Lem::ODE::Gronwall} in the way similar to the proof of Lemma \ref{Lem::ODE::LogLipODEReg}, we see that there is a $C_1'=C_1'(\phi,M)>0$ such that
\begin{equation*}
    \text{When }\sigma 2^{-\sigma}<\tfrac1{C_1'},\quad\text{then }|\gamma^\sigma(t)-\gamma(t)|\le e\cdot (C_1'\sigma 2^{-\sigma})^{e^{-C_1't}}.
\end{equation*}
Thus we can take $\tau_1>0$ such that $e^{-C_1'\tau_1}=1-\delta/2$ and $C_1>0$ be such that
\begin{equation*}
    C_1>2^{\sigma_0}\tau_1\max_{1\le j\le m}\|X_j\|_{L^\infty}+eC_1'\max_{\sigma\ge\sigma_0}\sigma 2^{-\frac\delta2\sigma},\quad\text{where }\sigma_0\text{ is the largest number such that }\sigma_02^{-\sigma_0}\ge\tfrac1{C_1'}.
\end{equation*}

We see that $\tau_1,C_1$ depends only on $\phi,\delta,M$ but not on the choice of $\gamma$. Thus \eqref{Eqn::ODE::ODEforSing::GammaEst} holds.

Let $\sigma\to\infty$ we see that $\gamma^\sigma(t)\to\gamma(t)$ uniformly for $t\in[0,\tau_1]$ and all such $\gamma$. Take $\gamma(t):=e^{tX_j}(p)$ or $\gamma(t):=e^{-tX_j}(p)$ we see that $\gamma^\sigma(t)=e^{\pm tX_j^\sigma}(p)$ and therefore $e^{tX_j^\sigma}(p)\to e^{tX_j}(p)$ uniformly for $t\in[-\tau_1,\tau_1]$ and $p\in\R^n$.
\end{proof}

Lemma \ref{Lem::ODE::IntCurveApprox} extends to the estimate of $\gamma(t)$ precomposing a log-Lipschitz function $F$:

\begin{lem}\label{Lem::ODE::GronApprox}
Let $m,n\ge1$. Let $(\phi_\sigma)_{\sigma=0}^\infty$ be a fixed dyadic resolution that defines $\{\Su_\sigma\}_{\sigma=0}^\infty$. For any $0<\delta<1$ and $M>1$ there is a $\tau_2>0$ and $C_2>0$ that satisfy the following:

Let $p\in\R^n$, let $X_1,\dots,X_m\in\Co^\LogL(\R^n;\R^{1\times n})$ be log-Lipschitz vector fields as well as row-vector valued functions. Let $F\in\Co^{\LogL}(\R^n;\R^m)$ be a vector-valued function and let $A_1,\dots,A_m\in\Co^{\LogL-1}(\R^n;\R^{m\times m})$ be matrix-valued distributions. Suppose
\begin{itemize}[parsep=-0.3ex]
    \item $F(p)=0\in\R^m$.
    \item $X_j\cdot\nabla F=A_j\cdot F$ as $\R^m$-valued distributions for $1\le j\le m$.
    \item $\sup_{1\le j\le m}\|X_j\|_{\Co^\LogL}+\|A_j\|_{\Co^{\LogL-1}}<M$.
\end{itemize}

For $\sigma\ge1$, let $F^\sigma:=\Su_\sigma F$ and let $\gamma,\gamma^\sigma$ be defined in \eqref{Eqn::ODE::IntCurveApprox::GammaDef}. Then
\begin{equation}\label{Eqn::ODE::GronApprox}
    \sup_{0\le t\le \tau_2}|F^\sigma(\gamma^\sigma(t))|\le C_22^{-(1-\delta)\sigma},\quad\forall\sigma\ge1.
\end{equation}
In particular $F(\gamma(t))\equiv0$ for $0\le t\le\tau_2$.
\end{lem}
\begin{proof}
    By assumption we have $(X_j\cdot\nabla F)^\top-F^\top\cdot A_j^\top=0\in\R^{m\times 1}$ where $X_j,F\in\Co^\alpha$ and $\nabla F,A_j\in\Co^{\alpha-1}$ for all $\frac12<\alpha<1$. Applying Theorem \ref{Thm::Hold::ApproxThm} to their components we see that there is a $C_1'=C_1'(m,n,M,\delta)>1$ such that
    \begin{equation*}
        \|X_j^\sigma\cdot\nabla F^\sigma-A_j^\sigma\cdot F^\sigma\|_{L^\infty(\R^n;\R^m)}\le C_1'2^{-(1-\delta/2)\sigma},\quad\forall\sigma\ge1,\quad1\le j\le m.\quad\text{Here }X_j^\sigma:=\Su_\sigma X_j,\quad A_j^\sigma:=\Su_\sigma A_j.
    \end{equation*}
    
    Since $A_j\in\Co^{\LogL-1}$, we have $\|A_j^\sigma\|_{L^\infty}\le \sigma\|A_j\|_{\Co^{\LogL-1}}$ by Definition \ref{Defn::Hold::LogL-1}. Thus by \eqref{Eqn::ODE::IntCurveApprox::GammaDef} we have: for every $\sigma\ge0$ and every $t$ in the domain of $\gamma$,
    \begin{equation}\label{Eqn::ODE::GronApprox::GronwallPrep}
        \begin{aligned}
        \bigg|\frac d{dt}F^\sigma\circ\gamma^\sigma(t)\bigg|&=\bigg|\sum_{j=1}^mu^j(t)(X_j^\sigma\cdot\nabla F^\sigma)(\gamma^\sigma(t))\bigg|\le\sum_{j=1}^m|u^j(t)||(A_j^\sigma\cdot F^\sigma)(\gamma^\sigma(t))|+C_1'2^{-(1-\delta/2)\sigma}
        \\
        &\le\sigma|F^\sigma(\gamma^\sigma(t))|\cdot\sup_{1\le j\le m}\|A_j\|_{\Co^{\LogL-1}}+C_1'2^{-(1-\delta/2)\sigma}\le C_2'\sigma|F^\sigma\circ\gamma^\sigma(t)|+C_2'2^{-(1-\delta/2)\sigma}.
    \end{aligned}
    \end{equation}
    Here $C_2':=\max_{1\le j\le m}\|A_j\|_{\Co^{\LogL-1}}+C_1'>1$ depends only on $m,n,\phi,\delta,M$.
    
    Therefore by Gronwall's inequality (or taking $\mu(t):=t$ in Lemma \ref{Lem::ODE::GroOsgInq}) we have 
    \begin{equation}\label{Eqn::ODE::GronApprox::Pf1}
        |F^\sigma\circ\gamma^\sigma(t)|\le|F^\sigma\circ\gamma^\sigma(0)|e^{C_2'\sigma t}+ \sigma^{-2}2^{-(1-\delta/2)\sigma}(e^{C_2'\sigma t}-1)\le |F^\sigma\circ\gamma^\sigma(0)|2^{\frac{C_2't}{\log2}\sigma }+2^{-(1-\frac\delta2-\frac{C_2't}{\log2})\sigma}.
    \end{equation}
    
    On the other hand we have $|\gamma^\sigma(0)|=|\gamma^\sigma(0)-\gamma(0)|\le C_12^{-(1-\delta/4)\sigma}$ where $C_1=C_1(\delta/4,M)>1$ is the constant in Lemma \ref{Lem::ODE::IntCurveApprox} depend only on $m,n,\phi,\delta,M$. Thus 
    \begin{align*}
        &|F^\sigma\circ\gamma^\sigma(0)|=|F^\sigma(\gamma^\sigma(0))-F(\gamma(0))|\le|F(\gamma^\sigma(0))-F(\gamma(0))|+|F^\sigma(\gamma^\sigma(0))-F(\gamma^\sigma(0))|
        \\
        \le&\|F\|_{\Co^{1-\frac\delta4}}|\gamma^\sigma(0)-\gamma(0)|^{1-\frac\delta4}+\|F^\sigma-F\|_{L^\infty}\le\|F\|_{\Co^{1-\frac\delta4}}C_1^{1-\frac\delta4}2^{-(1-\frac\delta4)^2\sigma}+\sum_{j=\sigma+1}^\infty\|\De_jF\|_{L^\infty}
        \\
        \lesssim&_{\delta}\|F\|_{\Co^{1-\delta/4}}2^{-(1-\frac\delta2)\sigma}+\|F\|_{\Co^{1-\delta/2}} 2^{-(1-\frac\delta2)\sigma}\lesssim_{\delta,M}2^{-(1-\frac\delta2)\sigma}. 
    \end{align*}
    We see that there is a $C'_3=C'_3(m,n,\phi,\delta,M)>0$ such that $|F^\sigma\circ\gamma^\sigma(0)|\le C_3'2^{-(1-\frac\delta2)\sigma}$. Plugging this to \eqref{Eqn::ODE::GronApprox::Pf1} we get
    \begin{equation*}
        |F^\sigma\circ\gamma^\sigma(t)|\le C_3'2^{-(1-\frac\delta2-\frac{C_2't}{\log2})\sigma}+2^{-(1-\frac\delta2-\frac{C_2't}{\log2})\sigma}.
    \end{equation*}
    
    Taking $\tau_2:=\min(\tau_1,\frac\delta2\cdot\frac{\log2}{C_2'})$ and $C_2:=C_3'+1$, where $\tau_1$ is the constant in  Lemma \ref{Lem::ODE::IntCurveApprox}, we obtain \eqref{Eqn::ODE::GronApprox}.
    
    Finally both convergence $F^\sigma\to F$ and $\gamma^\sigma\to\gamma$ are uniform, we conclude that $F(\gamma(t))=\lim_{\sigma\to\infty}F^\sigma(\gamma^\sigma(t))=0$ uniformly in $t\in[0,\tau_2]$.
\end{proof}

The key to prove the Frobenius theorem is to show that if two vector fields commute, then their ODE flows also commute. In order to prove the singular Frobenius theorem, Theorem \ref{MainThm::SingFro},  we need a stronger version to this property.
\begin{prop}[Strong flow commuting for log-Lipschitz vector fields]\label{Prop::ODE::StrFlowComm}
Let $\Mf$ be a $C^{1,1}$ manifold and let $S\Subset\Mf$ be a precompact set. Let $1\le r\le m$ and let $X_1,\dots,X_m$ be log-Lipschitz real vector fields on $\Mf$. Suppose we have
\begin{itemize}[parsep=-0.3ex]
    \item $X_1(p),\dots,X_r(p)\in T_p\Mf$ are linearly independent, and $X_{r+1}(p)=\dots=X_m(p)=0$ for all $p\in S$.
    
    \item There are $c_{ij}^k\in\Co^{\LogL-1}_\loc(\Mf)$ for $1\le i\le r$, $1\le j\le m$ and $r+1\le k\le m$, such that as distributions,
    \begin{equation*}
        [X_i,X_j]=\sum_{k=r+1}^mc_{ij}^kX_k,\quad1\le i\le r,\quad 1\le j\le m.
    \end{equation*}
\end{itemize}
Then there is a $\tau_0>0$ such that 
\begin{enumerate}[parsep=-0.3ex,label=(\roman*)]
    \item\label{Item::ODE::StrFlowComm::Xj=0} $X_j(e^{t^1X_1}\dots e^{t^rX_r}(p))\equiv0$ for all $r+1\le j\le m$ and $t^1,\dots,t^r\in(-\tau_0,\tau_0)$.
    \item\label{Item::ODE::StrFlowComm::Comm} the following equality holds with both sides are defined on $\Mf$:
\begin{equation}\label{Eqn::ODE::StrFlowComm::Result}
    e^{t^1X_1}\dots e^{t^rX_r}(p)=e^{t^{\sigma_1}X_{\sigma_1}}\dots e^{t^{\sigma_r}X_{\sigma_r}}(p),\quad\text{for }t^1,\dots,t^r\in(-\tau_0,\tau_0)\quad p\in S,
\end{equation}
for every permutation $\sigma=(\sigma_1,\dots,\sigma_r)$ of $\{1,\dots,r\}$. In particular
\begin{equation}\label{Eqn::ODE::StrFlowComm::ddtFlow}
    \Coorvec{t^j}e^{t^1X_1}\dots e^{t^rX_r}(p)=X_j\circ e^{t^1X_1}\dots e^{t^rX_r}(p),\quad\text{for }t^1,\dots,t^r\in(-\tau_0,\tau_0)\quad p\in S.
\end{equation}

\end{enumerate}

\end{prop}
\begin{remark}
Here unlike \eqref{Eqn::MainThm::SingFro::Inv} we do not need $[X_i,X_j]=\sum_kc_{ij}^kX_k$ to hold for $r+1\le i,j\le m$.

When $X_1,\dots,X_r$ pairwise commute, \eqref{Eqn::ODE::StrFlowComm::Result} holds automatically, as we can take $c_{ij}^k\equiv0$.
\end{remark}
\begin{proof}
By passing to a local, with possibly a partition of unity, we can assume $\Mf=\Omega\subseteq\R^n$ is an open subset.
We can take a precompact subset $\Omega'\Subset\Omega$ containing $\bar S$ and a cut-off $\chi\in C_c^\infty(\Omega)$ such that $\chi|_{\Omega'}\equiv1$. Let 
\begin{equation}\label{Eqn::ODE::StrFlowComm::DefT0}
    T_0:=(r\sup_{1\le i\le m}\|X_i\|_{L^\infty(\Omega')})^{-1}\cdot\dist(S,\partial\Omega'),
\end{equation}
we see that the following equality holds with both sides are defined. 
\begin{equation}\label{Eqn::ODE::StrFlowComm::Tmp1}
    e^{t^{i_1}X_{i_1}}\dots e^{t^{i_r}X_{i_r}}(p)=e^{t^{i_1}(\chi X_{i_1})}\dots e^{t^{i_r}(\chi X_{i_r})}(p)\in\Omega',\quad\forall i_1,\dots,i_r\in\{1,\dots,r\},\quad |t^{i_1}|,\dots,|t^{i_r}|\le T_0,\quad p\in S.
\end{equation}

First we note that it suffices find a $\tau_0\in(0,T_0)$ that satisfies the following: let $p\in S$ be a fixed point,
\begin{itemize}[parsep=-0.3ex]
    \item Let $\gamma\in C^{0,1}([0,r\tau_0];\Omega')$ that satisfies $\gamma(0)=p$ and $\dot\gamma(t)=\sum_{j=1}^ru^j(t)X_j(\gamma(t))$ for some $u=(u^1,\dots,u^r)\in L^\infty[0,r\tau_0]$ with $\|u\|_{L^\infty}\le1$, we have 
    \begin{equation}\label{Eqn::ODE::StrFlowComm::Claim}
        e^{t_iX_i}e^{t_jX_j}(\gamma(t))=e^{t_jX_j}e^{t_iX_i}(\gamma(t)),\quad\forall 0\le t\le r\tau_0,\quad t_i,t_j\in[-\tau_0,\tau_0],\quad1\le i,j\le r.
    \end{equation}
\end{itemize}

Once the \eqref{Eqn::ODE::StrFlowComm::Claim} is proved, by taking $\gamma$ such that $(u^1(t),\dots,u^r(t))\in\{e_1,\dots,e_n\}$ (the unit coordinate directions) for all $t\in[0,k\tau_0]$, we see that for all $2\le k\le r$, $$e^{t^{i_1}X_{i_1}}e^{t^{i_2}X_{i_2}}e^{t^{i_3}X_{i_3}}\dots e^{t^{i_k}X_{i_k}}(p)=e^{t^{i_2}X_{i_2}}e^{t^{i_1}X_{i_1}}e^{t^{i_3}X_{i_3}}\dots e^{t^{i_k}X_{i_k}}(p),\ 1\le i_1\dots,i_k\le r,\ t^{i_1},\dots,t^{i_k}\in(-\tau_0,\tau_0).$$
By decomposing a permutation into transpositions, \eqref{Eqn::ODE::StrFlowComm::Result} then follows.

\medskip To prove \eqref{Eqn::ODE::StrFlowComm::Claim}, let $\gamma:[0,T_0]\to\Omega'$ be a such fixed curve. We take approximations $X^\sigma_j:=\Su_\sigma(\chi X_j)$, $c_{ij}^{k\sigma}:=\Su_\sigma(\chi c_{ij}^k)$ for $\sigma\ge1$, $1\le i\le r$, $1\le j\le m$ and $r+1\le k\le m$. And we define $\gamma^\sigma:[0,T_0]\to\Omega'$ be the unique curve such that $\gamma^\sigma(0)=p$ and $\dot\gamma^\sigma(t)=\sum_{j=1}^ru^j(t)X_j^\sigma(\gamma^\sigma(t))$. 

By Lemma \ref{Lem::ODE::IntCurveApprox} there is a $0<\tau_1<1$ that depends on $X_1,\dots,X_r,\chi, T_0$ but not on $\gamma$, such that $e^{uX_i^\sigma}e^{vX_j^\sigma}(\gamma^\sigma(t))\xrightarrow{\sigma\to\infty}e^{uX_i}e^{vX_j}(\gamma(t))$ (since $\chi X_j\equiv X_j$ in $\Omega'$) uniformly for all $1\le i,j\le r$, $t\in[0,\tau_1]$ and $u,v\in[-\tau_1,\tau_1]$. On the other hand, using the commutator formula from \cite[Lemma 4.1]{RampazzoSussmanCommutators}, for such $u,v,t$,
\begin{equation}\label{Eqn::ODE::StrFlowComm:Pf::CommFormula}
\begin{aligned}
    &e^{-uX_i^\sigma}e^{-vX_j^\sigma}e^{uX_i^\sigma}e^{vX_j^\sigma}(\gamma^\sigma(t))-\gamma^\sigma(t)
    \\
    &\quad=\int_0^u\int_0^v\left([X_i^\sigma,X_j^\sigma]\cdot\nabla(e^{-uX_i^\sigma}e^{-v'X_j^\sigma}e^{-u' X_i^\sigma})\right)\circ\left(e^{(u-u')X_i^\sigma}e^{v' X_j^\sigma}(\gamma^\sigma(t))\right)du' dv'.
\end{aligned}
\end{equation}
Since $\tau_1\le T_0$ is small, by \eqref{Eqn::ODE::StrFlowComm::DefT0} we see that $e^{(u-u')X_i^\sigma}e^{v' X_j^\sigma}(\gamma^\sigma(t))\in\Omega'$ for $u,u',v'\in[-\tau_1,\tau_1]$ and $t\in[0,\tau_1]$. Applying Proposition \ref{Prop::Hold::InvVFApt} we see that
\begin{equation*}
   |[X_i^\sigma,X_j^\sigma](q)|\le \sum_{k=r+1}^m|c_{ij}^{k\sigma}(q)||X_k^\sigma(q)|+C_1'2^{-\frac23\sigma},\quad\forall\sigma\ge1,\quad1\le i,j\le m,\quad q\in\Omega'.
\end{equation*}
Here $C_1'=C_1'(X,\chi,(c_{ij}^k))>1$ does not depend on $\sigma$ and $\gamma$. Thus,
\begin{align}
    |e^{-uX_i^\sigma}e^{-vX_j^\sigma}e^{uX_i^\sigma}e^{vX_j^\sigma}(\gamma^\sigma(t))-\gamma^\sigma(t)|
    \le&\sum_{k=r+1}^m\|c_{ij}^{k\sigma}\|_{C^0(\Omega')}\bigg|\int_0^u\int_0^vX_k^\sigma\circ e^{(u-u')X_i^\sigma} e^{v' X_j^\sigma}(\gamma^\sigma(t))du' dv'\bigg|
    \label{Eqn::ODE::StrFlowComm:Pf::Est1}
    \\&+C_1'2^{-\frac23\sigma}\bigg|\int_0^u\int_0^v\bigg\|\nabla\left(e^{-uX_i^\sigma}e^{-v'X_j^\sigma}e^{-u' X_i^\sigma}\right)\bigg\|_{C^0(\Omega';\R^{n\times n})}du' dv'\bigg|.
    \label{Eqn::ODE::StrFlowComm:Pf::Est2}
\end{align}
For \eqref{Eqn::ODE::StrFlowComm:Pf::Est2}, by Gronwall's inequality $\|\nabla(e^{-uX_i^\sigma}e^{-v'X_j^\sigma}e^{-u' X_i^\sigma})\|_{C^0}\le e^{|u|\|\nabla X_i^\sigma\|_{C^0}+|v'|\|\nabla X_j^\sigma\|_{C^0}+|u'|\|\nabla X_i^\sigma\|_{C^0}}$. Since $\|\nabla X_l^\sigma\|_{C^0}\lesssim\sigma\|X_l\|_{\Co^{\LogL}}$ for all $1\le l\le r$, we can find a $C_2'>0$ that does not depend on $\sigma$ and $\gamma$, such that
\begin{align*}
    C_1'2^{-\frac23\sigma}\Big|\int_0^u\int_0^v\Big\|\nabla\left(e^{-uX_i^\sigma}e^{-v'X_j^\sigma}e^{-u' X_i^\sigma})\right)\Big\|_{C^0(\Omega')}du' dv'\Big|\le C_1'2^{-\frac23\sigma}e^{(2|u|+|v|)C_2'\sigma},\quad\forall \sigma\ge1,\quad u,v\in[-\tau_1,\tau_1].
\end{align*}
Taking $\tau_0':=(9C_2)^{-1}\log2$, we see that when $|u|,|v|\le\tau_0'$, \eqref{Eqn::ODE::StrFlowComm:Pf::Est2} is bounded by $C_1'2^{-\sigma/3}$, which goes to $0$ as $\sigma\to\infty$.

For \eqref{Eqn::ODE::StrFlowComm:Pf::Est1}, write $X_i=\sum_{l=1}^na_i^l\Coorvec{x^l}=[a_i^1,\dots,a_i^n]$ for $1\le i\le m$, we define $F:\R^n\to\R^{(m-r) n}$ to be the coefficients of $(\chi X_{r+1},\dots,\chi X_m)$, i.e. $F=(\chi a_j^l)_{r+1\le j\le m;1\le l\le n}$. The assumption $[X_i,X_j]=\sum_kc_{ij}^kX_k$ implies $[\chi X_i,\chi X_j]=\sum_{k=r+1}^m\chi c_{ij}^k\cdot\chi X_k+(X_i\chi)\cdot\chi X_j-(X_j\chi)\cdot\chi X_i$ for $1\le i\le r$, $r+1\le j\le m$. Therefore
\begin{equation*}
    \chi X_i\cdot\nabla(\chi a_j^l)=\sum_{k=r+1}^m\chi c_{ij}^k\cdot\chi a_k^l+(X_i\chi)\cdot\chi a_j^l-\sum_{q=1}^n(\chi a_i^l\partial_q\chi)\cdot \chi a_j^q.\quad1\le i\le r<j\le m,\quad1\le l\le n.
\end{equation*}

Therefore there is a $A\in\Co^{\LogL-1}_c(\Omega;\R^{(m-r)n\times (m-r)n})$ such that $\chi X_i\cdot\nabla F=A\cdot F$ as $\R^{(m-r)n}$-valued distributions in $\R^n$.

Now by the assumption $X_{r+1}(p),\dots,X_m(p)=0$ we have $F(p)=0$. Therefore applying Lemma \ref{Lem::ODE::GronApprox}, we get \ref{Item::ODE::StrFlowComm::Xj=0} and we can find a $0<\tau_2<\tau_0'$ and $C_3'>1$ that does not depend on $\gamma$ and $\sigma$, such that
\begin{equation*}
    |X_k^\sigma\circ e^{(u-u')X_i^\sigma} e^{v' X_j^\sigma}(\gamma^\sigma(t))|\le C_3'2^{-\frac12\sigma},\quad\forall1\le i,j\le r<k\le m,\quad u,u',v'\in[-\tau_2,\tau_2],\quad t\in[0,\tau_2],\quad \sigma\ge1.
\end{equation*}
Plugging this into \eqref{Eqn::ODE::StrFlowComm:Pf::Est1} and since $\|c_{ij}^{k\sigma}\|_{C^0}\lesssim\sigma$, we see that \eqref{Eqn::ODE::StrFlowComm:Pf::Est1} goes to $0$ as $\sigma\to\infty$. 

We now conclude that $e^{-uX_i^\sigma}e^{-vX_j^\sigma}e^{uX_i^\sigma}e^{vX_j^\sigma}(\gamma^\sigma(t))-\gamma^\sigma(t)\to0$ uniformly for $u,v\in[-\tau_2,\tau_2]$, $t\in[0,\tau_2]$ and all subunit $\gamma$ along $X_1,\dots,X_r$ such that $\gamma(0)=p$. Take $\tau_0:=\tau_2/r$ we get \eqref{Eqn::ODE::StrFlowComm::Claim}.

Finally, for $1\le j\le r$,
\begin{align*}
    &\textstyle\Coorvec{t^j}e^{t^1X_1}\dots e^{t^rX_r}(p)=\Coorvec{t^j}e^{t^jX_j}e^{t^1X_1}\dots e^{t^{j-1}X_{j-1}}e^{t^{j+1}X_{j+1}}\dots e^{t^rX_r}(p)
    \\
    =&X_j\circ e^{t^jX_j}e^{t^1X_1}\dots e^{t^{j-1}X_{j-1}}e^{t^{j+1}X_{j+1}}\dots e^{t^rX_r}(p)=X_j\circ e^{t^1X_1}\dots e^{t^rX_r}(p).
\end{align*}
This completes the proof of \ref{Item::ODE::StrFlowComm::Comm}.
\end{proof}

Now we see that for two log-Lipschitz vector fields $X$ and $Y$, if $[X,Y]=0$ in the sense of distributions, then $e^{tX}\circ e^{sY}=e^{sY}\circ e^{tX}$ holds locally for small $t$ and $s$. This is not true for large time in general.

\begin{example}
Consider the vector fields $X(x,y)=(e^{-x}\cos y,-e^{-x}\sin y)$ and $Y(x,y)=(e^{-x}\sin y,e^{-x}\cos y)$ defined in $\R^2$. By direct computation, $[X,Y]=0$ so $X$ and $Y$ commute.

One can see that we have the following ODE solutions, for $t,s\in\R$:
\begin{align*}
    e^{tX}(\ln\sqrt2,\tfrac54\pi)=(\tfrac12\ln(t^2+2t+2),\tfrac32\pi+\arctan(t-1)),&\quad e^{2X}(\ln\sqrt2,\tfrac54\pi)=(\ln\sqrt2,\tfrac74\pi);
    \\
    e^{sY}(\ln\sqrt2,\tfrac54\pi)=(\tfrac12\ln(s^2+2s+2),\pi-\arctan(s-1)),&\quad e^{2Y}(\ln\sqrt2,\tfrac54\pi)=(\ln\sqrt2,\tfrac34\pi);
    \\
    e^{sY}(\ln\sqrt2,\tfrac74\pi)=(\tfrac12\ln(s^2+2s+2),2\pi+\arctan(s-1)),&\quad e^{2Y}(\ln\sqrt2,\tfrac74\pi)=(\ln\sqrt2,\tfrac94\pi);
    \\
    e^{tX}(\ln\sqrt2,\tfrac34\pi)=(\tfrac12\ln(t^2+2t+2),\tfrac12\pi-\arctan(t-1)),&\quad e^{2X}(\ln\sqrt2,\tfrac34\pi)=(\ln\sqrt2,\tfrac14\pi).
\end{align*}
Therefore $e^{2X}e^{2Y}(\ln\sqrt2,\tfrac54\pi)=(\ln\sqrt2,\tfrac94\pi)$ and $e^{2Y}e^{2X}(\ln\sqrt2,\tfrac34\pi)=(\ln\sqrt2,\tfrac14\pi)$ are not equal.

A related fact is that for the holomorphic vector field $Z(z)=e^{-z}\Coorvec z$ defined on $\C^1$, we cannot talk about its holomoprhic ODE flow globally.
\end{example}

\subsection{The canonical coordinate system along vector fields}\label{Section::CanonicalCoordinates}

In the quantitative Frobenius theorem the canonical coordinates give an important initial normalization of the given collection of vector fields. Roughly speaking, if we have a collection of vector fields $X_1,\dots,X_n$ on a manifold that span the tangent space at every point, then under the pullback of the canonical coordinates the size of $X_1,\dots,X_n$ depends only on their intrinsic structures and lose at most 1 (intrinsic) derivative.

Recall the flow map defined in Definition \ref{Defn::ODE::FlowMap}.

\begin{defn}\label{Defn::ODE::MultiFlow}
    Let $\Mf$ be a $C^{1,1}$-manifold and let $X=(X_1,\dots,X_m)$ be a collection of log-Lipschitz vector fields. The multi-flow $\exp_X$ is a map $\exp_X:\Do_X\subseteq\R^m\times\Mf\to\Mf$ given by 
    \begin{equation*}
        \Do_X:=\{(t^1,\dots,t^m,p)\in\R^m\times\Mf:(1,p)\in\Do_{t^1X_1+\dots+t^mX_m}\},\quad\exp_X(t,p):=\exp_{t\cdot X}(1,p)=e^{t\cdot X}(p).
    \end{equation*}
\end{defn}

By definition if $(t,p)\in\Do_X$ then $(rt,p)\in\Do_X$ holds for $0\le r\le 1$, thus $\Do_X\cap(\R^m\times\{p\})\hookrightarrow\R^m$ is star-shaped convex at $0\in\R^m$. 

For a direct link to Proposition \ref{Prop::ODE::StrFlowComm}, we have
\begin{lem}\label{Lem::ODE::CommMultFlow}
Let $X_1,\dots,X_r$ be log-Lipschitz vector fields on $\Mf$ which are pairwise commutative. Let $U\Subset\Mf$ be a precompact open subset. There is a $\tau_0>0$ such that $e^{t^1X_1}\dots e^{t^rX_r}(p)=\exp_X(t,p)$ with both sides are defined for all $t\in B^m(0,\tau_0)$ and $p\in S$. In particular $\Coorvec{t^j}e^{t\cdot X}(p)=X_j(e^{t\cdot X}(p))$ holds for all $p\in U$, $t\in B^r(0,\tau_0)$ and $1\le j\le r$.
\end{lem}
\begin{proof}
Let $\tau_0$ be the same number from the result of Proposition \ref{Prop::ODE::StrFlowComm}. By \eqref{Eqn::ODE::StrFlowComm::Result} we have for $p\in U$, $t\in B^r(0,\tau_0)$ and $1\le j\le r$,
\begin{align*}
    &\textstyle\Coorvec{t^j}e^{t^1X_1}\dots e^{t^rX_r}(p)=\Coorvec{t^j}e^{t^jX_j}e^{t^1X_1}\dots e^{t^{j-1}X_{j-1}}e^{t^{j+1}X_{j+1}}\dots e^{t^rX_r}(p)
    \\
    =&X_j\circ e^{t^jX_j}e^{t^1X_1}\dots e^{t^{j-1}X_{j-1}}e^{t^{j+1}X_{j+1}}\dots e^{t^rX_r}(p)=X_j\circ e^{t^1X_1}\dots e^{t^rX_r}(p).
\end{align*}
Thus for a fixed $p\in U$ and $t\in B^r(0,\tau_0)$, by chain rule, 
\begin{equation*}
    \frac d{ds}e^{st^1X_1}\dots e^{st^rX_r}(p)=(t^1X_1+\dots+t^rX_r)\circ(e^{st^1X_1}\dots e^{st^rX_r})(p),\quad s\in[0,1].
\end{equation*}
In other words $s\mapsto e^{st^1X_1}\dots e^{st^rX_r}(p)$ is the solution to the log-Lipschitz ODE $\dot\gamma(s)=(t^1X_1+\dots+t^rX_r)(s)$, $\gamma(0)=p$ for $s\in[0,1]$. On the other hand, by definition $\gamma(s)=e^{s(t\cdot X)}(p)$ whenever the ODE is solvable in the domain, thus taking $s=1$ we get $e^{t^1X_1}\dots e^{t^rX_r}(p)=e^{t\cdot X}(p)$ and in particular $\Coorvec{t^j}e^{t\cdot X}(p)=X_j(e^{t\cdot X}(p))$ holds for all $t\in B^r(0,\tau_0)$ and $1\le j\le r$.
\end{proof}

\begin{defn}\label{Defn::ODE::CanCoor}
    Following Definition \ref{Defn::ODE::MultiFlow}, suppose $m=\dim\Mf$. Let $p\in\Mf$.  The \textit{canonical coordinates} associated with $X=(X_1,\dots,X_m)$ at $p$ is the map $\Phi_p:\Omega\subseteq\R^m\to\Mf$, $\Phi_p(t):=e^{t\cdot X}(p)$, where $\Omega$ is an open set containing $0\in\R^m$ such that $\Omega\times\{p\}\subseteq\Do_X$ and that $\Phi_p:\Omega\to\Mf$ is homeomorphism onto its image.
\end{defn}
Note that $\Phi_p$ is a canonical coordinates then necessarily we know $X_1(p),\dots,X_m(p)$ form a (real) basis of $T_p\Mf$.

\begin{remark}
For most case in the thesis, a parameterization is a map of the form $\Phi:\Omega\subseteq\R^m\to\Mf$, and a coordinate system is a map of the form $F:U\subseteq\Mf\to\R^m$. In this convention we should call such $\Phi_p$ ``the canonical parameterization''. The name ``canonical coordinates'' follows from \cite{CoordAdapted}. 
\end{remark}

For convenience we use the following terminology for the topological assumption on the canonical coordinates mentioned in Theorem \ref{MainThm::QuantFro}:
\begin{defn}\label{Defn::ODE::AlmostInjRad}
    Let $\Mf$ be a  $C^{1,1}$-manifold. Let $X_1,\dots,X_n$ be log-Lipschitz vector fields on $\Mf$. 
    
    Let $p\in\Mf$, the \textit{almost injective radius} associated with $X=(X_1,\dots,X_n)$ at $p$ is the number $R\in[0,\infty]$ which is maximal among the following conditions:
    \begin{itemize}[parsep=-0.3ex]
        \item $\exp_X(t,p)$ is defined for all $t\in B^n(0,R)$.
        \item For every $0<r\le R$ and every $q\in B_X(p,R)$ such that $B^n(0,r)\times \{q\}\subseteq\Do_X$ and $\exp_X(B^n(0,r)\times \{q\})\subseteq B_X(p,R)$, we have $e^{t\cdot X}(q)\neq q$ for all $t\in B^n(0,r)\backslash\{0\}$.
    \end{itemize}
\end{defn}

\begin{remark}
In \cite{CoordAdapted} the second condition is stated a bit differently: for such $r$ and $q$ we have, either  $e^{t\cdot X}(q)\neq q$ for all $0<|t|<r$, or $X_1(q),\dots,X_n(q)$ are linearly dependent. The linear independence assumption is unnecessary because by the singular Frobenius theorem (Theorem \ref{MainThm::SingFro}, or \cite{MontanariMorbidelliSingularFrobenius} when $X_1,\dots,X_n$ are Lipschitz), if $X_1(p),\dots,X_n(p)$ are already linearly independent, then $X_1(q),\dots,X_n(q)$ will always be linearly independent for all $q\in B_X(p,\infty)$.

Suppose in addition that $X_1,\dots,X_n$ satisfy $[X_i,X_j]=\sum_{k=1}^nc_{ij}^kX_k$ for some $c_{ij}^k\in C^0(\Mf)$. Then for any $R_0>0$ there is a $0<R_1\le R_0$ depending only $R_0$ and the upper bounded of $\sup_{i,j,k}\|c_{ij}^k\|_{C^0}$, such that, 
\begin{itemize}[parsep=-0.3ex]
    \item If $B(0,R_0)\times\{p\}\subseteq\Do_X$, then $B^n(R_1)\times B_X(p,R_1)\subseteq\Do_X$ (see Proposition \ref{Prop::ODE::InjLemma} \ref{Item::ODE::InjLemma::DistBall}).
    \item If $\Mf$ has almost injective radius $R_0$, then $\Phi_p:B^n(0,R_1)\to\Mf$ is a  $C^1$-embedding (see Corollary \ref{Cor::ODE::InjRad}).
\end{itemize}
\end{remark}

In the following discussion we assume that $X_1,\dots,X_n$ are all $C^1$. The definability and injectivity of the canonical coordinates between different points are related by the following. Recall that the  $C^{0,s}_X$ spaces are defined in \eqref{Eqn::Intro::CsXNorm} and \eqref{Eqn::Intro::CCDist}.

\begin{prop}\label{Prop::ODE::InjLemma}
    Let $X_1,\dots,X_n$ be $C^1$ vector fields on a manifold $\Mf$. Suppose there are $c_{ij}^k\in C^0(\Mf)$, $1\le i,j,k\le n$ such that $[X_i,X_j]=\sum_{i,j,k=1}^nc_{ij}^kX_k$ for $1\le i,j\le n$.
    
    Let $p_0\in\Mf$ be a fixed point. Let $\mu_0>0$ be such that $B^n(0,\mu_0)\times\{p_0\}\subseteq\Do_X$. Let  $\Phi_0=\Phi_{p_0}:B^n(0,\mu_0)\to\Mf$, $\Phi_0(t):=e^{t\cdot X}(p_0)$ be the canonical coordinates at $p_0$.
    
    Let $r_0:=\min(\mu_0,1/M_0)$ where 
    \begin{equation}\label{Eqn::ODE::InjLemma::DefM0}
        M_0:=\sup_{x\in\Mf}\sup_{\substack{u,v,w\in\R^n;\\|u|=|v|=|w|=1}}\bigg|\sum_{i,j,k=1}^nc_{ij}^k(x)u^iv^jw_k\bigg|.
    \end{equation}
    Then we have the following:
    \begin{enumerate}[parsep=-0.3ex,label=(\roman*)]

        \item\label{Item::ODE::InjLemma::Diffeo} $\Phi_0:B^n(0,r_0)\to\Mf$ is a $C^1$-locally diffeomorphism.
        \item\label{Item::ODE::InjLemma::A} Write $X=[X_1,\dots,X_n]^\top$ and $Y=\Phi_0^*X$ on $B^n(0,r_0)$. Then $Y$ is of the form $Y=(I+A)\nabla$ where $A\in C^0(B^n(0,r_0);\R^{n\times n})$ satisfies $|A(t)|\le M_0|t|$. Here for a matrix $B=(b_i^j)\in\R^{n\times n}$ we use the standard $\ell^2$ matrix operator norm $|B|:=\sup_{u,v\in\R^n;|u|=|v|=1}|\sum_{i,j=1}^nu^iv_jb_i^j|$.
        \item\label{Item::ODE::InjLemma::DistNorm} In particular we have $\Co^s(B^n(0,r_0/2))=C^{0,s}_Y(B^n(0,r_0/2))$ for all $0<s<1$ with equivalent norms $(\frac23)^s\|f\|_{C^{0,s}_Y}\|f\|_{\Co^s}\le2^s\|f\|_{C^{0,s}_Y}$ for all $f\in\Co^s(B^n(0,r_0/2))$.
        \item\label{Item::ODE::InjLemma::DistBall} For every $0<r\le r_0/2$ we have $B_X(p_0,r)\subseteq \Phi_0(B^n(0,2r))$. In particular $B^n(0,r_1)\times B_X(p_0,r_2)\subseteq\Phi_0(B^n(0,2r_1+2r_2))\subseteq\Do_X$ for all $r_1,r_2>0$ such that $r_1+r_2\le r_0/2$.
        
        \item\label{Item::ODE::InjLemma::HoldNorm}Let $0<s<1$. Suppose $c_{ij}^k\in C^{0,s}_X(\Mf)$, then $A\in \Co^s(B^n(0,r_0/4);\R^{n\times n})$. Moreover, we let $$M_1:=\sup_{u,v,w\in\R^n;|u|=|v|=|w|=1}\bigg\|\sum_{i,j,k=1}^nu^iv^jw_kc_{ij}^k\bigg\|_{C^{0,s}_X(\Mf)},$$ then $\|A\|_{\Co^s(B^n(0,r_0/4);\R^{n\times n})}\le5(1+M_0+r_0M_1)$.

    \end{enumerate}
\end{prop}
\begin{remark}
\begin{enumerate}[parsep=-0.3ex,label=(\roman*)]
    \item The constant $r_0$ is almost a ``0-admissible constant'', see \cite[Definitions 4.1]{CoordAdapted} except it depends on $\eta>0$ in \cite[Section 3.2]{CoordAdapted}. The constant $5(1+M_0+r_0M_1)$ is called a ``$\{s\}$-admissible constant'', see \cite[Definition 3.10]{CoordAdapted}. We do not know whether we still have the similar estimates if we relax $X_1,\dots,X_n\in C^1$ and $(c_{ij}^k)\in C^0$ to $X_1,\dots,X_n\in\Co^\LogL$ and $(c_{ij}^k)\in\Co^{\LogL-1}$.
    \item In \ref{Item::ODE::InjLemma::HoldNorm} we see that $\Phi_0^*X_1,\dots,\Phi_0^*X_n\in\Co^s$. While under the same assumption, in Theorem \ref{MainThm::QuantFro} we can find an different pullback $\Phi$ such that $\Phi^*X_1,\dots,\Phi^*X_n\in\Co^{s+1}$, which gain 1 more derivative. In the proof of Theorem \ref{MainThm::QuantFro} (see Section \ref{Section::PfQuantFro}), the $\Phi$ is given by taking an additional pre-composition to $\Phi_0$. The additional pre-composition is necessary: in Lemma \ref{Lem::Exmp::CanonicalCoords} we see that the regularity statement $\Phi_0^*X\in\Co^s$ in general is the best possible.
\end{enumerate}
\end{remark}
\begin{proof}[Proof of Proposition \ref{Prop::ODE::InjLemma}]
Define a matrix function $C:B^n(0,\mu_0)\to\R^{n\times n}$ as 
\begin{equation}\label{Eqn::ODE::InjLemma:Pf::DefforC}
    C(t)_i^j:=\sum_{k=1}^nt^k\cdot c_{ik}^j(\Phi_0(t)).
\end{equation}
Thus by definition $|C(t)|\le M_0|t|$ in the sense of the matrix operator norm. By \cite[Lemma 9.2 and Proposition 9.1]{CoordAdapted} we know if $\tilde A\in C^0(B^n(0,\mu_0);\R^{n\times n})$ be satisfied
\begin{equation}\label{Eqn::ODE::InjLemma:Pf::ODEforA}
    \Coorvec r(r\tilde A(r\theta))=-\tilde A(r\theta)^2-C(r\theta)\tilde A(r\theta)-C(r\theta),\quad\text{for }|r|<\mu_0\text{ and }\theta\in\mathbb S^{n-1},\quad\text{with }\tilde A(0)=0,
\end{equation}
then the tangent map $(d\Phi_0)_t((I+\tilde A(t))\nabla)=X(\Phi_0(t))$ holds for all $t\in B^n(0,\mu_0)$.

First we give estimate for $\tilde A$. Indeed rewriting \eqref{Eqn::ODE::InjLemma:Pf::ODEforA} we have
\begin{equation}\label{Eqn::ODE::InjLemma:Pf::ODEforA2}
    \tilde A(t)=-\int_0^1(\tilde A(\rho t)^2+C(\rho t)\tilde A(\rho t)+C(\rho t))d\rho.
\end{equation}
Denote the right hand side by $\Tc\tilde A(t)$. We see that if $A_1,A_2\in C^0(B^n(0,\mu_0);\R^{n\times n})$ satisfy $|A_1(t)|,|A_2(t)|\le M_0|t|$ for $|t|\le \mu_0$, then
\begin{align*}
    |\Tc A_i(t)|&\le\int_0^1(M_0^2|\rho t|^2+M_0^2|\rho t|^2+M_0|\rho t|)d\rho\le\tfrac{2M_0^2}3|t|^2+\tfrac {M_0}2|t|.
\\
    \tfrac{|\Tc A_1(t)-\Tc A_2(t)|}{|t|}&\le\int_0^1|A_1(\rho t)+A_2(\rho t)+C(\rho t)|\tfrac{|A_1(\rho t)-A_2(\rho t)|}{|\rho t|}\rho d\rho\le M_0|t|\sup_{0\le \rho\le 1}\tfrac{|A_1(\rho t)-A_2(\rho t)|}{|\rho t|}.
\end{align*}

For $0<R<\frac1{M_0}$, we define $\Xs_R=\{A\in C^0(B^n(0,R);\R^{n\times n}):|A(t)|\le M_0|t| \}$ with metric $d(A_1,A_2):=\sup_{|t|<R}|t|^{-1}|A_1(t)-A_2(t)|$, we see that $\Tc:\Xs_R\to\Xs_R$ is a contraction map. Thus the solution $\tilde A$ for \eqref{Eqn::ODE::InjLemma:Pf::ODEforA} is unique in $\Xs_R$ for all $0<R<r_0$, hence $\tilde A\in C^0(B^n(0,r_0);\R^{n\times n})$ and $|\tilde A(t)|\le M_0|t|$ holds for all $t\in B^n(0,r_0)$. Therefore we get \ref{Item::ODE::InjLemma::A} once we obtain \ref{Item::ODE::InjLemma::Diffeo}.

Now $|\tilde A(t)|<M_0|t|$, thus we have the following absolute series:
\begin{equation*}\label{Eqn::ODE::InjLemma:Pf::Sum(I+A)-1}
    |(I+\tilde A(t))^{-1}|=\Big|\sum_{k=0}^\infty(-1)^k\tilde A(t)^k\Big|\le\sum_{k=0}^\infty M_0^k|t|^k=\frac1{1-M_0|t|}<\infty,\quad t\in B^n(0,r_0).
\end{equation*}

Therefore $(d\Phi_0)_t\nabla=(I+\tilde A(t))^{-1}\cdot(X\circ\Phi_0(t))$ with both sides being collections of $n$ linearly independent vector fields. Since $r_0\le\frac1{M_0}$, we conclude that $A=\tilde A$ in $B^n(0,r_0)$ and $\Phi_0:B^n(0,r_0)\to\Mf$ is $C^1$-locally diffeomorphism, finishing the proof of \ref{Item::ODE::InjLemma::Diffeo}.

\medskip
Now in we have matrix norms $|I+A(t)|\le\frac32$ and $|(I+A(t))^{-1}|\le2$ for $t\in B^n(0,r_0/2)$, thus the distance $|t_1-t_2|=\dist_\nabla(t_1,t_2)$ satisfies $\frac12\dist_Y(t_1,t_2)\le|t_1-t_2|\le\frac32\dist_Y(t_1,t_2) $ for all $t_1,t_2\in B^n(0,r_0/2)$. In other words for any continuous function $f:B^n(0,r_0/2)\to\R$ we have $(\frac32)^s\frac{|f(t_1)-f(t_2)|}{\dist_Y(t_1,t_2)^s}\le\frac{|f(t_1)-f(t_2)|}{|t_1-t_2|^s}\le2^s\frac{|f(t_1)-f(t_2)|}{\dist_Y(t_1,t_2)^s}$ for all $t_1,t_2\in B^n(0,r_0/2)$. Taking supremum over $t_1,t_2$ we get \ref{Item::ODE::InjLemma::DistNorm}.

\medskip
Let $0<r<r_0/2$ and let $q\in B_X(p_0,r)$, by definition since $\dist_X(p_0,q)<r$, we can find a Lipschitz curve $\gamma:[0,r]\to\Mf$ such that $\gamma(0)=p_0$, $\gamma(r)=q$ and $\dot\gamma(\tau)=\sum_{j=1}^nb^j(\tau)X_j(\gamma(t))$ ($0\le \tau\le r$) for some $b\in L^\infty([0,r];\B^n)$. Thus we define a correspondent curve $\zeta$ in $B^n(0,r_0)$ as
\begin{equation*}
    \dot\zeta(\tau)=\sum_{j=1}^nb^j(\tau)Y_j(\zeta(\tau)),\quad \zeta(0)=0,
\end{equation*}
for $0\le \tau\le r$ in the existence domain. Since $Y=\Phi_0^*X$ we have $\gamma(\tau)=\Phi_0\circ\zeta(\tau)$ for $\tau$ in the domain.

By the Maximal Existence Theorem, if $\zeta$ is only defined in $[0,r']$ for some $r'<r$ then $\limsup_{\tau\to r'}|\zeta(\tau)|=r_0$.

On the other hand, since $Y=(I+A)\nabla$ with $|A(t)|\le M_0|t|$, we have $|\dot\zeta(\tau)|\le|Y(\zeta(\tau))|\le1+M_0|\zeta(\tau)|$ for all $\tau$ in the domain. Therefore by Gronwall's inequality we have $|\zeta(\tau)|\le \frac1{M_0}(e^{M_0\tau}-1)\le2\tau$ for all $\tau\in[0,r]\subset[0,r_0\log 2]$ in the domain, in particular $\zeta(\tau)\in B^n(0,2r)\subset B^n(0,r_0)$. We conclude that $\zeta(\tau)\in B^n(0,2r)$ for all $\tau$ in the domain, thus $\zeta:[0,r]\to B^n(0,2r)$ is defined. Using $\gamma=\Phi_0\circ\zeta$ we get $q=\gamma(r)\in \Phi_0(B^n(0,2r))$, finishing the proof of $B_X(p,r)\subseteq\Phi_0(B^n(0,2r))$ for $0<r<r_0/2$.

For every $r_1>0$ and every $q\in B_X(p_0,r_1)$, we have $e^{t\cdot X}(q)\in B_X(p_0,r_1+|t|)$ whenever it is defined. By \ref{Item::ODE::InjLemma::DistBall}, when $|t|<r_2$ we have $e^{t\cdot X}(q)\in B_X(p_0,r_1+r_2)\subset \Phi_0(B^n(0,2r_1+2r_2))$ which is in particular defined for $r_1+r_2<r_0/2$. This completes the proof of \ref{Item::ODE::InjLemma::DistBall}.

\medskip
When $t\in B^n(0,r_0/4)$, by \ref{Item::ODE::InjLemma::A} and \eqref{Eqn::ODE::InjLemma:Pf::DefforC} we have $\|A\|_{C^0(B^n(0,r_0/4))}\le\frac14$ and $\|C\|_{C^0(B^n(0,r_0/4))}\le\frac14$. Thus by Lemma \ref{Lem::Hold::Product} \ref{Item::Hold::Product::Hold2} and \eqref{Eqn::ODE::InjLemma:Pf::ODEforA2} (since $A=\tilde A$) we have for $\Co^s=\Co^s(B^n(0,r_0/4);\R^{n\times n})$ and $C^0=C^0(B^n(0,r_0/4);\R^{n\times n})$,
\begin{equation*}
    \|A\|_{\Co^s}\le2\|A\|_{\Co^s}\|A\|_{C^0}+\|A\|_{\Co^s}\|C\|_{C^0}+\|C\|_{\Co^s}\|A\|_{C^0}+\|C\|_{\Co^s}\le\tfrac34\|A\|_{\Co^s}+\tfrac54\|C\|_{\Co^s}\le 5\|C\|_{\Co^s}.
\end{equation*}

To prove \ref{Item::ODE::InjLemma::HoldNorm} it remains to show $\|C\|_{\Co^s(B^n(0,r_0/4);\R^{n\times n})}\le1+M_0+r_0M_1$.

Let $t\in B^n(0,\frac{r_0}4)$, since $\nabla\Phi_0(t)=(I+A(t))^{-1}\cdot X(\Phi_0(t))$, by \eqref{Eqn::ODE::InjLemma:Pf::Sum(I+A)-1} (since $r_0\le\frac1{M_0}$) we have $\dist_X(\Phi_0(t_1),\Phi_0(t_2))\le|t_1-t_2|\sup_{t\in B^n(0,r_0/4)}|(I+A(t))^{-1}|\le\frac43|t_1-t_2|$ for all $t_1,t_2\in B^n(0,\frac{r_0}4)$. Therefore
\begin{align*}
    &\|C\|_{\Co^s(B^n(0,\frac{r_0}4))}=\|C\|_{C^0}+\sup_{|t_1|,|t_2|<r_0/4}\frac{\big|\big(\sum_{k=1}^n(t_1^kc_{ik}^j(\Phi_0(t_1))-t_2^kc_{ik}^j(\Phi_0(t_2))\big)_{i,j}\big|_{\R^{n\times n}}}{|t_1-t_2|^s}
    \\
    \le&\frac14+\sup_{|t_1|,|t_2|<r_0/4}\frac{\big|\big(\sum_{k=1}^n(t_1^k-t_2^k)c_{ik}^j(\Phi_0(t_1))\big)_{i,j}\big|_{\R^{n\times n}}+\big|\big(\sum_{k=1}^nt_2^k(c_{ik}^j(\Phi_0(t_1))-c_{ik}^j(\Phi_0(t_2))\big)_{i,j}\big|_{\R^{n\times n}}}{|t_1-t_2|^s}
    \\
    \le&\tfrac14+\sup_{|t_1|,|t_2|<r_0/4}\frac{M_0|t_1-t_2|}{|t_1-t_2|^s}+\tfrac{r_0}4\sup_{|t_1|,|t_2|<r_0/4}\frac{\sup_{|u|,|v|,|w|=1}|\sum_{i,j,k=1}^nu^iv_jw^k(c_{ik}^j(\Phi_0(t_1))-c_{ik}^j(\Phi_0(t_2)))|}{|t_1-t_2|^s}
    \\
    \le&\tfrac14+M_0(\tfrac{r_0}4)^{1-s}+\tfrac{r_0}4M_1\sup_{|t_1|,|t_2|<r_0/4}\frac{\dist_X(\Phi_0(t_1),\Phi_0(t_2))^s}{|t_1-t_2|^s}\le\tfrac14+\tfrac{M_0^s}{4^{1-s}}+\tfrac{(4/3)^s}4r_0M_1<1+M_0+r_0M_1.
\end{align*}
This completes the proof of \ref{Item::ODE::InjLemma::HoldNorm}.
\end{proof}

As a corollary to Proposition \ref{Prop::ODE::InjLemma}, the almost injective radius gives a lower bounds the actual injective radius, which is also  almost a ``0-admissible constant'' (see \cite[Definitions 4.1]{CoordAdapted}) except it depends on the quantity $\eta\in\R_+$ in \cite[Section 3.2]{CoordAdapted}.
\begin{cor}[Injective radius]\label{Cor::ODE::InjRad}
Let $X_1,\dots,X_n$ be $C^1$ vector fields on a $C^2$-manifold $\Mf$ that span the tangent space at every point. Let $p\in\Mf$ and let $R_0$ be the almost injective radius at $p$ (see Definition \ref{Defn::ODE::AlmostInjRad}).

Suppose $[X_i,X_j]=\sum_{k=1}^nc_{ij}^kX_k$ for some bounded functions $(c_{ij}^k)\subset C^0(\Mf)$. Let $M_0$ be as in \eqref{Eqn::ODE::InjLemma::DefM0}. Then $\Phi_p(t):=e^{t\cdot X}(p)$ defines a $C^1$-regular parameterization $\Phi_p:B^n(0,\frac1{10}\min(R_0,\frac1{M_0}))\to\Mf$, i.e. $\Phi_p$ is $C^1$-diffeomorphism onto its image with injective radius at least $\frac1{10}\min(R_0,\frac1{M_0})$.
\end{cor}

\begin{proof}
Let $r_0:=\min(R_0,1/M_0)$. By Proposition \ref{Prop::ODE::InjLemma}, $\Phi_p:B^n(0,r_0)\to\Mf$ is a locally $C^1$-diffeomorphism.

Let $q_1,q_2\in\Phi_p(B^n(0,\frac1{10}r_0))$, clearly $\dist_X(p,q_1)\le\frac1{10}r_0$. By Proposition \ref{Prop::ODE::InjLemma} \ref{Item::ODE::InjLemma::DistBall} with $\mu_0=R_0$, $p_0=p$ and $r_2=\frac1{10}r_0$ in the assumption, we get that $\Phi_{q_1}(t)=e^{t\cdot X}(q_1)$ is defined for $|t|<\frac12(r_0-\frac15r_0)=\frac25r_0$ and satisfies $\Phi_{q_1}(B^n(0,\frac25r_0))\subseteq\Phi_p(B^n(0,r_0))$. Therefore, applying Proposition \ref{Prop::ODE::InjLemma} \ref{Item::ODE::InjLemma::DistBall} again with $\mu_0=\frac25r_0$ and $p_0=q_1$, we see that $B_X(q_1,\frac15r_0)\subseteq\Phi_{q_1}(B^n(0,\frac25r_0))$. And since $\dist_X(q_1,q_2)\le\dist_X(p,q_1)+\dist_X(p,q_2)\le\frac15r_0$, we get $q_2\in \Phi_{q_1}(B^n(0,\frac25r_0))$, i.e. there is a $|t_1|<\frac25r_0$ such that $q_2=e^{t_1\cdot X}(q_1)$. 

Now by Definition \ref{Defn::ODE::AlmostInjRad} with $r=\frac25r_0$ we see that $q_2=e^{t_1\cdot X}(q_1)\neq q_1$ unless $t_1=0$. Therefore $\Phi_p:B^n(0,\frac1{10}r_0)\to\Mf$ is injective. On the other hand from Proposition \ref{Prop::ODE::InjLemma} \ref{Item::ODE::InjLemma::Diffeo} we already know $\Phi_p:B^n(0,\frac1{10}r_0)\to\Mf$ is  locally a $C^1$-diffeomorphism, thus $\Phi_p$ is an actual $C^1$-diffeomorphism onto its image.
\end{proof}

When $X_1,\dots,X_n$ form a basis of the tangent space  at every point in $\Mf$, we can endow $\Mf$ with a Riemannian metric $g$ such that $X_1,\dots,X_n$ form an orthonormal frame.

\begin{lem}\label{Lem::ODE::FlowVsMetric}
    Let $\Mf$ be a connnected $C^{1,1}$-manifold. Let $X_1,\dots,X_n$ be log-Lipschitz vector fields that  form a basis of the tangent space at every point. Let $g$ be the associated metric described as above. Let $\dist_g:\Mf\times\Mf\to[0,\infty)$ be the correspondent distance with respect to $g$. For $p\in\Mf$ and $R\in(0,\infty]$ we denote $\bar B_g(p,R)=\{q\in\Mf:\dist_g(p,q)\le R\}$. 
    
    Then $\dist_X=\dist_g$. Moreover following are equivalent:
    \begin{enumerate}[parsep=-0.3ex,label=(\alph*)]
        \item\label{Item::ODE::FlowVsMetric::Metric} $\bar B_g(p,R)$ is a complete metric space. In other words, $\bar B_g(p,R)$ is a compact space.
        \item\label{Item::ODE::FlowVsMetric::Flow} For any $b=(b^1,\dots,b^n)\in L^\infty([0,R];\R^n)$ there exists a unique Lipschitz curve $\gamma:[0,R]\to\Mf$ such that $\gamma(0)=p$ and $\dot\gamma(t)=\sum_{j=1}^nb^j(t)X_j(\gamma(t))$ for $0\le t\le R$.
    \end{enumerate}
    
    In particular, $\Do_X=\R^n\times\Mf$ if and only if $(\Mf,g)$ is a complete Riemannian manifold.
\end{lem}
\begin{proof} Recall that $\dist_g(p,q)$ is given by the infimum of $T>0$ such that there is a piecewise $C^1$-curve $\gamma:[0,T]\to\Mf$ such that $\gamma(0)=p$, $\gamma(T)=q$ and $|\dot\gamma(t)|_g\le1$ for almost every $0\le t\le T$.

For such $\gamma$, we can write $\dot\gamma(t)=\sum_{j=1}^na^j(t)X_j(\gamma(t))$ for uniquely determined functions $a^1,\dots,a^n\in L^\infty[0,T]$. By assumption $X_1,\dots,X_n$ are orthonormal, thus $\gamma$ is a subunit Lipschitz path (see Definition \ref{Defn::ODE::SubUnitCurve}) with respect to $X_1,\dots,X_n$. Taking infimum over $T$ with such curves $\gamma$, we conclude that $\dist_X(p,q)\le\dist_g(p,q)$.

On the other hand, if $T>\dist_X(p,q)$, and if $\gamma:[0,T]\to\Mf$ is a subunit Lipschitz path joining $p$ and $q$, then $|\dot\gamma(t)|_g\le1$ for almost every $0\le t\le T$. Taking mollification we see that for any $\eps>0$ that is a $C^{1,1}$-curve $\gamma_\eps:[0,T+\eps]\to\Mf$ joining $p$ and $q$, which is also subunit as well, thus $\dist_g(p,q)\le T+\eps$. Taking infimum over all $T>\dist_X(p,q)$ and all $\eps>0$ we get $\dist_g(p,q)\le\dist_X(p,q)$, finishing the proof of $\dist_X=\dist_g$.

\smallskip \noindent\ref{Item::ODE::FlowVsMetric::Metric}$\Rightarrow$\ref{Item::ODE::FlowVsMetric::Flow}: By \cite[Theorem 3.2]{BahouriCheminDanchin} if such $\gamma$ in the assumption exists, then it must be unique.

Suppose \ref{Item::ODE::FlowVsMetric::Flow} is not true. There is a $b\in L^\infty([0,R];\R^n)$ such that $\|b\|_{L^\infty}\le1$, and the ODE $\dot\gamma=b\cdot( X\circ\gamma)$, $\gamma(0)=p$ only exists in $[0,\tau)$ for some $\tau\le R$. By the maximal ODE existence theorem, for every compact subset $K\Subset\Mf$ there is a $\delta>0$ such that $\gamma((\tau-\delta,\tau))\in\Mf\backslash K$. Take $K=\bar B_g(p,R)$, which is compact by assumption, we see that $\dist_g(p,\gamma(\tau-\frac\delta2))>R$. Thus $\dist_X(p,\gamma(\tau-\frac\delta2))>R$ as well.
However by Definition \ref{Defn::ODE::SubUnitCurve} and \eqref{Eqn::ODE::CCDistNew}, we see that $\dist_X(p,\gamma(t))\le t$ for all $0<t<\tau$. Taking $t=\tau-\frac\delta2<R$ we get a contradiction.



\smallskip\noindent\ref{Item::ODE::FlowVsMetric::Flow}$\Rightarrow$\ref{Item::ODE::FlowVsMetric::Metric}: We let $\tilde \Mf$ be the completion of $(\Mf,\dist_g)$, which is still a metric space but may not be a manifold. We denote by $\tilde\dist_g$ the completion of $\dist_g$ in $\Mf$.
Let $R_0\in(0,\infty]$ be the smallest number such that $\bar B_g(p,r)\subset\Mf$ is completed for all $0<r<R_0$. In other words the distance ball in $B_{\tilde \Mf}(p,r)\subseteq\tilde\Mf$ satisfies $B_{\tilde \Mf}(p,r)\subset\Mf$ for all $0<r<R_0$. Clearly such $R_0$ exists. Note that if $\bar B_g(p,r)\subset\Mf$ is completed, since $\Mf$ is locally compact, we see that there is a $\delta>0$ such that $B_{\tilde \Mf}(p,r+\delta)=B_g(p,r+\delta)\subseteq\Mf$ holds. Thus if $R_0<\infty$ then the closed ball $\bar B_{\tilde\Mf}(p,R_0)=\{q\in\tilde\Mf:\tilde\dist_g(p,q)\le R_0\}\not\subset\Mf$.

Therefore to prove \ref{Item::ODE::FlowVsMetric::Metric} it suffices to prove $R_0>R$.

Suppose otherwise $R_0\le R$, then there is a sequence $\{q_k\}_{k=1}^\infty\subset\bar B_g(p,R_0)(\subset\Mf)$ such that $q_k\to q_\infty$ for some $q_\infty\in\bar B_{\tilde\Mf}(p,R_0)\backslash\Mf$.
Since we have $\dist_g=\dist_X$, by assumption there are subunit Lipschitz curves $\gamma_k:[0,R+\frac1k]\to\Mf$ such that $\gamma_k(0)=p$ and $\gamma_k(R+\frac1k)=q_k$ for all $k$.
Now $\{\gamma_k\in \Co^\Lip([0,R_0];\Mf)\}_{k=1}^\infty$ are equicontinuous since they are all subunit. By the Arzel\`a–Ascoli theorem there is a convergence subsequence in $C^0([0,R];\tilde\Mf)$, whose limit we denote by $\gamma_\infty:[0,R_0]\to\tilde\Mf$. Since $\dist_g(q_k,\gamma_k(R))\to0$ and $q_k\to q_\infty$, we see that $\gamma_\infty(R_0)=q_\infty$.

Since $\|\dot\gamma_k\|_{L^\infty_g}\le1$ for all $k$, we know $\tilde\dist_g(\gamma_\infty(t_1),\gamma_\infty(t_2)\le|t_1-t_2|$ for all $t_1,t_2\in[0,R_0]$, in particular $\gamma_\infty(t)\in B_{\tilde \Mf}(p,t)(\subseteq\tilde\Mf)$ for all $0<t\le R_0$. By assumption of $R_0$, we have $\gamma_\infty(t)\in\Mf$ for all $0\le t<R_0$, thus $\gamma_\infty:[0,R_0)\to\Mf$ is subunit and we can write $\dot\gamma_\infty(t)=\sum_{j=1}^nu^j(t)X_j(\gamma_\infty(t))$ for some $u\in L^\infty([0,R_0);\R^n)$ such that $\|u\|_{L^\infty}\le1$.

Let $b(t):=u(t)$ for $t<R_0$ and $b(t):=0$ for $t\in[R_0,R]$, we see that $b\in L^\infty([0,R_0];\R^n)$ also satisfies $\|b\|_{L^\infty}\le1$. Therefore by assumption \ref{Item::ODE::FlowVsMetric::Flow}, we get $\gamma_\infty(R_0)\in\Mf$, contradict to the assumption $q_\infty=\gamma_\infty(R_0)\notin\Mf$. This completes the direction \ref{Item::ODE::FlowVsMetric::Flow}$\Rightarrow$\ref{Item::ODE::FlowVsMetric::Metric}.

\smallskip $\Do_X=\R^n\times\Mf$ if and only if \ref{Item::ODE::FlowVsMetric::Flow} is true for all $R>0$, if and only if \ref{Item::ODE::FlowVsMetric::Metric} is true for all $R>0$, if and only if $(\Mf,g)$ is complete.
\end{proof}
\begin{remark}\label{Rmk::ODE::FlowVsMetric}
By $\dist_X=\dist_g$ in Lemma \ref{Lem::ODE::FlowVsMetric}, we have $C^{0,\alpha}_X(\Mf)=\Co^\alpha(\Mf,g)$ for $0<\alpha<1$. Here the space $\Co^\alpha(\Mf,g)$ is defined by replacing $\dist_X$ in \eqref{Eqn::Intro::CsXNorm} with $\dist_g$, or by replacing $|x-y|$ in Definition \ref{Defn::Intro::DefofHold} \ref{Item::Intro::DefofHold::<1} with $\dist_g(x,y)$. In particular by passing to local, we obtain $C^{0,\alpha}_{X,\loc}(\Mf)=\Co^\alpha_\loc(\Mf)$.
\end{remark}

\subsection{The Lipschitz ODE flow and a holomorphic Frobenius theorem}\label{Section::HoloFro}

In the real Frobenius theorem (Theorem \ref{KeyThm::RealFro::BLFro}) and the holomorphic Frobenius theorem (Proposition \ref{Prop::ODE::ParaHolFro}) we need ODE flows to construct the parameterizations. We need the following regularity estimates for ODE flows:

\begin{lem}[Flow regularity]\label{Lem::ODE::ODEReg}
Let $\alpha\in\R_\Eb^+$ such that $\alpha\ge\Lip$.  Let $X$ be a $\Co^\alpha$-vector field on an open set $U_0\subseteq\R^n_x$. Then for any $U_1\Subset U_0$ and $0<\tau_0<\frac12\dist(U_1,\partial U_0)\cdot\|X\|_{C^0(U_0)}^{-1}$, the map $\exp_X$ is defined on $(-\tau_0,\tau_0)\times U_1$, $\exp_X\in\Co^{\alpha+1,\alpha}((-\tau_0,\tau_0),U_1;U_0)$ and $\nabla_x\exp_X\in\Co^{\alpha,\alpha-1}((-\tau_0,\tau_0),U_1;U_0)$.
\end{lem}
Here for $\alpha=\infty$ we use the standard convention $\alpha+1=\infty$.
\begin{proof}
It suffices to prove that case where $\alpha\in\{k+\Lip,k+1+\LogL:k=0,1,2,\dots\}\cup(1,\infty)$, since by Lemma \ref{Lem::Hold::LimitArgument}, for $\gamma\in(1,\infty)$, if $\exp_X\in\Co^{\gamma+1-\eps,\gamma-\eps}_{t,x}$ for all $\eps>0$ then $\exp_X\in\Co^{(\gamma+1)-,\gamma-}_{x,s}$, and the same argument holds for $\nabla_x\exp_X$.

By scaling we can assume that $U_0\subseteq\B^n$.

Let $\chi\in C_c^\infty(U_0)$ be such that $\chi\equiv1$ in a neighborhood of $\{x:\dist(x,U_1)\le\frac12\dist(U_1,\partial U_0)\}$. We see that $e^{tX(\cdot,s)}|_{U_1}=e^{t(\chi X(\cdot,s))}|_{U_1}$ whenever $|t|<\tau_0$ and $s\in V_1$. Therefore, in $U_1$ we can replace $X$ by $\chi X$ in the discussion. In other words, we can assume $X\in\Co^\alpha_c(\frac12\B^n;\R^n)$ and the $\exp_X$ is globally defined. Clearly in this case $\exp_X(t,x)\equiv x$ for $|x|\ge\frac12$ and $t\in\R$.

It suffices to show that there is a $\tau>0$ such that $\exp_X\in\Co^\alpha((-\tau,\tau)\times\B^n;\R^n)$. Indeed since $\exp_X(t+s,x)=\exp_X(t,\exp_X(s,x))$ for all $t,s\in\R$, by Proposition \ref{Prop::Hold::QComp} \ref{Item::Hold::QComp::>1} that $\exp_X\in\Co^\alpha((-\tau,\tau)\times\B^n;\R^n)$ implies $\exp_X\in\Co^\alpha_\loc(\R\times\B^n;\R^n)$. And since $\partial_t\exp_X=X\circ\exp_X\in\Co^\alpha_\loc(\R\times\B^n;\R^n)$, we get $\sup_x\|\exp_X(\cdot,x)\|_{\Co^{\alpha+1}(-T,T)}<\infty$ for all $T$, which concludes $\exp_X\in\Co^{\alpha+1,\alpha}_{t,x,\loc}$. Note that $\partial_t\nabla_x(e^{tX})=(\nabla X)(e^{tX})\cdot\nabla_x(e^{tX})$, thus $\nabla_x\exp_X\in\Co^{\alpha-1}_{t,x,\loc}$ and $(\nabla X)\circ\exp_X\in\Co^{\alpha-1}_{t,x,\loc}$ implies $\nabla_x\exp_X\in\Co^{\alpha,\alpha-1}_{t,x,\loc}$, which will complete the proof.

\medskip
First we see that $\phi(t,x):=\exp_X(t,x)$ is the unique fixed point $\phi=\Tc\phi$ where $\Tc: C^0_\loc(\R\times \B^n;\R^n)\to C^0_\loc(\R\times \B^n;\R^n)$ is given by
\begin{equation*}
    \Tc\phi(t,x):=x+\int_0^tX(\phi(r,x))dr,\quad t\in\R,\quad x\in\B^n.
\end{equation*}
We have 
\begin{equation*}
    (\Tc\phi_1-\Tc\phi_2)(t,x)=\int_0^t\int_0^1(\nabla X)(u\phi_1(r,x)+(1-u)\phi_2(r,x))\cdot(\phi_1(r,x)-\phi_2(r,x))duds.
\end{equation*}

Set $\Ic(t,x):=x$. Therefore if $\alpha\ge1+\Lip$, then by Proposition \ref{Prop::Hold::QComp} \ref{Item::Hold::QComp::>1} along with a scaling, we see that there is a $M=M(n,\alpha)>1$ such that
\begin{align*}
    \|\Tc\phi_1-\Ic\|_{\Co^{\alpha-1}((-\tau,\tau)\times\B^n;\R^n)}&\le M\tau\| X\|_{\Co^\alpha(\B^n;\R^n)},
    \\
    \|\Tc\phi_1-\Tc\phi_2\|_{\Co^{\alpha-1}((-\tau,\tau)\times\B^n;\R^n)}&\le M\tau\|\nabla X\|_{\Co^{\alpha-1}(\B^n;\R^{n\times n})}\|\phi_1-\phi_2\|_{\Co^{\alpha-1}((-\tau,\tau)\times\B^n;\R^n)},
\end{align*}
for every $\tau>0$ and every $\phi_1,\phi_2\in \Co^{\alpha-1}((-\tau,\tau)\times\B^n;\R^n)$ such that $\|\phi_1-\Ic\|_{\Co^{\alpha-1}},\|\phi_2-\Ic\|_{\Co^{\alpha-1}}\le1$.

Take $0<\tau_1<(2M\|X\|_{\Co^\alpha})^{-1}$, we see that $\Tc$ is a contraction map in $\{\phi\in \Co^{\alpha-1}((-\tau_1,\tau_1)\times\B^n;\R^n):\|\phi-\Ic\|_{\Co^{\alpha-1}}\le1\}$. By uniqueness of the fixed point we get that
\begin{equation}\label{Eqn::ODE::ODEReg::PfTmp1}
    \exp_X\in\Co^{\alpha-1}((-\tau_1,\tau_1)\times\B^n;\R^n),\text{ when }\alpha\ge1+\Lip.
\end{equation}

To get one derivative back, notice that $\nabla_x\exp_X-I_n:\R\times\B^n\to\R^{n\times n}$ is unique solution $\psi\in C^0_\loc(\R\times\B^n;\R^{n\times n})$ to the affine linear equation $\psi=I_n+\Se\psi$, where
\begin{equation*}
    \Se\psi(t,x):=\int_0^t(\nabla X)(\exp_X(r,x))\cdot\psi(r,x)dr.
\end{equation*}

By \cite[Proposition 3.10]{BahouriCheminDanchin}, since $X$ is Lipschitz, we see that $\exp_X-\id\in\Co^{\Lip}(\R,\R^n;\R^n)$. Combining with \eqref{Eqn::ODE::ODEReg::PfTmp1} and Proposition \ref{Prop::Hold::QComp} \ref{Item::Hold::QComp::>1} we see that $(\nabla X)\circ\exp_X\in\Co^{\alpha-1}$, when $\alpha\neq2$. Therefore when $\alpha\neq2$, for every $0<\tau\le\tau_1$,
\begin{equation}\label{Eqn::ODE::ODEReg::PfTmp2}
    \|\Se\psi\|_{\Co^{\alpha-1}((-\tau,\tau)\times\B^n;\R^{n\times n})}\le\tau\|(\nabla X)\circ\exp_X\|_{\Co^{\alpha-1}}\|\psi\|_{\Co^{\alpha-1}((-\tau,\tau)\times\B^n;\R^{n\times n})},\quad\psi\in\Co^{\alpha-1}((-\tau,\tau)\times\B^n;\R^{n\times n}).
\end{equation}
When $\alpha\neq2$, take $0<\tau_2<\tau_1$ such that $\tau_2<(2\|(\nabla X)\circ\exp_X\|_{\Co^{\alpha-1}})^{-1}$, we see that $\Se$ is a contraction linear map in $\Co^{\alpha-1}((-\tau_2,\tau_2)\times\B^n;\R^{n\times n}) $. Thus $\nabla\exp_X\in\Co^{\alpha-1}((-\tau_2,\tau_2)\times\B^n;\R^{n\times n})$, which implies $\exp_X\in\Co^\alpha((-\tau,\tau)\times\B^n;\R^n)$ for some small $\tau>0$. Therefore $\exp_X\in\Co^\alpha_\loc(\R\times\B^n;\R^n)$ when $\alpha\neq2$.

For the remaining case $\alpha=2$. Note that in the above discussion we have $\exp_X\in\Co^{2-\eps}_\loc(\R\times\B^n;\R^n)$ for all $\eps>0$. By Proposition \ref{Prop::Hold::QComp} \ref{Item::Hold::QComp::>1}, since $\Co^1\circ\Co^\frac32\subseteq\Co^1$, we see that $(\nabla X)\circ\exp_X\in\Co^1$ holds, thus \eqref{Eqn::ODE::ODEReg::PfTmp2} is still valid for $\alpha=1$ and $\tau>0$. By taking $0<\tau_2<(2\|(\nabla X)\circ\exp_X\|_{\Co^{\alpha-1}})^{-1}$, we see that  $\nabla\exp_X\in\Co^1((-\tau_2,\tau_2)\times\B^n;\R^{n\times n})$. This concludes that $\exp_X\in\Co^\alpha_\loc(\R\times\B^n;\R^n)$ when $\alpha=2$ and finishes the proof.
\end{proof}

In applications, we consider the multi-flow associated with commutative vector fields.

\begin{cor}\label{Cor::ODE::MultFlowReg} Let $\alpha\in\{\logl,k+\LogL,k+\Lip:k=0,1,2,\dots\}\cup[1,\infty]$ and let $X=[X_1,\dots,X_r]^\top$ be a collection of commutative $\Co^\alpha$-vector fields on an open set $U_0\subseteq\R^n_x$. Then for any $U_1\Subset U_0$ there is a $\tau_0>0$ such that the multi-flow $\exp_X(t,x):=e^{t^1X_1+\dots+t^rX_r}(x)$ is defined for $t\in B^r(0,\tau_0)$ and $x\in U_1$. Moreover
\begin{enumerate}[parsep=-0.3ex,label=(\roman*)]
    \item\label{Item::ODE::MultFlowReg::Reg} When $\alpha\ge\Lip$, $\exp_X\in\Co^{\alpha+1,\alpha}(B^r(0,\tau_0),U_1;U_0)$.
    
    When $\alpha=\logl$, $\exp_X\in\Co^{1+\LogL,1-}(B^r(0,\tau_0),U_1;U_0)$.
    
    When $\alpha=\LogL$, for every $\eps>0$ there is a $\tau_\eps\in(0,\tau_0)$ such that $\exp_X\in\Co^{1+\LogL,1-\eps}(B^r(0,\tau_\eps),U_1;U_0)$.
    
    When $\alpha=1$, then $\exp_X\in\Co^2L^\infty(B^r(0,\tau_0),U_1;U_0)$. In particular for the same $\tau_\eps\in(0,\tau_0)$ we have $\exp_X\in\Co^{2,1-\eps}(B^r(0,\tau_\eps),U_1;U_0)$.
    
    \item\label{Item::ODE::MultFlowReg::Comm} $\Coorvec{t^j}\exp_X(t,x)=X_j\circ\exp_X(t,x)$ for $j=1,\dots,r$.
\end{enumerate}
\end{cor}
\begin{proof}\ref{Item::ODE::MultFlowReg::Comm} follows from Lemma \ref{Lem::ODE::CommMultFlow}.

By Lemmas \ref{Lem::ODE::ODEReg} and \ref{Lem::ODE::LogLipODEReg} where we have regularity estimates of $\exp_{X_j}$ for each $1\le j\le r$. Taking compositions and applying Lemma \ref{Lem::Hold::CompofMixHold} \ref{Item::Hold::CompofMixHold::Comp} we get \ref{Item::ODE::MultFlowReg::Reg}.
\end{proof}


The ODE flows also work for holomorphic vector fields on complex manifolds. Here for a map $f:U\subseteq\C^m\to\C^n$, we recall the notations $\partial_zf=(\frac{\partial f}{\partial z^1},\dots,\frac{\partial f}{\partial z^n})$ and  $\partial_{\bar z}f=(\frac{\partial f}{\partial\bar z^1},\dots,\frac{\partial f}{\partial\bar z^n})$.

\begin{lem}[Holomorphic ODE flow with parameter]\label{Lem::ODE::HoloFlow}
    Let $U_0\subseteq\C^m$ and $V_0\subseteq\R^q$ be bounded open sets that are endowed with standard complex coordinate system $z=(z^1,\dots,z^m)$ and real coordinate system $s=(s^1,\dots,s^q)$ respectively. 
    
    Let $\beta\in\R_\Eb^+$, and let $Z=(Z^1,\dots,Z^m):U_0\times V_0\to\C^m$ be a $\Co^\beta$-map that is holomorphic in $z$-variable. Then
    \begin{enumerate}[parsep=-0.3ex,label=(\roman*)]
        \item\label{Item::ODE::HoloFlow::DefReg} For any $U_1\times V_1\Subset U_0\times V_0$, let $\tau_0:=\dist(U_1,\partial U_0)\cdot(4\|Z\|_{C^0(U_0\times V_0)})^{-1}$. There is a unique continuous map $\Phi^Z:\tau_0\B^2\times U_1\times V_1\subseteq\C^1_t\times\C^m_z\times\R^q_s\to U_0$ such that $\Phi^Z$ is holomorphic in $(t,z)$ and satisfies
    \begin{equation}\label{Eqn::ODE::HoloFlow::EqnDef}
        \textstyle\frac{\partial\Phi^Z}{\partial t}(t,z,s)=Z(\Phi^Z(t,z,s)),\quad\frac{\partial\Phi^Z}{\partial\bar t}(t,z,s)=0,\quad\Phi^Z(0,z,s)=z,\quad\forall t\in \tau_0\B^2,\ z\in U_1,\ s\in V_1.
    \end{equation}
    In addition, $\Phi^Z$ is holomorphic in $( t,z)$ and $\Phi^Z\in\Co^\infty\Co^\beta(\tau_0\B^2\times U_1,V_1;\C^m)$.
    \item\label{Item::ODE::HoloFlow::Comm} Moreover, suppose $W:U_0\times V_0\to\C^m$ is another such map as $Z$, such that two vector fields $\sum_{j=1}^mZ^j(z,s)\Coorvec{z^j}$ and $\sum_{j=1}^mW^j(z,s)\Coorvec{z^j}$  are commutative. Let $\Phi^Z$ and $\Phi^W$ be as in \eqref{Eqn::ODE::HoloFlow::EqnDef}. Then for any for any $U_1\times V_1\Subset U_0\times V_0$ there are an $\tau_0>0$ such that
    \begin{equation}\label{Eqn::ODE::HoloFlow::FlowComm}
        \Phi^Z( t_1,\Phi^W( t_2,z,s),s)=\Phi^W( t_2,\Phi^Z( t_1,z,s),s)=\Phi^{t_1Z+t_2 W}(1,z,s),s),\ \forall  t_1, t_2\in\tau_0\B^2,\ z\in U_1,\ s\in V_1.
    \end{equation}
    \end{enumerate}
\end{lem}
\begin{remark}
    If we forget the $s$-variable, the result is indeed the well-definedness of the holomorphic ODE flow $e^{ t Z}(z)$. In this way, we can write $\Phi^Z$ as 
    \begin{equation}\label{Eqn::ODE::HoloFlow::DefExpZ}
        \Phi^Z( t,z,s)=\exp_{Z(\cdot,s)}( t,z)=e^{ t Z(\cdot,s)}(z).
    \end{equation}
    
    And thus \eqref{Eqn::ODE::HoloFlow::FlowComm} is saying that $e^{ t_1Z(\cdot,s)}e^{ t_2W(\cdot,s)}=e^{ t_2W(\cdot,s)}e^{ t_1Z(\cdot,s)}=e^{t_1Z_1(\cdot,s)+t_2W(\cdot,s)}$.
\end{remark}

\begin{proof}[Proof of Lemma \ref{Lem::ODE::HoloFlow}]By Lemma \ref{Lem::Hold::NablaHarm} we see that $Z\in\Co^\infty_\loc\Co^\beta_\loc(U_0,V_0;\C^m)$ and $\partial_zZ\in\Co^\infty_\loc\Co^\beta_\loc(U_0,V_0;\C^m)$. By shrinking $U_0,V_0$ and expanding $U_1,V_1$, we can assume that $U_0,V_0,U_1,V_1$ are all smooth domains, and $Z,\partial_zZ$ both have bounded norms.

Write $(x,y):=(\re z,\im z)$. We define real vector fields $X$ and $Y$ on $U_0\subseteq\R^{2m}_{x,y}$ with parameter on $V_0\subseteq\R^{2m+q}_s$ as
\begin{equation}\label{Eqn::ODE::HoloFlow::PfXY}
    X(x,y,s)=\sum_{j=1}^m\re Z^j(z,s)\Coorvec{x^j}+\im Z^j(z,s)\Coorvec{y^j}, \quad Y(x,y,s)=\sum_{j=1}^m-\im Z^j(z,s)\Coorvec{x^j}+\re Z^j(z,s)\Coorvec{y^j}.
\end{equation}
By assumption, $Z$ is holomorphic in $z$, so $[X,Y]=0$ is commutative. Clearly for each $s\in V_1$, $\exp_{(X(\cdot,s),Y(\cdot,s))}:\tau_0\B^2\times U_1\to U_0$ is defined.

Let $\tilde\Phi=(\tilde\Phi^R,\tilde\Phi^I):=\exp_{(X,Y)}$,  where $\tilde\Phi^R,\tilde\Phi^I:\tau_0\B^2\times U_1\times V_1\to\R^m$ are the components of $x$ and $y$-coordinates respectively. Thus, by direct computation we see that $\Phi^Z:=\tilde\Phi^R+i\tilde\Phi^I$ solves \eqref{Eqn::ODE::HoloFlow::EqnDef} and is holomorphic in $(t,z)$. Conversely if $\Phi^Z$ solves \eqref{Eqn::ODE::HoloFlow::EqnDef} then we must have $\exp_{(X,Y)}=(\re\Phi^Z,\im\Phi^Z)$. Therefore we have the existence and uniqueness of $\Phi^Z$.

To show $\Phi^Z\in\Co^\infty_{t,z}\Co^\beta_s$ in the domain, by Lemma \ref{Lem::Hold::NablaHarm} it suffices to show $\Phi^Z\in\Co^{\beta+1,\beta}_{(t,z),s}$. By Definition \ref{Defn::Intro::DefofHold} \ref{Item::Intro::DefofHold::-} we only need to prove the case $\beta\in(0,\infty)\cup\{k+\LogL,k+\Lip:k=0,1,2,\dots\}$, since the domain $\tau_0\B^2\times U_0$ does not depend on $\beta$.

By the group law $\Phi^Z(t_1+t_2,z,s)=\Phi^Z(t_1,\Phi^Z(t_2,z,s),s)$ and Lemma \ref{Lem::Hold::CompofMixHold} \ref{Item::Hold::CompofMixHold::Comp}, it suffices to prove $\Phi^Z\in\Co^{\beta+1,\beta}_{(t,z),s}(\tau_0'\B^2,U_1;U_0)$ for some $\tau_0'\in(0,\tau_0]$.

Let $\tau_1:=\frac12\dist(U_1,\partial U_0)\cdot(1+\|Z\|_{C^0(U_0\times V_0)})^{-1}$, and let $\Ic( t,z,s)=z$. Similar to the proof of Lemma \ref{Lem::ODE::ODEReg} we see that the map $\Phi^Z$ is the $\phi$ for the equation $\phi=\Tc\phi$, where $\Tc:\{\phi\in C^0(\tau_1\B^2\times U_1\times V_1;\C^m):\partial_{\bar  t,\bar z}\phi=0,\|\phi-\Ic\|_{C^0}\le\tau_1\}\to C^0(\tau_1\B^2\times U_1\times V_1;\C^m)$ is given by
\begin{equation}\label{Eqn::ODE::HoloFlow::PfTmp1}
    \Tc\phi( t,z,s):=z+\int_0^ t Z(\phi(r,z,s),s)dr=z+\int_0^1 t\cdot Z(\phi(\rho t,z,s),s)d\rho.
\end{equation}
Here $\|\phi-\Ic\|_{C^0}\le\tau_1$ implies $\phi(U_1\times V_1)\subset U_0$, thus $Z(\phi(r,z,s),s)$ is defined. Since $\phi$ is holomorphic in $( t,z)$, the map $( t,z)\mapsto Z(\phi( t,z,s),s)$ is also holomorphic. Thus $\Tc\phi$ is holomorphic in $( t,z)$.

Clearly for $ t\in\tau_1\B^2$, $z\in U_1$ and $s\in V_1$,
\begin{equation}\label{Eqn::ODE::HoloFlow::PfTmp2}
    (\Tc\phi_1-\Tc\phi_2)( t,z,s):=\int_0^ t\int_0^1 \big((\partial_z Z)\circ(\lambda\phi_1(r,z,s)+(1-\lambda)\phi_2(r,z,s),s)\big)\cdot(\phi_1(r,z,s)-\phi_2(r,z,s) )d\lambda dr.
\end{equation}

Since $Z,\partial_zZ\in\Co^{\beta+1,\beta}_{z,s}$ and $\phi(\tau_1\B^2,B^m(z,\tau_1),s)\subseteq B^m(z,2\tau_1)\subseteq U_0$ for every $B^m(z,\tau_1)\subseteq U_1$ and $s\in V_1$. By Lemma \ref{Lem::Hold::QPComp} \ref{Item::Hold::QComp::>1} along with a scaling, we see that there is a $M=M(U_1,\tau_0,\beta,Z)>\|Z\|_{C^0}+\|\partial_zZ\|_{C^0}$ such that 
\begin{equation}\label{Eqn::ODE::HoloFlow::PfTmp3}
    \|\phi-\Ic\|_{\Co^\beta(\tau_1\B^2\times U_1\times V_1;\C^m)}\le\tau_1\Rightarrow\|Z(\phi,\id_s)\|_{\Co^\beta(\tau_1\B^2\times U_1\times V_1;\C^m)}+\|(\partial_zZ)(\phi,\id_s)\|_{\Co^\beta(\tau_1\B^2\times U_1\times V_1;\C^{n\times n})}\le M.
\end{equation}

For every $\tau\in(0,\tau_1]$, let 
\begin{align*}
    \Xs_\tau&:=\{\phi\in C^0(\tau\B^2\times U_1\times V_1;\C^m):\partial_{\bar t,\bar z}\phi=0, \|\phi-\Ic\|_{C^0}\le\tau_1\};\\
    \Ys_\tau&:=\{\phi\in \Co^\beta(\tau\B^2\times U_1\times V_1;\C^m):\partial_{\bar t,\bar z}\phi=0, \|\phi-\Ic\|_{\Co^\beta}\le\tau_1\}.
\end{align*}
Here recall that we assume $\beta\in(0,\infty)\cup\{k+\LogL,k+\Lip:k=0,1,2,\dots\}$, where $\Co^\beta(\tau\B^2\times U_1\times V_1;\C^m)$ is a Banach space.

 So $\Xs_\tau\subseteq C^0(\tau\B^2\times U_1\times V_1;\C^m)$ and $\Ys_\tau\subseteq\Co^\beta(\tau\B^2\times U_1\times V_1;\C^m)$ are both complete metric space, and $\Xs_\tau\supseteq\Ys_\tau$ for all $\tau$. By \eqref{Eqn::ODE::HoloFlow::PfTmp1}, \eqref{Eqn::ODE::HoloFlow::PfTmp2} and \eqref{Eqn::ODE::HoloFlow::PfTmp3} we get
\begin{align*}
    \|\Tc\phi_1-\Ic\|_{C^0}\le \tau\|Z\|_{C^0}\le\tau M,\quad
    \|\Tc\phi_1-\Tc\phi_2\|_{C^0}\le\tau M\|\phi_1-\phi_2\|_{C^0},&\quad\forall \phi_1,\phi_2\in\Xs_\tau;
    \\
    \|\Tc\phi_1-\Ic\|_{\Co^\beta}\le \tau M,\quad
    \|\Tc\phi_1-\Tc\phi_2\|_{\Co^\beta}\le\tau M\|\phi_1-\phi_2\|_{\Co^\beta},&\quad\forall \phi_1,\phi_2\in\Ys_\tau.
\end{align*}

Take $\tau_0':=\tau_1/(2M)$, we see that $\Tc$ is a contraction map in both $\Xs_{\tau_0'}$ and $\Ys_{\tau_0'}$. Thus by uniqueness of the fixed point $\Phi^Z\in C^0(\tau_0'\B^2\times U_1\times V_1;\C^m)$ is holomorphic in $(t,z)$ and satisfies $\Phi^Z\in \Co^\beta(\tau_0'\B^2\times U_1\times V_1;\C^m)$. By Lemma \ref{Lem::Hold::NablaHarm} again we get $\Phi^Z\in\Co^\infty_{t,z}\Co^\beta_s(\tau_0'\B^2\times U_1,V_1;\C^m)$, in particular $\Phi^Z\in\Co^{\beta+1,\beta}_{(t,z),s}$. This proves \ref{Item::ODE::HoloFlow::DefReg}.

\smallskip
\noindent \ref{Item::ODE::HoloFlow::Comm}:  We can fix a $s\in V_1$ in the discussion below, thus $Z$ and $W$ are holomorphic vector fields in $U_0$.

 Let $W=(W^1,\dots,W^m)$ be from the assumption, we define vector fields $X'$ and $Y'$ to be the ``real and imaginary part of $W$'' by replacing $Z^j$ with $W^j$ in \eqref{Eqn::ODE::HoloFlow::PfXY}. Similarly we get $\exp_{(X',Y')}=(\re\Phi^W,\im\Phi^W)$ in the domain. Now $[X,Y]=0$ and $[X',Y']=0$ since $Z$ and $W$ are holormorphic in $z$. The assumption $[Z,W]=0$ implies that $X,X',Y,Y'$ are all pairwise commutative. 

By Corollary \ref{Cor::ODE::MultFlowReg} \ref{Item::ODE::MultFlowReg::Comm} we have 
\begin{equation}\label{Eqn::ODE::HoloFlow::PfTmp4}
    e^{u_1X}e^{v_1Y}e^{u_2X'}e^{v_2Y'}=e^{u_2X'}e^{v_2Y'}e^{u_1X}e^{v_1Y}=e^{u_1X+v_1Y+u_2X'+v_2Y'}\quad\text{in }U_1,\text{ for }u_1,u_2,v_1,v_2\in\R\text{ small }.
\end{equation}By writing out the $\R^m_x$ and $\R^m_y$ components of \eqref{Eqn::ODE::HoloFlow::PfTmp4}, we get $\Phi^Z(t_1,\Phi^W(t_2,z))=\Phi^W(t_2,\Phi^Z(t_1,z))=\Phi^{t_1Z+t_2W}(1,z)$ for $z\in U_1$ and $t_1,t_2\in\C$ closed to $0$. By shrinking $\tau_0$, we get \eqref{Eqn::ODE::HoloFlow::FlowComm} and finish the proof.
\end{proof}

We can now prove a holomorphic Frobenius theorem with parameters.

We endow $\C^n$ and $\R^q$ with standard (real and complex) coordinate system $\zeta=(\zeta^1,\dots,\zeta^n)$ and $s=(s^1,\dots,s^q)$ respectively.

\begin{prop}[Holomorphic Frobenius theorem with parameter]\label{Prop::ODE::ParaHolFro}
    Let $\tilde U_0\subseteq\C^n$ and $V_0\subseteq\R^q$. Let $1\le m\le n$ and let $\beta\in\R_\Eb^+$ be a positive generalized index. And let $Z_j=\sum_{k=1}^nb_j^k(\zeta,s)\Coorvec{\zeta^k}$, $j=1,\dots,m$ be complex vector fields on $\tilde U_0\times V_0$ such that 
    \begin{itemize}[nolistsep]
        \item $Z_1,\dots,Z_m$ are linearly independent and pairwise commutative at every point in $\tilde U_0\times V_0$.
        \item $b_j^k\in\Co^\beta_\loc(\tilde U_0\times V_0;\C^n)$ are holomorphic in $\zeta$-variable.
    \end{itemize}
    
    Then for any $(\zeta_0,s_0)\in \tilde U_0\times V_0$ there is a neighborhood $\tilde U\times V\subseteq \tilde U_0\times V_0$ of $(\zeta_0,s_0)$ and a $\Co^\beta$-map $\widetilde G'':\tilde U\times V\to\C^{n-m}$, such that 
    \begin{enumerate}[parsep=-0.3ex,label=(\roman*)]
        \item\label{Item::ODE::ParaHolFro::0}$\widetilde G''(\zeta_0,s_0)=0$.
        \item\label{Item::ODE::ParaHolFro::Reg} $\widetilde G''\in\Co^\infty\Co^\beta(\tilde U,V;\C^{n-m})$ is holomorphic in $\zeta$-variable. In particular $\widetilde G''\in\Co^{\infty,\beta}_{\zeta,s,\loc}$ and $\nabla_\zeta \widetilde G''\in\Co^{\infty,\beta}_{\zeta,s,\loc}$. 
        \item\label{Item::ODE::ParaHolFro::Span} For each $s\in V$, the contangent subbundle $\Span(Z_1(\cdot,s),\dots,Z_m(\cdot,s))^\bot|_{\tilde U}\le \C T^*\tilde U$ has a holomorphic local basis $(d\widetilde G''^1(\cdot,s),\dots,d\widetilde G''^{n-m}(\cdot,s),d\bar \zeta^1,\dots,d\bar \zeta^n)$.
    \end{enumerate}
    
    Let  $U:=\tilde U\cap\R^n\subseteq\R^n$ be the open set in the real space endowed with the real coordinates $(\xi^1,\dots,\xi^n)$ where $\xi=\re\zeta$. Let $X_j:=\sum_{k=1}^nb_j^k(\xi+i0,s)\Coorvec{\xi^k}$ for $1\le j\le m$ be the ``real domain restriction'' of $Z_1,\dots,Z_m$. Then
    \begin{enumerate}[parsep=-0.3ex,label=(\roman*)]\setcounter{enumi}{3}
        \item\label{Item::ODE::ParaHolFro::SpanR} For each $s\in V$, the cotangent subbundle $\Span(X_1(\cdot,s),\dots,X_m(\cdot,s))^\bot\le \C T^*U$ has a real-analytic local basis $\big(d(\widetilde G''^1(\cdot,s)|_U),\dots,d(\widetilde G''^{n-m}(\cdot,s)|_U)\big)$.
    \end{enumerate}
\end{prop}
\begin{remark}When $\beta>\Lip$, the result \ref{Item::ODE::ParaHolFro::SpanR} is the same as to say that the $\Co^\beta$-subbundle $\Span(X_1,\dots,X_m)\le \C T^*(U\times V)$ has a local basis $(d(\widetilde G''^1|_{U\times V}),\dots,d(\widetilde G''^{n-m}|_{U\times V}),ds^1,\dots,ds^q)$. Similar result holds for \ref{Item::ODE::ParaHolFro::Span}.

We state \ref{Item::ODE::ParaHolFro::Span} and \ref{Item::ODE::ParaHolFro::SpanR} in this way because the pointwise span of $dG''$ or $d(G''|_{U\times V})$ does not make sense when $\beta\le1$. See also Convention \ref{Conv::EllipticPara::CoSpan} and Remark \ref{Rmk::EllipticPara::CoSpan}.
\end{remark}
\begin{proof}[Proof of Proposition \ref{Prop::ODE::ParaHolFro}]
Endow two complex spaces $\C^m$ and $\C^{n-m}$ with standard complex coordinate $z=(z^1,\dots,z^m)$ and $w=(w^1,\dots,w^{n-m})$.

Since $Z_1,\dots,Z_m$ are linearly independent, the span of $Z_1,\dots,Z_m$ is a rank $m$ complex subbundle on $\tilde U_0\times V_0$. Hence we can find an affine complex linear map $\Gamma:\C^{n-m}_w\to\C^n_\zeta$ such that $\Gamma(0)=\zeta_0$ and 
\begin{equation}\label{Eqn::ODE::PfParaHolFro::LinInd}
    \partial_{w^1}\Gamma,\dots,\partial_{w^{n-m}}\Gamma,Z_1|_{(\zeta_0,s_0)},\dots,Z_m|_{(\zeta_0,s_0)}\in\C T_{\zeta_0} \tilde U_0\text{ are $\C$-linearly independent}.
\end{equation}

We define $\Co^\beta$-maps $\Psi(z,w,s)$ to $\tilde U_0$ and $\widehat\Psi(z,w,s)$ to $\tilde U_0\times V_0$ as
\begin{equation*}
    \Psi(z,w,s):=e^{z\cdot Z}(\Gamma(w),s+s_0)=\big(e^{z^1 Z_1(\cdot,s+s_0)+\dots +z^m Z_m(\cdot,s+s_0)}(\Gamma(w))\big),\quad\widehat\Psi(z,w,s):=(\Psi(z,w,s),s+s_0).
\end{equation*}
Here $e^{z^jZ_j(\cdot,s)}(\zeta)=\Phi^{Z_j}(z_j,\zeta,s)$ is given by \eqref{Eqn::ODE::HoloFlow::EqnDef} and \eqref{Eqn::ODE::HoloFlow::DefExpZ}.

Clearly $\Psi(0,0,0)=(\zeta_0,s_0)$. By Lemma \ref{Lem::ODE::HoloFlow} \ref{Item::ODE::HoloFlow::DefReg} and taking compositions, $\Psi$ is defined in a neighborhood $\Omega'\times\Omega''\times\Omega'''\subseteq\C^m_z\times\C^m_w\times\R^q_s$ of $0$, such that $\Psi\in\Co^\infty\Co^\beta(\Omega'\times\Omega'',\Omega''';\tilde U_0)$ and  $\Psi$ is holomorphic in $(z,w)$.

Since $\Coorvec{z^j}e^{z^jZ_j}=Z_j\circ e^{z^jZ_j}$ from \eqref{Eqn::ODE::HoloFlow::EqnDef}, by Lemma \ref{Lem::ODE::HoloFlow} \ref{Item::ODE::HoloFlow::Comm} we have
\begin{equation}\label{Eqn::ODE::PfParaHolFro::DPsi}
    \textstyle\Coorvec{z^j}\widehat\Psi(z,w,s)=Z_j\circ\widehat\Psi(z,w,s),\quad 1\le j\le m,\quad (z,w,s)\in\Omega.
\end{equation}

By condition \eqref{Eqn::ODE::PfParaHolFro::LinInd} we see that $\nabla_{z,w}\Psi(\cdot,0)$ is of full rank in a neighborhood of $0$. 
Therefore, by Lemma \ref{Lem::Hold::CompofMixHold} \ref{Item::Hold::CompofMixHold::InvFun} and shrinking $\Omega'\times\Omega''\times\Omega'''$ if necessary, $\Psi(\cdot,s):\Omega'\times\Omega''\to\tilde U_0$ is locally biholomorphic for each $s\in\Omega'''$ and $\widehat\Psi^\Inv\in\Co^{\infty,\beta}_\loc(\tilde U,V;\C^n\times\R^q)$ for some neighborhood $\tilde U\times V\subseteq\tilde U_0\times V_0$ of $(\zeta_0,s_0)$. Thus we get the well-defined map $\widetilde G(\zeta,s):=\Psi(\cdot,s-s_0)^\Inv(\zeta)$ for $(\zeta,s)\in\tilde U\times V$ with $\widetilde G\in\Co^{\infty,\beta}_{\zeta,s}$.

Write $\widetilde G=(\widetilde G',\widetilde G'',\widetilde G''')$ with respect to the $z$, $w$, $s$-coordinate components. We are going to show that $\widetilde G'':\tilde U\times V\to\C^{n-m}$ is as desired.

Immediately $\widetilde G''(\zeta_0,s_0)=0$ because $\widehat\Psi^\Inv(\zeta,s)=(\widetilde G''(\zeta,s),s-s_0)$ and $\widehat\Psi(0)=(\zeta_0,s_0)$. We have \ref{Item::ODE::ParaHolFro::0}.

Since $\widetilde G\in\Co^{\infty,\beta}_{\zeta,s}(\tilde U,V;\C^n)$ and $\Psi$ is holomorphic in $z,w$, we know $\widetilde G$ is holomorphic in $\zeta$. By Lemma \ref{Lem::Hold::NablaHarm} we get $\widetilde G\in\Co^\infty_\zeta\Co^\beta_s$. In particular $\widetilde G''\in\Co^\infty_\zeta\Co^\beta_s$ is holomorphic in $\zeta$, giving \ref{Item::ODE::ParaHolFro::Reg}.

\medskip To prove \ref{Item::ODE::ParaHolFro::Span} and \ref{Item::ODE::ParaHolFro::SpanR} we can fix a $s\in V$. Thus without loss of generality $V=\{\ast\}$ is a singleton, and we can assume $Z(\zeta,s)=Z(\zeta)$ and $\widetilde G(\zeta,s)=\widetilde G(\zeta)$ in the following argument.

\smallskip\noindent
\textit{Proof of \ref{Item::ODE::ParaHolFro::Span}}: Now $\widetilde G:\tilde U\to\C^m$ is a complex (holomorphic) coordinate chart, so $d\widetilde G''^1,\dots,d\widetilde G''^{n-m},d\bar \zeta^1,\dots,d\bar \zeta^n$ are linearly independent differentials on $\tilde U$.

By \eqref{Eqn::ODE::PfParaHolFro::DPsi} we have $Z_j=\Psi_*\Coorvec{z^j}=\widetilde G^*\Coorvec{z^j}$, thus $Z_j\widetilde G''^k=(\widetilde G^*\Coorvec{z^j})(\widetilde G^*dw^k)=0$ all vanish for $1\le j\le m$ and $1\le k\le n-m$.
This finishes the proof of \ref{Item::ODE::ParaHolFro::Span}.

\smallskip\noindent\textit{Proof of \ref{Item::ODE::ParaHolFro::SpanR}}: By \ref{Item::ODE::ParaHolFro::Span}, we know the $n\times (n-m)$ matrix $(\partial_{\zeta^j}\widetilde G''^k)_{\substack{1\le j\le n\\1\le k\le n-m}}$ has full rank $n-m$ at every point in $\tilde U$.
By \ref{Item::ODE::ParaHolFro::Reg} $\widetilde G''$ is holomorphic in $\zeta$, so $\partial_\xi \widetilde G''=\partial_\zeta \widetilde G''$. Therefore, the matrix-valued function $(\partial_{\xi^j}\widetilde G''^k)_{n\times(n-m)}:U\to\C^{n\times (n-m)}$ still has rank $n-m$ at every point. Since $d(\widetilde G''^k|_{U})=\sum_{j=1}^n\frac{\partial \widetilde G''^k}{\partial\xi^j}d\xi^j$, we conclude that $d(\widetilde G''|_{U})$ is a collection of $(n-m)$-differentials that are linearly independent on $U$.

Since $\partial_\xi \widetilde G''=\partial_\zeta \widetilde G''$ and $Z_j\widetilde G''=0$ from above, we know $X_j(\widetilde G''|_{U})=0$ as well, for $1\le j\le m$. Therefore, $d(\widetilde G''|_{U})$ annihilates $X_1,\dots,X_m$ on $U$. This completes the proof of \ref{Item::ODE::ParaHolFro::SpanR}.
\end{proof}

\section{Holomorphic Extension of An Inverse Laplacian}\label{Section::SecHolLap}

In this section, for a holomorphic function $f(z)$ we use $\nabla f=\partial_zf=(\partial_1f,\dots,\partial_nf)$ (cf. \eqref{Eqn::Intro::ColumnNote}).

We use $\Sp^{n-1}=\partial\B^{n-1}$ for unit sphere.
We use $\Hb^n:=\{x+iy\in\C^n:|x|<1,4|y|<1-|x|\}$ for a fixed complex cone. We use $\Oh(\Hb^n)=\Oh_\loc(\Hb^n)$ for the space of holomorphic functions in $\Hb^n$, endowed with  compact-open topology.

Our goal in this section is to prove the following:
\begin{prop}[Holormorphic property of Laplacian]\label{Prop::HolLap}Let $\B^n\subset\R^n$ and $\Hb^n\subset\C^n$ be as above.
	\begin{enumerate}[parsep=-0.3ex,label={(\roman*)}]
		\item\label{Item::HolLap::P} There is a linear operator $\Pv$ acting on functions on $\B^n$ such that $\Delta \Pv=\id_{\B^n}$ and $\Pv:\Co^{\alpha-2}(\B^n)\to\Co^\alpha(\B^n)$ is bounded for all $-2<\alpha<1$.
		\item\label{Item::HolLap::TildeP} There is a linear operator  $\tilde \Pv$ acting on functions on $\Hb^n$, such that for every $-2<\alpha<1$:
		\begin{itemize}[nolistsep]
		    \item $\tilde \Pv:\Co^{\alpha-2}_\Oh(\Hb^n)\to \Co^\alpha_\Oh(\Hb^n)$ is bounded linear;
		    \item $\Delta \tilde \Pv u=u$ and $(\tilde \Pv u)|_{\B^n}=\Pv(u|_{\B^n})$ for all $u\in\Co^{\alpha-2}_\Oh(\Hb^n)$. Here $\Delta=\Delta_x=\sum_{j=1}^n\frac{\partial^2}{\partial x^j\partial x^j}$.
		\end{itemize}
		\item\label{Item::HolLap::E} If $\Delta u=0$ in $\B^n$, then $u$ extend holomorphically into $\Hb^n$, call it $\Ex u\in\Oh(\Hb^n)$. Moreover the holomorphic extension operator $\Ex :\Co^{\alpha}(\B^n;\C)\cap\ker\Delta\to \Co^{\alpha}_\Oh(\Hb^n)$ is bounded for all $-2<\alpha<1$.
	\end{enumerate}
\end{prop}

The idea follows from Morrey \cite{Analyticity}, which only consider the cases for positive $\alpha$. However in our cases we need $\alpha$ to be negative. So we include a complete proof.
By taking direct modification to our construction (see \eqref{Eqn::SecHolLap::FormulaForP}) we can show that there exist an inverse Laplacian $\Pv$ and its extension $\tilde \Pv$ such that $\Pv:\Co^{\alpha-2}(\B^n)\to\Co^\alpha(\B^n)$ and $\tilde\Pv:\Co^{\alpha-2}_\Oh(\Hb^n)\to\Co^\alpha(\Hb^n)$ is bounded for all $-M<\alpha<M$, where $M$ is an arbitrarily fixed large number.
\begin{remark}
    In \ref{Item::HolLap::TildeP} since $u\in\Co^{\alpha-2}_\Oh$ is holomorphic, the statement $\Delta_x \tilde \Pv u=u$ is equivalent as $\Delta_z\tilde\Pv u=u$ where $\Delta_z=\sum_{j=1}^n\frac{\partial^2}{\partial z^j\partial z^j}$. In this case $\Delta_z$ is NOT an elliptic operator on $\Hb^n$: it is not the complex Laplacian $\square=\sum_{j=1}^n\frac{\partial^2}{\partial z^j\partial\bar z^j}$.
\end{remark}

\begin{defn}\label{Defn::SecHolLap::HoloHoldSpace}
Let $\alpha<1$. We denote by $\Co^\alpha_\Oh(\Hb^n)=\Co^\alpha_\Oh(\Hb^n;\C)$  the space of holomorphic function $g\in\Oh(\Hb^n)$ such that 
\begin{equation*}
    \|g\|_{\Co^\alpha_\Oh(\Hb^n)}:=|g(0)|+\sup\limits_{z\in\Hb^n}\dist(z,\partial\Hb^n)^{1-\alpha}|\nabla g(z)|<\infty.
\end{equation*}
\end{defn}

We have in fact $\Co^\alpha_\Oh(\Hb^n)=\Co^\alpha(\Hb^n;\C)\cap\Oh(\Hb^n)$, see \cite[Chapter 7]{ZhuSpaceBall} for some illustrations. For completeness we prove some partial results in Section \ref{Section::SecHolLap::HLLem}.

We construct the inverse Laplacian $\Pv$ in Definition \ref{Defn::SecHolLap::DefofP} and its holomorphic extension $\tilde\Pv$ in Definition \ref{Defn::SecHolLap::DefofTildeP}. We prove Proposition \ref{Prop::HolLap} \ref{Item::HolLap::E} in Section \ref{Section::SecHolLap::HLLem}; we prove  Proposition \ref{Prop::HolLap} \ref{Item::HolLap::P} in Section \ref{Section::SecHolLap::DefP}; we show $\tilde \Pv$ preserve holomorphic functions in Section \ref{Section::SecHolLap::DefTildeP} and prove  Proposition \ref{Prop::HolLap} \ref{Item::HolLap::TildeP} in Section \ref{Section::SecHolLap::BddTildeP}.

Recall in Notation \ref{Note::Hold::Newtonian} we use $\Ga$ as the Newtonian potential in $\R^n$:
\begin{equation}\label{Eqn::SecHolLap::NewtonPotent}
    \Ga(x)=\begin{cases}\frac{|x|}2&n=1,\\\frac1{2\pi}\log|x|&n=2,\\
	-\frac{|x|^{2-n}}{(2-n)|\Sp^{n-1}|}&n\ge3,
	\end{cases}\quad x\in\R^n\backslash\{0\}.
\end{equation}

Note that $\R^n\backslash\{0\}$ is contained in $\{z\in\C^n:\re(z\cdot z)>0\}=\{z\in\C^n:|\re z|>|\im z|\}$, which is a connected set in $\C^n$. So $[x\in\R^n\backslash\{0\}\mapsto|x|]$ has unique holomorphic extension to $\{|\re z|>|\im z|\}$ as 
\begin{equation}\label{Eqn::SecHolLap::Eqn|z|}
    z\in\big\{|\re z|>|\im z|\big\}\mapsto(z\cdot z)^\frac12,\quad\text{with principle argument }-\tfrac\pi2<2\arg z<\tfrac\pi 2.
\end{equation}

We still use $\Ga$ as the holomorphic extended fundamental solution on $\{|\re z|>|\im z|\}$. 

Therefore, for each $\theta\in\Sp^{n-1}=\partial\B^n$, the function $[x\in\B^n\mapsto\Ga(x-\theta)]$ can be holomorphic extended to $[z\in\Hb^n\mapsto\Ga(z-\theta)]$. Moreover, for every multindex $\nu\neq0$, since $|\partial^\nu\Ga(x)|\lesssim_\nu|x|^{2-n-|\nu|}$, we have the following:
\begin{align}\label{Eqn::SecHolLap::DevGa1}
    |(\partial^\nu\Ga)(z-t)|\lesssim_\nu|\re z-t|^{2-n-|\nu|}&,\quad\forall t\in\R^n,z\in\C^n\text{ satisfying }|\re z-t|\ge2|\im z|.
\end{align}




\subsection{Holomorphic H\"older spaces and extension of harmonic functions}\label{Section::SecHolLap::HLLem}
For holomorphic functions in a complex domain whose H\"older norms are finite, we can use the blow-up speed of its derivative near the boundary.

\begin{lem}[Hardy-Littlewood type characterizations]\label{Lem::SecHolLap::HLLem}
\begin{enumerate}[parsep=-0.3ex,label=(\roman*)]
    \item\label{Item::SecHolLap::HLLem::Grad} Let $\alpha<1$. The differentiation $[f\mapsto\nabla f]:\Co^\alpha_\Oh(\Hb^n)\to\Co^{\alpha-1}_\Oh(\Hb^n;\C^n)$ is bounded linear.
    \item\label{Item::SecHolLap::HLLem::0Char} Let $\alpha<0$. Then $\Co^\alpha_\Oh(\Hb^n)$ has equivalent norm 
    $f\mapsto\sup_{z\in\Hb^n}\dist(z,\partial\Hb^n)^{-\alpha}|f(z)|$.
    
    \item\label{Item::SecHolLap::HLLem::Res} The restriction map $[f\mapsto f|_{\B^n}]:\Co^\alpha_\Oh(\Hb^n)\to\Co^\alpha(\B^n;\C)$ and the inclusion map $\Co^\alpha_\Oh(\Hb^n)\hookrightarrow L^\infty(\Hb^n;\C)$ are both bounded linear for $0<\alpha<1$.
    \item\label{Item::SecHolLap::HLLem::Harm} Let $\alpha<1$. If $g\in\Co^\alpha(\B^n)$ satisfies $\Delta g=0$, then $\sup_{x\in\B^n}\dist(x,\Sp^{n-1})^{1-\alpha}|\nabla g(x)|<\infty$. Moreover there is a $C=C(n,\alpha)>0$ such that $$\sup_{x\in\B^n}\dist(x,\Sp^{n-1})^{1-\alpha}|\nabla g(x)|\le C_3\|g\|_{\Co^\alpha(\B^n)},\quad\forall g\in\Co^\alpha(\B^n)\cap\ker\Delta.$$
\end{enumerate}
\end{lem}

\begin{proof}
\ref{Item::SecHolLap::HLLem::Grad}:  By assumption $|\partial_jf(z)|\le\|f\|_{\Co^\alpha_\Oh}\dist(z,\partial\Hb^n)^{\alpha-1}$. By Cauchy's integral formula,
\begin{equation}\label{Eqn::SecHolLap::HLLem::CauInt}
    \partial_jg(z)=\frac1{2\pi  r}\int_0^{2\pi}e^{-i\theta}g(z+re^{i\theta}\mathbf e_j)d\theta,\quad 1\le j\le n,\quad 0<r<\dist(z,\partial\Hb^n),\quad g\in\Oh(\Hb^n).
\end{equation}

Taking $r=\frac12\dist(z,\partial\Hb^n)$ and $g=\partial_kf$, we get $|\partial_{jk}f(z)|\le \frac1r\sup_{|w-z|=r}|\partial_k f(w)|\le2^{2-\alpha}\|f\|_{\Co^\alpha_\Oh}\dist(z,\partial\Hb^n)^{\alpha-2}$, for $1\le j,k\le n$, which means $\sup_{z\in\Hb^n}\dist(z,\partial\Hb^n)^{2-\alpha}|\nabla^2 f(z)|\lesssim\|f\|_{\Co^\alpha_\Oh}$. This proves \ref{Item::SecHolLap::HLLem::Grad}.
 
\medskip\noindent\ref{Item::SecHolLap::HLLem::0Char}: Now $\alpha<0$. Clearly $|f(0)|\le\dist(0,\partial\Hb^n)^\alpha\cdot \sup_z\dist(z,\partial\Hb^n)^{-\alpha}|f(z)|$. If $\sup_z\dist(z,\partial\Hb^n)^{-\alpha}|f(z)|<\infty$, then by \eqref{Eqn::SecHolLap::HLLem::CauInt} with $g=f$ and $r=\frac12\dist(z,\partial\Hb^n)$ we get $|\nabla f(z)|\le \frac1r\sup_{|w-z|=r}|f(w)|\lesssim\dist(z,\partial\Hb^n)^{\alpha}$. So $\|f\|_{\Co^\alpha_\Oh}\lesssim \sup_z\dist(z,\partial\Hb^n)^{-\alpha}|f(z)|$.

Conversely since $\Hb^n$ is a convex set, we have $f(z)=f(0)+\int_0^1z\cdot \nabla f(tz)dt$ for all $f\in\Oh(\Hb^n)$ and $z\in\Hb^n$. Note that $\dist(tz,\partial\Hb^n)\ge t\dist(z,\partial\Hb^n)+(1-t)\dist(0,\partial\Hb^n)$ for all $0\le t\le 1$. Therefore for $f\in\Co^\alpha_\Oh(\Hb^n)$,
\begin{align*}
    |f(z)|&\le|f(0)|+|z|\|f\|_{\Co^\alpha_\Oh}\int_0^1\big(t\dist(z,\partial\Hb^n)+(1-t)\dist(0,\partial\Hb^n)\big)^{\alpha-1}dt
    \\
    &\lesssim\|f\|_{\Co^\alpha_\Oh}+\|f\|_{\Co^\alpha_\Oh}\int_0^{1-\dist(z,\partial\Hb^n)}(1-t)^{\alpha-1}dt\lesssim_\alpha\dist(z,\partial\Hb^n)^\alpha\|f\|_{\Co^\alpha_\Oh}.
\end{align*}

So $\sup_z\dist(z,\partial\Hb^n)^{-\alpha}|f(z)|\lesssim\|f\|_{\Co^\alpha_\Oh} $, completing the proof of \ref{Item::SecHolLap::HLLem::0Char}.

\medskip\noindent\ref{Item::SecHolLap::HLLem::Res}: Now $0<\alpha<1$. It suffices to show $\Co^\alpha_\Oh(\Hb^n)\subset\Co^\alpha(\Hb^n;\C)$. The boundedness of the inclusion and restriction would then follow from the natural maps $\Co^\alpha(\Hb^n)\hookrightarrow L^\infty(\Hb^n)$ and $\Co^\alpha(\Hb^n)\twoheadrightarrow\Co^\alpha(\B^n)$ respectively.

Now for $z_0,z_1\in\Hb^n$, what we need is to show that $|f(z_0)-f(z_1)|\lesssim_{n,\alpha}|z_0-z_1|^\alpha\|f\|_{\Co^\alpha_\Oh(\Hb^n)}$.

Since $\Hb^n$ is convex, by its geometry we can find a point $w\in\Hb^n$ (depending on $z_0,z_1$) such that $|w-z_j|\le|z_0-z_1|$ and $\dist(tw+(1-t)z_j,\partial\Hb^n)\ge\frac1{100} t|z_0-z_1|$ for $j=0,1$ and $0\le t\le 1$. Therefore
\begin{align*}
    |f(z_0)-f(z_1)|&\le |f(z_0)-f(w)|+|f(z_1)-f(w)|\le \sum_{j=0}^1\int_0^1|z_j-w||\nabla f(tw+(1-t)z_j)|dt
    \\&\lesssim|z_0-z_1|\|f\|_{\Co^\alpha_\Oh}\int_0^1(|z_0-z_1|t)^{\alpha-1}dt\lesssim|z_0-z_1|^\alpha\|f\|_{\Co^\alpha_\Oh}.
\end{align*}

Thus $\Co^\alpha_\Oh(\Hb^n)\subset\Co^\alpha(\Hb^n;\C)$ and we finish the proof of \ref{Item::SecHolLap::HLLem::Res}.

\medskip\noindent\ref{Item::SecHolLap::HLLem::Harm}: Let $\rho\in C_c^\infty(\B^n)$ be a radial function such that $\int\rho=1$. Let $\rho_j(x):=2^{jn}\rho(2^jx)$ for $j\ge0$. Therefore, for a harmonic function $g$ in $\B^n$, by the ball average formula we have $g(x)=\rho_k\ast g(x)$ whenever $1-|x|>2^{-k}$.

Let $(\phi_j)_{j=0}^\infty$ be a dyadic resolution (see Definition \ref{Defn::Hold::DyadicResolution}) for $\R^n$. Since $\nabla\rho$ is Schwartz for every $M>0$ we have $\|\phi_j\ast \nabla\rho\|_{L^1}\lesssim_{\phi,\rho,M}2^{-Mj}$ for $j\ge0$, so by scaling, we get $\|\phi_j\ast\nabla \rho_k\|_{L^1}=2^k\|\phi_{j-k}\ast\nabla\rho\|_{L^1}\lesssim_{\phi,\rho,M}2^{k-M(j-k)}$ if $j\ge k$. Therefore $\|\phi_j\ast\nabla\rho_k\|_{L^1}\lesssim_{\phi,\rho,M}2^k\min(1,2^{-M(j-k)})$ for all $j,k\ge0$.

Let $\psi_j:=\sum_{k=j-1}^{j+1}\phi_j$ for $j\ge0$ (we use $\phi_{-1}:=0$). By \eqref{Eqn::Hold::RmkDyaSupp} we have $\psi_j\ast \phi_j=\phi_j$. Therefore
\begin{equation*}
    \|\psi_j\ast\nabla\rho_k\|_{L^1}\lesssim_M 2^k\min(1,2^{-M(j-k)}),\quad j,k\ge0.
\end{equation*}

Now for a harmonic function $g\in\Co^\alpha(\B^n)$, consider an extension $\tilde g:=E_xg\in\Co^\alpha(\R^n)$ where $E_x$ is in Lemma \ref{Lem::Hold::CommuteExt}. By Lemma \ref{Lem::Hold::HoldChar} \ref{Item::Hold::HoldChar::LPHoldChar} (if $0<\alpha<1$) and Definition \ref{Defn::Hold::NegHold} (if $\alpha\le0$) we have $\|\phi_k\ast\tilde g\|_{L^\infty}\lesssim_{\phi,E_x}\|g\|_{\Co^\alpha(\B^n)}2^{-k\alpha}$.
Thus, for any $x\in\B^n$, take $k\ge0$ to be the unique integer such that $1-2^{1-k}\le |x|<1-2^{-k}$, we get
\begin{align*}
    &|\nabla g(x)|=|\nabla\rho_k\ast g(x)|\le\|\nabla\rho_k\ast\tilde g\|_{L^\infty}\le\sum_{j=0}^\infty\|\phi_j\ast\psi_j\ast\nabla\rho_k\ast\tilde g\|_{L^\infty}=\sum_{j=0}^\infty\|\psi_j\ast\nabla\rho_k\ast\phi_j\ast\tilde g\|_{L^\infty}
    \\
    \le&\sum_{j=0}^\infty\|\psi_j\ast\nabla\rho_k\|_{L^1}\|\phi_j\ast\tilde g\|_{L^\infty(\R^n)}\lesssim_{\alpha,\phi}\sum_{j=0}^\infty2^k\min(1,2^{-(|\alpha|+1)(j-k)})2^{-k\alpha}\|\tilde g\|_{\Co^\alpha(\R^n)}\lesssim_{\alpha}2^{k(1-\alpha)}\|g\|_{\Co^\alpha(\B^n)}.
\end{align*}
We conclude that $|\nabla g(x)|\lesssim(1-|x|)^{\alpha-1}\|g\|_{\Co^\alpha(\B^n)}$ for $g\in\Co^\alpha(\B^n)$ and $x\in\B^n$, finishing the proof of \ref{Item::SecHolLap::HLLem::Harm}.
\end{proof}



We can now prove Proposition \ref{Prop::HolLap} \ref{Item::HolLap::E}. Recall that for a real analytic function $f$, we use $\Ex f$ to be the (unique) holomorphic extension to the complex domain where it is defined.

\begin{proof}[Proof of Proposition \ref{Prop::HolLap} \ref{Item::HolLap::E}]Let $f\in\Co^\alpha(\B^n)$ be such that $\Delta f=0$. By Lemma \ref{Lem::SecHolLap::HLLem} \ref{Item::SecHolLap::HLLem::Harm} $f\in C^0_\loc(\B^n)$, and by Poisson's formula (See \cite[Chapter 2.2 (45)]{Evans} for example)
$$ f(x)=\frac{r^2-|x-x_0|^2}{|r\Sp^{n-1}|}\int_{\partial B(x_0,r)}\frac{f(\theta)d\sigma(\theta)}{|x-\theta|^n},\quad\text{for }x_0\in \B^n,\quad 0<r<1-|x_0|,\quad x\in B^n(x_0,r).$$
Here $d\sigma$ is the standard spherical measure. 

Using \eqref{Eqn::SecHolLap::Eqn|z|} we have holomorphic extension,
$$\Ex f(z)=\frac{r^2-(z-x_0)^2}{|r\Sp^{n-1}|}\int_{\partial B(x_0,r)}\frac{f(\theta)d\sigma(\theta)}{((z-\theta)^2)^\frac n2},\quad x_0\in\B^n,\quad r<1-|x_0|,\quad|\im z|<r-|\re z-x_0|.$$
For $z\in\Hb^n$, take $x_0=\re z$ and $r=\frac{1-|x_0|}2$, we can write
\begin{equation}\label{Eqn::SecHolLap::EqnEx}
    \Ex f(x+iy)=\frac{(\frac{1-|x|}2)^2+y^2}{|\frac{1-|x|}2\cdot\Sp^{n-1}|}\int_{\frac{1-|x|}2\Sp^{n-1}}\frac{f(x+\theta)d\sigma(\theta)}{((iy-\theta)^2)^\frac n2},\quad x+iy\in\Hb^n.
\end{equation}

Since $|y|<\frac14(1-|x|)=\frac12|\theta|$ for $\theta\in\frac{1-|x|}2\Sp^{n-1}$, we have $\re((iy-\theta)^2)\approx|\theta|^2\approx(1-|x|)^2$.

Note that $\nabla f$ is also harmonic, which gives $\nabla_z\Ex f=\Ex(\nabla_xf)$. Therefore by Lemma \ref{Lem::SecHolLap::HLLem} \ref{Item::SecHolLap::HLLem::Harm},
$$|\nabla\Ex f(x+iy)|\lesssim\frac{(1-|x|)^2}{1-|x|}\int_{\frac{1-|x|}2\Sp^{n-1}}\frac{|\nabla f(x+\theta)|d\sigma(\theta)}{(1-|x|)^n}\lesssim\sup\limits_{\partial B^n(x,\frac{1-|x|}2)}|\nabla f|\lesssim(1-|x|)^{\alpha-1}\|f\|_{\Co^\alpha(\B^n)}.$$

Since $\dist(z,\partial\Hb^n)\le\dist(x,\partial\Hb^n)\approx 1-|x|$, we get $\sup_{z\in\Hb^n}\dist(z,\partial\Hb^n)^{1-\alpha}|\nabla\Ex f(z)|\lesssim\|f\|_{\Co^\alpha(\B^n)}$. 

Clearly $|\Ex f(0)|\lesssim\|f\|_{\Co^\alpha(\B^n)}$, so by definition of $\Co^\alpha_\Oh$-spaces we get $\|\Ex f\|_{\Co^\alpha_\Oh(\Hb^n)}\lesssim\|f\|_{\Co^\alpha(\B^n)}$, finishing the proof.
\end{proof}

\subsection{Construction and H\"older regularity for real inverse Laplacian $\Pv$}\label{Section::SecHolLap::DefP}

We are about to define $\Pv f=(\Ga\ast Ef)|_{\B^n}$ where $E$ is a concrete extension operator for $\B^n$ such that $E:\Co^\alpha(\B^n)\to\Co^\alpha(\R^n)$ is bounded linear for $-4<\alpha<1$. 
\begin{lem}[{\cite[Theorem 2.9.2]{Triebel1}}]\label{Lem::SecHolLap::ExtLem}
Let $M\in\Z_+$, let  $(a_j,b_j)_{j=1}^M\subset\R$ satisfy $b_j>0$, $1\le j\le M$ and $\sum_{j=1}^Ma_j(-b_j)^k=1$ for $-4\le k\le 4$. Let $(\tilde\chi_j)_{j=1}^M\subset C_c^\infty(\R)$ satisfies $\chi_j(t)\equiv1$ for $t$ closed to 0 and each $1\le j\le M$. Define
\begin{equation}\label{Eqn::HolLap::ExtHPOp}
    E^Hf(x',x_n):=\begin{cases}f(x',x_n)&x_n>0,\\\sum_{j=1}^M\tilde\chi_j(x_n)a_j\cdot f(x',-b_jx_n)&x_n<0.\end{cases}
\end{equation}
Then $E^H$ defines an extension operator on the half plane $\R^n_+=\{x_n>0\}$ such that $E^H:\Co^\alpha(\R^n_+)\to\Co^\alpha(\R^n)$ for all $-4<\alpha<1$.
\end{lem}

\begin{cor}
Let $(a_j,b_j)_{j=1}^9\subset\R$ satisfy $b_j>0$, $1\le j\le 9$ and $\sum_{j=1}^9a_j(-b_j)^{-k}=1$ for $-4\le k\le 4$. For $1\le j\le 9$, let $\tilde\chi_j\in C_c^\infty(\R)$ be such that $\tilde\chi_j\equiv1$ in an neighborhood of $0$.
Then 
\begin{equation}\label{Eqn::SecHolLap::ExtBallOp}
    Ef(e^{-\rho}\theta):=\begin{cases}f(e^{-\rho}\theta)&\rho>0,\\\sum_{j=1}^9a_j\tilde\chi_j(\rho)f(e^{\rho/b_j}\theta)&\rho<0,\end{cases}\quad\rho\in\R,\quad\theta\in\Sp^{n-1},
\end{equation}
defines an extension operator on $\B^n$ such that $E:\Co^\alpha(\B^n)\to\Co^\alpha(\R^n)$ is bounded linear for $-4<\alpha<1$.

Moreover its formal adjoint $E^*$ is given by
\begin{equation}\label{Eqn::SecHolLap::EqnE*}
    E^*g(x)=g(x)+\sum_{j=1}^9\frac{a_jb_j}{|x|^{n(b_j+1)}}\tilde\chi_j(b_j\log|x|)g\Big(\frac x{|x|^{b_j+1}}\Big),\quad x\in\B^n.
\end{equation}

\end{cor}
\begin{proof}
Clearly $(\rho,\theta)\mapsto e^{-\rho}\theta$ is a diffeomorphism from $\R\times\Sp^{n-1}$ to $\R^n\backslash\{0\}$ and maps $\{\rho>0\}$ onto $\B^n\backslash\{0\}$. By \cite[Theorem 2.10.2(i)]{Triebel1} we see that $\Co^\alpha=\Bs_{\infty\infty}^\alpha$ is preserved by bounded diffeomorphism. Therefore, by passing to coordinate covers of $\Sp^{n-1}$ if necessary (since we only need to worry about the place where $\rho$ is closed to $0\in\R$), by Lemma \ref{Lem::SecHolLap::ExtLem} with $b_j$ replacing by $\frac1{b_j}$, we know \eqref{Eqn::SecHolLap::ExtBallOp} is an extension operator for $\B^n$ and has $\Co^\alpha$-boundedness for $-4<\alpha<1$.

To get \eqref{Eqn::SecHolLap::EqnE*}, let $g:\R^n\to\R$, we have
\begin{align*}
    &\int_{\R^n}Ef(y)g(y)dt=\int_{\Sp^{n-1}} d\theta\int_\R Ef(e^{-\rho}\theta)g(e^{-\rho}\theta)e^{-n\rho}d\rho
   \\=&\int_{\Sp^{n-1}}\bigg(\int_{\R_+} f(e^{-\rho}\theta)g(e^{-\rho}\theta)e^{-n\rho}d\rho+\sum_{j=1}^9a_j\tilde\chi_j(\rho)\int_{\R_-}f(e^{\rho/b_j}\theta)g(e^{-\rho}\theta)e^{-n\rho}d\rho\bigg) d\theta
   \\=&\int_{\Sp^{n-1}} d\theta\int_{\R_+} f(e^{-r}\theta)g(e^{-r}\theta)e^{-nr}dr+\sum_{j=1}^9a_jb_j\int_{\Sp^{n-1}} d\theta\int_{\R_+}\tilde\chi_j(b_jr)f(e^{-r}\theta)g(e^{b_jr}\theta)e^{nr(b_j+1)}e^{-nr}dr
   \\=&\int_{\B^n} f(x)\bigg(g(x)+\sum_{j=1}^9\frac{a_jb_j\cdot\tilde\chi_j(b_j\log|x|)}{|x|^{n(b_j+1)}}g\Big(\frac x{|x|^{b_j+1}}\Big)\bigg)dx.
\end{align*}



This gives \eqref{Eqn::SecHolLap::EqnE*} and finishes the proof.
\end{proof}

We now give an concrete formula of $\Pv$ based on \eqref{Eqn::SecHolLap::EqnE*},

\begin{defn}\label{Defn::SecHolLap::DefofP} We define $\Pv$ to be
\begin{equation}\label{Eqn::SecHolLap::FormulaForP}
    \Pv f(x):=\int_{\B^n} f(t)\bigg(\Ga(x-t)+\chi(|t|)\sum_{j=1}^9\frac{a_jb_j}{|t|^{n(b_j+1)}}\Ga\Big(x-\frac{t}{|t|^{b_j+1}}\Big)\bigg)dt.
    \end{equation}
Here $\chi\in C_c^\infty(\frac13,\frac53)$ satisfies $\chi|_{[\frac12,\frac32]}\equiv1$, and $a_j\in\R,b_j>0$ satisfy $\sum_{j=1}^9a_j(-b_j)^k=1$ for $-4\le k\le 4$.
\end{defn}

\begin{proof}[Proof of Proposition \ref{Prop::HolLap} \ref{Item::HolLap::P}]
By \eqref{Eqn::SecHolLap::EqnE*} we see that $\Pv f$ is of the form $\Pv f(x)=\int_{\B^n}f(t)E^*(\Ga(x-\cdot))(t)dt$ for $f\in L^1(\B^n)$, where $E^*$ is the adjoint of some extension operator $E$ which is $\Co^\alpha$-bounded for $-4<\alpha<1$. Therefore $\Pv f(x)=(\Ga\ast Ef)(x)$ for $x\in\B^n$ and we get
$\Pv f=(\Ga\ast Ef)|_{\B^n}$. In particular $\Delta\Pv f=(\Delta\Ga\ast Ef)|_{\B^n}=f$ holds.

By construction, such $E$ has compact support, namely there is a bounded open set $U\subset\R^n$ (in fact we can take $U=B^n\big(0,e^{\max_j\frac1{3b_j}}\big)$) such that $\supp Ef\subset U$ has compact support for every $f\in\Co^{(-4)+}(\B^n)$. Therefore by Lemma \ref{Lem::Hold::GreensOp}, $\Pv=[f\mapsto(\Ga\ast Ef)|_{\B^n}]:\Co^\alpha(\B^n)\to\Co^{\alpha+2}(\B^n)$ is bounded linear for all $-4<\alpha<1$, in particular for $-4<\alpha<-1$.
\end{proof}



To prove Proposition \ref{Prop::HolLap} \ref{Item::HolLap::TildeP}, we need to decompose the integral in \eqref{Eqn::SecHolLap::FormulaForP} into the domains $\{t:|\re z-t|<\frac{1-|\re z|}2\}$ and $\{t:|\re z-t|>\frac{1-|\re z|}2\}$. We need the following estimate.
\begin{lem}\label{Lem::SecHolLap::AprioriBddforP}Let $a,b,\chi$ be given in Definition \ref{Defn::SecHolLap::DefofP}. There is a $C=C(n,a,b,\chi)>0$, such that for every $z\in\Hb^n$ and $t\in\B^n$ satisfying $|\re z-t|\ge\frac{1-|\re z|}2$, and for every multi-index $\nu$ such that $|\nu|\le1$,
\begin{equation}\label{Eqn::SecHolLap::AprioriBddforPEqn}
    \bigg|(\partial^\nu\Ga)(z-t)+\chi(|t|)\sum_{j=1}^9\frac{a_jb_j}{|t|^{n(b_j+1)}}(\partial^\nu\Ga)\Big(z-\frac{t}{|t|^{b_j+1}}\Big)\bigg|\le C|\re z-t|^{-n-2-|\nu|}(1-|t|)^4.
\end{equation}
\end{lem}
\begin{proof}
That $|\re z-t|\ge\frac{1-|\re z|}2$ implies $|\re z-t|\ge2|\im z|$, so by \eqref{Eqn::SecHolLap::DevGa1} the left hand side of \eqref{Eqn::SecHolLap::AprioriBddforPEqn} is finite.

When $1-|t|\ge\frac12|\re z-t|$, we have $|\re z-\frac t{|t|^{b_j+1}}|\ge \frac 1{b_j}|\re z-t|$. Thus by \eqref{Eqn::SecHolLap::DevGa1},
$$\textstyle\text{LHS of }\eqref{Eqn::SecHolLap::AprioriBddforPEqn}\lesssim|\re z-t|^{1-n}+|\re z-\frac t{|t|^{b_j+1}}|^{1-n}\lesssim_{n,b} |\re z-t|^{1-n}\lesssim|\re z-t|^{-3-n}(1-|t|)^4.$$

We then focus on the case $1-|t|<\frac12|\re z-t|$.

By assumption $\sum_{j=1}^9a_jb_j(-b_j)^k=-1$ for $0\le k\le3$. Therefore, by the Taylor's expansion of degree 3,  for any $\tilde \chi\in C_c^\infty(\R)$ that equals 1 in the neighborhood of 0, there is a $C'=C_{a,b,\tilde\chi}'>0$ such that
\begin{equation}\label{Eqn::SecHolLap::AprPEqn2}
    \textstyle\big|\psi(\rho)+\tilde\chi(\rho)\sum_{j=1}^9a_jb_j\psi(-b_j\rho)\big|\le C'\|\psi\|_{C^4}\cdot\rho^4,\quad\forall\ R>0,\  \psi\in C^4[-R\max_{1\le j\le 9}b_j,R],\  \rho\in[0,R].
\end{equation}

For fixed $t$ and $z$, we define $I_t:=\big[\log|t|\max\limits_{1\le j\le 9}b_j,-\log|t|\big]$ and $\psi_{z,t}^\nu(\rho):=e^{-n\rho}\cdot(\partial^\nu\Ga)(z-e^{-\rho}\frac t{|t|})$ for $\rho\in I_t$. Since for $\rho\in I_t$ we have $|\re z-e^{-\rho}\frac t{|t|}|\ge|\re z-t|$ and by assumption $|\re z-t|\ge2|\im z|$ holds in this case, using \eqref{Eqn::SecHolLap::DevGa1} we have $\|\psi_{z,t}^\nu\|_{C^4(I_t)}\lesssim|\re z-t|^{-n-2-|\nu|}$.

Thus taking $\tilde\chi(\rho):=\chi(e^{-\rho})$ and $\rho=-\log|t|(=R)$ in \eqref{Eqn::SecHolLap::AprPEqn2} with $\psi=\psi_{z,t}^\nu$ we get
$$\bigg||t|^n(\partial^\nu\Ga)(z-t)+\chi(|t|)\sum_{j=1}^9\frac{a_jb_j}{|t|^{nb_j}}(\partial^\nu\Ga)\Big(z-\frac{t}{|t|^{b_j+1}}\Big)\bigg|\lesssim |\re z-t|^{-n-2-|\nu|}(1-|t|)^4.$$

Note that $1-|t|<\frac12|\re z-t|$ implies $\frac13<|t|<1$, so by multiplying $\frac1{|t|^n}$ we get \eqref{Eqn::SecHolLap::AprioriBddforPEqn} for the case $1-|t|<\frac12|\re z-t|$. This finishes the proof.
\end{proof}

\begin{cor}\label{Cor::SecHolLap::IntofP}
If $(1-|t|)^4f(t)\in  L^1(\B^n)$, then the following integral converges locally uniformly for  $z\in\Hb^n$: 
\begin{equation}\label{Eqn::SecHolLap::IntofP::Tmp}
    \int_{\{t\in\B^n:|\re z-t|\ge\frac{1-|\re z|}2\}}f(t)\Big(\Ga(z-t)-\chi(|t|)\sum_{j=1}^9\frac{a_jb_j}{|t|^{n(b_j+1)}}\Ga\big(z-\frac{t}{|t|^{b_j+1}}\big)\Big)dt.
\end{equation}
\end{cor}
\begin{proof}By Lemma \ref{Lem::SecHolLap::AprioriBddforP} $\left|\Ga(z-t)-\sum_{j=1}^9\frac{a_jb_j}{|t|^{n(b_j+1)}}\Ga\big(z-\frac{t}{|t|^{b_j+1}}\big)\right|\le C|\re z-t|^{-n-2} (1-|t|)^4$ holds for $t\in\B^n\cap B^n(\re z,\frac{1-|\re z|}2)$, where $C=C_{n,a,b,\chi}$. Therefore \eqref{Eqn::SecHolLap::IntofP::Tmp} is bounded by
\begin{equation*}
    \|t\mapsto(1-|t|)^4f(t)\|_{L^1(\B^n)}\sup_{t\in\B^n:2|\re z-t|\ge1-|\re z|}C|\re z-t|^{-n-2}\le C(1-|\re z|)^{-n-2}.
\end{equation*}
The term $(1-|\re z|)^{-n-2}$ is locally bounded for $z\in\Hb^n$, thus we complete the proof.
\end{proof}

\subsection{Construction of the holomorphic extension $\tilde{\mathbf P}$}\label{Section::SecHolLap::DefTildeP}

We now extend $\Pv$ to the complex cone $\Hb^n$. 

In Sections \ref{Section::SecHolLap::DefTildeP} and \ref{Section::SecHolLap::BddTildeP}, we use a function\footnote{In \cite{Analyticity} he use $1-\lambda$ of our $\lambda$.} $\lambda:\B^n\times\B^n\to[0,1]$ as the following
\begin{equation}\label{Eqn::SecHolLap::DefLambda}
    \lambda(t,x):=\begin{cases}1-\frac{2|t-x|}{1-|x|},&\text{if }2|t-x|\le1-|x|\\0,&\text{if }2|t-x|\ge1-|x|\end{cases}.
\end{equation}
Immediately we have:
\begin{enumerate}[parsep=-0.3ex,label=($\Lambda$.\arabic*)]
    \item\label{Item::SecHolLap::LambdaPro1} $\lambda$ is locally Lipschitz, and the function $[(t,x+iy)\mapsto\lambda(t,x)\cdot y]:\B^n\times\Hb^n\to\R^n$ is bounded Lipschitz.
    \item\label{Item::SecHolLap::LambdaPro2} For every $t\in\B^n$ and $x+iy\in\Hb^n$, we have $4|\lambda(t,x)\cdot y|<1-|x|$ and $|(1-\lambda(t,x))y|\le\frac12|t-x|$.
\end{enumerate}

For each $z=x+iy\in\Hb^n$, let 
\begin{equation}\label{Eqn::SecHolLap::NotSzV}
    \textstyle S_z:=\{t+i\lambda(t,x) y:t\in\B^n\},\quad V(t,z):=1+iy\cdot\partial_t\lambda(t,x)=1+i\sum_{j=1}^ny^j(\partial_{t^j}\lambda)(t,x).
\end{equation}

By \ref{Item::SecHolLap::LambdaPro2} we see that $S_z\subset\Hb^n$, and $z-\zeta\in\{w\in\C^n:|\re w|>|\im w|\}$ for each $\zeta\in S_z$, so $\Ga(z-\zeta)$ is defined for $\zeta\in S_z$.

Fix $z\in\Hb^n$, by viewing $S_z$ as a set parameterized by $t$, we have a change of variable $\zeta=\zeta(t,z):=t+i\lambda(t,x)y$. Its Jacobian matrix is $\frac{\partial\zeta}{\partial t}(t,z)=1+iy\otimes\partial_t\lambda(t,x)$, so $\det\frac{\partial\zeta}{\partial t}(t,z)=1+iy\cdot\partial_t\lambda(t,x)$. We denote the ``volume form'' $V(t,z)$ as
\begin{equation*}
    V(t,x+iy):=1+iy\cdot\partial_t\lambda(t,x),\quad t\in\B^n,\quad x+iy\in\Hb^n.
\end{equation*}

In this way for an integrable function $g$ on $S_z$, we can formally write
$$\textstyle\int_{S_z}g(\zeta)d\zeta=\int_{\B^n}g(t+i\lambda(t,x)y)(1+iy\cdot \partial_t\lambda(t,x))dt.$$

\begin{defn}\label{Defn::SecHolLap::DefofTildeP}
Let $(a_j)_{j=1}^9,(b_j)_{j=1}^9,\chi$ to be as in Definition \ref{Defn::SecHolLap::DefofP}. Let $S_z$ be as in \eqref{Eqn::SecHolLap::NotSzV} for $z\in\Hb^n$. Define 
\begin{equation}\label{Eqn::SecHolLap::EqnHoloP0}
    \tilde\Pv f(z):=\int_{S_z}f(\zeta)\Ga(z-\zeta)d\zeta+\int_{\B^n}\chi(|t|)\sum_{j=1}^9\frac{a_jb_j}{|t|^{n(b_j+1)}}f(t)\Ga\Big(z-\frac{t}{|t|^{b_j+1}}\Big)dt,\quad\text{for }f\in\Co^{(-4)+}_\Oh(\Hb^n),
\end{equation}
in the sense that,
\begin{align}
    \tilde\Pv f(z)=\tilde\Pv f(x+iy)=
    &
    \label{Eqn::SecHolLap::EqnHoloP1}
    \int_{|t-x|<\frac{1-|x|}2}f(t+i\lambda(t,x)y)\Ga(x-t+i(1-\lambda(t,x))y)V(t,z)dt
    \\
    &\label{Eqn::SecHolLap::EqnHoloP2}
    +\int_{|t-x|<\frac{1-|x|}2}\sum_{j=1}^9\frac{a_jb_j\chi(|t|)}{|t|^{n(b_j+1)}}f(t)\Ga\Big(z-\frac{t}{|t|^{b_j+1}}\Big)dt\\
    &\label{Eqn::SecHolLap::EqnHoloP3}
    +\int_{t\in\B^n:|t-x|>\frac{1-|x|}2}f(t)\bigg(\Ga(z-t)-\sum_{j=1}^9\frac{a_jb_j\chi(|t|)}{|t|^{n(b_j+1)}}\Ga\Big(z-\frac{t}{|t|^{b_j+1}}\Big)\bigg)dt
    \\
    =:&\tilde\Pv_1 f(z)+\tilde\Pv_2 f(z)+\tilde\Pv_3 f(z)\label{Eqn::SecHolLap::EqnHoloPall}
\end{align}

\end{defn}

By \ref{Item::SecHolLap::LambdaPro2}, we see that \eqref{Eqn::SecHolLap::EqnHoloP1} is a classical Lebesgue integral pointwisely for $z\in\Hb^n$ and all $f\in\Oh(\Hb^n)$. The same holds for \eqref{Eqn::SecHolLap::EqnHoloP2}.

When $f\in\Co^{(-4)+}_\Oh(\Hb^n)$, by Lemma \ref{Lem::SecHolLap::HLLem} \ref{Item::SecHolLap::HLLem::0Char} we have $[t\mapsto (1-|t|)^4f(t)]\in L^\infty(\B^n;\C)$. By Corollary \ref{Cor::SecHolLap::IntofP},  \eqref{Eqn::SecHolLap::EqnHoloP3} is integrable as well.

\medskip
Note that $S_{x+i0}\equiv\B^n$. So comparing the expressions \eqref{Eqn::SecHolLap::FormulaForP} and \eqref{Eqn::SecHolLap::EqnHoloP0} we see that $\tilde \Pv f(x+i0)=\Pv (f(\cdot+i0))(x)$. Therefore to prove that $\tilde \Pv$ is the holomorphic extension of $\Pv$, what we need is the following:

\begin{prop}\label{Prop::SecHolLap::TildePisHolo}
Let $f\in\Co^{(-4)+}_\Oh(\Hb^n)$, then $\tilde \Pv f$ gives a holomorphic function defined on $\Hb^n$.

Moreover we have $(\tilde \Pv f)|_{\B^n}=\Pv[f|_{\B^n}]$ and $\Delta\tilde\Pv f=f$, where $\Delta=\sum_{j=1}^n\frac{\partial^2}{(\partial z^j)^2}$.
\end{prop}
We postpone its proof after Proposition \ref{Lem::SecHolLap::PfTildePHolo}. Essentially we can assume $f\in L^\infty(\Hb^n;\C)$ so that the two integrals in \eqref{Eqn::SecHolLap::EqnHoloP0} converge individually:
\begin{lem}\label{Lem::SecHolLap::PfTildePHolo}
    For $f\in L^\infty(\Hb^n;\C)\cap\Oh(\Hb^n)$, the function $z\in\Hb^n\mapsto\int_{S_z}f(\zeta)\Ga(z-\zeta)d\zeta$ is holomorphic. Moreover
    \begin{equation}
        \Coorvec{z^j}\int_{S_z}f(\zeta)\Ga(z-\zeta)d\zeta=\int_{\B^n}f(t+i\lambda(t,x)y)\Ga_{,j}(x-t+iy(1-\lambda(t,x)))V(t,z)dt,\quad 1\le j\le n,\quad z\in\Hb^n.
    \end{equation}
\end{lem}Here we use the comma notation $\Ga_{,j}(\zeta)=(\partial_{z^j}\Ga)(\zeta)$ for $\zeta$ in the domain.
\begin{proof} Using the comma notation we can write $(\partial_{z^j}f)(\zeta)=f_{,x^j}(\zeta)=-if_{,y^j}(\zeta)=f_{,j}(\zeta)$.

In the following computations, we use abbreviations $f=f(t+i\lambda(t,x)y)$, $\Ga=\Ga(x-t+i(1-\lambda(t,x))y)$ and $V=V(t,z)=1+iy\cdot\lambda_{,t}(t,x)$. First we claim to have the following (cf. \cite[(4.12)]{Analyticity}):
\begin{gather}\label{Eqn::SecHolLap::PfTildePHolo::Tmp1}
    \sum_{k=1}^n\Coorvec{t^k}(y^kf\Ga\nabla\lambda)=\sum_{k=1}^n(y^k(f_{,k}\Ga-f\Ga_{,k})V\nabla\lambda+y^kf\Ga(\nabla\lambda)_{,t^k}).
    \\\label{Eqn::SecHolLap::PfTildePHolo::Tmp2}
    \Coorvec{t^j}(\lambda f\Ga V)-i\sum_{k=1}^n\Coorvec{t^k}(y^k\lambda\lambda_{,t^j}f\Ga)=f\Ga\lambda_{,t^j}+\lambda(f_{,j}\Ga-f\Ga_{,j})V\quad\text{for }j=1,\dots,n.
\end{gather}
Here $\nabla\lambda=(\lambda_{,x^1},\dots,\lambda_{,x^n},\lambda_{,t^1},\dots,\lambda_{,t^n})$.

Indeed by the chain rule, \eqref{Eqn::SecHolLap::PfTildePHolo::Tmp1} is done by the following,
    \begin{align*}
    &\textstyle\sum_{k=1}^n\Coorvec{t^k}(y^kf\Ga\nabla\lambda)
    \\
    =&\textstyle\sum_{k=1}^ny^kf_{,k}\Ga\nabla\lambda+\sum_{k,l=1}^niy^kf_{,l}y^l\lambda_{,t^k}\Ga\nabla\lambda-\sum_{k=1}^ny^kf\Ga_{,k}\nabla\lambda+\sum_{k,l=1}^niy^k\Ga_{,l}y^l\lambda_{,t^k}\nabla\lambda+\sum_{k=1}^ny^kf\Ga\nabla\lambda_{,t^k}
    \\
    =&\textstyle\sum_{k,l=1}^ny^k(1+iy^l\lambda_{,t^l})(f_{,k}\Ga-f\Ga_{,k})\nabla\lambda+\sum_{k=1}^ny^kf\Ga\nabla\lambda_{,t^k}
    \\
    =&\textstyle\sum_{k=1}^ny^k(f_{,k}\Ga-f\Ga_{,k})V\nabla\lambda+\sum_{k=1}^ny^kf\Ga\nabla\lambda_{,t^k}.
    \end{align*}
    
Since $\partial_t\lambda$ is the component of $\nabla\lambda$, using \eqref{Eqn::SecHolLap::PfTildePHolo::Tmp1} we get \eqref{Eqn::SecHolLap::PfTildePHolo::Tmp2} by the following
\begin{align*}
    &\textstyle\Coorvec{t^j}(\lambda f\Ga V)-i\sum_{k=1}^n\Coorvec{t^k}(y^k\lambda\lambda_{,t^j}f\Ga)
    \\
    =&\textstyle\lambda_{,t^j}f\Ga V+\lambda(f_{,j}\Ga-f\Ga_{,j})V+i\sum_{k=1}^ny^k\lambda\lambda_{,t^j}(f_{,k}\Ga-f\Ga_{,k})V+i\lambda f\Ga\sum_{k=1}^ny^k\lambda_{,t^jt^k}
    \\
    &\textstyle-i\sum_{k=1}^ny^k\lambda_{,t^k}\lambda_{,t^j}f\Ga-i\lambda\sum_{k=1}^n\Coorvec{t^k}(y^kf\Ga\lambda_{,t^j})
    \\
    =&\textstyle\lambda_{,t^j}f\Ga \big(1+i\sum_{k=1}^ny^k\lambda_{,t^k}\big)+\lambda(f_{,j}\Ga-f\Ga_{,j})V-i\sum_{k=1}^ny^k\lambda_{,t^k}\lambda_{,t^j}f\Ga
    \\
    =&f\Ga\lambda_{,t^j}+\lambda(f_{,j}\Ga-f\Ga_{,j})V.
\end{align*}

Next we compute $\Coorvec{x^j}\int_{\B^n}f\Ga Vdt$ and $\Coorvec{y^j}\int_{\B^n}f\Ga Vdt$. Note that $\lambda=0$ and $\nabla\lambda=0$ near $t\in\Sp^{n-1}$, so both \eqref{Eqn::SecHolLap::PfTildePHolo::Tmp1} and \eqref{Eqn::SecHolLap::PfTildePHolo::Tmp2} have integrals zero on $t\in\B^n$. Therefore,
\begin{align*}
    \textstyle\Coorvec{x^j}\int_{\B^n}f\Ga Vdt=&\textstyle\int_{\B^n}f\Ga_{,j}V+\sum_{k=1}^n\big(if_{,k}y^k\lambda_{,x^j}\Ga V-if\Ga_{,k}y^k\lambda_{,x^j}V+if\Ga y^k\lambda_{,t^kx^k}\big)dt
    \\
    =&\textstyle\int_{\B^n}f\Ga_{,j}Vdt+i\sum_{k=1}^n\int_{\B^n}\big(y^k(f_{,k}\Ga-f\Ga_{,k})V\lambda_{,x^j}+y^kf\Ga\lambda_{,x^jt^k}\big)dt
    \\
    =&\textstyle\int_{\B^n}f\Ga_{,j}Vdt+i\sum_{k=1}^n\int_{\B^n}\Coorvec{t^k}(y^kf\Ga\lambda_{,x^j})dt&(\text{by }\eqref{Eqn::SecHolLap::PfTildePHolo::Tmp1})
    \\
    =&\textstyle\int_{\B^n}f\Ga_{,j}Vdt+i\sum_{k=1}^n\int_{\Sp^{n-1}}y^kf\Ga\lambda_{,x^j}d\sigma_k(t)
    \\
    =&\textstyle\int_{\B^n}f\Ga_{,j}Vdt.&(\nabla\lambda|_{\Sp^{n-1}}=0)
\end{align*}
\begin{align*}
    \textstyle\Coorvec{y^j}\int_{\B^n}f\Ga Vdt=&\textstyle i\int_{\B^n}\big(f_{,j}\lambda\Ga V+f\Ga_{,j}(1-\lambda)V+f\Ga\lambda_{,t^j}\big)dt&(f_{,j}=-if_{,y^j})
    \\
    =&\textstyle i\int_{\B^n}f\Ga_{,j}Vdt+i\int_{\B^n}\big(f\Ga\lambda_{,t^j}+\lambda(f_{,j}\Ga-f\Ga_{,j})V\big)dt
    \\
    =&\textstyle i\int_{\B^n}f\Ga_{,j}Vdt+i\int_{\B^n}\Coorvec{t^j}(\lambda f\Ga V)dt+\sum_{k=1}^n\int_{\B^n}\Coorvec{t^k}(y^k\lambda\lambda_{,t^j}f\Ga)dt&(\text{by }\eqref{Eqn::SecHolLap::PfTildePHolo::Tmp2})
    \\
    =&\textstyle i\int_{\B^n}f\Ga_{,j}Vdt+i\int_{\Sp^{n-1}}\lambda f\Ga Vd\sigma_j(t)+\sum_{k=1}^n\int_{\Sp^{n-1}}y^k\lambda\lambda_{,t^j}f\Ga d\sigma_k(t)
    \\=&\textstyle i\int_{\B^n}f\Ga_{,j}Vdt.&(\lambda|_{\Sp^{n-1}}=0)
\end{align*}
Here we use $d\sigma_j$ as the $j$-th outer normal direction of $\Sp^{n-1}=\partial\B^n$.

Therefore $\Coorvec{x^j}\int_{\B^n}f\Ga Vdt=i\Coorvec{y^j}\int_{\B^n}f\Ga Vdt=\int_{\B^n}f\Ga_{,j}Vdt$, finishing the proof.
\end{proof}

\begin{proof}[Proof of Proposition \ref{Prop::SecHolLap::TildePisHolo}]
When $f\in L^\infty(\Hb^n;\C)$, two integrals in \eqref{Eqn::SecHolLap::EqnHoloP0} converge individually. 

By Proposition \ref{Lem::SecHolLap::PfTildePHolo} the first integral $z\mapsto\int_{S_z}f(\zeta)\Ga(z-\zeta)d\zeta$ is holomorphic in $\Hb^n$. Clearly the second integral in \eqref{Eqn::SecHolLap::EqnHoloP0} is holomorphic since $z\mapsto\Ga(z-t/|t|^{b_j+1})$ is a holomorphic in $\Hb^n$ as well. So $\tilde \Pv f\in\Oh(\Hb^n)$ when $f\in L^\infty(\Hb^n;\C)\cap\Oh(\Hb^n)$.

For general $f\in\Co^{(-4)+}_\Oh(\Hb^n)$, by Lemma \ref{Lem::SecHolLap::HLLem} \ref{Item::SecHolLap::HLLem::0Char} we have $f|_{\B^n}\in L^1(\B^n,(1-|x|)^4dx;\C)$. Take $f_\eps(z):=f((1-\eps)z)$ for $\eps>0$, we have $f_\eps\in L^\infty(\Hb^n;\C)$ for each $\eps$, and $f_\eps|_{\B^n}\xrightarrow{\eps\to0}f|_{\B^n}$ converges in $L^1(\B^n,(1-|x|)^4dx;\C)$. Therefore in the notation \eqref{Eqn::SecHolLap::EqnHoloPall} we have $\lim_{\eps\to0}\tilde\Pv_3f_\eps(z)=\tilde \Pv_3f(z)$ locally uniformly in $z$.

Since the convergence $f_\eps\xrightarrow{\eps\to0}f$ holds automatically in $\Oh(\Hb^n)$, we see that $\lim_{\eps\to0}\tilde\Pv_1f_\eps(z)=\tilde \Pv_1f(z)$ and $\lim_{\eps\to0}\tilde\Pv_2f_\eps(z)=\tilde \Pv_2f(z)$ locally uniformly in $z$ as well.
Therefore $\tilde\Pv f_\eps(z)\xrightarrow{\eps\to0}\tilde\Pv f(z)$ locally uniformly in $z\in\Hb^n$. By Proposition \ref{Lem::SecHolLap::PfTildePHolo} $\tilde\Pv f_\eps(z)$ is holomorphic, we conclude that $\tilde\Pv f(z)$ is also holomorphic.

For $f\in\Co^{(-4)+}_\Oh(\Hb^n)$, by Lemma \ref{Lem::SecHolLap::HLLem} \ref{Item::SecHolLap::HLLem::Res}  $f|_{\B^n}\in\Co^{(-4)+}(\B^n;\C)$ so $\Pv[f|_{\B^n}]\in\Co^{(-2)+}(\B^n;\C)$ is defined. 

By \eqref{Eqn::SecHolLap::NotSzV} we have $S_{x+i0}=\B^n$ for all $x\in\B^n$. By comparing \eqref{Eqn::SecHolLap::FormulaForP} with \eqref{Eqn::SecHolLap::EqnHoloP1}, \eqref{Eqn::SecHolLap::EqnHoloP2} and \eqref{Eqn::SecHolLap::EqnHoloP3}, we get $(\tilde \Pv f)|_{\B^n}=\Pv[f|_{\B^n}]$ for $f\in\Co^{(-4)+}_\Oh(\Hb^n)$.

And since $\B^n\subset\C^n$ is a totally real submanifold with full dimension $n$, by uniqueness of the holomorphic extension we have $\tilde \Pv f=\Ex\big[\Pv[f|_{\B^n}]\big]$ for all $f\in\Co^{(-4)+}_\Oh(\Hb^n)$. Therefore $\Delta_z\tilde \Pv f=\Ex\big[\Delta_x\Pv[f|_{\B^n}]\big]=\Ex[f|_{\B^n}]=f$ holds for $f\in\Co^{(-4)+}_\Oh(\Hb^n)$ as well.
\end{proof}

\subsection{H\"older-Zygmund regularities for $\tilde{\mathbf P}$}\label{Section::SecHolLap::BddTildeP}

Recall $\tilde\Pv$ in Definition \ref{Defn::SecHolLap::DefofTildeP}. In this subsection we are going to prove $\tilde \Pv:\Co^\alpha_\Oh(\Hb)\to\Co^{\alpha+2}_\Oh(\Hb)$ for $-4<\alpha<-1$. By Proposition \ref{Prop::SecHolLap::TildePisHolo} we already know $\tilde\Pv:\Co^{(-4)+}_\Oh(\Hb^n)\to\Oh(\Hb^n)$ is the holomorphic extension of $\Pv$. 

\begin{proof}[Proof of Proposition \ref{Prop::HolLap} \ref{Item::HolLap::TildeP}]
What we remain to do is to show that for $-4<\alpha<-1$, there is a $C=C_{n,\alpha}>0$ such that $$|\nabla_z\tilde\Pv f(z)|\le C_{n,\alpha}\|f\|_{\Co^\alpha_\Oh(\Hb^n)}\dist(z,\partial\Hb^n)^{\alpha+1},\quad f\in\Co^\alpha_\Oh(\Hb^n).$$

 By Proposition \ref{Lem::SecHolLap::PfTildePHolo} we have for $k=1,\dots,n$ and $z=x+iy\in\Hb^n$,
\begin{align}
    \partial_{z^k}\tilde\Pv f(z)=
    &\label{Eqn::SecHolLap::EstHoloP1}
    \int_{|t-x|<\frac{1-|x|}2}f(t+i\lambda(t,x)y)\Ga_{,k}(x-t+i(1-\lambda(t,x))y)(1+iy\cdot\partial_t\lambda(t,x))dt
    \\\label{Eqn::SecHolLap::EstHoloP2}
    &+\int_{|t-x|<\frac{1-|x|}2}f(t)\sum_{j=1}^9\frac{a_jb_j\chi(|t|)}{|t|^{n(b_j+1)}}\Ga_{,k}\Big(z-\frac t{|t|^{b_j+1}}\Big)dt
    \\
    &\label{Eqn::SecHolLap::EstHoloP3}
    +\int_{t\in\B^n:|t-x|>\frac{1-|x|}2}f(t)\bigg(\Ga_{,k}(z-t)-\sum_{j=1}^9\frac{a_jb_j\chi(|t|)}{|t|^{n(b_j+1)}}\Ga_{,k}\Big(z-\frac{t}{|t|^{b_j+1}}\Big)\bigg)dt.
\end{align}

We need to show the absolute value of \eqref{Eqn::SecHolLap::EstHoloP1}, \eqref{Eqn::SecHolLap::EstHoloP2} and \eqref{Eqn::SecHolLap::EstHoloP3} are all bounded by a constant times $\|f\|_{\Co^\alpha_\Oh}\dist(z,\partial\Hb^n)^{\alpha+1}$.
 
We first estimate \eqref{Eqn::SecHolLap::EstHoloP1}. By Lemma \ref{Lem::SecHolLap::HLLem} \ref{Item::SecHolLap::HLLem::0Char} we have $|f(z)|\lesssim\|f\|_{\Co^\alpha_\Oh}\dist(z,\partial\Hb^n)^\alpha$; by \eqref{Eqn::SecHolLap::DevGa1} and \eqref{Eqn::SecHolLap::DefLambda} we have $|F_{,k}(z-t+i(1-\lambda)y)|\lesssim|x-t|^{1-n}$; by property \ref{Item::SecHolLap::LambdaPro1} we have $|1+iy\cdot\partial_t\lambda|\lesssim1$. Therefore
\begin{equation}\label{Eqn::SecHolLap::PfEstHolo::Tmp1}
    \Big|\int_{|t-x|<\frac{1-|x|}2}f\Ga_{,k}Vdt\Big|\lesssim_{n,\alpha}\|f\|_{\Co^\alpha_\Oh}\int_{|t-x|<\frac{1-|x|}2}\dist(t+i\lambda y,\partial\Hb^n)^\alpha\cdot |x-t|^{1-n}dt 
\end{equation}


Since $t+i\lambda y$ lays in a line segment connecting $z$ and a point in $\partial B^n(x,\frac{1-|x|}2)$, by convexity of $\Hb^n$ we get $\dist(t+i\lambda(t,x)y,\partial\Hb^n)\ge\lambda(t,x)\dist(z,\partial\Hb^n)+(1-\lambda(t,x))\dist(\partial B^n(x,\frac{1-|x|}2),\partial\Hb^n)$. Since $\dist(z,\partial\Hb^n)\approx\frac{1-|x|}2-2|y|$ and $\dist\big(\partial B^n(x,\frac{1-|x|}2),\partial\Hb^n\big)\approx \frac{1-|x|}2$, by plugging in the expression of $\lambda(t,x)$ we get (recall $\alpha+1<0$)
\begin{align*}
    &\int_{|t-x|<\frac{1-|x|}2}\dist(t+i\lambda y,\partial\Hb^n)^\alpha\cdot |x-t|^{1-n}dt
    \\
    \le& \int_{|t-x|<\frac{1-|x|}2}\big(\lambda\dist(z,\partial\Hb^n)+(1-\lambda)\dist\big(\partial B^n(x,\tfrac{1-|x|}2),\partial\Hb^n\big)\big)^\alpha\cdot |x-t|^{1-n}dt 
    \\
    \lesssim&\int_{B^n(0,\frac{1-|x|}2)}\Big(\big(1-\tfrac{2|s|}{1-|x|}\big)\big(\tfrac{1-|x|}2-2|y|\big)+\tfrac{2|s|}{1-|x|}\tfrac{1-|x|}2\Big)^\alpha|s|^{1-n}ds
    \\
    \lesssim&\int_0^{\frac{1-|x|}2}\Big(\tfrac{1-|x|}2-2|y|+\tfrac{4|y|r}{1-|x|}\Big)^\alpha dr\lesssim_\alpha\Big(\tfrac{1-|x|}2-2|y|\Big)^{\alpha+1}\approx_\alpha\dist(z,\partial\Hb^n)^{\alpha+1}.
\end{align*}
Putting this in \eqref{Eqn::SecHolLap::PfEstHolo::Tmp1} we control \eqref{Eqn::SecHolLap::EstHoloP1}.

\medskip
For \eqref{Eqn::SecHolLap::EstHoloP2} and \eqref{Eqn::SecHolLap::EstHoloP3}, we use $\dist(t,\partial\Hb^n)\approx\dist(t,\partial\B^n)=1-|t|$ for $t\in\B^n$. By Lemma \ref{Lem::SecHolLap::HLLem} \ref{Item::SecHolLap::HLLem::0Char} $|f(t)|\lesssim\|f\|_{\Co^\alpha_\Oh}(1-|t|)^\alpha$.

By \eqref{Eqn::SecHolLap::DevGa1} we have $|\Ga_{,k}(z-\frac t{|t|^{b_j+1}})|\lesssim|x-\frac t{|t|^{b_j+1}}|^{1-n}\le (1-|x|)^{1-n}$ since $\frac t{|t|^{b_j+1}}\notin\B^n$ when $t\in\B^n$. So we bound \eqref{Eqn::SecHolLap::EstHoloP2} by the following
\begin{align*}
    &\Big|\int_{|t-x|<\frac{1-|x|}2}f(t)\sum_{j=1}^9\frac{a_jb_j\chi(|t|)}{|t|^{n(b_j+1)}}\Ga_{,k}\Big(z-\frac t{|t|^{b_j+1}}\Big)dt\Big|
    \lesssim\|f\|_{\Co^\alpha_\Oh}\int_{|t-x|<\frac{1-|x|}2}(1-|t|)^\alpha(1-|x|)^{1-n}dt
    \\\lesssim&\|f\|_{\Co^\alpha_\Oh}(1-|x|)^{\alpha+1}\lesssim\|f\|_{\Co^\alpha_\Oh}\dist(z,\partial\Hb^n)^{\alpha+1}.
\end{align*}

By Lemma \ref{Lem::SecHolLap::AprioriBddforP} we bound \eqref{Eqn::SecHolLap::EstHoloP3} by the following
\begin{align*}
    &\bigg|\int_{t\in\B^n:|t-x|>\frac{1-|x|}2}f(t)\bigg(\Ga_{,k}(z-t)-\sum_{j=1}^9\frac{a_jb_j\chi(|t|)}{|t|^{n(b_j+1)}}\Ga_{,k}\Big(z-\frac{t}{|t|^{b_j+1}}\Big)\bigg)dt\bigg|
    \\
    \lesssim&\|f\|_{\Co^\alpha_\Oh}\int_{t\in\B^n:|t-x|>\frac{1-|x|}2}(1-|t|)^\alpha|x-t|^{-n-3}(1-|t|)^4dt
    \\
    \lesssim&\|f\|_{\Co^\alpha_\Oh}\bigg(\int_{|t|<\frac{3|x|-1}2}+\int_{\substack{\frac{3|x|-1}2<|t|<1\\|t-x|>\frac{1-|x|}2}}\bigg)(1-|t|)^{4+\alpha}|x-t|^{-n-3}dt
    \\
    \lesssim&\|f\|_{\Co^\alpha_\Oh}\int_0^{\frac{3|x|-1}2}(1-r)^{\alpha+4}r^{n-1}dr\int_{|x|-r}^1\rho^{-n-3}\rho^{n-2}d\rho+\|f\|_{\Co^\alpha_\Oh}\int_{\frac{3|x|-1}2}^1(1-r)^{\alpha+4}dr\int_{\frac{1-|x|}2}^1\rho^{-n-3}\rho^{n-2}d\rho
    \\
    \lesssim&\|f\|_{\Co^\alpha_\Oh}\Big(\int_0^{|x|}(1-r)^{\alpha}dr+(1-|x|)^{\alpha+5}(1-|x|)^{-4}\Big)\approx\|f\|_{\Co^\alpha_\Oh}(1-|x|)^{\alpha+1}\lesssim\|f\|_{\Co^\alpha_\Oh}\dist(z,\partial\Hb)^{\alpha+1}.
\end{align*}
Here we use $\rho$ as the radial parameter of polar coordinates on $\Sp^{n-1}$ centered at $\frac x{|x|}$.

Thus \eqref{Eqn::SecHolLap::EstHoloP1}, \eqref{Eqn::SecHolLap::EstHoloP2} and \eqref{Eqn::SecHolLap::EstHoloP3} are all bounded by $\|f\|_{\Co^\alpha}\dist(z,\partial\Hb^n)^{\alpha+1}$. The proof is now complete.
\end{proof}

\chapter{On Singular Subbundles, Rough Involutivity and Integrability}\label{Chapter::DisInv}

In many literature people use ``distributions'' for subbundles or generalized subbundles. In the thesis we keep the word strictly referring to \textit{generalized functions}. We still use ``subbundles'' for \textit{regular tangent subbundles}. 

In Sections \ref{Section::DisInv::CharSec} and \ref{Section::DisInv::CharInv} we only discuss the characterization of the regular subbundles. In Section \ref{Section::DisInv::Sing} we give a brief discussion on singular subbundles. 
\begin{defn}\label{Defn::DisInv::GenSub}
A \textbf{generalized complex tangential subbundle} on a $C^1$ manifold $\Mf$ is a subset $\Gc\subseteq\C T\Mf$ such that for every $p\in \Mf$, $\Gc_p:=\Gc\cap\C T_p\Mf$ is a complex linear subspace of $\C T_p\Mf$.

Let $\alpha\in\R_\Eb^+$ be a positive generalized index (recall Definition \ref{Defn::Hold::ExtIndex}) and assume that $\Mf$ is at least $\Co^{\alpha+1}$. We say the generalized subbundle $\Gc\subseteq \C T\Mf$ is $\Co^\alpha$, if for any $p\in \Mf$, there are a neighborhood $U\subseteq \Mf$ of $p$ and finitely many $\Co^\alpha$ complex vector fields $X_1,\dots, X_m\in\Co^\alpha(U;\C TU)$, such that the fiber $\Gc_q:=\Gc\cap\C T_q\Mf$ is a linear subspace generated by $X_1|_q,\dots,X_m|_q$  for every $q\in U$.

We say a generalized subbundle $\Gc$ is a \textbf{regular subbundle}, or just a subbundle, if its fibres have constant rank, i.e. the map $p\in \Mf\mapsto\rank(\Gc\cap\C T_p\Mf)$ is constant.

We say a generalized subbundle $\Gc$ is \textbf{singular}, if it is not regular.
\end{defn}

The definition is equivalent if we remove the finiteness assumption of the vector fields. See Proposition \ref{Prop::DisInv::SerreSwanGenSub}.

To discuss the involutivity of generalized subbundles, it is more convenient to consider this property on module of sections.
\begin{defn}\label{Defn::DisInv::ModInv}
Let $\alpha\in\R_\Eb$ be a generalized index such that $\alpha>\frac12$. Let $\beta\in\R$ such that $\alpha\ge\beta-1$ and $\alpha>-\beta$. Let $\Fs\subseteq\Co^\alpha_\loc(\Mf;\C T\Mf)$ be a submodule over the ring $\Co^\alpha_\loc(\Mf;\C)$.

We say $\Fs$ is $\Co^\beta$-involutive, if for every vector fields $X,Y\in\Fs$, the Lie bracket $[X,Y]\in\Co^\beta_\loc(\Mf;\C)\otimes_{\Co^\alpha_\loc(\Mf;\C)}\Fs$.
\end{defn}
One can also consider the characterization on sheaves of $\Co^\alpha$ rings. But sheaves are more useful than modules only when we discuss real-analytic structures, or other function structures that have rigidity. In our discussion everything are at most smooth, and we can use partition of unity to patch the local constructions into a global one. 

When $\Se$ is a regular subbundle, in Proposition \ref{Prop::DisInv::CharInv1} we see that the involutivity of the submodule $\Co^\alpha_\loc(\Mf;\Se)\subseteq\Co^\alpha_\loc(\Mf;\C T\Mf)$ does not depend on the choice of $\beta$ in Definition \ref{Defn::DisInv::ModInv}. However for a general $\Co^\alpha$ submodules $\Fs$ it is possible that $\Fs$ is not $\Co^{\alpha-1}$-involutive but is $\Co^\beta$-involutive for some $\beta<\alpha-1$. See Examples \ref{Exmp::DisInv::ModInv1} and \ref{Exmp::DisInv::ModInv2}.

\section{Distributional sections of regular subbundles}\label{Section::DisInv::CharSec}
In this section we give some equivalent characterizations of sections in the sense of distributions. We focus on complex bundles, since for real bundles all properties and characterizations below work and the proof are simpler.

Using Lemmas \ref{Lem::Hold::PushForwardFuncSpaces} and \ref{Lem::Hold::LogL-1Comp}, we can define objects on manifolds via patchings, which are  analogous to \cite[Definition 6.3.3]{Hormander}.

\begin{defn}\label{Defn::DisInv::DefFunVF}
 Let $\kappa>1$ and let $\Mf$ be a $n$ dimensional $\Co^\kappa$-manifold with maximal $\Co^\kappa$-atlas $\Af=\{\psi:U_\psi\subseteq \Mf\to\R^n\}_\psi$. Namely, each $\psi\in\Af$ is a homeomorphism $\psi=(\psi^1,\dots,\psi^n):U_\psi\xrightarrow{\sim}\psi(U_\psi)\subseteq\R^n$; we have $\psi\circ\phi^\Inv\in \Co^\kappa_\loc(\phi(U_\phi\cap U_\psi);\R^n)$ whenever $\psi,\phi\in\Af$ satisfy $U_\phi\cap U_\psi\neq \varnothing$; and $\Af$ is maximal with these properties.

\begin{enumerate}[parsep=-0.3ex,label=(\roman*)]
    \item  Let $\beta\in(1-\kappa,\kappa]\cup\{\LogL-1\}$. A (locally) $\Co^\beta$-complex function $f$ is a collection $\{f_\psi\in\Co^\beta_\loc(\psi(U_\psi);\C):\psi\in\Af\}$ such that
    \begin{equation}\label{Eqn::DisInv::DefFunVF::TransFun}
        f_\psi=f_\phi\circ(\phi\circ\psi^\Inv),\quad\text{on }\psi(U_\psi\cap U_\phi)\subseteq\R^n\text{ whenever }U_\psi\cap U_\phi\neq\varnothing.
    \end{equation}
    
    We denote by $\Co^\beta_\loc(\Mf;\C)$ the space of all $\Co^\beta$-functions on $\Mf$.
    \item\label{Item::DisInv::DefFunVF::VF} Let $\beta\in(1-\kappa,\kappa-1]\cup\{\LogL-1\}$. A (locally) $\Co^\beta$-complex vector field $X$ is a collection $\{X_\psi=(X_\psi^1,\dots,X_\psi^n)\in\Co^\beta_\loc(\psi(U_\psi);\C^n):\psi\in\Af\}$ such that
    \begin{equation}\label{Eqn::DisInv::DefFunVF::TransVF}
        X_\psi^i=\sum_{j=1}^n\big(X_\phi^j\circ (\phi\circ \psi^\Inv)\big)\cdot\Big(\frac{\partial \psi^i }{\partial \phi^j}\circ \psi^\Inv\Big),\quad\text{on }\psi(U_\psi\cap U_\phi)\subseteq\R^n\text{ whenever }U_\psi\cap U_\phi\neq\varnothing.
    \end{equation}
    We denote by $\Co^\beta_\loc(\Mf;\C T\Mf)$ the space of all $\Co^\beta$ complex vector fields on $\Mf$. We denote by $\Co^\beta_\loc(\Mf; T\Mf)$  the space of all $\Co^\beta$ real vector fields on $\Mf$.
    \item\label{Item::DisInv::DefFunVF::Form} Let $\beta\in(1-\kappa,\kappa-1]\cup\{\LogL-1\}$. Let $k=1,\dots,n$. A (locally) $\Co^\beta$ complex $k$-form $\omega$ is a collection $\{\omega^\psi=(\omega^\psi_{i_1\dots i_k})_{ i_1,\dots,i_k=1}^n\in\Co^\beta_\loc(\psi(U_\psi);\C^{n^k}):\psi\in\Af\}$ such that $\omega^\psi_{\dots i_ji_{j+1}\dots}=-\omega^\psi_{\dots i_{j+1}i_j\dots}$ for every $1\le j\le k-1$, $1\le i_j,i_{j+1}\le n$, $\psi\in\Af$. And for every $\phi,\psi\in\Af$ such that $U_\psi\cap U_\phi\neq\varnothing$,
    \begin{equation}\label{Eqn::DisInv::DefFunVF::TransForm}
        \omega^\psi_{i_1\dots i_k}=\sum_{ j_1,\dots,j_k=1}^n\big(\omega^\psi_{j_1\dots j_k}\circ (\phi\circ \psi^\Inv)\big)\cdot\Big(\frac{\partial \phi^{j_1} }{\partial \psi^{i_1}}\circ \psi^\Inv\Big)\dots\Big(\frac{\partial \phi^{j_k} }{\partial \psi^{i_k}}\circ \psi^\Inv\Big),\quad\text{on }\psi(U_\psi\cap U_\phi)\subseteq\R^n.
    \end{equation}
    We denote by $\Co^\beta_\loc(\Mf;\wedge^k\C T^*\Mf)$ the space of all $\Co^\beta$ complex $k$-form on $\Mf$. We denote by $\Co^\beta_\loc(\Mf;\wedge^k T^*\Mf)$ the space of all $\Co^\beta$ real $k$-form on $\Mf$.
    \item\label{Item::DisInv::DefFunVF::dForm} If $2-\kappa<\beta\le\kappa-1$ (provided $\kappa>\frac32$), for a $\Co^\beta$ $k$-form $\omega=\{\omega_{i_1\dots i_k}^\psi\}$ on $\Mf$, the exterior differential $d\omega=\{(d\omega)_{i_1\dots i_{k+1}}^\psi\}$ is defined in the way that
    \begin{equation*}
        \sum_{1\le i_1<\dots<i_{k+1}\le n}(d\omega)_{i_1\dots i_{k+1}}^\psi dx^{i_1}\wedge\dots\wedge dx^{i_{k+1}}:=\sum_{1\le j_1<\dots<j_k\le n}d(\omega_{j_1\dots j_k}^\psi)\wedge dx^{j_1}\wedge\dots\wedge dx^{j_{k}}\quad\text{on }\psi(U_\psi)\subseteq\R^n_x.
    \end{equation*}
\end{enumerate}
\end{defn}

\begin{remark}\label{Rmk::DisInv::RmkMfldObj}
\begin{enumerate}[parsep=-0.3ex,label=(\alph*)]
    \item Here \eqref{Eqn::DisInv::DefFunVF::TransFun}, \eqref{Eqn::DisInv::DefFunVF::TransVF} and \eqref{Eqn::DisInv::DefFunVF::TransForm} are given by the pullback equations $f_\psi=(\phi\circ\psi^\Inv)^*f_\phi$, $X_\psi=(\phi\circ\psi^\Inv)^*X_\phi$ and $\omega_\psi=(\phi\circ\psi^\Inv)^*\omega_\phi$. By Lemma \ref{Lem::Hold::PushForwardFuncSpaces} if $\beta\in\R_-$ and Lemma \ref{Lem::Hold::LogL-1Comp} if $\beta=\LogL-1$, the right hand sides of \eqref{Eqn::DisInv::DefFunVF::TransFun}, \eqref{Eqn::DisInv::DefFunVF::TransVF} and \eqref{Eqn::DisInv::DefFunVF::TransForm} are well-defined $\Co^{\beta}_\loc$ objects. 
    \item\label{Item::DisInv::RmkMfldObj::CheckCover} In particular, to show a function $f$ (or a vector fields $X$ or a form $\omega$) is $\Co^{\beta}$, it suffices to show $f_{\psi_j}\in \Co^{\beta}_\loc(\psi_j(U_{\psi_j});\C)$ (or $X_{\psi_j}\in \Co^{\beta}_\loc(\psi_j(U_{\psi_j});\C^n)$ or $\omega^{\psi_j}\in \Co^{\beta}_\loc(\psi_j(U_{\psi_j});\C^{n^k})$) on a coordinate cover $\{\psi_j\}_j$, namely a subcollection of atlas $\{\psi_j:U_{\psi_j}\to\R^n\}_j$ such that $\bigcup_jU_{\psi_j}=\Mf$. In this way to study the local problems we can use one coordinate chart near a fixed point and work on subsets in Euclidean spaces.
    \item We do not talk about ``bounded $\Co^\beta$ functions or vector fields'' on a general manifold, since a manifold may not have structures like metric to control the size of functions or vector fields.
\end{enumerate}
\end{remark}

\begin{remark}
The notations $\Co^\beta(\Mf;T\Mf)$ and $\Co^\beta(\Mf;\wedge^kT^*\Mf)$ are different from taking $\Nf=T\Mf$ or $\Nf=\wedge^kT^*\Mf$ in Notation \ref{Note::ODE::DiffMap}. By $\Co^\beta(\Mf;T\Mf)$ we mean the space of vector fields, i.e. all $\Co^\beta$ maps $X:\Mf\to T\Mf$ that satisfies $X(p)\in T_p\Mf$. 

In practice there are $\Co^\beta$ maps $F:\Mf\to T\Mf$ we use in the proof such that $F(p)\notin T_p\Mf$. For example, the map $\Coorvec t\exp_X(t,\cdot)=X\circ\exp_X(t,\cdot):\Mf\to T\Mf$ sends every $p\in\Mf$ to a tangent vector in $T_{\exp(t,p)}\Mf$ which in general is not $T_p\Mf$. We do not use the word ``$F\in\Co^\beta(\Mf;T\Mf)$'' for such map $F$ in the thesis.
\end{remark}

\begin{remark}\label{Rmk::DisInv::Functional}
A $\Co^\beta$-function on manifold can be identified as a linear functional by the way of \cite[Chapter 3.1]{GeomDist}. We sketch the construction as the following: 

For a (smooth or non-smooth) $n$-dimensional manifold $\Mf$ with atlas $\Af=\{\psi:U_\psi\subseteq \Mf\to\R^n\}$, we define the \textit{volume bundle} $\Vol_\Mf$ to be a real rank 1 vector bundle on $\Mf$ as
\begin{gather}\label{Eqn::DisInv::DefVolBundle}
    \Vol_\Mf:=\Big(\coprod_{\psi\in\Af}\psi(U_\psi)\times\R\Big)\Big/\left\{(x,t)\sim(\phi\circ\psi^\Inv(x),g_{\phi\psi}(x)t):x\in \psi(U_\phi\cap U_\psi),t\in\R\right\},
    \\
    \text{where }g_{\phi\psi}(x):=|\det \nabla(\phi\circ\psi^\Inv)(x)|,\quad x\in \psi(U_\phi\cap U_\psi).\notag
\end{gather}

When $(\Mf,\Af)$ is a $\Co^\kappa$-manifold, then $\Vol_\Mf$ is a $\Co^{\kappa-1}$ vector bundle. Similar to Definition \ref{Defn::DisInv::DefFunVF} a $\Co^{\kappa-1}$-section of $\Vol_\Mf$ can be defined by a collection $s=\{s_\psi\in\Co^{\kappa-1}_\loc(\psi(U_\psi)):\psi\in\Af\}$ such that $s_\psi=(s_\phi\circ(\phi\circ\psi^\Inv))\cdot|\det \nabla(\phi\circ\psi^\Inv)|$.

In this case take a $\Co^\kappa$ coordinate cover $\{\psi_i:U_{\psi_i}\subset \Mf\to\R^n\}_{i\in I}$ of $\Mf$ and take a $\Co^\kappa$ partition of unity $\{\chi_i\in \Co_\loc^\kappa(U_{\psi_i})\}_{i\in I}$, we have a pairing: for $f=(f_\psi)_{\psi\in\Af}\in C^0_\loc(\Mf)$ and $s=(s_\psi)_{\psi\in\Af}\in\Co^{\kappa-1}_c(\Mf)$,
\begin{equation*}
    \langle f,s\rangle:=\sum_{i\in I}\int_{\psi_i(U_{\psi_i})}\chi_i(\psi_i^\Inv(x))f_{\psi_i}(x)s_{\psi_i}(x)dx.
\end{equation*}
The pairing is still defined for $f\in \Co^{\beta}_\loc(\Mf)$, thus we have a map $\Co^{\beta}_\loc(\Mf)\hookrightarrow \Co_c^\kappa(\Mf;\Vol_\Mf)'$ for $\beta\in(1-\kappa,\kappa]$. The pairing is canonical in the sense that it does not depend on the choice of coordinates and partition of unity. Therefore we can say $\Co^{\beta}_\loc(\Mf)\subset \Co_c^\kappa(\Mf;\Vol_\Mf)'$ for all $0<\beta<1$. 

Similarly we have the embedding of the space of vector fields $\Co^{\beta}_\loc(\Mf; T\Mf)\subset \Co_c^\kappa(\Mf;T^*\Mf\otimes_\R\Vol_\Mf)'$ and the space of differential forms $\Co^{\beta}_\loc(\Mf;\wedge^kT^*\Mf)\subset \Co_c^\kappa(\Mf;(\wedge^k T\Mf)\otimes_\R\Vol_\Mf)'$ for $1\le k\le n$. For complex vector fields and complex forms the almost same embedding hold.

If $\Mf$ is a oriented manifold, then $\Vol_\Mf$ is canonically isomorphic to $\wedge^nT^*\Mf$, the top degree alternating tensor of cotangent bundle. In this case a distribution is indeed an $n$-current, see \cite[Example 3.1.17]{GeomDist}. We leave the details to the reader.
\end{remark}

To make the result more general we can consider general vector bundles and their subbundles.
\begin{defn}\label{Defn::DisInv::DefVBSec}
Let $n,m\in\Z_+$, $\kappa>1$, $\alpha\in(0,\kappa]$ and let $\Mf$ be a $n$-dimensional $\Co^\kappa$ manifold with maximal $\Co^\kappa$-atlas $\Af$. A rank $m$ $\Co^\alpha$ complex vector bundle $\E$ over $\Mf$ is a $(n+m)$-dimensional topological manifold with a surjective topological submersion $\pi:\E\twoheadrightarrow\Mf$ and an atlas $\Bf$ such that
\begin{itemize}[parsep=-0.3ex]
    \item Elements in $\Bf$ are homeomorphisms $\Phi=(\psi,\phi):\pi^{-1}(U_\psi)\xrightarrow{\sim}\psi(U_\psi)\times\C^m(\subseteq\R^n\times\C^m)$ where $\psi\in\Af$.
    \item For every $\Phi_1=(\psi_1,\phi_1),\Phi_2=(\psi_2,\phi_2)\in\Bf$ such that $U_{\psi_1}\cap U_{\psi_2}\neq\varnothing$, the map $\phi_1\circ\Phi_2^\Inv:\psi_2(U_{\psi_1}\cap U_{\psi_2})\times\C^m\to\C^m$ is a $\Co^\alpha$-map such that $\phi_1\circ\Phi_2^\Inv:\C^m\to\C^m$ is linear for each $x\in \psi_2(U_{\psi_2})$. In other words we can write $g_{\Phi_1\Phi_2}\in\Co^\alpha(\psi_2(U_{\psi_1}\cap U_{\psi_2});\C^{m\times m})$ such that
    \begin{equation*}
        \Phi_1\circ\Phi_2^\Inv(x,v)=(\psi_1\circ\psi_2^\Inv(x),g_{\Phi_1\Phi_2}(x)\cdot v),\quad\forall x\in \psi_2(U_{\psi_1}\cap U_{\psi_2}),\quad v\in\C^m.
    \end{equation*}
    \item $\Bf$ is maximal among the above two conditions.
\end{itemize}

Let $\beta\in(\max(1-\kappa,-\alpha),\alpha]\cup\{\LogL-1\}$. A $\Co^\beta$-section of $\E$ is a collection $\xi=\{\xi_{(\psi,\phi)}\in\Co^\beta_\loc(\psi(U_\psi);\C^m):(\psi,\phi)\in\Bf\}$ such that 
\begin{equation*}
    \xi_{(\psi',\phi')}=g_{(\psi,\phi)(\psi',\phi')}\cdot(\xi_{(\psi,\phi)}\circ\psi\circ\psi'^\Inv),\quad\text{on }\psi'(U_{\psi'}\cap U_\psi)\subseteq\R^n\text{ whenever }U_{\psi'}\cap U_\psi\neq\varnothing.
\end{equation*}

We denote by $\Co^\beta_\loc(\Mf;\E)$ the space of all (complex) $\Co^\beta$-sections of $\E$.
\end{defn}
For real vector bundles the definition is the same except we replace $\C^m$ and $\C^{m\times m}$ by $\R^m$ and $\R^{m\times m}$.

\begin{remark}\label{Rmk::DisInv::RmkSecVB}
\begin{enumerate}[parsep=-0.3ex,label=(\alph*)]
    \item Just like manifolds, a coordinate cover $\{\Phi_j\}\subset\Bf$ of $\E$ uniquely determines its bundle atlas.
    \item As the direct correspondences to Remark \ref{Rmk::DisInv::RmkMfldObj}, we see that for each $\Phi,\Phi'\in\Bf$ we can say $\xi_\Phi=(\Phi'\circ\Phi^\Inv)^*\xi_{\Phi'}$. The well-definedness of the pullback is guaranteed by Lemmas \ref{Lem::Hold::PushForwardFuncSpaces} and \ref{Lem::Hold::LogL-1Comp}. And to check a $\xi$ is $\Co^\beta$ it suffices to check $\xi_{(\psi_j,\phi_j)}\in\Co^\beta_\loc(\psi_j(U_{\psi_j});\C^m)$ where $\{(\psi_j,\phi_j)\}\subset\Bf$ such that $\{\psi_j\}$ is a  coordinate cover of $\Mf$. We leave the proof to readers.
    \item\label{Item::DisInv::RmkSecVB::DefSec} Let $\Se\le\C T\Mf$ be  a complex tangent subbundle  and let $\beta\le0$. By Define \ref{Defn::DisInv::DefVBSec} we define the $\Co^\beta$-sections of $\Se$ by considering $\Se$ as an abstract vector bundle. The definition is different from \cite[Definition 1.12]{YaoLLFro}, where we define $\Co^\beta_\loc(\Mf;\Se)$ to be the subset of space of vector fields $\Co^\beta_\loc(\Mf;\C T\Mf)$. Nevertheless Definition \ref{Defn::DisInv::DefVBSec} for $\Se$ and \cite[Definition 1.12]{YaoLLFro} coincide by Proposition \ref{Prop::DisInv::CharSec}.
    
\end{enumerate}
\end{remark}

\begin{lem}\label{Lem::DisInv::DisSecVB}
    Let $\kappa\in(1,\infty]$, $\alpha\in(0,\kappa]$. Let $\Mf$ be a $\Co^\kappa$-manifold, and let $\E$ be a $\Co^{\alpha}$ complex vector bundle over $\Mf$. Let $\gamma\in(0,\alpha]$ and $\beta\in (\max(1-\kappa,\gamma),\gamma]\cup\{\LogL-1\}$.
    \begin{enumerate}[parsep=-0.3ex,label=(\roman*)]
        \item\label{Item::DisInv::DisSecVB::FiniteGen} The space of sections $\Co^{\gamma}_\loc(\Mf;\E)$ is a finitely generated module over the ring of function $\Co^{\gamma}_\loc(\Mf;\C)$. In particular if $\xi_1,\dots,\xi_N\in \Co^{\gamma}_\loc(\Mf;\E)$ generate the module, then $\xi_1(p),\dots,\xi_N(p)$ span the fibre $\E_p$ for every $p\in \Mf$.
        \item\label{Item::DisInv::DisSecVB::Loc}Sections has local property. More precisely, let $\{U_j\}_{j\in I}$ be an open cover of $\Mf$, then $\xi\in\Co^\beta_\loc(\Mf;\E)$ if and only if $\xi|_{U_j}\in \Co^\beta_\loc(U_j;\E|_{U_j})$ for all $j$.
        \item\label{Item::DisInv::DisSecVB::Ten} We have equality of tensor products, or the so-called based change formula,
        \begin{equation}\label{Eqn::DisInv::DisSecVB::Ten}
            \Co^\beta_\loc(\Mf;\E)=\Co^\beta_\loc(\Mf;\C)\otimes_{\Co^\alpha_\loc(\Mf;\C)}\Co^{\alpha}_\loc(\Mf;\E)=\Co^\beta_\loc(\Mf;\C)\otimes_{\Co^{\gamma}_\loc(\Mf;\C)}\Co^{\gamma}_\loc(\Mf;\E).
        \end{equation}
        
    \end{enumerate}
\end{lem}
\begin{remark}\label{Rmk::DisInv::RmkDisSecVB}
\begin{enumerate}[parsep=-0.3ex,label=(\alph*)]
    \item Algebraically tensor products can only be equal in the sense of a unique isomorphism. Here we can make them equal set theoretically by  embedding both sides of \eqref{Eqn::DisInv::DisSecVB::Ten} in a larger ambient spaces. For example, we can view them as linear functionals on $\Co^{\gamma}_c(\Mf;\E^*\otimes\Vol_\Mf)$ (see Remark \ref{Rmk::DisInv::Functional}).
    \item\label{Item::DisInv::RmkDisSecVB::NoDepReg} We can use either side of \eqref{Eqn::DisInv::DisSecVB::Ten} to define the space of $\Co^{\beta}$-sections of $\E$. We see that the definition does not depend on the regularity assumption $\E\in \Co^{\alpha}$: the space is the same if we view $\E$ as a $\Co^{\gamma}$-vector bundle. In fact the space $\Co^{\beta}_\loc(\Mf;\E)$ is defined as long as $\E\in \Co^{(-\beta)+}$.
\end{enumerate}
\end{remark}
\begin{proof}[Proof of Lemma \ref{Lem::DisInv::DisSecVB}] We assume that $\dim\Mf=n$ and $\rank\E=m$.

\noindent\ref{Item::DisInv::DisSecVB::FiniteGen}: This is the non-smooth version of the Swan's Theorem. One can find the proof in \cite[Theorem 5.3]{Lewis}. In fact one can take $N=(\dim\Mf+1)\cdot\rank\E$ many generators.

\medskip
\noindent\ref{Item::DisInv::DisSecVB::Loc}: Recall that $\xi$ is defined to be a collection $\{\xi_\Phi:\Phi\in\Bf\}$ of distributions.

Let $(\psi,\phi)\in\Bf=\Bf^\Mf$ be a $\Co^\alpha$ bundle chart of $\Mf$. By assumption $\{U_j\cap U_\psi\}$ covers $U_\psi\subseteq\Mf$, the domain of $(\psi,\phi)$. Therefore $\xi_\psi\in\Co^\beta(\psi(U_\psi);\C^m)$ if and only if $\xi_\psi|_{\psi(U_\psi\cap U_j)}\in\Co^\beta(\psi(U_\psi\cap U_j);\C^m)$ for all $j\in I$, whenever $U_\psi\cap U_j\neq\varnothing$. Since $(\psi,\phi)|_{U_j\cap U_\psi}\in\Bf^{U_j}$ is a  $\Co^\alpha$ bundle chart of $U_j$ for each $j$, we finish the proof.

\medskip
\noindent\ref{Item::DisInv::DisSecVB::Ten}: Take finitely many $\Co^{\alpha}$-sections $\xi_1,\dots,\xi_N$ that generates $\Co^{\alpha}_\loc(\Mf;\E)$. Clearly a $\Co^{\gamma}$-section of $\E$ can be written as some $\Co^{\gamma}$-linear combinations of $\xi_1,\dots,\xi_N$, which gives a based change formula
\begin{equation}\label{Eqn::DisInv::DisSecVB:Pf::Ten>0}
    \Co^{\gamma}_\loc(\Mf;\E)=\Co^{\gamma}_\loc(\Mf;\C)\otimes_{\Co^{\alpha}_\loc(\Mf;\C)}\Co^{\alpha}_\loc(\Mf;\E)=\Co^{\gamma}_\loc(\Mf)\otimes_{\Co^{\alpha}_\loc(\Mf)}\Co^{\alpha}_\loc(\Mf;\E).
\end{equation}

On the other hand, if $R$ is a commutative ring, $S$ is an $R$-commutative algebra, $P$ is a $S$-module and $Q$ is a $R$-module, we have
\begin{equation*}
    P\otimes_RQ=P\otimes_S(S\otimes_RQ).
\end{equation*}
Taking $R=\Co^{\alpha}_\loc(\Mf;\C)$, $S=\Co^{\gamma}_\loc(\Mf;\C)$, $P=\Co^{\beta-1}_\loc(\Mf;\C)$ and $Q=\Co^{\alpha}_\loc(\Mf;\E)$; or $R=\Co^{\alpha}_\loc(\Mf)$, $S=\Co^{\gamma}_\loc(\Mf)$, $P=\Co^{\beta-1}_\loc(\Mf)$ and $Q=\Co^{\alpha}_\loc(\Mf;\E)$, we get the right equality of \eqref{Eqn::DisInv::DisSecVB::Ten}.

To prove the left equality of \eqref{Eqn::DisInv::DisSecVB::Ten}, using \ref{Item::DisInv::DisSecVB::Loc} it suffices to show that for every $p\in\Mf$ there is a neighborhood $U\subseteq\Mf$ of $p$ such that $\Co^\beta_\loc(U;\E|_U)=\Co^\beta_\loc(U;\C)\otimes_{\Co^\alpha_\loc(U;\C)}\Co^{\alpha}_\loc(U;\E|_U)$. 

Indeed, by re-ordering the indices we can assume that $\xi_1(p),\dots,\xi_m(p)$ form a basis of $\E_p$. So $\xi_1,\dots,\xi_m$ for a local basis of some neighborhood $U\subseteq\Mf$ of $p$ and for each $m+1\le j\le N$ we can write $\xi_j=\sum_{k=1}^mf_j^k\xi_k$ for some $f_j^k\in\Co^\alpha_\loc(U)$.

By shrinking $U$ we can find a coordinate chart $\psi\in\Af$ of $\Mf$ whose domain is $U$. Let $\phi:\E|_U\to \C^m$ be given by $\phi(q,\sum_{j=1}^,v^j\xi_j):=(\psi(q),(v^1,\dots,v^m))$. We see that $(\psi,\phi)$ is a $\Co^\alpha$ bundle chart of $\E$.

Clearly for every $\eta\in\Co^\beta_\loc(U;\E|_U)$ we have $\eta|_U=\sum_{j=1}^m(\eta_{(\psi,\phi)}^j\circ\psi)\cdot\xi_j\in \Co^\beta_\loc(U;\C)\otimes_{\Co^\alpha_\loc(U;\C)}\Co^{\alpha}_\loc(U;\E|_U)$. Conversely for $\xi\in\Co^\beta_\loc(U;\C)\otimes_{\Co^\alpha_\loc(U;\C)}\Co^{\alpha}_\loc(U;\E|_U)$, we can write $\eta=\sum_{j=1}^Ng^j\xi_j$ where $g^j\in\Co^\beta_\loc(U;\C)$. By linear dependence of $\xi_{m+1},\dots,\xi_N$ we can write $\eta=\sum_{j=1}^N\sum_{k=1}^mg^jf_j^k\xi_k$. Therefore by definition of $(\psi,\phi)$ we have $\eta_{(\psi,\phi)}=\big(\sum_{j=1}^N\sum_{k=1}^m(g^jf_j^k)\circ\psi^\Inv\big)_{k=1}^m\in\Co^\beta_\loc(\psi(U);\C^m)$ i.e. $\eta\in \Co^\beta_\loc(U;\E|_U)$. This completes the proof.
\end{proof}
\begin{remark}\label{Rmk::DisInv::DisSubSec}
Let $\E\twoheadrightarrow\Mf$ be a $\Co^\alpha$ vector bundle and let $\V\le\E$ be a $\Co^\alpha$ subbundle. By viewing $\V$ itself as an abstract vector bundle we can define $\Co^\beta_\loc(\Mf;\V)$ for $\beta\in(\max(1-\kappa,-\alpha),\alpha]$ as well.

For pointwise defined section we can identify the collection $\{\xi_\Phi:\Phi\in\Bf^\E\}$ with a map $\xi:\Mf\to\E$. In this way we get the classical inclusion  $\Co^\alpha_\loc(\Mf;\V)\hookrightarrow \Co^\alpha_\loc(\Mf;\E)$. Therefore by Lemma \ref{Lem::DisInv::DisSecVB} \ref{Item::DisInv::DisSecVB::Ten}, we obtain $\Co^\beta_\loc(\Mf;\V)\hookrightarrow \Co^\beta_\loc(\Mf;\E)$ by the following:
\begin{equation}\label{Eqn::DisInv::DisSubSec}
    \Co^\beta_\loc(\Mf;\V)=\Co^\beta_\loc(\Mf;\C)\otimes_{\Co^\alpha_\loc(\Mf;\C)}\Co^{\alpha}_\loc(\Mf;\V)\hookrightarrow \Co^\beta_\loc(\Mf;\C)\otimes_{\Co^\alpha_\loc(\Mf;\C)}\Co^{\alpha}_\loc(\Mf;\E)=\Co^\beta_\loc(\Mf;\E).
\end{equation}
\end{remark}

Now we can give some properties of distributional sections of a vector subbundle. In the following if we take $\E=\C T\Mf$ then we get the characterizations of distributional sections of a tangent subbundle. 

\begin{prop}[Characterization of distributional sections]\label{Prop::DisInv::CharSec}
Let $\kappa>1$, $0\le\alpha\in(0,\kappa]$ and $\beta\in(\max(1-\kappa,-\alpha),\alpha]$, let $\Mf$ be an $n$-dimensional $\Co^\kappa$-manifold, let $\E$ be a $\Co^{\alpha}$ vector bundle over $\Mf$ and let $\V\le\E$ be a $\Co^{\alpha}$ vector bundle.
Let $X\in \Co^{\beta}_\loc(\Mf;\E)$. Let $\max(-\beta,0)<\gamma\le\alpha$ and $\max(1-\kappa,-\alpha)<\delta\le\beta$. The following are equivalent:
\begin{enumerate}[parsep=-0.3ex,label=(S.\arabic*)]
    \item\label{Item::DisInv::CharSec::Ten1} $X\in\Co^\beta_\loc(\Mf;\V)$ in the sense of \eqref{Eqn::DisInv::DisSubSec}. That is $X\in\Co^\beta_\loc(\Mf;\C)\otimes_{\Co^\alpha_\loc(\Mf;\C)}\Co^{\alpha}_\loc(\Mf;\V)$.
    \item\label{Item::DisInv::CharSec::Ten2}  $X\in \Co^\beta_\loc(\Mf;\C)\otimes_{\Co^{\gamma}_\loc(\Mf;\C)}\Co^{\gamma}_\loc(\Mf;\V)$.
    \item\label{Item::DisInv::CharSec::Pullback} $X\in \Co^{\delta}_\loc(\Mf;\V)$, or equivalently $X\in \Co^{\delta}_\loc(\Mf;\V)\cap \Co^{\beta}_\loc(\Mf;\E)$,  in the sense of \eqref{Eqn::DisInv::DisSubSec}.
    \item\label{Item::DisInv::CharSec::Dual} $\langle\theta,X\rangle=0\in\Co^\beta_\loc(\Mf;\C)$ for every dual section $\theta\in\Co^{\alpha}_\loc(\Mf;\V^\bot_\E)$. Here $\V^\bot_\E\le\E^*$ is the dual bundle of $\V$ in $\E$, given by 
    \begin{equation*}
        \V^\bot_\E:=\{(p,l)\in\E^*:\langle l,v\rangle=0\in\C,\ \forall v\in\V_p\}.
    \end{equation*}
\end{enumerate}
\end{prop}
\begin{remark}
By the equivalence of \ref{Item::DisInv::CharSec::Ten1} and \ref{Item::DisInv::CharSec::Ten2}, we see that definition of $\Co^{\beta}$-sections does not depend on the regularity assumption that $\V$ is $\Co^{\alpha}$. Also see Remark \ref{Rmk::DisInv::RmkDisSecVB} \ref{Item::DisInv::RmkDisSecVB::NoDepReg}.
\end{remark}
\begin{proof}[Proof of Proposition \ref{Prop::DisInv::CharSec}] The equivalence of \ref{Item::DisInv::CharSec::Ten1} and \ref{Item::DisInv::CharSec::Ten2} follows from Lemma \ref{Lem::DisInv::DisSecVB} \ref{Item::DisInv::DisSecVB::Ten}.
The equivalence of \ref{Item::DisInv::CharSec::Ten1} and \ref{Item::DisInv::CharSec::Pullback} comes from the standard computation:
\begin{align*}
    \Co^\delta_\loc(\Mf;\V)\cap\Co^\beta_\loc(\Mf;\C)=&(\Co^\delta_\loc(\Mf;\C)\otimes_{\Co^\alpha_\loc(\Mf;\C)}\Co^\alpha_\loc(\Mf;\V))\cap(\Co^\beta_\loc(\Mf;\C)\otimes_{\Co^\alpha_\loc(\Mf;\C)}\Co^\alpha_\loc(\Mf;\E))
    \\
    =&\Co^\beta_\loc(\Mf;\C)\otimes_{\Co^\alpha_\loc(\Mf;\C)}\Co^\alpha_\loc(\Mf;\V)=\Co^\beta_\loc(\Mf;\V).
\end{align*}

To see \ref{Item::DisInv::CharSec::Pullback} $\Rightarrow$ \ref{Item::DisInv::CharSec::Dual}, we write $X=\sum_{i=1}^Nf^iY_i$ where $f^i\in\Co^\delta_\loc(\Mf;\C)$ and $Y_i\in\Co^\alpha_\loc(\Mf;\C)$, then
\begin{equation*}
    \langle \theta,X\rangle=\sum_{i=1}^N\langle\theta,f^iY_i\rangle=\sum_{i=1}^Nf^i\langle\theta,Y_i\rangle=0,\quad\text{for every }\Co^{\alpha}\text{ 1-form }\theta\text{ on }\V^\bot.
\end{equation*}

It remains to show  \ref{Item::DisInv::CharSec::Dual} $\Rightarrow$ \ref{Item::DisInv::CharSec::Ten1}. 

\medskip
By Lemma \ref{Lem::DisInv::DisSecVB} \ref{Item::DisInv::DisSecVB::FiniteGen} we can find finitely many $\Co^{\alpha}$ sections $Y_1,\dots,Y_N$ of $\E$ that generates the module  $\Co^{\alpha}_\loc(\Mf;\V)$ over the ring $\Co^{\alpha}_\loc(\Mf;\C)$. In particular $Y_1(p),\dots,Y_N(p)$ spans $\V_p\le\E_p$ for all $p\in \Mf$.

We can find an open cover $\{U_\zeta\}_\zeta$, associated with $r$-element subsets $I_\zeta\subseteq\{1,\dots,N\}$ and vector fields $Z^\zeta_1,\dots,Z^\zeta_{n-r}\in \Co^{\alpha}_\loc(U_\zeta)$ for each $\zeta$, such that $\{Y_i\}_{i\in I}\cup\{Z^\zeta_1,\dots,Z^\zeta_{n-r}\}$ form a local basis for $\E|_{U_\zeta}$. So for each $\zeta$, $(Y_i)_{i\in I_\zeta}$ are linearly independent and span $\V|_{U_\zeta}$.

Let $(\mu_\zeta^1,\dots,\mu_\zeta^r,\theta_\zeta^1,\dots,\theta_\zeta^{n-r})$ be  $\Co^{\alpha}$ 1-forms defined on $U_\zeta$ that form the dual basis to vector fields $(Y_{\sigma_\zeta(1)},\dots,Y_{\sigma_\zeta(r)},Z^\zeta_1,\dots,Z^\zeta_{n-r})$ where we take $\sigma_\zeta:\{1,\dots,r\}\to I_\zeta$ such that $\sigma_\zeta(1)<\dots<\sigma_\zeta(r)$. Therefore every $\tilde X\in \Co^{\beta}_\loc(U_\zeta;\E|_{U_\zeta})$ can be written as
\begin{equation*}
    \tilde X=\sum_{j=1}^r\langle\mu_\zeta^j,\tilde X\rangle Y_{\sigma_\zeta(j)}+\sum_{k=1}^{n-r}\langle\theta_\zeta^k,\tilde X\rangle Z^\zeta_k,\quad\text{on }U_\zeta\text{ for each }\zeta.
\end{equation*}

Let $\{\chi_\zeta\in\Co^\kappa_c(U_\zeta)\}_\zeta$ be a $\Co^\kappa$ partition of unity to the open cover $\{U_\zeta\}_\zeta$, thus we have
\begin{equation*}
    X=\sum_{\zeta}\Big(\sum_{j=1}^r\langle\chi_\zeta\mu_\zeta^j, X\rangle\cdot Y_{\sigma_\zeta(j)}+\sum_{k=1}^{n-r}\langle\chi_\zeta\theta_\zeta^k,X\rangle \cdot Z^\zeta_k\Big).
\end{equation*}
Clearly $\chi_\zeta\theta_\zeta^1,\dots,\chi_\zeta\theta_\zeta^{n-r}\in \Co^{\alpha}(\Mf;\V^\bot)$ are globally defined, since $\theta_\zeta^1,\dots,\theta_\zeta^{n-r}$ are $\Co^{\alpha}$-sections of $V^\bot|_{U_\zeta}$.

By assumption $X\in \Co^{\beta}_\loc(\Mf;\X)\cap C^{-1,\delta}_\loc(\Mf;\V)$, we have $\langle\chi_\zeta\theta_\zeta^k,X\rangle=0 \in C^{-1,\delta}_\loc(\Mf)$, so 
\begin{equation}\label{Eqn::DisInv::CharSecPf::Tmp}
    X=\sum_{\zeta}\sum_{j=1}^r\langle\chi_\zeta\mu_\zeta^j, X\rangle\cdot Y_{\sigma_\zeta(j)})=\sum_{i=1}^N\Big(\sum_\zeta\langle\chi_\zeta\mu_\zeta^{\sigma_\zeta^{-1}(j)}, X\rangle\Big)\cdot Y_j.
\end{equation}
Here we let $\mu^{\sigma_\zeta^{-1}(j)}_\zeta=0$ if $j\notin I_\zeta(=\sigma_\zeta(\{1,\dots,r\}))$.

Since $\langle\chi_\zeta\mu_\zeta^j, X\rangle\in \Co^{\beta}_\loc(\Mf;\C)$ and the sums in \eqref{Eqn::DisInv::CharSecPf::Tmp} are locally finite, we know $X$ is the linear combinations of $Y_1,\dots,Y_N\in \Co^{\alpha}_\loc(\Mf;\V)$ with coefficients in $\Co^{\beta}_\loc(\Mf;\C)$. Therefore $X\in \Co^{\beta}_\loc(\Mf;\C)\otimes_{\Co^{\alpha}_\loc(\Mf;\C)}\Co^{\alpha}_\loc(\Mf;\V)$, which is the condition \ref{Item::DisInv::CharSec::Ten1}. 
\end{proof}

\section{Distributional involutivity for regular subbundles}\label{Section::DisInv::CharInv}
In this part we provide some equivalent characterizations of distributional involutivity for a tangent subbundle. Recall the definition of  distributional involutivity in Definition \ref{Defn::Intro::DisInv}. See also in Section \ref{Section::RealFro::HisRmk} for some historical remark.

Recall that in the smooth setting, when $\Se\le\C T\Mf$ is a smooth complex tangent subbundle over a smooth manifold $\Mf$, we have the following equivalent characterizations of involutivity:
\begin{enumerate}[parsep=-0.3ex,label=(\ref{Chapter::DisInv}.\Alph*)]
    \item\label{Item::DisInv::SmoothChar::LieAlg} $C^\infty_\loc(\Mf;\Se)$ is a Lie subalgebra of $C^\infty_\loc(\Mf;\C T\Mf)$, i.e. if $X$ and $Y$ are smooth sections of $\Se$ then $[X,Y]$ is also a smooth section of $\Se$.
    \item\label{Item::DisInv::SmoothChar::LieDer}  If $X\in C^\infty_\loc(\Mf;\Se)$ and $\theta\in C^\infty_\loc(\Mf;\Se^\bot )$, then the Lie derivative $\Lie X\theta\in C^\infty_\loc(\Mf;\Se^\bot)$.
    \item\label{Item::DisInv::SmoothChar::DiffIdeal}  On the exterior algebra $C^\infty_\loc(\Mf;\wedge^\bullet)=\bigoplus_{k=0}^{\dim \Mf}C^\infty_\loc(\Mf;\wedge^k\C T^*\Mf)$, the (two-sided) ideal $C^\infty_\loc(\Mf;\Se^\bot)\wedge C^\infty_\loc(\Mf;\wedge^\bullet \C T^*\Mf)$ is closed under differential.
    \item\label{Item::DisInv::SmoothChar::2Form} On space of 2-forms we have $d(C^\infty_\loc(\Mf;\Se^\bot))\subseteq C^\infty_\loc(\Mf;\Se^\bot)\wedge C^\infty_\loc(\Mf;\C T^*\Mf)$.
\end{enumerate}

Here we recall that $\langle\Lie X\theta,Y\rangle=X\langle\theta,Y\rangle-\langle\theta,[X,Y]\rangle$ for all smooth vector fields $X,Y$ and 1-form $\theta$. The equality is still true for rough vector fields and 1-forms as long as the products are defined.

These characterizations can all be generalized to non-smooth subbundles.

We first describe some of the characterizations involving Lie brackets. Recall the involutivity definition of modules in Definition \ref{Defn::DisInv::ModInv}.


\begin{prop}[Characterizations of distributional involutivity I]\label{Prop::DisInv::CharInv1}
Let $\alpha\in(\frac12,\infty)$, let $\Mf$ be a $n$-dimensional $\Co^{\alpha+1}$-manifold and let $\Se$ be a rank $r$ $\Co^{\alpha}$-subbundle of $\C T\Mf$. 

Let $\beta\in(\frac12,\alpha]$ and $\gamma\in(-\beta,\beta-1]$, the following conditions for $\Se$ are equivalent:
\begin{enumerate}[parsep=-0.3ex,label=(I.\arabic*)]
    \item\label{Item::DisInv::CharInv1::Sec1} $\Co^\alpha_\loc(\Mf;\Se)$ is $\Co^{\alpha-1}$-involutive. That is, for every $X,Y\in \Co^{\alpha}_\loc(\Mf;\Se)$ we have $[X,Y]\in \Co^{\alpha-1}_\loc(\Mf;\Se)$.
    \item\label{Item::DisInv::CharInv1::Sec2} $\Co^\beta_\loc(\Mf;\Se)$ is $\Co^\gamma$-involutive.
    \item\label{Item::DisInv::CharInv1::Pair0}$\Se$ is distributional involutive (by Definition \ref{Defn::Intro::DisInv}). That is, for every $X,Y\in \Co^{\frac12+}_\loc(\Mf;\Se)$ and $\theta\in \Co^{\frac12+}_\loc(\Mf;\Se^\bot)$, we have $\langle\theta,[X,Y]\rangle=0\in \Co^{(-\frac12)+}_\loc(\Mf;\C)$.
    \item\label{Item::DisInv::CharInv1::Pair1} For every $X,Y\in \Co^{\alpha}_\loc(\Mf;\Se)$ and $\theta\in \Co^{\alpha}_\loc(\Mf;\Se^\bot)$, we have $\langle\theta,[X,Y]\rangle=0\in \Co^{\alpha-1}_\loc(\Mf;\C)$.
    \item\label{Item::DisInv::CharInv1::Pair2} For every $X,Y\in \Co^{\alpha}_\loc(\Mf;\Se)$ and $\theta\in \Co^{\alpha}_\loc(\Mf;\Se^\bot)$, we have $d\theta(X,Y)=0\in \Co^{\alpha-1}_\loc(\Mf;\C)$.
    \item\label{Item::DisInv::CharInv1::Pair3} For every $X\in \Co^{\alpha}_\loc(\Mf;\Se) $ and $\theta\in \Co^{\alpha}_\loc(\Mf;\Se^\bot)$, we have $\Lie X\theta\in \Co^{\alpha-1}_\loc(\Mf;\Se^\bot)$.
    \item\label{Item::DisInv::CharInv1::Gen1} For every $p\in \Mf$ there is a neighborhood $U\subseteq \Mf$ of $p$ and a $\Co^{\alpha}$-local basis $(X_1,\dots,X_r)$ for $\Se|_U$ such that $X_1,\dots,X_r$ are pairwise commutative in the sense of distributions.
    \item\label{Item::DisInv::CharInv1::Gen2} If $X_1,\dots,X_{N}\in \Co^{\alpha}_\loc(\Mf;\Se)$ are global generators for $\Se$ (that is $X_1,\dots,X_N$ span $\Se$ at every point in $\Mf$), then there are $c_{ij}^k\in \Co^{\alpha-1}_\loc(\Mf;\C)$ for $1\le i,j,k\le N$ such that 
    \begin{equation}\label{Eqn::DisInv::CharInv1::Gen2cijk}
        [X_i,X_j]=\sum_{k=1}^Nc_{ij}^kX_k\quad\text{on }\Mf,\quad \forall\ 1\le i,j\le N.
    \end{equation}
\end{enumerate}
\end{prop}
\begin{remark}
\begin{enumerate}[parsep=-0.3ex,label=(\alph*)]
    \item In the statements, \ref{Item::DisInv::CharInv1::Sec1} and \ref{Item::DisInv::CharInv1::Pair3} are the direct generalizations to \ref{Item::DisInv::SmoothChar::LieAlg} and \ref{Item::DisInv::SmoothChar::LieDer} respectively.
    \item By the equivalence of \ref{Item::DisInv::CharInv1::Sec1} and \ref{Item::DisInv::CharInv1::Sec2} or \ref{Item::DisInv::CharInv1::Pair0} and \ref{Item::DisInv::CharInv1::Pair1}, we see that the definition of distributional involutivity does not depend on the regularity assumption $\Se\in \Co^{\alpha}$.
    
    \smallskip
    In this way the characterizations are still equivalent if we replace every $\alpha$ by $\beta$ in \ref{Item::DisInv::CharInv1::Pair1} - \ref{Item::DisInv::CharInv1::Gen2}.

    \item The above characterizations  do not use the information of $\dim \Mf$ and $\rank \Se$. That is why we can extend Definition \ref{Defn::Intro::DisInv}, the involutivity of a subset (subbundle) of $\Se$, to Definition \ref{Defn::DisInv::ModInv}, the involutivity of a submodule $\Co^\alpha_\loc(\Mf;\Se)$.
    
    \item In particular from \ref{Item::DisInv::CharInv1::Gen2} we can focus on the specific choice of $X_1,\dots,X_N$ rather than on the ambient space $\Se$. In the setting of regular subbundles, we only consider the module spanned by $X_1,\dots,X_N$, where $\Span(X_1(p),\dots,X_N(p))\le\C T_p\Mf$ have constant rank among all $p\in\Mf$. 
    
\end{enumerate}
\end{remark}

\begin{proof}[Proof of Proposition \ref{Prop::DisInv::CharInv1}]We consider the directions \ref{Item::DisInv::CharInv1::Sec2} $\Rightarrow$ \ref{Item::DisInv::CharInv1::Sec1} $\Leftrightarrow$ \ref{Item::DisInv::CharInv1::Pair1} $\Leftrightarrow$ \ref{Item::DisInv::CharInv1::Pair2} $\Leftrightarrow$ \ref{Item::DisInv::CharInv1::Pair3}; \ref{Item::DisInv::CharInv1::Pair1} $\Rightarrow$ \ref{Item::DisInv::CharInv1::Gen1} $\Rightarrow$ \ref{Item::DisInv::CharInv1::Sec2}, \ref{Item::DisInv::CharInv1::Sec1} $\Leftrightarrow$ \ref{Item::DisInv::CharInv1::Gen2}, and \ref{Item::DisInv::CharInv1::Sec1} $\Leftrightarrow$ \ref{Item::DisInv::CharInv1::Pair0}.

When $X,Y\in \Co^{\alpha}_\loc(\Mf;\C T\Mf)$ and $\theta\in \Co^{\alpha}_\loc(\Mf;\C T^*\Mf)$, by Lemma \ref{Lem::Hold::MultLoc} \ref{Item::Hold::MultLoc::WellDef} we have $[X,Y]\in \Co^{\alpha-1}_\loc(\Mf;\C T\Mf)$, $\langle\theta,[X,Y]\rangle,d\theta(X,Y)\in \Co^{\alpha-1}_\loc(\Mf;\C)$ and $\Lie X\theta\in \Co^{\alpha-1}_\loc(\Mf;T^*\Mf)$.

\noindent\ref{Item::DisInv::CharInv1::Sec2} $\Rightarrow$ \ref{Item::DisInv::CharInv1::Sec1}: Let $X,Y\in\Co^\alpha_\loc(\Mf;\Se)$. By condition \ref{Item::DisInv::CharInv1::Sec2} $[X,Y]\in\Co^\gamma_\loc(\Mf;\Se)$. Since we have $[X,Y]\in\Co^{\alpha-1}_\loc(\Mf;\C T\Mf)$, by Proposition \ref{Prop::DisInv::CharSec} \ref{Item::DisInv::CharSec::Pullback} $\Rightarrow$ \ref{Item::DisInv::CharSec::Ten1}, $[X,Y]\in\Co^{\gamma}_\loc(\Mf;\Se)\cap \Co^{\alpha-1}_\loc(\Mf;\C T\Mf)=\Co^{\alpha-1}_\loc(\Mf;\Se)$.

The direction \ref{Item::DisInv::CharInv1::Sec1} $\Leftrightarrow$ \ref{Item::DisInv::CharInv1::Pair1} follows from Proposition \ref{Prop::DisInv::CharSec} \ref{Item::DisInv::CharSec::Ten1} $\Leftrightarrow$ \ref{Item::DisInv::CharSec::Dual}, where $[X,Y]$ is a section of $\Se$ if and only if $\langle\theta,[X,Y]\rangle=0$ for all $\theta\in \Co^{\alpha}_\loc(\Mf;\Se^\bot)$.

\medskip
\noindent\ref{Item::DisInv::CharInv1::Pair1} $\Leftrightarrow$ \ref{Item::DisInv::CharInv1::Pair2} $\Leftrightarrow$ \ref{Item::DisInv::CharInv1::Pair3}: For $X,Y\in \Co^{\alpha}_\loc(\Mf;\Se)$ and $\theta\in \Co^{\alpha}_\loc(\Mf;\Se^\bot)$, in the sense of distributions we have $d\theta(X,Y)=X\langle \theta, Y\rangle-Y\langle\theta,X\rangle-\langle\theta,[X,Y]\rangle=-\langle\theta,[X,Y]\rangle$ and $\langle\Lie X\theta,Y\rangle=X\langle\theta,Y\rangle-\langle\theta,[X,Y]\rangle$. Thus
\begin{equation}\label{Eqn::DisInv::CharInv1::Tmp}
    \langle\theta,[X,Y]\rangle=-d\theta(X,Y)=\langle\Lie X\theta,Y\rangle.
\end{equation}
If one term in \eqref{Eqn::DisInv::CharInv1::Tmp} is zero then the other two vanish as well. Thus we get \ref{Item::DisInv::CharInv1::Pair1} $\Leftrightarrow$ \ref{Item::DisInv::CharInv1::Pair2}.
Note that \ref{Item::DisInv::CharInv1::Pair3} is equivalent as saying $\langle\Lie X\theta,Y\rangle=0$ for all $X,Y\in \Co^{\alpha}_\loc(\Mf;\Se)$ and $\theta\in \Co^{\alpha}_\loc(\Mf;\Se^\bot)$, thus \ref{Item::DisInv::CharInv1::Pair3} is equivalent to  \ref{Item::DisInv::CharInv1::Pair1} and \ref{Item::DisInv::CharInv1::Pair2}.
    
    \medskip
\noindent \ref{Item::DisInv::CharInv1::Pair1} $\Rightarrow$ \ref{Item::DisInv::CharInv1::Gen1}: This is done by Lemma \ref{Lem::ODE::GoodGen} \ref{Item::ODE::GoodGen::InvComm}.

\medskip
\noindent \ref{Item::DisInv::CharInv1::Gen1} $\Rightarrow$ \ref{Item::DisInv::CharInv1::Sec2}: Let $Y_1,Y_2\in \Co^{\beta}_\loc(\Mf;\Se)$, by locality of distributions (see Lemma \ref{Lem::DisInv::DisSecVB} \ref{Item::DisInv::DisSecVB::Loc}) it suffices to show that for any $p\in \Mf$ there is a neighborhood $U\subseteq \Mf$ of $p$ such that $[Y_1,Y_2]\in \Co^{\beta-1}_\loc(U;\Se|_U)$. Fix $p\in \Mf$, let $U\subseteq \Mf$ and  $X_1,\dots,X_r\in \Co^{\alpha}_\loc(U;\Se|_U)$ be as in \ref{Item::DisInv::CharInv1::Gen1}. By Proposition \ref{Prop::DisInv::CharSec} \ref{Item::DisInv::CharSec::Dual} $\Rightarrow$ \ref{Item::DisInv::CharSec::Ten1} we can write $Y_i=\sum_{j=1}^rf_i^jX_j$ for some $f_i^j\in \Co^{\beta-1}_\loc(U)$, $i=1,2$ and $1\le j\le r$. Therefore
\begin{equation*}
    [Y_1,Y_2]=\sum_{i,j=1}^r[f_1^iX_i,f_2^jX_j]=\sum_{i,j=1}^rf_1^i(X_if_2^j)X_j-f_2^j(X_jf_1^i)X_i+f_1^if_2^j[X_i,X_j]=\sum_{i,j=1}^r\big(f_1^j(X_jf_2^i)-f_2^j(X_jf_1^i)\big)X_i.
\end{equation*}

We see that $[Y_1,Y_2]$ is the $\Co^{\beta-1}_\loc$ linear combinations of $X_1,\dots,X_r$, thus $[Y_1,Y_2]\in \Co^{\beta-1}_\loc(U)\otimes_{\Co^{\alpha}_\loc(U)}\Co^{\alpha}(U;\Se|_U)$. By Proposition \ref{Prop::DisInv::CharSec} \ref{Item::DisInv::CharSec::Ten1}$\Rightarrow$\ref{Item::DisInv::CharSec::Dual} $[Y_1,Y_2]\in \Co^{\beta-1}_\loc(U;\Se|_U)$. Since $p\in \Mf$ is arbitrary, we get $[Y_1,Y_2]\in \Co^{\beta-1}_\loc(\Mf;\Se)(\subseteq\Co^\gamma_\loc(\Mf;\Se))$, proving \ref{Item::DisInv::CharInv1::Sec2}.

\medskip The proof of \ref{Item::DisInv::CharInv1::Gen2} $\Rightarrow$ \ref{Item::DisInv::CharInv1::Sec1} is identical to the proof of \ref{Item::DisInv::CharInv1::Gen1} $\Rightarrow$ \ref{Item::DisInv::CharInv1::Sec2} except we do not need to localize the objects on an open subset. We omit the details to the readers.

\medskip
\noindent \ref{Item::DisInv::CharInv1::Sec1} $\Rightarrow$ \ref{Item::DisInv::CharInv1::Gen2}: Note that by Lemma \ref{Lem::DisInv::DisSecVB} \ref{Item::DisInv::DisSecVB::FiniteGen} we can find finitely many $X_1,\dots,X_N\in \Co^{\alpha}_\loc(\Mf;\Se)$ that spans $\Se$ at every point.

By assumption \ref{Item::DisInv::CharInv1::Sec1}, $[X_i,X_j]\in \Co^{\alpha-1}(\Mf;\Se)$ for all $1\le i,j\le N$.
So by Proposition \ref{Prop::DisInv::CharSec} \ref{Item::DisInv::CharSec::Dual}$\Rightarrow$\ref{Item::DisInv::CharSec::Ten1}, $[X_i,X_j]\in \Co^{\alpha-1}_\loc(\Mf)\otimes_{\Co^{\alpha}_\loc(\Mf)}\Co^{\alpha}(\Mf;\Se)$.
Since $X_1,\dots,X_N$ generates $\Co^{\alpha}(\Mf;\Se)$ as a module over $\Co^{\alpha}_\loc(\Mf)$, we can write $[X_i,X_j]=\sum_{k=1}^Nc_{ij}^kX_k$ for some $c_{ij}^k\in \Co^{\alpha-1}_\loc(\Mf)$. This gives the condition \ref{Item::DisInv::CharInv1::Gen2}.

 Now we have all the done all the equivalence of each conditions except \ref{Item::DisInv::CharInv1::Pair0}.
 
\medskip \noindent \ref{Item::DisInv::CharInv1::Sec1} $\Leftrightarrow$ \ref{Item::DisInv::CharInv1::Pair0}: Clearly \ref{Item::DisInv::CharInv1::Pair0} implies \ref{Item::DisInv::CharInv1::Pair1}. We proved \ref{Item::DisInv::CharInv1::Pair1} $\Rightarrow$ \ref{Item::DisInv::CharInv1::Sec1} so we have \ref{Item::DisInv::CharInv1::Pair0} $\Rightarrow$ \ref{Item::DisInv::CharInv1::Sec1}.

Conversely, for $X,Y,\theta\in\Co^{\frac12+}$ in the assumption of \ref{Item::DisInv::CharInv1::Pair0}, we can assume that $X,Y,\theta\in\Co^{\frac12+\eps_0}$ for some $\eps_0>0$. Using \ref{Item::DisInv::CharInv1::Sec1} $\Leftrightarrow$ \ref{Item::DisInv::CharInv1::Sec2} for arbitrary $\frac12<\beta\le\alpha$ and $-\beta<\gamma\le\beta-1$, we see that \ref{Item::DisInv::CharInv1::Sec1} implies that $\Co^{\frac12+\eps}_\loc(\Mf;\Se)$ is $\Co^{\eps-\frac12}$-involutive, for every $\eps>0$. Therefore $[X,Y]\in\Co^{\eps_0-\frac12}_\loc(\Mf;\Se)$.  Thus by Proposition \ref{Prop::DisInv::CharSec} \ref{Item::DisInv::CharSec::Ten1} $\Rightarrow$ \ref{Item::DisInv::CharSec::Dual} with $\beta=\eps_0-\frac12$, we see that $\langle\theta,[X,Y]\rangle=0\in\Co^{\frac12-\eps_0}_\loc(\Mf;\C)$. This gives the condition \ref{Item::DisInv::CharInv1::Pair0} and completes the proof.
\end{proof}

In Proposition \ref{Prop::DisInv::CharInv1} we mostly focus on the side of vector fields. We can consider the characterizations using differential forms which are related to the sections of $\Se^\bot$.

\begin{prop}[Characterizations of distributional involutivity II]\label{Prop::DisInv::CharInv2}
Let $\alpha\in(\frac12,\infty)$, let $\Mf$ be a $n$-dimensional $\Co^{\alpha+1}$-manifold and let $\Se\le\C T\Mf$ be a rank $r$ $\Co^{\alpha}$-subbundle. 

Then $\Se$ is distributional involutive if and only if  either of the following condition holds:
\begin{enumerate}[parsep=-0.3ex,label=(I.\arabic*)]\setcounter{enumi}{8}
    
        \item\label{Item::DisInv::CharInv2::Ann} Let $1\le k\le n$ and $\omega\in \Co^{\alpha}_\loc(\Mf;\wedge^k\C T^*\Mf)$ satisfy $\omega(Y_1,\dots,Y_k)=0$ for every $Y_1,\dots,Y_k\in \Co^{\alpha}_\loc(\Mf;\Se)$. Then we have $d\omega(Y_1,\dots,Y_{k+1})=0$ as distributions for every $Y_1,\dots,Y_{k+1}\in \Co^{\alpha}_\loc(\Mf;\Se)$.   
    
    \item\label{Item::DisInv::CharInv2::Diff} $d(\Co^{\alpha}_\loc(\Mf;\Se^\bot))\subseteq (\Co^{\alpha}_\loc(\Mf;\Se^\bot)\wedge \Co^{\alpha}_\loc(\Mf;\C T^*\Mf))\otimes_{\Co^{\alpha}_\loc(\Mf;\C)}\Co^{\alpha-1}_\loc(\Mf;\C)$. That is, for every $\lambda\in \Co^{\alpha}_\loc(\Mf;\Se^\bot)$, $d\lambda\in \Co^{\alpha-1}_\loc(\Mf;\wedge^2\C T^*\Mf)$ can be written as a finite sum
    \begin{equation*}
        d\lambda=\sum_{i=1}^Nf_i\cdot\theta^i\wedge\mu^i,\quad\text{for some }N\ge1,\quad f_i\in \Co^{\alpha-1}_\loc(\Mf;\C),\quad \theta^i\in \Co^{\alpha}_\loc(\Mf;\Se^\bot),\quad \mu^i\in \Co^{\alpha}_\loc(\Mf;\C T^*\Mf).
    \end{equation*}

    \item\label{Item::DisInv::CharInv2::coSec} For every $\lambda^1,\dots,\lambda^{n-r+1}\in \Co^{\alpha}_\loc(\Mf;\Se^\bot)$ we have 
    \begin{equation}\label{Eqn::DisInv::CharInv2::coSecEqn}
        \lambda^1\wedge\dots\wedge\lambda^{n-r}\wedge d\lambda^{n-r+1}=0\in \Co^{\alpha-1}_\loc(\Mf;\wedge^{n-r+2}\C T^*\Mf).
    \end{equation}
\end{enumerate}
\end{prop}

Here \ref{Item::DisInv::CharInv2::Diff} is the direct generalization to \ref{Item::DisInv::SmoothChar::2Form}.

We shall see that \ref{Item::DisInv::CharInv2::Ann} generalizes \ref{Item::DisInv::SmoothChar::DiffIdeal}, because the space of all differential forms that annihilate sections of $\Se$ is indeed an ideal of the exterior algebra $\Co^{\alpha}_\loc(\Mf;\wedge^\bullet)$.

To prove the equivalence of \ref{Item::DisInv::CharInv2::Diff} and \ref{Item::DisInv::CharInv2::Ann} we need the following:

\begin{lem}\label{Lem::DisInv::IdealChar}
    Let $\alpha>0$ and let $\Se\le\C T\Mf$ be a $\Co^\alpha$ complex subbundle. Then for $1\le k\le n$,
    \begin{equation}\label{Eqn::DisInv::IdealChar::Eqn}
        \Co^{\alpha}_\loc(\Mf;\Se^\bot)\wedge \Co^{\alpha}_\loc(\Mf;\wedge^{k-1}\C T^*\Mf)=\{\omega\in \Co^{\alpha}_\loc(\Mf;\Lambda^k):\omega(X_1,\dots,X_k)=0\ \forall X_1,\dots,X_k\in \Co^{\alpha}_\loc(\Mf;\Se)\}.
    \end{equation}
\end{lem}
\begin{proof}

Let $W$ be a $n$-dimensional (complex) vector space and $V\le W$ is a rank $r$ (complex) subspace. By results in standard linear algebra, as the subspaces of alternating tensors we have
\begin{equation}\label{Eqn::DisInv::IdealChar::PointwiseEqn}
    V_W^\bot\wedge \textstyle\bigwedge^{k-1}W^*=\big\{\phi\in\textstyle\bigwedge^{k}W^*:\phi(v_1,\dots,v_k)=0\text{ for all }v_1,\dots,v_k\in V\big\}.
\end{equation}
Here $V_W^\bot=\{l\in W^*:\langle l,v\rangle=0,\ \forall v\in V\}$.

For each $p\in \Mf$, by taking $W=\C T_p\Mf$ and $V=\Se_p$ in \eqref{Eqn::DisInv::IdealChar::PointwiseEqn}, we get \eqref{Eqn::DisInv::IdealChar::Eqn} and complete the proof.
\end{proof}

\begin{proof}[Proof of Proposition \ref{Prop::DisInv::CharInv2}]

By Lemma \ref{Lem::DisInv::IdealChar} we see that 
\ref{Item::DisInv::CharInv2::Diff} is the special case to \ref{Item::DisInv::CharInv2::Ann} when $k=1$. So \ref{Item::DisInv::CharInv2::Ann} $\Rightarrow$ \ref{Item::DisInv::CharInv2::Diff}.


\medskip
We are going to prove \ref{Item::DisInv::CharInv1::Gen2} $\Rightarrow$  \ref{Item::DisInv::CharInv2::Ann} and \ref{Item::DisInv::CharInv2::Diff} $\Rightarrow$ \ref{Item::DisInv::CharInv2::coSec} $\Rightarrow$ \ref{Item::DisInv::CharInv1::Gen1}. This would complete the proof since by Proposition \ref{Prop::DisInv::CharInv1}, \ref{Item::DisInv::CharInv1::Gen2} and \ref{Item::DisInv::CharInv1::Gen1} are both equivalent characterizations of involutivity.

\medskip
\noindent\ref{Item::DisInv::CharInv1::Gen2} $\Rightarrow$ \ref{Item::DisInv::CharInv2::Ann}: Let $X_1,\dots,X_N$ satisfy the condition \ref{Item::DisInv::CharInv1::Gen2}, and let $\omega$ be a $k$-form annihilating all sections of $\Se$. In particular $\omega(X_{i_1},\dots,X_{i_k})=0$ for all $\{i_1,\dots,i_k\}\subseteq\{1,\dots,N\}$. By passing to $\Co^{\alpha}$-linear combinations it suffices to show that $\omega(Y_1,\dots,Y_{k+1})=0$ for $\{Y_1,\dots,Y_{k+1}\}\subseteq\{X_1,\dots,X_N\}$. 

Indeed by direct computation along with \eqref{Eqn::DisInv::CharInv1::Gen2cijk} we have
\begin{align*}
    d\omega(X_{i_0},\dots,X_{i_{k}})=&\sum_{0\le u<v\le k}(-1)^{u+v}\omega([X_{i_u},X_{i_v}],X_{i_0},\dots , X_{i_{u-1}},X_{i_{u+1}}\dots ,X_{i_{v-1}},X_{i_{v+1}},\dots ,X_{i_k})
    \\
    =&\sum_{0\le u<v\le k}(-1)^{u+v}\sum_{l=1}^Nc_{i_ui_v}^l\cdot\omega(X_l,X_{i_0},\dots , X_{i_{u-1}},X_{i_{u+1}}\dots ,X_{i_{v-1}},X_{i_{v+1}},\dots ,X_{i_k})=0.
\end{align*}
This gives \ref{Item::DisInv::CharInv2::Ann}.

\medskip
\noindent\ref{Item::DisInv::CharInv2::Diff} $\Rightarrow$ \ref{Item::DisInv::CharInv2::coSec}:
Let $\lambda^1,\dots,\lambda^{n-r+1}\in \Co^{\alpha}_\loc(\Mf;\Se^\bot)$, by assumption \ref{Item::DisInv::CharInv2::Diff} we can write $d\lambda^{n-r+1}=\sum_{i=1}^Nf_i\theta^i\wedge\mu^i$ for some $\theta^i\in \Co^{\alpha}_\loc(\Mf;\Se^\bot)$. So
\begin{equation*}
    \lambda^1\wedge\dots\wedge\lambda^{n-r}\wedge d\lambda^{n-r+1}=\sum_{i=1}^N \lambda^1\wedge\dots\wedge\lambda^{n-r}\wedge \theta^i\wedge (f_i\mu^i)
\end{equation*}

On the other hand $\Se^\bot $ is a cotangent subbundle with rank $(n-r)$, we see that $\lambda^1,\dots,\lambda^{n-r},\theta^i$ are linear dependent at every point in $\Mf$, thus $\lambda^1\wedge\dots\wedge\lambda^{n-r}\wedge \theta^i\equiv0$. We conclude that $\lambda^1\wedge\dots\wedge\lambda^{n-r}\wedge d\lambda^{n-r+1}=0$ as a distributional differential form.

\medskip
\noindent\ref{Item::DisInv::CharInv2::coSec} $\Rightarrow$ \ref{Item::DisInv::CharInv1::Gen1}: Fixed a point $p\in \Mf$, by Lemma \ref{Lem::ODE::GoodGen} we can find a $\Co^{\alpha+1}$-chart $(x^1,\dots,x^r,y^1,\dots,y^{n-r}):U\subseteq \Mf\to\R^n$ near $p$ and $\Co^{\alpha}$-vector fields $X_1,\dots,X_r$ on $U$ that span $\Se|_U$ and have the form $X_j=\Coorvec{x^j}+\sum_{k=1}^{n-r}b_j^k\Coorvec{y^k}$. What we need is to show that $[X_j,X_k]=0$ for all $1\le j<k\le r$.

Clearly $(X_1,\dots,X_r,\Coorvec{y^1},\dots,\Coorvec{y^{n-r}})$ span the tangent space at every point in $U$. The dual basis has the form $(dx^1,\dots,dx^r,\theta^1,\dots,\theta^{n-r})$ where $\theta^l=dy^l-\sum_{j=1}^rb_j^ldx^j$, $l=1,\dots,n-r$. Thus for $1\le l\le n-r$
\begin{align*}
    d\theta^l&=\sum_{i,j=1}^r\frac{\partial b_i^l}{\partial x^j}dx^i\wedge dx^j+\sum_{i=1}^r\sum_{k=1}^{n-r}\frac{\partial b_i^l}{\partial y^k}dx^i\wedge dy^k
    \\
    &=\sum_{1\le i<j\le r}\Big(\frac{\partial b_i^l}{\partial x^j}-\frac{\partial b_j^l}{\partial x^i}+\sum_{k=1}^{n-r}\big(b_j^k\frac{\partial b_i^l}{\partial y^k}-b_i^k\frac{\partial b_i^k}{\partial y^l}\big)\Big)dx^i\wedge dx^j+\sum_{i=1}^r\sum_{k=1}^{n-r}\frac{\partial b_j^l}{\partial y^k}dx^i\wedge\theta^k
    \\
    &\equiv\sum_{1\le i<j\le r}\langle\theta^l,[X_i,X_j]\rangle dx^i\wedge dx^j\pmod{\theta^1,\dots,\theta^r}.
\end{align*}
Thus $$\theta^1\wedge\dots\wedge\theta^{n-r}\wedge d\theta^l=\sum_{1\le i<j\le r}\langle\theta^l,[X_i,X_j]\rangle \theta^1\wedge\dots\wedge\theta^{n-r}\wedge dx^i\wedge dx^j,\quad 1\le l\le n-r.$$
By shrinking $U$ if necessary we can find global sections $\lambda^1,\dots,\lambda^r\in \Co^{\alpha}_\loc(\Mf;\Se^\bot)$ such that $\lambda^l|_U=\theta^l$. By assumption \ref{Item::DisInv::CharInv2::coSec} we see that $\langle\theta^l,[X_i,X_j]\rangle=0\in \Co^{\alpha-1}_\loc(U)$. Since $[X_i,X_j]=\sum_{k=1}^{n-r}\langle\theta^l,[X_i,X_j]\rangle\Coorvec{y^l}$ we know $[X_i,X_j]=0$ are commutative, proving the condition \ref{Item::DisInv::CharInv1::Gen1}.
\end{proof}
\section{Remarks on Integrability and Problems in Regular Subbundles}\label{Section::DisInv::Integrable}

Unlike involutivity, it is more difficult to discuss integrability in the low regularity setting.

On real tangent subbundles, from discussion of the real Frobenius theorem we get the following definition:
\begin{defn}\label{Defn::DisInv::LocInvPara}
Let $\V$ be a continuous rank $r$ real tangent subbundle over a $n$-dim $C^1$-manifold $\Mf$.

We say $\V$ is \textbf{locally integrable via parameterizations} if for every point $p\in \Mf$ there exists a neighborhood $U\subset\R^r_u\times\R^{n-r}_v$ of $(0,0)$ and a map $\Phi(u,v)$ that is homeomorphism onto its image such that 
\begin{itemize}[nolistsep]
    \item $\Phi(0,0)=p$.
    \item $\frac{\partial\Phi}{\partial u^1},\dots,\frac{\partial\Phi}{\partial u^r}:U\to T\Mf$ are continuous maps.
    \item For every $(u,v)\in U$, $\frac{\partial\Phi}{\partial u^1}(u,v),\dots,\frac{\partial\Phi}{\partial u^r}(u,v)$ form a basis for $\V_{\Phi(u,v)}$.
\end{itemize}
\end{defn}

In this setting, Theorem \ref{MainThm::LogFro} implies the following:
\begin{prop}\label{Prop::DisInv::LogFroReStated}
An involutive log-Lipschitz real tangent subbundle must be locally integrable via parameterizations.
\end{prop}

It is not difficult to see that a H\"older involutive real subbundles may not be  locally integrable via parameterizations since it may not admit a foliation. See Example \ref{Exmp::DisInv::ModInv2}.





We can talk about integrability pointwisely. For convenience we define them on generalized subbundles (recall Definition \ref{Defn::DisInv::GenSub}).

\begin{defn}\label{Defn::DisInv::PointwiseInt}
Let $\Gc$ be a continuous generalized subbundle over a $n$-dim $C^1$-manifold $\Mf$.

\begin{itemize}[parsep=-0.3ex]
    \item We say $\Gc$ is \textbf{pointwisely integrable}, if for every $p\in \Mf$ there is a $C^1$-submanifold $\Sf\subset \Mf$ containing $p$, such that $\Se|_\Sf=T\Sf$, i.e. $\Gc\cap T_q\Mf=T_q\Sf$ for all $q\in \Sf$. We call $\Sf$ an \textbf{integral submanifold} for $\Gc$ near $p$.
    \item We say $\Gc$ is \textbf{uniquely integrable}, if in addition to pointwisely integrability the integral submanifold is unique: If both $\Sf_1$ and $\Sf_2$ are integral submanifolds for $\Se$ near $p\in \Mf$, then there is a neighborhood $U\subseteq \Mf$ of $p$ such that $\Sf_1\cap U=\Sf_2\cap U$.
\end{itemize}
\end{defn}

On the H\"older regular subbundles, it is easy to see that pointwise integrability does not implies unique integrability. For example, for $0<\alpha<1$, the subbundle $\V:=\Span(\Coorvec x+|y|^\alpha\Coorvec y)$ in $\R^2_{x,y}$ is pointwise integrable because it is of rank 1 and every continuous vector field has ODE existence. $\V$ is not uniquely integrable because on $(0,0)$,
\begin{equation*}
    \gamma_0(t):=(t,0),\quad \gamma_1(t):=(t,((1-\alpha)\max(0,t))^\frac1{1-\alpha}),\quad t\in\R,
\end{equation*}
give two integral curves for $\V$ that has different images.

For regular subbundles that have Osgood regularity, the integrability conditions in Definitions \ref{Defn::DisInv::LocInvPara} and \ref{Defn::DisInv::PointwiseInt} are all coincide.
\begin{defn}
Let $\mu$ be a modulus of continuity, that is, $\mu:[0,1]\to[0,1]$ is a strictly increasing map  such that $\mu(0)=0$. We say $\mu$ is a Osgood modulus, if $\int_0^1dr/\mu(r)=\infty$.

For an open subset $\Omega\subseteq\R^n$, we define $\Co^\mu(\Omega)$ be the same as in \eqref{Eqn::Hold::CMuSpace}.

We say a function $f:\Omega\subseteq\R^n\to\R$ is locally Osgood continuous, if for every precompact open $\Omega'\Subset\Omega$, there is a Osgood modulus $\mu$ such that $f\in\Co^\mu(\Omega')$. We say a subbundle is Osgood continuous, if it is spanned by vector fields whose coefficients are (locally) Osgood continuous.
\end{defn}

\begin{lem}\label{Lem::DisInv::CharofInt}
Let $\V\le T\Mf$ be a regular subbundle which is Osgood continuous. Then the following are equivalent.
\begin{enumerate}[parsep=-1ex,label=(\roman*)]
    \item\label{Item::DisInv::CharofInt::0} $\V$ is locally integrable via parameterizations.
    \item\label{Item::DisInv::CharofInt::1} $\V$ is pointwise integrable.
    \item\label{Item::DisInv::CharofInt::2} $\V$ is uniquely integrable.
\end{enumerate}
\end{lem}
\begin{proof}Clearly
\ref{Item::DisInv::CharofInt::0} implies \ref{Item::DisInv::CharofInt::1}, and \ref{Item::DisInv::CharofInt::2} implies \ref{Item::DisInv::CharofInt::1}.

We assume that $\dim\Mf=n$ and $\rank\V=r$.

By passing to a local coordinate chart we can assume that $\Mf\subseteq\R^n$ is an open subset. Applying Lemma \ref{Lem::ODE::GoodGen} we can find a linear coordinates $(x^1,\dots,x^r,y^1,\dots,y^{n-r})$ and Osgood vector fields $X_j=\Coorvec{x^j}+\sum_{k=1}^{n-r}a_j^k\Coorvec{y^j}$, $1\le j\le r$ that span $\V$. By \cite[Theorem 3.2]{BahouriCheminDanchin} the multi-flow $\exp_X$ is locally defined near $(0,0)\in\R^r\times\R^n$. By \cite[Corollary 3.6]{BahouriCheminDanchin} (where we can consider a single vector field $\tilde X(t,x):=(t\cdot X(x),0)$ in $\R^{n+r}$) we see that $\exp_X$ is a continuous map.

\medskip\noindent
\ref{Item::DisInv::CharofInt::1} $\Rightarrow$ \ref{Item::DisInv::CharofInt::2} and \ref{Item::DisInv::CharofInt::0}: It suffices to consider the case $p=0\in\R^n$. By taking scaling we can assume that $X_1,\dots,X_r$ are defined on $3\B^n$, and $\|(a_j^k)\|_{C^0(3\B^n;\R^{(n-r)\times n}}\le1$. Thus $\exp_X:\B^r\times\B^n\to3\B^n$ is defined.

Now for every $y\in\B^{n-r}$ there is a $r$-dimensional integral submanifold containing $(0,y)$. Since $X_1,\dots,X_r$ spans $\V$ and their $\Coorvec x$-components do not vanish, we can find a regular $C^1$-parameterization $\phi^y:\B^r\to3\B^n$ such that $\phi^y(0)=(0,y)$ and $\phi^y(t)=(t,\ast)$ for every $t\in\B^n$. 

By assumption $\frac{\partial\phi^y}{\partial t^1}(t),\dots,\frac{\partial\phi^y}{\partial t^r}(t)\in\V_{\phi^y(t)}$, thus $(t\cdot\Coorvec t)\phi^y(t)\in\V_{\phi^y(t)}$. However $(t\cdot\Coorvec t)\phi^y(t)=(t,\ast)$, we must have $(t\cdot\Coorvec t)\phi^y(t)=(t\cdot X)(\phi^y(t))$, for every $t\in\B^n$. We conclude that $\phi^y(t)=\exp_X(t,(0,y))$.

Therefore $\phi^y$ is uniquely determined by $X_1,\dots,X_r$ and $y\in\B^{n-r}$, thus $\phi^0(\B^r)\subseteq\R^n$ is the unique integral manifold containing $0$, giving the condition \ref{Item::DisInv::CharofInt::2}.

Let $\Phi(t,s):=\exp_X(t,(0,s))$, we see that $\Phi$ is continuous and $\Phi(t,s)=\phi^s(t)$. Since $\phi^s$ is $C^1$ whose derivatives span $\V$, we conclude that $\Phi$ is $C^1$ in $t$ and $\nabla_t\Phi(t,s)$ span $\V_{\Phi(t,s)}$.

Notice that $F(x,y)=(x,e^{-x\cdot X}y)$ is a continuous map which is the inverse map of $\Phi$, we conclude that $\Phi$ is homeomorphism onto its image. This finishes the proof of \ref{Item::DisInv::CharofInt::0}.
\end{proof}

    

\begin{remark}
It is not known to the author whether \ref{Item::DisInv::CharofInt::1} $\Rightarrow$ \ref{Item::DisInv::CharofInt::2} in Lemma \ref{Lem::DisInv::CharofInt} remained true for general continuous real subbundles, or even for subbundles with H\"older continuity. It is also not known to the author whether local integrability (via parameterizations) always implies involutivity when the subbundles are $\Co^{\frac12+}$.
\end{remark}

On complex subbundles there is a different type of integrability. The following definition can be founded in \cite{HypoAnalytic} and \cite{Involutive}.
\begin{defn}\label{Defn::DisInv::LocInvDiff}
Let $\Mf$ be a $n$-dimensional $C^2$ manifold and let $\Se\le\C T\Mf$ be a $C^1$ complex tangent subbundle of rank $r$. 

We say $\Se$ is \textbf{locally integrable via differentials}, if for any $p\in \Mf$ there are neighborhood $U\subseteq\Mf$ of $p$ and $C^1$ functions $f^1,\dots,f^{n-r}:U\to\C$, such that $df^1(q),\dots,df^{n-r}(q)\in\C T_q^*\Mf$ are linearly independent for all $q\in U$, and $\Se=\Span(df^1,\dots,df^{n-r})^\bot$.
\end{defn}

We cannot use parameterizations since an involutive subbundle may not necessarily be spanned by coordinate vector fields unless it is a complex Frobenius structure.
\begin{example}
In Remark \ref{Rmk::ODE::RmkInvStr} \ref{Item::ODE::RmkInvStr::NotBund}, on $\R^2_{x,t}$ we see that the subbundle generated by Mizohata operator $X=\partial_t+it\partial_x$ cannot be the span of coordinate vector fields. On the other hand $(\Span X)^\bot=\Span(dx-itdt)=\Span(d(x-it^2/2))$, so $\Span X$ is locally integrable via differentials.
\end{example}

It is known that real-analytic involutive structures are always locally integrable, while there exists a smooth involutive structure which is not locally integrable, see \cite[Section I.16]{Involutive}.

It is difficult to generalize Definition \ref{Defn::DisInv::LocInvDiff} for non $C^1$ complex subbundles since it does not make sense to talk about linearly independence of differentials pointwisely. On one hand, without linearly independence, it is possible that the set $\Span(df^1,\dots,df^{n-r})^\bot\subseteq\C T\Mf$ has fibre whose rank is larger than $r$. On the other hand, in the non $C^1$-setting, it is also possible to construct a span of differentials whose dual bundle has rank less than $r$:

\begin{lem}\label{Example::DisInv::BadDfSpan}
For a $\Co^{\frac12+}$ function $f:\R^n\to\R$, we associate a set $\V^f\subseteq T\Mf$ given by the following:
\begin{equation*}
    \V^f:=\{(p,X(p))\in T\Mf:X\in\Co^{\frac12+}(U;\R^n)\text{ for some open }U\ni p\text{ such that }\langle df,X\rangle=Xf=0\text{ as distributions}\}.
\end{equation*}
Then there is a Lipschitz function $f$ such that $\V^f=\{0\}$.
\end{lem}
We see that if $f$ is $\Co^{\frac32+}$, then $\V^f=(\Span df)^\bot$ is the classical pointwise dual set.

\begin{proof}Take $f_0(x):=|x|$, we claim that 
\begin{equation}\label{Eqn::DisInv::BadDfSpan::Vf1}
    \V^{f_0}=\{(p,v)\in\R^{2n}:p\neq0,\ p\cdot v=0\}\cup\{(0,0)\}
\end{equation}
Indeed $f_0$ is regular for $x\neq0$, so $\V^{f_0}_p=(\Span df_0(p))^\bot$ for $p\neq0$. At $p=0$, if $X$ is a continuous vector fields near $0$ such that $Xf=0$, then $\frac x{|x|}\cdot X(x)=0$. Let $x\to 0$ along each ray, we see that $\theta\cdot X(0)=0$ for all $\theta\in\Sp^{n-1}$, i.e. $X(0)=0$. Therefore $\V^{f_0}_0=\{0\}$.

Now let $f=\sum_{j=1}^\infty3^{-j}|x-q_j|$ where $\{q_j\}_{j=1}^\infty=\mathbb Q^n\subset\R^n$ is an enumeration of rational point. We claim that $\V^f=\R^n\times\{0\}$ is the zero subbundle as desired. 

Indeed for each $j=1,2,3,\dots$, let $f_j(x):=\sum_{k\neq j}^\infty3^{-k}|x-q_k|$. We see that $\nabla f_j(x)=\sum_{k\neq j}3^{-k}(x-q_k)/|x-q_k|$ for all $x\notin\{q_k:k\neq j\}$, and therefore the following limit exists (in the calculus setting):
\begin{equation*}
    \lim_{x\to q_j;x\notin\{q_k:k\neq j\}}\nabla f_j(x)=\nabla f_j(q_j)=\sum_{k\neq j}3^{-k}\frac{q_j-q_k}{|q_j-q_k|}.
\end{equation*}

Let $U_j\ni q_j$ be an open set and let $X\in C^0(U_j;\R^n)$ be a continuous vector field  such that $Xf_j=\nabla f_j\cdot X=0$ almost everywhere. Since $X$ is continuous, we see that $\lim_{r\to0}\|\nabla f(x)\cdot X(q_j)\|_{L^\infty(B(q_j,r))}=0$. On the other hand, for each irrational direction $\theta\in\Sp^{n-1}$ (i.e. $(\theta\cdot\R)\cap\mathbb Q^n=\{0\}$), we see that in the calculus setting
\begin{equation*}
    \lim_{r\to0}\nabla f(q_j+r\theta)\cdot X(q_j)=\theta\cdot X(q_j)+\lim_{r\to0}\nabla f_j(q_j+r\theta)\cdot X(q_j)=(\theta+\nabla f_j(q_j))\cdot X(q_j).
\end{equation*}

Therefore $(\theta+\nabla f_j(q_j))\cdot X(q_j)=0$ for almost every $\theta\in\Sp^{n-1}$. Since any positive (surface) measurable subset of $\nabla f_j(q_j)+\Sp^{n-1}$ contains $n$-linearly independent vectors, we conclude that $X(q_j)=0$.

Since $\{q_j\}$ is dense in $\R^n$, if $X:U\subseteq\R^n\to\V^f$ is any continuous vector field, then we must have $X(q_j)=0$ for all $q_j\in U$. By continuity of $X$ and density of $\{q_j\}$, we see that $X\equiv0$. This is also true for $\Co^{\frac12+}$-vector fields as well since $\Co^{\frac12+}\subsetneq C^0$.
Therefore $\V^f$ is a zero subbundle.
\end{proof}
\begin{remark}
If we consider $\Span(df)$ as a ``$L^\infty$-cotangent subbundle'', which is an almost everywhere defined section to the Grassmannian spaces, then $\Span(df)^\bot$ is still a rank $(n-1)$ tangent subbundle though only almost everywhere defined. The above example shows that there is no nonzero continuous vector field which is almost everywhere contained in $\Span(df)^\bot$.
\end{remark}

\begin{remark}
Example \ref{Example::DisInv::BadDfSpan} does not have issue with Convention \ref{Conv::EllipticPara::CoSpan}. In Convention \ref{Conv::EllipticPara::CoSpan} we have a product manifold $\Mf\times\Nf$ endowed with a product coordinate chart $(\mu,\nu)=(\mu^1,\dots,\mu^n,\nu^1,\dots,\nu^q)$. We consider a subbundle $\Se\le (\C T\Mf)\times\Nf$ and functions $f^1,\dots,f^{n-r}$ which are $C^1$ in $\mu$, such that $\frac{\partial f}{\partial \mu}d\mu$ is a collection of $(n-r)$ linearly independent differentials and $\Se=\Span(\frac{\partial f}{\partial \mu}d\mu,d\nu)^\bot$. 

Even though $f^1,\dots,f^{n-r}$ are not $C^1$, on each leaf $\Mf\times\{q_0\}$, the restriction $\Se|_{\Mf\times\{q_0\}}\le \C T\Mf$ is locally integrable via differential. Since we have $(\Se|_{\Mf\times\{q_0\}})^\bot=\Span d(f(\cdot,q_0))$, where $d(f^1(\cdot,q_0)),\dots,d(f^{n-r}(\cdot,q_0))$ are continuous linearly independent differentials on an open subset of $\Mf$.
\end{remark}

\medskip
We can naturally to ask whether the Frobenius theorem holds in the general Osgood setting:
\begin{prob}\label{Prob::OsgFro}
Let $\V\le T\Mf$ be an Osgood subbundle which is involutive in the sense of distributions. Does $\V$ always pointwise integrable? If so, does $\V$ always locally integrable via parameterizations.
\end{prob}

We know that if such parameterization $\Phi(t,s):\Omega'\times\Omega''\to\Mf$ exists, then applying \cite[Lemma 3.4]{BahouriCheminDanchin}, the regularity estimate $\Phi\in C^1_t(\Omega';C^0_s(\Omega'';\Mf))$ is the best possible we can get. 

Problem \ref{Prob::OsgFro} is reduced to the flow commuting problem, when vector fields are merely Osgood:
\begin{prob}\label{Prob::ConjOsgoodFlowCom}
Let $X,Y$ be two Osgood vector fields in $\R^n$ satisfying $[X,Y]=0$ in the sense of distributions. Do we must have $e^{tX}\circ e^{sY}=e^{sY}\circ e^{tX}$ locally for small $t$ and $s$?
\end{prob}

Unfortunately we do not know how to prove the Problem \ref{Prob::ConjOsgoodFlowCom} for non log-Lipschitz vector fields in general. One attempt to deal with the Problem \ref{Prob::ConjOsgoodFlowCom} is to show whether an Osgood subbundle is asymptotically involutive (see  \cite[Definition 1.17]{ContFro}). Based on this, we can ask an even weaker question:
\begin{prob}
Is there a modulus of continuity $\mu$ such that the following hold?
\begin{itemize}[parsep=-0.3ex]
    \item $\lim\limits_{r\to0}\mu(r)^{-1}\cdot r\log\frac1r=0$.
    \item If $\V$ is a tangent subbundle with modulus of continuity $\mu$, and if $\V$ is involutive in the sense of distributions, then $\V$  is asymptotically involutive.
\end{itemize}
\end{prob}

A major difficulty in understanding Problem \ref{Prob::ConjOsgoodFlowCom} is to verify \eqref{Eqn::Hold::CorAsyInv::Eqn1}. But the exponential term in \eqref{Eqn::Hold::CorAsyInv::Eqn2} is somehow subtle.

\section{Remarks and Examples on Nonsmooth Singular Subbundles}\label{Section::DisInv::Sing}

We recall the Definition \ref{Defn::DisInv::GenSub} for generalized tangent subbundles and Definition \ref{Defn::DisInv::ModInv} for involutivity for modules of sections. In this section we focus only on real singular subbundles since we do not have a good analogy for the singular version of complex Frobenius theorem.

For completeness we define integrability on submodules:
\begin{defn}
Let $\alpha>0$, let $\Mf$ be a $\Co^{\alpha+1}$-manifold, and let $\Fs\subseteq\Co^\alpha_\loc(\Mf;T\Mf)$ be a submodule over the ring $\Co^\alpha_\loc(\Mf)$. We say $\Fs$ is pointwise integrable, if for every $p\in\Mf$ there is a neighborhood $U\subseteq\Mf$ of $p$ and a $C^1$-submanifold $\Sf\subseteq U$ containing $p$, such that $T_q\Sf=\Span\{X(q):X\in\Fs\}$ for every $q\in\Sf$.
\end{defn}

Similar to Proposition \ref{Prop::DisInv::LogFroReStated}, Theorem \ref{MainThm::SingFro} can be restated as the following:
\begin{prop}\label{Prop::DisInv::SingFroReStated}
Let $\Mf$ be a $C^{1,1}$-manifold and let $\Fs\subseteq\Co^\LogL_\loc(\Mf;T\Mf)$ be a finitely generated module over the ring $\Co^\LogL_\loc(\Mf)$. If $\Fs$ is $\Co^{\LogL-1}$-involutive, then $\Fs$ is pointwise integrable.
\end{prop}

Given a generalized subbundle, we can still talk about $\Co^\beta$-sections for $\beta>0$. To make it more general, we define sections for generalized subbundle over general vector bundles.
\begin{defn}\label{Defn::DisInv::GenSubofVB}

    Let $\alpha\in(0,\infty)\cup\{\infty\}$. Let  $\Mf$ be a $C^1$ manifold which is at least $\Co^\alpha$. Let $\E\twoheadrightarrow\Mf$ be a $\Co^\alpha$ real vector bundle. A $\Co^\alpha$ generalized subbundle of $\E$ is a subset $\Gc\subset\E$, such that $\Gc_p\subseteq\E_p$ is a (real) vector subspace for every $p\in\Mf$, and for every $p\in\Mf$ there is a neighborhood $U\subseteq\Mf$ of $p$ and finite sections $\xi_1,\dots,\xi_m\in\Co^\alpha_\loc(U;\E|_U)$ such that $\Gc_q=\Span(\xi_1(q),\dots,\xi_m(q))(\le\E_q)$ for every $q\in U$.
    
    We define $\rank\Gc:=\max\{\rank \Gc_p:p\in\Mf\}$ and $\Co^\alpha_\loc(\Mf;\Gc):=\{\xi\in\Co^\alpha_\loc(\Mf;\E):\xi(p)\in\E_p,\ \forall p\in\E \}$.
    
\end{defn}

The definition for complex vector bundles and complex generalized subbundles are the same and we omit the details.

We can remove the finiteness assumption of local generators for $\Gc$ in Definitions \ref{Defn::DisInv::GenSub} and \ref{Defn::DisInv::GenSubofVB} due to the following:
\begin{prop}[Serre-Swan for generalized subbundles]\label{Prop::DisInv::SerreSwanGenSub}
Let $\alpha\in(0,\infty)\cup\{\infty\}$ and $n\ge1$. Let $\Mf$ be a $n$-dimensional manifold (which is at least $C^1$ and $\Co^\alpha$) and let $\E\twoheadrightarrow\Mf$ be a $\Co^\alpha$ real vector bundle. 

Let $\Gc\subset\E$ be a generalized subbundle and has rank $r$. Then the following are equivalent:
\begin{enumerate}[parsep=-0.3ex,label=(\roman*)]
    \item\label{Item::DisInv::SerreSwanGenSub::GlobalAll} $\Gc_p=\Span\{X(p):X\in\Co^\alpha_\loc(\Mf;\E)\}$ for every $p\in\Mf$.
    \item\label{Item::DisInv::SerreSwanGenSub::Swan} There are at most $r(n+1)$ many (global) sections $\xi_1,\dots,\xi_{r(n+1)}\in\Co^\alpha_\loc(\Mf;\E)$ such that $\Gc_p$ is spanned by $\xi_1(p),\dots,\xi_{r(n+1)}(p))$ for all $p\in\Mf$. 
    \item\label{Item::DisInv::SerreSwanGenSub::LocalAll} For every $p\in\Gc$ there is a $U\subseteq\Mf$ of $p$ such that  $\Gc_q=\Span\{X(q):X\in\Co^\alpha_\loc(U;\E|_U)\}$ for every $q\in U$.
    \item\label{Item::DisInv::SerreSwanGenSub::LocalFinite} $\Gc$ is a $\Co^\alpha$ generalized subbundle, from Definition \ref{Defn::DisInv::GenSubofVB}.
    
\end{enumerate}
\end{prop}

Here we need Ostrand's theorem on colored dimension:
\begin{lem}\label{Lem::DisInv::ColorDim}
Let $\Mf$ be a $n$-dimensional topological manifold and let $\Pf\subseteq\Mf$ be a closed subset. Let $\Us=\{U_i:i\in I\}$ be an open cover of $\Pf$. Then there are $(n+1)$-collections $\Vs_1,\dots,\Vs_{n+1}$ of distinct open sets on $\Mf$ such that:
\begin{enumerate}[parsep=-0.3ex,label=(\alph*)]
    \item\label{Item::DisInv::ColorDim::Disjoint} For each $1\le j\le n+1$, $\Vs_j$ is a collection of disjoint open sets of . 
    \item For each $1\le j\le n+1$ and each $V\in\Vs_j$, there exists an $U\in\Us$ such that $V\subseteq U$.
    \item $\bigcup_{j=1}^{n+1}\bigcup_{V\in\Vs_j}V\supseteq\Pf$ is an open cover.
\end{enumerate}
\end{lem}
\begin{proof}
A closed subset of $n$-dimensional topological manifold has covering dimension at most $n$, see \cite[Exercise 50.9 and Theorem 50.1]{MunkresTopology}. The result then follows from \cite[Theorem 1]{ColorDim}.
\end{proof}
\begin{proof}[Proof of Proposition \ref{Prop::DisInv::SerreSwanGenSub}]
Clearly \ref{Item::DisInv::SerreSwanGenSub::Swan} $\Rightarrow$ \ref{Item::DisInv::SerreSwanGenSub::LocalFinite} $\Rightarrow$ \ref{Item::DisInv::SerreSwanGenSub::LocalAll}. By the standard partition of unity argument we have \ref{Item::DisInv::SerreSwanGenSub::LocalAll} $\Rightarrow$ \ref{Item::DisInv::SerreSwanGenSub::GlobalAll}.

Our goal is to prove \ref{Item::DisInv::SerreSwanGenSub::GlobalAll} $\Rightarrow$ \ref{Item::DisInv::SerreSwanGenSub::Swan}. We proceed by induction on $k=0,1,\dots,r$ of the following
\begin{equation}\label{Eqn::DisInv::SerreSwanGenSub::Induction}
    \exists \xi_1,\dots,\xi_{k(n+1)}\in\Co^\alpha_\loc(\Mf;\Gc),\text{ such that }\rank\big(\Span\{\xi_1(p),\dots,\xi_{k(n+1)}(p)\}\big)\ge\min(\rank\Gc_p,k),\ \forall p\in\Mf.
\end{equation}
The base case $k=0$ is trivial as there is no section to take.

Let $k\ge1$ and suppose we have the case $k-1$. We define a set $\Nf_k:=\{p\in\Mf:\rank\Gc_p\ge k\}$. For $p\in\Nf_k$ by assumption \ref{Item::DisInv::SerreSwanGenSub::GlobalAll} we can find sections $\eta_1,\dots,\eta_k\in\Co^\alpha_\loc(\Mf;\Gc)$ such that $\eta_1(p),\dots,\eta_k(p)\in\E_p$ are linearly independent. Thus $\eta_1,\dots,\eta_k$ are linearly independent in a neighborhood of $p$, which means $p$ is an interior point of $\Nf_k$. We conclude that $\Nf_k$ is an open set.

By induction hypothesis there are global $\Co^\alpha$ sections $\xi_1,\dots,\xi_{(k-1)(n+1)}$ of $\Gc$ such that their span at each $p\in\Mf$ has rank $\ge\min(\rank\Gc_p,k-1)$. Let $p\in\Nf_k$, since $\Nf_k$ is open, we can find a precompact neighborhood $U^p\Subset\Nf_k$ of $p$ and a section $\eta^p\in\Co^\alpha_\loc(U^p;\Gc|_{U^p})$ such that 
\begin{equation}\label{Eqn::DisInv::SerreSwanGenSub::AssumeSpan}
    \rank\Span\{\xi_1(q),\dots,\xi_{(k-1)(n+1)}(q),\eta^p(q)\}\ge k,\quad\forall q\in U^p.
\end{equation}

Applying Lemma \ref{Lem::DisInv::ColorDim} to the open cover $\{U^p:p\in\Nf_k\}$ of $\Nf_k$, there are collections $\Vs_1=\Vs_1^k,\dots,\Vs_{n+1}=\Vs_{n+1}^k$ of open sets, associated with an index function $\tau:\Vs_1\cup\dots\cup\Vs_{n+1}\to\Mf$, such that
\begin{enumerate}[nolistsep,label=(\alph*)]
    \item\label{Item::DisInv::SerreSwanGenSub:Color::Disjoint} $\Vs_j$ are disjoint collections, for each $j=1,\dots,n+1$.
    \item $V\subseteq U^{\tau (V)}$ holds for every $V\in\Vs_1\cup\dots\cup\Vs_{n+1}$, where $U^{\tau(V)}\Subset\Nf_k$ is the neighborhood of $\tau(V)\in\Nf_k$.
    \item $\Vs_1\cup\dots\cup\Vs_{n+1}$ is an open cover of $\Nf_k$.
\end{enumerate}

We can assume that $\Mf$ is connected, hence $\Nf_k$ is second counterable, which means that $\Vs_1\cup\dots\cup\Vs_{n+1}$ is a collection of at most counterable elements.
By reordering the index we can write $\Vs_j=\{V_{j,1},V_{j,2},V_{j,3},\dots\}$. For each $j$, if $\Vs_j$ is a finite set, then we define $V_{j,l}:=\varnothing$  for $l\ge\#\Vs_j+1$. 

We endow $\E$ with a $\Co^\alpha$ bundle metric so that the $\Co^\beta$-norms are defined on the space of sections of $\E$, for every $0<\beta\le\alpha$ such that $\beta<\infty$.

Since each open subset in $\Vs_1\cup\dots\cup\Vs_{n+1}=\{V_{j,l}:1\le j\le n+1,\ l\ge0\}$ is precompact in $\Nf_k$, we can take a partition of unity $\{\chi_{j,l}:1\le j\le n+1,\ l\ge0\}\subset \Co_c^\alpha(\Nf_k)$ associated with $\{V_{j,l}\}_{j,l}$, where we emphasize that each $\chi_{j,l}$ has compact support in $\Nf_k$.


We then define for $j=1,\dots,n+1$,
\begin{equation}\label{Eqn::DisInv::SerreSwanGenSub::DefXi}
    \xi_{(k-1)(n-1)+j}:=\sum_{l=1}^\infty 2^{-l}\frac{\chi_{j,l}\cdot\eta^{\tau(V_{j,l})}}{\|\chi_{j,l}\cdot\eta^{\tau(V_{j,l})}\|_{\Co^{\min(\alpha,l)}}}.
\end{equation}
Here the $\Co^{\min(\alpha,l)}$ norm is only technical. One can use $\Co^\alpha$-norm if $\alpha<\infty$, and $\Co^l$-norm if $\alpha=\infty$.

Clearly the sum in \eqref{Eqn::DisInv::SerreSwanGenSub::DefXi} converges in the space $\Co^\alpha(\Mf;\E)$ with respect to the fixed bundle metric. In particular we know $\xi_{(k-1)(n-1)+j}\in\Co^\alpha_\loc(\Mf;\Gc)$.

For each $p\in\Nf_k$, there is a $1\le j_0\le n+1$ and $V_{j_0,l}\in\Vs_{j_0}$ such that $\chi_{j_0,l}(p)\neq0$.
By  property \ref{Item::DisInv::SerreSwanGenSub:Color::Disjoint}, the sum in \eqref{Eqn::DisInv::SerreSwanGenSub::DefXi} has at most 1 nonzero at each point. Thus for such $j$, $\xi_{(k-1)(n+1)+j_0}(p)=c_{j_0}(p)\eta^{\tau(V_{j_0,l})}(p)$ for some constant $c_{j_0}(p)>0$ that comes from \eqref{Eqn::DisInv::SerreSwanGenSub::DefXi}. 

Thus by \eqref{Eqn::DisInv::SerreSwanGenSub::AssumeSpan}, for such $j_0=j_0(p)$, 
\begin{align*}
    \rank\Span(\xi_1(p),\dots,\xi_{k(n+1)}(p))\ge&\rank\Span(\xi_1(p),\dots,\xi_{(k-1)(n+1)}(p),\xi_{(k-1)(n+1)+j_0}(p))
    \\
    =&\rank\Span(\xi_1(p),\dots,\xi_{(k-1)(n+1)}(p),\eta^{\tau(V_{j_0,l})}(p))\ge k.
\end{align*}
Therefore \eqref{Eqn::DisInv::SerreSwanGenSub::Induction} is satisfied for all $p\in\Nf_k$.

On the other hand, if $p\notin\Nf_k$, then we have $\rank\Gc_p\le k-1$. By induction hypothesis \eqref{Eqn::DisInv::SerreSwanGenSub::Induction} is already satisfied. Therefore $\xi_1,\dots,\xi_{k(n+1)}$ are the desired sections for the case $k$ in \eqref{Eqn::Hold::PfQPIFT::Induction}. 

By reaching the case $k=r$, we get $r(n+1)$ many sections and complete the proof.
\end{proof}

Unlike Definition \ref{Defn::ODE::CpxSubbd}, on a singular subbundle $\Gc$ of rank $r$, locally we may not be able to find $r$-many sections that span $\Gc$ in a neighborhood. See the following example.

\begin{example}
On $\R^3$, we let $\Gc_1:=\{(0,0)\}\cup\{(x,v):x\in\R^3\backslash\{0\},\ x\cdot v=0\}\subset T\R^3$. So $\Gc_1=\{(0,0)\}\cup\coprod_{x\neq0}T_{\frac x{|x|}}\Sp^2$. Clearly $\rank\Gc_1=2$.

Take finitely many vector fields $X_1,\dots,X_m$ (in fact we can take $m=3$) that span $T\Sp^2$, we see that $\Gc_1$ is spanned by smooth vector fields $\big [x\mapsto e^{-1/|x|}X_j(\frac x{|x|})\big]$, $1\le j\le m$. Therefore $\Gc_1$ is a smooth singular subbundle.

Suppose we can find continuous vector fields $Y_1,Y_2:U\subseteq\R^3\to\R^3$ near $0\in\R^3$ such that $\Gc_1|_U$ is spanned by $Y_1$ and $Y_2$. Since $U\ni 0$, we have $B^3(0,\delta)\Subset U$ for some $\delta>0$. Thus on the boundary $\partial B^3(0,\delta)=\delta\Sp^2$,  $Y_1|_{\delta\Sp^2}$ and $Y_2|_{\delta\Sp^2}$ span $\Gc_1|_{\delta\Sp^2}\simeq T\Sp^2$. Since $T\Sp^2$ is a regular vector bundle over $\Sp^2$ with rank $2$, we see that $Y_1|_{\delta\Sp^2}$ and $Y_2|_{\delta\Sp^2}$ are linearly independent at every point in $\delta\Sp^2$.  However by the Hairy Ball Theorem $Y_1|_{\delta\Sp^2}$ must have zero points, contradiction.

Therefore we need at least 3 vector fields that span $\Gc_1$ in a neighborhood of $0$.
\end{example}

Unlike the regular case it does not make sense to talk about sections of singular subbundle which are not pointwisely defined.
\begin{example}
Consider $\Gc_2:=((\R\backslash\{0\})\times\R^1)\cup\{(0,0)\}\subset T\R$. It is a smooth singular subbundle since it is generated by $x\Coorvec x$. 

Note that $\Gc_2(x)=\R^1$ for almost every $x\in\R$ since $\{0\}\subset\R$ is a zero measured set. If we define a ``$L^\infty$-section'' of $\Gc_2$ to be a $L^\infty$-vector field that lays in $\Gc_2$ for almost every point, then $L^\infty(\R;\Gc_2)=L^\infty(\R;T\R)$. There is no distinguishment between $\Gc_2$ and $T\R$ in this setting.
\end{example}

Another reason why we cannot define distributional sections of $\Gc$ is that we do not have \eqref{Eqn::DisInv::DisSecVB:Pf::Ten>0} for $\alpha,\gamma>0$.

\begin{example}
Consider $\Gc_3:=((0,\infty)\times\R^1)\cup((-\infty,0]\times\{0\})\subset T\R$, this is a smooth singular subbundle generated by the smooth vector field $e^{-\frac1x}\mathbf1_{\R_+}(x)\Coorvec x$.

For every non-integer $\gamma>0$, we see that $X\in\Co^\gamma_\loc(\R;\Gc_3)$ implies $|X(x)|\le C|x|^\gamma$ for $x>0$ closed to 0. (In fact $|X(x)|\le C|x|^\gamma\log\frac1{|x|}$ for $x>0$ closed to $0$, if $\gamma=1,2,3,\dots$.)

On the other hand, suppose $X\in\Co^\gamma_\loc(\R)\otimes_{\Co^\alpha_\loc(\R)}\Co^\alpha_\loc(\R;\Gc_3)$ for some non-integer $\alpha>\gamma$, we can write $X=\sum_{j=1}^Nf^jY_j$ where $Y_j\in\Co^\alpha_\loc(\R;\Gc_3)$ and $f^j\in \Co^\alpha_\loc(\R)$. Clearly $X\in \Co^\gamma_\loc(\R;\Gc_3)$, $f^j$ are bounded near $0$ and $|Y_j(x)|\le C_j|x|^\alpha$ for $x>0$ near $0$. Therefore $|X(x)|\le C\max(0,x)^\alpha\ll\max(0,x)^\gamma$ for $x$ near $0$.

To summarize, for non-integers $0<\gamma<\alpha$, we have $\Co^\gamma_\loc(\R)\otimes_{\Co^\alpha_\loc(\R)}\Co^\alpha_\loc(\R;\Gc_3)\subsetneq\Co^\gamma_\loc(\R;\Gc_3)$ and the inclusion is strict. (This is still true if $\alpha$ or $\gamma$ is a positive integer.)
\end{example}

The following example shows that the $\Co^\beta$-involutivity of modules are not equivalent if we change the index $\beta$.
\begin{example}[{\cite[Example 8.3]{Lewis}}]\label{Exmp::DisInv::ModInv1}
Let $\alpha>0$. Let $\Fs_3\subset\Co^\alpha_\loc(\R^2;T\R^2)$ be the $\Co^\alpha_\loc(\R^2)$-submodule  generated by $X(x,y)=(x^2+y^2)\Coorvec x$ and $Y(x,y)=(x^4+y^4)\Coorvec y$.

By direct computation
\begin{equation*}
    \textstyle[X,Y](x,y)=-\frac{2y(x^4+y^4)}{x^2+y^2}X+\frac{4x^3(x^2+y^2)}{x^4+y^4}Y.
\end{equation*}
Note that the coefficients $\frac{2y(x^4+y^4)}{x^2+y^2}$ and $\frac{4x^3(x^2+y^2)}{x^4+y^4}$ are uniquely determined, because $X$ and $Y$ are linearly independent almost everywhere (they are linearly independent on $\R^2\backslash\{0\}$).

We see that $\Fs_3$ is $\Co^\Lip$-involutive but not $\Co^\beta$-involutive for any $\beta>1$. In particular if $\alpha>2$, then $\Fs_3$ is not $\Co^{\alpha-1}$-involutive.

Nevertheless, $\Gc_4:=\Span\Fs_3$ is still pointwise integrable. If we choose the submodule $\Fs_3^\alpha$ to be the $\Co^\alpha$ span of $X$ and $Y':=(x^2+y^2)\Coorvec y$, then $[X,Y']=-2yY'+2xX$, which means $\Fs_3^\alpha$ is always $\Co^{\alpha-1}$-involutive.
\end{example}

In the classical setting, we know if a smooth singular subbundle $\Gc\subseteq T\Mf$ is (pointwise) integrable, then $\Co^\infty_\loc(\Mf;\Gc)$ is (smoothly) involutive, and the converse direction is false, see Example \ref{Exmp::DisInv::ModInv2} below. 

In the smooth setting, the Hermann's theorem \cite{HermannFoliations} (see also \cite[Section 8.3]{Lewis} and \cite[Theorem 3]{ShortGuideFrobenius}) says the following, which is the smooth version to Proposition \ref{Prop::DisInv::LogFroReStated}:
\begin{prop}\label{Prop::DisInv::Hermann}
Let $\Gc$ be a $\Co^\infty$-generalized subbundle of a smooth manifold $\Mf$. Suppose there exists a finitely generated $\Co^\infty$-module $\Fs\subset\Co^\infty_\loc(\Mf;\Gc)$ that spans $\Gc$ (i.e. $\Gc_p=\Span\{X(p):X\in\Fs\}$ for all $p\in\Mf$) such that $\Fs$ is $\Co^\infty$-involutive. Then $\Gc$ is pointwise integrable.
\end{prop}

In other words the $\Gc_5$ in Example \ref{Exmp::DisInv::ModInv2} below cannot be spanned by a finitely generated $\Co^\infty$-module.

However, if we allow to use certain distributional involutivity on singular subbundles, then the finitely generated involutivity may not imply integrability.

\begin{example}[{\cite[Example 8.2]{Lewis}}]\label{Exmp::DisInv::ModInv2}
Let $\Gc_5:=\Span(\Coorvec x,e^{-\frac1x}\mathbf1_{\R_+}(x)\Coorvec y)\subset T\R^2$. We see that $\Gc_5$ is a smooth singular subbundle, $\Gc_5(x,y)=\R_x\times\{0\}$ if $x\le 0$ and $\Gc_5(x,y)=\R^2$ if $x>0$.

For every $y\in\R$, $\Gc_5(0,y)=\R\times\{0\}$ has rank 1. If $\gamma(t)=(x(t),y(t))$ is a $C^1$ integral curve passing $(0,y)$, say $x(t)>0$ when $t>0$ small, then we must have $\rank\Gc_5(\gamma(t))=2>\rank\Span(\dot\gamma(t))=1$ when $t>0$ small. Thus $\Gc_5$ is not pointwise integrable.

On the other hand, for every $\alpha>1$, let $X:=\Coorvec x$ and $Y^\alpha:=\max(0,x)^\alpha\Coorvec y$, we see that $\Fs^\alpha_5:=\{X,Y^\alpha\}\cdot\Co^\alpha_\loc(\R^2;T\R^2)$ is a $\Co^\alpha_\loc(\R^2)$-submodule of $\Co^\alpha$-vector fields, and
\begin{equation*}
    \textstyle[X,Y^\alpha]=\alpha\max(0,x)^{\alpha-1}\Coorvec y=(\pv\frac\alpha x)\cdot Y^\alpha.
\end{equation*}

Here $\pv\frac1x=\frac d{dx}\log|x|$, which is a well-defined $\Co^{-1}$-function on $\R^2_{x,y}$ by Lemma \ref{Lem::Hold::ZygExample}. 

Therefore $\Fs^\alpha_5$ is finitely generated and $\Co^{-1}$-involutive, but for the ground space of $\Fs^\alpha$, the singular subbundle $\Gc_5$ is not (pointwise) integrable.
\end{example}

For more examples and discussions of the singular Frobenius theorem, we refer \cite{Lewis} and \cite{ShortGuideFrobenius} to readers.

\begin{remark}
Based on Theorem \ref{MainThm::SingFro} or Proposition \ref{Prop::DisInv::SingFroReStated}, we know $\Gc_5$ cannot support a $\Co^\LogL$-module which is $\Co^{\LogL-1}$-involutive. 

However, we do not know whether Proposition \ref{Prop::DisInv::SingFroReStated} is still true if we replace the $\Co^{\LogL-1}$-involutivity by $\Co^\beta$-involutivity for $\beta<0$, even when the module is smooth. More precisely we can ask the following problem.
\end{remark}
\begin{prob}\label{Prob::SingHoldInv}
Let $\Mf$ be a smooth manifold and let $\Fs\subseteq\Co^\infty_\loc(\Mf;T\Mf)$ be a submodule over $\Co^\infty_\loc(\Mf)$. Let $-1<\beta<0$.

Suppose $\Fs$ is $\Co^\beta$-involutive. Does $\Fs$ always pointwise integrable?
\end{prob}

The major difficulty of Problem \ref{Prob::SingHoldInv} is somehow similar to that for Problem \ref{Prob::ConjOsgoodFlowCom}. If we consider $A_1,\dots,A_m\in\Co^\beta$ in the proof of Lemma \ref{Lem::ODE::GronApprox}, then for the right hand side of  \eqref{Eqn::ODE::GronApprox::GronwallPrep} we only have $|\frac d{dt}F^\sigma\circ\gamma^\sigma(t)|\le C_2'(2^{|\beta|\sigma}|F^\sigma\circ\gamma^\sigma(t)|+2^{-(1-\delta/2)\sigma})$, which blows up too fast and cannot tends to 0 if we try to use Gronwall's inequality.








\chapter{The Variants of Malgrange's Method}\label{Chapter::Malgrange}
In this chapter we give two key estimates, one for proving the complex Frobebius theorem, Theorems \ref{MainThm::CpxFro}, \ref{MainThm::RoughFro1} and \ref{MainThm::RoughFro2},  one for proving the quantitative Frobenius theorem, Theorem \ref{MainThm::QuantFro}. These two estimates are inspired by Malgrange's proof
of the Newlander-Nirenberg Theorem \cite{Malgrange}, see Sections \ref{Section::EllipticPara::Overview} and \ref{Section::Rough1FormOV} for their overviews.


\section{Estimates of Elliptic Structures with Parameters}\label{Section::EllipticPara}

Let $\Mf$  and $\Nf$ be two smooth manifolds. In this section we consider the product manifold $\Mf\times\Nf$ as the ambient space, so $(\C T\Mf)\times\Nf=\coprod_{q\in\Nf}\C T\Mf$ is a subbundle of $\C T(\Mf\times\Nf)$ whose (complex) rank equals to the (real) dimension of $\Mf$. 

We consider a $\Co^{\alpha,\beta}$ involutive complex tangential subbundle $\Se$ of $\Mf\times\Nf$ such that $\Se\le(\C T\Mf)\times\Nf$ and $\Se+\bar\Se=(\C T\Mf)\times \Nf$. Note that $(\C T\Mf)\times\Nf$ is automatically involutive, so $\Se$ is a complex Frobenius structure.

We recall the space $\Co^{\alpha,\beta}_\loc(\Mf,\Nf)$ ($\alpha>\frac12$, $\beta>0$) in Definition \ref{Defn::ODE::MixHoldMaps}, the $\Co^{\alpha,\beta}$-subbundle $\Se$ in Definition \ref{Defn::ODE::CpxPaSubbd}, and the involutivity of $\Se\le(\C T\Mf)\times\Nf$ in Definition \ref{Defn::ODE::InvMix} (also see Remark \ref{Rmk::ODE::InvMix}). We endow $\R^r$, $\C^m$ and $\R^q$ with standard (real and complex) coordinate system $t=(t^1,\dots,t^r)$, $z=(z^1,\dots,z^m)$ and $s=(s^1,\dots,s^q)$ respectively. 

To shorten the indices, we set for $\gamma\in\R_+$,
\begin{equation}\label{Eqn::EllipticPara::Betas}
\begin{aligned}
    \gamma^{\sim\alpha}:=\min\left(\gamma,\alpha+1\right),\quad \gamma^{\wedge\alpha}:=\min\left(\gamma-,\alpha+1\right),&\quad\text{for }\alpha\in(1,\infty);\\
    \gamma^{\sim\alpha}=\gamma^{\wedge\alpha}:=\min\left(\gamma-,(2-\tfrac1\alpha)\gamma,\alpha+1\right),&\quad\text{for }\alpha\in(\tfrac12,1];
    \\
    (\gamma-)^{\sim\alpha}=(\gamma-)^{\wedge\alpha}:=\min\left(\gamma-,(2-\tfrac1\alpha)\gamma-,\alpha+1\right),&\quad\text{for }\alpha\in(\tfrac12,\infty).
\end{aligned}
\end{equation}
Here for $\beta\in\{\gamma,\gamma-:\gamma\in\R_+\}$, $\beta^{\sim\alpha}$ stands for the parameter regularity of the coordinate charts, and $\beta^{\wedge\alpha}$ stands for the parameter regularity of the coordinate vector fields.

\begin{keythm}[Elliptic structures with parameter]\label{KeyThm::EllipticPara}
Let $r,m,q\ge0$, let $\Mf$ and $\Nf$ be two smooth manifolds with dimensions $(r+2m)$ and $q$ respectively. Let $\alpha,\beta\in\R_\Eb$ be two indices such that either
\begin{itemize}[nolistsep]
    \item $\alpha\in(\frac12,1)\cup(1,\infty)$ and $\beta\in\R_+$; or
    \item $\alpha\in(\frac12,1]$ and $\beta\in\R_+^-(=\{\gamma-:\gamma\in(0,\infty)\})$.
\end{itemize}
Let $\Se\le (\C T\Mf)\times \Nf$ be a $\Co^{\alpha,\beta}$-involutive subbundle of rank $r+m$ such that $\Se+\bar\Se=(\C T\Mf)\times \Nf$.

Then for any $u_0\in \Mf$ and $v_0\in\Nf$, there are neighborhoods $U\subseteq\Mf$ of $u_0$, $V\subseteq\Nf$ of $v_0$ and a map $F=(F',F'',F'''):U\times V\to\R^r_t\times\C^m_z\times\R^q_s$ such that 
\begin{enumerate}[parsep=-0.3ex,label=(\arabic*)]
    \item\label{Item::EllipticPara::F'''}$F''':V\to\R^q_s$ is a smooth coordinate chart of $\Nf$ that is independent of $U\subseteq\Mf$, such that $F'''(v_0)=0$; $F':U\to\R^r_t$ is a smooth map that is independent of $V\subseteq\Nf$, such that $F'(u_0)=0$.
    \item\label{Item::EllipticPara::FBase}$F:U\times V\to F(U\times V)$ is homeomorphism. And for each $v\in V$, $[u\mapsto (F'(u),F''(u,v))]:U\to\R^r_t\times\C^m_z$ is a $\Co^{\alpha+1}$-coordinate chart.
    \item\label{Item::EllipticPara::F''Reg} $F''(u_0,v_0)=0$, $F''\in\Co^{\alpha+1,\beta^{\sim\alpha}}_\loc(U,V;\C^m)$ and $XF''\in\Co^{\alpha,\beta^{\wedge\alpha}}_\loc(U,V;\C^m)$ for every $C^\infty$ vector field $X$ on $\Mf$.
    \item\label{Item::EllipticPara::Span} $F^*\Coorvec{t^1},\dots,F^*\Coorvec{t^r},F^*\Coorvec{z^1},\dots,F^*\Coorvec{z^m}$ are well-defined $\Co^{\beta^{\wedge\alpha}}$-vector fields and span $\Se|_{U\times V}$.
\end{enumerate}

For the parameterization side, let $\Phi=(\Phi',\Phi''):F(U\times V)\to\Mf\times\Nf$ be the inverse map of $F$, and let $\Omega:=\Omega'\times\Omega''\times\Omega'''\subseteq\R^r_t\times\C^m_z\times\R^q_s$ be any neighborhood of $(0,0,0)$ which is contained in $F(U\times V)$. Then
\begin{enumerate}[parsep=-0.3ex,label=(\arabic*)]\setcounter{enumi}{4}
    \item\label{Item::EllipticPara::Phi''} $\Phi'':V\to \Nf$ is a smooth regular parameterization. Moreover $\Phi''=(F''')^\Inv$.
    \item\label{Item::EllipticPara::Phi0}  $\Phi:\Omega\to\Mf\times\Nf$ is homeomorphic onto its image,  and $\Phi(0,0,0)=(u_0,v_0)$. For each $s\in\Omega''$, $\Phi'(\cdot,s):\Omega'\times\Omega''\to\Mf$ is a $\Co^{\alpha+1}$-regular parameterization.
    \item\label{Item::EllipticPara::PhiReg}  $\Phi'\in\Co^{\alpha+1,\beta^{\sim\alpha}}_\loc(\Omega'\times\Omega'',\Omega''';\Mf)$ and $\frac{\partial\Phi'}{\partial t^j},\frac{\partial\Phi'}{\partial z^k},\frac{\partial\Phi'}{\partial \bar z^k}\in\Co^{\alpha,\beta^{\wedge\alpha}}_\loc(\Omega'\times\Omega'',\Omega''';\C T\Mf)$ for $1\le j\le r$ and $1\le k\le m$.
    \item\label{Item::EllipticPara::PhiSpan}  $\Se_{\Phi(t,z,s)}\le\C T_{\Phi'(t,z,s)}\Mf$ is spanned by $\frac{\partial\Phi}{\partial t^1}(t,z,s),\dots,\frac{\partial\Phi}{\partial t^r}(t,z,s),\frac{\partial\Phi}{\partial z^1}(t,z,s),\dots,\frac{\partial\Phi}{\partial z^m}(t,z,s)$, for every $(t,z,s)\in\Omega$.
\end{enumerate}
In particular,
\begin{enumerate}[parsep=-0.3ex,label=(\arabic*)]\setcounter{enumi}{8}
    \item\label{Item::EllipticPara::>1} When $\beta^{\sim\alpha}>1$, i.e. $\min(\beta,(2-\frac1\alpha)\beta,\alpha+1)>1$, $F$ is a $\Co^{\min(\beta,(2-\frac1\alpha)\beta,\alpha+1)}$-coordinate chart and $\Phi$ is a $\Co^{\min(\beta,(2-\frac1\alpha)\beta,\alpha+1)}$-regular parameterization.
\end{enumerate}

\end{keythm}
\begin{remark}\label{Rmk::EllipticPara::KeyFirst}
\begin{enumerate}[parsep=-0.3ex,label=(\roman*)]

\item When $\alpha>1$ we shall see in Section \ref{Section::Sharpddz} that the index ``$\beta-$'' for  $XF''\in\Co^{\alpha,\beta-}$ and $F^*\Coorvec z\in\Co^{\min(\alpha,\beta-)}$ are both sharp.

\item Theorem \ref{KeyThm::EllipticPara} implies Street's result \cite[Theorem 1.1]{SharpElliptic} for the sharp estimate on elliptic structure by considering $\alpha=\beta>1$ and forgetting the $s$-parameter. (Recall $\Co^{\alpha,\alpha}(U,V)=\Co^\alpha(U\times V)$ from Lemma \ref{Lem::Hold::CharMixHold}.)
\end{enumerate}
\end{remark}


\begin{remark}\label{Rmk::EllipticPara::KeySpecial}
For $\alpha>1$, our results allow the range $\beta\in(0,1]$ and $\beta\in(\alpha,\alpha+1]$. 

For the case $\beta\le1$, it does not make sense to talk about $\Co^\beta$-coordinate charts nor $\Co^\beta$-parameterizations. Nevertheless by \ref{Item::EllipticPara::FBase} we have 
$F^*\Coorvec t|_{(u,v)}=(F',F''(\cdot,v))^*\Coorvec t|_u $ and $F^*\Coorvec z|_{(u,v)}=(F',F''(\cdot,v))^*\Coorvec z|_u $, since $(F',F''(\cdot,v))$ is a coordinate chart for a fix $v\in V$. Alternatively by \ref{Item::EllipticPara::Phi''} $F^*\Coorvec t=\frac{\partial\Phi}{\partial t}\circ F$ and $F^*\Coorvec z=\frac{\partial\Phi}{\partial t}\circ F$ are both continuous. Either way $F^*\Coorvec t,F^*\Coorvec z$ are collections of pointwise defined vector fields.

For the case $\beta>\alpha$, the parameter part is actually more regular than the part of spatial variable. In this case $F^*\Coorvec t,F^*\Coorvec z$ are $\Co^\alpha$-vector fields on $U\times V\subseteq\Mf\times\Nf$.
\end{remark}
\begin{remark}\label{Rmk::EllipticPara::KeyImpliesGong}
If we only consider a family of elliptic structures rather than a family of complex Frobenius structures, and assume $\alpha,\beta$ are non-integers, then Theorem \ref{KeyThm::EllipticPara} is much stronger than \cite[Theorem 1.3]{Gong}:
\begin{itemize}[parsep=-0.3ex]
    \item We consider all $(\alpha,\beta)$ such that $\alpha>1$ and $0<\beta\le\alpha+1$. In \cite[Theorem 1.3]{Gong} only the case $\alpha>\beta+3\ge4$ and the case $\alpha=\beta$ is considered.
    \item For the regularity of the coordinate charts we have the sharp estimate $F\in\Co^{\alpha+1,\beta}(U,V)$. \cite[Theorem 1.3]{Gong} only gives $F\in\Co^{\alpha-,\beta-}(U,V)$.
\end{itemize}
\end{remark}



We give a detailed overview to the idea of its proof in Section \ref{Section::EllipticPara::Overview} using Malgrange's factorization.

\subsection{Overview of the estimate}\label{Section::EllipticPara::Overview}
 For vector fields $X_1,\dots,X_{r+m}$, we write $X=[X_1,\dots,X_{r+m}]^\top$ following from the convention in Section \ref{Section::Convention}.

The most important part is to find the map $F'':U\times V\to\C^m$ that has the desired regularity property. We factorize $F''$ to the following composition of maps:
\begin{equation}\label{Eqn::EllipticPara::FacterizationofF}
    F''=\overline{G''}\circ H\circ L.
\end{equation}

Endow $\R^r,\C^m,\R^q$ with standard coordinates $\tau=(\tau^1,\dots,\tau^r)$, $w=(w^1,\dots,w^m)$ and $s=(s^1,\dots,s^q)$ respectively.
There are smooth coordinate charts $(L',L''):U_0\subseteq\Mf\xrightarrow{\sim}\B^{r+2m}\subset\R^r_\tau\times\C^m_w$ and $L''':V_0\subseteq\Nf\xrightarrow{\sim}\B^q\subset\R^q_s$, such that by taking $L=(L',L'',L'''):U_0\times V_0\to\B^{r+2m}\times\B^q$,  $L_*\Se$ has a $\Co^{\alpha,\beta}_{(\tau,w),s}$-local basis on $\B^{r+2m}\times\B^q$ with the form 
\begin{equation}\label{Eqn::EllipticPara::DefofA}
    X=\begin{pmatrix}X'\\X''\end{pmatrix}=\begin{pmatrix}I_r&0&A'\\0&I_m&A''\end{pmatrix}\begin{pmatrix}\partial_\tau\\\partial_w\\\partial_{\bar w}\end{pmatrix},\quad A:=\begin{pmatrix}
    A'\\A''\end{pmatrix},
\end{equation}
where $\|A\|_{\Co^{\alpha,\beta}(\B^{r+2m},\B^q;\C^{(r+m)\times m})}$ is small enough. This can be done by some scaling argument, see Lemma \ref{Lem::EllipticPara::IniNorm}.

For simplicity we only illustrate the case $\alpha>1$, where we can assume $0<\beta\le\alpha+1$. When $\frac12<\alpha<1$ the argument is the same except we replace the indices $\beta$ and $\beta-$ by $(2-\frac1\alpha)\beta$ and the space $\Co^{\alpha+1}L^\infty\cap\Co^1\Co^\beta$ by $\Co^{\alpha+1}L^\infty\cap\Co^1\Co^\beta+\Co^{\alpha+1}\Co^{(2-\frac1\alpha)\beta}$.

There is a ``$\Co^{\alpha+1,\beta}_{(\tau,w),s}$-coordinate change'' $(\sigma,\zeta,s)=H(\tau,w,s)$ where $H:\B^{r+2m}_{\tau,w}\times\B^q_s\to\B^{r+2m}_{\sigma,\zeta}\times\B^q_s$ is of the form $H(\tau,w,s)=(\tau,H''(\tau,w,s),s)$, such that the local generators for $(H\circ L)_*\Se$ are $\Co^\beta$ and  real-analytic in $(\sigma,\zeta)$-variable. More precisely, the canonical local basis of $(H\circ L)_*\Se$ has the the form
\begin{equation}\label{Eqn::EllipticPara::DefofB}
\begin{pmatrix}T\\Z\end{pmatrix}=\begin{pmatrix}I_r&0&B'\\0&I_m&B''\end{pmatrix}\begin{pmatrix}\partial_\sigma\\\partial_\zeta\\\partial_{\bar \zeta}\end{pmatrix},\ B= \left.\begin{pmatrix}B'\\B''\end{pmatrix}\right|_{(\sigma,\zeta,s)}:=\left.\begin{pmatrix}I&H''_\tau+ A'H''_{\bar w}\\&H''_ w+ A''H''_{\bar w}\end{pmatrix}^{-1}\!\!\begin{pmatrix}\bar H''_\tau+ A'\bar H''_{\bar w}\\\bar H''_ w+ A''\bar H''_{\bar w}\end{pmatrix}\right|_{H^\Inv(\sigma,\zeta,s)}.
\end{equation}
Here $B'$ and $B''$ are defined on $\B^{r+2m}_{\sigma,\zeta}\times\B^q_s\subset(\R^r\times\C^m)\times\R^q$ and admit $(\sigma,\zeta)$-variables holomorphic extensions to the domain $\frac12\Hb^{r+2m}\times\frac12\B^q\subset\C^{r+2m}\times\R^q$.

In this case $H\in\Co^{\alpha+1}_{\tau,w}L^\infty_s\cap\Co^1_{\tau,w}\Co^\beta_s$, note that $H''_\tau,H''_{\bar w}\in\Co^{\alpha}_{\tau,w}L^\infty_s\cap\Co^0_{\tau,w}\Co^\alpha_s\subset\Co^{\alpha,\beta-}_{(\tau,w),s}$ has a regularity loss on the parameter, while $T,Z\in\Co^\infty_{\sigma,\zeta}\Co^\beta_s\subset\Co^{\infty,\beta}_{(\sigma,\zeta),s}$ retrieve the full $\Co^\beta$ regularity on the parameter. 

We denote by $\Tf$ and $\Zf$ the  collections of holomorphic vector fields corresponding to $T$ and $Z$. Finally, applying Proposition \ref{Prop::ODE::ParaHolFro}, the holomorphic Frobenius theorem with parameter, on $\Tf,\Zf$ we can find a $\Co^{\infty}_{\sigma,\zeta}\Co^\beta_s$-map $G''(\sigma,\zeta,s)$ such that $(H\circ L)_*\Se^\bot$ is spanned by $dG'',d\bar\zeta,ds$. See Lemma \ref{Lem::EllipticPara::FinalModf}.

In this way, the composition map $\overline{G''}\circ H\circ L$ is as desired.

\medskip To make $B(\sigma,\zeta,s)$ real-analytic, we impose extra conditions on $H$: we require $H''\in\Co^{\alpha+1}_{\tau,w}L^\infty_s\cap\Co^1_{\tau,w}\Co^\beta_s$, and $H$ satisfies the following condition:
\begin{equation}\label{Eqn::EllipticPara::KeyEqn1}
\sum_{j=1}^r\frac{\partial T_j}{\partial \sigma^j}+\sum_{k=1}^m\frac{\partial Z_k}{\partial \bar \zeta^k}=0.
\end{equation}

Rewriting \eqref{Eqn::EllipticPara::KeyEqn1} in the $(\tau,w,s)$-space using the coordinate change $H$, we get: for $l=1,\dots,m$, $(\tau,w)\in\B^{r+2m}$ and $s\in\B^q$,
\begin{equation}\label{Eqn::EllipticPara::ExistenceH}\begin{aligned}
		&\sum_{j=1}^r\left(\begin{pmatrix}I&H''_\tau&\bar H''_\tau\\&H''_ w&\bar H''_ w\\&H''_{\bar w}&\bar H''_{\bar w}\end{pmatrix}^{-1}\!\!\begin{pmatrix}\partial_{\tau}\\\partial_{ w}\\\partial_{\bar w}\end{pmatrix}\right)_j\left(\begin{pmatrix}I&H''_\tau+ A'H''_{\bar w}\\&H''_ w+ A''H''_{\bar w}\end{pmatrix}^{-1}\!\!\begin{pmatrix}\bar H''_\tau+ A'\bar H''_{\bar w}\\\bar H''_ w+ A''\bar H''_{\bar w}\end{pmatrix}\right)_j^l\Bigg|_{(\tau, w;s)}
		\\
		+&\sum_{k=1}^m\left(\begin{pmatrix}I&H''_\tau&\bar H''_\tau\\&H''_ w&\bar H''_ w\\&H''_{\bar w}&\bar H''_{\bar w}\end{pmatrix}^{-1}\!\!\begin{pmatrix}\partial_{\tau}\\\partial_{ w}\\\partial_{\bar w}\end{pmatrix}\right)_{k+r+m}\!\!\left(\begin{pmatrix}I&H''_\tau+ A'H''_{\bar w}\\&H''_ w+ A''H''_{\bar w}\end{pmatrix}^{-1}\!\!\begin{pmatrix}\bar H''_\tau+ A'\bar H''_{\bar w}\\\bar H''_ w+ A''\bar H''_{\bar w}\end{pmatrix}\right)_{k+r}^l\Bigg|_{(\tau, w;s)}=0.
		\end{aligned}
\end{equation}
See Lemma \ref{Lem::EllipticPara::NewGen} for this deduction and Proposition \ref{Prop::EllipticPara::ExistPDE} for the existence of such $H$.

We write $B'=(b_j^l)_{\substack{1\le j\le r\\1\le l\le m}}$ and $B''=(b_{r+k}^l)_{\substack{1\le k\le m\\1\le l\le m}}$.
By Lemma \ref{Lem::ODE::GoodGen} \ref{Item::ODE::GoodGen::InvComm} the involutivity of $\Se$ implies that $T_1,\dots,T_r,Z_1,\dots,Z_m$ are commutative. Writing out the conditions $[T_j,T_{j'}]=0$, $[T_j,Z_k]=0$ and $[Z_k,Z_{k'}]=0$ in terms of $B'$ and $B''$ and combining \eqref{Eqn::EllipticPara::KeyEqn1}, we have the following:
\begin{equation}\label{Eqn::EllipticPara::KeyEqnB}
\begin{aligned}
\frac{\partial b_{j'}^l}{\partial \sigma^j}-\frac{\partial b_j^l}{\partial \sigma^{j'}}&=\sum_{q=1}^m\Big(b_{j'}^q\frac{\partial b_j^l}{\partial\bar  \zeta^q}-b_j^q\frac{\partial b_{j'}^l}{\partial\bar \zeta^q}\Big),&1\le j<j'\le r,1\le l\le m;\\
\frac{\partial b_{r+k}^l}{\partial \sigma^j}-\frac{\partial b_j^l}{\partial \zeta^k}&=\sum_{q=1}^m\Big(b_{r+k}^q\frac{\partial b_j^l}{\partial\bar \zeta^q}-b_j^q\frac{\partial b_{r+k}^l}{\partial\bar \zeta^q}\Big),&1\le j\le r,1\le k,l\le m;\\
\frac{\partial b_{r+k'}^l}{\partial \zeta^k}-\frac{\partial b_{r+k}^l}{\partial \zeta^{k'}}&=\sum_{q=1}^m\Big(b_{r+k'}^q\frac{\partial b_{r+k}^l}{\partial\bar \zeta^q}-b_{r+k}^q\frac{\partial b_{r+k'}^l}{\partial\bar \zeta^q}\Big),&1\le k<{k'}\le m,1\le l\le m;\\
\sum_{j=1}^r\frac{\partial b_j^l}{\partial\sigma^j}+\sum_{k=1}^m\frac{\partial b_{r+k}^l}{\partial\bar\zeta^k}&=0,&1\le l\le m.
\end{aligned}\end{equation}

\eqref{Eqn::EllipticPara::KeyEqnB} is a quasilinear PDE system of the form $DB=\Theta(B,\nabla_{\sigma,\zeta}B)$, where $D$ is a first order differential operator with constant (vector-valued) coefficients and $\Theta$ is a product map with constant (vector-valued) coefficients. We can check that $D^*D=\Delta_\sigma+\Box_\zeta$ is the real and complex Laplacian acting on the components. So \eqref{Eqn::EllipticPara::KeyEqnB} is an elliptic PDE system and in Proposition \ref{Prop::EllipticPara::AnalyticPDE} we will see that $B$ admits a $(\sigma,\re\zeta,\im\zeta)$-variable holomorphic extension. See Lemma \ref{Lem::EllipticPara::BSatAna}.

\subsection{The existence proposition}\label{Section::EllipticPara::ExistPDE}
In this subsection we use $\tau=(\tau^1,\dots,\tau^r)$, $w=(w^1,\dots,w^m)$, $s=(s^1,\dots,s^q)$ as the standard (real or complex) coordinate system for $\R^r$, $\C^m$ and $\R^q$ respectively. We consider the unit balls $\B^{r+2m}\subset\R^r_\tau\times\C^m_w\simeq\R^{r+2m}$ and $\B^q\subset\R^q_s$.

\begin{conv}\label{Conv::EllipticPara::ConvofExtPDE} 

We denote $\Ic:\B^{r+2m}\times\B^q\to\B^{2m}$, $\Ic(\tau,w,s):=w$ as the projection map. 


For a matrix $B=(b_j^k)$, we use $b_j^k$ to denote the $j$-th row $k$-th column entry of $B$.

We use $A=\begin{pmatrix}	A'\\A''	\end{pmatrix}:\B^{r+2m}_{\tau,w}\times \B^q_s\to\C^{(r+m)\times m}$ as a matrix map, $H''=(H''^1,\dots,H''^m):\B^{r+2m}_{\tau,w}\times \B^q_s\to \C^m$ as a vector-valued map whose image is set to be a $m$-dimensional row vector. We use 
\begin{equation}\label{Eqn::EllipticPara::GradH}
    H''_\tau=\begin{pmatrix}\frac{\partial H''^1}{\partial \tau^1}&\cdots&\frac{\partial H''^m}{\partial \tau^1}\\\vdots&\ddots&\vdots\\\frac{\partial H''^1}{\partial \tau^r}&\cdots&\frac{\partial H''^m}{\partial \tau^r}\end{pmatrix},\quad H''_w=\begin{pmatrix}\frac{\partial H''^1}{\partial w^1}&\cdots&\frac{\partial H''^m}{\partial w^1}\\\vdots&\ddots&\vdots\\\frac{\partial H''^1}{\partial w^m}&\cdots&\frac{\partial H''^m}{\partial w^m}\end{pmatrix},\quad\nabla_{\tau,w} H''=\begin{pmatrix}H''_\tau\\H''_w\\H''_{\bar w}\end{pmatrix}
\end{equation}

Let $\Psi,\Theta,\Lambda$ be matrix-valued functions of $H''$ and $A$ as follow:
\begin{equation}\label{Eqn::EllipticPara::matrixfun}
\Psi[H'']:=\begin{pmatrix}I&H''_\tau&\bar H''_\tau\\&H''_ w&\bar H''_ w\\&H''_{\bar w}&\bar H''_{\bar w}\end{pmatrix},\ \Theta[A;H'']:=\begin{pmatrix}I&H''_\tau+ A'H''_{\bar w}\\&H''_ w+ A''H''_{\bar w}\end{pmatrix},\  \Lambda[A;H'']:=\Theta[A;H'']^{-1}\begin{pmatrix}\bar H''_\tau+ A'\bar H''_{\bar w}\\\bar H''_ w+ A''\bar H''_{\bar w}\end{pmatrix}.\end{equation}

\end{conv}

\begin{remark}One can see that $\Psi,\Theta,\Lambda$ are all rational functions in the components of $A$ and $\nabla_{\tau,w}H''$.
\end{remark} 

We are going to consider the following PDE system, whose deduction is in Lemma \ref{Lem::EllipticPara::NewGen}.

Recall $\partial_\tau=[\partial_{\tau^1},\dots,\partial_{\tau^r}]^\top$ in Section \ref{Section::Convention}. In this subsection we prove the following: 
\begin{prop}\label{Prop::EllipticPara::ExistPDE}
	Let $\Xs^1,\Xs^0,\Xs^{-1}$ be bi-parameter H\"older-Zygmund structures such that:
	\begin{enumerate}[parsep=-0.3ex,label=(X.\arabic*)]
	\item\label{Item::EllipticPara::XsAssume::C0} $\Xs^1(\B^{r+2m},\B^q)\subseteq\Xs^0(\B^{r+2m},\B^q)\subseteq C^0(\B^{r+2m},\B^q)\cap\Xs^{-1}(\B^{r+2m},\B^q)$.
    \item\label{Item::EllipticPara::XsAssume::Grad} $\nabla_{\tau,w}:\Xs^1(\B^{r+2m},\B^q)\to\Xs^0(\B^{r+2m},\B^q;\R^{r+2m})$ and $\nabla_{\tau,w}:\Xs^0(\B^{r+2m},\B^q)\to\Xs^{-1}(\B^{r+2m},\B^q;\R^{r+2m})$ are bounded.
    \item\label{Item::EllipticPara::XsAssume::Delta} $\Delta_\tau+\Box_w=\sum_{j=1}^r\partial_{\tau^j}^2+\sum_{k=1}^m\partial_{w^k\bar w^k}^2:\{f\in\Xs^1(\B^{r+2m},\B^q):f|_{\partial\B^{r+2m}\times\B^q}=0\}\to\Xs^{-1}(\B^{r+2m},\B^q)$ is invertible.
    \item\label{Item::EllipticPara::XsAssume::Mult} The product maps $\Xs^0(\B^{r+2m},\B^q)\times\Xs^{-1}(\B^{r+2m},\B^q)\to\Xs^{-1}(\B^{r+2m},\B^q)$ and $\Xs^0(\B^{r+2m},\B^q)\times\Xs^0(\B^{r+2m},\B^q)\to\Xs^0(\B^{r+2m},\B^q)$ are bilinearly bounded.
	\end{enumerate}
	
	Using Convention \ref{Conv::EllipticPara::ConvofExtPDE}, there exists a $\delta_0=\delta_0(\Xs,r,m,q)>0$ and $C_0=C_0(\Xs,r,m,q)>1$, such that if $A\in\Xs^0(\B^{r+2m},\B^q;\C^{(r+m)\times m})$ satisfies $\|A\|_{\Xs^0}<\delta_0$, then there is a unique $H''\in\Xs^1(\B^{r+2m}, \B^q;\C^m)$ that solves the following PDE system \eqref{Eqn::EllipticPara::ExistenceH} and satisfies the boundary condition $H''|_{(\partial\B^{r+2m})\times\B^q}=\Ic$, and satisfies $\|H''-\Ic\|_{\Xs^1}<C_0\|A\|_{\Xs^0}$.
\end{prop}

Before we prove Proposition \ref{Prop::EllipticPara::ExistPDE}, we need to show the PDE \eqref{Eqn::EllipticPara::ExistenceH} is well-defined, in the sense that $\Psi[H''](\tau,w,s)$ and $\Theta[A;H''](\tau,w,s)$ are invertible matrices for all $(\tau,w,s)\in\B^{r+2m}\times \B^q$. 

For $\Xs^1,\Xs^0,\Xs^{-1}$ that satisfies the assumptions of Proposition \ref{Prop::EllipticPara::ExistPDE}, and for $\eps>0$, we define spaces
	\begin{equation}\label{Eqn::EllipticPara::MatisInv::SpaceUV}
	    \begin{aligned}
    &\U_{\Xs,\eps}:=\{A\in\Xs^0(\B^{r+2m},\B^q;\C^{(r+m)\times m}):\|A\|_{\Xs^0}<\eps\},\\
    &\V_{\Xs,\eps}:=\{H''\in\Xs^1(\B^{r+2m},\B^q;\C^m):H''|_{(\partial\B^{r+2m})\times\B^q}\equiv\Ic,\ \|H''-\Ic\|_{\Xs^1}<\eps\},
        \end{aligned}
	\end{equation}
	with metrics induced by the norm structures of their respective ambient Banach spaces.
\begin{lem}\label{Lem::EllipticPara::MatisInv}
    There is a $\eps_1=\eps_1(r,m,q,\Xs)>0$ such that
    \begin{enumerate}[parsep=-0.3ex,label=(\roman*)]
        \item\label{Item::EllipticPara::MatisInv::Inv}For every $A\in\U_{\Xs,\eps_1}$, $H''\in \V_{\Xs,\eps_1}$, $(\tau,w)\in\B^{r+2m}$ and $s\in\B^q$, the matrices $\Psi[H''](\tau,w,s)$ and $\Theta[A;H''](\tau,w,s)$ given in \eqref{Eqn::EllipticPara::matrixfun} are both invertible.
        \item\label{Item::EllipticPara::MatisInv::BddLambda} For every $A\in\U_{\Xs,\eps_1}$ and $H''\in\V_{\Xs,\eps_1}$,
        \begin{equation}\label{Eqn::EllipticPara::MatisInv::BddLambda}
	    \|\Lambda[A;H'']\|_{\Xs^0}\le\eps_1^{-1}(\|H''-\Ic\|_{\Xs^1}+\|A\|_{\Xs^0}).
	\end{equation}
	
	    \item\label{Item::EllipticPara::MatisInv::Diff} Denote by $T[A;H'']^l$ the left hand side of \eqref{Eqn::EllipticPara::ExistenceH} for $l=1,\dots,m$, and $T[A;H'']:=(T[A;H'']^l)_{l=1}^m$ as the vector-valued function. Then $T:\U_{\Xs,\eps_1}\times\V_{\Xs,\eps_1}\to \Xs^{-1}(\B^{r+2m},\B^q;\C^m)$ is a real-analytic map, in particular second order Fr\'echet differentiable. Moreover
	    \begin{equation}\label{Eqn::EllipticPara::MatisInv::TangT}
	        \frac{\partial T}{\partial H''}[A;H'']\Big|_{A=0;H''=\Ic}=\Delta_\tau+\Box_w:\{F\in\Xs^1(\B^{r+2m},\B^q;\C^m):F|_{(\partial\B^{r+2m})\times\B^q}=0\}\to\Xs^{-1}(\B^{r+2m},\B^q;\C^m).
	    \end{equation}

	    Here $\Delta_\tau+\Box_w=\sum_{j=1}^r\frac{\partial^2}{\partial\tau^j\partial\tau^j}+\sum_{k=1}^m\frac{\partial^2}{\partial w^k\partial\bar w^k}$ is the mixed real and complex Laplacian.
    \end{enumerate}
\end{lem}
\begin{proof} We let $\eps_1$ be a small constant which may change from line to line.

Note that when $H''=\Ic$ and $A=0$, we have
$$\nabla_{\tau,w}\Ic=\begin{pmatrix}0_{r\times m}\\I_m\\0_{m\times m}\end{pmatrix},\quad\Psi[\Ic]=\begin{pmatrix}I_r\\&I_m\\&&I_m\end{pmatrix}=I_{r+2m},\quad\Theta[0;\Ic]=\begin{pmatrix}I_r+ 0_{r\times r}&0_{r\times m}\\0_{m\times r}&I_m+0_{m\times m}\end{pmatrix}=I_{r+m}.$$

By cofactor representations of the matrices, $\Psi^{-1}$, $\Theta^{-1}$ and $\Lambda$ are rational functions in the components of $A$, $\nabla_{\tau,w}H''$ and $\nabla_{\tau,w}\bar H''$, which can be expressed as power series expansions in the components of $A$, $\nabla_{\tau,w}H''$ and $\nabla_{\tau,w}\bar H''$ at the point $A=0$ and $\nabla_{\tau,w}H''=\nabla_{\tau,w}\Ic$ and $\nabla_{\tau,w}\bar H''=\nabla_{\tau,w}\Ic$.

By Lemma \ref{Lem::Hold::CramerMixedPara}, inverting a  $\Xs^0$-matrix is a real-analytic map near the identity matrix. Therefore there is a $\tilde\eps_1>0$ such that if $\|A\|_{\Co^\alpha L^\infty\cap\Co^0\Co^\beta}<\tilde\eps_1$ and $\|\nabla_{\tau,w}(H''-\Ic)\|_{\Co^\alpha L^\infty\cap\Co^0\Co^\beta}<\tilde\eps_1$, then $\Psi[H'']$ and $\Theta[A;H'']$ are invertible matrices at every point in $\B^{r+2m}\times\B^q$, and moreover as maps between $\Xs^0$-spaces,
\begin{align*}
&\Psi^{-1}:\{\|\nabla_{\tau,w}(H''-\Ic)\|_{\Xs^0}<\tilde\eps_1\}\to\Xs^0(\B^{r+2m},\B^q;\C^{(r+2m)\times(r+2m)}),
\\
    &\Lambda:\{\|A\|_{\Xs^0}<\tilde\eps_1\}\times\{\|\nabla_{\tau,w}(H''-\Ic)\|_{\Xs^0}<\tilde\eps_1\}\to\Xs^0(\B^{r+2m},\B^q;\C^{(r+m)\times m}),
\end{align*}
are both real-analytic.

By the assumption \ref{Item::EllipticPara::XsAssume::Grad} there is a  $\tilde\eps_2<C^{-1}\tilde\eps_1$ where $C=C(\Xs,r,m,q)>0$ is a large constant, such that $H''\in\V_{\Xs,\tilde\eps_2}$ implies $\|\nabla_{\tau,w}(H''-\Ic)\|_{\Xs^0}<\tilde\eps_1$.

Thus by \ref{Item::EllipticPara::XsAssume::C0}, when $\eps\le\tilde\eps_2$, $\Psi[H'']$ and $\Theta[A;H'']$ are pointwise invertible in $\B^{r+2m}\times\B^q$. This finishes the proof of \ref{Item::EllipticPara::MatisInv::Inv}. Moreover we have the following real-analytic maps:
\begin{equation}\label{Eqn::EllipticPara::MatisInv::PsiLambdaDiff}
    \begin{aligned}
    \Psi^{-1}:\U_{\Xs,\tilde\eps_2}\to\Xs^0(\B^{r+2m},\B^q;\C^{(r+2m)\times (r+2m)}),\quad
    \Lambda:\U_{\Xs,\tilde\eps_2}\times\V_{\Xs,\tilde\eps_2}\to\Xs^0(\B^{r+2m},\B^q;\C^{(r+2m)\times m}).
\end{aligned}
\end{equation}

Since $\Psi[\Ic]=I_{r+2m}$ and $\Lambda[0;\Ic]=0$, taking the first order terms of the power expansion we have
\begin{gather}\label{Eqn::EllipticPara::MatisInv::PfBddLambda}
\|\Psi[H'']^{-1}-I_{r+2m}\|_{\Xs^0}\lesssim\|\nabla_{\tau,w}H''-\nabla_{\tau,w}\Ic\|_{\Xs^0},\quad
    \|\Lambda[A;H'']\|_{\Xs^0}\lesssim\|A\|_{\Xs^0}+\|\nabla_{\tau,w}H''-\nabla_{\tau,w}\Ic\|_{\Xs^0},
\end{gather}
whenever $A\in\U_{\Xs,\tilde\eps_2}$ and $H''\in\V_{\Xs,\tilde\eps_2}$. Therefore by \eqref{Eqn::EllipticPara::MatisInv::PfBddLambda} with possibly shrinking $\eps_1$, we get \eqref{Eqn::EllipticPara::MatisInv::BddLambda} and prove \ref{Item::EllipticPara::MatisInv::BddLambda}.

\medskip
We now prove \ref{Item::EllipticPara::MatisInv::Diff}. From the expression \eqref{Eqn::EllipticPara::ExistenceH}, $T[A;H'']$ is the linear combinations to the components of $\Psi[H'']^{-1}\otimes\nabla_{\tau,w}\Lambda[A;H'']$. As the maps in \eqref{Eqn::EllipticPara::MatisInv::PsiLambdaDiff}, $\Psi^{-1}$ and $\Lambda$ are both real-analytic and hence Fr\'echet differentiable. Therefore 
$$\big[(A,H'')\mapsto \nabla_{\tau,w}\Lambda[A;H'']\big]:\U_{\Xs,\eps_1}\times\V_{\Xs,\eps_1}\to\Xs^{-1}(\B^{r+2m},\B^q;\C^{(r+m)\times m\times(r+2m)}),$$ is also real-analytic and differentiable.

By \ref{Item::EllipticPara::XsAssume::Mult} that the product map $\Xs^0\times \Xs^{-1}\to \Xs^{-1} $ is bounded hence is differentiable, we see that
$$\big[(A,H'')\mapsto\Psi[H'']^{-1}\otimes \nabla_{\tau,w}\Lambda[A;H'']\big]:\U_{\Xs,\eps_1}\times\V_{\Xs,\eps_1}\to\Xs^{-1}(\B^{r+2m},\B^q;\C^{(r+2m)^2\times m(r+m)(r+2m)}),$$ is also differentiable. Taking linear combinations we get that $T:\U_{\Xs,\eps_1}\times \V_{\Xs,\eps_1}\to\Xs^{-1}$ is Fr\'echet differentiable.

Since $H''\mapsto\Lambda[0;H'']$ is a real-analytic function to the components of $\nabla_{\tau,w}H''$ and $\nabla_{\tau,w}\bar H''$, taking second order power expansion at $\nabla_{\tau,w}H''=\nabla_{\tau,w}\bar H''=\nabla_{\tau,w}\Ic$ we have
\begin{equation*}
    \Lambda[0;H'']=\begin{pmatrix}I_r&\\&I_m\end{pmatrix}\begin{pmatrix}\bar H''_\tau\\\bar H''_w\end{pmatrix}+O(|\nabla_{\tau,w}H''-\nabla_{\tau,w}\Ic|^2).
\end{equation*}

By  Corollary \ref{Cor::Hold::CorMult} \ref{Item::Hold::CorMult::Prin} again, since $\eps_1$ is small, we have
\begin{equation}\label{Eqn::EllipticPara::MatisInv::PfBddLambda2}
    \bigg\|\Lambda[0;H'']-\begin{pmatrix}\bar H''_\tau\\\bar H''_w\end{pmatrix}\bigg\|_{\Xs^0}\lesssim\|\nabla_{\tau,w}(H''-\Ic)\|_{\Xs^0}^2\lesssim\|H''-\Ic\|_{\Xs^1}^2.
\end{equation}

Taking $\nabla_{\tau,w}$ on $\Lambda[0;H'']$, since by \eqref{Eqn::EllipticPara::MatisInv::PfBddLambda} that $\Psi[H'']^{-1}=I_{r+2m}+O(|\nabla_{\tau,w}(H''-\Ic)|)$ near $H''=\Ic$, we have
\begin{align*}
    T[0;H'']^l=&\sum_{j=1}^r\left((I_{r+2m}+O(\nabla_{\tau,w}(H''-\Ic))\begin{pmatrix}\partial_{\tau}\\\partial_{ w}\\\partial_{\bar w}\end{pmatrix}\right)_j\left(\begin{pmatrix}\bar H''_\tau\\\bar H''_ w\end{pmatrix}+O(\nabla_{\tau,w}(H''-\Ic))^2\right)_j^l
		\\
		&+\sum_{k=1}^m\left((I_{r+2m}+O(\nabla_{\tau,w}(H''-\Ic))\begin{pmatrix}\partial_{\tau}\\\partial_{ w}\\\partial_{\bar w}\end{pmatrix}\right)_{k+r+m}\left(\begin{pmatrix}\bar H''_\tau\\\bar H''_ w\end{pmatrix}+O(\nabla_{\tau,w}(H''-\Ic))^2\right)_{k+r}^l
		\\
		=&\sum_{j=1}^r\Coorvec{\tau^j}\frac{\partial\bar H''^l}{\partial\tau^j}+\sum_{k=1}^m\Coorvec{w^k}\frac{\partial\bar H''^l}{\partial \bar w^k}+O(\nabla_{\tau,w}(H''-\Ic))^2+O(\nabla_{\tau,w}(H''-\Ic))O(\nabla^2_{\tau,w}(H''-\Ic))
		\\
		=&(\Delta_\tau+\Box_w)\bar H''+O(\nabla_{\tau,w}(H''-\Ic))^2+O(\nabla_{\tau,w}(H''-\Ic))O(\nabla^2_{\tau,w}(H''-\Ic)).
\end{align*}

Combining \eqref{Eqn::EllipticPara::MatisInv::PfBddLambda}, \eqref{Eqn::EllipticPara::MatisInv::PfBddLambda2}, we get
\begin{equation*}
    \|T[0;H'']-(\Delta_\tau+\Box_w)\bar H''\|_{\Xs^{-1}(\B^{r+2m},\B^q;\C^m)}\lesssim\|H''-\Ic\|_{\Xs^1(\B^{r+2m},\B^q;\C^m)}^2,\quad\forall H''\in\V_{\Xs,\eps_1}.
\end{equation*}
Since $T:\U_{\Xs,\eps_1}\times\V_{\Xs,\eps_1}\to\Xs^{-1}$ is real-analytic, we see that $\|T[0;H'']-T[A;H'']\|_{\Xs^{-1}}\lesssim\|A\|_{\Xs^0}$. Therefore $\|T[A;H'']-(\Delta_\tau+\Box_w)\bar H''\|_{\Xs^{-1}}\lesssim\|A\|_{\Xs^0}+\|H''-\Ic\|_{\Xs^1}$. By taking $\eps_1$ smaller, we  get \eqref{Eqn::EllipticPara::MatisInv::TangT}, finishing the proof of \ref{Item::EllipticPara::MatisInv::Diff}.
\end{proof}


\begin{proof}[Proof of Proposition \ref{Prop::EllipticPara::ExistPDE}]
	By Lemma \ref{Lem::EllipticPara::MatisInv} \ref{Item::EllipticPara::MatisInv::Diff}, $T:\U_{\Xs,\eps_1}\times\V_{\Xs,\eps_1}\to\Xs^{-1}(\B^{r+2m},\B^q,\C^m)$ is a well-defined real-analytic map for some small $\eps_1>0$. 
	
By \ref{Item::EllipticPara::XsAssume::Delta}, \eqref{Eqn::EllipticPara::MatisInv::TangT} is invertible, whose inverse we denote by $\tilde\Pc:\Xs^{-1}\to\Xs^1$. 

For a given $A\in\U_{\Xs,\eps_1}$, we define $\Se^A:\V_{\Xs,\eps_1}\to\Xs^1(\B^{r+2m},\B^q,\C^m)$ by $\Se^A[H'']:=H''-\tilde\Pc T[A;H'']$. Clearly $\Se^0[\Ic]=\Ic$. By Lemma \ref{Lem::EllipticPara::MatisInv} \ref{Item::EllipticPara::MatisInv::Diff}, since $T$ is second order differentiable, $T[0;\Ic]=0$ and $(\Delta_\tau+\Box_w)\Ic=0$, we have, for $A\in\U_{\Xs,\eps_1}$ near $0$ and $H''\in\V_{\Xs,\eps_1}$ near $\Ic$,
\begin{equation*}
    T[A;H'']=0+(\Delta_\tau+\Box_w)(H''-\Ic)+O_{\Xs^1\to\Xs^{-1}}((H''-\Ic)\otimes (H''-\Ic))+O_{\Xs^0}(A).
\end{equation*}

Therefore there is a $C_1=C_1(r,m,q,\Xs,\eps_1)>2$ that does not depend on $A,H''$ such that
\begin{equation}\label{Eqn::EllipticPara::ExistPDE::SeACont}
    \begin{aligned}
    \|\Se^A[H'']-\Ic\|_{\Xs^1}&\le C_1(\|H''-\Ic\|_{\Xs^1}^2+\|A\|_{\Xs^0}),&\forall A\in\U_{\Xs,\eps_1},\ H''\in\V_{\Xs,\eps_1};
    \\
    \|\Se^A[H''_1]-\Se^A[H''_2]\|_{\Xs^1}&\le C_1(\|H_1-\Ic\|_{\Xs^1}+\|H_2-\Ic\|_{\Xs^1})\|H_1-H_2\|_{\Xs^1},&\forall A\in\U_{\Xs,\eps_1},\ H''_1,H''_2\in\V_{\Xs,\eps_1}.
\end{aligned}
\end{equation}

Take $\delta_0:=\min(\frac12\eps_1,\frac1{6C_1^2})$. When $A\in\U_{\Xs,\delta_0}$, by \eqref{Eqn::EllipticPara::ExistPDE::SeACont} $\Se^A:\V_{\Xs,\eta}\to\V_{\Xs,\eta}$ is a contraction map for all $\frac43C_1\|A\|_{\Xs^0}<\eta\le\frac43C_1{\delta_0}$. Therefore there exists a unique fixed point $H''\in\V_{\Xs,\frac43C_1{\delta_0}}$ such that $H''=\Se^A[H'']$. Therefore $H''=H''-\tilde\Pc T[A;H'']$ and thus $T[A;H'']=0$. By uniqueness of the fixed point we see that $H''\in \V_{\Xs,\frac43C_1\|A\|_{\Xs^0}}$, which is $\|H''-\Ic\|_{\Xs^1}<\frac43C_1\|A\|_{\Xs^0}$. Taking $C_0:=\frac43 C_1$, we complete the proof.
\end{proof}

In application to Theorem \ref{KeyThm::EllipticPara} we have the following:
\begin{cor}\label{Cor::EllipticPara::EstofH}
    Let $\alpha\in(\frac12,\infty)\backslash\{1\}$ and $\beta>0$. Then for any $\eps>0$ there is a $\delta=\delta(\eps,\alpha,\beta,r,m,q)>0$ such that: for any $A\in\Co^{\alpha,\beta}(\B^{r+2m},\B^q;\C^{(r+m)\times m})$ that satisfies $\|A\|_{\Co^{\alpha,\beta}}<\delta$, there is a $H''$ that solves \eqref{Eqn::EllipticPara::ExistenceH} with $H|_{(\partial\B^{r+2m})\times\B^q}=\Ic$. Moreover
    \begin{enumerate}[parsep=-0.3ex,label=(\roman*)]
        \item\label{Item::EllipticPara::EstofH::HReg}When $\alpha>1$, $H''\in\Co^{\alpha+1,\beta}(\B^{r+2m},\B^q;\C^m)$ and $\nabla_{\tau,w} H''\in\Co^{\alpha,\beta-}(\B^{r+2m},\B^q;\C^{(r+2m)\times m})$;
        
        When $\frac12<\alpha<1$, $H''\in\Co^{\alpha+1,(2-\frac1\alpha)\beta}(\B^{r+2m},\B^q;\C^m)$ and $\nabla_{\tau,w} H''\in\Co^{\alpha,(2-\frac1\alpha)\beta}(\B^{r+2m},\B^q;\C^{(r+2m)\times m})$.
        \item\label{Item::EllipticPara::EstofH::HisDiffeo} Let $H(\tau,w,s):=(\tau,H''(\tau,w,s),s)$. Then $H:\B^{r+2m}\times\B^q\to\B^{r+2m}\times\B^q$ is homeomorphism, and $H(\frac14\B^{r+2m}\times\frac14\B^q)\subseteq\frac12\B^{r+2m}\times\frac12\B^q$.
        \item\label{Item::EllipticPara::EstofH::Phi} Endow $\R^r$ and $\C^m$ with another (real or complex) coordinate system $\sigma=(\sigma^1,\dots,\sigma^r)$ and $\zeta=(\zeta^1,\dots,\zeta^m)$. Let $\tilde\Phi:=H^\Inv:\B^{r+2m}_{\sigma,\zeta}\times\B^q_s\to\B^{r+2m}_{\tau,w}\times\B^q_s$ be the inverse map. Then the matrix map $\Lambda[A;H'']$ in \eqref{Eqn::EllipticPara::matrixfun} we have $\Lambda[A;H'']\circ\tilde\Phi\in\Co^\alpha_{\sigma,\zeta} L^\infty_s\cap\Co^{-1}_{\sigma,\zeta}\Co^{\min(\alpha+1,\beta,(2-\frac1\alpha)\beta)}_s(\B^{r+2m},\B^q;\C^{(r+m)\times m})$ with
        \begin{equation}\label{Eqn::EllipticPara::EstofH::BddLambda}
	        \|\Lambda[A;H'']\circ\tilde\Phi\|_{\Co^\alpha L^\infty \cap\Co^{-1} \Co^{\min(\alpha+1,\beta,(2-1/\alpha)\beta)}(\B^{r+2m},\B^q;\C^{(r+m)\times m})}<\eps.
	    \end{equation}
    \end{enumerate}
\end{cor}
\begin{proof}
We set for $k=1,0,-1$,
\begin{equation}\label{Eqn::EllipticPara::EstofH::DefXs}
    \begin{aligned}
        \Xs^k&:=\Co^{\alpha+k}_{\tau,w}L^\infty_s\cap\Co^k_{\tau,w}\Co^\beta_s,&&\text{when }\alpha>1
        \\
        \Xs^k&:=\Co^{\alpha+k}_{\tau,w}L^\infty_s\cap\Co^k_{\tau,w}\Co^\beta_s+\Co^{\alpha+k}_{\tau,w}\Co^{(2-1/\alpha)\beta}_s,&&\text{when }\tfrac12<\alpha<1.
    \end{aligned}
\end{equation}

Clearly \ref{Item::EllipticPara::XsAssume::C0} and \ref{Item::EllipticPara::XsAssume::Grad} hold. By Lemma \ref{Lem::Hold::LapInvBdd} along with a scaling on the $w$-variable we get \ref{Item::EllipticPara::XsAssume::Delta}. By Corollary \ref{Cor::Hold::CorMult} \ref{Item::Hold::CorMult::Prin} \ref{Item::Hold::CorMult::0} and Proposition \ref{Prop::Hold::MultLow} we get \ref{Item::EllipticPara::XsAssume::Mult}. Thus 
$\Xs^1,\Xs^0,\Xs^{-1}$ satisfy the assumptions in Proposition \ref{Prop::EllipticPara::ExistPDE}.

And since $\Co^{\alpha,\beta}\subset\Xs^0$ is an embedding. Therefore by Proposition \ref{Prop::EllipticPara::ExistPDE} for a suitable small $\delta_0>0$ and for $A$ such that $\|A\|_{\Co^{\alpha,\beta}}<\delta_0$, there is a $H''\in\Xs^1(\B^{r+2m},\B^q;\C^m)$ that solves \eqref{Eqn::EllipticPara::ExistenceH} with $H|_{(\partial\B^{r+2m})\times\B^q}=\Ic$.

From \eqref{Eqn::EllipticPara::EstofH::DefXs} and \eqref{Eqn::EllipticPara::Betas} we see that $H''\in\Co^{\alpha+1,\beta^{\sim\alpha}}(\B^{r+2m},\B^q;\C^m)$. By Remark \ref{Rmk::Hold::BiHoldInterpo} we have $\Co^\alpha L^\infty\cap\Co^0\Co^\beta\subset\Co^\alpha L^\infty\bigcap_{\theta>0}\Co^{\theta\alpha}\Co^{(1-\theta)\beta}\subset\Co^{\alpha,\beta-}$, thus $\nabla_{\tau,w}H''\in\Co^{\alpha,\beta-}$ when $\alpha>1$ and $\nabla_{\tau,w}H''\in\Co^{\alpha,\beta-}+\Co^\alpha\Co^{(2-1/\alpha)\beta}\subset\Co^{\alpha,(2-1/\alpha)\beta}$ when $\frac12<\alpha<1$. Thus $\nabla_{\tau,w}H''\in\Co^{\alpha,\beta^{\wedge\alpha}}(\B^{r+2m},\B^q;\C^{(r+2m)\times m})$ and we obtain \ref{Item::EllipticPara::EstofH::HReg}.

By taking $\eps<K_5(r+2m,q,\alpha,\beta^{\sim\alpha})^{-1}$ where $K_5$ is the constant in Proposition \ref{Prop::Hold::QPIFT}, we can find a $\delta_1\in(0,\delta_0)$ such that when $\|A\|_{\Co^{\alpha,\beta}}<\delta_1$ we have that $H:\B^{r+2m}\times\B^q\to\B^{r+2m}\times\B^q$ is a bijective $\Co^{\alpha,\beta^{\sim\alpha}}$-map and $\|\tilde \Phi\|_{\Co^{\alpha+1,\beta^{\sim\alpha}}(\B^{r+2m},\B^q;\R^{r+2m+q})}\lesssim1$. In particular $H$ is a homeomorphism. The result $H(\frac14\B^{r+2m}\times\frac14\B^q)\subseteq\frac12\B^{r+2m}\times\frac12\B^q$ can be done whenever $\|H-\id\|_{C^0(\B^{r+2m}\times\B^q)}=\|H''-\Ic\|_{C^0(\B^{r+2m}\times\B^q)}<\frac14$, which is achieved when $\delta_1$ is suitably small.  This proves \ref{Item::EllipticPara::EstofH::HisDiffeo}.

By Lemma \ref{Lem::EllipticPara::MatisInv} \ref{Item::EllipticPara::MatisInv::BddLambda} we have $\|\Lambda[A;H'']\|_{\Xs^0}\lesssim\|H''-\Ic\|_{\Xs^1}+\|A\|_{\Xs^0}\lesssim\|H''-\Ic\|_{\Xs^1}+\|A\|_{\Co^{\alpha,\beta}}$. So $\alpha>1$, \ref{Item::EllipticPara::EstofH::Phi} then follows from Proposition \ref{Prop::Hold::CompThm}.

When $\frac12<\alpha<1$ we have $\Xs^0\subset\Co^{\alpha,\beta^{\sim\alpha}}$, thus by Remark \ref{Rmk::Hold::BiHoldInterpo} and Lemma \ref{Lem::Hold::QPComp}, $$\|\Lambda\circ\tilde\Phi\|_{\Co^\alpha L^\infty\cap\Co^{-1}\Co^{\beta^{\sim\alpha}}}\lesssim\|\Lambda\circ\tilde\Phi\|_{\Co^{\alpha,\beta^{\sim\alpha}}}\lesssim\|\Lambda\|_{\Co^{\alpha,\beta^{\sim\alpha}}}\lesssim\|\Lambda\|_{\Xs^0},$$ where the implied constant depends on $\|\tilde \Phi\|_{\Co^{\alpha+1,\beta^{\sim\alpha}}(\B^{r+2m},\B^q;\R^{r+2m+q})}\lesssim1$. Therefore by keeping track of the constant we get \ref{Item::EllipticPara::EstofH::Phi} for $\frac12<\alpha<1$, finishing the proof.
\end{proof}


\subsection{A scaling argument on parameters}\label{Section::EllipticPara::Scaling}

In this part we improve Corollary \ref{Cor::EllipticPara::EstofH} to the case where $A\in\Co^{\alpha,\beta-}$. In particular we include the discussion of $\alpha=1$.

\begin{prop}\label{Prop::EllipticPara::ImprovedEstofH}
Let $\alpha\in(\frac12,1]$. Following the Convention \ref{Conv::EllipticPara::ConvofExtPDE}, for any $\eps>0$ there is a $\delta=\delta(\eps,\alpha,r,m,q)>0$ such that: 

For every $\beta\in\R_+$ such that $(2-\frac1\alpha)\beta\le\alpha+1$, and for every $A\in\Co^{\alpha,\beta-}(\B^{r+2m},\B^q;\C^{(r+m)\times m})$ that satisfies $\|A\|_{\Co^\alpha L^\infty(\B^{r+2m},\B^q)}<\delta$, there is a $H''$ that solves \eqref{Eqn::EllipticPara::ExistenceH} with $H|_{(\partial\B^{r+2m})\times\B^q}=\Ic$. Moreover
    \begin{enumerate}[parsep=-0.3ex,label=(\roman*)]
        \item\label{Item::EllipticPara::ImprovedEstofH::HReg}$H''\in\Co^{\alpha+1,(2-\frac1\alpha)\beta-}_\loc(\B^{r+2m},\B^q;\C^m)$ and $\nabla_{\tau,w} H''\in\Co^{\alpha,(2-\frac1\alpha)\beta-}_\loc(\B^{r+2m},\B^q;\C^{(r+2m)\times m})$.
        \item\label{Item::EllipticPara::ImprovedEstofH::HisDiffeo} Let $H(\tau,w,s):=(\tau,H''(\tau,w,s),s)$. Then $H:\B^{r+2m}\times\B^q\to\B^{r+2m}\times\B^q$ is homeomorphism, and $H(\frac14\B^{r+2m}\times\frac14\B^q)\subseteq\frac12\B^{r+2m}\times\frac12\B^q$.
        \item\label{Item::EllipticPara::ImprovedEstofH::Phi} Endow $\R^r$ and $\C^m$ with another (real or complex) coordinate system $\sigma=(\sigma^1,\dots,\sigma^r)$ and $\zeta=(\zeta^1,\dots,\zeta^m)$. Let $\tilde\Phi:=H^\Inv:\B^{r+2m}_{\sigma,\zeta}\times\B^q_s\to\B^{r+2m}_{\tau,w}\times\B^q_s$ be the inverse map. Then the matrix map $\Lambda[A;H'']$ in \eqref{Eqn::EllipticPara::matrixfun} we have $\Lambda[A;H'']\circ\tilde\Phi\in\Co^{\alpha,(2-\frac1\alpha)\beta-}_\loc(\B^{r+2m},\B^q;\C^{(r+m)\times m})$ with \begin{equation}\label{Eqn::EllipticPara::ImprovedEstofH::BddLambda}
            \|\Lambda[A;H'']\circ\tilde\Phi\|_{\Co^\alpha L^\infty}<\eps
        \end{equation}
    \end{enumerate}
\end{prop}

\begin{remark}
\begin{enumerate}[parsep=-0.3ex,label=(\roman*)]
    \item The results are also true for $\alpha>1$, where we replace $(2-\frac1\alpha)\beta- $ by $\beta-$. We do not state the result for $\alpha>1$ because in application to Theorem \ref{MainThm::RoughFro1} and \ref{MainThm::RoughFro2} we only use the case $\alpha\le1$.
    \item If we start with a $A\in\Co^{\alpha,\beta}$ from the technics in the proof we only obtain $H''\in\Co^{\alpha,(2-\frac1\alpha)\beta-}$. Comparing to Corollary \ref{Cor::EllipticPara::EstofH}, where we can get $H''\in\Co^{\alpha,(2-\frac1\alpha)\beta}$, there is a slight regularity loss in the parameter. By doing a more careful Schauder's estimate in the proof, it is possible to get the estimate $H''\in\Co^{\alpha,(2-\frac1\alpha)\beta}$ back.
\end{enumerate}
\end{remark}

\begin{proof}[Proof of Proposition \ref{Prop::EllipticPara::ImprovedEstofH}]
Let $\tilde C=\tilde C(r,m,\alpha)>0$ be maximum of the operator norms of the inclusion maps $\Co^\alpha(\B^{r+2m};\C^{(r+m)\times m})\hookrightarrow\Co^{\frac{2\alpha+1}4}(\B^{r+2m};\C^{(r+m)\times m})$ and $\Co^{\alpha+1}(\B^{r+2m};\C^m)\hookrightarrow\Co^{\frac{2\alpha+1}4+1}(\B^{r+2m};\C^m)$. 

We still use the notation $T[A;H'']$ for the left hand side of \eqref{Eqn::EllipticPara::ExistenceH}, following from Lemma \ref{Lem::EllipticPara::MatisInv} \ref{Item::EllipticPara::MatisInv::Diff}.

By taking $\Xs^j\in\{\Co^{\alpha+j}L^\infty,\Co^{\frac{2\alpha+1}4+j} L^\infty\}$ for $j=1,0,-1$ in Proposition \ref{Prop::EllipticPara::ExistPDE}, there is a $\delta_0=\delta_0(r,m,q,\alpha)>0$ and $C_0=C_0(r,m,q,\alpha)>0$ such that
\begin{enumerate}[parsep=-0.3ex,label=(T.\alph*)]
    \item\label{Item::EllipticPara::ImprovedEstofH:Pf::AlphaInfty} For every $A\in\Co^\alpha L^\infty(\B^{r+2m},\B^q;\C^{(r+m)\times m})$ such that $\|A\|_{\Co^\alpha L^\infty}<\delta_0$ there is a unique $H''\in\Co^{\alpha+1} L^\infty(\B^{r+2m},\B^q;\C^m)$ such that $\|H''-\Ic\|_{\Co^{\alpha+1} L^\infty}<C_0\delta_0$, $H''|_{(\partial\B^{r+2m})\times\B^q}=\Ic$ and $T[A;H]=0$. Moreover $\|H''-\Ic\|_{\Co^{\alpha+1}L^\infty}\le C_0\|A\|_{\Co^\alpha L^\infty}$.
    \item\label{Item::EllipticPara::ImprovedEstofH:Pf::Alpha} In particular for every $A_0\in\Co^\alpha(\B^{r+2m};\C^{(r+m)\times m})$ such that $\|A_0\|_{\Co^\alpha}<\delta_0$ there is a unique $H''_0\in\Co^{\alpha+1}(\B^{r+2m};\C^m)$ such that $\|H''_0-\Ic\|_{\Co^{\alpha+1}}<C_0\delta_0$, $H''_0|_{\partial\B^{r+2m}}=\Ic$ and $T[A_0;H_0]=0$.
    \item\label{Item::EllipticPara::ImprovedEstofH:Pf::SmallInfty} For every $A\in\Co^\frac{2\alpha+1}4 L^\infty(\B^{r+2m},\B^q;\C^{(r+m)\times m})$ such that $\|A\|_{\Co^\frac{2\alpha+1}4 L^\infty}<\tilde C\delta_0$ there is a unique $H''\in\Co^\frac{2\alpha+5}4L^\infty(\B^{r+2m},\B^q;\C^m)$ such that $\|H''-\Ic\|_{\Co^\frac{2\alpha+5}4L^\infty}<C_0\tilde C\delta_0$, $H''|_{(\partial\B^{r+2m})\times\B^q}=\Ic$ and $T[A;H]=0$. Moreover $\|H''-\Ic\|_{\Co^\frac{2\alpha+5}4L^\infty}\le C_0\|A\|_{\Co^\frac{2\alpha+1}4 L^\infty}$.
    \item\label{Item::EllipticPara::ImprovedEstofH:Pf::Small}  In particular for every $A_0\in\Co^\frac{2\alpha+1}4(\B^{r+2m};\C^{(r+m)\times m})$ such that $\|A_0\|_{\Co^\frac{2\alpha+1}4}<\tilde C\delta_0$ there is a unique $H''_0\in\Co^{\frac{2\alpha+5}4}(\B^{r+2m};\C^m)$ such that $\|H''_0-\Ic\|_{\Co^{\frac{2\alpha+1}4}}<C_0\tilde C\delta_0$, $H''_0|_{\partial\B^{r+2m}}=\Ic$ and $T[A_0;H_0]=0$.
\end{enumerate}
Note that by the assumption of $\tilde C$, $\|A_0\|_{\Co^\alpha}<\delta_0$ implies $\|A_0\|_{\Co^\frac{2\alpha+1}4}<\tilde C\delta_0$.


By shrinking $\delta_0$, we can assume that $C_0\delta_0<K_2(n,\alpha,\alpha)^{-1}$ where $K_2$ is the constant in Proposition \ref{Prop::Hold::QIFT}. Therefore by Proposition \ref{Prop::Hold::QIFT} \ref{Item::Hold::QIFT::PhiEst},
\begin{enumerate}[parsep=-0.3ex,label=(T.\alph*)]\setcounter{enumi}{4}
    \item \label{Item::EllipticPara::ImprovedEstofH:Pf::HInv}For $H''$ from the consequence of \ref{Item::EllipticPara::ImprovedEstofH:Pf::AlphaInfty} (since $\|H''\|_{\Co^{\alpha+1}L^\infty}<K_2^{-1}$), $H''(\cdot,s)$ has inverse map $\Phi(\cdot,s)$ with $\|\Phi(\cdot,s)\|_{\Co^{\alpha+1}}<K_2$, for each $s\in\B^q$. 
\end{enumerate}

Moreover by Proposition \ref{Prop::Hold::QComp} \ref{Item::Hold::QComp::>1},
\begin{equation}\label{Eqn::EllipticPara::ImprovedEstofH:Pf::BddLambda}
    \|\Lambda[A;H'']\circ\tilde\Phi\|_{\Co^\alpha L^\infty}=\sup_{s\in\B^q}\|\Lambda[A;H''](\cdot,s)\circ\Phi(\cdot,s)\|_{\Co^\alpha}\le K_1(m,r+2m,\alpha+1,\alpha,K_2)\|\Lambda[A;H'']\|_{\Co^\alpha L^\infty},
\end{equation}
where $K_1$ is the constant in Proposition \ref{Prop::Hold::QComp}.

On the other hand by Lemma \ref{Lem::EllipticPara::MatisInv} \ref{Item::EllipticPara::MatisInv::BddLambda}, $\|\Lambda[A;H'']\|_{\Co^\alpha L^\infty}\lesssim_{r,m,q,\alpha,\delta_0}\|H''-\Ic\|_{\Co^{\alpha+1}L^\infty}+\|A\|_{\Co^\alpha L^\infty}$. Therefore there exist a $C_1=C_1(r,m,q,\alpha,\delta_0)>0$ such that
\begin{equation}\label{Eqn::EllipticPara::ImprovedEstofH:Pf::LambdaBdd}
    \|\Lambda[A;H'']\circ\tilde\Phi\|_{\Co^\alpha L^\infty}\le C_1\|A\|_{\Co^\alpha L^\infty},\quad\forall A\in\Co^{\alpha}L^\infty(\B^{r+2m},\B^q;\C^{(r+m)\times m})\text{ with }\|A\|_{\Co^\alpha L^\infty}<\delta_0.
\end{equation}
This gives \eqref{Eqn::EllipticPara::ImprovedEstofH::BddLambda} once $\|A\|_{\Co^\alpha L^\infty}$ is small enough, depending on the $\eps>0$ in the assumption.

\medskip 
By direct computation, for each $A$ and $H''$, the tangent map $\frac{\partial T}{\partial H''}[A;H'']$ is the second differential operator that has the following form: for a function $(F^1,\dots,F^m):\B^{r+2m}\times\B^q\to\C^m$ and $i=1,\dots,m$
\begin{equation*}
    \lim_{h\to0}\frac{T[A;H''+hF]^i-T[A;H'']^i}h=(\Delta_\tau+\Box_w)F^i+\sum_{j,k=1}^{r+2m}\sum_{u,v=1}^m\big(B^{1,iu}_{jk}\partial_j(B^{2,uv}_{jk}\partial_kF^v)+B^{3,iu}_{jk}(\partial_jB^{4,uv}_{jk})\partial_kF^u\big).
\end{equation*}
Here $(\partial_1,\dots,\partial_{r+2m})=(\partial_{\tau^1},\dots,\partial_{\tau^1},\partial_{\re w^1},\partial_{\im w^1}\dots,\partial_{\re w^m},\partial_{\im w^m})$. Here $B^{\sigma,uv}_{jk}[A;H'']$, $1\le \sigma\le 4$, $1\le u,v\le m$, $1\le j,k\le r+2m$ are all rational functions of the components of $A$, $\nabla_{\tau,w}H''$ and $\nabla_{\tau,w}\bar H''$, such that $B^{\sigma,iu}_{jk}(A,\nabla_{\tau,w} H'',\nabla_{\tau,w}\bar H'')|_{A=0;\nabla_{\tau,w} H''=\nabla_{\tau,w}\bar H''=\nabla_{\tau,w}\Ic}=0$.

Taking the first order Taylor expansion of $B^{\sigma,uv}_{jk}(A,\nabla_{\tau,w} H'',\nabla_{\tau,w}\bar H'')$ we see that there is a $C_2=C_2(r,m,\alpha,\delta_0)>0$ such that for every $\|A\|_{\Co^\alpha L^\infty}<\delta_0$ and $\|H''-\Ic\|_{\Co^{\alpha+1}L^\infty}<\tilde C\delta_0$,
\begin{equation}\label{Eqn::EllipticPara::ImprovedEstofH:Pf::BddB}
    \sum_{\sigma=1}^4\sum_{u,v=1}^m\sum_{j,k=1}^{r+2m}\|B^{\sigma,uv}_{jk}(A,\nabla_{\tau,w} H'',\nabla_{\tau,w}\bar H'')\|_{\Co^\alpha L^\infty}\le C_2(\|A\|_{\Co^\alpha L^\infty}+\|H''-\Ic\|_{\Co^{\alpha+1} L^\infty}).
\end{equation}

By Proposition \ref{Prop::Hold::StrongLapInv} with a scaling so that $\Delta_\tau+\Box_w$ becomes a real Laplacian, there is a $\delta_1=\delta_1(r,m,\alpha,\delta_1)\in(0,\delta_0)$ such that
\begin{enumerate}[parsep=-0.3ex,label=(T.\alph*)]\setcounter{enumi}{5}
    \item\label{Item::EllipticPara::ImprovedEstofH:Pf::TangTInv} If $\|A_0\|_{\Co^\alpha}<\delta_1$ and $\|H''_0-\Ic\|_{\Co^{\alpha+1}}<C_0\delta_1$, then the following tangent map is invertible:
    \begin{equation}\label{Eqn::EllipticPara::ImprovedEstofH:Pf::TangTInv}
        \frac{\partial T}{\partial H''}[A_0;H''_0]:\{F\in\Co^{\gamma+1}(\B^{r+2m};\C^m):F|_{\partial\B^{r+2m}}=0\}\to\Co^{\gamma-1}(\B^{r+2m};\C^m),\ \forall\gamma\in(1-\alpha,\alpha].
    \end{equation}
\end{enumerate}
This can be done by taking $\delta_1$ such that $C_2(\delta_1+\tilde C\delta_1)\lesssim_mK_3^{-1}$ where $K_3=K_3(\{|\tau|^2+\frac14|w|^2<1\},m,\alpha)>0$ is the constant in Proposition \ref{Prop::Hold::StrongLapInv}.

\medskip
Let $0<\eps<\frac{2\alpha-1}4$, we consider the spaces
\begin{equation*}
    \Xs_\eps^j:=\Co^{\alpha-\eps+j}\Co^\frac{\eps\beta}{2\alpha}\cap\Co^{1-\alpha+2\eps+j}\Co^{(2-\frac{1+3\eps}\alpha)\beta},\quad j=1,0,-1.
\end{equation*}

We have $\Xs_\eps^0\subset\Co^\frac{2\alpha+1}4 L^\infty$. And by Remark \ref{Rmk::Hold::BiHoldInterpo} we see that $\Xs_\eps^0\subset\Co^{\alpha,\frac{2\alpha-1-2\eps}{2\alpha-1-3\eps}\beta}\subset\Co^{\alpha,\beta-}$, for every $0<\eps<\frac{2\alpha-1}4$.

We have the product map $\Co^{\alpha-\eps}\Co^\frac{\eps\beta}{2\alpha}\times \Co^{\alpha-\eps-1}\Co^\frac{\eps\beta}{2\alpha}\to\Co^{\alpha-\eps-1}\Co^\frac{\eps\beta}{2\alpha}$. By Proposition \ref{Prop::Hold::Mult}, we have product map $\Co^{\alpha-\eps}\Co^\frac{\eps\beta}{2\alpha}\cap\Co^{1-\alpha+2\eps}\Co^{(2-\frac{1+3\eps}\alpha)\beta}\times \Co^{\alpha-\eps-1}\Co^\frac{\eps\beta}{2\alpha}\cap\Co^{2\eps-\alpha}\Co^{(2-\frac{1+3\eps}\alpha)\beta}\to \Co^{\alpha-\eps-1}\Co^\frac{\eps\beta}{2\alpha}\cap\Co^{2\eps-\alpha}\Co^{(2-\frac{1+3\eps}\alpha)\beta}$. Taking intersections we have $\Xs_\eps^0\times\Xs_\eps^{-1}\to\Xs_\eps^{-1}$. Therefore $(\Xs_\eps^1,\Xs_\eps^0,\Xs_\eps^{-1})$ satisfy conditions \ref{Item::EllipticPara::XsAssume::C0} - \ref{Item::EllipticPara::XsAssume::Delta} in Proposition \ref{Prop::EllipticPara::ExistPDE}.

By taking $\gamma\in\{\alpha-\eps,1-\alpha+2\eps\}$ in \eqref{Eqn::EllipticPara::ImprovedEstofH:Pf::TangTInv} and applying Lemma \ref{Lem::Hold::OperatorExtension}, we see that
\begin{enumerate}[parsep=-0.3ex,label=(T.\alph*)]\setcounter{enumi}{6}
    \item\label{Item::EllipticPara::ImprovedEstofH:Pf::XsInv}Let $A_0\in\Co^\alpha(\B^{r+2m};\C^{(r+m)\times m})$ and $H''_0\in \Co^{\alpha+1}(\B^{r+2m};\C^m)$ such that $\|A_0\|_{\Co^\alpha}<\delta_1$, $H''_0|_{\partial\B^{r+2m}}=\Ic$, $\|H''_0-\Ic\|_{\Co^{\alpha+1}}<C_0\delta_1$ and $T[A_0;H_0]=0$. Then the tangent map 
\begin{equation*}
    \frac{\partial T}{\partial H''}[A_0;H''_0]:\{F\in\Xs^1(\B^{r+2m},\B^q;\C^m):F|_{(\partial\B^{r+2m})\times\B^q}=0\}\to\Xs^{-1}(\B^{r+2m},\B^q;\C^m)\text{ is invertible}.
\end{equation*}
\end{enumerate}

Therefore by the standard Implicit Function Theorem on Banach Spaces (see for example \cite[Page 417]{Edwards}), or by doing the contraction mapping argument similar to \eqref{Eqn::EllipticPara::ExistPDE::SeACont}, we see that
\begin{enumerate}[parsep=-0.3ex,label=(T.\alph*)]\setcounter{enumi}{7}
    \item\label{Item::EllipticPara::ImprovedEstofH:Pf::XsEps}Let $A_0$ and $H''_0$ be as in \ref{Item::EllipticPara::ImprovedEstofH:Pf::XsInv}. There is a $\delta'_{\eps,A_0,H''_0}=\delta'(r,m,q,\alpha,\eps,A_0,H''_0)>0$ and a continuous map\footnote{In fact $\mathfrak h$ does not depend on $\eps,A_0,H''_0$ but we do not need this property explicitly in the proof.} 
    $$\hspace{-0.4in}\mathfrak h^{A_0,H_0}_\eps:\{A\in\Xs_\eps^0(\B^{r+2m};\C^{(r+m)\times m}):\|A-A_0\|_{\Xs_\eps^0}<\delta'\}\to\{H''\in\Xs_\eps^1(\B^{r+2m},\B^q;\C^m):H''|_{(\partial\B^{r+2m})\times\B^q}=\Ic\},$$ such that $\mathfrak h^{A_0,H_0}_\eps[A_0]=H_0$ and $T[A,\mathfrak h^{A_0,H_0}_\eps[A]]=0$ for every $\|A-A_0\|_{\Xs_\eps^0}<\delta' $.

\end{enumerate}

Since $\Xs_\eps^j\subset \Co^{\frac{2\alpha+1}4+j} L^\infty$, we can shrink $\delta'$ in \ref{Item::EllipticPara::ImprovedEstofH:Pf::XsEps}, which still depends only on $r,m,q,\alpha,\eps,A_0,H''_0$, such that 
\begin{equation}\label{Eqn::EllipticPara::ImprovedEstofH:Pf::AssumpDelta1}
    \|A_0\|_{\Co^\alpha}<\delta_1\text{ and }\|A-A_0\|_{\Xs^0}<\delta'_{\eps,A_0,H''_0}\quad\Rightarrow\|A\|_{ \Co^\frac{2\alpha+1}4 L^\infty}<C_0\delta_0\text{ and }\|\mathfrak h^{A_0,H_0}_\eps[A]-\Ic\|_{\Co^\frac{2\alpha+5}4 L^\infty}\le C_0\tilde C\delta_0.
\end{equation}

We are now ready to prove \ref{Item::EllipticPara::ImprovedEstofH::HReg}, \ref{Item::EllipticPara::ImprovedEstofH::HisDiffeo} and \ref{Item::EllipticPara::ImprovedEstofH::Phi}.

\medskip Now let $A\in\Co^{\alpha,\beta-}(\B^{r+2m},\B^q;\C^{(r+m)\times m})$ be satisfied $\|A\|_{\Co^\alpha L^\infty}<\delta_1$, where $\delta_1$ is the constant in \ref{Item::EllipticPara::ImprovedEstofH:Pf::XsInv}. By \ref{Item::EllipticPara::ImprovedEstofH:Pf::AlphaInfty}, since $\delta_1\le\delta_0$, there is a unique $H''\in\Co^{\alpha+1}L^\infty(\B^{r+2m},\B^q;\C^m)$ that satisfies the consequences in \ref{Item::EllipticPara::ImprovedEstofH:Pf::AlphaInfty}. 

Since $\Xs_\eps^j\subset\Co^{\alpha+j,\frac{2\alpha-1-2\eps}{2\alpha-1-3\eps}\beta}$ and $\Co^{\alpha+j,\beta-}=\bigcap_{0<\eps<\frac{2\alpha-1}4}\Co^{\alpha+j,\frac{2\alpha-1-2\eps}{2\alpha-1-3\eps}\beta}$ for $j=0,1$, to prove \ref{Item::EllipticPara::ImprovedEstofH::HReg} it suffices to that, for every $s_0\in\B^q$ and $\eps>0$, there is a $0<\mu<1-|s_0|$ such that $H''\in\Xs_\eps^1(\B^{r+2m},B^q(s_0,\mu);\C^m)$.

Set $A_0:=A(\cdot,s_0)$ and $H''_0:=H''(\cdot,s_0)$. We see that $A_0$ and $H''_0$ satisfy the assumptions of \ref{Item::EllipticPara::ImprovedEstofH:Pf::XsInv}. 

For $\mu\in(0,1-|s_0|)$, we let $A^{s_0,\mu}(\tau,w,s):=A(\tau,w,\mu(s-s_0))$. By Lemma \ref{Lem::Hold::ScalingLem}, we see that
\begin{align*}
    \|A^{s_0,\mu}-A_0\|_{\Co^{\alpha-\eps}\Co^\frac{\eps\beta}{2\alpha}(\B^{r+2m},\B^q)}\lesssim&_{\eps,\alpha,\beta}\mu^{\min(\frac12,\frac{\eps\beta}{2\alpha})}\|A\|_{\Co^{\alpha-\eps}\Co^\frac{\eps\beta}{2\alpha}}\lesssim_{\eps,\alpha}\mu^{\min(\frac12,\frac{\eps\beta}{2\alpha})}\delta_1;
    \\
    \|A^{s_0,\mu}-A_0\|_{\Co^{1-\alpha+2\eps}\Co^{(2-\frac{1+3\eps}\alpha)\beta}(\B^{r+2m},\B^q)}\lesssim&_{\eps,\alpha,\beta}\mu^{1-\alpha+2\eps}\|A\|_{\Co^{1-\alpha+2\eps}\Co^{(2-\frac{1+3\eps}\alpha)\beta}}\lesssim_{\eps,\alpha}\mu^{1-\alpha+2\eps}\delta_1.
\end{align*}
In particular there is a small $\mu>0$, such that
\begin{equation*}
    \|A^{s_0,\mu}-A_0\|_{\Xs_\eps^0(\B^{r+2m},\B^q)}<\delta'_{\eps,A_0,H''_0},
\end{equation*}
where $\delta'_{\eps,A_0,H''_0}$ is the small constant in \ref{Item::EllipticPara::ImprovedEstofH:Pf::XsEps} and satisfies \eqref{Eqn::EllipticPara::ImprovedEstofH:Pf::AssumpDelta1}. Therefore by \ref{Item::EllipticPara::ImprovedEstofH:Pf::XsEps}, there is a map $\mathfrak h^{A_0,H_0}_\eps$ such that $T[A^{s_0,\mu};\mathfrak h^{A_0,H_0}_\eps[A^{s_0,\mu}]]=0$.

On the other hand we have $T[A;H'']=0$. By scaling on $s\in\B^q$ we see that $H''^{s_0,\mu}(\tau,w,s):=H''(\tau,w,\mu(s-s_0))$ satisfies $T[A^{s_0,\mu};H'']=0$. Applying the uniqueness of $H''$ in \ref{Item::EllipticPara::ImprovedEstofH:Pf::SmallInfty} and the assumption \eqref{Eqn::EllipticPara::ImprovedEstofH:Pf::AssumpDelta1}, we see that $\mathfrak h^{A_0,H_0}_\eps[A^{s_0,\mu}]\in H''^{s_0,\mu}$.

Therefore $H''^{s_0,\mu}\in\Xs_\eps^1(\B^{r+2m},\B^q;\C^m)$. Taking scaling back we see that $H''\in\Xs_\eps^1(\B^{r+2m},B^q(s_0,\mu);\C^m)$. Since $s_0\in\B^q$ is arbitrary, we get $H''\in\Xs_{\eps,\loc}^1(\B^{r+2m},\B^q;\C^m)$. Since $\eps>0$ is arbitrary, we get $H''\in\Co^{\alpha+1,(2-\frac1\alpha)\beta-}_\loc(\B^{r+2m},\B^q;\C^m)$ and $\nabla_{\tau,w}H''\in\Co^{\alpha,(2-\frac1\alpha)\beta-}_\loc(\B^{r+2m},\B^q;\C^{(r+2m)\times m})$. Therefore for every $\delta<\delta_1$, we have the result \ref{Item::EllipticPara::ImprovedEstofH::HReg}.

By \ref{Item::EllipticPara::ImprovedEstofH:Pf::HInv}, we see that $H:\B^{r+2m}\times\B^q\to\B^{r+2m}\times\B^q$ is bijective, hence $\tilde\Phi$ is defined.
By Lemma \ref{Lem::Hold::CompofMixHold} \ref{Item::Hold::CompofMixHold::InvFun}, since $(2-\frac1\alpha)\beta\le\alpha+1$, we see that $\tilde\Phi$ is locally $\Co^{\alpha+1,(2-\frac1\alpha)\beta-\eps}$, for every $\eps$ near every $(\tau,w,s)\in\B^{r+2m}\times\B^q$. Therefore $\tilde \Phi\in\Co^{\alpha+1,(2-\frac1\alpha)\beta-}_\loc(\B^{r+2m},\B^q;\B^{r+2m}\times\B^q)$, in particular $H$ is homeomorphism.

We see that $H(\frac14\B^{r+2m}\times\frac14\B^q)\subseteq\frac12\B^{r+2m}\times\frac12\B^q$ whenever $\|H''-\Ic\|_{\Co^{\delta+1}L^\infty}<\frac14$. Thus by shrinking $\delta_1$ such that $\tilde C\delta_1<\frac14$, we see that when $\|A\|_{\Co^\alpha L^\infty}<\delta_1$ the result \ref{Item::EllipticPara::ImprovedEstofH::HisDiffeo} is satisfied.

Finally for every $\eps>0$, by taking $\delta=\min(\delta_1,\eps/C_1)$ where $C_1$ is the constant in \eqref{Eqn::EllipticPara::ImprovedEstofH:Pf::BddLambda}, we get \eqref{Eqn::EllipticPara::ImprovedEstofH::BddLambda} from \eqref{Eqn::EllipticPara::ImprovedEstofH:Pf::BddLambda} and we complete the proof.
\end{proof}

\begin{remark}[An alternative method via Generalized H\"older-Zygmund spaces]
Here the proof of Proposition \ref{Prop::EllipticPara::ImprovedEstofH} is a kind of Schauder's estimate. Alternatively we can have a different proof use generalized H\"older-Zygmund spaces.

 Note that this method does not require the scaling argument Lemma \ref{Lem::Hold::ScalingLem}, and we can also discuss the case $A\in\Co^{\alpha-,\beta-}$.

Let $(\phi_j)_{j=0}^\infty$ be a dyadic resolutions for $\R^n$. In Lemma \ref{Lem::Hold::HoldChar} we see that $f\in\Co^\beta(\R^n)$ if and only if $\|\phi_j\ast f\|_{L^\infty}\lesssim 2^{-j\beta}$. We can generalize the decay/growth rate $2^{-j\beta}$ to $\mu(2^{-j})$ where $\mu:(0,1]\to\R_+$ is a certain modulus of smoothness. We define $\Co^\mu(\R^n)$ to be the space of tempered distributions $f$ such that $\|\phi_j\ast f\|_{L^\infty}\lesssim \mu(2^{-j})$. For example we can pick $\mu(2^{-j})=2^{-j\beta+\sqrt j}$, where we have $\Co^\beta\subsetneq\Co^\mu\subsetneq\Co^{\beta-}$.  The $\Co^\mu$-spaces can naturally be extended to the bi-parameter case.

We assume $A\in\Co^{\alpha,\beta-}$ and we wish to find a $H''\in\Co^{\alpha+1,\beta-}$ that solves \eqref{Eqn::EllipticPara::ExistenceH}. 

By assumption we can find modulus of smoothness $\mu$ such that $A\in\Co^{\alpha,\mu}$ where $\mu:(0,1]\to\R_+$ satisfies $\lim_{t\to0}\frac{\log\mu(t)}{\log t}=\beta$, along with some other good properties. 

Then we look for spaces $\Xs^1,\Xs^0,\Xs^{-1}$ that satisfy \ref{Item::EllipticPara::XsAssume::C0} - \ref{Item::EllipticPara::XsAssume::Mult} in Proposition \ref{Prop::EllipticPara::ExistPDE}, such that $\Co^{\alpha,\mu}\subseteq\Xs^0\subseteq\Co^{\alpha,(\beta^{\sim\alpha})-}$ and $\Co^{\alpha+1,\mu}\subseteq\Xs^1\subseteq\Co^{\alpha+1,(\beta^{\sim\alpha})-}$. To achieve this, especially to check the condition \ref{Item::EllipticPara::XsAssume::Mult}, we need to run again the arguments in Section \ref{Section::BiHoldMult}. For example one needs to prove the paraproduct operators $$\Pf,\Rf:\Co^\mu\times L^\infty\to\Co^\mu,$$ are bounded linear.

Once this is done we know $H\in\Co^{\alpha+1,\mu}\subseteq\Co^{\alpha+1,(\beta^{\sim\alpha})-}$ and $\nabla_{\tau,w}H\in\Co^{\alpha+1,\mu}\subseteq \Co^{\alpha,(\beta^{\sim\alpha})-}$. This completes the sketch of the proof.


For more details of generalized H\"older-Zygmund spaces, we refer \cite{KreitThesisGeneralizedHolder} to readers.
\end{remark}

\subsection{The regularity proposition}
\label{Section::EllipticPara::AnalPDE}
In this part on a complex space we denote the fixed open cone
\begin{equation}\label{Eqn::HCone}
    \Hb^n=\{x+iy:x,y\in\B^n,\ 4|y|<1-|x|\}\subset\C^n_z.
\end{equation}
Clearly $\Hb^n_z\cap\R^n_x=\B^n$ is a ``complex extension'' of the unit ball $\B^n\subset\R^n$ in the real domain. 

For a holomorphic function $f\in\Oh(\Hb^n)$, we denote $\partial f=\partial_zf=(\frac{\partial f}{\partial z^1},\dots,\frac{\partial f}{\partial z^n})$ as a $n$-dimensional vector valued function.

In this part we give the proof of the following regularity result for a quasilinear elliptic equation system:
\begin{prop}\label{Prop::EllipticPara::AnalyticPDE}Let $m,n,p,q\in\Z_+$, let $\Theta:\C^m\times\C^{n\times m}\to\C^p$ be a complex bilinear map, and let $L=\sum_{j=1}^na^j\Coorvec{x^j}$ be a vector-valued first order differential operator on $\R^n$ with constant matrix coefficients $a^j\in\C^{m\times p}$.

Assume $\alpha>\frac12$, there is a $\eps_0=\eps_0(m,n,p,q,\alpha,\Theta,L)>0$ that satisfies the following:

Let $\tilde\beta\in\{\gamma,\gamma-:\gamma\in\R_+\}$. Suppose $f\in \Co^\alpha_xL^\infty_s\cap\Co^{-1}_x\Co^{\tilde\beta}_s(\B^n,\B^q;\C^m)$ satisfies $\|f\|_{\Co^\alpha_xL^\infty_s}<\eps_0$ and $\Delta_xf=L_x\Theta(f,\partial_x f)$, then $f$ admits extension $\mathfrak f:\Hb^n\times \B^q\subset\C^n_z\times \R^q_s\to\C^m$ such that $\mathfrak f$ is holomorphic in $z$ and $\mathfrak f(x+i0;s)=f(x;s)$ for $x\in\B^n,s\in \B^q$. 
Moreover $\mathfrak f\in\Co^\infty_\loc\Co^{\tilde\beta}_\loc(\Hb^n,\B^q;\C^m)$.
\end{prop}

In application we take $\tilde\beta=\beta^{\sim\alpha}=\min(\beta,(2-1/\alpha)\beta,\alpha+1)$. We note that $\eps_0$ can be choose independent of $\tilde\beta$.

The proof of Proposition \ref{Prop::EllipticPara::AnalyticPDE} is based on Proposition \ref{Prop::HolLap}.  Here  we need the parameter version of the spaces $\Co^\alpha_\Oh(\Hb^n)$, which are the natural generalizations to Definition \ref{Defn::SecHolLap::HoloHoldSpace}:


\begin{defn}\label{Defn::EllipticPara::ParaHoloHoldSpace}
Let $\alpha<1$.
We define $\Co^\alpha_{\Oh}L^\infty(\Hb^n,\B^q)=\Co^\alpha_{\Oh}L^\infty(\Hb^n,\B^q;\C)$ as the set of $f\in L^\infty (\Hb^n,\B^q;\C)$ such that $f(\cdot,s)\in \Co^\alpha_\Oh(\Hb^n)$ for almost every $s\in\B^q$ and $\|f\|_{\Co^\alpha_{\Oh}L^\infty}:=\essup_{s\in\B^q}\|f(\cdot,s)\|_{\Co^\alpha_\Oh}<\infty$.

For $\beta>0$ we define $\Co^\alpha_{\Oh}\Co^\beta(\Hb^n,\B^q):=\Co^\beta(\B^q;\Co^\alpha_{\Oh}(\Hb^n))$. 

We define the space $\Co^\alpha_{\Oh}L^\infty\cap\Co^{-\eta}_{\Oh}\Co^\beta(\Hb^n,\B^q)$ with norm $\|f\|_{\Co^\alpha_{\Oh}L^\infty\cap\Co^{-\eta}_{\Oh}\Co^\beta(\Hb^n,\B^q)}:=\|f\|_{\Co^\alpha_{\Oh}L^\infty(\Hb^n,\B^q)}+\|f\|_{\Co^{-\eta}_{\Oh}\Co^\beta(\Hb^n,\B^q)}$, for $\eta<1$.
\end{defn}

\begin{remark}\label{Rmk::EllipticPara::HolLapPara}
    By Lemmas \ref{Lem::Hold::OperatorExtension} and \ref{Lem::Hold::CharCalphaInfty1} \ref{Item::Hold::CharCalphaInfty1::OptExt}, we see that $\Pv,\tilde\Pv,\Ex$ in Proposition \ref{Prop::HolLap} have natural extension to the parameter case. We still denote them as $\Pv,\tilde\Pv,\Ex$: let $V\subseteq\R^q$ be a bounded smooth subset, we have boundedness:
    \begin{itemize}[nolistsep]
        \item $\Pv:\Co^{\alpha-2}\Xs(\B^n,V)\to\Co^\alpha\Xs(\B^n,V)$ for $-2<\alpha<1$ and $\Xs\in\{L^\infty,\Co^\beta:\beta>0\}$.
        \item $\tilde\Pv:\Co^{\alpha-2}_\Oh\Xs(\Hb^n,V)\to\Co^\alpha_\Oh\Xs(\Hb^n,V)$ for $-2<\alpha<1$ and $\Xs\in\{L^\infty,\Co^\beta:\beta>0\}$.
        \item $\Ex:\Co^{\alpha}\Xs(\B^n,V)\cap\ker\Delta_x\to\Co^\alpha_\Oh\Xs(\Hb^n,V)$ for $-2<\alpha<1$ and $\Xs\in\{L^\infty,\Co^\beta:\beta>0\}$.
    \end{itemize}
    Moreover, their operator norms equal to $\|\Pv\|_{\Co^{\alpha-2}(\B^n)\to\Co^\alpha(\B^n)}$, $\|\tilde\Pv\|_{\Co^{\alpha-2}_\Oh(\Hb^n)\to\Co^\alpha_\Oh(\Hb^n)}$ and $\|\Ex\|_{\Co^{\alpha}(\B^n)\cap\ker\Delta\to\Co^\alpha(\Hb^n)}$, respectively.
\end{remark}

The following lemma is useful in the proof of Proposition \ref{Prop::EllipticPara::AnalyticPDE}.
\begin{lem}\label{Lem::EllipticPara::HoloParaLem}
    Let $\beta>0$ and let $\Xs\in\{L^\infty,\Co^\beta\}$.
    \begin{enumerate}[parsep=-0.3ex,label=(\roman*)]
        \item\label{Eqn::EllipticPara::HoloParaLem::Grad} The differentiation $[\tilde f\mapsto\partial_z\tilde f]:\Co^{\alpha}_{\Oh}\Xs(\Hb^n,\B^q)\to\Co^{\alpha-1}_{\Oh}\Xs(\Hb^n,\B^q;\C^n)$ is bounded linear for all $\alpha<1$.
        \item\label{Item::EllipticPara::HoloParaLem::EqvNorm} For $\alpha<0$, $\Co^\alpha_{\Oh}\Xs(\Hb^n,\B^q)$ has an equivalent norm
        \begin{equation}\label{Eqn::EllipticPara::HoloParaLem::EqvNorm}
            \textstyle\tilde f\mapsto\sup_{z\in\Hb^n}\dist(z,\partial\Hb^n)^{-\alpha}\|\tilde f(z,\cdot)\|_{\Xs(\B^q;\C)}.
        \end{equation}
        
        \item\label{Item::EllipticPara::HoloParaLem::Res} The restriction map $[\tilde f\mapsto \tilde f|_{\B^n\times\B^q}]:\Co^\alpha_\Oh\Xs(\Hb^n,\B^q)\to\Co^\alpha_x\Xs(\B^n,\B^q;\C)$ is bounded linear for $0<\alpha<1$.
        \item\label{Item::EllipticPara::HoloParaLem::Inc} We have inclusion map $\Co^\alpha_\Oh\Xs(\Hb^n,\B^q)\hookrightarrow L^\infty\Xs(\Hb^n,\B^q;\C)$ for $0<\alpha<1$.
    \end{enumerate}
    Moreover all operator norms above are not depend on $\beta$.
\end{lem}
\begin{proof}These results follow from Lemma \ref{Lem::SecHolLap::HLLem}, the parameter free cases, along with  Lemma \ref{Lem::Hold::OperatorExtension} when $\Xs=\Co^\beta$, and the definition when $\Xs=L^\infty$. 
\end{proof}

We need to show that in $\Hb^n\times \B^q$, $\mathfrak f\mapsto \Theta(\mathfrak f,\partial_z\mathfrak f)$ maps $\mathfrak f$ into the desired function space.

\begin{lem}\label{Lem::EllipticPara::HoloMult}
	Let $0<\alpha<1$ and $\beta>0$. 
For $f\in\Co^\alpha_\Oh L^\infty\cap\Co^{-1}_\Oh\Co^\beta(\Hb^n,\B^q)$ and $g\in\Co^{\alpha-1}_\Oh L^\infty\cap\Co^{-2}_\Oh\Co^\beta(\Hb^n,\B^q)$, we have $fg\in\Co^{\alpha-1}_\Oh L^\infty\cap\Co^{-2}_\Oh\Co^\beta(\Hb^n,\B^q)$. Moreover there is a $C=C(n,q,\alpha)>0$ that does not depend on $\beta$ and $f,g$ such that
\begin{equation}\label{Eqn::EllipticPara::HoloMult}
    \|fg\|_{\Co^{\alpha-1}_\Oh L^\infty\cap\Co^{-2}_\Oh\Co^\beta(\Hb^n,\B^q)}\le C\|f\|_{\Co^{\alpha}_\Oh L^\infty\cap\Co^{-1}_\Oh\Co^\beta(\Hb^n,\B^q)}\|g\|_{\Co^{\alpha-1}_\Oh L^\infty\cap\Co^{-2}_\Oh\Co^\beta(\Hb^n,\B^q)}.
\end{equation}
\end{lem}

\begin{proof}Clearly the product $f(z,s)g(z,s)$ is locally continuous and holomorphic in $z$.
    
    By Lemma \ref{Lem::EllipticPara::HoloParaLem} \ref{Item::EllipticPara::HoloParaLem::EqvNorm} and \ref{Item::EllipticPara::HoloParaLem::Inc}, we have almost every $s\in\B^q$,
    \begin{align*}
        |(fg)(z,s)|\lesssim_{n,\alpha}\|f\|_{L^\infty}\dist(z,\partial\Hb^n)^{\alpha-1}\|g\|_{\Co^{\alpha-1}_\Oh L^\infty}\lesssim_{n,\alpha}\|f\|_{\Co^\alpha_\Oh L^\infty}\dist(z,\partial\Hb^n)^{\alpha-1}\|g\|_{\Co^{\alpha-1}_\Oh L^\infty}.
    \end{align*}
    Dividing $\dist(z,\partial\Hb^n)^{\alpha-1}$ on both side, taking essential sup over $z$ and $s$, and using Lemma \ref{Lem::EllipticPara::HoloParaLem} \ref{Item::EllipticPara::HoloParaLem::EqvNorm} we get the first control $\|fg\|_{\Co^{\alpha-1}_\Oh L^\infty}\lesssim\|f\|_{\Co^{\alpha}_\Oh L^\infty\cap\Co^{-1}_\Oh\Co^\beta}\|g\|_{\Co^{\alpha-1}_\Oh L^\infty\cap\Co^{-2}_\Oh\Co^\beta}$.
    
    To control $\|fg\|_{\Co^{-2}_\Oh\Co^\beta}$. By Lemmas \ref{Lem::Hold::Product} \ref{Item::Hold::Product::Hold2} and \ref{Lem::EllipticPara::HoloParaLem} \ref{Item::EllipticPara::HoloParaLem::EqvNorm}, for every $z\in\Hb^n$,
    \begin{align*}
        &\|fg(z,\cdot)\|_{\Co^\beta(\B^q)}\le\|f(z,\cdot)\|_{L^\infty}\|g(z,\cdot)\|_{\Co^\beta}+\|f(z,\cdot)\|_{\Co^\beta}\|g(z,\cdot)\|_{L^\infty}
        \\
        \le&\|f\|_{\Co^\alpha_\Oh L^\infty}\|g\|_{\Co^{-2}_\Oh\Co^\beta}\dist(z,\partial\Hb^n)^{-2}+\|f\|_{\Co^{-1}_\Oh\Co^\beta}\|g\|_{\Co^{\alpha-1}_\Oh L^\infty}\dist(z,\partial\Hb^n)^{-1}\cdot\dist(z,\partial\Hb^n)^{\alpha-1}
        \\
        \le&\|f\|_{\Co^{\alpha}_\Oh L^\infty\cap\Co^{-1}_\Oh\Co^\beta}\|g\|_{\Co^{\alpha-1}_\Oh L^\infty\cap\Co^{-2}_\Oh\Co^\beta}\dist(z,\partial\Hb^n)^{-2}.
    \end{align*}
	
	Dividing $\dist(z,\partial\Hb^n)^{-2}$ on both side, we get $\|fg\|_{\Co^{-2}_\Oh \Co^\beta}\lesssim_\alpha\|f\|_{\Co^{\alpha}_\Oh L^\infty\cap\Co^{-1}_\Oh\Co^\beta}\|g\|_{\Co^{\alpha-1}_\Oh L^\infty\cap\Co^{-2}_\Oh\Co^\beta}$, completing the proof.
\end{proof}
\begin{remark}
    By considering higher derivatives we can define $\Co^\alpha_\Oh$ for $\alpha\ge1$ naturally, where Lemma \ref{Lem::EllipticPara::HoloMult} still holds. On the other hand by Example \ref{Exam::Hold::MultUnBdd}, we see that $\Co^\alpha_\Oh$ in \eqref{Eqn::EllipticPara::HoloMult} cannot be replaced by $\Co^\alpha$ unless $\alpha>2$. This is partly because the following: for $\Co^\alpha_\Oh$-functions the high frequency information only occur near the boundary, whereas for $\Co^\alpha$-functions the singularity are allowed to be ``evenly spread''.
\end{remark}

\begin{proof}[Proof of Proposition \ref{Prop::EllipticPara::AnalyticPDE}]It suffices to consider $\frac12<\alpha<1$ and $\tilde\beta\in\R_+$. 

First we claim that by possibly shrinking $\eps_0>0$ we can assume that $\|f\|_{\Co^\alpha L^\infty\cap\Co^{-1}\Co^{\tilde\beta}}<\eps_0$.

Let $s_0\in\B^q$ be a fixed point, to prove the proposition it suffices to show that there exists a $\delta>0$ such that $\mathfrak f\in\Co^\infty_\loc\Co^{\tilde\beta}(\B^n,B^q(s_0,\delta);\C^m)$.

For $0<\mu\le1-|s_0|$, let $f_{s_0,\mu}(x,s):=f(x,\mu(s-s_0))$ be a function defined on $(x,s)\in\B^n\times\B^q$. By Lemma \ref{Lem::Hold::ScalingLem} we see that
\begin{equation*}
    \lim\limits_{\mu\to0}\|f_{s_0,\mu}\|_{\Co^\alpha L^\infty\cap\Co^{-1}\Co^{\tilde\beta}(\B^n,\B^q;\C^m)}\le\|f\|_{\Co^\alpha L^\infty(\B^n,\B^q;\C^m)}+\|f(\cdot,s_0)\|_{\Co^{-1}(\B^n;\C^m)}.
\end{equation*}

Since $\Co^\alpha L^\infty\cap\Co^{-1}\Co^{\tilde\beta}\subset C^0$, we have $\|f(\cdot,s_0)\|_{\Co^{-1}(\B^n)}\lesssim_\alpha \|f(\cdot,s_0)\|_{\Co^\alpha(\B^n)}\le\|f\|_{\Co^\alpha L^\infty(\B^n,\B^q)}$. Therefore, if $\|f\|_{\Co^\alpha L^\infty}<\eps_0$, then $\lim\limits_{\mu\to0}\|f_{s_0,\mu}\|_{\Co^\alpha L^\infty\cap\Co^{-1}\Co^{\tilde\beta}(\B^n,\B^q)}<C_\alpha\eps_0$ where $C_\alpha$ does not depend on $\eps_0$. 

Thus  we can find a $\mu\in(0,1-|s_0|]$ such that $\|f_{s_0,\mu}\|_{\Co^\alpha L^\infty\cap\Co^{-1}\Co^{\tilde\beta}(\B^n,\B^q)}<2C_\alpha\eps_0$. By replacing $\eps_0$ with $2C_\alpha \eps_0$, we have a $\mathfrak f_{s_0,\mu}\in\Co^\infty_\loc(\Hb^n,\B^q;\C^m)$ such that $\mathfrak f_{\mu,s_0}|_{\B^n\times\B^q}=f_{\mu,s_0}$. Taking scaling back we get $\mathfrak f\in\Co^\infty_\loc(\Hb^n, B^q(s_0,\mu);\C^m)$, completing the proof of the claim.


\medskip
Now we assume $\|f\|_{\Co^\alpha L^\infty\cap\Co^{-1}\Co^{\tilde\beta}}<\eps_0$, where $\eps_0>0$ is a small number to be determined.

We use $\Theta[u]=\Theta(u,\nabla_xu)$ for function $u$ either on $\B^n\times\B^q$ or on $\Hb^n\times\B^q$. Note that for $\tilde u\in\Co^{\alpha}_\Oh L^\infty\cap\Co^{-1}_\Oh\Co^{\tilde\beta} $ we have $\nabla_x\tilde u=\partial_z\tilde u$.
	
	The proof uses contraction mappings twice. The first one take place in $\Co^\alpha L^\infty(\B^n,\B^q;\C^m)$ and shows that $\Pv\Delta f$ has some particular properties. The second one take place in $\Co^\alpha_\Oh L^\infty\cap\Co^{-1}_\Oh\Co^{\tilde\beta}(\Hb^n,\B^q;\C^m)$, showing that $\Pv\Delta f$ is a restriction of a function in $\Hb^n\times\B^q$ which is holomorphic in $z$. 
	
	\medskip
	For $u\in\Co^\alpha  L^\infty(\B^n,\B^q;\C^m)$, by Lemma \ref{Lem::Hold::Product} (since $\alpha>\frac12$) we have $u\otimes\nabla_x u(\cdot,s)\in\Co^{\alpha-1}_x$ uniformly for almost every $s$, so $\Theta[u]\in\Co^{\alpha-1} L^\infty(\B^n,\B^q;\C^p)$ and thus $L\Theta[u]\in \Co^{\alpha-2}  L^\infty(\B^n,\B^q;\C^m)$. Moreover we have for every $u,v\in\Co^\alpha  L^\infty(\B^n,\B^q;\C^m)$
{\small\begin{equation}\label{Eqn::EllipticPara::AnaPDE::ContMap1}
\begin{aligned}
    \|L\Theta[u]\|_{\Co^{\alpha-2}L^\infty}&\lesssim \|\Theta(u,\nabla_x u)\|_{\Co^{\alpha-1}L^\infty}\lesssim\|u\otimes\nabla_x u\|_{\Co^{\alpha-1}L^\infty}\lesssim\|u\|_{\Co^{\alpha-2}L^\infty}^2,
    \\
        \|L\Theta[u]-L\Theta[v]\|_{\Co^{\alpha-2}L^\infty}&\lesssim\|\Theta(u-v,\nabla_x u)\|_{\Co^{\alpha-1}L^\infty}+\|\Theta(v,\nabla_x (u-v))\|_{\Co^{\alpha-1}L^\infty}\lesssim(\|u\|_{\Co^{\alpha}L^\infty}+\|v\|_{\Co^{\alpha}L^\infty})\|u-v\|_{\Co^{\alpha}L^\infty}.
\end{aligned}
\end{equation}}
Here the implicit constants depend only on $m,n,q,\alpha,L,\Theta$ but not on $u$ or $v$.

    For $\tilde u\in\Co^\alpha_\Oh L^\infty\cap\Co^{-1}_\Oh\Co^{\tilde\beta}(\Hb^n,\B^q;\C^m)$,  $L\Theta[\tilde u]$ is a function on $\Hb^n\times\B^q$ that is holormorphic in $z$. By Lemma \ref{Lem::EllipticPara::HoloMult}, $\tilde u\otimes\partial_z \tilde u\in \Co^{\alpha-1}_\Oh L^\infty\cap\Co^{-2}_\Oh\Co^{\tilde\beta}(\Hb^n,\B^q;\C^{m^2 n})$, we know $L\Theta[\tilde u]\in \Co^{\alpha-2}_\Oh L^\infty\cap\Co^{-3}_\Oh\Co^{\tilde\beta}(\Hb^n,\B^q;\C^m)$. Moreover, similar to \eqref{Eqn::EllipticPara::AnaPDE::ContMap2} we have for $\tilde u,\tilde v\in\Co^{\alpha}_\Oh L^\infty\cap\Co^{-1}_\Oh\Co^{\tilde\beta}(\Hb^n,\B^q;\C^m)$:
	\begin{equation}\label{Eqn::EllipticPara::AnaPDE::ContMap2}
	\begin{aligned}
    \|L\Theta[\tilde u]\|_{\Co^{\alpha-2}_{\Oh}L^\infty\cap\Co^{-3}_{\Oh}\Co^{\tilde\beta}}&\lesssim_\alpha\|\tilde u\|_{\Co^{\alpha-1}_{\Oh}L^\infty\cap\Co^{-1}_{\Oh}\Co^{\tilde\beta}}^2,\\
       \|L\Theta[\tilde u]-L\Theta[\tilde v]\|_{\Co^{\alpha-2}_{\Oh}L^\infty\cap\Co^{-3}_{\Oh}\Co^{\tilde\beta}}&\lesssim_\alpha(\|\tilde u\|_{\Co^{\alpha-1}_{\Oh}L^\infty\cap\Co^{-1}_{\Oh}\Co^{\tilde\beta}}+\|\tilde v\|_{\Co^{\alpha-1}_{\Oh}L^\infty\cap\Co^{-1}_{\Oh}\Co^{\tilde\beta}})\|\tilde u-\tilde v\|_{\Co^{\alpha-1}_{\Oh}L^\infty\cap\Co^{-1}_{\Oh}\Co^{\tilde\beta}}.
    \end{aligned}
\end{equation}
Here the implicit constants depend only on $m,n,q,\alpha,L,\Theta$ but not on $\tilde\beta$, $\tilde u$ or $\tilde v$.

	Let $\Pv:\Co^{\alpha-2}(\B^n)\to\Co^\alpha(\B^n)$ be as in Proposition \ref{Prop::HolLap}. So the assumption $\Delta f=L\Theta[f]$ implies that 
	\begin{equation}\label{Eqn::EllipticPara::AnaPDE::FixPtf}
	    f=f-\Pv\Delta f+\Pv\Delta f=(f-\Pv\Delta f)+\Pv L\Theta[\Pv\Delta f+(f-\Pv\Delta f)].
	\end{equation}
	
	Define $T_f$ on $\Co^\alpha  L^\infty(\B^n,\B^q;\C^m)$ as 
	$$T_f[u]:=(f-\Pv\Delta f)+\Pv L\Theta[u+(f-\Pv\Delta f)],\quad u\in\Co^\alpha L^\infty(\B^n,\B^q;\C^m).$$
	
	Since $\|f-\Pv\Delta f\|_{\Co^\alpha L^\infty(\B^n,\B^q;\C^m)}\lesssim\|f\|_{\Co^\alpha L^\infty(\B^n,\B^q;\C^m)}$. By \eqref{Eqn::EllipticPara::AnaPDE::ContMap1} and Remark \ref{Rmk::EllipticPara::HolLapPara}, we can find a $C_1=C_1(m,n,q,\alpha,L,\Theta,\Pv)>1$ that does not depend on $f$, that bounds \eqref{Eqn::EllipticPara::AnaPDE::ContMap1} and satisfies $C_1>\|\id-\Pv\Delta \|_{\Co^\alpha L^\infty}$. More precisely,  for every $f,u,v\in \Co^\alpha L^\infty(\B^n,\B^q,\C^m)$, we have
\begin{equation*}
    \begin{aligned}
    \|T_f[u]\|_{\Co^\alpha L^\infty}&\le C_1\|f\|_{\Co^\alpha L^\infty}+ C_1\big(\|u\|_{\Co^\alpha L^\infty}+C_1\|f\|_{\Co^\alpha L^\infty}\big)^2;
    \\
    \|T_f[u]-T_f[v]\|_{\Co^\alpha L^\infty}&\le C_1\left(\|u\|_{\Co^\alpha L^\infty}+\|v\|_{\Co^\alpha L^\infty}+C_1\|f\|_{\Co^\alpha L^\infty}\right)\|u-v\|_{\Co^\alpha L^\infty};
    \\
    \|f-\Pv\Delta f\|_{\Co^\alpha L^\infty}&\le C_1\|f\|_{\Co^\alpha L^\infty}.
    \end{aligned}
\end{equation*}

Take $\eps_1>0$ such that $9\eps_1C_1+10\eps_1 C_1^2<1$. If $\|f\|_{\Co^\alpha  L^\infty}<\eps_1$ and $\|u\|_{\Co^\alpha  L^\infty},\|v\|_{\Co^\alpha  L^\infty}\le 2C_1\eps_1$, then 
\begin{gather*}
    \|T_f[u]\|_{\Co^\alpha L^\infty}\le C_1\eps_1+C_1(3C_1\eps_1)^2<2C_1\eps_1,\\
    \|T_f[u]-T_f[v]\|_{\Co^\alpha L^\infty}<C_1(5C_1\eps_1)\|u-v\|_{\Co^\alpha L^\infty}<\tfrac12\|u-v\|_{\Co^\alpha L^\infty}.
\end{gather*}

Therefore when $\|f\|_{\Co^\alpha  L^\infty}<\eps_1$, $T_f$ is a contraction mapping in $\{u\in\Co^\alpha   L^\infty(\B^n,\B^q;\C^m):\|u\|_{\Co^\alpha L^\infty}\le2C_1\eps_1\}$. In such case, there is a unique $u$ satisfying $\|u\|_{\Co^\alpha L^\infty}\le2C_1\eps_1$ such that $T_f[u]=u$. On the other hand, $f\in\{u\in\Co^\alpha   L^\infty(\B^n,\B^q;\C^m):\|u\|_{\Co^\alpha L^\infty}\le2C_1\eps_1\}$ satisfies \eqref{Eqn::EllipticPara::AnaPDE::FixPtf}. Thus $f$ is that unique fixed point.

\medskip
Similarly, by Proposition \ref{Prop::HolLap} \ref{Item::HolLap::P} (also see Remark \ref{Rmk::EllipticPara::HolLapPara}), we have boundedness $\Pv:\Co^{-3}(\B^n)\to\Co^{-1}(\B^n)$, so $\|f-\Pv\Delta f\|_{\Co^\alpha L^\infty\cap\Co^{-1} \Co^{\tilde\beta}}\lesssim\|f\|_{\Co^\alpha L^\infty\cap\Co^{-1} \Co^{\tilde\beta}}$. Note that $f-\Pv\Delta f$ is harmonic function in $x$-variable, so by Proposition \ref{Prop::HolLap}  \ref{Item::HolLap::E} and Remark \ref{Rmk::EllipticPara::HolLapPara} we have $\Ex(f-\Pv\Delta f)\in \Co^\alpha_\Oh L^\infty\cap\Co^{-1}_\Oh\Co^{\tilde\beta}(\Hb^n,\B^q;\C^m)$ with 
\begin{equation}\label{Eqn::EllipticPara::AnaPDE::BddExf}
    \|\Ex(f-\Pv\Delta f)\|_{\Co^\alpha_\Oh L^\infty\cap\Co^{-1}_\Oh \Co^{\tilde\beta}(\Hb^n,\B^q;\C^m)}\lesssim_{\alpha}\|f\|_{\Co^\alpha L^\infty\cap\Co^{-1}\Co^{\tilde\beta}(\B^n,\B^q;\C^m)}.
\end{equation}

For a $f\in \Co^\alpha L^\infty\cap\Co^{-1} \Co^{\tilde\beta}(\B^n,\B^q;\C^m)$, define $\tilde T_f$ on $\Co^\alpha_\Oh L^\infty\cap\Co^{-1}_\Oh\Co^{\tilde\beta}(\Hb^n,\B^q;\C^m)$ as 
$$\tilde T_f[\tilde u]=\Ex(f-\Pv\Delta f)+\tilde\Pv L\Theta[\tilde u+\Ex(f-\Pv\Delta f)],\quad\tilde u\in \Co^\alpha_\Oh L^\infty\cap\Co^{-1}_\Oh\Co^{\tilde\beta}(\Hb^n,\B^q;\C^m).$$
Here $\tilde\Pv$ is as in Proposition \ref{Prop::HolLap} \ref{Item::HolLap::TildeP}.

Take $C_2=C_2(m,n,q,\alpha,L,\Theta,\tilde\Pv,\Ex)>1$ that bounds \eqref{Eqn::EllipticPara::AnaPDE::ContMap2} and \eqref{Eqn::EllipticPara::AnaPDE::BddExf}. Therefore, for every $\tilde u,\tilde v\in\Co^\alpha_\Oh L^\infty\cap\Co^{-1}_\Oh\Co^{\tilde\beta}(\Hb^n,\B^q;\C^m)$,
{\small\begin{align*}
    \|\tilde T_f[\tilde u]\|_{\Co^\alpha_{\Oh} L^\infty\cap\Co^{-1}_{\Oh}\Co^{\tilde\beta}}&\le C_2\|f\|_{\Co^\alpha L^\infty\cap\Co^{-1}\Co^{\tilde\beta}}+ C_2\big(\|\tilde u\|_{\Co^\alpha_{\Oh} L^\infty\cap\Co^{-1}_{\Oh}\Co^{\tilde\beta}}+C_2\|f\|_{\Co^\alpha L^\infty\cap\Co^{-1}\Co^{\tilde\beta}}\big)^2;
    \\
    \|\tilde T_f[\tilde u]-\tilde T_f[\tilde v]\|_{\Co^\alpha_{\Oh} L^\infty\cap\Co^{-1}_{\Oh}\Co^{\tilde\beta}}&\le C_2\big(\|\tilde u\|_{\Co^\alpha_\Oh L^\infty\cap\Co^{-1}_\Oh\Co^{\tilde\beta}}+\|\tilde v\|_{\Co^\alpha_\Oh L^\infty\cap\Co^{-1}_\Oh\Co^{\tilde\beta}}+C_2\|f\|_{\Co^\alpha L^\infty\cap\Co^{-1}\Co^{\tilde\beta}}\big)\|\tilde u-\tilde v\|_{\Co^\alpha_\Oh L^\infty\cap\Co^{-1}_\Oh\Co^{\tilde\beta}}.
\end{align*}}

By Lemma \ref{Lem::EllipticPara::HoloParaLem} \ref{Item::EllipticPara::HoloParaLem::Res} we can find a $C_3=C_3(m,n,q,\alpha)>0$ such that,
\begin{equation}\label{Eqn::EllipticPara::AnaPDE::RestrictionTmp}
    \|\tilde u|_{\B^n\times\B^q}\|_{\Co^\alpha  L^\infty(\B^n,\B^q;\C^m)}\le C_3\|\tilde u\|_{\Co^\alpha_\Oh L^\infty\cap\Co^{-1}_\Oh\Co^{\tilde\beta}(\Hb^n,\B^q;\C^m)},\quad \forall u\in\Co^\alpha_\Oh L^\infty\cap\Co^{-1}_\Oh\Co^{\tilde\beta}(\Hb^n,\B^q;\C^m).
\end{equation}
Take $\eps_0>0$ such that $C_3C_2\eps_0 <\eps_1$ and $9\eps_0 C_2+10\eps_0 C_2^2<1$. So $\|f\|_{\Co^\alpha L^\infty\cap\Co^{-1}\Co^{\tilde\beta}}<\eps_0$ implies $\|f\|_{\Co^\alpha L^\infty}<\eps_1$.
For such $f$, we know $\tilde T_f$ is a contraction mapping in $\{\tilde u\in\Co^\alpha_{\Oh} L^\infty\cap\Co^{-1}_{\Oh}\Co^{\tilde\beta}(\Hb^n,\B^q;\C^m):\|\tilde u\|_{\Co^\alpha_{\Oh} L^\infty\cap\Co^{-1}_{\Oh}\Co^{\tilde\beta}}\le2C_2\eps_0\}$.
Thus there is a unique $\tilde u\in\Co^\alpha_\Oh L^\infty\cap\Co^{-1}_\Oh\Co^{\tilde\beta}(\Hb^n,\B^q;\C^m) $ such that $\|\tilde u\|_{\Co^\alpha_\Oh}\le2C_2\eps_0$ and $\tilde u=\tilde T_f[\tilde u]$. 

By Proposition \ref{Prop::HolLap} \ref{Item::HolLap::TildeP}, 
$$\tilde u|_{\B^n\times \B^q}=\big(\Ex(f-\Pv\Delta f)+\tilde\Pv L\Theta[\tilde u+\Ex(f-\Pv\Delta f)]\big)\big|_{\B^n\times \B^q}=f-\Pv\Delta f+\Pv L\Theta[\tilde u|_{\B^n\times \B^q}+f-\Pv\Delta f]=T_f[\tilde u|_{\B^n\times \B^q}].$$
By \eqref{Eqn::EllipticPara::AnaPDE::RestrictionTmp} and the assumption $C_3C_2\eps_0<\eps_1 <2C_1\eps_1$, we have  $\|\tilde u|_{\B^n\times \B^q}\|_{\Co^\alpha  L^\infty(\B^n,\B^q;\C^m)}\le 2C_1\eps_1$. By the uniqueness of the fixed point for $T_f$, we get $f=\tilde u|_{\B^n\times \B^q}$. 

So $\mathfrak f=\tilde u\in\Co^{-1}_\Oh\Co^{\tilde\beta}(\Hb^n,\B^q;\C^m)$ is the desired extension of $f$ to $\Hb^n\times\B^q$.

Finally, by Lemma \ref{Lem::Hold::NablaHarm} we have $\mathfrak f\in\Co^\infty_\loc \Co^{\tilde\beta}(\Hb^n,\B^q;\C^m)\subset \Co^\infty_\loc \Co^{\tilde\beta}_\loc(\Hb^n,\B^q;\C^m)$, completing the proof.
\end{proof}

\subsection{Some transition lemmas}\label{Section::KeyLems}
We first construct the coordinate chart $L$.
\begin{lem}\label{Lem::EllipticPara::IniNorm}
Endow $\R^r,\C^m,\R^q$ with standard coordinate systems $\tau=(\tau^1,\dots,\tau^r)$, $w=(w^1,\dots,w^m)$, $s=(s^1,\dots,s^q)$ respectively.

Let $\alpha\in(\frac12,\infty)$, $\beta\in(0,\infty)$. Let $\Se\le(\C T\Mf)\times\Nf$ and $(u_0,v_0)\in\Mf\times\Nf$ be as in the assumption of Theorem \ref{KeyThm::EllipticPara}. Then for any $\delta>0$ there are neighborhoods $U_0\subseteq\Mf$ of $p_0$, $V_0\subseteq\Nf$ of $q_0$ and smooth coordinate charts $(L',L''):U_0\xrightarrow{\sim}\B^{r+2m}\subset\R^r_\tau\times\C^m_w$, $L''':V_0\xrightarrow{\sim}\B^q_s$, such that for the product chart $L=(L',L'',L''')$, $L_*\Se$ has a local basis $X=[X_1,\dots,X_{r+m}]^\top$ with the form \eqref{Eqn::EllipticPara::DefofA} defined on $\B^{r+2m}_{\tau,w}\times\B^q_s$, and
\begin{equation}\label{Eqn::EllipticPara::IniNorm::A}
    \|A\|_{\Co^{\alpha,\beta}(\B^{r+2m},\B^q;\C^{(r+m)\times m})}<\delta.
\end{equation}
\end{lem}
\begin{proof}Let $\tilde L''':\tilde V_0\subseteq\Nf\to\R^q_s$ be an arbitrary smooth coordinate chart near $v_0$ such that $\tilde L'''(v_0)=0$. Since $\rank\Se=r+m$ and $\rank(\Se+\bar\Se)=r+2m$, by standard linear algebra argument (see also \cite[Page 18]{Involutive}) we can find a smooth coordinate chart $(\tilde L',\tilde L''):\tilde U_0\subseteq\Mf\to\R^r_\tau\times\C^m_w$ such that
\begin{equation}\label{Eqn::EllipticPara::PfIniNorm::GenAssum1}
    \textstyle(\tilde L',\tilde L'')(u_0)=(0,0)\text{ and }\Se_{(u_0,v_0)}=(\tilde L',\tilde L'')^*\Span\big(\Coorvec{\tau^1},\dots,\Coorvec{\tau^r},\Coorvec{w^1},\dots,\Coorvec{w^m}\big)\big|_{u_0}\le \C T_{u_0}\Mf.
\end{equation}
Therefore, by considering $\tilde L:=(\tilde L',\tilde L'',\tilde L''')$ and working on $\Se_{(u_0,v_0)}\le\C T_{(u_0,v_0)}(\Mf\times\Nf)$, we have
\begin{equation}\label{Eqn::EllipticPara::PfIniNorm::GenAssum2}
    \textstyle(\tilde L_*\Se)_{(0,0,0)}\oplus\Span\big(\Coorvec{\bar w^1}\big|_{(0,0,0)},\dots,\Coorvec{\bar w^m}\big|_{(0,0,0)},\Coorvec{s^1}\big|_{(0,0,0)},\dots\Coorvec{s^q}\big|_{(0,0,0)}\big)=\R^r\times\C^m\times\R^q.
\end{equation}

Applying Lemma \ref{Lem::ODE::GoodGen}, since we have \eqref{Eqn::EllipticPara::PfIniNorm::GenAssum2}, there is a $\lambda_0>0$ such that $\lambda_0\B^{r+2m}_{\tau,w}\times\lambda_0\B^q_s\subseteq\tilde L(U_0\times V_0)$ and we can find a $\Co^\alpha$-local basis $\tilde X=[\tilde X_1,\dots,\tilde X_{r+m}]^\top$ of $\tilde L_*\Se$ that is defined on $\lambda_0\B^{r+2m}_{\tau,w}\times\lambda_0\B^q_s$ and has the form,
\begin{equation*}
    \tilde X=\begin{pmatrix}\tilde X'\\\tilde X''\end{pmatrix}=\begin{pmatrix}I_r&&\tilde A'&\tilde A'''\\&I_m&\tilde A''&\tilde A''''\end{pmatrix}\begin{pmatrix}\partial_\tau\\\partial_w\\\partial_{\bar w}\\\partial_s\end{pmatrix},\quad\text{on }\lambda_0\B^{r+2m}_{\tau,w}\times\lambda_0\B^q_s.
\end{equation*}
Clearly $\tilde A'''=0$ and $\tilde A''''=0$ since $\Se\le(\C T\Mf)\times\Nf$. By \eqref{Eqn::EllipticPara::PfIniNorm::GenAssum1} we have $\tilde A'(0)=0$ and $\tilde A''(0)=0$.

Write $\tilde A=\begin{pmatrix}\tilde A'\\\tilde A''\end{pmatrix}$. By Lemma \ref{Lem::Hold::ScalingLem} there is a  $C=C(r,m,q,\lambda_0,\alpha,\beta)>0$ such that 
\begin{equation}\label{Eqn::EllipticPara::IniNorm::Scaling}
    \|\tilde A(\lambda\cdot)\|_{\Co^{\alpha,\beta}(\B^{r+2m},\B^q;\C^{(r+m)\times m})}\le C\lambda^{\min(\beta,\frac12)}\|\tilde A\|_{\Co^{\alpha,\beta}(\lambda_0\B^{r+2m},\lambda_0\B^q;\C^{(r+m)\times m})},\quad\forall \lambda\in(0,\lambda_0].
\end{equation}

Take $\lambda=(C\|\tilde A\|_{\Co^{\alpha,\beta}(\lambda_0\B^{r+2m},\lambda_0\B^q)})^{-2}\cdot(\frac\delta2)^2$. Then the right hand side of \eqref{Eqn::EllipticPara::IniNorm::Scaling} is less than $\delta$. 

Take $L:=\lambda\tilde L$, we see that $L_*\Se$ has a $\Co^{\alpha,\beta}$-local basis $X=[X_1,\dots,X_{r+m}]^\top$ with the form
\begin{equation*}
    X=\begin{pmatrix}I_r&&A'\\&I_m& A''\end{pmatrix}\begin{pmatrix}\partial_\tau\\\partial_w\\\partial_{\bar w}\end{pmatrix},\quad\text{where }\begin{pmatrix}
    A'(\tau,w,s)\\A''(\tau,w,s)
    \end{pmatrix}=\begin{pmatrix}
    \tilde A'(\lambda\tau,\lambda w,\lambda s)\\\tilde A''(\lambda\tau,\lambda w,\lambda s)
    \end{pmatrix}\text{ for }(\tau,w)\in\B^{r+2m},s\in\B^q.
\end{equation*}

So we have $\|A\|_{\Co^{\alpha,\beta}(\B^{r+2m},\B^q)}=\|\tilde A(\lambda\cdot)\|_{\Co^{\alpha,\beta}(\lambda_0\B^{r+2m},\lambda_0\B^q)}<\delta$ as desired.
\end{proof}

Next we show that \eqref{Eqn::EllipticPara::KeyEqn1} and \eqref{Eqn::EllipticPara::ExistenceH} are associated. Here we endow $\R^r\times\C^m$ with two standard coordinate systems $(\tau,w)=(\tau^1,\dots,\tau^r,w^1,\dots,w^m)$ and $(\sigma,\zeta)=(\sigma^1,\dots,\sigma^r,\zeta^1,\dots,\zeta^m)$, and we still use $s=(s^1,\dots,s^q)$ as the standard coordinate system for $\R^q$.
\begin{lem}\label{Lem::EllipticPara::NewGen}Let $\eps,\delta>0$, $A:\B^{r+2m}_{\tau,w}\times\B^q_s\to\C^{(r+m)\times m}$ and $H:\B^{r+2m}_{\tau,w}\times\B^q_s\to\B^{r+2m}_{\sigma,\zeta}\times\B^q_s$ be as in (the assumptions and consequence) of Proposition \ref{Prop::EllipticPara::ExistPDE}. Then,
\begin{enumerate}[parsep=-0.3ex,label=(\roman*)]
    \item\label{Item::EllipticPara::NewGen::B} If $L_*\Se$ has a local basis with the form \eqref{Eqn::EllipticPara::DefofA}, then $(H\circ L)_*\Se$ has a local basis with the form \eqref{Eqn::EllipticPara::DefofB}.
    \item\label{Item::EllipticPara::NewGen::Eqn} Suppose $\Se$ is involutive and $H''$ is a solution to \eqref{Eqn::EllipticPara::ExistenceH}, then the coefficient matrix $B=[B',B'']^\top$ in \eqref{Eqn::EllipticPara::DefofB} is a solution to \eqref{Eqn::EllipticPara::KeyEqnB}.
\end{enumerate} 
\end{lem}
\begin{proof}
Following Convention \ref{Conv::EllipticPara::ConvofExtPDE}, by assumption $H(\tau,w,s)=(\tau,H''(\tau,w,s),s)$, we have
\begin{equation}\label{Eqn::EllipticPara::NewGen::H}
    H^*(d\sigma,d\zeta,d\bar\zeta,ds)=(d\tau,dw,d\bar w,ds)\begin{pmatrix}I&H''_\tau&\bar H''_\tau\\&H''_w&\bar H''_w\\&H''_{\bar w}&\bar H''_{\bar w}\\& H''_s&\bar H''_s&I\end{pmatrix}\ \Rightarrow\ \begin{pmatrix}
    \Coorvec\tau\\\Coorvec w\\\Coorvec{\bar w}\end{pmatrix}=\begin{pmatrix}I&H''_\tau&\bar H''_\tau\\&H''_w&\bar H''_w\\&H''_{\bar w}&\bar H''_{\bar w}\end{pmatrix}\cdot \begin{pmatrix}
    H^*\Coorvec\sigma\\H^*\Coorvec\zeta\\H^*\Coorvec{\bar\zeta}\end{pmatrix}.
\end{equation}
Therefore taking pushforward of \eqref{Eqn::EllipticPara::DefofA} under $H$ we have, \begin{equation*}
    H_*X=H_*\left(\begin{pmatrix}I&&A'\\&I&A''\end{pmatrix}\begin{pmatrix}I&H''_\tau&\bar H''_\tau\\&H''_w&\bar H''_w\\&H''_{\bar w}&\bar H''_{\bar w}\end{pmatrix}\right)\begin{pmatrix}
    \Coorvec\sigma\\\Coorvec\zeta\\\Coorvec{\bar\zeta}\end{pmatrix}=\left(\begin{pmatrix}I&H''_\tau+A'H''_{\bar w}&\bar H''_\tau+A'\bar H''_{\bar w}\\&H''_w+A''\bar H''_{\bar w}&\bar H''_w+A''\bar H''_{\bar w}\end{pmatrix}\circ H^\Inv\right)\cdot\begin{pmatrix}
    \Coorvec\sigma\\\Coorvec\zeta\\\Coorvec{\bar\zeta}\end{pmatrix}.
\end{equation*}
The matrix $\begin{pmatrix}I&H''_\tau+A'H''_{\bar w}\\&H''_w+A''\bar H''_{\bar w}\end{pmatrix}$ is invertible (see Lemma \ref{Lem::EllipticPara::MatisInv} \ref{Item::EllipticPara::MatisInv::Inv}), we see that $H_*(L_*\Se)=(H\circ L)_*\Se$ has a local basis
\begin{equation*}
    H_*\left(\begin{pmatrix}I&H''_\tau+A'H''_{\bar w}\\&H''_w+A''\bar H''_{\bar w}\end{pmatrix}^{-1}\begin{pmatrix}X'\\X''\end{pmatrix}\right)=\left(\left(I_{r+m}+\begin{pmatrix}I&H''_\tau+ A'H''_{\bar w}\\&H''_ w+ A''H''_{\bar w}\end{pmatrix}^{-1}\begin{pmatrix}\bar H''_\tau+ A'\bar H''_{\bar w}\\\bar H''_ w+ A''\bar H''_{\bar w}\end{pmatrix}\right)\circ H^\Inv\right)\cdot\begin{pmatrix}
    \Coorvec\sigma\\\Coorvec\zeta\\\Coorvec{\bar\zeta}\end{pmatrix}.
\end{equation*}
This gives \eqref{Eqn::EllipticPara::DefofB} and finishes the proof of \ref{Item::EllipticPara::NewGen::B}.

\medskip Recall the notations $B'=(b_j^l)_{\substack{1\le j\le r;1\le l\le m}}$ and $B''=(b_{r+k}^l)_{\substack{1\le k,l\le m}}$. By \eqref{Eqn::EllipticPara::DefofB} we have,
\begin{equation}\label{Eqn::EllipticPara::NewGen::Tmp1}
    b_j^l\circ H=\left(\begin{pmatrix}I&H''_\tau+ A'H''_{\bar w}\\&H''_ w+ A''H''_{\bar w}\end{pmatrix}^{-1}\!\!\begin{pmatrix}\bar H''_\tau+ A'\bar H''_{\bar w}\\\bar H''_ w+ A''\bar H''_{\bar w}\end{pmatrix}\right)_j^l,\quad 1\le j\le m+r,\quad 1\le l\le m.
\end{equation}

By \eqref{Eqn::EllipticPara::NewGen::H} we have, for $1\le j\le r$ and $1\le k\le m$,
\begin{equation}\label{Eqn::EllipticPara::NewGen::Tmp2}
    \left.\Coorvec{\sigma^j}\right|_{H(\tau,w,s)}=\left(\begin{pmatrix}I&H''_\tau&\bar H''_\tau\\&H''_w&\bar H''_w\\&H''_{\bar w}&\bar H''_{\bar w}\end{pmatrix}^{-1}\;\begin{pmatrix}
    \Coorvec\tau\\\Coorvec w\\\Coorvec{\bar w}\end{pmatrix}\right)_j,
    \
    \left.\Coorvec{\bar\zeta^k}\right|_{H(\tau,w,s)}=\left(\begin{pmatrix}I&H''_\tau&\bar H''_\tau\\&H''_w&\bar H''_w\\&H''_{\bar w}&\bar H''_{\bar w}\end{pmatrix}^{-1}\;\begin{pmatrix}
    \Coorvec\tau\\\Coorvec w\\\Coorvec{\bar w}\end{pmatrix}\right)_{r+m+k}.
\end{equation}

Plugging \eqref{Eqn::EllipticPara::NewGen::Tmp1} and \eqref{Eqn::EllipticPara::NewGen::Tmp2} into \eqref{Eqn::EllipticPara::ExistenceH} we get $\sum_{j=1}^r\frac{\partial b_j^l}{\partial\sigma^j}+\sum_{k=1}^m\frac{\partial b_{r+k}^l}{\partial\bar\zeta^k}=0$ for $1\le l\le m$, which is the last row of \eqref{Eqn::EllipticPara::KeyEqnB}.

Since $\Se$ is involutive, $(H\circ L)_*\Se$ is also involutive, by Lemma \ref{Lem::ODE::GoodGen} \ref{Item::ODE::GoodGen::InvComm}, $T_1,\dots,T_r,Z_1,\dots,Z_m$ given in \eqref{Eqn::EllipticPara::DefofB} are pairwise commutative. Writing out the condition $[T_j,T_{j'}]=0$ $(1\le j<j'\le r)$, $[T_j,Z_k]=0$ $(1\le j\le r,1\le k\le m)$ and $[Z_k,Z_{k'}]=0$ $(1\le k<k'\le m)$ in the coordinate components we get the first three rows of \eqref{Eqn::EllipticPara::KeyEqnB}, finishing the proof of \ref{Item::EllipticPara::NewGen::Eqn}.
\end{proof}

Now we have the PDE system \eqref{Eqn::EllipticPara::KeyEqnB}. We show that it satisfies the assumption of Proposition \ref{Prop::EllipticPara::AnalyticPDE}.
\begin{lem}\label{Lem::EllipticPara::BSatAna}
Suppose $B:\B^{r+2m}_{\sigma,\zeta}\times\B^q_s\to\C^{(r+m)\times m}$ solves \eqref{Eqn::EllipticPara::KeyEqnB}, then $f(\sigma,\xi,\eta,s):=B(\sigma,\frac12(\xi+i\eta),s)$ is defined on $\B^{r+2m}_{\sigma,\xi,\eta}\times\B^q_s$ and $f$ solves a PDE with the form $Lf=\Delta_{\sigma,\xi,\eta}\Theta(f,\nabla_{\sigma,\xi,\eta}f)$ where $L$ and $\Theta$ satisfy the assumptions of Proposition \ref{Prop::EllipticPara::AnalyticPDE}.
\end{lem}
\begin{proof}
We can write \eqref{Eqn::EllipticPara::KeyEqnB} as the form $DB=\tilde\Theta(B,\nabla_{\sigma,\xi,\eta}B)$ where the coefficients of $D$ and $\Theta$ take value in $\C^{m\times\frac{(r+m)(r+m-1)}2+m}$.

By direct computation (see also \cite[Lemma B.5]{SharpElliptic}) we have $D^*D=\sum_{j=1}^r\frac{\partial^2}{\partial\sigma^j\partial\sigma^j}+\sum_{k=1}^m\frac{\partial^2}{\partial\zeta^k\partial\bar\zeta^k}=\Delta_\sigma+\Box_\zeta$, acting on the components of the vector value functions. 

Using the real coordinates $\xi=\re\zeta$, $\eta=\im\zeta$, we see that $\Box_\zeta=\frac14\sum_{k=1}^m(\frac{\partial^2}{\partial\xi^k\partial\xi^k}+\frac{\partial^2}{\partial\eta^k\partial\eta^k})=\frac14\Delta_{\xi,\eta}$. So for $f(\sigma,\xi,\eta,s)=B(\sigma,\frac{\xi+i\eta}2,s)$ we have $(D^*DB)(\sigma,\frac{\xi+i\eta}2,s)=(\Delta_\sigma+\Delta_{\xi,\eta})f(\sigma,\xi,\eta,s)$ and $\nabla_{\xi,\eta}f(\sigma,\xi,\eta,s)=\frac12(\nabla_{\xi,\eta}B)(\sigma,\frac{\xi+i\eta}2,s)$. Therefore $f$ satisfies $\Delta_{\sigma,\xi,\eta}f=D^*\tilde \Theta(f,(\nabla_\sigma f,2\nabla_{\xi,\eta}f))$. Taking $L:=D^*$ and $\Theta(f,\nabla_{\sigma,\xi,\eta}f):=\tilde \Theta(f,(\nabla_\sigma f,2\nabla_{\xi,\eta}f))$ we finish the proof.
\end{proof}

Finally we find the map $G''$ in \eqref{Eqn::EllipticPara::FacterizationofF} using Proposition \ref{Prop::ODE::ParaHolFro}. We work on $\frac12\B^{r+2m}\subset\R^{r+2m}_{\sigma,\xi,\eta}$ and its ``complex extension'' $\frac12\Hb^{r+2m}\subset\C^{r+2m}$.

We endow $\C^{r+2m}$ with standard complex coordinates $(\boldsymbol\sigma,\boldsymbol\xi,\boldsymbol\eta)=(\boldsymbol\sigma^1,\dots,\boldsymbol\sigma^r,\boldsymbol\xi^1,\dots,\boldsymbol\xi^m,\boldsymbol\eta^1,\dots,\boldsymbol\eta^m)$ in the way that $\sigma=\re\boldsymbol\sigma$, $\xi=\re\boldsymbol\xi$ and $\eta=\re\boldsymbol\eta$. We denote $\Coorvec{\boldsymbol\zeta}:=\frac12(\Coorvec{\boldsymbol\xi}-i\Coorvec{\boldsymbol\eta})=[\frac12(\Coorvec{\boldsymbol\xi^1}-i\Coorvec{\boldsymbol\eta^1}),\dots,\frac12(\Coorvec{\boldsymbol\xi^m}-i\Coorvec{\boldsymbol\eta^m})]^\top$. 

Note that in this place, we temporarily forget the original complex structure of $\C^m_\zeta$. The objects $\boldsymbol\zeta=\boldsymbol\xi+i\boldsymbol\eta$ and $\Coorvec{\boldsymbol\zeta}$ are just conventions.

In this setting, $T_1,\dots,T_r,Z_1,\dots,Z_m$ given in \eqref{Eqn::EllipticPara::DefofB} are complex vector fields defined on $\frac12\B^{r+2m}\times\frac12\B^q\subset\R^{r+2m}_{\sigma,\xi,\eta}\times\R^q_s$.
\begin{lem}\label{Lem::EllipticPara::FinalModf}
Let $\tilde\beta\in\{\gamma,\gamma-:\gamma\in\R_+\}$. Suppose $T_1,\dots,T_r,Z_1,\dots,Z_m$ in \eqref{Eqn::EllipticPara::DefofB} are commutative, and the coefficient matrix function $B:\frac12\B^{r+2m}_{\sigma,\xi,\eta}\times\frac12\B^q_s\to\C^{(r+m)\times m}$ has a $(\sigma,\xi,\eta)$-variable holomorphic extension $\Bf\in\Co^\infty_\loc\Co^{\tilde\beta}(\frac12\Hb^{r+2m},\frac12\B^q;\C^{(r+m)\times m})$. Then for any $(\sigma_0,\xi_0,\eta_0,s_0)\in \frac12\B^{r+2m}\times\frac12\B^q$ there is a neighborhood $U_1\times V_1\subset\frac12\B^{r+2m}\times\frac12\B^q$ of $(\sigma_0,\xi_0,\eta_0,s_0)$ and a map $G''\in\Co^\infty\Co^{\tilde\beta}(U_1,V_1;\C^m)$ such that $G''(\sigma_0,\xi_0,\eta_0,s_0)=0$ and
\begin{equation*}
    (\Span(T_1(\cdot,s),\dots,T_r(\cdot,s),Z_1(\cdot,s),\dots,Z_m(\cdot,s)))^\bot|_{U_1}=\Span(dG''(\cdot,s)^1,\dots,dG''(\cdot,s)^m)\le\C T^*U_1,\quad\forall\ s\in V_1.
\end{equation*}
\end{lem}
Recall the notation of dual bundle in \eqref{Eqn::ODE::DualConjBundle}.
\begin{proof}
We define collections of vector fields $\Tf=[\Tf_1,\dots,\Tf_r]^\top$ and $\Zf=[\Zf_1,\dots,\Zf_m]^\top$  as
\begin{equation*}
    \begin{pmatrix}\Tf\\\Zf\end{pmatrix}:=\begin{pmatrix}I_r&&\Bf'\\&I_m&\Bf''\end{pmatrix}\begin{pmatrix}\partial_{\boldsymbol\sigma}\\\partial_{\boldsymbol\zeta}\\\partial_{\bar{\boldsymbol\zeta}}\end{pmatrix}.
\end{equation*}
Here $\Bf'$ is the upper $r\times m$ block of $\Bf$, and $\Bf''$ is the lower $m\times m$ block of $\Bf$.

Thus $(T_1,\dots,T_r,Z_1,\dots,Z_m)$  is the ``real domain restriction'' of $(\Tf_1,\dots,\Tf_r,\Zf_1,\dots,\Zf_m)$ in the sense of Proposition \ref{Prop::ODE::ParaHolFro}. We see that the coefficients of $[\Tf_j,\Tf_{j'}]$, $[\Tf_j,\Zf_k]$ and $[\Zf_k,\Zf_{k'}]$ are the holomorphic extension those of $[T_j,T_{j'}]$, $[T_j,Z_k]$ and $[Z_k,Z_{k'}]$. 

By assumption $T_1,\dots,T_r,Z_1,\dots,Z_m$ are commutative, so $[T_j,T_{j'}]=[T_j,Z_k]=[Z_k,Z_{k'}]=0$. Since a holomorphic function in $\Hb^{r+2m}$ that vanishes in $\B^{r+2m}$ must be a zero function, we know $[\Tf_j,\Tf_{j'}]=[\Tf_j,\Zf_k]=[\Zf_k,\Zf_{k'}]=0$ as well. Thus $\Tf_1,\dots,\Tf_r,\Zf_1,\dots,\Zf_m$ are commutative.

Now the assumptions of Proposition \ref{Prop::ODE::ParaHolFro} are satisfied. By Proposition \ref{Prop::ODE::ParaHolFro} we can find a neighborhood $\tilde U_1\times V_1\subseteq\frac12\Hb^{r+2m}\times\frac12\B^q$ of $(\sigma_0,\xi_0,\eta_0,s_0)$ and a map $\widetilde G''\in\Co^\infty\Co^{\tilde\beta}(\tilde U_1,V_1;\C^m)$ such that $\widetilde G''(\sigma_0,\xi_0,\eta_0,s_0)=0$ and for $U_1:=\tilde U_1\cap\R^{r+2m}$ we have for each $s\in V_1$,
$$(\Span(T_1(\cdot,s),\dots,T_r(\cdot,s),Z_1(\cdot,s),\dots,Z_m(\cdot,s)))^\bot|_{U_1}=\Span(d(\widetilde G''(\cdot,s)^1|_{U_1}),\dots,d(\widetilde G''(\cdot,s)^m|_{U_1}))\le \C T^*U_1.$$

Clearly $\widetilde G''|_{U_1\times V_1}\in\Co^\infty\Co^{\tilde\beta}(U_1,V_1;\C^m)$ holds. Thus, $G'':=\widetilde G''|_{U_1\times V_1}$ is as desired.
\end{proof}

\subsection{Proof of Theorem \ref{KeyThm::EllipticPara}}\label{Section::KeyPf}

In the proof we often consider the local bases of cotangent subbundles. 
\begin{conv}\label{Conv::EllipticPara::CoSpan}
Let $M$ and $N$ be two smooth manifolds, let $\Tc\le (\C TM)\times N$ be a continuous subbundle. Let $\mu=(\mu^1,\dots,\mu^n):M\to\R^n$ and $\nu=(\nu^1,\dots,\nu^q):N\to\R^q$ be two smooth coordinate charts and let $f=(f^1,\dots,f^{n-r}):M\times N\to\C^{n-r}$ be a continuous map which is $C^1$ in $M$. We say $\Tc^\bot$ is spanned by $df^1,\dots,df^{n-r},d\nu^1,\dots,d\nu^q$, denoted by $\Tc^\bot=\Span(df,d\nu)$, if
\begin{equation}\label{Eqn::EllipticPara::CoSpan}
    \Span\Big(\sum_{j=1}^n\frac{\partial f^i}{\partial\mu^j}(x,s)d\mu^j|_x\Big)_{i=1}^r=\{\theta\in T^*_xM:\langle\theta,X\rangle=0,\ \forall X\in\Tc_{(x,s)}\},\quad\forall\ u\in M,\ s\in N.
\end{equation}
\end{conv}

\begin{remark}\label{Rmk::EllipticPara::CoSpan}
When $\tilde\beta>1$ and when $f\in\Co^{\tilde\beta}$, \eqref{Eqn::EllipticPara::CoSpan} coincides with the classical notion $\Tc^\bot=\Span(df,d\nu)$ since $\Span(df,d\nu)=\Span(\frac{\partial f}{\partial\mu}d\mu+\frac{\partial f}{\partial\nu}d\nu,d\nu)=\Span(\frac{\partial f}{\partial\mu}d\mu,d\nu)$. 

The convention is used specifically for the case $\tilde\beta\le1$, since when $f\in\Co^{\alpha+1,\tilde\beta}$, the collection of differentials $df$ is no longer pointwise defined, while $\frac{\partial f}{\partial\mu}\cdot d\mu\equiv df \pmod{d\nu}$ and the left hand side is still pointwise.
\end{remark}

We now begin the proof.
On the space $\R^r_\sigma\times\C^m_\zeta$, we write $\xi=\re\zeta$ and $\eta=\im\zeta$. Recall the notation $\beta^{\sim\alpha}$ and $\beta^{\wedge\alpha}$ in \eqref{Eqn::EllipticPara::Betas}. 

\bigskip
\noindent\textit{Proof of Theorem \ref{KeyThm::EllipticPara}.} By Lemma \ref{Lem::EllipticPara::BSatAna} and Proposition \ref{Prop::EllipticPara::AnalyticPDE}, along with a scaling, we can find a $\eps_0=\eps_0(r,m,q,\alpha)>0$ such that,
\begin{enumerate}[parsep=-0.3ex,label=(B.1)]
    \item\label{Item::EllipticPara::PfKey::B} If $B\in\Co^\alpha L^\infty\cap\Co^{-1}\Co^{\beta^{\sim\alpha}}(\B^{r+2m}_{\sigma,\eta,\xi},\B^q_s;\C^{(r+m)\times m})$ satisfies \eqref{Eqn::EllipticPara::KeyEqnB} with 
    \begin{equation}\label{Eqn::EllipticPara::PfKey::BddB}
        \|B\|_{\Co^\alpha L^\infty(\B^{r+2m},\B^q)}<\eps_0,
    \end{equation}then $B(\sigma,\frac{\xi+i\eta}2,s)$ has a   $(\sigma,\xi,\eta)$-holomorphic extension $\Bf(\boldsymbol\sigma,\frac12\boldsymbol\xi,\frac12\boldsymbol\eta,s)\in\Co^\infty_\loc\Co^{\beta^{\sim\alpha}}(\Hb^{r+2m}_{\boldsymbol\sigma,\boldsymbol\xi,\boldsymbol\eta},\B^q_s;\C^{(r+m)\times m})$.
\end{enumerate}
For such $\eps_0>0$, we can find a $\delta_0>0$ such that the results of Proposition \ref{Prop::EllipticPara::ExistPDE} are satisfied with $\eps=\eps_0$, $\delta=\delta_0$.

For such $\delta_0>0$, by Lemma \ref{Lem::EllipticPara::IniNorm} there is a  smooth coordinate chart $L:=((L',L''),L'''):U_0\times V_0\subseteq\Mf\times\Nf\xrightarrow{\sim}\B^{r+2m}_{\tau,w}\times\B^q_s$ near $(u_0,v_0)$ such that,
\begin{enumerate}[parsep=-0.3ex,label=(L.\arabic*)]
    \item\label{Item::EllipticPara::PfKey::LBasic} $(L',L''):U_0\subseteq\Mf\xrightarrow{\sim}\B^{r+2m}_{\tau,w}$ and $L''':V_0\subseteq\Nf\xrightarrow{\sim}\B^q_s$ are both smooth coordinate charts, and $L(u_0,v_0)=(0,0)$.
    \item $L_*\Se$ has a $\Co^{\alpha,\beta}_{(\tau,w),s}$-local basis with form \eqref{Eqn::EllipticPara::DefofA}, such that:
    \begin{itemize}[nolistsep]
        \item When $\beta\in\R_+$, \eqref{Eqn::EllipticPara::IniNorm::A} is satisfied with $\delta=\delta_0$, i.e. $\|A\|_{\Co^{\alpha,\beta}(\B^{r+2m},\B^q;\C^{(r+m)\times m})}<\delta_0$.
        \item When $\beta\in\R_+^-(=\{\gamma-:\gamma\in\R_+\})$, $\|A\|_{\Co^{\alpha}L^\infty(\B^{r+2m},\B^q;\C^{(r+m)\times m})}<\delta_0$.
    \end{itemize}
\end{enumerate}

By Corollary \ref{Cor::EllipticPara::EstofH} if $\beta\in\R_+$ and $\alpha\neq1$, and Proposition \ref{Prop::EllipticPara::ImprovedEstofH} if $\beta\in\R_+^-$ or $\alpha=1$, we can find a $\Co^{\alpha+1,\beta^{\sim\alpha}}$-map $H:\B^{r+2m}_{\tau,w}\times\B^q_s\to\R^r_\sigma\times\C^m_\zeta\times\B^q_s$ that has the form $H(\tau,w,s)=(\tau,H''(\tau,w,s),s)$ such that
\begin{enumerate}[parsep=-0.3ex,label=(H.\arabic*)]
    \item\label{Item::EllipticPara::PfKey::HDiffeo} $H:\B^{r+2m}_{\tau,w}\times\B^q_s\to\B^{r+2m}_{\sigma,\zeta}\times\B^q_s$ is homeomorphism, and $H,H^\Inv\in\Co^{\alpha+1,\beta}(\B^{r+2m},\B^q;\R^r\times\C^m)$.
    \item\label{Item::EllipticPara::PfKey::HTmp} $H(\frac14\B^{r+2m}\times\frac14\B^q)\subseteq\frac12\B^{r+2m}\times\frac12\B^q$.
    \item\label{Item::EllipticPara::PfKey::HReg}$H''\in\Co^{\alpha+1,\beta^{\sim\alpha}}(\B^{r+2m},\B^q;\C^m)$ and $\nabla_{\tau,w}H''\in\Co^{\alpha,\beta^{\wedge\alpha}}(\B^{r+2m},\B^q;\C^{(r+2m)\times m})$.
    \item\label{Item::EllipticPara::PfKey::HPDE} $H''$ solves \eqref{Eqn::EllipticPara::ExistenceH}.
    \item\label{Item::EllipticPara::PfKey::BddLambda} When $\beta\in\R_+$ and $\alpha\neq1$, \eqref{Eqn::EllipticPara::EstofH::BddLambda} is satisfied with $\eps=\eps_0$.
    
    When $\beta\in\R_+^-$ or $\alpha=1$, \eqref{Eqn::EllipticPara::ImprovedEstofH::BddLambda} is satisfied with $\eps=\eps_0$.
\end{enumerate}

Here \ref{Item::EllipticPara::PfKey::HDiffeo} shows that $(H\circ L)_*\Se$ is defined on $\B^{r+2m}_{\sigma,\zeta}\times\B^q_s$.

By Lemma \ref{Lem::EllipticPara::NewGen} \ref{Item::EllipticPara::NewGen::Eqn},  \ref{Item::EllipticPara::PfKey::HPDE} implies that $B$ solves \eqref{Eqn::EllipticPara::KeyEqnB}. By comparing \eqref{Eqn::EllipticPara::matrixfun} and \eqref{Eqn::EllipticPara::DefofB} we see that $B=\Lambda[A,H'']\circ H^\Inv$, so \ref{Item::EllipticPara::PfKey::BddLambda} implies \eqref{Eqn::EllipticPara::PfKey::BddB} with $\eps=\eps_0$. Thus we get the result \ref{Item::EllipticPara::PfKey::B}: a matrix-valued function $\Bf\in\Co^\infty_\loc\Co^{\beta^{\sim\alpha}}(\frac12\Hb^{r+2m},\frac12\B^q;\C^{(r+m)\times m})$.

Applying Lemma \ref{Lem::EllipticPara::FinalModf} on $H\circ L(u_0,v_0)= H(0,0)\in\frac12\B^{r+2m}\times\frac12\B^q$, we can find a neighborhood $U_1\times V_1\subseteq \frac12\B^{r+2m}\times\frac12\B^q$ of $H(0,0)$ and a $G'':U_1\times V_1\to\C^m$ such that

\begin{enumerate}[parsep=-0.3ex,label=(G.\arabic*)]
    \item\label{Item::EllipticPara::PfKey::GReg} $G''\in\Co^\infty\Co^{\beta^{\sim\alpha}}(U_1,V_1;\C^m)$ and $G''(H(0,0))=0$. In particular $G''\in\Co^{\infty,\beta^{\sim\alpha}}_{(\sigma,\zeta),s}$ and $\nabla_{\sigma,\zeta}G''\in\Co^{\infty,\beta^{\sim\alpha}}_{(\sigma,\zeta),s}$.
    \item\label{Item::EllipticPara::PfKey::GSpan} For each $s\in V_1$, $\Span d(G''(\cdot,s))=\Span(T(\cdot,s),Z(\cdot,s))|_{U_1}^\bot=(H\circ L)_*\Se|_{U_1\times \{s\}}^\bot$. Equivalently we have $\Span(dG'',ds)=(H\circ L)_*\Se|_{U_1\times V_1}^\bot $ in the sense of Convention \ref{Conv::EllipticPara::CoSpan}.
\end{enumerate}

Therefore, $(\Span d(G''\circ H),ds)|_{H^\Inv(U_1\times V_1)}=L_*\Se|_{H^\Inv(U_1\times V_1)}$ in the sense of Convention \ref{Conv::EllipticPara::CoSpan}, which means
\begin{equation}\label{Eqn::EllipticPara::PfKey::Span}
    \Span(d(G''\circ H\circ L),dL''')|_{(H\circ L)^\Inv(U_1\times V_1)}=\Se|_{(H\circ L)^\Inv(U_1\times V_1)}^\bot\text{ in the sense of Convention \ref{Conv::EllipticPara::CoSpan}}.
\end{equation}

Next, we claim that 
\begin{enumerate}[parsep=-0.3ex,label=(L.3)]
    \item\label{Item::EllipticPara::PfKey::Claim} $\big(dL',d(\re G''\circ H\circ L(\cdot,v_0)),d(\im G''\circ H\circ L(\cdot,v_0))\big)$ is a collection of $(r+2m)$-real differentials that are linearly independent at $u_0\in\Mf$.
\end{enumerate} 

Indeed, we have $\Span(d(\overline{G''}\circ H),ds)=L_*\bar\Se^\bot$ in the sense of Convention \ref{Conv::EllipticPara::CoSpan}, which means for each $s\in V_1$,
$$\Span(d(\re G''\circ H(\cdot,s)),d(\im G''\circ H(\cdot,s))^\bot=\Span(d(G''\circ H(\cdot,s)),d(\overline{G''}\circ H(\cdot,s)))^\bot=L_*(\Se\cap\bar \Se)|_{H^\Inv(U_1\times\{s\})},$$ has rank $r=\dim\Mf-(r+2m)$. Therefore, $(d(\re G''\circ H(\cdot,s),d(\im G''\circ H(\cdot,s))$ is a linearly independent collection in the domain for all $s\in V_1$.

On the other hand, by \eqref{Eqn::EllipticPara::DefofA} we have 
$L_*\Se=\Span(\Coorvec \tau+A'\Coorvec w,\Coorvec w+A''\Coorvec{\bar w})$, so $L_*(\Se\cap\bar\Se)=\Span(\Coorvec\tau+\frac12A'\Coorvec w+\frac12\bar A'\Coorvec{\bar w}$, which means $L_*(\Se\cap\bar\Se)^\bot\cap\Span d\tau=\{0\}$. We conclude that $(d\tau,d(\re G''\circ H(\cdot,s)),d(\im G''\circ H(\cdot,s))$ is a linearly independent collection in the domain for all $s\in V_1$. Since $dL'=L^*d\tau$ and $dL'''=L^*ds$, taking compositions with $L$ we obtain \ref{Item::EllipticPara::PfKey::Claim}.

\medskip
Now endow $\R^r\times\C^m\times\R^q$ with standard coordinates $(t,z,s)=(t^1,\dots,t^r,z^1,\dots,z^m,s^1,\dots,s^q)$. We take
\begin{equation}\label{Eqn::EllipticPara::PfKey::DefF}
    F=(F',F'',F'''):=(L',\overline{G''}\circ H\circ L,L'''):(H\circ L)^\Inv(U_1\times V_1)\subseteq\Mf\times\Nf\to\R^r_t\times\C^m_z\times\R^q_s.
\end{equation}

By \ref{Item::EllipticPara::PfKey::LBasic}, \ref{Item::EllipticPara::PfKey::HReg}, \ref{Item::EllipticPara::PfKey::GReg} and Lemma \ref{Lem::Hold::CompofMixHold} \ref{Item::Hold::CompofMixHold::Comp} we have 

\begin{enumerate}[parsep=-0.3ex,label=(F.1)]
    \item\label{Item::EllipticPara::PfKey::FReg1} $F'(u_0)=0$, $F''(u_0,v_0)=0$, $F'''(v_0)=0$. And for every product neighborhood $U\times V\subseteq (H\circ L)^\Inv(U_1\times V_1)$ of $(u_0,v_0)$, we have $F'\in\Co^\infty_\loc(U;\R^r)$, $F''\in\Co^{\alpha+1,\beta^{\sim\alpha}}_\loc(U,V;\C^m)$ and $F'''\in\Co^\infty_\loc(V;\R^q)$.
\end{enumerate}
In particular \ref{Item::EllipticPara::F'''} is satisfied, since $L'''$ is already a smooth coordinate chart.

By \ref{Item::EllipticPara::PfKey::FReg1} $(F',F'')\in\Co^{\alpha+1,\beta^{\sim\alpha}}(U,V;\R^r\times\C^m)$ for every $U\times V\subseteq (H\circ L)^\Inv(U_1\times V_1) $. By \ref{Item::EllipticPara::PfKey::Claim}, $(dF',dF''(\cdot,v_0))$ has full rank $r+2m$ in the domain.
Therefore, by Lemma \ref{Lem::Hold::CompofMixHold} \ref{Item::Hold::CompofMixHold::InvFun}, we can find a small enough neighborhood $U\times V\subseteq (H\circ L)^\Inv(U_1\times V_1)$ of $(u_0,v_0)$, such that $F:U\times V\to F(U\times V)\subseteq\R^r_t\times\C^m_z\times\R^q_s$ is homeomorphism, and 
\begin{enumerate}[parsep=-0.3ex,label=(F.2)]
    \item\label{Item::EllipticPara::PfKey::FReg2} $F^\Inv\in\Co^{\alpha+1,\beta}_{(t,z),s}(\Omega'\times\Omega'',\Omega''';\Mf\times\Nf)$ whenever $\Omega'_t\times\Omega''_z\times\Omega'''_s\subseteq F(U\times V)$.
\end{enumerate}

Define $\Phi:=F^\Inv:\Omega\to U\times V$ where $\Omega=\Omega'_s\times\Omega''_z\times\Omega'''_s$ is an arbitrary product neighborhood of $(0,0,0)$ in the assumption of Theorem \ref{KeyThm::EllipticPara}. We are going to show that $F$ and $\Phi$ are as desired.

\medskip
Firstly, by \ref{Item::EllipticPara::PfKey::FReg1}, we get $\Phi''=(L''')^\Inv=(F''')^\Inv$ is a smooth parameterization, which is the result \ref{Item::EllipticPara::Phi''}.

By \ref{Item::EllipticPara::PfKey::FReg1} and \ref{Item::EllipticPara::PfKey::FReg2} we know $F\in\Co^{\beta}$ and $\Phi\in\Co^{\beta}$. Thus \ref{Item::EllipticPara::>1} holds and $F$, $\Phi$ are both homeomorphic to their images respectively. Since $F(\cdot,v)\in\Co^{\alpha+1}$, $\Phi(\cdot,s)\in\Co^{\alpha+1}$ for each $v$ and $s$, we know $(F',F''(\cdot,v))$ is a $\Co^{\alpha+1}$-chart and $\Phi'(\cdot,s)$ is a $\Co^{\alpha+1}$-parameterization, which gives \ref{Item::EllipticPara::FBase} and \ref{Item::EllipticPara::Phi0}.

\medskip
To prove \ref{Item::EllipticPara::F''Reg} and \ref{Item::EllipticPara::PhiReg}, we use $\widehat U:=(L',L'')(U)\subseteq\R^r_\tau\times\C^m_w$, $\widehat V:=L'''(V)\subseteq\R^q_s$, $\widehat F'':=F''\circ L^\Inv$ and $\widehat\Phi:=L\circ\Phi$.
Since $(L',L'')$ and $L'''$ are smooth charts, it is equivalent but more convenient to work on the chart $\widehat F:\widehat U\times\widehat V\to\R^r_t\times\C^m_z\times\R^q_s$ and the parameterization $\widehat\Phi:\Omega\to\widehat U\times\widehat V$:
 it suffices to show the following
\begin{equation}\label{Eqn::EllipticPara::PfKey::AllReg}
    \begin{gathered}
    \widehat F''\in\Co^{\alpha+1,\beta^{\sim\alpha}}_{(\tau,w),s}(\widehat U,\widehat V;\C^m), \quad\nabla_{\tau,w}\widehat F''\in \Co^{\alpha,\beta^{\wedge\alpha}}(\widehat U,\widehat V;\C^{(r+2m)\times m}),
    \\
    \widehat\Phi\in\Co^{\alpha+1,\beta^{\sim\alpha}}_{(t,z),s}(\Omega'\times\Omega'',\Omega''';\widehat U\times \widehat V),
    \quad\nabla_{t,z}\widehat\Phi\in\Co^{\alpha,\beta^{\wedge\alpha}}(\Omega'\times\Omega'',\Omega''';\C^{(r+2m)\times(r+2m)}).
\end{gathered}
\end{equation}

Here for every vector field $X\in\Co^\infty_\loc(\Mf;T\Mf)$, $(L',L'')_*X$ must be the smooth linear combinations of $\Coorvec\tau,\Coorvec w,\Coorvec{\bar w}$ on $\widehat U$. Thus $\nabla_{\tau,w}\widehat F''\in \Co^{\alpha,\beta^{\wedge\alpha}}_{(\tau,w),s}$ implies $XF''\in\Co^{\alpha,\beta^{\wedge\alpha}}_\loc$.

\medskip\noindent\textit{Proof of \eqref{Eqn::EllipticPara::PfKey::AllReg}}: By \eqref{Eqn::EllipticPara::PfKey::DefF} we have $\widehat F''(\tau,w,s)=\overline{G''}(H(\tau,w,s))=\overline{G''}(\tau,H''(\tau,w,s),s)$.
By \ref{Item::EllipticPara::PfKey::HReg} and \ref{Item::EllipticPara::PfKey::GReg} we have $G''\in\Co^{\infty,\beta^{\sim\alpha}}_{(\sigma,\zeta),s}$, $H\in\Co^{\alpha+1,\beta^{\sim\alpha}}_{(\tau,w),s}$, $\nabla_{\sigma,\zeta}G''\in\Co^{\infty,\beta^{\sim\alpha}}_{(\sigma,\zeta),s}$ and $\nabla_{\tau,w}H\in\Co^{\alpha,\beta^{\wedge\alpha}}_{(\tau,w),s}$. Therefore, by Lemma \ref{Lem::Hold::CompofMixHold} \ref{Item::Hold::CompofMixHold::Comp} we get $\widehat F''\in \Co^{\alpha+1,\beta^{\sim\alpha}}_{(\tau,w),s}$ and $\nabla_{\tau,w}\widehat F''=((\nabla_{\sigma,\zeta}\overline{G''})\circ H)\cdot\nabla_{\tau,w}H\in\Co^{\alpha+1,\beta^{\sim\alpha}}_{(\tau,w),s}\cdot \Co^{\alpha,\beta^{\wedge\alpha}}_{(\tau,w),s}\subseteq \Co^{\alpha,\beta^{\wedge\alpha}}_{(\tau,w),s}$.

Since $F\circ (L^\Inv)(\tau,w,s)=(\tau,\widehat F''(\tau,w,s),s)$, we see that $\widehat\Phi$ has the form
\begin{equation*}
    \widehat\Phi(t,z,s)=(t,\widehat\Psi(t,z,s),s),\quad\text{where }\widehat\Psi:\Omega'_t\times\Omega''_z\times\Omega'''_s\to\C^m_w,\quad\widehat\Psi(t,z,s):=\widehat F''(t,\cdot,s)^\Inv(z).
\end{equation*}

Since $\widehat F''\in\Co^{\alpha+1,\beta^{\sim\alpha}}_{(\tau,w),s}$, applying Lemma \ref{Lem::Hold::CompofMixHold} \ref{Item::Hold::CompofMixHold::InvFun} to the map $[(\tau,w),s)\mapsto(\tau,\widehat F''(\tau,w,s))]$ we get $\widehat\Psi\in\Co^{\alpha+1,\beta^{\sim\alpha}}_{(t,z),s}$. 

Using the chain rules on $z=\widehat F''(\widehat\Phi(t,z,s))=\widehat F''(t,\widehat\Psi(t,z,s),s)$ we have (recall $\nabla_w=[\Coorvec w,\Coorvec{\bar w}]^\top$),
\begin{align*}
    I_{2m}&=((\nabla_w\widehat F'')\circ\widehat\Phi)\cdot\nabla_z\widehat\Psi,& 0_{r\times m}&=\nabla_\tau \widehat F''\circ\widehat\Phi+((\nabla_w\widehat F'')\circ\widehat\Phi)\cdot\nabla_t\widehat\Psi;
    \\
    \Longrightarrow\quad \nabla_z\widehat\Psi&=(\nabla_w\widehat F'')^{-1}\circ\widehat\Phi,&\nabla_t\widehat\Psi&=((\nabla_w\widehat F'')^{-1}\cdot\nabla_\tau \widehat F'')\circ\widehat\Phi.
\end{align*}

Since $\nabla_{\tau,w}\widehat F''\in\Co^{\alpha,\beta^{\wedge\alpha}}_{(\tau,w),s}$, by Lemma \ref{Lem::Hold::CompofMixHold} \ref{Item::Hold::CompofMixHold::InvMat} we have $(\nabla_w\widehat F'')^{-1}\in\Co^{\alpha,\beta^{\wedge\alpha}}_{(\tau,w),s}$ thus $(\nabla_w\widehat F'')^{-1}\cdot\nabla_\tau \widehat F''\in\Co^{\alpha,\beta^{\wedge\alpha}}_{(\tau,w),s}$. Applying Lemma \ref{Lem::Hold::CompofMixHold} \ref{Item::Hold::CompofMixHold::Comp} with $\widehat\Phi\in\Co^{\alpha+1,\beta^{\sim\alpha}}_{(t,z),s}$, we get $\nabla_{t,z}\widehat\Psi\in\Co^{\alpha,\beta^{\wedge\alpha}}_{(t,z),s}$. Therefore $\nabla_{t,z}\widehat\Phi\in\Co^{\alpha,\beta^{\wedge\alpha}}_{(t,z),s}$ finishing the proof of \eqref{Eqn::EllipticPara::PfKey::AllReg} and hence the proof of \ref{Item::EllipticPara::F''Reg} and \ref{Item::EllipticPara::PhiReg}.

\medskip
Now by \ref{Item::EllipticPara::PfKey::FReg1} and \ref{Item::EllipticPara::PhiReg} we see that $F\in\Co^{\alpha+1,\beta^{\sim\alpha}}_{u,v}$ and $\nabla_{t,z}\Phi\in\Co^{\alpha,\beta^{\wedge\alpha}}_{(t,z),s}$. By Lemma \ref{Lem::Hold::CompofMixHold} \ref{Item::Hold::CompofMixHold::Comp} we see that $\frac{\partial\Phi}{\partial t^j}\circ F\in\Co^{\alpha,\beta^{\wedge\alpha}}_\loc(U,V;T\Mf)$ and $\frac{\partial\Phi}{\partial z^k}\circ F\in\Co^{\alpha,\beta^{\wedge\alpha}}_\loc(U,V;\C T\Mf)$ for $1\le j\le r$ and $1\le k\le m$. Thus $F^*\Coorvec{t^j}=\frac{\partial\Phi}{\partial t^j}\circ F$, $F^*\Coorvec{z^k}=\frac{\partial\Phi}{\partial z^k}\circ F$ are both well-defined vector fields and have regularity $\Co^\alpha\cap\Co^{\beta^{\wedge\alpha}}=\Co^{\min(\alpha,\beta^{\wedge\alpha})}$.

By definition, $dF'=F^*dt$, $dF'''=F^*ds$,  $dF''=F^*dz$ and $d\overline{F''}=F^*d\bar z$. So by \eqref{Eqn::EllipticPara::PfKey::Span} and \eqref{Eqn::EllipticPara::PfKey::DefF}, in the sense of Convention \ref{Conv::EllipticPara::CoSpan},
\begin{equation*}
    \textstyle\Span(F^*\Coorvec t,F^*\Coorvec z)=\Span(F^*d\bar z,F^*ds)^\bot=\Span(d\overline{F''},dF''')^\bot=\Se|_{U\times V}.
\end{equation*}

This proves \ref{Item::EllipticPara::Span}. Immediately \ref{Item::EllipticPara::PhiSpan} follows.
Now the whole proof is complete.\qed

\section{Improving the Regularity of Rough 1-forms}\label{Section::Rough1Form}



Let $\alpha>0$ and $\beta\in [\alpha,\alpha+1]$.  Suppose $\lambda^1,\ldots, \lambda^{n}$
are $\Co^{\alpha}$ $1$-forms on an open set $U\subseteq \R^n$ which span
the cotangent space at every point in $U$.  If we know that $d\lambda^j\in \Co^{\beta-1}_{\loc}$ for each $j$, it is not necessarily true that $\lambda^j\in \Co^{\beta}_{\loc}$. However, one can always
change coordinates so that the forms are in $\Co^{\beta}_{\loc}$
(see also Corollary \ref{Cor::Rough1Form::BasicCaseforMainThm}).

The next result is a special case of this idea, where we present an intital setting where we may find a $\Co^{\alpha+1}$-diffemorphism such that $F_{*}\lambda^j\in \Co^{\beta}_{\loc}$.

\begin{keythm}[Improving regularity of rough 1-forms]\label{KeyThm::Rough1Form}
Let $\alpha>0$ and $\beta\in [\alpha,\alpha+1]$. Let $x=(x^1,\dots,x^n)$ and $y=(y^1,\dots,y^n)$ be two coordinate systems for $\R^n$.
There exists $c=c(n,\alpha,\beta)>0$ such that the following holds.

Suppose $\lambda^i\in \Co^{\alpha}(\B^n; T^{*}\B^n)$, $i=1,\dots,n$ are 1-forms on $\B^n$ such that $\supp (\lambda^i-dx^i)\subsetneq\frac{1}{2}\B^n$ for each $i$, and
\begin{equation}\label{Eqn::Rough1Form::Assumption}
    \sum_{i=1}^n(\|\lambda^i-dx^i\|_{\Co^{\alpha}(\B^n;T^{*}\B^n)}+\|d\lambda^i\|_{\Co^{\beta-1}(\B^n;\wedge^2 T^{*}\B^n)})\le c.
\end{equation}

Then, there exists a $\Co^{\alpha+1}$-diffeomorphism $F:\B^n_x\xrightarrow{\sim}\B^n_y$, such that $B^n(F(0),\frac16)\subseteq F(\frac13\B^n)\cap\frac34\B^n$, $F_*\lambda^i\in\Co^\beta(\B^n;T^*\B^n)$ for $i=1,\dots,n$, and moreover
\begin{equation}\label{Eqn::Rough1Form::Conclusion}
\begin{aligned}
    &\|F-\id\|_{\Co^{\alpha+1}(\B^n;\R^n)}+\sum_{i=1}^n\|F_*\lambda^i-dy^i\|_{\Co^\beta\mleft(\B^n;T^*\B^n\mright)}
    \le c^{-1}\sum_{i=1}^n\left(\|\lambda^i-dx^i\|_{\Co^\alpha\mleft(\B^n;T^*\B^n\mright)}+\|d\lambda^i\|_{\Co^{\beta-1}\mleft(\B^n;\wedge^2T^*\B^n\mright)}\right).
\end{aligned}
\end{equation}
\end{keythm}

\begin{remark}\label{Rmk::Rough1From::RmkofThm}
Notice that if $\tilde c(n,\alpha,\beta)>0$ is a small constant that satisfies the results in Theorem \ref{KeyThm::Rough1Form} with $c=\tilde c$, then any constant $c<\tilde c$ would satisfies the results as well.
Therefore by shrinking $c$ if possible, we can make the left hand side of \eqref{Eqn::Rough1Form::Assumption} bounded by, say $\frac1{12}$.
\end{remark}

\subsection{Outline of the proof: the dual Malgrange method}\label{Section::Rough1FormOV}
The proof of Theorem \ref{KeyThm::Rough1Form}  is inspired by Malgrange's proof
of the Newlander-Nirenberg Theorem \cite{Malgrange}.

Let $x=(x^1,\dots,x^n)$ and $y=(y^1,\dots,y^n)$ be two coordinate systems on the unit ball $\B^n\subset\R^n$. In this section, we write 1-forms $\lambda^1,\dots,\lambda^n$, $\eta^1,\dots,\eta^n$ as
\begin{equation}\label{Eqn::Rough1Form::LambdaEta}
    \lambda^i=dx^i+\sum_{j=1}^na_j^idx^j,\quad\eta^i=dy^i+\sum_{j=1}^nb_j^idy^j,\quad i=1,\dots,n,
\end{equation}
and define coefficient matrices
\begin{equation}\label{Eqn::Rough1Form::LambdaEtaMatrix}
    A:=(a_j^i)_{n\times n}=\begin{pmatrix}a_1^1&\cdots&a_n^1\\\vdots&\ddots&\vdots\\a_1^n&\cdots&a_n^n\end{pmatrix},\quad B:=(b_j^i)_{n\times n}=\begin{pmatrix}b_1^1&\cdots&b_n^1\\\vdots&\ddots&\vdots\\b_1^n&\cdots&b_n^n\end{pmatrix}.
\end{equation}

In this section, $\lambda^1,\ldots, \lambda^n$ are given $\Co^{\alpha}$
$1$-forms on $\B^n\subset \R^n$ which span the cotangent space at every point.
And $\eta^i:=F_{*}\lambda^i$ are the push-forward $1$-forms by the unknown
$\Co^{\alpha+1}$-diffeomorphism $F:\B^n\xrightarrow{\sim}\B^n$,
which we are solving for.  Thus, $\eta^1,\ldots, \eta^n$ are also $\Co^{\alpha}$ $1$-forms defined on $\B^n$ which span the cotangent
space at every point.


As in Malgrange's work \cite{Malgrange},
the main idea is to choose $F$ so that the matrix $B$ satisfies a nonlinear
elliptic PDE.  That $\eta=F_{*}\lambda\in \Co^{\beta}$ will follow
from the classical interior regularity of elliptic PDEs.  We will show
such an $F$ exists by showing that it suffices for $F$ to satisfy
a different nonlinear elliptic PDE, whose solution is guaranteed by
classical elliptic theory.

    

Given collections $(\lambda^1,\dots,\lambda^n)$ and $(\eta^1,\dots,\eta^n)$ of 1-forms on $\B^1$, as above, that both span their co-tangent spaces at every point, we define Riemannian metrics $g$ and $h$ by
\begin{equation}\label{Eqn::Rough1Form::Riemannianmetric1}
\begin{split}
    &g=\sum_{i,j=1}^ng_{ij}dx^idx^j:=\sum_{i,j,k=1}^n(\delta_i^k+a_i^k)(\delta_j^k+a_j^k)dx^idx^j,\\ &h=\sum_{i,j=1}^nh_{ij}dy^idy^j:=\sum_{i,j,k=1}^n(\delta_i^k+b_i^k)(\delta_j^k+b_j^k)dy^idy^j.
\end{split}
\end{equation}
Here, $\delta_i^j$, $\delta_{ij}$, $\delta^{ij}$ are the Kronecker delta functions:
\begin{equation}\label{Eqn::Rough1Form::Krnoecker}
    \delta_i^j=\delta_{ij}=\delta^{ij}=\begin{cases}1,&i=j,\\0,&i\neq j.\end{cases}
\end{equation}

We use the following notations from classical Riemannian differential geometry:
\begin{equation}\label{Eqn::Rough1Form::Riemannianmetric2}
\begin{split}
    &g^{ij}:=g(dx^i,dx^j),\quad \sqrt{\det g}:=\left|\frac{\lambda^1\wedge\dots\wedge\lambda^n}{dx^1\wedge\dots\wedge dx^n}\right|,
    \\&h^{ij}:=h(dy^i,dy^j),\quad \sqrt{\det h}:=\left|\frac{\eta^1\wedge\dots\wedge\eta^n}{dy^1\wedge\dots\wedge dy^n}\right|.
\end{split}
\end{equation}

\begin{remark}\label{Rmk::Rough1Form::RmkRiemMetric}
\begin{enumerate}[label=(\alph*)]
    \item We can write $h=\sum_{i=1}^n\eta^i\cdot\eta^i$. It is non-degenerate since $\eta^1,\dots,\eta^n$ span the cotangent space at every point. Moreover $\eta^1,\dots,\eta^n$ form an orthogonal basis with respect to this metric $h$. Similar remarks hold for $g=\sum_{i=1}^n\lambda^i\cdot\lambda^i$.
    
    \item\label{Item::Rough1Form::RmkRiemMetric::RationalFunction} Using matrix notations in \eqref{Eqn::Rough1Form::LambdaEta} we have $(h_{ij})_{n\times n}=(I+B)^\top(I+B)$ and so we know $(h^{ij})_{n\times n}=(h_{ij})^{-1}=(I+B)^{-1}\left((I+B)^{-1}\right)^\top$ and $\sqrt{\det h}=\det(I+B)$. Similarly, we  have $(g_{ij})=(I+A)^\top(I+A)=(g^{ij})^{-1}$ and $\sqrt{\det g}=\det(I+A)$. More importantly, we have the following lemma.
\end{enumerate}
\end{remark}
\begin{lem}\label{Lem::Rough1Form::TaylorExpansionofRiemMetric}
Let $B$, $h^{ij}$ and $\sqrt{\det h}$ be as above. Then $h^{ij}$ and $\sqrt{\det h}$ are rational functions of the components of $B$. Moreover for every $\gamma>0$, there is a $c_{n,\gamma}>0$ such that
\begin{equation*}
    h^{ij},\sqrt{\det h}:\mleft\{B\in\Co^\gamma(\B^n;\R^{n\times n}):\|B\|_{\Co^\gamma}<c_{n,\gamma}\mright\}\to \Co^\gamma(\B^n),\quad 1\le i,j\le n,
\end{equation*}
are norm continuous maps, with 
$$\sum_{i,j=1}^n\|h^{ij}-\delta^{ij}\|_{\Co^\gamma(\B^n)}\leq c_{n,\gamma}^{-1}\|B\|_{\Co^\gamma(\B^n;\R^{n\times n})},\quad\|\sqrt{\det h}-1\|_{\Co^\gamma(\B^n)}\leq c_{n,\gamma}^{-1}\|B\|_{\Co^\gamma(\B^n;\R^{n\times n})}.$$
\end{lem}
\begin{remark}\label{Rmk::Rough1Form::TaylorExpansionofRiemMetric}
The same results hold for $g^{ij}$ and $\sqrt{\det g}$. Namely,  $\|g^{ij}-\delta^{ij}\|_{\Co^\alpha}+\|\sqrt{\det g}-1\|_{\Co^\alpha}\leq c_{n,\alpha}^{-1} \|A\|_{\Co^\alpha}$ holds with the same constant $c_{n,\alpha}>0$.
\end{remark}
\begin{proof}
By Lemma \ref{Lem::Hold::Product}, the space $\Co^\gamma(\B^n;\R^{n\times n})$ is closed under matrix multiplication. By Remark \ref{Rmk::Rough1Form::RmkRiemMetric} \ref{Item::Rough1Form::RmkRiemMetric::RationalFunction} $\sqrt{\det h}=\det(I+B)$ is a polynomial in the components of $B$, so in particular is a norm continuous function on $\Co^\gamma(\B^n;\R^{n\times n})$. Note that $\det(I+0)=1$ so we have 
$\|\sqrt{\det h}-1\|_{\Co^\gamma(\B^n)}\lesssim_\gamma\|B\|_{\Co^\gamma(\B^n;\R^{n\times n})}$ when the right hand side is small.

By Lemma \ref{Lem::Hold::CramerMixedPara} \ref{Item::Hold::CramerMixedPara::SingHold}, choosing $c_{n,\gamma}<\tilde c(\B^n,\alpha,n)$ where $\tilde c$ is the constant in Lemma \ref{Lem::Hold::CramerMixedPara} \ref{Item::Hold::CramerMixedPara::SingHold}, we see that the map $B\mapsto(I+B)^{-1}$ is $\Co^\gamma$-norm continuous with $\|(I+B)^{-1}-I\|_{\Co^\gamma(\B^n;\R^{n\times n})}\le2\|B\|_{\Co^\gamma(\B^n;\R^{n\times n})}$. 

Thus, in the domain $\|B\|_{\Co^\gamma}<c_{n,\gamma}$, the map $B\mapsto(I+B)^{-1}\left((I+B)^{-1}\right)^\top$ is also continuous and satisfies $\|(I+B)^{-1}\left((I+B)^{-1}\right)^\top-I\|_{\Co^\gamma}\lesssim\|(I+B)^\top(I+B)-I\|_{\Co^\gamma}\lesssim\|B\|_{\Co^\gamma}$. It follows with $(h^{ij})_{n\times n}=(h_{ij})^{-1}_{n\times n}=(I+B)^{-1}\left((I+B)^{-1}\right)^\top$, we have that $h^{ij}$ are norm continuous and $\|h^{ij}-\delta^{ij}\|_{\Co^\gamma}\lesssim\|B\|_{\Co^\gamma}$.

By possibly shrinking $c_{n,\gamma}$ we get $\|\sqrt{\det h}-1\|_{\Co^\gamma}\le c_{n,\gamma}^{-1}\|B\|_{\Co^\gamma}$ and $\sum_{i,j=1}^n\|h^{ij}-\delta^{ij}\|_{\Co^\gamma}\le c_{n,\gamma}^{-1}\|B\|_{\Co^\gamma}$.
\end{proof}




\begin{conv}\label{Conv::Rough1Form::CoDiff}
Given a Riemannian metric $g$, we use the co-differential $\codiff_g$ as the adjoint of differential with respect to $h$. That is, for any $k$-form $\phi$ and any compactly supported $(k-1)$-form $\psi$, 
\begin{equation*}
    (\codiff_g\phi,\psi)_g:=(\phi,d\psi)_g=\int g(\phi,d\psi)\ d\Vol_g,
\end{equation*}
where $d\Vol_g$ is the Riemannian volume density induced by $g$.  In local coordinates, $d\Vol_g=\sqrt{\det g}\ d\Vol_{\R^n}$, where $d\Vol_{\R^n}$ is the usual Lebesgue density on $\R^n$.
We write $\codiff_{\R^n}=\codiff$ for the usual co-differential with respect to the flat metric $\sum_{i,j=1}^n\delta_{ij}dy^idy^j$ on $\R^n$.
\end{conv}


\begin{lem}\label{Lem::Rough1Form::TransitionBetweenPDEs}
Let $\alpha>0$, and let $\lambda^1,\dots,\lambda^n$ be $\Co^\alpha$ 1-forms defined on $\B^n$ that span the cotangent space at every point, with $(a_i^j)$, $g$, $(g^{ij})$, and $\sqrt{\det g}$ given in \eqref{Eqn::Rough1Form::LambdaEtaMatrix}, \eqref{Eqn::Rough1Form::Riemannianmetric1}, and \eqref{Eqn::Rough1Form::Riemannianmetric2}. 

Suppose $F=\id+R:\B^n_x\to\B^n_y$ is a $\Co^{\alpha+1}$-diffeomorphism that satisfies
\begin{equation}\label{Eqn::Rough1Form::PDEforR}
        \sum_{i,j=1}^n\Coorvec{x^j}\Big(\sqrt{\det g}\cdot g^{ij}\frac{\partial R^k}{\partial x^i}\Big)=\sum_{i,j=1}^n\Coorvec{x^j}\big(\sqrt{\det g}\cdot g^{ij}a_i^k\big),\quad\text{in }\B^n_x,\quad k=1,\dots,n.
\end{equation}
Then for the pushforward 1-forms $\eta^k=F_*\lambda^k$, $k=1,\dots,n$, the coefficients $(b_i^j)$ defined in \eqref{Eqn::Rough1Form::LambdaEtaMatrix} satisfy
\begin{equation}\label{Eqn::Rough1Form::PDEforB}
    \sum_{i,j=1}^n\Coorvec{y^i}\left(\sqrt{\det h}h^{ij}b_j^k\right)=0,\quad\text{in }\B^n_y,\quad k=1,\dots,n.
\end{equation}
\end{lem}
\begin{proof}
    Note that the composition of a $C^1$-function and a $\Co^{\alpha+1}$-diffeomorphism is still $C^1$, and being compactly supported is preserved under homeomorphism, so we have
    \begin{equation}\label{Eqn::Rough1Form::TransitionBetweenPDEs::Proof}
        C_c^1(\B_x^n)=\{v\circ F:v\in C_c^1(\B_y^n)\}.
    \end{equation}
    
    By assumption $R\in\Co^{\alpha+1}(\B^n;\R^n)$, and $a_i^k,g^{ij},\sqrt{\det g}\in\Co^\alpha(\B^n)$, so \eqref{Eqn::Rough1Form::PDEforR} makes sense in $\Co^{\alpha-1}_\loc(\B^n_x)\subsetneq C_c^1(\B^n_x)'$ and the equality
    can be viewed as elements of the dual of $C_c^1(\B^n)$. 
    
    For any $u\in C_c^1(\B^n)$, integrating by parts, we obtain
    \begin{align*}
        0&=\bigg\langle\sum_{i,j=1}^n\Coorvec{x^j}\bigg(\sqrt{\det g}\cdot g^{ij}\Big(\frac{\partial R^k}{\partial x^i}-a_i^k\Big)\bigg),u\bigg\rangle_{\B^n_x}=-\int_{\B^n}\sum_{i,j=1}^ng^{ij}\Big(\frac{\partial R^k}{\partial x^i}-a_i^k\Big)\frac{\partial u}{\partial x^j}\sqrt{\det g}dx
        \\
        &=-\langle dR^k-(\lambda^k-dx^k),du\rangle_{\B^n_x;g}=-\langle dF^k-\lambda^k,du\rangle_{\B_x^n;g}
    \end{align*}
    Here $\langle\cdot,\cdot\rangle_{\B_x^n;g}$ are the dual pairs for linear functionals and test functions induced by $g$. Namely, for $u,v\in C^0(\B^n)$ and $\phi,\psi\in C^0(\B^n;T^*\B^n)$, $\langle u,v\rangle_{\B^n;g}=\int_{\B^n}uv\sqrt{\det g}dx$ and $\langle \phi,\psi\rangle_{\B^n;g}=\int_{\B^n}g(\phi,\psi)\sqrt{\det g}dx$.

    Using \eqref{Eqn::Rough1Form::PDEforR}, we get $\langle dF^k-\lambda^k,du\rangle_{\B^n_x;g}=0$ for all $u\in C_c^1(\B^n)$. By \eqref{Eqn::Rough1Form::TransitionBetweenPDEs::Proof} we have $\langle dF^k-\lambda^k,d(v\circ F)\rangle_{\B^n_x;g}=0$ for all $v\in C_c^1(\B^n_y)$.
    
    Note that $F_*(g(\phi,\psi))=(F_*g)(F_*\phi,F_*\psi)=h(F_*\phi,F_*\psi)$ for all $\phi,\psi\in C^0_\loc(\B^n;T^*\B^n)$, and $F_*\sqrt{\det g}=\sqrt{\det F_*g}=\sqrt{\det h}$, so 
    we have, for every $v\in C_c^1(\B^n_y)$,
    \begin{equation*}
        0=\langle dF^k-\lambda^k,d(v\circ F)\rangle_{\B^n_x;g}=\langle F_*(dF^k-\lambda^k),F_*F^*dv\rangle_{\B^n_y;h}=\langle dy^k-\eta^k,dv\rangle_{\B^n_y;h}=\int_{\B^n}\sum_{i,j=1}^nh^{ij}b_j^k\frac{\partial v}{\partial y^i}\sqrt{\det h}dy.
    \end{equation*}
    Integrating by parts, we obtain \eqref{Eqn::Rough1Form::PDEforB}.
\end{proof}

We will choose a coordinate chart $F$ so that \eqref{Eqn::Rough1Form::PDEforR}
is satisfied, and therefore \eqref{Eqn::Rough1Form::PDEforB}
will be satisfied as well.

To  prove Theorem \ref{KeyThm::Rough1Form}  we will
prove the following:
\begin{itemize}
    \item There exists a $R\in\Co^{\alpha+1}(\B^n;\R^n)$ satisfying \eqref{Eqn::Rough1Form::PDEforR} with boundary condition $R\big|_{\partial\B^n}=0$. Moreover,
    we can choose $R$ with $\|R\|_{\Co^{\alpha+1}}\lesssim_{\alpha,\beta}\|A\|_{\Co^\alpha}$.
    Thus, by taking $c>0$ small, we may take $\|R\|_{\Co^{\alpha+1}}$ small.
    
    \item When $\|R\|_{\Co^{\alpha+1}}$ is small, $F=\id+R$ is a $\Co^{\alpha+1}$-diffeomorphism of $\B^n$.
    And under the assumption $\supp A\subsetneq\frac12\B^n$, we have $\|B\|_{\Co^\beta(\partial\B^n)}\lesssim_{\alpha,\beta}\|A\|_{\Co^\alpha}$ and $\|d\eta\|_{\Co^{\beta-1}}\lesssim_{\alpha,\beta}\|d\lambda\|_{\Co^{\beta-1}}$.
    In particular, by taking $c>0$ small, we may take 
    $\|B\|_{\Co^\beta(\partial\B^n)}+\|d\eta\|_{\Co^{\beta-1}} $ small.
    
    \item Using that $B\in\Co^\alpha(\partial\B^n;\R^{n\times n})$ satisfies \eqref{Eqn::Rough1Form::PDEforB}, if $\|B\|_{\Co^\beta(\partial\B^n)}+\|d\eta\|_{\Co^{\beta-1}}$ is small, we will show $B\in\Co^\beta(\B^n;\R^{n\times n})$ and  $\|B\|_{\Co^\beta}\lesssim_{\alpha,\beta}\|A\|_{\Co^\alpha}$.
\end{itemize}

\subsection{The existence proposition}
In this section, we show that there exists a $\Co^{\alpha+1}$-diffeomorphism $F=\id+R$ solving \eqref{Eqn::Rough1Form::PDEforR} and which satisfies
good estimates.

\begin{prop}\label{Prop::Rough1Form::ExistPDE}
Let $\alpha>0$ and let $\beta\in[\alpha,\alpha+1]$. There is a $c_1=c_1(n,\alpha,\beta)\in(0,1)$ such that, if $A=(a_j^k)_{n\times n}:\B^n\to\R^{n\times n}$ satisfies
\begin{itemize}[parsep=-0.3ex]
    \item[-] $A\in\Co^\alpha_c(\frac12\B^n;\R^{n\times n})$ and  $\|A\|_{\Co^\alpha}<c_1$,
\end{itemize}
then the matrix $(I+A)(x)$ is invertible for every $x\in\B^n$, and there is a  $\Co^{\alpha+1}$-map $ F=\id+R$ on $\B^n$ such that
\begin{enumerate}[parsep=-0.3ex,label=(\roman*)]
    \item\label{Item::Rough1Form::ExistPDE::1} $R$ solves the equation
    \begin{equation}\tag{\ref{Eqn::Rough1Form::PDEforR}}
        \sum_{i,j=1}^n\Coorvec{x^j}\Big(\sqrt{\det g}\cdot g^{ij}\frac{\partial R^k}{\partial x^i}\Big)=\sum_{i,j=1}^n\Coorvec{x^j}\big(\sqrt{\det g}\cdot g^{ij}a_i^k\big),\quad\text{in }\B^n_x,\quad k=1,\dots,n.
    \end{equation}
     with boundary condition $R\big|_{\partial\B^n}=0$, and we have
    \begin{equation}\label{Eqn::Rough1Form::ExistPDEQuantControl1}
        \textstyle\|R\|_{\Co^{\alpha+1}(\B^n; \R^n)}+\|\nabla R\|_{\Co^{\beta}(\B^n\backslash\frac34\B^n;\R^n)}\le c_1^{-1}\|A\|_{\Co^\alpha}.
    \end{equation}
    \item\label{Item::Rough1Form::ExistPDE::1.5} $F:\B^n_x\to\B^n_y$ is a $\Co^{\alpha+1}$-diffeomorphism such that $B^n(F(0),\frac16)\subseteq F(\frac13\B^n)\cap\frac34\B^n $.  
    \item\label{Item::Rough1Form::ExistPDE::2}Let $\Phi=F^\Inv:\B^n_y\to\B^n_x$ be its
    inverse map,  then 
    \begin{equation}\label{Eqn::Rough1Form::ExistPDEQuantControl2}
        \|\nabla \Phi-I\|_{\Co^{\alpha}(\B^n;\R^{n\times n})}+\|\nabla \Phi-I\|_{\Co^{\beta}(\partial\B^n;\R^{n\times n})}\le c_1^{-1}\|A\|_{\Co^\alpha}.
    \end{equation}
    
    In particular $\|\Phi\|_{\Co^{\alpha+1}(\B^n;\R^n)}\le c_1^{-1}$.
\end{enumerate}
\end{prop}

\begin{remark}
The map $F$ in Proposition \ref{Prop::Rough1Form::ExistPDE} is uniquely determined by $A$. This is due to the well-posedness of the Dirichlet problem for the second order elliptic equations, since $R$ satisfies \eqref{Eqn::Rough1Form::PDEforR} with $R\big|_{\partial\B^n}=0$.
\end{remark}

\begin{remark}\label{Rmk::Rough1Form::ExistPDE:FisNeeded}
As we will see in the proof of Theorem \ref{KeyThm::Rough1Form}, the map $F$
from Proposition \ref{Prop::Rough1Form::ExistPDE}
is the map of the same name in Theorem \ref{KeyThm::Rough1Form}.
\end{remark}

\begin{proof}
We let $c_1$ be a small constant which may change from line to line.
     Note that if \eqref{Eqn::Rough1Form::ExistPDEQuantControl1} and \eqref{Eqn::Rough1Form::ExistPDEQuantControl2} are valid for some $\tilde c_1$, then they are also valid for any $0<c_1\le\tilde c_1$. 
    
    First pick $c_1<\tilde c(\B^n,\alpha,n)$ where $\tilde c$ is the constant in Lemma \ref{Lem::Hold::CramerMixedPara} \ref{Item::Hold::CramerMixedPara::SingHold}. By Lemma \ref{Lem::Hold::CramerMixedPara} \ref{Item::Hold::CramerMixedPara::SingHold}, the assumption $\|A\|_{\Co^\alpha}<c_1(<\tilde c_{\B^n,\alpha,n})$ implies that $I+A$ is invertible at every point and
    $(I+A)^{-1}\in \Co^{\alpha}(\B^n; \R^{n\times n})$. Therefore,  $g$ given in \eqref{Eqn::Rough1Form::Riemannianmetric1} is indeed a $\Co^\alpha$-Riemannian metric.

    By the assumption $\supp A\subsetneq\frac12\B^n$, we have $\sqrt{\det g}\cdot g^{ij}\big|_{\B^n\backslash\frac12\B^n}=\delta^{ij}$. {Since $(\sqrt{\det g}g^{ij}(x))_{n\times n}$ is  an invertible matrix for $x\in\frac12\B^n$}, the second order operator $\sum_{i,j=1}^n\partial_{x^j}(\sqrt{\det g}\cdot g^{ij}\partial_{x^i})$ is uniformly elliptic on $\B^n$. 
    Classical existence theorems (for example, \cite[Theorem 8.3]{GilbargTrudinger}) show that for each $k=1,\dots,n$ there exists\footnote{Here $H^1(\B^n)$ stands for the classical $L^2$-Sobolev space of order $1$, and $H^{-1}(\B^n)=H_0^1(\B^n)^*$ is the $L^2$-Sobolev space of order $-1$.} a $R^k\in H^1(\B^n)$ that satisfies \eqref{Eqn::Rough1Form::PDEforR} with Dirichlet boundary condition $R^k\big|_{\partial\B^n}=0$, since $\sum_{i,j=1}^n\Coorvec{x^j}\left(\sqrt{\det g}g^{ij}a_i^k\right)\in\Co^{\alpha-1}(\B^n)\subset H^{-1}(\B^n)$. 
    By a classical regularity estimate (see \cite[Theorem 8.34]{GilbargTrudinger} or  \cite[Theorem 15]{DirichletBoundedness}), we know $R^k\in\Co^{\alpha+1}(\B^n; \R^n)$.
    
    To show $\|R\|_{\Co^{\alpha+1}}\lesssim\|A\|_{\Co^\alpha}$, we  write \eqref{Eqn::Rough1Form::PDEforR} as
    \begin{equation}\label{Eqn::Rough1Form::ExistPDE::Proof1}
        \Delta R^k=\sum_{i,j=1}^n\Coorvec{x^j}\left(\left(\delta^{ij}-\sqrt{\det g}g^{ij}\right)\frac{\partial R^k}{\partial x^i}\right)+\sum_{i,j=1}^n\Coorvec{x^j}\left(\sqrt{\det g}g^{ij}a_i^k\right),\quad\text{in }\B^n_x,\quad k=1,\dots,n.
    \end{equation}
    
    By  Remark \ref{Rmk::Rough1Form::TaylorExpansionofRiemMetric} (see also Lemma \ref{Lem::Rough1Form::TaylorExpansionofRiemMetric}), we see that $\|\delta^{ij}-\sqrt{\det g}g^{ij}\|_{\Co^\alpha(\B^n)}\lesssim_\alpha\|A\|_{\Co^\alpha}$.
    Therefore
    \begin{equation}\label{Eqn::Rough1Form::EstimateLapR}
        \|\Delta R\|_{\Co^{\alpha-1}}\lesssim_\alpha\sum_{i,j,k=1}^n\left(\|\delta^{ij}-\sqrt{\det g}g^{ij}\|_{\Co^\alpha}\|\partial_{x^i}R^k\|_{\Co^\alpha}+\|\sqrt{\det g}g^{ij}a_i^k\|_{\Co^\alpha}\right)\lesssim
        \|A\|_{\Co^\alpha}\|R\|_{\Co^{\alpha+1}}+\|A\|_{\Co^\alpha},
    \end{equation}
    where the implicit constants depend only on $n$ and $\alpha$ but not $A$ or $R$.
    
    The assumption $R\big|_{\partial\B^n}=0$ implies that $R=\Pc_0^{\B^n}\Delta R$, where $\Pc_0^{\B^n}$ is the zero Dirichlet boundary solution operator given in Notation \ref{Note::Hold::DiriSol}. Since, by Lemma \ref{Lem::Hold::DiriSol}, $\Pc_0^{\B^n}:\Co^{\alpha-1}(\B^n)\to\Co^{\alpha+1}(\B^n)$ is bounded, \eqref{Eqn::Rough1Form::EstimateLapR} implies
    \begin{equation}\label{Eqn::Rough1Form::ExistPDE::Tmp1}
        \| R\|_{\Co^{\alpha+1}(\B^n;\R^n)}\le\tilde C_1\|A\|_{\Co^\alpha(\B^n;\R^{n\times n})}\|R\|_{\Co^{\alpha+1}(\B^n;\R^n)}+\tilde C_1\|A\|_{\Co^\alpha(\B^n;\R^{n\times n})},
    \end{equation}
    where $\tilde C_1=\tilde C_1(n,\alpha)>1$ is a constant  depending only on $n$ and $\alpha$ but not $A$ or $R$.
    
    Choosing $c_1$ small enough so that $c_1\tilde C_1\le\frac13$, then we get $\|R\|_{\Co^{\alpha+1}}\le\frac13\|R\|_{\Co^{\alpha+1}}+\tilde C_1\|A\|_{\Co^\alpha}$ when $A$ satisfies the assumption $\|A\|_{\Co^\alpha}<c_1$. Therefore
    \begin{equation}\label{Eqn::Rough1Form::ExistPDE::Tmp2}
        \textstyle\|R\|_{\Co^{\alpha+1}(\B^n;\R^n)}\le \frac32\tilde C_1\|A\|_{\Co^\alpha}\le\frac12 c_1^{-1}\|A\|_{\Co^\alpha},\quad\text{when }\|A\|_{\Co^\alpha}<c_1\le(3\tilde C_1)^{-1}.
    \end{equation}
     This is part of the estimate in \eqref{Eqn::Rough1Form::ExistPDEQuantControl1}.

    \medskip
    Next we show $\|R\|_{\Co^{\beta+1}(\B^n\backslash\frac34\B^n)}\lesssim\|A\|_{\Co^\alpha}$. Note that by the support assumption {$\supp A\subsetneq\frac12\B^n$} we have $\sqrt{\det g}g^{ij}\big|_{\B^n\backslash\frac12\B^n}=\delta^{ij}$ and $a_i^k\big|_{\B^n\backslash\frac12\B^n}=0$, so the right hand side of \eqref{Eqn::Rough1Form::ExistPDE::Proof1} is zero in $\B^n\backslash\frac12\B^n$. Therefore each $R^k$ is a harmonic functions in the domain $\B^n\backslash\frac12\B^n$.
    
    The estimate $\|R\|_{\Co^{\alpha+1}(\B^n)}\lesssim_\alpha\|A\|_{\Co^\alpha}$ implies $\|R\|_{\Co^{\alpha+1}(\partial(\frac12\B^n))}\lesssim_\alpha\|A\|_{\Co^\alpha}$ since the trace map $\Co^{\alpha+1}(\B^n)\to\Co^{\alpha+1}(\partial(\frac12\B^n))$ is bounded (see Remark \ref{Rmk::Hold::ZygmundonSphere}). By classical interior estimates of harmonic functions (for example, \cite[Theorem 2.10]{GilbargTrudinger}) since $\Delta R\big|_{\B^n\backslash\frac12\B^n}=0$, we have 
    \begin{equation}\label{Eqn::Rough1Form::ExistPDE::TmpBoundEstR}
        \|R\|_{\Co^{\beta+1}(\partial\frac34\B^n)}\lesssim\|R\|_{C^{\lceil\beta\rceil+1}(\partial\frac34\B^n)}\lesssim\|R\|_{C^0(\B^n\backslash\frac23\B^n)}\lesssim\|A\|_{\Co^\alpha}.
    \end{equation}
    
    Therefore, along with the fact that $R\big|_{\partial \B^n}=0$, 
    $R\big|_{\partial(\B^n\backslash\frac34\B^n)}=R\big|_{\partial\frac34\B^n}\cup R\big|_{\partial\B^n}$ has $\Co^{\beta+1}$ norm bounded by a constant times $\|A\|_{\Co^\alpha}$. By classical regularity estimates of harmonic functions (also see Lemma \ref{Lem::Hold::DiriSol}) on $\B^n\backslash\frac34\B^n$ we know
    $$\|R\|_{\Co^{\beta+1}(\B^n\backslash\frac34\B^n)}\lesssim_\beta\|R\|_{\Co^{\beta+1}(\partial(\B^n\backslash\frac34\B^n))}\lesssim_{\alpha,\beta}\|A\|_{\Co^\alpha}.$$
    In particular, there is a $\tilde C_2=\tilde C_2(n,\alpha,\beta)>0$ that depends on neither $A$ nor $R$ such that
    \begin{equation}\label{Eqn::Rough1Form::ExistPDE::Tmp3}
        \|\nabla R\|_{\Co^{\beta}(\B^n\backslash\frac34\B^n)}\le\tilde C_2\|A\|_{\Co^\alpha}.
    \end{equation}
    
    Taking $c_1<\frac 12\tilde C_2^{-1}$ we have $\|\nabla R\|_{\Co^{\beta}(\B^n\backslash\frac34\B^n)}\le \frac12c_1^{-1}\|A\|_{\Co^\alpha}$. Combining this with \eqref{Eqn::Rough1Form::ExistPDE::Tmp2}, completes the proof of \eqref{Eqn::Rough1Form::ExistPDEQuantControl1}.

    \medskip
     We can take $c_1>0$ possibly smaller so that $c_1<\frac13(K_2\tilde C_1)^{-1}$, where $K_2=K_2(n,\alpha,\beta)$ is the constant in Proposition \ref{Prop::Hold::QIFT} and $\tilde C_1$ is the constant in \eqref{Eqn::Rough1Form::ExistPDE::Tmp1}. By \eqref{Eqn::Rough1Form::ExistPDE::Tmp2} we know $\|R\|_{\Co^{\alpha+1}}\le K_2^{-1}$. So by Proposition \ref{Prop::Hold::QIFT}, the map $F=\id+R$ has $\Co^{\alpha+1}$-inverse.  We conclude $\Phi=F^\Inv\in\Co^{\alpha+1}(\B^n;\R^n)$.
    
    Since $\|R\|_{C^1(\B^n;\R^{n\times n})}\lesssim_\alpha\|R\|_{\Co^{\alpha+1}(\B^n;\R^n)}$, by possibly shrinking $c_1$ we can ensure $\|R\|_{C^0}+\|\nabla R\|_{C^0}\le\frac12$. So $F(0)=R(0)\in B^n(0,\frac12)$, which implies $B^n(F(0),\frac16)\subseteq B^n(0,\frac12+\frac16)\subset\frac34\B^n$, and $|F(x_1)-F(x_2)|\ge|x_1-x_2|-|R(x_1)-R(x_2)|\ge\frac12|x_1-x_2|$ for $x_1,x_2\in\B^n$. 
    Thus, if  $|F(x)-F(0)|<\frac16$ then $|x-0|<\frac13$; i.e., $B^n(F(0),\frac16)\subseteq F(\frac13\B^n)$. So $B^n(F(0),\frac16)\subseteq F(\frac13\B^n)\cap\frac34\B^n$ finishing the proof of \ref{Item::Rough1Form::ExistPDE::1.5}.
    
    
    \medskip
    Finally we prove \ref{Item::Rough1Form::ExistPDE::2}. Note that by Proposition \ref{Prop::Hold::QIFT} \ref{Item::Hold::QIFT::PhiEst}, \eqref{Eqn::Rough1Form::ExistPDEQuantControl2} gives $\|\Phi\|_{\Co^{\alpha+1}(\B^n;\R^n)}\lesssim1$ and $ $, which is $\|\Phi\|_{\Co^{\alpha+1}}\le c_1^{-1}$ by choosing $c_1$ small.
    
    By Proposition \ref{Prop::Hold::QIFT} \ref{Item::Hold::QIFT::PhiEst} and using that $c_1<\frac13(C_0\tilde C_1)^{-1}$, we get $\|\nabla\Phi-I\|_{\Co^\alpha}\le K_2\|R\|_{\Co^{\alpha+1}}\le K_2\tilde C_1\|A\|_{\Co^\alpha}\le\frac12c_1^{-1}\|A\|_{\Co^\alpha}$, which proves half of \eqref{Eqn::Rough1Form::ExistPDEQuantControl2}.


    To show the second half  of \eqref{Eqn::Rough1Form::ExistPDEQuantControl2}, we need to show $\|\nabla\Phi-I\|_{\Co^\beta(\partial\B^n)}\lesssim\|A\|_{\Co^\alpha}$. 
    
    The assumption $R\big|_{\partial\B^n}=0$ implies $F\big|_{\partial\B^n}=\id\big|_{\partial\B^n}=\Phi\big|_{\partial\B^n}$ and therefore 
    \begin{equation*}
        (\nabla\Phi-I)\big|_{\partial\B^n}=((\nabla\Phi)\circ\Phi^\Inv)\big|_{\partial\B^n}-I=(\nabla F)^{-1}\big|_{\partial\B^n}-I.
    \end{equation*}

    Fix $\chi\in C_c^\infty(2\B^n\backslash\frac34\B^n)$ such that $\chi\equiv1$ in a neighborhood of $\partial\B^n$, so $\nabla R(x)=\chi(x)\nabla R(x)$ for $x\in\B^n$ near $\partial \B^n$.
    
    We shrink $c_1>0$ so that $c_1<\tilde c_{\B^n,\beta,n}\cdot(2\tilde C_2\|\chi\|_{\Co^\beta})^{-1}$, where $\tilde  C_2$ is in \eqref{Eqn::Rough1Form::ExistPDE::Tmp3}, 
    and $\tilde c_{\B^n,\beta,n}$ is in Lemma \ref{Lem::Hold::CramerMixedPara} \ref{Item::Hold::CramerMixedPara::SingHold}. Then by Lemma \ref{Lem::Hold::Product} \ref{Item::Hold::Product::Hold2} the assumption $\|A\|_{\Co^\alpha}<c_1$ implies
    \begin{equation*}
        \|\chi\nabla R\|_{\Co^\beta(\B^n)}\le 2\|\chi\|_{\Co^\beta_c(\B^n\backslash\frac34\B^n)}\|\nabla R\|_{\Co^\beta(\B^n\backslash\frac34\B^n)}\le2\tilde C_2\|A\|_{\Co^\alpha_c(\frac12\B^n)}<\tilde c_{\B^n,\beta}.
    \end{equation*}
    Therefore we can apply Lemma \ref{Lem::Hold::CramerMixedPara} \ref{Item::Hold::CramerMixedPara::SingHold} to $\chi\nabla R\in\Co^\beta(\B^n;\R^{n\times n})$ to obtain $\|(I+\chi\nabla R)^{-1}-I\|_{\Co^\beta}\le 2\|\chi\nabla R\|_{\Co^\beta}$. Hence, along with \eqref{Eqn::Hold::HZNormforSphere},
    $$\|(I+\nabla R)^{-1}-I\|_{\Co^\beta(\partial\B^n)}\lesssim\|(I+\chi\nabla R)^{-1}-I\|_{\Co^\beta(\B^n)}\lesssim\|\chi\nabla R\|_{\Co^\beta(\B^n)}\lesssim_\chi\|R\|_{\Co^{\beta+1}(\B^n\backslash\frac34\B^n)}\lesssim\|A\|_{\Co^\alpha}.$$
    So by possibly shrinking $c_1>0$, we get $\|\nabla \Phi-I\|_{\Co^\beta(\partial\B^n)}\le \frac12 c_1^{-1}\|A\|_{\Co^\alpha}$, which completes the second half of \eqref{Eqn::Rough1Form::ExistPDEQuantControl2}.
    %
\end{proof}

We now have pushforward 1-forms $\eta^1=F_*\lambda^1,\dots,\eta^n=F_*\lambda^n$. Their norms admit some control, as the next lemma shows. 

\begin{lem}\label{Lem::Rough1Form::AtoB}
Let $\alpha>0$ and let $\beta\in[\alpha,\alpha+1]$. There is a $c_2=c_2(n,\alpha,\beta)>0$ such that the following holds.
Let  $A\in\Co^\alpha_c(\frac12\B^n,\R^{n\times n})$ be the coefficient matrix for $\lambda^1,\dots,\lambda^n$ (see \eqref{Eqn::Rough1Form::LambdaEta}) satisfying the assumptions of Proposition \ref{Prop::Rough1Form::ExistPDE} and also satisfying
\begin{enumerate}[parsep=-0.3ex,label=(\alph*)]
    \item $\|A\|_{\Co^\alpha(\B^n;\R^{n\times n})}<c_2$.
    \item For $k=1,\dots,n$, $d\lambda^k\in\Co^{\beta-1}\mleft(\B^n;\wedge^2 T^*\B^n\mright)$  with $\sum_{k=1}^n\|d\lambda^k\|_{\Co^{\beta-1}(\B^n;\wedge^2T^*\B^n)}<c_2$.
\end{enumerate}

Suppose $\Phi =F^\Inv:\B^n_y\to\B^n_x$ satisfies the conclusions of Proposition \ref{Prop::Rough1Form::ExistPDE}.  Then, for the 1-forms $\eta^k=\Phi^*\lambda^k$ ($k=1,\dots,n$) with coefficient matrix $B=(b^i_j)_{n\times n}$ (see \eqref{Eqn::Rough1Form::LambdaEta}), we have:
\begin{enumerate}[parsep=-0.3ex,label=(\roman*)]
    \item\label{Item::Rough1Form::AtoB::1} $B$ satisfies the PDE system \eqref{Eqn::Rough1Form::PDEforB}.
    \item\label{Item::Rough1Form::AtoB::2} $\|B\|_{\Co^\alpha(\B^n;\R^{n\times n})}+\|B\|_{\Co^\beta(\partial\B^n;\R^{n\times n})}<c_2^{-1}\|A\|_{\Co^\alpha}$.
    \item\label{Item::Rough1Form::AtoB::3} $d\eta^k\in\Co^{\beta-1}\mleft(\B^n;\wedge^2T^*\B^n\mright)$ for $k=1,\dots,n$, with $\|d\eta^k\|_{\Co^{\beta-1}(\B^n;\wedge^2T^*\B^n)}<c_2^{-1}\|d\lambda^k\|_{\Co^{\beta-1}}$.
\end{enumerate}
\end{lem}
\begin{proof}Part \ref{Item::Rough1Form::AtoB::1} is obtained in Lemma \ref{Lem::Rough1Form::TransitionBetweenPDEs}.

For part \ref{Item::Rough1Form::AtoB::2}, write $\Phi=(\phi^1,\dots,\phi^n)$, where $\phi^k\in\Co^{\alpha+1}(\B^n)$, $k=1,\dots,n$. Therefore 
\begin{equation}\label{Eqn::Rough1Form::AtoB::Proof1}
    \eta^k=\Phi^*\Big(dx^k+\sum_{i=1}^na_i^kdx^i\Big)=d\phi^k+\sum_{i=1}^n(a_i^k\circ\Phi)d\phi^i,\quad b_j^k=\frac{\partial(\phi^k-y^k)}{\partial y^j}+\sum_{i=1}^n(a_i^k\circ\Phi)\frac{\partial\phi^i}{\partial y^j},\quad 1\le j,k\le n.
\end{equation}

From \eqref{Eqn::Rough1Form::AtoB::Proof1} we know that $\|B\|_{\Co^\alpha}\lesssim\|\nabla\Phi-I\|_{\Co^\alpha}+\|A\circ\Phi\|_{\Co^\alpha}\|\nabla\Phi\|_{\Co^\alpha}$. By Proposition \ref{Prop::Rough1Form::ExistPDE} \ref{Item::Rough1Form::ExistPDE::2} we know $\|\Phi\|_{\Co^{\alpha+1}}\lesssim1$ and $\|\nabla\Phi-I\|_{\Co^\alpha}\lesssim\|A\|_{\Co^\alpha}$. By 
Proposition \ref{Prop::Hold::QComp} \ref{Item::Hold::QComp::>1} we get $\|A\circ\Phi\|_{\Co^\alpha}\lesssim\|A\|_{\Co^\alpha}$. Combining these we get 
\begin{equation}\label{Eqn::Rough1Form::AtoB::Tmp1}
    \|B\|_{\Co^\alpha}\lesssim\|\nabla\Phi-I\|_{\Co^\alpha}+\|A\circ\Phi\|_{\Co^\alpha}\|\nabla\Phi\|_{\Co^\alpha}\lesssim\|A\|_{\Co^\alpha}.
\end{equation}

Since we have $A\equiv 0$ outside $\frac{1}{2}\B^n$ {in particular $A\big|_{\partial\B^n}=0$}, it follows that $\eta^k=d\phi^k$ on $\partial\B^n$. 
Therefore $\|B\|_{\Co^\beta(\partial\B^n)}=\|\nabla\Phi-I\|_{\Co^\beta(\partial\B^n)}$.
So by Proposition \ref{Prop::Rough1Form::ExistPDE} \ref{Item::Rough1Form::ExistPDE::2}
\begin{equation}\label{Eqn::Rough1Form::AtoB::Tmp2}
    \|B\|_{\Co^\beta(\partial\B^n)}=\|\nabla\Phi-I\|_{\Co^\beta(\partial\B^n)}\lesssim\|A\|_{\Co^\alpha}
\end{equation}

By choosing $c_2>0$ small, \eqref{Eqn::Rough1Form::AtoB::Tmp1} and \eqref{Eqn::Rough1Form::AtoB::Tmp2} complete the proof of \ref{Item::Rough1Form::AtoB::2}.

\medskip
Finally, for \ref{Item::Rough1Form::AtoB::3}, we apply 
Proposition \ref{Prop::Hold::QIFT} \ref{Item::Hold::QIFT::dFormEst}
with $\lambda=\lambda^k$,
for each $k=1,\dots,n$. Since $d(F_*\lambda^k)=d\eta^k$, by \eqref{Eqn::Hold::QIFT::dFormEst} we get $\|d\eta^k\|_{\Co^{\beta-1}(\B^n;\wedge^2T^*\B^n)}\lesssim\|d\lambda^k\|_{\Co^{\beta-1}}$. Taking $c_2$ smaller, we complete the proof.
\end{proof}
\subsection{The regularity proposition}

In this part, we show that the  1-forms $\eta^1,\dots,\eta^n$ are indeed $\Co^\beta$, by using the interior regularity theory for elliptic PDEs.

\begin{prop}\label{Prop::Rough1Form::RegularityPDE}
{Let $\alpha>0$ and $\beta\in[\alpha,\alpha+1]$.
There is a $c_3=c_3(n,\alpha,\beta)>0$}, such that if $\eta^1,\dots,\eta^n\in\Co^\alpha(\B^n;T^*\B^n)$ with coefficient matrix $B\in\Co^\alpha(\B^n;\R^{n\times n})$ (see \eqref{Eqn::Rough1Form::LambdaEta}) such that $B$ solves the PDE \eqref{Eqn::Rough1Form::PDEforB}, $B\big|_{\partial \B^n}\in \Co^{\beta}(\partial\B^n;\R^{n\times n})$ with
\begin{equation}\label{Eqn::Rough1Form::RegularityPDE::Assumption}
    \|B\|_{\Co^\alpha(\B^n;\R^{n\times n})}+\|B\|_{\Co^\beta(\partial\B^n;\R^{n\times n})}+\sum_{l=1}^n\|d\eta^l\|_{\Co^{\beta-1}\mleft(\B^n;\wedge^2T^*\B^n\mright)}<c_3,
\end{equation} then $B\in\Co^\beta(\B^n;\R^{n\times n})$. Moreover
\begin{equation}\label{Eqn::Rough1Form::RegularityPDE::Conclusion}
    \|B\|_{\Co^\beta(\B^n;\R^{n\times n})}\le c_3^{-1}\Big(\|B\|_{\Co^\alpha(\B^n;\R^{n\times n})}+\|B\|_{\Co^\beta(\partial\B^n;\R^{n\times n})}+\sum_{l=1}^n\|d\eta^l\|_{\Co^{\beta-1}(\B^n;\wedge^2T^*\B^n)}\Big).
\end{equation}
\end{prop}
\begin{proof}
We can write $B=\Pc_0\Delta B+(B-\Pc_0\Delta B)$, where $\Pc_0=\Pc_0^{\B^n}$ is given in Lemma \ref{Lem::Hold::LapInvBdd}, which is the zero Dirichlet boundary solution operator on the unit ball.

Note that $B-\Pc_0\Delta B$ is the harmonic function whose boundary value equals to $B\big|_{\partial\B^n}$ (which might not be zero). By Lemma \ref{Lem::Hold::DiriSol} using the assumption $B\big|_{\partial\B^n}\in\Co^\beta(\partial\B^n;\R^{n\times n})$, we get $B-\Pc_0\Delta B\in\Co^\beta(\B^n;\R^{n\times n})$ and
\begin{equation}\label{Eqn::Rough1Form::RegularityPDE::Proof0}
    \|B-\Pc_0\Delta B\|_{\Co^\beta(\B^n)}\lesssim\|B\|_{\Co^\beta(\partial\B^n)}.
\end{equation}

We can rewrite \eqref{Eqn::Rough1Form::PDEforB} as
\begin{equation}\label{Eqn::Rough1Form::RegularityPDE::Proof1}
    -\sum_{i=1}^n\Coorvec{y^i}b_i^k=\sum_{i,j=1}^n\Coorvec{y^i}\left(\big(\sqrt{\det h}h^{ij}-\delta^{ij}\big)b_j^k\right),\quad\text{in }\B^n_y,\quad k=1,\dots,n.
\end{equation}
The left hand side of \eqref{Eqn::Rough1Form::RegularityPDE::Proof1} is $\codiff_{\R^n}\eta^k$. By Lemma \ref{Lem::Rough1Form::TaylorExpansionofRiemMetric}, the right hand side of \eqref{Eqn::Rough1Form::RegularityPDE::Proof1} is the derivatives of rational functions of the components of $B$, which vanish to second order at $B=0$. More precisely, using Lemma \ref{Lem::Rough1Form::TaylorExpansionofRiemMetric}, we can rewrite \eqref{Eqn::Rough1Form::RegularityPDE::Proof1} as
\begin{equation}\label{Eqn::Rough1Form::RegularityPDE::EqnasCalR}
    \codiff_{\R^n}\eta^k=-\sum_{i=1}^n\Coorvec{y^i}b_i^k=-\sum_{i=1}^n\Coorvec{y^i}\Rc_i^k(B),\quad\text{in }\B^n_y,\quad k=1,\dots,n.
\end{equation}
Here $\Rc_i^k$ are rational functions (see Lemma \ref{Lem::Rough1Form::TaylorExpansionofRiemMetric}) defined in a neighborhood of origin in $\R^{n\times n}$ with $|\Rc_i^k(u)|\lesssim|u|_{\R^{n\times n}}^2$ for suitably small matrices $u\in\R^{n\times n}$, and we have 
\begin{equation}\label{Eqn::Rough1Form::RegularityPDE::ProofTmp1}
    |\Rc_i^k(u_1)-\Rc_i^k(u_2)|\lesssim(|u_1|_{\R^{n\times n}}+|u_2|_{\R^{n\times n}})|u_1-u_2|_{\R^{n\times n}},\quad\text{when }u_1,u_2\in\R^{n\times n}\text{ small}.
\end{equation}

We can pass this fact from matrices to matrix-valued functions. Indeed, $\Rc_i^k$ has convergent power expansion in a neighborhood of $0$ as
\begin{equation}\label{Eqn::Rough1Form::RegularityPDE::PowerSeriesforR}
    \Rc_i^k(u)=\sum_{r=2}^\infty\sum_{j_1,\dots,j_r,l_1,\dots,l_r=1}^n a_{i,r;l_1\dots l_r}^{k;j_1\dots j_r}u_{j_1}^{l_1}\dots u_{j_r}^{l_r},\quad\text{converging when }|u|_{\R^{n\times n}}\text{ is small.}
\end{equation}
Here $a_{i,r;l_1\dots l_r}^{k;j_1\dots j_r}\in\R$. The power expansion starts at $r=2$ since the zero and the first order terms  all vanish.

By Lemma  \ref{Lem::Hold::Product}, we can replace $u\in\R^{n\times n}$ in \eqref{Eqn::Rough1Form::RegularityPDE::PowerSeriesforR} by $f\in\Co^\gamma(\B^n;\R^{n\times n})$,  
and as in \eqref{Eqn::Rough1Form::RegularityPDE::ProofTmp1}, for $\gamma>0$ there is a $\tilde C_{\mathcal R,\gamma}>0$, such that when $\|f_1\|_{\Co^\gamma(\B^n;\R^{n\times n})}+\|f_2\|_{\Co^\gamma(\B^n;\R^{n\times n})}\le\tilde C_{\mathcal R,\gamma}^{-1}$,
\begin{equation}\label{Eqn::Rough1Form::RegularityPDE::Proof1.5}
    \|\Rc_i^k(f_1)-\Rc_i^k(f_2)\|_{\Co^\gamma(\B^n)}\le \tilde C_{\mathcal R,\gamma}(\|f_1\|_{\Co^\gamma(\B^n;\R^{n\times n})}+\|f_2\|_{\Co^\gamma(\B^n;\R^{n\times n})})\|f_1-f_2\|_{\Co^\gamma(\B^n;\R^{n\times n})},\quad1\le i,k\le n.
\end{equation}

Using the fact that $-\Delta\eta^k=d\codiff_{\R^n}\eta^k+\codiff_{\R^n}d\eta^k$, we further have
$$\Delta\eta^k=d\sum_{i=1}^n\Coorvec{y^i}\Rc_i^k(B)-\codiff_{\R^n} d\eta^k=\sum_{i,j=1}^n\frac{\partial^2}{\partial y^j\partial y^i}\Rc_i^k(B)dy^j-\sum_{j=1}^n\mleft\langle \codiff_{\R^n} d\eta^k,\Coorvec{y^j}\mright\rangle dy^j.$$
Here $\langle \cdot, \cdot \rangle$ denotes the pairing between
$1$ forms and vector fields.

On the other hand, $\Delta\eta^k=\sum_{j=1}^n\Delta(\delta_j^k+b_j^k)dy^j=\sum_{j=1}^n\Delta b_j^kdy^j$, therefore
\begin{equation}\label{Eqn::Rough1Form::RegularityPDE::Proof2}
    \Delta b_j^k=\mleft\langle\Delta\eta^k,\Coorvec{y^j}\mright\rangle=\sum_{i,j=1}^n\frac{\partial^2}{\partial y^j\partial y^i}\Rc_i^k(B)-\mleft\langle \codiff_{\R^n} d\eta^k,\Coorvec{y^j}\mright\rangle,\quad\text{in }\B^n_y,\quad k=1,\dots,n.
\end{equation}


Let $\tilde\xi_0:=\min(\tilde C_{\Rc,\alpha}^{-1},\tilde C_{\Rc,\beta}^{-1})$ where $\tilde C_{\Rc,\gamma}$ is the constant in \eqref{Eqn::Rough1Form::RegularityPDE::Proof1.5}. Let $\xi=\xi_B\in(0,\tilde\xi_0]$ to be determined. We define metric spaces $\Xs_{\gamma,\xi}$ and an operator $\Tc_B:\Xs_{\alpha,\tilde\xi_0}\to \Co^\alpha(\B^n;\R^{n\times n})$ by
\begin{equation}\label{Eqn::Rough1Form::RegularityPDE::DefofContMap}
\begin{gathered}
    \Xs_{\gamma,\xi}:=\{f\in\Co^\gamma(\B^n;\R^{n\times n}):\|f\|_{\Co^\gamma}\le \xi\}\subset\Co^\gamma(\B^n;\R^{n\times n}),\quad\text{for }\gamma\in\{\alpha,\beta\}\text{ and }\xi\in(0,\tilde\xi_0].
    \\
    \Tc_B[f]^k_j:=b_j^k-\Pc_0\Delta b_j^k-\Big\langle\Pc_0(\codiff_{\R^n} d\eta^k),\Coorvec{y^j}\Big\rangle+\sum_{i=1}^n\Pc_0\Big(\frac{\partial^2\Rc_i^k(f)}{\partial y^j\partial y^i}\Big),\quad 1\le j,k\le n.
\end{gathered}
\end{equation}
We endow $\Xs_{\gamma,\xi}$ with the metric induced by the norm $\|\cdot\|_{\Co^\gamma}$, which makes $\Xs_{\gamma,\xi}$ a complete metric space.

Note that from \eqref{Eqn::Rough1Form::RegularityPDE::DefofContMap} and \eqref{Eqn::Rough1Form::RegularityPDE::Proof2} we have $B=\Tc_B[B]$. Our goal is to show that when $c_3$ and $\xi$ are both suitably small, we have $B\in\Xs_{\alpha,\xi}$ and that  $\Tc_B$ is a contraction mapping on both $\Xs_{\alpha,\xi}$ and $\Xs_{\beta,\xi}$, thus by uniqueness of the fixed point we conclude that $B$ is a $\Co^\beta$-matrix and $\|B\|_{\Co^\beta}<\xi$.
\medskip


{By Lemma \ref{Lem::Hold::DiriSol}, $\Pc_0:\Co^{\gamma-2}(\B^n;\R^{n\times n})\to\Co^\gamma(\B^n;\R^{n\times n})$ is bounded for $\gamma\in\{\alpha,\beta\}$}. By \eqref{Eqn::Rough1Form::RegularityPDE::Proof1.5}, we know for every $f_1,f_2\in\Xs_{\gamma,\xi}$,
\begin{equation}\label{Eqn::KeyRegPDE::Tmp1}
    \|\Tc_B[f_1]-\Tc_B[f_2]\|_{\Co^\gamma}\le \|\Pc_0\|_{\Co^{\gamma-2}\to\Co^\gamma}\|\nabla^2\|_{\Co^\gamma\to\Co^{\gamma-2}}\sum_{k=1}^n\|\Rc_k^l(f_1)-\Rc_k^l(f_2)\|_{\Co^\gamma}\le \xi C'_{\Rc,\gamma}\|f_1-f_2\|_{\Co^\gamma}.
\end{equation}
where $C'_{\Rc,\gamma}>1$ is a constant that only depends on $n,(\Rc_j^k),\gamma$ but not on $B,\xi,f_1,f_2$.

On the other hand $\Tc_B[0]_j^k=b_j^k-\Pc_0\Delta b_j^k-\big\langle\Pc_0(\codiff_{\R^n} d\eta^k),\Coorvec{y^j}\big\rangle$. By \eqref{Eqn::Rough1Form::RegularityPDE::Proof0}, for $\gamma\in \{\alpha,\beta\}$,
\begin{equation}\label{Eqn::KeyRegPDE::Tmp2}
\|\Tc_B[0]\|_{\Co^\gamma}\le\|\Tc_B[0]\|_{\Co^\beta}\le\|B-\Pc_0\Delta B\|_{\Co^\beta(\B^n)}+\sum_{k=1}^n\|\Pc_0(\codiff_{\R^n} d\eta^k)\|_{\Co^\beta}\lesssim\|B\|_{\Co^\beta(\partial\B^n)}+\sum_{k=1}^n\|d\eta^k\|_{\Co^{\beta-1}}. 
\end{equation}
So by possibly increasing $C'_{\Rc,\gamma}$, we have, for {$f_1\in\Xs_{\gamma,\xi}$},
using \eqref{Eqn::KeyRegPDE::Tmp1} and \eqref{Eqn::KeyRegPDE::Tmp2},
\begin{equation}\label{Eqn::KeyRegPDE::Tmp3}
\|\Tc_B[f_1]\|_{\Co^\gamma}
\leq \|\Tc_B[f_1]-\Tc_B[0]\|_{\Co^\gamma}+\|\Tc_B[0]\|_{\Co^\gamma}
\le C'_{\Rc,\gamma}\Big(\|B\|_{\Co^\beta(\partial\B^n)}+\sum_{l=1}^n\|d\eta^l\|_{\Co^{\beta-1}}+\xi\|f_1\|_{\Co^\gamma}\Big).
\end{equation}

Take $c_3>0$ satisfying $c_3<\frac14\max(1,C'_{\Rc,\alpha},C'_{\Rc,\beta})^{-2}$, and take 
\begin{equation}\label{Eqn::Rough1Form::RegularityPDE::Proof3}
    \xi=\xi_B:=2\max(C'_{\Rc,\alpha},C'_{\Rc,\beta})\left(\|B\|_{\Co^\alpha(\B^n;\R^{n\times n})}+\|B\|_{\Co^\beta(\partial\B^n)}+\sum_{k=1}^n\|d\eta^k\|_{\Co^{\beta-1}}\right).
\end{equation}

By the assumption \eqref{Eqn::Rough1Form::RegularityPDE::Assumption}, $\xi_B\le\frac12\max(C'_{\Rc,\alpha},C'_{\Rc,\beta})^{-1}<\tilde\xi_0$, so $\Tc_B$ is defined on $\Xs_{\alpha,\tilde\xi_0}$ and by \eqref{Eqn::KeyRegPDE::Tmp3} $\Tc_B$ maps $\Xs_{\gamma,\xi_B}$ into $\Xs_{\gamma,\xi_B}$ for $\gamma\in\{\alpha,\beta\}$. 

Since $\xi_B C'_{\Rc,\gamma}<\frac12$ for $\gamma\in\{\alpha,\beta\}$, using \eqref{Eqn::KeyRegPDE::Tmp1},
$\Tc_B$ is a contraction mapping on the domain $\Xs_{\gamma,\xi_B}$, for $\gamma\in\{\alpha,\beta\}$.

Note that $\xi_B\ge\|B\|_{\Co^\alpha}$, and so $B\in \{f\in\Co^\alpha(\B^n;\R^{n\times n}):\|f\|_{\Co^\alpha}\le \xi_B\}=\Xs_{\alpha,\xi_B}$. Therefore, $B$ is a fixed point for $\Tc_B$ in $\Xs_{\alpha,\xi_B}$, which is unique since $\Tc_B$ is a contraction mapping on $\Xs_{\alpha,\xi_B}$.

On the other hand $\Tc_B$ also has a unique fixed point in $\Xs_{\beta,\xi_B}\subsetneq\Xs_{\alpha,\xi_B}$. Therefore, by uniqueness, $B\in\Xs_{\beta,\xi_B}=\{f\in\Co^\beta(\B^n;\R^{n\times n}):\|f\|_{\Co^\beta}\le \xi_B\}$.
In particular, $\|B\|_{\Co^\beta}\le\xi_B$. Thus by \eqref{Eqn::Rough1Form::RegularityPDE::Proof3}, 
$$\|B\|_{\Co^\beta(\B^n)}\le\xi_B\lesssim_{\alpha,\beta} \|B\|_{\Co^\alpha(\B^n)}+\|B\|_{\Co^\beta(\partial\B^n)}+\sum_{l=1}^n\|d\eta^l\|_{\Co^{\beta-1}}.$$

Thus, we have established \eqref{Eqn::Rough1Form::RegularityPDE::Conclusion}
which completes the proof.
\end{proof}

\subsection{The proof of Theorem \ref{KeyThm::Rough1Form} and a qualitative regularity improving}\label{Section::Rough1FormSum}
Using Propositions \ref{Prop::Rough1Form::ExistPDE} and \ref{Prop::Rough1Form::RegularityPDE} we can prove Theorem \ref{KeyThm::Rough1Form}.

\begin{proof}[Proof of Theorem \ref{KeyThm::Rough1Form}]
Let $c_1,c_2,c_3>0$ be the small constants in Proposition \ref{Prop::Rough1Form::ExistPDE}, Lemma \ref{Lem::Rough1Form::AtoB}, and Proposition \ref{Prop::Rough1Form::RegularityPDE}. We take $c=\frac1{2n^2}\min(c_1,c_2c_3)$ in the assumption of Theorem \ref{KeyThm::Rough1Form}.

Let $F$, $R=F-\id$, $\Phi=F^\Inv$, $A$, $B$ and $\eta^i=F_*\lambda^i$ be as in Proposition \ref{Prop::Rough1Form::ExistPDE}. {Recall $\eta^i$ and $B=(b_i^j)$ are given in \eqref{Eqn::Rough1Form::LambdaEta} and \eqref{Eqn::Rough1Form::LambdaEtaMatrix}.} 

When the assumption \eqref{Eqn::Rough1Form::Assumption} is satisfied, by Proposition \ref{Prop::Rough1Form::ExistPDE} \ref{Item::Rough1Form::ExistPDE::1.5} we have that $F$ is $\Co^{\alpha+1}$-diffeomorphism and satisfies $B^n(F(0),\frac16)\subseteq F(\frac13\B^n)\cap\frac34\B^n$. And by \eqref{Eqn::Rough1Form::ExistPDEQuantControl1},  we have $\|F-\id\|_{\Co^{\alpha+1}}=\|R\|_{\Co^{\alpha+1}}\le c_1^{-1}\|A\|_{\Co^\alpha}\le \frac1{2n} c^{-1}\sum_{i=1}^n\|\lambda^i-dx^i\|_{\Co^\alpha}$. This implies half of the estimate \eqref{Eqn::Rough1Form::Conclusion}.

By Lemma \ref{Lem::Rough1Form::AtoB} \ref{Item::Rough1Form::AtoB::2} and \ref{Item::Rough1Form::AtoB::3}, we have $\|B\|_{\Co^\alpha(\B^n)}+\|B\|_{\Co^\beta(\partial\B^n)}<c_2^{-1}\|A\|_{\Co^\alpha}$ and $\|d\eta^k\|_{\Co^{\beta-1}}<c_2^{-1}\|d\lambda^k\|_{\Co^{\beta-1}}$, $k=1,\dots,n$. 

Thus, $\|B\|_{\Co^\alpha(\B^n)}+\|B\|_{\Co^\beta(\partial\B^n)}+\sum_{k=1}^n\|d\eta^k\|_{\Co^{\beta-1}}<c_2^{-1}\|A\|_{\Co^\alpha}+c_2^{-1}\sum_{k=1}^n\|d\lambda^k\|_{\Co^{\beta-1}}<2c_2^{-1}c<c_3$.
By \eqref{Eqn::Rough1Form::RegularityPDE::Conclusion} in Proposition \ref{Prop::Rough1Form::RegularityPDE}, we get 
\begin{align*}
    &\sum_{k=1}^n\|\eta^k-dy^k\|_{\Co^\beta}\le n\|B\|_{\Co^\beta(\B^n)}\le n c_3^{-1}\mleft(\|B\|_{\Co^\alpha}+\|B\|_{\Co^\beta(\partial\B^n)}+\sum_{k=1}^n\|d\eta^k\|_{\Co^{\beta-1}}\mright)\\
    &\le n c_3^{-1}\cdot c_2^{-1}\mleft(\|A\|_{\Co^\alpha}+\sum_{k=1}^n\|d\lambda^k\|_{\Co^{\beta-1}}\mright)\le n^2(c_2c_3)^{-1}\sum_{k=1}^n\mleft(\|\lambda^k-dx^k\|_{\Co^\alpha}+\|d\lambda^k\|_{\Co^{\beta-1}}\mright).
\end{align*}
This gives the second half of the estimate \eqref{Eqn::Rough1Form::Conclusion} since $n^2(c_2c_3)^{-1}\le\frac12 c^{-1}$.
\end{proof}


In Theorem \ref{KeyThm::Rough1Form}, we assumed \eqref{Eqn::Rough1Form::Assumption} which is a smallness assumption.  When \eqref{Eqn::Rough1Form::Assumption}
is not satisfied, we may use a scaling argument to transfer to a setting where it is satisfied, as the next result shows.

\begin{prop}[The scaling argument]\label{Prop::Rough1Form::Scaling}Let $\alpha>0$, $\beta\in[\alpha,\alpha+1]$
and let $\mu_0,\tilde c,M>0$. There exists a $\kappa_0=\kappa_0(n,\alpha,\beta,\mu_0,\tilde c,M)\in(0,\mu_0]$ that satisfies the following:

Suppose $\xi^1,\dots,\xi^n\in\Co^\alpha(\mu_0\B^n;T^*\R^n)$ such that $\xi^i\big|_0=dx^i\big|_0$ for $i=1,\dots,n$ and $d\xi^1,\dots,d\xi^n\in\Co^{\beta-1}\mleft(\mu_0\B^n;\wedge^2T^*\R^n\mright)$ with estimate 
\begin{equation}\label{Eqn::Rough1Form::Scaling::AssumptionOnTheta}
    \sum_{i=1}^n\|\xi^i\|_{\Co^\alpha(\mu_0\B^n;T^*\R^n)}+\|d\xi^i\|_{\Co^{\beta-1}(\mu_0\B^n;\wedge^2T^*\R^n)}<M.
\end{equation}

Then there are 1-forms $\lambda^1,\dots,\lambda^n\in\Co^\alpha(\B^n;T^*\B^n)$ such that
\begin{enumerate}[parsep=-0.3ex,label=(\roman*)]
    \item\label{Item::Rough1Form::Scaling::FormEqual} $\lambda^i\big|_{\frac13\B^n}=\frac1{\kappa_0}\cdot(\phi_{\kappa_0}^*\xi^i)\big|_{\frac13\B^n}$, $i=1,\dots,n$, where $\phi_{\kappa_0}(x):=\kappa_0\cdot x$ is the scaling map.
    \item\label{Item::Rough1Form::Scaling::Est} $\lambda^1,\dots,\lambda^n$ satisfy the assumptions of Theorem \ref{KeyThm::Rough1Form} with the constant $c=\tilde c$. That is,
    \begin{itemize}
        \item $\lambda^1,\dots,\lambda^n$ span the cotangent space at every point in $\B^n$.
        \item $\supp(\lambda^i-dx^i)\subsetneq\frac12\B^n $.
        \item $\sum_{i=1}^n(\|\lambda^i-dx^i\|_{\Co^\alpha}+\|d\lambda^i\|_{\Co^{\beta-1}})\le \tilde c$.
    \end{itemize}
\end{enumerate}
\end{prop}

\begin{proof}
First we construct 1-forms $\rho^1,\dots,\rho^n\in\Co^\beta(\mu_0\B^n;T^*\R^n)$ such that for $i=1,\dots,n$, 
\begin{enumerate}[parsep=-0.3ex,label=(\alph*)]
    \item\label{Item::Rough1Form::Scaling::RhoFormEqual} $\rho^i\big|_0=0$ and $d\rho^i\big|_{\frac{\mu_0}2\B^n}=d\xi^i\big|_{\frac{\mu_0}2\B^n}$.
    \item\label{Item::Rough1Form::Scaling::RhoFormEst} There is a $C_0=C_0(n,\alpha,\beta,\mu_0)>0$ that does not depend on $\xi^i$, such that 
    \begin{equation}\label{Eqn::Rough1Form::Scaling::RhoFormEst}
        \|\rho^i\|_{\Co^\beta(\mu_0\B^n;T^*\R^n)}\le C_0(\|\xi^i\|_{\Co^\alpha(\mu_0\B^n;T^*\R^n)}+\|d\xi^i\|_{\Co^{\beta-1}(\mu_0\B^n;\wedge^2T^*\R^n)}).
    \end{equation}
\end{enumerate}

Take a $\chi_0\in C_c^\infty(\mu_0\B^n)$ such that $\chi_0\big|_{\frac{\mu_0}2\B^n}\equiv1$. Define
\begin{equation}
    \tilde\rho^i:=-\Ga\ast\codiff d(\chi_0\xi^i)=-\Ga\ast\codiff (\chi_0\cdot d\xi^i+d\chi_0\wedge \xi^i),\quad\rho^i:=\tilde\rho^i-(\tilde\rho^i\big|_0),\quad i=1,\dots,n.
\end{equation}
Recall $\codiff$ is the codifferential from Notation \ref{Note::Hold::Codiff}, and $\Ga$ is the fundamental solution of Laplacian as in \eqref{Eqn::Hold::Newtonian}. 
The convolution is defined in $\R^n$ using Lemma \ref{Lem::Hold::GreensOp}, since the support $\supp \codiff d(\chi_0\xi^i)\subseteq\supp\chi_0\Subset\mu_0\B^n$ is compact.

Clearly $\rho^i\big|_0=0$. Similar to the proof of Proposition \ref{Prop::Hold::QIFT}, since $\chi_0\big|_{\frac{\mu_0}2\B^n}\equiv1$, we have $$d\xi^i\big|_{\frac{\mu_0}2\B^n}=d(\chi_0\xi^i)\big|_{\frac{\mu_0}2\B^n}=-(\codiff d+d\codiff)(\Ga\ast d(\chi_0\xi^i))\big|_{\frac{\mu_0}2\B^n}=-(\Ga\ast d\codiff d(\chi_0\xi^i))\big|_{\frac{\mu_0}2\B^n}=d\tilde\rho^i\big|_{\frac{\mu_0}2\B^n}=d\rho^i\big|_{\frac{\mu_0}2\B^n}.$$
So condition \ref{Item::Rough1Form::Scaling::RhoFormEqual} is satisfied.

By Lemma \ref{Lem::Hold::GreensOp} we have, for every $\mu>0$,
\begin{equation}\label{Eqn::Rough1Form::ScalingProof3}
    \|\Gc\ast\codiff\omega\|_{\Co^{\beta}(\mu\B^n;\wedge^2T^*\R^n)}\lesssim_{\beta,\mu}\|\omega\|_{\Co^{\beta-1}(\mu\B^n;T^*\R^n)},\quad\forall \omega\in\Co^{\beta-1}_c(\mu\B^n;T^*\R^n).
\end{equation}

Take $\omega=\codiff (\chi_0\cdot d\xi^i+d\chi_0\wedge \xi^i)$, $\mu=\mu_0$ in \eqref{Eqn::Rough1Form::ScalingProof3} and by Lemma \ref{Lem::Hold::Product}, we have
\begin{equation}\label{Eqn::Rough1Form::Scaling::RhoFormEstCompute}
    \begin{aligned}
    &\|\rho^i\|_{\Co^\beta(\mu_0\B^n;T^*\R^n)}\le \|\tilde\rho^i\|_{\Co^\beta(\mu_0\B^n;T^*\R^n)}+|\tilde\rho^i(0)|\le 2\|\tilde\rho^i\|_{\Co^\beta(\mu_0\B^n;T^*\R^n)}
    \\
    &\lesssim_{\beta,\mu_0}\|\codiff(\chi_0\cdot d\xi^i+d\chi_0\wedge \xi^i)\|_{\Co^{\beta-2}(\mu_0\B^n;\wedge^2T^*\R^n)}\lesssim_{\beta,\mu_0}\|\chi_0\|_{\Co^{\alpha+1}}\|d\xi^i\|_{\Co^{\beta-1}}+\|d\chi_0\|_{\Co^{\alpha}}\|\xi^i\|_{\Co^{\beta-1}}
    \\
    &\lesssim_{\alpha,\beta,\mu_0,\chi_0}\|d\xi^i\|_{\Co^{\beta-1}(\mu_0\B^n;\wedge^2T^*\R^n)}+\|\xi^i\|_{\Co^\alpha(\mu_0\B^n;T^*\R^n)}.
\end{aligned}
\end{equation}

\eqref{Eqn::Rough1Form::Scaling::RhoFormEstCompute} gives us the $C_0$ for condition \ref{Item::Rough1Form::Scaling::RhoFormEst}. This complete the proof of \ref{Item::Rough1Form::Scaling::RhoFormEqual} and \ref{Item::Rough1Form::Scaling::RhoFormEst} and we get $\rho^1,\dots,\rho^n$ as desired.

\medskip
Fix $\chi_1\in C_c^\infty(\frac12\B^n)$ such that $\chi_1\big|_{\frac13\B^n}\equiv1$. For $\kappa>0$, let $\phi_\kappa(x):=\kappa\cdot x$, so $\phi_\kappa$ maps $\frac12\B^n$ into $\frac{\mu_0}2\B^n$ when $\kappa\in(0,\mu_0]$. For $\kappa\in(0,\mu_0]$, we define 1-forms $\lambda_\kappa^1,\dots,\lambda_\kappa^n$ and $\tau_\kappa^1,\dots,\tau_\kappa^n$ by
\begin{equation}\label{Eqn::Rough1Form::Scaling::DefLambdaTau}
    \lambda_\kappa^i:=dx^i+\tfrac1\kappa\chi_1\cdot\phi_\kappa^*(\xi^i-dx^i),\quad\tau_\kappa^i:=\tfrac1\kappa\chi_1\cdot(\phi_\kappa^*\rho^i)+\tfrac1\kappa\Gc\ast\codiff\left(d\chi_1\wedge \phi_\kappa^*(\xi^i-dx^i)\right),\quad i=1,\dots,n.
\end{equation}

Since $\xi^i\in\Co^\alpha$ and $\rho^i\in\Co^\beta$, we have $\lambda_\kappa^i\in\Co^\alpha(\B^n;T^*\R^n)$, $\tau_\kappa^i\in\Co^\beta(\B^n;T^*\R^n)$ (by \eqref{Eqn::Rough1Form::ScalingProof3}) and $\supp(\lambda_\kappa^i-dx^i)\subseteq\supp\chi_1\Subset\frac12\B^n$. And since $\chi_1\big|_{\frac13\B^n}\equiv1$ and $\phi_\kappa^*dx=\kappa dx$,
\begin{equation}\label{Eqn::Rough1Form::Scaling::FormEqualProof}
    \lambda_\kappa^i\big|_{\frac13\B^n}=dx^i+\tfrac1\kappa\cdot\phi_\kappa^*(\xi^i-dx^i)\big|_{\frac13\B^n}=\tfrac1\kappa(\phi_\kappa^*\xi^i)\big|_{\frac13\B^n}.
\end{equation}

We write $\xi^i$ and $\rho^i$, $i=1,\dots,n$ as
\begin{equation}\label{Eqn::Rough1Form::Scaling::ThetaRhoAB}
    \xi^i=dx^i+\sum_{j=1}^na^i_j(x)dx^j,\quad\rho^i=\sum_{j=1}^nb^i_j(x)dx^j,\quad\text{where }a^i_j\in\Co^\alpha(\mu_0\B^n),\ b^i_j\in\Co^\beta(\mu_0\B^n).
\end{equation}
By assumption $\xi^i\big|_0=dx^i\big|_0$ and $\rho^i\big|_0=0$ for $i=1,\dots,n$, so $a^i_j(0)=b^i_j(0)=0$ for all $1\le i,j\le n$. And we have
\begin{equation}\label{Eqn::Rough1Form::ScalingProofLambdaCoordinate}
    \lambda_\kappa^i=dx^i+\sum_{j=1}^n\chi_1(x)a^i_j(\kappa x)dx^j,\quad\tau_\kappa^i=\sum_{j=1}^n\Big(\chi_1(x) b^i_j(\kappa x)dx^j+\Gc\ast\codiff\big(a^i_j(\kappa x)d\chi_1\wedge dx^j\big)\Big),\quad i=1,\dots,n.
\end{equation}

Since $\phi_\kappa(\frac12\B^n)\subseteq\frac{\mu_0}2\B^n$ and $\supp\chi_1\Subset\frac12\B^n$, by condition \ref{Item::Rough1Form::Scaling::RhoFormEqual}  we have $\chi_1\cdot\phi_\kappa^*d\rho^i=\chi_1\cdot\phi_\kappa^*d\xi^i$ for $i=1,\dots,n$. It follows that  $d\lambda_\kappa^i=d\tau_\kappa^i$ for $i=1,\dots,n$; indeed,
\begin{equation}\label{Eqn::Rough1Form::ScalingTmp2}
    d\lambda_\kappa^i-d\tau_\kappa^i
    =\tfrac1\kappa d\chi_1\wedge \phi_\kappa^*(\xi^i-dx^i)
    +\tfrac1\kappa\chi_1\cdot\phi_\kappa^*d\xi^i
    -\tfrac1\kappa d\chi_1\wedge\phi_\kappa^*(\xi^i-dx^i)
    -\tfrac1\kappa\chi_1\cdot\phi_\kappa^*d\rho^i
    =\tfrac1\kappa \chi_1\cdot\phi_\kappa^*(d\xi^i-d\rho^i)=0.
\end{equation}

Applying Lemmas \ref{Lem::Hold::Product} and \ref{Lem::Hold::ScalingLem} to $a^i_j$  we have
\begin{equation}\label{Eqn::Rough1Form::Scaling::BddA}
    \|\chi_1\cdot a^i_j(\kappa\cdot)\|_{\Co^\alpha(\B^n)}\le2\|\chi_1\|_{\Co^\alpha}\|a^i_j(\kappa\cdot)\|_{\Co^\alpha(\B^n)}\lesssim_{\alpha,\mu_0}\kappa^{\min(\alpha,\frac12)}\|\chi_1\|_{\Co^\alpha}\|a^i_j\|_{\Co^\alpha(\mu_0\B^n)},\quad\forall\kappa\in(0,\mu_0].
\end{equation}

Using \eqref{Eqn::Rough1Form::ScalingProofLambdaCoordinate} and \eqref{Eqn::Rough1Form::Scaling::BddA}, we deduce that

\begin{equation}\label{Eqn::Rough1Form::ScalingProof2.1}
\sum_{i=1}^n\|\lambda_\kappa^i-dx^i\|_{\Co^\alpha}\le\sum_{i,j=1}^n\|\chi_1(x)a^i_j(\kappa x)\|_{\Co^\alpha}\lesssim_{\alpha,\mu_0}\kappa^{\min(\alpha,\frac12)}\|\chi_1\|_{\Co^\alpha}\sum_{i,j=1}^n\| a^i_j\|_{\Co^\alpha(\mu_0\B^n)}.
\end{equation}

By \eqref{Eqn::Rough1Form::ScalingTmp2} we have $d\lambda_\kappa^i=d\tau_\kappa^i$, and therefore 
\begin{equation}\label{Eqn::Rough1Form::ScalingTmpDlambdaDkappa}
    \|d\lambda_\kappa^i\|_{\Co^{\beta-1}(\B^n;\wedge^2T^*\B^n)}=\|d\tau_\kappa^i\|_{\Co^{\beta-1}(\B^n;\wedge^2T^*\B^n)}\lesssim_\beta\|\tau_\kappa^i\|_{\Co^{\beta}(\B^n;T^*\B^n)},\quad i=1,\dots,n.
\end{equation}

Applying Lemma \ref{Lem::Hold::ScalingLem} to $ a^i_j$ and $ b^i_j$ we get that for $0<\kappa<\mu_0$,
\begin{gather}\label{Eqn::Rough1Form::ScalingTmp2.1}
    \|\chi_1\cdot b^i_j(\kappa x)\|_{\Co^\beta(\B^n)}\lesssim_\beta\|\chi_1\|_{\Co^\beta}\|b^i_j(\kappa x)\|_{\Co^\beta(\B^n)}\lesssim_\beta\kappa^{\min(\beta,\frac12)}\|\chi_1\|_{\Co^\beta}\|b^i_j\|_{\Co^\beta}\le \kappa^{\min(\alpha,\frac12)}\|\chi_1\|_{\Co^\beta}\|b^i_j\|_{\Co^\beta(\mu_0\B^n)}.
    \\\label{Eqn::Rough1Form::ScalingTmp2.2}
    \|a^i_j(\kappa x)d\chi_1\|_{\Co^\alpha(\B^n;T^*\R^n)}\lesssim_\alpha\|a^i_j(\kappa x)\|_{\Co^\alpha(\B^n)}\|\chi_1\|_{\Co^{\alpha+1}}\lesssim_\alpha\kappa^{\min(\alpha,\frac12)}\|\chi_1\|_{\Co^{\alpha+1}}\|a^i_j\|_{\Co^\alpha(\mu_0\B^n)}.
\end{gather}

Letting $\omega=a_j^i(\kappa x)d\chi_1\wedge dx^j\in\Co^{\alpha}(\B^n;T^*\B^n)\subseteq\Co^{\beta-1}(\B^n;T^*\B^n)$ and $\mu=1$ in \eqref{Eqn::Rough1Form::ScalingProof3}, we see that  
\begin{equation}
\label{Eqn::Rough1Form::ScalingProof2.2}
\begin{aligned}
    &\sum_{i=1}^n\|d\lambda_\kappa^i\|_{\Co^{\beta-1}(\B^n;\wedge^2T^*\B^n)}=\sum_{i=1}^n\|d\tau_\kappa^i\|_{\Co^{\beta-1}(\B^n;\wedge^2T^*\B^n)}\lesssim_\beta\sum_{i=1}^n\|\tau_\kappa^i\|_{\Co^\beta(\B^n;T^*\B^n)}&\text{by \eqref{Eqn::Rough1Form::ScalingTmpDlambdaDkappa}}
    \\
    &\le\sum_{i,j=1}^n\big(\|\chi_1(x)b_j^i(\kappa x)dx^j\|_{\Co^\beta}+\big\|\Gc\ast\codiff\big(a^i_j(\kappa x)d\chi_1\wedge dx^j\big)\big\|_{\Co^\beta}\big)&\text{by \eqref{Eqn::Rough1Form::ScalingProofLambdaCoordinate}}
    \\
    &\lesssim_\beta\sum_{i,j=1}^n\big(\|\chi_1(x)b_j^i(\kappa x)\|_{\Co^\beta(\B^n)}+\|a^i_j(\kappa x)d\chi_1\wedge dx^j\big\|_{\Co^{\beta-1}(\B^n;\wedge^2T^*\B^n)}\big)&\text{by \eqref{Eqn::Rough1Form::ScalingProof3}}
    \\
    &\lesssim_{\alpha,\beta}\|\chi_1\|_{\Co^{\alpha+1}}\sum_{i,j=1}^n\big(\|b_j^i(\kappa x)\|_{\Co^\beta(\B^n)}+\|a^i_j(\kappa x)d\chi_1\big\|_{\Co^{\alpha}(\B^n;T^*\B^n)}\big)&\text{since }\alpha\ge\beta-1
    \\
    &\lesssim_{\alpha,\beta}\kappa^{\min(\alpha,\frac12)}\|\chi_1\|_{\Co^{\alpha+1}}\sum_{i,j=1}^n\Big(\| a^i_j\|_{\Co^\alpha(\mu_0\B^n)}+\|b^i_j\|_{\Co^\beta(\mu_0\B^n)}\Big)&\text{by \eqref{Eqn::Rough1Form::ScalingTmp2.1} and \eqref{Eqn::Rough1Form::ScalingTmp2.2}}.
\end{aligned}\end{equation}

Note that $\chi_1$ is a fixed cut-off function whose $\Co^\alpha$ and $\Co^{\alpha+1}$-norms depend only on $n,\alpha$. So combining \eqref{Eqn::Rough1Form::ScalingProof2.1} and \eqref{Eqn::Rough1Form::ScalingProof2.2} we have
\begin{equation}\label{Eqn::Rough1Form::Scaling::FinalBdd1}
    \sum_{i=1}^n(\|\lambda_\kappa^i-dx^i\|_{\Co^\alpha}+\|d\lambda_\kappa^i\|_{\Co^{\beta-1}})\lesssim_{\alpha,\beta,\mu_0}\kappa^{\min(\alpha,\frac12)}\sum_{i,j=1}^n\big(\|a^i_j\|_{\Co^\alpha(\mu_0\B^n;T^*\R^n)}+\|b^i_j\|_{\Co^{\beta}(\mu_0\B^n;\wedge^2T^*\R^n)}\big).
\end{equation}

By \eqref{Eqn::Rough1Form::Scaling::ThetaRhoAB} we have $\sum_{i,j=1}^n\|a^i_j\|_{\Co^\alpha}\lesssim 1+\sum_{i=1}^n\|\xi^i\|_{\Co^\alpha}$ and $\sum_{i,j=1}^n\|b^i_j\|_{\Co^\beta}\lesssim\sum_{i=1}^n\|\rho^i\|_{\Co^\beta}$. And combining \eqref{Eqn::Rough1Form::Scaling::FinalBdd1} with \eqref{Eqn::Rough1Form::Scaling::RhoFormEst}, we can find a $C_1=C_1(n,\alpha,\beta,\mu_0)>0$ that does not depend on the other quantities, such that
\begin{equation}\label{Eqn::Rough1Form::Scaling::FinalBdd2}
    \sum_{i=1}^n\|\lambda_\kappa^i-dx^i\|_{\Co^\alpha}+\|d\lambda_\kappa^i\|_{\Co^{\beta-1}}\le C_1\cdot\kappa^{\min(\alpha,\frac12)}\sum_{i=1}^n\Big(1+\|\xi^i\|_{\Co^\alpha(\mu_0\B^n;T^*\R^n)}+\|d\xi^i\|_{\Co^{\beta-1}(\mu_0\B^n;\wedge^2T^*\R^n)}\Big),\quad\forall \kappa\in(0,\mu_0].
\end{equation}

Now applying assumption \eqref{Eqn::Rough1Form::Scaling::AssumptionOnTheta} to \eqref{Eqn::Rough1Form::Scaling::FinalBdd2} we have $$\sum_{i=1}^n\|\lambda_\kappa^i-dx^i\|_{\Co^\alpha}+\|d\lambda_\kappa^i\|_{\Co^{\beta-1}}\le \kappa^{\min(\alpha,\frac12)}C_1\cdot (M+n).$$ 

Since $\|(\langle\lambda_\kappa^i-dx^i,\Coorvec{x^j}\rangle)_{n\times n}\|_{C^0(\B^n;\R^{n\times n})}\lesssim_\alpha\sum_{i=1}^n\|\lambda_\kappa^i-dx^i\|_{\Co^\alpha(\B^n)}$, we can find a $\tilde c'=\tilde c'(n,\alpha)>0$ such that 
$$\sum_{i=1}^n\|\lambda_\kappa^i-dx^i\|_{\Co^\alpha}\le \tilde c'\quad\text{implies}\quad\textstyle\|(\langle\lambda_\kappa^i-dx^i,\Coorvec{x^j}\rangle)_{n\times n}\|_{C^0(\B^n;\R^{n\times n})}\le\frac12.$$ In particular $\mleft(\langle\lambda_\kappa^i,\Coorvec{x^j}\rangle\mright)_{n\times n}=I+\mleft(\langle\lambda_\kappa^i-dx^i,\Coorvec{x^j}\rangle\mright)_{n\times n}$ is invertible at every point in $\B^n$, which means $(\lambda_\kappa^1,\dots,\lambda_\kappa^n)$  span the tangent space at every point in $\B^n$.

We take $\kappa_0=\kappa_0(n,\alpha,\beta,\mu_0,M,\tilde c)>0$ such that
$$0<\kappa_0<\mu_0\quad\text{and}\quad\kappa_0^{\min(\alpha,\frac12)}C_1\cdot (M+n)\le \min(\tilde c,\tilde c').$$ 

Take $\lambda^i=\lambda_{\kappa_0}^i$ for $i=1,\dots,n$. We have $\sum_{i=1}^n\|\lambda^i-dx^i\|_{\Co^\alpha}+\|d\lambda^i\|_{\Co^{\beta-1}}\le\min(\tilde c,\tilde c')$. {By our assumption on $\tilde c'$,} $\lambda^1,\dots,\lambda^n$ span the tangent space at every point in $\B^n$. Note that by \eqref{Eqn::Rough1Form::ScalingProofLambdaCoordinate} $\supp(\lambda_{\kappa_0}^i-dx^i)\subseteq\supp\chi_1\Subset\frac12\B^n$. This shows conclusion \ref{Item::Rough1Form::Scaling::Est} is satisfied.

By \eqref{Eqn::Rough1Form::Scaling::FormEqualProof} we get the conclusion \ref{Item::Rough1Form::Scaling::FormEqual}, finishing the proof.
\end{proof}


We can give a sufficient and necessary condition for some 1-forms to have better regularity in some special coordinate system:
\begin{cor}\label{Cor::Rough1Form::BasicCaseforMainThm}
Let $\alpha>0$ and let $\beta\in[\alpha,\alpha+1]$ be such that $\alpha+\beta>1$. Let $\lambda^1,\dots,\lambda^n$ be $\Co^\alpha$ 1-forms on a smooth manifold $\Mf$ of dimension $n$, which span the cotangent space at every point. Then following are equivalent:
\begin{enumerate}[parsep=-0.3ex,label=(\arabic*)]
    \item\label{Item::Rough1Form::BasicCaseforMainThm::a} For every $p\in\Mf$  there exist a neighborhood $U\subseteq\Mf$ of $p$ and a $\Co^{\alpha+1}_\loc$-diffeomorphism $\Phi:\B^n\xrightarrow{\sim}U\subseteq \Mf$ 
    such that $\Phi(0)=p$ and $\Phi^*\lambda^1,\dots,\Phi^*\lambda^n\in\Co^\beta(\B^n;T^*\B^n)$.
    \item\label{Item::Rough1Form::BasicCaseforMainThm::b} $d\lambda^1,\dots,d\lambda^n\in\Co^{\beta-1}_\loc(\Mf;\wedge^2T^*\Mf)$. (see Definition \ref{Defn::DisInv::DefFunVF} \ref{Item::DisInv::DefFunVF::Form}).
\end{enumerate}
\end{cor}

\begin{remark}
\begin{enumerate}[parsep=-0.3ex,label=(\roman*)]
    \item When $X_1,\dots,X_n$ are $\Co^\alpha$ vector fields that span the cotangent space at every point, we can apply Corollary \ref{Cor::Rough1Form::BasicCaseforMainThm} by taking $(\lambda^1,\dots,\lambda^n)$ to be the dual basis of $(X_1,\dots,X_n)$, and we see that $\Phi^*\lambda^1,\dots,\Phi^*\lambda^n$ are $\Co^\beta$ if and only if $\Phi^*X_1,\dots,\Phi^*X_n$ are $\Co^\beta$.
    \item When $\beta>\alpha+1$, the result is still true if we replace $d\lambda^1,\dots,d\lambda^n\in\Co^{\beta-1}_\loc$ by $d\lambda^1,\dots,d\lambda^n\in\Co^{\beta-1}_{X,\loc}$ in \ref{Item::Rough1Form::BasicCaseforMainThm::b}. See \cite[Definition 2.15]{StreetYaoVectorFields}.
    \item The assumption $\alpha+\beta>1$ can be removed. When $\alpha+\beta\le1 $, in particular $0<\alpha\le\beta<1$, we cannot pullback a $\Co^{\beta-1}$-form under a $\Co^{\alpha+1}$-map. However
    the statement $d\lambda^1,\dots,d\lambda^n\in\Co^{\beta-1}_\loc(\Mf;\wedge^2T^*\Mf)$ still make sense because for a 1-form $\xi\in\Co^{0+}$ and a $\Co^{\alpha+1}$-diffeomorphism $\phi$, we have $d\xi\in\Co^{\beta-1}_\loc$ if and only if $d(\phi^*\xi)\in\Co^{\beta-1}_\loc$, see \cite[Proposition 2.8]{StreetYaoVectorFields}.
\end{enumerate}
\end{remark}

\begin{proof}
By passing to a local coordinate system, we can assume $\Mf$ to be an open subset of $\R^n$. 

\ref{Item::Rough1Form::BasicCaseforMainThm::a}$\Rightarrow$\ref{Item::Rough1Form::BasicCaseforMainThm::b}: For such $\Phi$, we have $d(\Phi^*\lambda^i)\in\Co^{\beta-1}_\loc(\B^n;T^*\B^n)$ for $i=1,\dots,n$. By Lemma \ref{Lem::Hold::Product} \ref{Item::Hold::Product::Hold1} along with direct computation (see \cite[Lemma 4.4]{StreetYaoVectorFields}) we see that $\Phi^*d\lambda^i=d(\Phi^*\lambda^i)$.
So $\Phi^\Inv:U\subseteq\Mf\to\R^n$ is the desired coordinate chart that shows $d\lambda^1,\dots,d\lambda^n\in\Co^{\beta-1}_\loc$.

\ref{Item::Rough1Form::BasicCaseforMainThm::b}$\Rightarrow$\ref{Item::Rough1Form::BasicCaseforMainThm::a}: Let $p\in\Mf$. 

By passing to local coordinate system and applying an invertible linear transformation we can find a $\mu_0>0$ and a $\Co^{\alpha+1}$-coordinate chart $F_0:U_0\subseteq\Mf\xrightarrow{\sim} \mu_0\B^n_x$ such that
\begin{itemize}[nolistsep]
    \item $F_0(p)=0$.
    \item $((F_0)_*\lambda^i)\big|_0=dx^i\big|_0$ for $i=1,\dots,n$.
\end{itemize}

Take $c>0$ be the small constant in Theorem \ref{KeyThm::Rough1Form}. By Proposition \ref{Prop::Rough1Form::Scaling} with $\tilde c=c$, we can find a $\kappa_0\in(0,\mu_0]$ and 1-forms $\tilde \lambda^1,\dots,\tilde\lambda^n\in\Co^\alpha(\B^n;T^*\B^n)$ such that for scaling map $\phi_{\kappa_0}:\B^n\to\mu_0\B^n$, $\phi_{\kappa_0}(x)=\kappa_0x$, we have
\begin{enumerate}[parsep=-0.3ex,label=(\alph*)]
    \item $\tilde\lambda^1,\dots,\tilde\lambda^n$ span the cotangent space at every point.
    \item $\supp(\tilde\lambda^i-dx^i)\subsetneq\frac12\B^n$ for $i=1,\dots,n$.
    \item\label{Item::Rough1Form::BasicCaseforMainThm::Tmp} $(\tilde\lambda^1,\dots,\tilde \lambda^n)\big|_{\frac13\B^n}=\frac1{\kappa_0}\cdot((F_0^{-1}\circ\phi_{\kappa_0})^*\lambda^1,\dots,(F_0^{-1}\circ\phi_{\kappa_0})^*\lambda^n)\big|_{\frac13\B^n}$.
    \item $\sum_{i=1}^n\|\tilde\lambda^i-dx^i\|_{\Co^\alpha}+\|d\tilde\lambda^i\|_{\Co^{\beta-1}}<c$.
\end{enumerate}
We set $F_1:=\phi_{\kappa_0}^\Inv$. 

Applying Theorem \ref{KeyThm::Rough1Form} to $\tilde\lambda^1,\dots,\tilde\lambda^n$ we obtain a $\Co^{\alpha+1}$-chart $F_2:\B^n\xrightarrow{\sim}\B^n$ such that $F_2(\frac13\B^n)\supseteq B^n(F_2(0),\frac16)$, $\|F_2-\id\|_{\Co^{\alpha+1}}<c$ and $(F_2)_*\tilde\lambda^1,\dots,(F_2)_*\tilde\lambda^n\in\Co^\beta(\B^n;T^*\B^n)$.

By \ref{Item::Rough1Form::BasicCaseforMainThm::Tmp}, $\left((F_2\circ F_1\circ F_0)_*\lambda^i\right)\big|_{F_2(\frac13\B^n)}={\tfrac1{\kappa_0}}\cdot(F_2)_*\tilde\lambda^i\big|_{F_2(\frac13\B^n)}\in\Co^\beta$ for $i=1,\dots,n$. We can take an
affine linear transformation
$F_3:\R^n\to\R^n$ such that $F_3(F_2(0))=0$ and $F_3(B^n(F_2(0),\frac16))\supseteq\B^n$. 


{Now we have $F_0(U_0)\supseteq\mu_0\B^n$, $F_1(\mu_0\B^n)\supseteq\B^n$, $F_2(\frac13\B^n)\supseteq B^n(F_2(0),\frac16)$ and $F_3(B^n(F_2(0),\frac16)\supseteq\B^n$.} Take $\Phi:=(F_3\circ F_2\circ F_1\circ F_0)^\Inv:\B^n\to\Mf$. Since $F_0,F_1,F_2,F_3$ are all $\Co^{\alpha+1}$-diffeomorphism onto their images, we know know $\Phi:\B^n\xrightarrow{\sim}\Phi(\B^n)\subseteq\Mf$ is a $\Co^{\alpha+1}$ diffeomorphism. Moreover, we have $\Phi(0)=p$ because $F_0(p)=0$, $F_1(0)=0$, $F_3(F_2(0))=0$. Thus, $\Phi$ is the diffeomorphism we desire with $U:=\Phi(\B^n)$, completing the proof.
\end{proof}

\subsection{Remarks in harmonic coordinates}\label{Section::Rough1FormRmkCoor}
Given a non-smooth Riemannian metric $g$ on a manifold, $\Mf$, DeTurck and Kazdan showed that $g$ has optimal regularity
in harmonic coordinates \cite[Lemma 1.2]{DeTurckKazdan} (in the Zygmund-H\"older sense),
but may not have optimal regularity in geodesic normal coordinates \cite[Example 2.3]{DeTurckKazdan} (in fact, the regularity
of $g$ in geodesic normal coordinates may be two derivatives worse than the regularity in harmonic coordinates).  


We present analogous results for vector fields.
 Let $X_1,\dots,X_n$ be $C^1_\loc$-vector fields on a $C^2$-manifold $\Mf$ that form a local basis for the tangent space at every point. In Section \ref{Section::CanCoor-1Reg} we give an example showing that the canonical coordinate system with respect to $X_1,\dots,X_n$ may not have optimal regularity
 
 In Section \ref{Section::Rough1FormOV} we defined a Riemannian metric $g=\sum_{i=1}^n\lambda^i\cdot\lambda^i$ where $(\lambda^1,\dots,\lambda^n)$ is the dual basis of $(X_1,\dots,X_n)$ (see Remark \ref{Rmk::Rough1Form::RmkRiemMetric}). 
 With respect to this metric, $X_1,\dots,X_n$ form an orthonormal basis at every point. Since $X_1,\dots,X_n\in C^1$, 
 we can talk about the metric Laplacian $\Delta_g$ with respect to $g$.

\begin{prop}\label{Prop::Rough1Form::RmkCoor::HarmOptimal}
In harmonic coordinates with respect to $g$, $X_1,\dots,X_n$ have optimal regularity. 

More precisely, let $X_1,\dots,X_n$ and $g$ be as above, and let $\beta>1$. Suppose there is a $\Co^{\beta+1}$-atlas $\Af$ which is compatible with the $C^2$-atlas of $\Mf$, such that $X_1,\dots,X_n$ are $\Co^\beta$ on $\Af$, let $\psi:U\subseteq\Mf\to V\subseteq\R^n$ be a harmonic coordinate chart\footnote{Such a harmonic coordinate chart $\psi$ always exists locally when $\beta>1$, see also \cite[Lemma 1.2]{DeTurckKazdan}.}, then $\psi_*X_1,\dots,\psi_*X_n\in\Co^\beta_\loc(V;\R^n)$.
\end{prop}
\begin{proof}
It suffices to show that $\psi$ is a $\Co_{\loc}^{\beta+1}$-map with respect to $(\Mf,\Af)$. Once this is done, using that $X_1,\dots,X_n$ are $\Co^\beta_\loc$ with respect to $(\Mf,\Af)$, we get that $\psi_*X_1,\dots,\psi_*X_n\in\Co^\beta_\loc$.

Since the statement $\psi\in\Co_{\loc}^{\beta+1}(U;V)$ is local, 
we may without loss of generality, shrink $U$.
By doing so, the hypotheses of the proposition
imply that there is a $\Co^{\beta+1}$-coordinate chart $x=(x^1,\dots,x^n):U\subseteq\Mf\to\R^n$ on $U$ (respect to $\Af$). In this coordinate chart we can write $\Delta_g=\frac1{\sqrt{\det g}}\sum_{i,j=1}^n\Coorvec{x^i}(\sqrt{\det g}g^{ij}\Coorvec{x^j})$ where $g^{ij}$ and $\sqrt{\det g}$ are as in \eqref{Eqn::Rough1Form::Riemannianmetric2}.

By assumption $\Delta_g \psi^k=0$ for $k=1,\dots,n$. Note that on $(x^1,\dots,x^n)$, $\Delta_g$ is a second order divergent form elliptic operator with $\Co^\beta$-coefficients. By a classical elliptic estimates (for example, \cite[Proposition 4.1]{TaylorPDE3}) we have that $\psi^k$ are all $\Co^{\beta+1}_\loc$ with respect to $\Af$, completing the proof.
\end{proof}
\begin{remark}


In fact the coordinate chart we construct in Proposition \ref{Prop::Rough1Form::ExistPDE} (also see \eqref{Eqn::Rough1Form::PDEforR}) is closely related to harmonic coordinates.
\end{remark}
\begin{remark}
While Proposition \ref{Prop::Rough1Form::RmkCoor::HarmOptimal} shows that $X_1,\dots,X_n$ have optimal regularity with respect to harmonic coordinates, this fact along does not give a practical test for what the optimal regularity is. Corollary \ref{Cor::Rough1Form::BasicCaseforMainThm}, on the other hand, provides such a test.
\end{remark}

\begin{remark}
Proposition \ref{Prop::Rough1Form::RmkCoor::HarmOptimal} shows that harmonic coordinates induces a $\Co^{\beta+1}$-atlas with respect to which $X_1,\dots,X_n$ are $\Co^\beta_\loc$. It is possible that the harmonic coordinates induces some $\Co^{\gamma+1}$-atlas for some $\gamma>\beta$ while $X_1,\dots,X_n$ are only $\Co^\beta$ with respect to this atlas; see Example \ref{Exmp::Rough1Form::RmkCoor::CanonicalMighNotBeOptimal} below.
\end{remark}

\begin{example}\label{Exmp::Rough1Form::RmkCoor::CanonicalMighNotBeOptimal}
Endow $\R^2$ with standard coordinates $(x,y)$, and let $\theta\in C^1(\R^2)$ be a function which is not smooth. Set $X:=\cos(\theta(x,y))\Coorvec x+\sin(\theta(x,y))\Coorvec y$ and $Y:=-\sin(\theta(x,y))\Coorvec x+\cos(\theta(x,y))\Coorvec x$. The corresponding metric is $g=dx^2+dy^2$ since $X,Y$ form an orthonormal basis with respect to the standard Euclidean metric, thus $\Delta_g=\Delta=\partial_x^2+\partial_y^2$.

Therefore the singleton $\{(x,y):\R^2\to\R^2\}$ is an atlas of harmonic coordinates for $\R^2$, and since harmonic functions are real-analytic, we know the collection of harmonic coordinates with respect to $\Delta_g$ defines an real-analytic structure for $\R^2$ (which coincides with the standard real-analytic structure). Even though the differential structure induced by the harmonic coordinates is real analytic, $X$ and $Y$ cannot be smooth under any coordinate system (since they are not smooth with respect to these
harmonic coordinates).
\end{example}





\chapter{The Theorems and The Proofs}
\section{The Real Frobenius Theorems}\label{Section::RealFro}
\subsection{History remark for rough involutivity}\label{Section::RealFro::HisRmk}
For vector fields $X$ and $Y$ which are not $C^1$, the commutator $[X,Y]$ may not be defined pointwise. Because of this, any Frobenius type theorem for subbundles of less than $C^1$-regularity requires a suitable notion of involutivity. This leads to the following:

\begin{ques}\label{Ques::RealFro::Intro::Sec}
Let $\V$ be a tangent subbundle which is not $C^1$, and let vector fields $X,Y$ be sections of $\V$.
\begin{itemize}[parsep=-0.3ex]
    \item What is the good definition of the Lie bracket $[X,Y]$?
    \item With this definition how can we say whether $[X,Y]$ is a section of $\V$? 
\end{itemize}
\end{ques}

When $X$ and $Y$ are Lipschitz vector fields, the Lie bracket $[X,Y]$ makes sense as a vector field with $L^\infty$-coefficients. And we can talk about sections and involutivity almost everywhere. In \cite{Simic}, Simi\'c gives the following definition:
\begin{defn}[Almost everywhere characterization]\label{Defn::RealFro::Intro::AEsection}Let $\Mf$ be a $C^{1,1}$-manifold, let  $\V$ be a Lipschitz subbundle of $T\Mf$, and let $X,Y$ be two Lipschitz sections of $\V$.
We say $[X,Y]$ is an \textit{almost everywhere section} of $\V$, if as a vector field with $L^\infty$-coefficients, $[X,Y](p)\in\V_p$ holds for almost every $p\in \Mf$.

We say $\V$ is \textit{almost everywhere involutive}, if for any Lipschitz sections $X,Y$ of $\V$, the Lie bracket $[X,Y]$ is an almost everywhere section of $\V$.
\end{defn}

 Later Rampazzo \cite{Rampazzo} gave another characterizations using set-valued Lie brackets:
\begin{defn}[Set valued characterization]\label{Defn::RealFro::Intro::SetInv}
Let $\Mf$ be a $C^{1,1}$-manifold, let  $\V$ be a Lipschitz subbundle of $T\Mf$, and let $X,Y$ be two Lipschitz sections of $\V$. The \textit{set-valued Lie bracket} $[X,Y]_{set}\subseteq\V$ is defined as the convex hall $[X,Y]_{set}(p):=\operatorname{conv}([X,Y]_{set}^+(p))\subseteq T_p\Mf$, $p\in \Mf$, where
$$[X,Y]_{set}^+:=\overline{\{(q,[X,Y](q))\in T\Mf:X,Y\text{ are differentiable at }q\}}\subseteq T\Mf.$$
And we say $\V$ is \textit{set-valuedly involutive}, if $[X,Y]_{set}\subseteq\V$ for every Lipschitz sections $X,Y$ of $\V$.
\end{defn}

Definitions \ref{Defn::RealFro::Intro::AEsection} and \ref{Defn::RealFro::Intro::SetInv} are equivalent, see \cite[Section 2]{Rampazzo} and \cite[Theorem 4.11 (I) and (III)]{Rampazzo}.

An advantage of set valued bracket is that it is defined pointwise everywhere.
Note that a Lipschitz function is differentiable almost everywhere, so $[X,Y]_{set}(p)\neq\varnothing$ for all $p\in M$.

\begin{prop}[Lipschitz Frobenius theorem, \cite{Simic,Rampazzo}]
Let $\Mf$ be a $n$-dimensional $C^{1,1}$ manifold, and let $\V\le T\Mf$ be a Lipschitz tangent subbundle of rank $r$. If $\V$ is either almost everywhere involutive or set-valuedly involutive, then for any $p\in \Mf$, there is a neighborhood $\Omega\subseteq\R^r_t\times\R^{n-r}_s$ of $0$ and a bi-Lipschitz parameterization $\Phi(t,s):\Omega\to \Mf$ satisfying the following:
\begin{enumerate}[nolistsep,label=(\roman*)]
    \item $\Phi(0)=p$.
    \item \label{LipFroConse2} $\frac{\partial\Phi}{\partial t^1},\dots,\frac{\partial\Phi}{\partial t^r}:\Omega\to T\Mf$ are Lipschitz continuous maps.
    \item For any $(t,s)\in \Omega$, $\frac{\partial\Phi}{\partial t^1}(t,s),\dots,\frac{\partial\Phi}{\partial t^r}(t,s)\in T_{\Phi(t,s)}\Mf$ form a basis of the linear subspace $\V_{\Phi(t,s)}$.
\end{enumerate}
\end{prop}

Both \cite{Simic} and \cite{Rampazzo} only deal with Lipschitz vector fields and Lipschitz subbundles. What happen if they are non-Lipschitz? Indeed, $[X,Y]$ is not necessarily $L^1_\loc$ so we cannot use an almost everywhere characterization. A non-Lipschitz section can be nowhere differentiable so we cannot use set-valued bracket either. 

To resolve this problem, we use generalized functions. When the vector fields $X,Y$ are log-Lipschitz, or more generally when $X,Y$ are H\"older $\Co^{\alpha}$ for some $\frac12<\alpha<1$, by Corollary \ref{Cor::Hold::[X,Y]WellDef}, $[X,Y]$ is a vector field with $\Co^{\alpha-1}$-coefficients, see Lemma \ref{Lem::Hold::PushForwardFuncSpaces} and Definition \ref{Defn::DisInv::DefFunVF}.

By imitating Definition \ref{Defn::Intro::DisInv} we can give a definition for the distributional sections of $\V$ as follows:
\begin{defn*}[Distributional sections of subbundles]
Let $\alpha>0$ and $\beta>-\alpha$. Let $\V$ be a $\Co^{\alpha}$ tangent subbundle of a manifold $\Mf$. We say a $\Co^{\beta}$-vector field $Z=\sum_{i=1}^nZ^i\Coorvec{x^i}$ is a distributional section of $\V$, if for any log-Lipschitz section $\theta$ of $\V^\bot$, written as $\theta=\sum_{j=1}^n\theta_jdx^j$, the generalized function\footnote{Using Lemma \ref{Lem::Hold::MultLoc} \ref{Item::Hold::MultLoc::WellDef}, the right hand side of \eqref{Eqn::RealFro::Intro::DisSec::DefEqn} makes sense as a distribution in $\Co^{\beta}_\loc(\Mf)$. Also see Definition \ref{Defn::DisInv::DefFunVF}.}
\begin{equation}\label{Eqn::RealFro::Intro::DisSec::DefEqn}
    \langle\theta,Z\rangle=\sum_{j=1}^nZ^j\theta_j\in\D'(\Mf),
\end{equation}
is identically zero in the sense of distributions.
\end{defn*}

One can see that this is the characterization in Proposition \ref{Prop::DisInv::CharSec} \ref{Item::DisInv::CharSec::Dual}. See also Remark \ref{Rmk::DisInv::RmkSecVB} \ref{Item::DisInv::RmkSecVB::DefSec}.

The characterization of distributional sections and distributional involutivity recover the above characterizations when $\V$ is Lipschitz. Indeed, for Lipschitz vector fields $X,Y$ that are sections of $\V$, $[X,Y]=0$ lays in $\V$ almost everywhere if and only if $\langle \theta,[X,Y]\rangle=0$, so Definitions \ref{Defn::Intro::DisInv} and \ref{Defn::RealFro::Intro::AEsection} are equivalent. 

\medskip
Recently in \cite{ContFro}, a version of involutivity was introduced using approximations by smooth functions. \cite{ContFro} works on differential forms (see \cite[Definition 1.12]{ContFro}). An equivalent characterization using vector fields is as follows:
\begin{defn}\label{Defn::RealFro::Intro::AsyInv}
Let $\V$ be a continuous rank-$r$ tangent subbundle over $\R^n$.
We say $\V$ is \textit{strongly asymptotically involutive}, if for any $p\in \R^n$ there is a neighborhood $U\subseteq \R^n$ of $p$, and sequence of $C^1$-vector fields $\{X_1^\nu,\dots,X^\nu_r\}_{\nu=1}^\infty$ on $U$, such that
\begin{enumerate}[parsep=-0.3ex,label=(\roman*)]
    \item For each $j=1,\dots,r$, $\{X_j^\nu\}_{\nu=1}^\infty$ uniformly converges to a continuous vector field $X_j$, such that for each $q\in U$, $X_1(q),\dots,X_r(q)\in \R^n$ form a linear basis of $\V_q$.
    \item There exists $\{c_{jk}^{l\nu}\}_{\nu=1}^\infty\subset C^0(\bar U)$ for $1\le j,k,l\le r$, such that 
    \begin{equation}\label{Eqn::RealFro::Intro::AsyInv1}\exists t_0>0,\
        \lim\limits_{\nu\to\infty}\bigg(\max\limits_{1\le j,k\le r}\Big\|[X_j^\nu,X_k^\nu]-\sum_{l=1}^rc_{jk}^{l\nu}X_l^\nu\Big\|_{C^0(U;\R^n)}\bigg)\cdot\exp\big(t_0\cdot\max\limits_{1\le l\le r}\|\nabla X_l^\nu\|_{C^0(U;\R^{n\times n})}\big)=0.
    \end{equation}
\end{enumerate}
\end{defn}

In \cite{ContFro} it is shown that if $\V$ is strongly asymptotically involutive, then locally there exists a topological parameterization $\Phi(t,s)\in C^1_t(C^0_s)$ such that $\Phi_*\Coorvec {t^1},\dots,\Phi_*\Coorvec{t^r}$ span $\V$. However, given a continuous subbundle, we do not know an algorithm to check whether such an approximation exists. In contrast, Definition \ref{Defn::Intro::DisInv} can be directly checked using coordinate systems.

In the log-Lipschitz setting we prove that distributional involutivity implies strongly asymptotically involutivity. See Corollary \ref{Cor::Hold::CorAsyInv} and Remark \ref{Rmk::Hold::CorAsyInv}.

\subsection{The Frobenius theorem for log-Lipschitz subbundles: proof of Theorem \ref{MainThm::LogFro}}\label{Section::RealFro::PfLogFro}

For application to the proof of Theorems \ref{MainThm::RoughFro1} and \ref{MainThm::RoughFro2}, we give a refinement of estimates in the results of Theorem \ref{MainThm::LogFro}.
\begin{keythm}\label{KeyThm::RealFro::ImprovedLogFro}
    Keep the assumptions of Theorem \ref{MainThm::LogFro}, by shrinking the neighborhood $\Omega$ that has the form $\Omega=\Omega'\times\Omega''\subseteq\R^r_t\times\R^{n-r}_s$ and write the inverse map $F:=\Phi^\Inv$ as $F=(F',F'')$ where $F':\Phi(\Omega)\to\R^r$ and $F'':\Phi(\Omega)\to\R^{n-r}$, then $F'\in C^{1,1}_\loc(\Phi(\Omega);\R^r)$. Moreover
    \begin{enumerate}[parsep=-0.3ex,label=(\roman*)]
        \item[\mylabel{Item::RealFro::ImprovedLogFro::BigLReg}{(5)}]For any $\eps>0$ there is an open subset $\Omega'_\eps\subseteq\Omega'$ of $0$, such that $\Phi\in\Co_\loc^{1+\LogL,1-\eps}(\Omega'_\eps,\Omega'';\Mf)$ and $\frac{\partial\Phi}{\partial t^1},\dots,\frac{\partial\Phi}{\partial t^r}\in\Co_\loc^{\LogL,1-\eps}(\Omega'_\eps,\Omega'';T\Mf)$. Moreover by taking $\Omega_\eps=\Omega'_\eps\times\Omega''$, and $F''\in\Co^{1-\eps}_\loc(\Phi(\Omega_\eps);\R^{n-r})$.
        \item[\mylabel{Item::RealFro::ImprovedLogFro::Zyg1Reg}{(5')}]Suppose $\V\in\Co^1$. Then for  the same $\Omega'_\eps\subseteq\Omega'$ in \ref{Item::RealFro::ImprovedLogFro::BigLReg}, automatically $\Phi\in\Co_\loc^{2,1-\eps}(\Omega'_\eps,\Omega'';\Mf)$ and $\frac{\partial\Phi}{\partial t^1},\dots,\frac{\partial\Phi}{\partial t^r}\in\Co_\loc^{1,1-\eps}(\Omega'_\eps,\Omega'';T\Mf)$.
        \item[\mylabel{Item::RealFro::ImprovedLogFro::LittleLReg}{(5'')}]Suppose $\V$ is little log-Lipschitz. Then automatically $\Phi\in\Co^{1+\LogL,1-}(\Omega',\Omega'';\Mf)$, $\frac{\partial\Phi}{\partial t^1},\dots,\frac{\partial\Phi}{\partial t^r}\in\Co^{\LogL,1-}(\Omega',\Omega'';T\Mf)$. and $F''\in\Co^{1-}_\loc(\Phi(\Omega);\R^{n-r})$.
        \item[\mylabel{Item::RealFro::ImprovedLogFro::C1Reg}{(5'\!'\!')}] Let $\alpha\in\{\Lip\}\cup(1,\infty]\cup\{k+\Lip,k+\LogL:k=1,2,3,\dots\}$. Suppose $\Mf$ is smooth and $\V\in\Co^\alpha$ for some $\alpha>1$, then automatically $\Phi\in\Co_\loc^{\alpha+1,\alpha}(\Omega',\Omega'';\Mf)$, $\frac{\partial\Phi}{\partial t^1},\dots,\frac{\partial\Phi}{\partial t^r}\in\Co_\loc^{\alpha}(\Omega;T\Mf)$, $F'\in \Co^\infty_\loc(\Phi(\Omega);\R^r)$ and $F''\in\Co^\alpha_\loc(\Phi(\Omega);\R^{n-r})$.
    \end{enumerate}
\end{keythm}
Note that \ref{Item::RealFro::ImprovedLogFro::BigLReg}, \ref{Item::RealFro::ImprovedLogFro::LittleLReg} and \ref{Item::RealFro::ImprovedLogFro::C1Reg} in Theorem \ref{KeyThm::RealFro::ImprovedLogFro} are stronger statement than \ref{Item::MainThm::LogFro::BigLReg}, \ref{Item::MainThm::LogFro::LittleLReg} and \ref{Item::MainThm::LogFro::C1Reg} in Theorem \ref{MainThm::LogFro}, respectively. 

\begin{proof}[Proof of Theorems \ref{MainThm::LogFro} and \ref{KeyThm::RealFro::ImprovedLogFro}]
By Remark \ref{Rmk::DisInv::RmkMfldObj} and by passing to a local coordinate chart of the manifold $\Mf$, we can assume that $\Mf\subseteq\R^n$ is an open subset with the base point $p=0\in\R^n$.

By Lemma \ref{Lem::ODE::GoodGen} we can find a neighborhood $U\subseteq \Mf$ of $p=0$, a linear coordinate system $(x,y)=(x^1,\dots,x^r,y^1,\dots,y^{n-r})$ for $\R^n$, and commutative log-Lipschitz vector fields $X_1,\dots,X_r$ defined on $U$ that  span $\V$ at every point in $U$ and has the form \eqref{Eqn::ODE::GoodGen::GoodGenFormula}, i.e. for some $a_j^k\in \Co^\LogL(U)$,
\begin{equation}\label{Eqn::RealFro::PfLogThm::FormXj}
    X_j=\Coorvec{x^j}+\sum_{k=1}^{n-r}a_j^k\Coorvec{y^k},\quad 1\le j\le r.
\end{equation}

Let $\Omega''\subseteq\R^{n-r}$ be an open neighborhood of $0\in\R^{n-r}$ such that $\{(0,s)\in\R^r\times\R^{n-r}:v\in\Omega''\}\Subset U$. We can take a small enough number $\tau_0>0$ such that the following map $\Phi(t,s):B^r(0,\tau_0)\times\Omega''\subset\R^r\times\R^{n-r}\to U$ is defined:
\begin{equation}\label{Eqn::RealFro::PfLogThm::DefofPhi}
    \Phi(t,s):=\Phi(t^1,\dots,t^r,s):=e^{t^1X_1}\dots e^{t^rX_r}(0,s),\quad |t|<\tau_0,\quad s\in\Omega''.
\end{equation}
Clearly $\Phi(0,0)=0$ which gives \ref{Item::MainThm::LogFro::0}.

\medskip
Applying Lemma \ref{Lem::ODE::CommMultFlow} with $S=\{0\}\times\Omega''$, by choosing $\tau_0>0$ small enough, $\Phi$ is well-defined, differentiable in $t$, and satisfies
\begin{equation*}
    \Coorvec{t^j}\Phi(t,s)=X_j(\Phi(t,s)),\quad1\le j\le r,\quad|t|<\tau_0,\quad s\in\Omega''.
\end{equation*}

Since $X_1,\dots,X_r$ are linearly independent and span $\V$ in the domain, we see that $\frac{\partial\Phi}{\partial t^1}(t,s),\dots,\frac{\partial\Phi}{\partial t^r}(t,s)$ spans $\V_{\Phi(t,s)}$ for $(t,s)$ in the domain $B^r(0,\tau_0)\times\Omega''$. This gives \ref{Item::MainThm::LogFro::Span}.

Next we prove that $\Phi$ is homeomorphism onto its image by constructing the inverse map. We still use $(x,y)$ as the standard coordinate system for $\Mf\subseteq\R^n$. Thus we have the natural coordinate maps $ x^j(q^1,\dots,q^n)=q^j$ for $1\le j\le r$ and $ y^k(q^1,\dots,q^n)=q^{r+k}$ for $1\le k\le n-r$ and $q\in\R^n$. Define $F'(q):=x(q)$ and $F''=(\lambda^1,\dots\lambda^{n-r}):\Phi(\Omega)\subseteq\R^n_q\to\R^{n-r}$ as 
\begin{equation}\label{Eqn::RealFro::PfLogThm::lambda}
    \lambda^j(q):=y^j(e^{-q^rX_r}\dots e^{-q^1X_1}(q)),\quad 1\le j\le n-r.
\end{equation}

We are going to show $\Phi\circ(F',F'')=(F',F'')\circ\Phi=\id$. This implies  $F=(F',F'')=\Phi^\Inv$ and thus $\Phi$ is homeomorphic onto its image.

Indeed, since $X_j$ are of the form \eqref{Eqn::RealFro::PfLogThm::FormXj}, we see that $\Phi$ can be written as
$$\Phi(t,s)=(t,\phi(t,s)),\quad\text{where }\phi:B^r(0,\tau_0)\times\Omega''\to \R^{n-r}.$$
Also we have
\begin{equation}\label{Eqn::RealFro::PfLogThm::Tmp}
    e^{-q^rX_r}\dots e^{-q^1X_1}(q)=(0, y(e^{-q^rX_r}\dots e^{-q^1X_1}(q))),\quad  x(e^{-q^rX_r}\dots e^{-q^1X_1}(q))=0,\quad q\in\Phi(B^r(0,\tau_0)\times\Omega'').
\end{equation}

Therefore for $(t,s)\in B^r(0,\tau_0)\times\Omega''$, using $\Phi(t,s)=(t,\phi(t,s))$ we have
\begin{equation}\label{Eqn::RealFro::PfLogThm::PfInv1}
    \begin{aligned}
    (F',F'')(\Phi(t,s))=&(x(t,\phi(t,s)),y(e^{-t^rX_r}\dots e^{-t^1X_1}(\Phi(t,s))))
    \\
    =&(t, y(e^{-t^rX_r}\dots e^{-t^1X_1}e^{t^1X_1}\dots e^{t^rX_r}(0,s)))
    \\
    =&(t, y(0,s))=(t,s).
\end{aligned}
\end{equation}
For $q\in\Phi(B^r(0,\tau_0)\times\Omega'')$, by \eqref{Eqn::RealFro::PfLogThm::Tmp} we have
\begin{equation}\label{Eqn::RealFro::PfLogThm::PfInv2}
    \begin{aligned}
    \Phi((F'(q),F''(q))=&\Phi(x(q),y(e^{-q^rX_r}\dots e^{-q^1X_1}(q)))
    \\
    =&e^{q^1X_1}\dots e^{q^rX_r}(0, y(e^{-q^rX_r}\dots e^{-q^1X_1}(q)))
    \\
    =&e^{q^1X_1}\dots e^{q^rX_r}(e^{-q^rX_r}\dots e^{-q^1X_1}(q))=q.
\end{aligned}
\end{equation}
We get $\Phi\circ(F',F'')=\id$ and $(F',F'')\circ\Phi=\id$. Since $B^r(0,\tau_0)\times\Omega''$ is open, we know $\Phi(B^r(0,\tau_0)\times\Omega'')=(x,\lambda)^{-1}(B^r(0,\tau_0)\times\Omega'')$ is also open.
This completes the proof of \ref{Item::MainThm::LogFro::Top} whenever $\Omega$ is an open neighborhood of $0$ satisfying $\Omega\subseteq B^r(0,\tau_0)\times\Omega''$.

Now $\Phi^\Inv=(F',F'')$ and it remains to prove the additional results in Theorem \ref{KeyThm::RealFro::ImprovedLogFro}, which implies  \ref{Item::MainThm::LogFro::BigLReg} and \ref{Item::MainThm::LogFro::LittleLReg} and \ref{Item::MainThm::LogFro::C1Reg}.

Since $F'(q)=x(q)$ arise from the coordinate chart, the regularity of $F'$ is the same as the regularity assumption for $\Mf$. 

By Lemma \ref{Lem::ODE::GoodGen} \ref{Item::ODE::GoodGen::Uniqueness}, since the log-Lipschitz vector fields $X_1,\dots,X_r$ in \eqref{Eqn::RealFro::PfLogThm::FormXj} are uniquely determined by $\V$ and the coordinate system $(x^1,\dots,x^r,y^1,\dots,y^{n-r})$, we see that
\begin{itemize}[parsep=-0.3ex]
    \item If $\V$ is $\Co^\alpha$ where $\alpha\in\{\logl,\Lip\}\cup[1,\infty]\cup\{k+\LogL,k+\Lip:k=1,2,3,\dots\}$, then $X_1,\dots,X_r$ are automatically $\Co^\alpha$ as well.
\end{itemize}

By Corollary \ref{Cor::ODE::MultFlowReg} we see that
\begin{itemize}[parsep=-0.3ex]
    \item For every $\eps>0$ there exists $\tau_\eps\in(0,\tau_0)$ such that $\Phi\in\Co^{1+\LogL,1-\frac\eps2}(B^r(0,\tau_\eps),\Omega'';U)$ and $F''\in\Co^{1-\frac\eps2}(\Phi(B^r(0,\tau_\eps)\times\Omega'');\R^{n-r})$. In particular by Lemma \ref{Lem::Hold::CompofMixHold} \ref{Item::Hold::CompofMixHold::Comp}, $\frac{\partial\Phi}{\partial t^j}=X_j\circ\Phi\in\Co^{\LogL,1-\eps}(B^r(0,\tau_\eps),\Omega'';U)$ for $1\le j\le r$.
    \item If $X_1,\dots,X_r\in\Co^1$, then automatically $\Phi\in\Co^2L^\infty(B^r(0,\tau_0),\Omega'';U)$ and $\frac{\partial\Phi}{\partial t^j}\in\Co^{1}L^\infty(B^r(0,\tau_0),\Omega'';U)$ for $1\le j\le r$. In particular $\Phi\in\Co^{2,1-\frac\eps2}(B^r(0,\tau_\eps),\Omega'';U)$ and $\frac{\partial\Phi}{\partial t^1},\dots,\frac{\partial\Phi}{\partial t^r}\in\Co^{1,1-\eps}(B^r(0,\tau_\eps),\Omega'';U)$.
    \item If $X_1,\dots,X_r\in\Co^\logl$, then automatically $\Phi\in\Co^{1+\LogL,1-}(B^r(0,\tau_0),\Omega'';U)$ and $F''\in\Co^{1-}(\Phi(B^r(0,\tau_0)\times\Omega'');\R^{n-r})$.  In particular $\frac{\partial\Phi}{\partial t^j}=X_j\circ\Phi\in\Co^{\LogL,1-}(B^r(0,\tau_0),\Omega'';U)$ for $1\le j\le r$.
\end{itemize}
By Lemmas \ref{Lem::ODE::ODEReg} and \ref{Lem::Hold::CompofMixHold} \ref{Item::Hold::CompofMixHold::Comp} and Corollary \ref{Cor::Hold::CompOp} \ref{Item::Hold::CompOp::Lip} we see that
\begin{itemize}[parsep=-0.3ex]
\item If $X_1,\dots,X_r$ are all $\Co^\alpha$ where $\alpha\in\{\Lip\}\cup(1,\infty]\cup\{k+\LogL,k+\Lip:k=1,2,3,\dots\}$, then $\Phi\in\Co^{\alpha+1,\alpha}(B^r(0,\tau_0),\Omega'';U)$ and $F''\in\Co^{\alpha}(\Phi(B^r(0,\tau_0)\times\Omega'');\R^{n-r})$. In particular $\frac{\partial\Phi}{\partial t^j}=X_j\circ\Phi\in\Co^{\alpha}(B^r(0,\tau_0)\times\Omega'';U)$ for $1\le j\le r$.
\end{itemize}
Take $\Omega':=B^r(0,\tau_0)$ and $\Omega'_\eps:=B^r(0,\tau_\eps)$. The above results give \ref{Item::RealFro::ImprovedLogFro::BigLReg}, \ref{Item::RealFro::ImprovedLogFro::Zyg1Reg}, \ref{Item::RealFro::ImprovedLogFro::LittleLReg} and  \ref{Item::RealFro::ImprovedLogFro::C1Reg} respectively, finishing the whole proof.
\end{proof}

\subsection{A PDE counterpart of Theorem \ref{MainThm::LogFro}}\label{Section::RealFro::PDEPartLogFro}

We can interpret Theorem \ref{MainThm::LogFro} in terms of the first-order PDE system following from \cite[Section 5]{Rampazzo}:

\begin{keythm}\label{KeyThm::RealFro::Further::PDEThm}Let $\Mf$ be a $n$-dimensional $C^{1,1}$-manifold. Let $L_1,\dots,L_r$ be some first order differential operators on $\Mf$ with log-Lipschitz coefficients, such that
\begin{equation}\label{Eqn::RealFro::Further::PDEAssumption1}
    L_jL_k-L_kL_j=\sum_{l=1}^rc_{jk}^l\cdot L_l,
\end{equation}
for some $c_{jk}^l\in \Co^{\LogL-1}_\loc(\Mf)$, $1\le j,k,l\le r$, in the sense of distributions.

Suppose there is a $(n-r)$-dimensional $C^1$-submanifold  $\Sf\subset \Mf$ such that 
\begin{equation}\label{Eqn::RealFro::Further::PDEAssumption2}
    \Span (L_1|_q\dots,L_r|_q)\oplus T_q\Sf=T_q \Mf,\qquad\forall q\in \Sf.
\end{equation}

Let $0<\beta\le1$ and $h\in \Co^{\beta}_\loc(\Sf)$. Then for any $p\in \Sf$ and any $0<\eps<\beta$, there is a neighborhood $U\subseteq \Mf$ of $p$, such that there exists a unique solution $f\in \Co^{\beta-\eps}_\loc(U)$ to the Cauchy problem
\begin{equation}\label{Eqn::RealFro::Further::CountPDE}
    \begin{cases}L_jf=0,&1\le j\le r,\\f|_{\Sf\cap U}=h.\end{cases}
\end{equation}
\end{keythm}Here $L_jf=0$ holds in the sense of distributions.

Applying Lemma \ref{Lem::Hold::MultLoc} \ref{Item::Hold::MultLoc::WellDef} with $L_j\in \Co^{1-}_\loc(U;T\Mf)$ and $df\in \Co^{\beta-1-\eps}_\loc(U;T^*\Mf)$, we know that $L_jf\in C^{\beta-1-\eps}_\loc(U)$ is defined as a distribution.

If in addition $L_1,\dots,L_r$ have little log-Lipschitz coefficients, we can show that $f$ given \eqref{Eqn::RealFro::Further::CountPDE} is indeed $\Co^{\beta-}$. We omit the proof to reader.

The proof is based on the construction of $\Phi$ in \eqref{Eqn::RealFro::PfLogThm::DefofPhi}, along with the following proposition on inverse functions. 
\begin{lem}\label{Lem::RealFro::Further::Mu}
Let $(x,y)=(x^1,\dots,x^r,y^1,\dots,y^{n-r})$ be the standard coordinate system of $\R^r\times\R^{n-r}$. Let $U\subseteq\R^n$ be an open neighborhood of $0$. Let $X_1,\dots,X_r$ be commutative log-Lipschitz vector fields on $U$ that have the form $X_j=\Coorvec{x^j}+\sum_{k=1}^{n-r}a_j^k\Coorvec{y^k}$ for $1\le j\le r$.

Let $\Gamma=(\Gamma^1,\dots,\Gamma^{n-r}):\R^{n-r}\to\R^r$ be a $C^1$-map such that $\Gamma(0)=0$ and $\nabla \Gamma(0)=0$. There is a neighborhood $U_0''\subseteq\R^{n-r}$ such that
\begin{equation*}
    e^{\Gamma(s)\cdot X}(0,\mu(s))=(\Gamma(s),s),
\end{equation*}
defines a unique map $\mu:U_0''\to\R^{n-r}$. Moreover,
\begin{itemize}[nolistsep]
    \item $\mu$ is hemeomorphism onto its image. In particular $\mu(U_0'')\subseteq\R^{n-r}$ is open.
    \item For any $\eps>0$ there is a smaller neighborhood $U_\eps''\subseteq U_0''$ such that $\mu\in\Co^{1-\eps}(U_\eps'';\R^{n-r})$ and $\mu^\Inv\in\Co^{1-\eps}(\mu(U_\eps'');\R^{n-r})$.
\end{itemize}

\end{lem}

\begin{proof}By assumption $X_j=\Coorvec{x^j}+\sum_{k=1}^{n-r}a_j^k\Coorvec{y^k}$, we get $\mu(s)=y(e^{-\Gamma(s)\cdot X}(\Gamma(s),s))$. Thus $\mu$ is uniquely defined and is a continuous map in the domain and satisfies $\mu(0)=0$.
By Lemma \ref{Lem::ODE::CommMultFlow} (since $X_1,\dots,X_r$ are commutative) we see that for every $0<\eps<1$ there is a neighborhood $\tilde U_\eps''$ of $0$ such that $\mu\in\Co^{1-\eps}(\tilde U_\eps'';\R^{n-r})$.

It remains to show that for every $0<\eps<1$ there are a $c_\eps>0$ and a neighborhood $U''_\eps\subseteq\tilde U_\eps''$ of $0$, such that
\begin{equation}\label{Eqn::RealFro::Further::MuInj}
    |e^{-\Gamma(s_1)\cdot X}(\Gamma(s_1),s_1)-e^{-\Gamma(s_2)\cdot X}(\Gamma(s_2),s_2)|\ge c|s_1-s_2|^{\frac1{1-\eps}},\quad\forall s_1,s_2\in U_\eps''.
\end{equation}

Once \eqref{Eqn::RealFro::Further::MuInj} is done, using $e^{-\Gamma(s)\cdot X}(\Gamma(s),s)=(0,\mu(s))$ we know $\mu$ is injective in $U_\eps''$ and $\mu^\Inv:\mu(U_\eps'')\to\R^{n-r}$ is $\Co^{1-\eps}$. In particular this implies $\mu$ and $\mu^\Inv$ are both continuous. Taking $U_0'':=U_{1/2}''$ then we finish the proof.

\medskip
To prove \eqref{Eqn::RealFro::Further::MuInj}, for any open neighborhood $V_1\subseteq \R^{n-r}$ of $0$ such that  $\{e^{u\cdot X}(\Gamma(s),s):s\in V_1,|u|\le|\Gamma(s)|\}\Subset U$, we have
\begin{equation*}
    |e^{(\Gamma(s_2)-\Gamma(s_1))\cdot X}(\Gamma(s_1),s_1)-(\Gamma(s_2),s_2)|\ge |s_1-s_2|-\|\nabla \Gamma\|_{C^0(V_1;\R^n)}\|X\|_{C^0(U;\R^n)}|s_2-s_1|,\quad s_1,s_2\in V_1.
\end{equation*}
Since $\nabla\Gamma(0)=0$, by continuity we can choose $V_1$ small enough so that $\|\nabla \Gamma\|_{C^0(V_1)}\|X\|_{C^0(U)}\le\frac12$. Thus
\begin{equation}\label{Eqn::RealFro::Further::PfMuInj::Tmp1}
    |e^{(\Gamma(s_2)-\Gamma(s_1))\cdot X}(\Gamma(s_1),s_1)-(\Gamma(s_2),s_2)|\ge\tfrac12|s_1-s_2|,\quad\forall s_1,s_2\in V_1.
\end{equation}

Let $W_1\Subset U$ be an open set containing $ \{(\Gamma(s),s):s\in V_1\}$. By Lemma \ref{Lem::ODE::CommMultFlow}, for each $\eps$ we can find $\delta>0$ such that $e^{t\cdot X}\in \Co^{1-\eps}(W_1;U)$ uniformly for $t\in B^r(0,\delta)$. Take $U'':=\{s\in V_1:|\Gamma(s)|<\frac12\delta\}$, we have for every $s_1,s_2\in U''$,
\begin{equation}\label{Eqn::RealFro::Further::PfMuInj::Tmp2}
    |e^{(\Gamma(s_2)-\Gamma(s_1))\cdot X}(\Gamma(s_1),s_1)-(\Gamma(s_2),s_2)|\le \sup_{|u|<\delta}\|e^{u\cdot X}\|_{C^{0,1-\eps}(W_1;\R^n)} |e^{-\Gamma(s_1)\cdot X}(\Gamma(s_1),s_1)-e^{-\Gamma(s_2)}(\Gamma(s_2),s_2)|^{1-\eps}.
\end{equation}
Combining \eqref{Eqn::RealFro::Further::PfMuInj::Tmp1} and \eqref{Eqn::RealFro::Further::PfMuInj::Tmp2} we get \eqref{Eqn::RealFro::Further::MuInj} with $c:=(2\sup_{|u|<\delta}\|e^{u\cdot X}\|_{\Co^{1-\eps}(W_1)})^{-\frac1{1-\eps}}>0$, finishing the proof.
\end{proof}

\begin{proof}[Proof of Theorem \ref{KeyThm::RealFro::Further::PDEThm}]
By \eqref{Eqn::RealFro::Further::PDEAssumption2} and that $\dim \Sf=n-r$, we know that $L_1,\dots,L_r$ are linearly independent at every point in $\Sf$. By continuity there is a neighborhood $\tilde U\subseteq \Mf$ of $\Sf$ such that they are still  linearly independent.

We can define a tangent subbundle $\V\le T\Mf|_{\tilde U}$ which is spanned by $L_1,\dots,L_r$. By construction $\V$ is a rank $r$ log-Lipschitz subbundle. By Proposition \ref{Prop::DisInv::CharInv1} \ref{Item::DisInv::CharInv1::Gen2}$\Rightarrow$\ref{Item::DisInv::CharInv1::Pair1}, $\V$ is involutive in the sense of distributions.

By passing to an invertible linear transform we can assume that $\Coorvec{y^1}|_p,\dots,\Coorvec{y^{n-r}}|_p$ span $T_p\Sf\le T_p \Mf$. Applying Lemma \ref{Lem::ODE::GoodGen} and shrinking $\tilde U$ if necessary, we can find a coordinate system $(x^1,\dots,x^r,y^1,\dots,y^{n-r})$ and a log-Lipschitz local basis $(X_1,\dots,X_r)$ for $\V$ on $\tilde U$ such that $[X_j,X_k]=0$ for $1\le j,k\le r$ and have expression $X_j=\Coorvec{x^j}+\sum_{k=1}^{n-r}a_j^k\Coorvec{y^k}$. Thus by shrinking $\tilde U$, in the coordinate $(x^1,\dots,x^r,y^1,\dots,y^{n-r})$ we can write
\begin{equation}\label{Eqn::RealFro::Further::PfPDEThm::Tmp0}
    \Sf\cap\tilde U=\{(\Gamma(s),s)\},\quad s\in U_0'',
\end{equation}
where $U_0''\subseteq\R^{n-r}$ is a small neighborhood of $0$ and $\Gamma:U_0''\to\R^r$ is a $C^1$-map.

By taking a translation we can assume $p=0\in\R^n$, thus $\Gamma$ satisfies $\Gamma(0)=0$ and $\nabla\Gamma(0)=0$.

Take $\tilde\eps>0$ be such that $\beta(1-\tilde\eps)^2\ge\beta-\eps$. Recall in the proof of Theorem \ref{MainThm::LogFro} in Section \ref{Section::RealFro::PfLogFro}, the map $\Phi(t,s)=e^{t\cdot X}(0,s)$ is defined in a neighborhood of $(t,s)=(0,0)$, such that
\begin{itemize}[parsep=-0.3ex]
    \item  $\Phi$ is a locally homeomorphism and $\Phi^\Inv$ is of the form $\Phi^\Inv=(x,\lambda)$, where $\lambda$ is given in \eqref{Eqn::RealFro::PfLogThm::lambda}.
    \item There is a neighborhood $U_1\subseteq \Mf$ of $0$, such that $\lambda:U_1\to\R^{n-r}$ is $\Co^{1-\tilde\eps}$.
    \item $\frac{\partial\Phi}{\partial t^j}=X_j\circ\Phi$ hold for $j=1,\dots,r$.
\end{itemize}

If $f$ is a H\"older function satisfying $X_jf=0$ for $j=1,\dots,r$, then $\Coorvec{t^j}(f\circ\Phi)(t,s)=(X_jf)\circ\Phi(t,s)=0$. Therefore for $(t,s)$ in a small enough neighborhood of $0$,
\begin{equation}\label{Eqn::RealFro::Further::PfPDEThm::Tmp1}
    f(\Phi(t,s))=f(\Phi(\Gamma(\mu^\Inv(s)),s))=f(\Gamma(\mu^\Inv(s)),\mu^\Inv(s)).
\end{equation}
Here $\mu$ is a map given in Lemma \ref{Lem::RealFro::Further::Mu} which is homeomorphism onto its image, and there is a neighborhood $U_2''$ of $0\in\R^{n-r}$ such that $\mu^\Inv\in \Co^{1-\tilde\eps}(U_2'';\R^{n-r})$.

Thus if $f|_{U\cap \Sf}=h$, then in a neighborhood of $0$ we must have
\begin{align*}
    f(x,y)&=f(\Phi(x,\lambda(x,y)))=f(\Gamma(\mu^\Inv(\lambda(x,y))),\mu^\Inv(\lambda(x,y)))=h(\Gamma(\mu^\Inv(\lambda(x,y))),\mu^\Inv(\lambda(x,y))).
\end{align*}
Thus $f$ is completely determined by $h$. This completes the proof of the uniqueness.

\medskip
Let  $U:=\mu^\Inv(U''_2)\cap U_1$. By composition we have $\mu^\Inv\circ\lambda\in \Co^{(1-\tilde\eps)^2}(U;\R^{n-r})$. Since $h\in \Co^{\beta}(\Sf)$, taking compositions we get $f\in \Co^{\beta(1-\tilde\eps)^2}(U)\subseteq \Co^{\beta-\eps}(U)$, finishing the proof of the existence and hence the whole proof.
\end{proof}

\subsection{The bi-layer Frobenius theorem}\label{Section::RealFro::BLFro}

In this subsection we apply Theorem \ref{MainThm::LogFro} to give a 2-step Frobenius theorem which is used in the proof of Theorems \ref{MainThm::CpxFro}, \ref{MainThm::RoughFro1} and \ref{MainThm::RoughFro2}.


In the proof of Theorems \ref{MainThm::CpxFro}, \ref{MainThm::RoughFro1} and \ref{MainThm::RoughFro2}, we need to apply the real Frobenius theorem on $\Se\cap\bar\Se$ and on $\Se+\bar\Se$ simultaneously. In order to achieve the best possible regularity, we need a slightly different proof. 



In the following, we endow $\R^r$, $\R^m$ and $\R^q$ with standard coordinate system $t=(t^1,\dots,t^r)$, $x=(x^1,\dots,x^m)$ and $s=(s^1,\dots,s^q)$ respectively. For the case $\kappa=\infty$, we use the convention $\kappa-1=\infty$.

\begin{keythm}[Bi-layer Frobenius theorem]\label{KeyThm::RealFro::BLFro}
	Let $\kappa\in(2,\infty]$,  Let $r,m,q$ be non-negative integers. Let $\Mf$ be a $(r+m+q)$-dimensional $\Co^\kappa$-manifold. Let $\V\le\E\le T\Mf$ be two log-Lipschitz involutive real tangent subbundles such that $\V\in \Co^{1+}$  has rank $r$, and $\E$ has rank $r+m$.

	Then for any $p\in \Mf$, there is a neighborhood $\Omega=\Omega'\times\Omega''\times\Omega'''\subseteq\R^r\times\R^m\times\R^q$ of $(0,0,0)$ and a topological parameterization $\Phi(t,x,s):\Omega\to \Mf$, such that by denoting $U:=\Phi(\Omega)$ and $F=(F',F'',F'''):=\Phi^\Inv:U\to\R^r\times\R^m\times\R^q$:
	\begin{enumerate}[parsep=-0.3ex,label=(\arabic*)]
		\item\label{Item::RealFro::BLFro::0} $\Phi(0)=p$ and $F(p)=(0,0,0)$. 
		\item\label{Item::RealFro::BLFro::1}$\frac{\partial\Phi}{\partial t^j},\frac{\partial\Phi}{\partial x^k}:\Omega\to T\Mf$ are all continuous maps for $1\le j\le r$ and $1\le k\le m$. In other words $F^*\Coorvec{t^j},F^*\Coorvec{x^k}$ are all defined continuous vector fields on $U$ for $1\le j\le r$ and $1\le k\le m$.
		\item\label{Item::RealFro::BLFro::Span} For each $\omega\in\Omega$, $\frac{\partial\Phi}{\partial t^1}(\omega),\dots,\frac{\partial\Phi}{\partial t^r}(\omega)\in T_{\Phi(\omega)}\Mf$ span $\V_{\Phi(\omega)}$, and $\frac{\partial\Phi}{\partial t^1}(\omega),\dots,\frac{\partial\Phi}{\partial t^r}(\omega)$, $\frac{\partial\Phi}{\partial x^1}(\omega),\dots,\frac{\partial\Phi}{\partial x^m}(\omega)\in T_{\Phi(\omega)}\Mf$ span $\E_{\Phi(\omega)}$.
		In other words $\V|_U$ is spanned by $F^*\Coorvec{t^1},\dots,F^*\Coorvec{t^r}$, and $\E|_U$ is spanned by $F^*\Coorvec{t^1},\dots,F^*\Coorvec{t^r},F^*\Coorvec{x^1},\dots,F^*\Coorvec{x^m}$.
	\end{enumerate}

	Let $\gamma\in(1,\infty]\cup\{k+\LogL,k+\Lip:k=1,2,\dots\}$ such that $\gamma\le\kappa-1$. Suppose $\V\in\Co^\gamma$, then automatically
	\begin{enumerate}[parsep=-0.3ex,label=(\arabic*)]\setcounter{enumi}{3}
	    \item\label{Item::RealFro::BLFro::FirstReg} $F'\in\Co_\loc^\kappa(U;\R^r)$, $F''\in\Co_\loc^\gamma(U;\R^m)$ and $F^*\Coorvec{t^j}\in\Co_\loc^\gamma(U;T\Mf|_U)$.
	    
	    \item\label{Item::RealFro::BLFro::RegF*dx} $F^*\Coorvec{x^1},\dots,F^*\Coorvec{x^m}\in\Co_\loc^{\min(\gamma-1,\LogL)}(U;T\Mf|_U)$.
	\end{enumerate}	
	
	Moreover for every $0<\eps<1$ there is a neighborhood $\Omega_\eps''\subseteq\Omega''(\subseteq\R^m)$, such that by writing $\Omega_\eps:=\Omega'\times\Omega_\eps''\times\Omega'''$ and $U_\eps:=\Phi(\Omega_\eps)$,
	\begin{enumerate}[parsep=-0.3ex,label=(\arabic*)]\setcounter{enumi}{5}
	    \item\label{Item::RealFro::BLFro::RegLLPhiF} $\Phi\in\Co_\loc^{\gamma+1,\min(\gamma,1+\LogL),1-\eps}(\Omega',\Omega_\eps'',\Omega''';\Mf)$ and $F'''\in\Co_\loc^{1-\eps}(U_\eps;\R^q)$.
	    \item\label{Item::RealFro::BLFro::RegLLDPhi}
	    $\frac{\partial\Phi}{\partial t^j}\in\Co_\loc^{\gamma,\min(\gamma,1+\LogL),1-\eps}(\Omega',\Omega_\eps'',\Omega''';T\Mf)$ and $\frac{\partial\Phi}{\partial x^k}\in\Co_\loc^{\gamma,\min(\gamma-1,\LogL),(\gamma-1)\circ(1-\eps)}(\Omega',\Omega_\eps'',\Omega''';T\Mf)$ for $1\le j\le r$ and $1\le k\le m$.
	\end{enumerate}
	
    Let $\beta\in[1,\infty]\cup\{\logl,k+\LogL,k+\Lip:k=0,1,2,\dots\}$ and we assume in addition that $\E\in\Co^\beta$. Then we have the following regularity estimates:
\begin{enumerate}[parsep=-0.3ex,label=(\arabic*')]\setcounter{enumi}{4}
    \item \label{Item::RealFro::BLFro::RegImprovedF*dx} $F^*\Coorvec{x^1},\dots,F^*\Coorvec{x^m}\in\Co_\loc^{\min(\gamma-1,\beta)}(U;T\Mf|_U)$.
\end{enumerate}	

	Moreover, when $\beta=1$,  automatically
	\begin{enumerate}[parsep=-0.3ex,label=(\arabic*')]\setcounter{enumi}{5}
	    \item\label{Item::RealFro::BLFro::RegZygPhi1} $\Phi\in\Co_\loc^{\min(\gamma,2)}L^\infty(\Omega'',\Omega'\times\Omega''';\Mf)$, $\frac{\partial\Phi}{\partial t^j}\in\Co_\loc^{1}L^\infty(\Omega'',\Omega'\times\Omega''';T\Mf)$ and $\frac{\partial\Phi}{\partial x^k}\in\Co_\loc^{\min(\gamma-1,1)}L^\infty(\Omega'',\Omega'\times\Omega''';T\Mf)$ for $1\le j\le r$ and $1\le k\le m$.
	    
	    \item\label{Item::RealFro::BLFro::RegZygPhi2} In particular let $\Omega_\eps''$ and $U_\eps=\Phi(\Omega_\eps)$ be as above, then $\Phi\in\Co_\loc^{\gamma+1,\min(\gamma,2),1-\eps}(\Omega',\Omega_\eps'',\Omega''';\Mf)$,
	    $\frac{\partial\Phi}{\partial t^1},\dots,\frac{\partial\Phi}{\partial t^r}\in\Co_\loc^{\gamma,\min(\gamma,2),1-\eps}(\Omega',\Omega_\eps'',\Omega''';T\Mf)$ and $\frac{\partial\Phi}{\partial x^1},\dots,\frac{\partial\Phi}{\partial x^m}\in\Co_\loc^{\gamma,\min(\gamma-1,1),(\gamma-1)\circ(1-\eps)}(\Omega',\Omega_\eps'',\Omega''';T\Mf)$.
	\end{enumerate}	
	
	When $\beta=\logl$, automatically
	\begin{enumerate}[parsep=-0.3ex,label=(\arabic*'')]\setcounter{enumi}{5}
	    \item\label{Item::RealFro::BLFro::RegLittlePhiF} $\Phi\in\Co_\loc^{\gamma+1,\min(\gamma,1+\LogL),1-}(\Omega',\Omega'',\Omega''';\Mf)$ and $F'''\in\Co_\loc^{1-}(U;\R^q)$.
	    
	    \item\label{Item::RealFro::BLFro::RegLittleDPhi} $\frac{\partial\Phi}{\partial t^j}\in\Co_\loc^{\gamma,\min(\gamma,1+\LogL),1-}(\Omega',\Omega'',\Omega''';T\Mf)$ and $\frac{\partial\Phi}{\partial x^k}\in\Co_\loc^{\gamma,\min(\gamma-1,\LogL),(\gamma-1)\circ(1-)}(\Omega',\Omega'',\Omega''';T\Mf)$ for $1\le j\le r$ and $1\le k\le m$.
	\end{enumerate}
	
	When $\beta\ge\Lip$, automatically
	\begin{enumerate}[parsep=-0.3ex,label=(\arabic*'\!'\!')]\setcounter{enumi}{5}
	    \item\label{Item::RealFro::BLFro::RegLipPhiF} $\Phi\in\Co_\loc^{\gamma+1,\min(\gamma,\beta+1),\min(\gamma,\beta)}(\Omega',\Omega'',\Omega''';\Mf)$ and $F'''\in\Co_\loc^\beta(U;\R^q)$.
	    
	    \item\label{Item::RealFro::BLFro::RegLipDPhi}
	    $\frac{\partial\Phi}{\partial t^j}\in\Co_\loc^{\gamma,\min(\gamma,\beta+1),\min(\gamma,\beta)}(\Omega',\Omega'',\Omega''';T\Mf)$ and $\frac{\partial\Phi}{\partial x^k}\in\Co_\loc^{\gamma,\min(\gamma-1,\beta),(\gamma-1)\circ\beta}(\Omega',\Omega'',\Omega''';T\Mf)$ for $1\le j\le r$ and $1\le k\le m$.
	\end{enumerate}

\end{keythm}

For definition of $\Co^{\alpha,\beta,\gamma}$-type maps, see Definitions \ref{Defn::Hold::MoreMixHold} \ref{Item::Hold::MoreMixHold::3Mix} and \ref{Defn::ODE::MixHoldMaps}. Recall in Definition \ref{Defn::Hold::ExtIndex} that we have orders $1-\eps<1-<\LogL<1<\Lip<1+\eps$ for all $\eps>0$. 

Here the regularity results for $F$ and those for $\Phi$ do not directly imply each other.

\begin{remark}\label{Rmk::RealFro::CompOpRecap}
    For the regularity estimates of $\frac{\partial\Phi}{\partial x^k}$ we use the notation of index composition in Definition \ref{Defn::Hold::CompIndex}: by Corollary \ref{Cor::Hold::CompOp} we have, for every $0<\eps<1$ and $\beta>\Lip$,
\begin{itemize}[nolistsep]
    \item $\alpha\in(0,1)$: $\alpha\circ(1-\eps)=\alpha(1-\eps)$, $\alpha\circ(1-)=\alpha\circ\LogL=\alpha\circ1=\alpha-$ and $\alpha\circ\Lip=\alpha\circ\beta=\alpha$.
    \item  $\alpha\in\{\LogL,1\}$: $\alpha\circ(1-\eps)=(1-\eps)-$, $\alpha\circ(1-)=\alpha\circ\LogL=\alpha\circ1=1-$, $\alpha\circ\Lip=\LogL$ and $\alpha\circ\beta=\alpha$.
    \item  $\alpha=\Lip$: $\alpha\circ\delta=\delta$ for $\delta\in\{1-\eps,1-,\LogL,\Lip\}$, $\alpha\circ1=\LogL$ and $\alpha\circ\beta=\alpha$.
    \item $\alpha>\Lip$: $\alpha\circ\delta=\min(\delta,\alpha)$ for all $\delta\in\{1-\eps,1-,\LogL,1,\Lip,\beta\}$.
\end{itemize}
\end{remark}


\begin{remark}Theorem \ref{KeyThm::RealFro::BLFro} is necessary in the proof of Theorems \ref{MainThm::RoughFro1} \ref{Item::MainThm::RoughFro1::Reg2} and \ref{MainThm::RoughFro2} \ref{Item::MainThm::RoughFro1::Reg2}, where we need a double foliations for $\Se+\bar\Se$ and $\Se\cap\bar\Se$. Say $\Se+\bar \Se\in\Co^\beta$ and $\Se\cap\bar\Se\in\Co^\gamma$ for some  $\Lip\le\beta<\gamma<2$, only using ``single'' Frobenius theorem then we will not achieve the optimal regularity estimate: if we foliate $\Se+\bar\Se$ first, then the pullback of  $\Se\cap\bar\Se$ is only $\Co^\beta$ which is not $\Co^\gamma$; if we foliate $\Se\cap\bar\Se$ first, then the pullback of $\Se+\bar\Se$ is $\Co^{\gamma-1}$, which is worse than log-Lipschitz, where real Frobenius theorem cannot apply.
    
\end{remark}


In Theorem \ref{KeyThm::RealFro::BLFro}, if we take either $m=0$ i.e. $\V=\E$ or $r=0$, $\gamma=\kappa-1$, i.e. $\V=\{0\}$, we get the H\"older estimate immediately for the classical real Frobenius theorem (see also Theorem \ref{MainThm::LogFro} \ref{Item::MainThm::LogFro::C1Reg} and Theorem \ref{KeyThm::RealFro::ImprovedLogFro} \ref{Item::RealFro::ImprovedLogFro::C1Reg}).

    


\begin{proof}[Proof of Theorem \ref{KeyThm::RealFro::BLFro}]We can pick an initial $\Co^\kappa$-coordinate system near $p$
\begin{equation}\label{Eqn::RealFro::BLFro::H}
    (u,v,w)=(u^1,\dots,u^r,v^1,\dots,v^m,w^1,\dots,w^q):U_0\subseteq \Mf\to\R^r\times\R^m\times\R^q,
\end{equation}
such that
\begin{itemize}[parsep=-0.3ex]
    \item $u(p)=0\in\R^r$, $v(p)=0\in\R^m$, $w(p)=0\in\R^q$;
    \item $\Coorvec{u^1}|_p,\dots,\Coorvec{u^r}|_p$ form a real basis for $\V_p\le T_p\Mf$;
    \item $\Coorvec{u^1}|_p,\dots,\Coorvec{u^r}|_p,\Coorvec{v^1}|_p,\dots,\Coorvec{v^m}|_p$  form a real basis for $\E_p\le T_p\Mf$.
\end{itemize}

Thus, we can identify $\V,\E$ with their pushforwards on the open neighborhood of  $(0,0,0)\in\R^r\times\R^m\times\R^q$. In other words, we can assume $\Mf\subseteq\R^r_u\times\R^m_v\times\R^q_w$ and $p=(0,0,0)$.

Applying Lemma \ref{Lem::ODE::GoodGen}, on $\V$ and on $\E$ respectively, we can find a smaller neighborhood $U_1\subseteq U_0$ of $p=0$, a $C^1$-local basis $T=[T_1,\dots,T_r]^\top$ for $\V$ on $U_1$, and a $\Co^\LogL$-local basis $\tilde X=[X',X'']^\top=[X_1,\dots,X_{r+m}]^\top$ for $\E$ on $U_1$, that have the following form:
\begin{equation}\label{Eqn::RealFro::BLFro::GenforTX}
    T=\Coorvec u+A'\Coorvec v+A''\Coorvec w,\quad \tilde X=\begin{pmatrix}X'\\X''\end{pmatrix}=\begin{pmatrix}I_r&&B'\\&I_m&B''\end{pmatrix}\begin{pmatrix}\Coorvec u\\\Coorvec v\\\Coorvec w\end{pmatrix},\quad\text{on }U_1\subseteq\R^r_u\times\R^m_v\times\R^q_w.
\end{equation}
Here $A'\in\Co^{1+}(U_1;\R^{r\times m})$, $A''\in\Co^{1+}(U_1;\R^{r\times n})$, $B'\in\Co^\LogL(U_1;\R^{r\times n})$ and $B''\in\Co^\LogL(U_1;\R^{m\times n})$.

For $\beta,\gamma$ in the assumptions, by Lemma \ref{Lem::ODE::GoodGen} \ref{Item::ODE::GoodGen::Uniqueness}, if $\V\in\Co^\gamma$ and $\E\in\Co^\beta$, then automatically $A',A''\in\Co^\gamma$ and $B',B''\in\Co^\beta$.

By \eqref{Eqn::RealFro::BLFro::GenforTX} we see that $T_1,\dots,T_r,X_1,\dots,X_m$ are linearly independent, so they form a $\Co^\LogL$-local basis for $\E$ on $U_1$.


We can now define our map $\Phi:\Omega\subseteq\R^r_t\times\R^m_x\times\R^q_s\to \Mf$ for a small enough $\Omega$, as
\begin{equation}\label{Eqn::RealFro::BLFro::DefPhi}
    \Phi(t,x,s):=\exp_T(t,\exp_{X''}(x,(0,0,s)))=e^{tT}e^{xX''}((0,0,s)).
\end{equation}

The notation $\exp_T,\exp_{X''}$ are given in Definition \ref{Defn::ODE::MultiFlow}. By taking $\tau_0$ in Lemma \ref{Lem::ODE::ODEReg} small enough we see that $\Phi$ is defined is when $\Omega$ is small enough.

\ref{Item::RealFro::BLFro::0} is immediately since $\Phi(0)=(0,0,0)=p$.

We first show that $\Phi$ is homeomorphisms onto its image by constructing the inverse map.

Denote $\tilde u:U_1\to\R^r_u$, $\tilde v:U_1\to\R^m_v$ and $\tilde w:U_1\to\R^q_w$ by $\tilde u(u,v,w)=u$, $\tilde v(u,v,w):=v$ and $\tilde w(u,v,w):=w$ as the natural projections\footnote{We abuse of notation here. It is more convenient to use $u,v,w$ as variables in the proof. Thus we use $\tilde u,\tilde v,\tilde w$ for $u,v,w$ in \eqref{Eqn::RealFro::BLFro::H}.}. We define maps $F=(F',F'',F''')$ in a neighborhood of $0\in\R^r_u\times\R^m_v\times\R^q_w$ as
\begin{equation}\label{Eqn::RealFro::BLFro::DefFG}
\begin{gathered}
    F':=\tilde u,\quad F''(u,v,w):=\tilde v(e^{-u\cdot T}(u,v,w)),\quad
    F'''(u,v,w):=\tilde w(e^{-F''(u,v,w)\cdot X''}e^{-u\cdot T}(u,v,w)). 
\end{gathered}
\end{equation}

Similar to the arguments in \eqref{Eqn::RealFro::PfLogThm::PfInv1} and \eqref{Eqn::RealFro::PfLogThm::PfInv2}, using \eqref{Eqn::RealFro::BLFro::GenforTX} we have
\begin{align*}
    F'(\Phi(t,x,s))&=\tilde u(e^{tT}e^{xX''}(0,0,s))=\tilde u(t,\ast,\ast)=t;
    \\
    F''(\Phi(t,x,s))&=\tilde v(e^{-t T}(e^{tT}e^{xX''}(0,0,s)))=\tilde v(e^{xX''}(0,0,s))=\tilde v(0,x,\ast)=x;
    \\
    F'''(\Phi(t,x,s))&=\tilde w(e^{-F''(\Phi(t,x,s))\cdot X''}e^{xX''}(0,0,s)))=\tilde w(e^{-xX''}e^{xX''}(0,0,s))=\tilde w(0,0,s)=s;
    \\
    \Phi(F(u,v,w))&=e^{uT}e^{F''(u,v,w)X''}(0,0,F'''(u,v,w))=e^{uT}e^{F''\cdot X''}(e^{-F''\cdot X''}e^{-u T}(u,v,w))=(u,v,w).
\end{align*}

Since $F$ and $\Phi$ are both continuous in their domains, we conclude that by choosing $\Omega$ small, $\Phi$ is homeomorphism onto its image, hence a topological parameterization. In particular, in what follows we assume that the domain $\Omega=\Omega'\times\Omega''\times\Omega''$ of $\Phi$, the domain $U:=\Phi(\Omega)$ of $F$, and the domain $U_1$ of $(T,X)$ satisfy the following:
\begin{equation}\label{Eqn::RealFro::BLFro::DomEqn}
    \exp_T(\Omega',U)\subseteq U_1,\quad\exp_X(\Omega'\times\Omega'',U)\subseteq U_1.
\end{equation}
Note that $F'=\tilde u$ implies $F'(U)=F'(\Phi(\Omega))=\Omega'$.


Now by construction \eqref{Eqn::RealFro::BLFro::DefPhi} we have for $1\le j\le r$ and $1\le k\le m$, in $\Omega$,
\begin{equation}\label{Eqn::RealFro::BLFro::DPhi}
    \textstyle\frac{\partial\Phi}{\partial t^j}(t,x,s)=T_j\circ\Phi(t,x,s),\quad\frac{\partial\Phi}{\partial x^k}(t,x,s)=((e^{tT})_*X_k)\circ\Phi(t,x,s)=(X_k\cdot\nabla_{u,v,w}(e^{tT}))\circ(e^{xX''}(0,0,s)).
\end{equation}
In other words
\begin{equation}\label{Eqn::RealFro::BLFro::PullbackVF}
    \textstyle F^*\Coorvec{t^j}=T_j,\quad F^*\Coorvec{x^k}(u,v,w)=((e^{uT})_*X_k)(u,v,w)=(X_k\cdot\nabla_{u,v,w}(e^{uT}))\circ(e^{-uT}(u,v,w)).
\end{equation}

Since $T_j,X_k$ are continuous, and $\exp_T$ is $C^1$, we see that $\frac{\partial\Phi}{\partial t^j},\frac{\partial\Phi}{\partial x^k}$ are all continuous, finishing the proof of \ref{Item::RealFro::BLFro::1}.

\smallskip
By construction \eqref{Eqn::RealFro::BLFro::H} and \eqref{Eqn::RealFro::BLFro::DefPhi}, we know 
$\frac{\partial\Phi}{\partial t^j}(0)=T_j(0)=\Coorvec{u^j}\big|_p$ and $\frac{\partial\Phi}{\partial x^k}(0)=X_{r+k}(0)=\Coorvec{v^k}\big|_p$, $1\le j\le r$ and $1\le k\le m$, are all linearly independent. By continuity we see that $\frac{\partial\Phi}{\partial t^1}(\omega),\dots,\frac{\partial\Phi}{\partial t^r}(\omega)$, $\frac{\partial\Phi}{\partial x^1}(\omega),\dots,\frac{\partial\Phi}{\partial x^m}(\omega)$ are still linearly independent for $\omega\in\R^{r+m+q}$ closed 0. Since $T_1=\Phi_*\Coorvec{t^1},\dots,T_r=\Phi_*\Coorvec{t^r}$ are sections of $\V$, we see that $\frac{\partial\Phi}{\partial t^1}(\omega),\dots,\frac{\partial\Phi}{\partial t^r}(\omega)$ spans $\V_{\Phi(\omega)}$ when $\omega$ small. Thus to prove \ref{Item::RealFro::BLFro::Span} it suffices to prove that $\frac{\partial\Phi}{\partial x^j}(\omega)\in\E_{\Phi(\omega)}$ for each $\omega$ in the domain.

Applying Theorem \ref{MainThm::LogFro} (or Theorem \ref{MainThm::SingFro}) to the subbundle $\E$ at the point $\Phi(\omega)\in\Mf$, we see that there is a $C^1$ (in fact $\Co^{\beta+1}$) submanifold $\Sf\subset\Mf$ containing $p$ such that $\E|_\Sf=T\Sf$. Since $\V\le\E$, we know $T_1|_\Sf,\dots,T_r|_\Sf$ are all ($C^1$) sections of $T\Sf$. Therefore by shrinking $\Omega$ if necessary, we have $e^{tT}|_{\Omega\cap\Sf}:\Omega\cap\Sf\to\Sf$. Since $e^{tT}$ is $C^1$, by restricted on $\Sf$, $(e^{tT})_*X_j$ are also sections of $T\Sf=\E|_\Sf$. In particular $\frac{\partial\Phi}{\partial x^j}(\omega)=((e^{tT})_*X_j)\circ\Phi(\omega)\in\E_{\Phi(\omega)}$, finishing the proof of \ref{Item::RealFro::BLFro::Span}.

\smallskip
We now begin the regularity estimate and we assume $\V\in\Co^\gamma$, $\E\in\Co^\beta$. Thus $T\in\Co^\gamma$ and $X',X''\in\Co^\beta$ in their domains.

Clearly $F'=\tilde u\in\Co^\kappa$. By Corollary \ref{Cor::ODE::MultFlowReg} we have $\exp_T\in\Co^{\gamma}$, thus by \eqref{Eqn::RealFro::BLFro::DefFG}  we get $F''\in\Co^\gamma$. By \eqref{Eqn::RealFro::BLFro::PullbackVF} and \eqref{Eqn::RealFro::BLFro::GenforTX} we have $F^*\Coorvec{t^j}=T_j\in\Co^\gamma$, finishing the proof of \ref{Item::RealFro::BLFro::FirstReg}.

For the rest of the regularity estimates, by Corollary \ref{Cor::ODE::MultFlowReg} we have
\begin{itemize}[parsep=-0.3ex]
    \item $\exp_T\in\Co^{\gamma+1,\gamma}(\Omega',U;U_1)$ and $\nabla_{u,v,w}\exp_T\in\Co^{\gamma,\gamma-1}(\Omega',U;\R^{r+m+q})$.
    \item $\exp_{X''}\in\Co^{\beta+1}L^\infty(\Omega'',U;U_1)$. Here we use $1+\beta=1+\LogL$ when $\beta=\logl$.
    \item When $\beta\in\{\LogL,1\}$, for every $0<\eps<1$ there is $\Omega''_\eps\subseteq\Omega''$ such that $\exp_{X''}\in\Co^{\beta+1,1-\eps/2}(\Omega''_\eps,U;U_1)$.
    \item When $\beta=\logl$, $\exp_{X''}\in\Co^{\beta+1,1-}(\Omega'',U;U_1)$.
    \item When $\beta\ge\Lip$, $\exp_{X''}\in\Co^{\beta+1,\beta}(\Omega'',U;U_1)$.
\end{itemize}

For convenience and unification of the discussion, we set for $0<\eps<1$,
\begin{equation}\label{Eqn::RealFro::BLFro::BetaEps}
    \textstyle\beta^{\vee\eps}:=1-\eps/2,\text{ if }\beta\in\{\LogL,1\};\quad\beta^{\vee\eps}:=1-,\text{ if }\beta=\logl;\quad\beta^{\vee\eps}:=\beta,\text{ if }\beta\ge\Lip.
\end{equation}
Here $\beta^{\vee\eps}$ stands for the regularity index of $\exp_{X''}$, on $\Omega''_\eps$ if $\beta\in\{\LogL,1\}$, and on $\Omega''$ if $\beta=\logl$ or $\beta\ge\Lip$.

Also we use the convention $\logl+1:=1+\LogL$. Therefore $\exp_{X''}\in\Co^{\beta+1,\beta^{\vee\eps}}$ on $\Omega''_\eps\times U$ if $\beta\in\{\LogL,1\}$, and on $\Omega''\times U$ if $\beta=\logl$ or $\beta\ge\Lip$.

Therefore applying Lemma \ref{Lem::Hold::CompofMixHold} \ref{Item::Hold::CompofMixHold::Comp} to $\Phi$ in \eqref{Eqn::RealFro::BLFro::DefPhi}, $F'''$ in \eqref{Eqn::RealFro::BLFro::DefFG}, $\frac{\partial\Phi}{\partial x^k} $ in \eqref{Eqn::RealFro::BLFro::DefPhi} and $F^*\Coorvec{x^k}$ in \eqref{Eqn::RealFro::BLFro::PullbackVF} we get
\begin{itemize}[parsep=-0.3ex]
    \item $\Phi\in\Co^{\gamma+1,\gamma\circ(\beta+1),\gamma\circ\beta^{\vee\eps}}_{t,x,s}=\Co^{\gamma+1,\min(\gamma,\beta+1),\min(\gamma,\beta^{\vee\eps})}_{t,x,s}$. This gives the estimate of $\Phi$ in \ref{Item::RealFro::BLFro::RegLLPhiF}, \ref{Item::RealFro::BLFro::RegZygPhi2}, \ref{Item::RealFro::BLFro::RegLittlePhiF} and \ref{Item::RealFro::BLFro::RegLipPhiF}.
    \item $F'''\in\Co^{\min((\beta+1)\circ\gamma,\beta^{\vee\eps}\circ\gamma)}_{u,v,w}=\Co^{\min(\gamma,\beta^{\vee\eps})}$ (since $F''\in\Co^\gamma$). This gives the estimate of $F$ in \ref{Item::RealFro::BLFro::RegLLPhiF} and \ref{Item::RealFro::BLFro::RegLittlePhiF}, and \ref{Item::RealFro::BLFro::RegLipPhiF} if $\beta\le\gamma$.
    \item $\frac{\partial\Phi}{\partial t^j}=T_j\circ\Phi\in \Co^{\gamma\circ(\gamma+1),\gamma\circ\min(\gamma,\beta+1),\gamma\circ\beta^{\vee\eps}}_{t,x,s}= \Co^{\gamma,\min(\gamma,\beta+1),\min(\gamma,\beta^{\vee\eps})}_{t,x,s}$. This gives the estimate of $\frac{\partial\Phi}{\partial t^j}$ in \ref{Item::RealFro::BLFro::RegLLDPhi}, \ref{Item::RealFro::BLFro::RegZygPhi2}, \ref{Item::RealFro::BLFro::RegLittleDPhi} and \ref{Item::RealFro::BLFro::RegLipDPhi}.
    \item $\frac{\partial\Phi}{\partial x^k}\in \Co^{\beta\circ(\beta+1),\beta\circ\beta^{\vee\eps}}_{x,s}\cdot\Co^{\gamma,(\gamma-1)\circ(\beta+1),(\gamma-1)\circ\beta^{\vee\eps}}_{t,x,s}\subseteq\Co^{\gamma,\min(\beta,\gamma-1)\circ(\beta+1),\min(\beta,\gamma-1)\circ\beta^{\vee\eps}}_{t,x,s}\subseteq\Co^{\gamma,\min(\beta,\gamma-1),(\gamma-1)\circ\beta^{\vee2\eps}}_{t,x,s}$. This gives the estimate of $\frac{\partial\Phi}{\partial x^k}$ in \ref{Item::RealFro::BLFro::RegLLDPhi}, \ref{Item::RealFro::BLFro::RegZygPhi2}, \ref{Item::RealFro::BLFro::RegLittleDPhi} and \ref{Item::RealFro::BLFro::RegLipDPhi}.
    \item $F^*\Coorvec{x^k}\in\Co^{\beta\circ\gamma}_{u,v,w}\cdot\Co^{(\gamma-1)\circ\gamma}_{u,v,w}\subseteq\Co^{\min(\beta,\gamma-1)}$. This gives \ref{Item::RealFro::BLFro::RegF*dx} and \ref{Item::RealFro::BLFro::RegImprovedF*dx}.
\end{itemize}

Also when $\beta=1$, by Corollary \ref{Cor::Hold::CompOp} we see that $\{\Phi(t,\cdot,s):t\in\Omega',s\in\Omega'''\}\subset\Co^{\gamma\circ2}_x=\Co^{\min(\gamma,2)}_x$, $\{\frac{\partial\Phi}{\partial t}(t,\cdot,s):t\in\Omega',s\in\Omega'''\}\subset\Co^{\gamma\circ1}_x=\Co^1_x$ and $\{\frac{\partial\Phi}{\partial x}(t,\cdot,s):t\in\Omega',s\in\Omega'''\}\subset\Co^{1\circ2}_x\cdot\Co^{(\gamma-1)\circ2}_x\subseteq\Co^{\min(\gamma-1,1)}_x$ are all bounded sets. Hence $\Phi\in\Co^{\min(\gamma,2)}_xL^\infty_{t,s}$ and $\frac{\partial\Phi}{\partial t^1},\dots,\frac{\partial\Phi}{\partial t^r}\in\Co^1_xL^\infty_{t,s}$ and $\frac{\partial\Phi}{\partial x^1},\dots,\frac{\partial\Phi}{\partial x^m}\in\Co^{\min(\gamma-1,1)}_xL^\infty_{t,s}$, finishing the proof of \ref{Item::RealFro::BLFro::RegZygPhi1}.

\medskip
The only thing remains to prove is that $F'''\in\Co^\gamma$ when $\beta>\gamma(>1)$, since the above composition argument only shows $F'''\in\Co^{\min(\beta,\gamma)}$ when $\beta\ge\Lip$.

To get $F'''\in\Co^\gamma$, we need to show the following
\begin{equation}\label{Eqn::RealFro::BLFro::F'''=G''}
    F'''(u,v,w)=\tilde w(e^{-u\cdot X'-v\cdot X''}(u,v,w)),\quad\text{ in a neighborhood of }p=(0,0,0)\in\Mf.
\end{equation}
Once \eqref{Eqn::RealFro::BLFro::F'''=G''} is known, by the property $\exp_X\in\Co^{\beta+1,\beta}$ (since $\beta>1$) from Lemma \ref{Lem::ODE::ODEReg}, we get $F'''\in\Co^\gamma$ and complete the proof of \ref{Item::RealFro::BLFro::RegLipPhiF}.

To prove \eqref{Eqn::RealFro::BLFro::F'''=G''}, we fix a point $(u,v,w)\in U$ in the domain of $F$ and let $\Sf:=\{(u',v',w')\in U:F'''(u',v',w')=F'''(u,v,w)\}$. By shrinking $U$ if necessary, since we only care about a neighborhood of $p$, we can assume that $\Sf$ is a connected set. 

Since $\min(\gamma,\beta)>1$, by \eqref{Eqn::RealFro::BLFro::DefFG} we know $F$ is a $C^1$-coordinate chart, thus $dF'''$ is a collection of $q$ differentials which are linearly independent at every point, therefore the level sets of $F'''$ are all $q$-dimensional $C^1$-submanifolds, in particular $\Sf$ is $C^1$. By the result \ref{Item::RealFro::BLFro::Span} we know that $\E=\Span(dF''')^\bot$, therefore $\E|_{\Sf}=T\Sf$.

On the other hand, $T_1,\dots,T_r,X_1,\dots,X_{r+m}$ are all sections of $\E$, thus by restricting to $\Sf$ they are all tangent vector fields on $\Sf$. By uniqueness of $C^1$ ODE flows we know both $e^{-F''(u,v,w)\cdot X''}e^{-u\cdot T}(u,v,w)$ and $e^{-u\cdot X'-v\cdot X''}(u,v,w)$ are in $\Sf$.

On the other hand by \eqref{Eqn::RealFro::BLFro::GenforTX} we see that $e^{-F''(u,v,w)\cdot X''}e^{-u\cdot T}(u,v,w)$ and $e^{-u\cdot X'-v\cdot X''}(u,v,w)$ lay in $\{(0,0)\}\times\R^q_w$. Thus to prove \eqref{Eqn::RealFro::BLFro::F'''=G''} it remains to show $\Sf\cap(\{(0,0)\}\times\R^q)$ is a singleton.

Suppose otherwise there are two distinct points $q_1,q_2\in \Sf\cap(\{(0,0)\}\times\R^q)$, by connectedness of $\Sf$ the distance $\dist_\Sf(q_1,q_2)$ is finite. Therefore the set $S:=\{q\in\Sf:\dist(q,q_1)^2+\dist(q,q_2)^2<\dist(q_1,q_2)^2\}$ is an open subset of $\Sf$ whose boundary contains $q_1,q_2$. By mean value theorem, since $\dim S=r+m$, there is a point $q_3\in S$ such that the displacement $q_1-q_2\in T_{q_3}S$. But $\Sf$ is spanned by $X_1|_\Sf,\dots,X_{r+m}|_\Sf$ whose $(u,v)$-coordinate components are all non-vanishing, contradiction.

Now $\Sf\cap(\{(0,0)\}\times\R^q)$ is singleton and we conclude that $e^{-F''(u,v,w)\cdot X''}e^{-u\cdot T}(u,v,w)=e^{-u\cdot X'-v\cdot X''}(u,v,w)$, which gives \eqref{Eqn::RealFro::BLFro::F'''=G''}. By composition we get $F'''\in\Co^\beta$, finishing the proof of \ref{Item::RealFro::BLFro::RegLipPhiF} and hence the whole proof.
\end{proof}


\begin{remark}
    It would be natural to generalize the result to the multi-layers case, where we have a flag of involutive tangent subbundles $0\le\V_1\le\dots\le \V_k\le T\Mf$. And the parameterization $\Phi(t^{(1)},\dots,t^{(k)},s)=e^{t^{(1)}T_{(1)}}\dots e^{t^{(k)}T_{(k)}}(\Gamma(s))$ gives the parameterization with desired regularities, where $T_{(j)}$ are collections of vector fields such that $T_{(1)},\dots,T_{(k)}$ are linearly independent and $\V_j=\Span_\R(T_{(1)},\dots,T_{(j)})$ for each $j$. We will not use this generalization in the paper.
\end{remark}

\subsection{A singular Frobenius theorem on log-Lipschitz vector fields: proof of Theorem \ref{MainThm::SingFro}}\label{Section::RealFro::PfSingFro}

First we give an analog of Lemma \ref{Lem::ODE::GoodGen} in the linear dependent setting.

\begin{lem}\label{Lem::SingFro::SingGoodGen}
Let $\Mf$ be a $C^{1,1}$-manifold with a fixed point $p\in\Mf$. Let $X_1,\dots,X_m$ be log-Lipschitz (real) vector fields on  $\Mf$. Let $r$ be the rank of $\Span(X_1(p),\dots,X_m(p))\in T_p\Mf$. Then there is a neighborhood $U\subseteq\Mf$ of $p$, a coordinate system $(x,y)=(x^1,\dots,x^r,y^1,\dots,y^{n-r}):U\to\R^r\times\R^{n-m}$, and vector fields $Y_1,\dots,Y_m$ such that
\begin{enumerate}[parsep=-0.3ex,label=(\roman*)]
    \item $(X_1,\dots,X_m)$ and $(Y_1,\dots,Y_m)$ are $\Co^\LogL$ linear combinations each other. In particular for every $q\in U$, $\Span(X_1(q),\dots,X_m(q))=\Span(Y_1(q),\dots,Y_m(q))$.
    \item $Y_1,\dots,Y_m$ have the form
    \begin{equation}\label{Eqn::SingFro::SingGoodGen}
        Y_j=\Coorvec{x^j}+\sum_{l=1}^{n-r}b_j^l\Coorvec{y^l},\quad Y_k=\sum_{l=1}^{n-r}b_k^l\Coorvec{y^l},\quad 1\le j\le r<k\le m,
    \end{equation}
    where $b_j^k\in\Co^\LogL_\loc(U)$ satisfy $b_j^k(0)=0$ for all $1\le j\le m$, $1\le k\le n-r$.
\end{enumerate}
\end{lem}
\begin{proof}
The argument is almost the same to the proof of Lemma \ref{Lem::ODE::GoodGen}. We choose a coordinate system $(x,y)$ near $p$ such that $\Span(X_1(p),\dots,X_m(p))=\Span(\Coorvec{x^1}|_p,\dots,\Coorvec{x^r}|_p)$. By reordering the index if necessary we can assume that $X_1,\dots,X_r$ are linearly independent at $p$. By continuity $X_1,\dots,X_r$ are linearly independent in a neighborhood $U\subseteq\Mf$ of $p$. By shrinking $U$ if necessary we can assume that $U$ is contained in the domain of $(x,y)$.

Write $X'=[X_1,\dots,X_r]^\top$ and $X''=[X_{r+1},\dots,X_m]^\top$ as two row collections of vector fields. By assumption of $(x,y)$ we can write $X'=A_1\Coorvec x+A_2\Coorvec y$ and $X''=A_3\Coorvec x+A_4\Coorvec y$ where $A_1\in\Co^\LogL_\loc(U;\R^{r\times r})$ is invertible at every point in the domain, $A_2\in\Co^\LogL_\loc(U;\R^{r\times (n-r)})$, $A_3\in\Co^\LogL_\loc(U;\R^{(m-r)\times r})$ and $A_4\in\Co^\LogL_\loc(U;\R^{(m-r)\times (n-r)})$ satisfy $A_2(p)=0$ and $A_4(p)=0$.

By cofactor representation of inverse matrices or Lemma \ref{Lem::Hold::CramerMixed}, we see that $A_1^{-1}\in\Co^\LogL_\loc(U;\R^{r\times r})$. We define $Y'=[Y_1,\dots,Y_r]^\top$ and $Y''=[Y_{r+1},\dots,Y_m]^\top$ as
\begin{equation*}
    Y':=A_1^{-1}X',\quad Y'':=X''-A_3Y'=X''-A_3A_1^{-1}X'.
\end{equation*}
Clearly $X''=A_1Y'$ and $X''=Y''+A_3Y'$, we see that $(X_1,\dots,X_m)$ and $(Y_1,\dots,Y_m)$ are $\Co^\LogL_\loc$ linear combinations each other.

By construction we have $Y'=\Coorvec x+B_1\Coorvec y$ and $Y''=B_2\Coorvec y$ where $B_1:=A_1^{-1}A_2\in\Co^\LogL_\loc(U;\R^{r\times(n-r)})$ and $B_2:=A_4-A_3A_1^{-1}A_2\in\Co^\LogL_\loc(U;\R^{(m-r)\times(n-r)})$. Since $A_2(p)=0$ and $A_4(p)=0$, we get $B_1(p)=0$ and $B_2(p)=0$ as well, finishing the proof.
\end{proof}

We can now start the proof by using new generators $Y_1,\dots,Y_m$ from Lemma \ref{Lem::SingFro::SingGoodGen}.
\begin{proof}[Proof of Theorem \ref{MainThm::SingFro}]
By Remark \ref{Rmk::DisInv::RmkMfldObj} and by passing to a local coordinate chart of the manifold $\Mf$, we can assume that $\Mf\subseteq\R^n$ is an open subset with the base point $p=0\in\R^n$.

By Lemma \ref{Lem::SingFro::SingGoodGen} we can find a coordinate system $(x,y):U\subseteq\Mf\to\R^r\times\R^{n-r}$ of $p$ and log-Lipschitz vector fields $Y_1,\dots,Y_m$ on $U$ which are $\Co^\LogL$ invertible linearly combinations of $X_1,\dots,X_m$ and have the form \eqref{Eqn::SingFro::SingGoodGen}.

We can write $Y_j=\sum_{k=1}^mf_j^kX_k$ and $X_j=\sum_{k=1}^mg_j^kY_k$ for $1\le j\le m$ where $f_j^k,g_j^k\in\Co^\LogL(U)$. Thus by \eqref{Eqn::MainThm::SingFro::Inv}, for $1\le j,k\le m$,
\begin{align*}
    [Y_j,Y_k]=&\sum_{j',k'=1}^m[f_j^{j'}X_{j'},f_k^{k'}Y_{k'}]=\sum_{j',k'=1}^mf_j^{j'}(X_{j'}f_k^{k'})X_{k'}-f_k^{k'}(X_{k'}f_j^l)X_{j'}+f_j^{j'}f_k^{k'}[X_{j'},Y_{k'}]
    \\
    =&\sum_{j',k',l=1}^m\bigg(f_j^{j'}(X_{j'}f_k^{k'})g_{k'}^l-f_k^{k'}(X_{k'}f_j^l)g_{j'}^l+\sum_{l'=1}^mf_j^{j'}f_k^{k'}c_{j'k'}^{l'}g_{l'}^l\bigg)Y_l=:\sum_{l=1}\tilde c_{jk}^lY_l.
\end{align*}
Here $\tilde c_{jk}^l:=\sum_{j',k'=1}^m\big(f_j^{j'}(X_{j'}f_k^{k'})g_{k'}^l-f_k^{k'}(X_{k'}f_j^l)g_{j'}^l+\sum_{l'=1}^mf_j^{j'}f_k^{k'}c_{j'k'}^{l'}g_{l'}^l\big)$ for $1\le j,k,l\le m$. By Lemma \ref{Lem::Hold::CharLogL-1} \ref{Item::Hold::CharLogL-1::Grad} we see that $\tilde c_{jk}^l\in\Co^{\LogL-1}(U)$.

On the other hand by direct computation (cf. \eqref{Eqn::ODE::GoodGen::LieBraofX}) we see that $[Y_j,Y_k]$ lays in the span of $\Coorvec{y^1},\dots,\Coorvec{y^{n-r}}$, thus $\tilde c_{jk}^l=0$ for all $1\le l\le r$.

By \eqref{Eqn::SingFro::SingGoodGen} we see that $Y_1(0),\dots,Y_r(0)$ are linearly independent and $Y_{r+1}(0)=\dots=Y_m(0)=0$. Therefore $Y_1,\dots,Y_m$ satisfy the assumptions of Proposition \ref{Prop::ODE::StrFlowComm}, and by \eqref{Eqn::ODE::StrFlowComm::ddtFlow} we see that the map $\Phi(t):=e^{t^1Y_1}\dots e^{t^rY_r}(p)$ is $C^1$ near $t=0$ and satisfies $\frac{\partial\Phi}{\partial t^j}=Y_j\circ\Phi(t)$ for $1\le j\le r$.

Since $Y_1,\dots,Y_r$ are linearly independent, we know $\frac{\partial\Phi}{\partial t^1}(t),\dots,\frac{\partial\Phi}{\partial t^r}(t)$ are linearly independent for every $t$ in the domain. Therefore by shrinking the domain if necessary, $\Phi$ gives a $C^1$-parameterization of a $r$-dimensional submanifold $\Sf$, such that $\Sf\ni p$ and
\begin{equation*}
    \textstyle T\Sf=\Span(\Phi_*\Coorvec{t^1},\dots,\Phi_*\Coorvec{t^r})=\Span(Y_1|_\Sf,\dots,Y_r|_\Sf).
\end{equation*}

On the other hand, $Y_j(\Phi(0))=Y_j(0)=0$ for all $r+1\le j\le m$, so by Proposition \ref{Prop::ODE::StrFlowComm} \ref{Item::ODE::StrFlowComm::Xj=0} we see that $Y_j(\Phi(t))\equiv0$ for all $r+1\le j\le m$ as well. We conclude that $T\Sf=\Span(Y_1,\dots,Y_m)|_\Sf=\Span(X_1,\dots,X_m)|_\Sf$ is the desired submanifold.

Now for every $q\in B_X(p,\infty)$, let $r'=r'(q)$ be the rank of $\Span(X_1(q),\dots,X_m(q))$. Repeating the argument on the point $q$, we can locally find a $r'$-dimensional $C^1$ submanifold containing $q$, whose tangent bundle is spanned by the restrictions of $X_1,\dots,X_m$. Now the set $B_X(p,\infty)$ is path connected, so $q\mapsto r'(q)$ must be a constant function. Therefore $r'(q)\equiv r$ and we conclude that $B_X(p,\infty)$ is a $r$-dimensional immersed submanifold.
\end{proof}

\section{H\"older Estimates for Nirenberg's Complex Frobenius Theorem}\label{Section::CpxFro}
Let $\Mf$ be a (smooth) manifold. We recall the following terminologies:

    We say a $\Co^{\frac12+}$ complex tangent subbundle $\Se\subseteq\C T\Mf$ is \textit{involutive}, if for every $\Co^{\frac12+}$-vector fields $X, Y$ that are sections of $\Se$, the Lie bracket $[X,Y]$ is a (distributional) section of $\Se$. (See Definition \ref{Defn::Intro::DisInv}.)
\begin{itemize}[parsep=-0.3ex]
    \item An \textbf{elliptic structure} on $\Mf$ is an involutive complex subbundle $\Se\subseteq\C T\Mf$ such that $\C T\Mf=\Se+\bar\Se$. 

    \item A \textbf{complex structure}\footnote{In this paper we use ``complex structure'' to refer the almost complex structure satisfying involutivity condition. See also \cite[Section I.8]{Involutive}.} on $\Mf$ is an involutive complex subbundle $\Se\subseteq\C T\Mf$ such that $\C T\Mf=\Se\oplus\bar\Se$.
    \item An \textbf{essentially real structure} is a complex subbundle $\Se\subseteq\C T\Mf$ such that $\Se=\bar\Se$.
\end{itemize}

Here $\bar\Se=\{(p,u-iv):p\in \Mf,\ u,v\in T_p\Mf,\ u+iv\in\Se_p\}$ is the complex conjugate subbundle of $\Se$.

The real involutive subbundles are one-one correspondent to essentially real structures in the way that, given an involutive $\V\le T\Mf$, $\Se:=\V\otimes\C\le\C T\Mf$ is essentially real, and converse given an essential real $\Se\le\C T\Mf$, $\V:=\Se\cap T\Mf$ is an involutive real subbundle.

\subsection{History remarks and ideas of the proof}\label{Section::CpxFro::Overview}
The complex Frobenius theorem was first introduced by Louis Nirenberg. It is a generalization of the real Frobenius theorem and the Newlander-Nirenberg theorem.

\begin{prop}[L. Nirenberg \cite{Nirenberg}]\label{Prop::Intro::Nirenberg}
Let $r,m,N$ be nonnegative integers satisfying $r+2m\le N$. Let $\Mf$ be an $N$-dimensional smooth manifold and let $\Se$ be a smooth complex Frobenius structure on $\Mf$, such that
$\rank(\Se\cap\bar\Se)=r$ and $\rank\Se=r+m$.
Then for any $p\in \Mf$ there is a smooth coordinate system $F=(F',F'',F'''):U\subseteq \Mf\to\R^r_t\times\C^m_z\times\R^{N-r-2m}_s$ near $p$, such that $\Se|_{U}$ is spanned by the smooth vector fields $F^*\Coorvec{t^1},\dots,F^*\Coorvec{t^r},F^*\Coorvec{z^1},\dots,F^*\Coorvec{z^m}$.
\end{prop}

There are many previous discussions of complex Frobenius theorem in the nonsmooth setting.
Nijenhuis \& Woolf \cite{NijenhuisWoolf} studied Newlander-Nirenberg theorem with  nonsmooth parameters, which is  a special case of  nonsmooth complex Frobenius theorem. In their assumptions, all structures are $\Co^\alpha$ for some $\alpha>1$.

Later Hill \& Taylor \cite{RoughComplex} studied nonsmooth complex Frobenius structures. They assumed the complex Frobenius structure to be at least $C^{1,1}$ in order to achieve the condition called ``Hypothesis V'' in their paper.

Although they did not say explicitly, a corollary they proved is the following:

\begin{prop}[\cite{RoughComplex}]\label{Prop::CpxFro::HillTaylerCpxThm}Let $\gamma>\frac12$, $\beta>0$ and let $\Mf,\Nf$ be two manifolds. Let $\Se\le(\C T\Mf)\times\Nf$ be a complex subbundle that admits generators with coefficients in $\Co^\gamma\Co^\beta(\Mf,\Nf)$, such that $\Se|_{\Mf\times\{s\}}\le\C T\Mf$ is a complex structure for all $s\in\Nf$. Then locally there is a map $F:U\times V\subseteq\Mf\times\Nf\to\C^m_z$, such that $F(\cdot,s)$ is a $C^1$ coordinate chart and $F(\cdot,s)^*\Coorvec{z^1},\dots,F(\cdot,s)^*\Coorvec{z^m}$ spans $\Se|_{\Mf\times\{s\}}$. Moreover $F$ has regularity\footnote{For the mixed regularity notation $\Co^\gamma\Co^\beta$, see Definition \ref{Defn::Hold::BiHold}.} $F\in\Co^{\gamma+1}\Co^\beta(U,V;\C^m)$.
\end{prop}

\begin{remark}
If we are begin with a such $\Se\in\Co^\alpha$ for some $\alpha>\frac12$, our regularity result in Theorem \ref{KeyThm::EllipticPara} (where we consider $\Se\in\Co^{\alpha,\alpha}$) is stronger than Proposition \ref{Prop::CpxFro::HillTaylerCpxThm}. Indeed, by interpolation the generators of $\Se$ belong to $\Co^{\theta\alpha}\Co^{(1-\theta)\alpha}(\Mf,\Nf;\C T\Mf)$ for every $0<\theta<1$. By Proposition \ref{Prop::CpxFro::HillTaylerCpxThm} for every $\theta\in(0,1)$ such that $\theta\alpha>\frac12$, we can find a map $F\in\Co^{\theta\alpha+1}\Co^{(1-\theta)\alpha}$ representing $\Se$. For overall regularity of $F$, by taking $\theta$ such that $\theta\alpha>\frac12$ and $\min(\theta\alpha+1,(1-\theta)\alpha)$ maximize:
\begin{itemize}[parsep=-0.3ex]
    \item When $\alpha>2$, by Theorem \ref{KeyThm::EllipticPara} we can find a $F\in\Co^\alpha$, while for Proposition \ref{Prop::CpxFro::HillTaylerCpxThm}, by taking $\theta=\frac{\alpha-1}{2\alpha}$, $F\in\Co^\frac{\alpha+1}2$ is the best possible.
    \item When $\frac12<\alpha\le 2$ and $\alpha\neq1$, by Theorem \ref{KeyThm::EllipticPara} we can find a $F\in\Co^{\min(\alpha,2\alpha-1)}$, while for Proposition \ref{Prop::CpxFro::HillTaylerCpxThm}, for every $0<\eps<1-\frac1{2\alpha}$ we can find a $F\in\Co^{(1-\eps)\alpha-\frac12}$ by taking $\theta=\frac1{2\alpha}+\eps$.
\end{itemize}

Notice that for both cases Theorem \ref{KeyThm::EllipticPara} gives a better estimate.
\end{remark}



Recently, Gong \cite{Gong} proved a better estimate for complex Frobenius structure. He showed that for a non-integer $\alpha>1$ and a $\Co^\alpha$-complex Frobenius structure $\Se\le \C T\Mf$, locally one can find a $\Co^{\alpha-}$-coordinate chart $F:U\subseteq\Mf\to\R^r\times\C^m\times\R^{N-r-2m}$ such that $F^*\Coorvec{t^1},\dots,F^*\Coorvec{t^r},F^*\Coorvec{z^1},\dots,F^*\Coorvec{z^m}$ span $\Se$. Though he did not have the sharp result, the statement $F\in\Co^{\alpha-}$ is endpoint sharp (up to a loss of $\eps>0$ derivatives), see Remark \ref{Rmk::EllipticPara::KeyImpliesGong}.
The proof technique in his paper is quite different from ours. He used homotopy formulae with Nash-Moser iterations.

Our result is stronger because we have not only the optimal regularity of the coordinate chart that $F\in\Co^\alpha$, but also the optimal regularity for other ingredients like the components of the coordinates and the coordinate vector fields. See Section \ref{Section::ExampleOverview}. 

\medskip
As the origin of the complex Frobenius theorem, the sharp estimate to Newlander-Nirenberg is due to Malgrange \cite{Malgrange}. He proved that for $\alpha>1$  one can always obtain a $\Co^{\alpha+1}$ local complex coordinate system for a $\Co^\alpha$ integrable almost complex structure. \cite{RoughNN} showed that his argument also works for $\frac12<\alpha\le1$. We remark that we cannot replace $(\Co^\alpha,\Co^{\alpha+1})$ by $(C^k,C^{k+1})$, see Theorem \ref{MainThm::CounterNN} and Section \ref{Section::CountNN}.

Later Street used Malgrange's method and obtained a sharp estimate for elliptic structure.
\begin{prop}[\cite{SharpElliptic}, Sharp estimate for elliptic structures]\label{Prop::Intro::SharpE}
	Let $\alpha>1$, let $\Mf$ be a $\Co^{\alpha+1}$-manifold, and let $r,m$ be integers satisfying $\dim \Mf=r+2m$. Endow $\R^r\times\C^m$ with standard coordinate system $(t,z)=(t^1,\dots,t^r,z^1,\dots,z^m)$.
	
	Let $\Se\le\C T\Mf$ be a $\Co^\alpha$-elliptic structure of rank $r+m$. Then for any $p\in \Mf$ there exists a $\Co^{\alpha+1}$-coordinate system $F:U\subseteq \Mf\to\R^r\times\C^m$ near $p$, such that $\Se|_U$ is spanned by $F^*\Coorvec{t^1},\dots,F^*\Coorvec{t^r},F^*\Coorvec{z^1},\dots,F^*\Coorvec{z^m}$.
\end{prop}

Note that in Proposition \ref{Prop::Intro::SharpE} we have $F^*\Coorvec{t^1},\dots,F^*\Coorvec{t^r},F^*\Coorvec{z^1},\dots,F^*\Coorvec{z^m}\in\Co^\alpha$ because $F\in\Co^{\alpha+1}$. However in a more general setting (Theorem \ref{MainThm::CpxFro}), due to the influence of the parameters, we only have $F^*\Coorvec{t^1},\dots,F^*\Coorvec{t^r}$, $F^*\Coorvec{z^1},\dots,F^*\Coorvec{z^m}\in\Co^{\alpha-}$. Both results for $F$ and $F^*\Coorvec{z^1},\dots,F^*\Coorvec{z^m}$ are sharp.

\bigskip
To prove Theorems \ref{MainThm::CpxFro}, \ref{MainThm::RoughFro1} and \ref{MainThm::RoughFro2}, we separate the discussion from the case $\Se\cap\bar\Se\in\Co^\alpha$ and $\Se\cap\bar\Se\in\Co^\gamma$ for some $\gamma>\alpha$. Notice that both $\Se+\bar\Se$ and $\Se\cap\bar\Se$ are an essentially real structure.

For the case $\Se\cap\bar\Se\in\Co^\alpha$, we only foliate $\Se+\bar\Se$. By applying the real Frobenius theorem, we get a foliation of $\Co^{\alpha+1}$-manifolds, such that by restricted to each leaf $\Se$ becomes an elliptic structure. The theorems are then reduced to the estimate of a family of elliptic structures, which is done in Theorem \ref{KeyThm::EllipticPara}. Theorem \ref{KeyThm::EllipticPara} implies Proposition \ref{Prop::Intro::SharpE} which can be viewed as the ``parameter free'' case. See also Remarks \ref{Rmk::EllipticPara::KeyFirst}, \ref{Rmk::EllipticPara::KeySpecial} and \ref{Rmk::EllipticPara::KeyImpliesGong}.

The idea is based on Malgrange's factorization method for the sharp Newlander-Nirenberg theorem, whose sketch can be found in, for example \cite[Page 47]{Involutive} and \cite[Section 3]{RoughNN}. Malgrange's idea is to find an intermediate coordinate change $H$ such that the generators for $H_*\Se$ satisfy an elliptic PDE system. Using the PDE we can show that $H_*\Se$ is a $\Co^\alpha$-family of real-analytic elliptic structures. Finally we construct another coordinate change $G$ using holomorphic Frobenius theorem, now with parameter, so that $G_*H_*\Se$ is the span of $\Coorvec{t^1},\dots,\Coorvec{t^r},\Coorvec{z^1},\dots,\Coorvec{z^m}$. Thus $F=G\circ H$ is the desired coordinate chart. Also see Section \ref{Section::EllipticPara::Overview} for a detailed sketch to the proof of Theorem \ref{KeyThm::EllipticPara}.

Note that in Theorem \ref{KeyThm::EllipticPara} we work on the elliptic structures rather than the complex structures. In this way we do not have to pick generators for $\Se$ that satisfy the ``Hypothesis V'' in \cite{RoughComplex}.

\medskip
For the case $\Se\cap\bar\Se\in\Co^\gamma$ where $\gamma\ge\alpha+1$, we follow Nirenberg's original proof by applying real Frobenius theorem twice on the essentially real structures $\Se\cap\bar\Se$ and $\Se+\bar\Se$, see \cite[the proof of Theorem I.12.1]{Involutive}. This ends up with a family of complex structures. We complete the proof by applying Theorem \ref{KeyThm::EllipticPara} with the special case $r=0$.

To get the optimal H\"older-Zygmund regularity during this reduction, we combine the two real Frobenius theorems into one step, which is done in Theorem \ref{KeyThm::RealFro::BLFro}.

\medskip
When $\alpha>1$, unlike the sharp Newlander-Nirenberg theorem or Street's result where $F^*\Coorvec{z^1},\dots,F^*\Coorvec{z^m}$ are $\Co^\alpha$-vector fields, in our result we show that $F^*\Coorvec{z^j}\in\Co^{\alpha-}$ is optimal. 

One reason for the regularity loss comes from solving elliptic pdes. Say we have a $\Co^\alpha$-subbundle for some $\alpha>0$. In the first step of Malgrange's proof of the sharp Newlander-Nirenberg theorem we need to solve an elliptic pde system for $f=f(x)$ that has the following form:

\begin{equation}\label{Eqn::Intro::ElliPDE}
    \Delta f=\phi(f,\nabla f,\nabla^2f,a,\nabla a).
\end{equation}
Here $a=a(x)\in\Co^\alpha$ is the vector-valued coefficient that comes from the subbundle, and $\phi$ is a rational function which is defined and small when $f$ and $a$ are both small.

Since Laplace's equation gains two derivatives in H\"older-Zygmund spaces, we can find a $f\in\Co^{\alpha+1}$ solving \eqref{Eqn::Intro::ElliPDE}. This is why we can find a $\Co^{\alpha+1}$-coordinate system for a $\Co^\alpha$-integrable almost complex structure. In particular the coordinate vector fields are $\Co^\alpha$ because $\nabla f\in\Co^\alpha$.

When we do estimate on the complex Frobenius structure, we need to solve a pde system that has similar expression to \eqref{Eqn::Intro::ElliPDE} but with a parameter: we need $f=f(x,s)$ for the system $\Delta_xf=\phi(f,\nabla_xf,a,\nabla_xa)$ with $\Co^\alpha$-coefficient $a=a(x,s)$.  See \eqref{Eqn::EllipticPara::ExistenceH} for the actual formula.

We need to introduce some suitable function spaces that capture the property that the pde system ``gains 2 derivatives on $x$-variable''. Note that by comparing to $a(x,s)$, the function $\phi(f,\nabla_xf,a,\nabla_xa)$ lose one derivative in $x$, and by taking inverse Laplacian, $f$ should gain 1 derivative in $x$. Thus we can ask what is the good function space $\Xs=\Xs_{x,s}$ to work on, such that if $a\in\Xs_{x,s}$ then $\nabla_xf\in\Xs_{x,s}$.

Unfortunately the space $\Co^\alpha_{x,s}$ does not satisfy our need. When $a\in\Co^\alpha_{x,s}$ we may not find a solution $f$ such that $\nabla f\in\Co^\alpha_{x,s}$. This is roughly because the $x$-variable Riesz transform (say $x\in\R^n$ and $s\in V\subseteq\R^m$)
$$R_jf(x,s):=\frac{\Gamma(\frac{n+1}2)}{\pi^\frac{n+1}2}\lim\limits_{\eps\to 0}\int_{\R^n\backslash B^n(0,\eps)}\frac{y_j-x_j}{|y-x|^{n+1}}f(y,s)dy, \quad j=1,\dots,n,$$
does not map $\Co^\alpha(\R^n\times V)$ into itself. That is why $F^*\Coorvec{z^j}\in\Co^{\alpha-}$ is the best one we can get.

To see that why $R_j:\Co^\alpha_{x,s}\not\to \Co^\alpha_{x,s}$ holds, we consider the bi-parameter decomposition  $\Co^\alpha_{x,s}=\Co^\alpha_xL^\infty_s\cap L^\infty_x\Co^\alpha_s$ (see Lemma \ref{Lem::Hold::CharMixHold}): we have $\|f\|_{\Co^\alpha_{x,s}}\approx\sup_s\|f(\cdot,s)\|_{\Co^\alpha_x}+\sup_x\|f(x,\cdot)\|_{\Co^\alpha_s}$. We see that $R_j:\Co^\alpha_xL^\infty_s\to \Co^\alpha_xL^\infty_s$ is still bounded because Riesz transforms are bounded in $\Co^\alpha$-spaces. However $R_j:L^\infty_x\Co^\alpha_s\not\to L^\infty_x\Co^\alpha_s$ is not bounded because Riesz transforms are unbounded on $L^\infty$-spaces.

In application we use the H\"older-Zygmund space $\Xs_{x,s}=\Co^\alpha_xL^\infty\cap\Co^0_x\Co^\alpha_s$, thus $f\in \Co^{\alpha+1}_xL^\infty\cap\Co^1_x\Co^\alpha_s\subset\Co^\alpha_{x,s}$ and $\nabla f\in \Co^{\alpha-}_{x,s}$, both of which are optimal.

To construct the example of complex Frobenius structure such that $F^*\Coorvec z\notin\Co^\alpha$, we follow a similar but more sophisticated construction to \cite{LidingCounterNN}. See Propositions \ref{Prop::Exmp::SharpddzRed} and \ref{Prop::Exmp::Expa} in Section \ref{Section::Sharpddz}.

When $\frac12<\alpha<1$, in Theorem \ref{KeyThm::EllipticPara} we are only able to get $F\in\Co^{2\alpha-1}$ instead of $F\in\Co^\alpha$. Because in solving the existence of \eqref{Eqn::Intro::ElliPDE}, the image of the product map $(\Co^\alpha_xL^\infty\cap\Co^0_x\Co^\alpha_s)\times (\Co^{\alpha-1}_xL^\infty\cap\Co^{-1}_x\Co^\alpha_s)$ is not contained in $\Co^{\alpha-1}_xL^\infty\cap\Co^{-1}_x\Co^\alpha_s$: there is an extra gain in $\Co^0_x\Co^{2\alpha-1}_s$, which contributes to the additional regularity loss in parameters. We do not know whether the $\Co^{2\alpha-1}$-smoothness is sharp or not.

\subsection{The general versions of Theorems \ref{MainThm::CpxFro}, \ref{MainThm::RoughFro1}, \ref{MainThm::RoughFro2}, and their proofs}\label{Section::PfCpxFro}

Now we give more general statements for Theorems \ref{MainThm::CpxFro}, \ref{MainThm::RoughFro1} and \ref{MainThm::RoughFro2} with the estimate of mixed regularity for the parameterizations. In the proof we combine the cases $\alpha<1$ and $\alpha>1$ together, but separate the discussion between $\Se\cap\bar\Se\in\Co^\alpha$ and $\Se\cap\bar\Se\in\Co^\gamma$ for some $\gamma\ge\alpha+1$.

Recall  Definitions \ref{Defn::Hold::MixHold}, \ref{Defn::Hold::MoreMixHold} and \ref{Defn::ODE::MixHoldMaps} for the $\Co^{\alpha,\beta}$-spaces, Definition \ref{Defn::ODE::CpxSubbd} for a local basis, and \eqref{Eqn::Intro::ColumnNote} for the column and matrix convention of Jacobian matrices. Recall in Definition \ref{Defn::Hold::ExtIndex} we use order $1-<\LogL<1<\Lip<1+\eps$ for generalized indices, and we do not include $\logl$ in the index set $\R_\Eb$.

We endow $\R^r\times\C^m\times\R^q$ with standard coordinate system $(t,z,s)=(t^1,\dots,t^r,z^1,\dots,z^m,s^1,\dots,s^q)$.

When $\Se\cap\bar\Se$ is as rough as $\Se$, we have the following:
\begin{keythm}\label{KeyThm::CpxFro1} 

Let $\kappa\in(2,\infty]$, $\alpha\in(\frac12,\kappa-1]$ and $\beta\in\{\LogL,\logl,\Lip\}\cup(1,\kappa-1]$ satisfy $\alpha\le\beta$ if $\beta\neq\logl$ (see Definition \ref{Defn::Hold::ExtIndex}) and $\alpha\le1$ if $\beta=\logl$. Let $r,m,q$ be nonnegative integers. Let $\Mf$ be a $(r+2m+q)$-dimensional $\Co^\kappa$-manifold. 

Let $\Se$ be a $\Co^\alpha$-complex Frobenius structure over $\Mf$ that has complex rank $r+m$, such that  $\Se+\bar\Se$ is a $\Co^\beta$-subbundle and has complex rank $r+2m$.
Then for any $p_0\in \Mf$ and $0<\eps<2\alpha-1$ {\normalfont($\eps$ is only required for the case $\beta=\LogL$)} there is an open neighborhood $\Omega=\Omega'\times\Omega''\times\Omega'''\subseteq\R^r_t\times\C^m_z\times\R^q_s$ of $(0,0,0)$ and a topological parameterization $\Phi:\Omega\to\Mf$, such that by denoting $U:=\Phi(\Omega)$ and $F=(F',F'',F'''):=\Phi^\Inv:U\to\R^r\times\C^m\times\R^q$:
\begin{enumerate}[parsep=-0.3ex,label=(\arabic*)]
    \item\label{Item::CpxFro::0}  $\Phi(0)=p_0$ and $F(p_0)=(0,0,0)$.
    \item\label{Item::CpxFro::1}  $\frac{\partial\Phi}{\partial t^j},\frac{\partial\Phi}{\partial z^k},\frac{\partial\Phi}{\partial\bar z^k}:\Omega\to \C T\Mf$ are all continuous maps for $1\le j\le r$ and $1\le k\le m$. In other words $F^*\Coorvec{t^j},F^*\Coorvec{z^k},F^*\Coorvec{\bar z^k}$ are all defined continuous vector fields on $U$ for $1\le j\le r$ and $1\le k\le m$.
    \item\label{Item::CpxFro::Span} For any $\omega\in\Omega$, the complex subspace $\Se_{\Phi(\omega)}$ is spanned by $\frac{\partial\Phi}{\partial t^1}(\omega),\cdots,\frac{\partial\Phi}{\partial t^r}(\omega),\frac{\partial\Phi}{\partial z^1}(\omega),\dots,\frac{\partial\Phi}{\partial z^m}(\omega)\in\C T_{\Phi(u)}\Mf$. In other words $\Se|_U$ is spanned by $F^*\Coorvec{t^1},\dots,F^*\Coorvec{t^r},F^*\Coorvec{z^1},\dots,F^*\Coorvec{z^m}$.
\end{enumerate}
In particular,
\begin{enumerate}[parsep=-0.3ex,label=(\arabic*)]\setcounter{enumi}{3}
    \item\label{Item::CpxFro::SCapSpan} $(\Se\cap\bar\Se)|_U$ has a local basis $(F^*\Coorvec{t^1},\dots,F^*\Coorvec{t^r})$.
    \item\label{Item::CpxFro::S+Span} $(\Se+\bar\Se)|_U$ has a local basis $(F^*\Coorvec{t^1},\dots,F^*\Coorvec{t^r},F^*\Coorvec{z^1},\dots,F^*\Coorvec{z^m},F^*\Coorvec{\bar z^1},\dots,F^*\Coorvec{\bar z^m})$.
\end{enumerate}

Moreover for the regularity estimates, when $\beta=\LogL$ (in particular $\frac12<\alpha\le1$):
\begin{enumerate}[parsep=-0.3ex,label=(\arabic*)]\setcounter{enumi}{5}
    \item\label{Item::CpxFro::PhiReg1} $\Phi\in\Co^{\alpha+1,2\alpha-1-\eps}_\loc(\Omega'\times\Omega'',\Omega''';\Mf)$, $\frac{\partial\Phi}{\partial t^j}\in\Co^{\alpha,2\alpha-1-\eps}_\loc(\Omega'\times\Omega'',\Omega''';T\Mf)$ for $1\le j\le r$, and $\frac{\partial\Phi}{\partial z^k},\frac{\partial\Phi}{\partial\bar z^k}\in\Co^{\alpha,2\alpha-1-\eps}_\loc(\Omega'\times\Omega'',\Omega''';\C T\Mf)$ for $1\le k\le m$.
    \item\label{Item::CpxFro::FReg1} $F'\in\Co^\kappa_\loc(U;\R^r)$, $F''\in\Co^{2\alpha-1-\eps}_\loc(U;\C^m)$ and $F'''\in\Co^{1-\eps}_\loc(U;\R^q)$.
    \item\label{Item::CpxFro::FDReg1}$F^*\Coorvec{t^j},F^*\Coorvec{z^k},F^*\Coorvec{\bar z^k}$ are $\Co^{2\alpha-1-\eps}$ vector fields on $U$ for $1\le j\le r$ and $1\le k\le m$.
\end{enumerate}

When $\beta=\logl$ (in particular $\frac12<\alpha\le1$), or $\beta\ge\Lip$ and $\alpha=1$:
\begin{enumerate}[parsep=-0.3ex,label=(\arabic*')]\setcounter{enumi}{5}
    \item\label{Item::CpxFro::PhiReg1.5} $\Phi\in\Co^{\alpha+1,(2\alpha-1)-}_\loc(\Omega'\times\Omega'',\Omega''';\Mf)$, $\frac{\partial\Phi}{\partial t^j}\in\Co^{\alpha,(2\alpha-1)-}_\loc(\Omega'\times\Omega'',\Omega''';T\Mf)$ for $1\le j\le r$, and $\frac{\partial\Phi}{\partial z^k},\frac{\partial\Phi}{\partial\bar z^k}\in\Co^{\alpha,(2\alpha-1)-}_\loc(\Omega'\times\Omega'',\Omega''';\C T\Mf)$ for $1\le k\le m$.
    \item\label{Item::CpxFro::FReg1.5} $F'\in\Co^\kappa_\loc(U;\R^r)$, $F''\in\Co^{(2\alpha-1)-}_\loc(U;\C^m)$ and $F'''\in\Co^{1-}_\loc(U;\R^q)$.
    \item\label{Item::CpxFro::FDReg1.5}$F^*\Coorvec{t^j},F^*\Coorvec{z^k},F^*\Coorvec{\bar z^k}$ are $\Co^{(2\alpha-1)-}$ vector fields on $U$ for $1\le j\le r$ and $1\le k\le m$.
\end{enumerate}

When $\beta\ge\Lip$ and $\frac12<\alpha<1$:
\begin{enumerate}[parsep=-0.3ex,label=(\arabic*'')]\setcounter{enumi}{5}
    \item\label{Item::CpxFro::PhiReg2} $\Phi\in\Co^{\alpha+1,2\alpha-1}_\loc(\Omega'\times\Omega'',\Omega''';\Mf)$, $\frac{\partial\Phi}{\partial t^j}\in\Co^{\alpha,2\alpha-1}_\loc(\Omega'\times\Omega'',\Omega''';T\Mf)$, and $\frac{\partial\Phi}{\partial z^k},\frac{\partial\Phi}{\partial\bar z^k}\in\Co^{\alpha,2\alpha-1}_\loc(\Omega'\times\Omega'',\Omega''';\C T\Mf)$ for $1\le j\le r$ and $1\le k\le m$.
    \item\label{Item::CpxFro::FReg2} $F'\in\Co^\kappa_\loc(U;\R^r)$, $F''\in\Co^{2\alpha-1}_\loc(U;\C^m)$ and $F'''\in\Co^{\beta}_\loc(U;\R^q)$.
    \item\label{Item::CpxFro::FDReg2}$F^*\Coorvec{t^j},F^*\Coorvec{z^k},F^*\Coorvec{\bar z^k}$ are $\Co^{2\alpha-1}$ vector fields on $U$ for $1\le j\le r$ and $1\le k\le m$.
\end{enumerate}

When $\beta>1$ and $1<\alpha\le\beta$ {\normalfont($\alpha$ is a real number)}:
\begin{enumerate}[parsep=-0.3ex,label=(\arabic*'\!'\!')]\setcounter{enumi}{5}
    \item\label{Item::CpxFro::PhiReg3} $\Phi\in\Co^{\alpha+1,\alpha}_\loc(\Omega'\times\Omega'',\Omega''';\Mf)$, $\frac{\partial\Phi}{\partial t^j}\in\Co^{\alpha,\alpha-}_\loc(\Omega'\times\Omega'',\Omega''';T\Mf)$, and $\frac{\partial\Phi}{\partial z^k},\frac{\partial\Phi}{\partial\bar z^k}\in\Co^{\alpha,\alpha-}_\loc(\Omega'\times\Omega'',\Omega''';\C T\Mf)$ for $1\le j\le r$ and $1\le k\le m$.
    \item\label{Item::CpxFro::FReg3} $F'\in\Co^\kappa_\loc(U;\R^r)$, $F''\in\Co^{\alpha}_\loc(U;\C^m)$ and $F'''\in\Co^{\beta}_\loc(U;\R^q)$.
    \item\label{Item::CpxFro::FDReg3}$F^*\Coorvec{t^j},F^*\Coorvec{z^k},F^*\Coorvec{\bar z^k}$ are $\Co^{\alpha-}$ vector fields on $U$ for $1\le j\le r$ and $1\le k\le m$.
\end{enumerate}

As a corollary, when $\alpha>1$,
\begin{enumerate}[parsep=-0.3ex,label=(\arabic*)]\setcounter{enumi}{8}
    \item\label{Item::CpxFro::ParaCoor} $\Phi$ is a $\Co^\alpha$ regular parameterization, and $F$ is a $\Co^\alpha$ coordinate chart.
    \item\label{Item::CpxFro::SBotSpan} $\Se^\bot|_U$ has a local basis $(F^*d\bar z^1,\dots,F^*d\bar z^m,F^*ds^1,\dots,F^*ds^q)$. 
\end{enumerate}
\end{keythm}
When $\Se\cap\bar\Se$ is at least one order more regular than $\Se$, we have the following:
\begin{keythm}\label{KeyThm::CpxFro2}
Let $\alpha,\beta,\kappa$, $r,m,q\ge0$, $\Mf$ and $\Se\le\C T\Mf$ be the same as in Theorem \ref{KeyThm::CpxFro1}.

Suppose in addition to the assumptions of Theorem \ref{KeyThm::CpxFro1}, there is a $\gamma\in(1,\kappa-1]$ such that $\gamma\ge\alpha$ and $\Se\cap\bar\Se$ is a  $\Co^\gamma$-subbundle.

Then for any $p_0\in \Mf$ and $0<\eps<2\alpha-1$ {\normalfont($\eps$ is only required for the case $\beta=\LogL$)} there are an open neighborhood $\Omega=\Omega'\times\Omega''\times\Omega'''\subseteq\R^r_t\times\C^m_z\times\R^q_s$ of $(0,0,0)$, and a map $\Phi:\Omega\to\Mf$, such that $\Phi$ satisfies all the consequences in Theorem \ref{KeyThm::CpxFro1} except  \ref{Item::CpxFro::PhiReg1}, \ref{Item::CpxFro::PhiReg1.5}, \ref{Item::CpxFro::PhiReg2}, \ref{Item::CpxFro::PhiReg3}. And we have the following:
\begin{enumerate}[parsep=-0.3ex,label=(\arabic*)]\setcounter{enumi}{10}
    \item\label{Item::CpxFro::ImprovePhi}
    When $\beta=\LogL$: $\Phi\in\Co^{\gamma+1,\alpha+1,2\alpha-1-\eps}_\loc(\Omega',\Omega'',\Omega''';\Mf)$ and $\frac{\partial\Phi}{\partial t^1},\dots,\frac{\partial\Phi}{\partial t^r}\in\Co^{\gamma,\alpha+1,2\alpha-1-\eps}_\loc(\Omega',\Omega'',\Omega''';T\Mf)$.
    
    When $\beta=\logl$, or $\beta\ge\Lip$ and $\alpha=1$: $\Phi\in\Co^{\gamma+1,\alpha+1,(2\alpha-1)-}_{t,z,s,\loc}$ and $\frac{\partial\Phi}{\partial t^1},\dots,\frac{\partial\Phi}{\partial t^r}\in\Co^{\gamma,\alpha+1,(2\alpha-1)-}_{t,z,s,\loc}$.
    
    When $\beta\ge\Lip$ and $\alpha\neq1$: $\Phi\in\Co^{\gamma+1,\alpha+1,\min(\alpha,2\alpha-1)}_{t,z,s,\loc}$ and  $\frac{\partial\Phi}{\partial t^1},\dots,\frac{\partial\Phi}{\partial t^r}\in\Co^{\gamma,\alpha+1,\min(\alpha,2\alpha-1)}_{t,z,s,\loc}$.
    
    \item\label{Item::CpxFro::ImproveDPhiDz}
    When $\beta=\LogL$: $\frac{\partial\Phi}{\partial z^1},\frac{\partial\Phi}{\partial\bar z^1},\dots,\frac{\partial\Phi}{\partial z^m},\frac{\partial\Phi}{\partial\bar z^m} \in\Co^{\gamma,\alpha,(2\alpha-1-\eps)\cdot\min(\gamma-1,1)}_\loc(\Omega',\Omega'',\Omega''';\C T\Mf)$.
    
    When $\beta=\logl$:  $\frac{\partial\Phi}{\partial z^1},\frac{\partial\Phi}{\partial\bar z^1},\dots,\frac{\partial\Phi}{\partial z^m},\frac{\partial\Phi}{\partial\bar z^m} \in\Co^{\gamma,\alpha,(2\alpha-1)\cdot\min(\gamma-1,1)-}_\loc(\Omega',\Omega'',\Omega''';\C T\Mf)$.
    
    When $\beta\ge\Lip$: $\frac{\partial\Phi}{\partial z^1},\frac{\partial\Phi}{\partial\bar z^1},\dots,\frac{\partial\Phi}{\partial z^m},\frac{\partial\Phi}{\partial\bar z^m} \in\Co^{\gamma,\alpha,(\gamma-1)\circ\min(\alpha-,2\alpha-1)}_\loc(\Omega',\Omega'',\Omega''';\C T\Mf)$.
    
    
    
    \item\label{Item::CpxFro::ImproveDDt}For $F=\Phi^\Inv$, the vector fields $F^*\Coorvec{t^1},\dots,F^*\Coorvec{t^r}$ defined on $U=\Phi(\Omega)$ are all $\Co^\gamma$.
    
    
\end{enumerate}
\end{keythm}
\begin{remark}\label{Rmk::CpxFro::RmkofCpxFro1and2}

For $(\gamma-1)\circ\min(\alpha-,2\alpha-1)$ in \ref{Item::CpxFro::ImproveDPhiDz} we recall the composition notation in Definition \ref{Defn::Hold::CompIndex} and Remark \ref{Rmk::RealFro::CompOpRecap}. 

Thus \ref{Item::CpxFro::ImprovePhi} is stronger than the estimates of $\Phi$ and $\nabla_t\Phi$ in \ref{Item::CpxFro::PhiReg1}, \ref{Item::CpxFro::PhiReg1.5}, \ref{Item::CpxFro::PhiReg2} and \ref{Item::CpxFro::PhiReg3}, and \ref{Item::CpxFro::ImproveDDt} is stronger than the estimates of $F^*\Coorvec t$ in \ref{Item::CpxFro::FDReg1}, \ref{Item::CpxFro::FDReg1.5}, \ref{Item::CpxFro::FDReg2} and \ref{Item::CpxFro::FDReg3}.

When $\gamma>2$ or $(\beta,\gamma)\in\{(\LogL,2),(\logl,2)\}$, the results in \ref{Item::CpxFro::ImproveDPhiDz}  are stronger than the estimates of $\nabla_z\Phi$ in \ref{Item::CpxFro::PhiReg1}, \ref{Item::CpxFro::PhiReg1.5}, \ref{Item::CpxFro::PhiReg2} and \ref{Item::CpxFro::PhiReg3}. But when $\frac32<\gamma\le2$ and  $(\beta,\gamma)\notin\{(\LogL,2),(\logl,2)\}$, in particular $\frac12<\alpha\le\gamma-1(\le1)$, the $s$-variable regularity in \ref{Item::CpxFro::ImproveDPhiDz} is weaker than those of $\nabla_z\Phi$ in \ref{Item::CpxFro::PhiReg1}, \ref{Item::CpxFro::PhiReg1.5} and \ref{Item::CpxFro::PhiReg2}.

\end{remark}

Some of the results are directly implied by others, we mark them as below.
\begin{remark}\label{Rmk::CpxFro::RmkofTrueThm}
\begin{enumerate}[parsep=-0.3ex,label=(\roman*)]
    \item\label{Item::RmkofCpxFro::Infty} The case $\alpha=\infty$ (which implies $\beta=\gamma=\kappa=\infty$) is the classical result \cite{Nirenberg} (see Proposition \ref{Prop::Intro::Nirenberg}).
    \item\label{Item::RmkofCpxFro::FDReg} \ref{Item::CpxFro::FDReg3} is implied by \ref{Item::CpxFro::PhiReg3} and \ref{Item::CpxFro::FReg3}.
\end{enumerate}

\smallskip
We have $F^*\Coorvec{t^j}=\frac{\partial\Phi}{\partial t^j}\circ F$ and $F^*\Coorvec{z^k}=\frac{\partial\Phi}{\partial z^k}\circ F$ for $1\le j\le r$, $1\le k\le m$. \ref{Item::CpxFro::PhiReg3} shows that $\frac{\partial\Phi}{\partial t^j},\frac{\partial\Phi}{\partial z^k}\in\Co^{\alpha-}$ and \ref{Item::CpxFro::FReg3} shows that $F\in\Co^{\min(\kappa,\alpha,\beta)}=\Co^\alpha$. Taking compositions we get $\frac{\partial\Phi}{\partial t^j}\circ F,\frac{\partial\Phi}{\partial z^k}\circ F\in\Co^{\alpha-}$.\hfill\qedsymbol

\begin{enumerate}[parsep=-0.3ex,label=(\roman*)]\setcounter{enumi}{3}
    \item\label{Item::RmkofCpxFro::OtherSpan} \ref{Item::CpxFro::SBotSpan}, \ref{Item::CpxFro::SCapSpan} and \ref{Item::CpxFro::S+Span} are all implied by \ref{Item::CpxFro::Span}.
\end{enumerate} 

For the complex conjugate of \ref{Item::CpxFro::Span}, $\bar\Se|_U$ has local basis $F^*\Coorvec t,F^*\Coorvec{\bar z}$. Taking intersection  we get that $\Se\cap\bar\Se|_U$ is spanned by $F^*\Coorvec t$. Taking union we get that $\Se+\bar\Se|_U$ is spanned by $F^*\Coorvec t,F^*\Coorvec z,F^*\Coorvec{\bar z}$.

When $\alpha>1$, since $(F^*dt,F^*dz,F^*d\bar z,F^*ds)$ is the dual basis of $(F^*\Coorvec t,F^*\Coorvec z,F^*\Coorvec{\bar z},F^*\Coorvec s)$,  and by \ref{Item::CpxFro::Span} $F^*\Coorvec t,F^*\Coorvec z$ span $\Se|_U$. Taking the dual bundle, we get that $F^*d\bar z,F^*ds$ span $\Se^\bot|_U$.\hfill\qedsymbol
\end{remark}

The proofs of Theorems \ref{KeyThm::CpxFro1} and \ref{KeyThm::CpxFro2} are combinations of the real Frobenius theorem, the Theorems \ref{KeyThm::RealFro::BLFro} or  \ref{KeyThm::RealFro::ImprovedLogFro}, and the estimate of the parameterized elliptic structures, the Theorem \ref{KeyThm::EllipticPara}. 

\begin{proof}[Proof of Theorem \ref{KeyThm::CpxFro1}]By Remark \ref{Rmk::CpxFro::RmkofTrueThm} \ref{Item::RmkofCpxFro::Infty} only the case $\alpha<\infty$ needs a proof.

Throughout the proof we choose $\eps>0$ such that $0<\eps<\alpha-\frac12$.

By Lemma \ref{Lem::ODE::S+IsBundle} since $\Se\in\Co^\alpha$, we know $\Se+\bar\Se\in\Co^\alpha$ as well. This is useful when $\alpha=1$ and $\beta\in\{\LogL,\logl\}$.

We fix the point $p_0\in\Mf$.
First we apply real Frobenius theorem, the Theorem \ref{KeyThm::RealFro::ImprovedLogFro} to $\Se+\bar\Se$. We endow $\R^{r+2m}\times\R^q$ with standard coordinate system $(\nu,\lambda)=(\nu^1,\dots,\nu^{r+2m},\lambda^1,\dots,\lambda^q)$.

Since $\Se+\bar\Se$ is at least log-Lipschitz and has rank $r+2m$, by Theorem \ref{KeyThm::RealFro::ImprovedLogFro}, we can find a neighborhood $\Omega_1=\Omega'_1\times\Omega'''_1\subseteq\R^{r+2m}_\nu\times\R^q_\lambda$ of $(0,0)$ and a topological parameterization $\Phi_1:\Omega_1\to\Mf$, which does not depend on $\beta$, such that by taking $U_1:=\Phi_1(\Omega_1)$ and $F_1=(F_1',F_1'''):U_1\to\R^{r+2m}\times\R^q$, $F_1:=\Phi_1^\Inv$, we have
\begin{enumerate}[parsep=-0.3ex,label=(S1.\arabic*)]
    \item\label{Item::CpxFroPf1::Phi10} $\Phi_1(0,0)=p_0$ and $\frac{\partial\Phi_1}{\partial\nu^1},\dots,\frac{\partial\Phi_1}{\partial\nu^{r+2m}}:\Omega\to T\Mf$ are continuous. 
    \item\label{Item::CpxFroPf1::RFroSpan} $(\Se+\bar\Se)|_{U_1}$ has a $\Co^\beta$ local basis $(F_1^*\Coorvec{\nu^1},\dots,F_1^*\Coorvec{\nu^{r+2m}})$.
    
    \item\label{Item::CpxFroPf1::DF1Reg}When $\Se+\bar\Se\in\Co^\beta$, automatically $F_1^*\Coorvec{\nu^1},\dots,F_1^*\Coorvec{\nu^{r+2m}}\in\Co^\beta(U_1;T\Mf|_{U_1})$.
    \item\label{Item::CpxFroPf1::Phi1Reg}
    When $\beta=\LogL$:  for every $0<\eps<1$ there is a neighborhood $\Omega_{1,\eps}'\subseteq\Omega_1'$ of $0$ such that 
    
    \quad $\Phi_1\in\Co^{2-,1-\eps/2}_{\nu,\lambda,\loc}(\Omega_{1,\eps}',\Omega_1''';\Mf)$ and $\nabla_\nu\Phi_1\in\Co^{1-,1-\eps/2}_\loc(\Omega_{1,\eps}',\Omega_1''';T\Mf\otimes\R^{r+2m})$.
    
    When $\alpha=1$: automatically $\Phi_1\in\Co^{2,1-\eps/2}_{\nu,\lambda,\loc}(\Omega_{1,\eps}',\Omega_1''';\Mf)$ and $\nabla_\nu\Phi_1\in\Co^{1,1-\eps/2}_\loc(\Omega_{1,\eps}',\Omega_1''';T\Mf\otimes\R^{r+2m})$.
    
    When $\beta=\logl$: automatically  $\Phi_1\in\Co^{2-,1-}_{\nu,\lambda,\loc}(\Omega_1',\Omega_1''';\Mf)$ and $\nabla_\nu\Phi_1\in\Co^{1-}_\loc(\Omega_1;T\Mf\otimes\R^{r+2m})$.
    
    When $\beta\ge\Lip$: automatically 
    $\Phi_1\in\Co^{\beta+1,\beta}_{\nu,\lambda,\loc}(\Omega_1',\Omega_1''';\Mf)$ and $\nabla_\nu\Phi_1\in\Co^\beta_\loc(\Omega_1;T\Mf\otimes\R^{r+2m})$.
    
    \item\label{Item::CpxFroPf1::F1Reg}When $\beta=\LogL$: for $U_{1,\eps}:=\Phi(\Omega_{1,\eps}'\times\Omega_1''')$, $F_1'\in\Co^\kappa_\loc(U_1;\R^{r+2m})$ and $F_1'''\in\Co^{1-\eps/2}_\loc(U_{1,\eps};\R^q)$.
    
    When $\beta=\logl$: automatically  $F_1'\in\Co^\kappa_\loc(U_1;\R^{r+2m})$ and $F_1'''\in\Co^{1-}_\loc(U_1;\R^q)$.
    
    When $\beta\ge\Lip$: automatically  $F_1'\in\Co^\kappa_\loc(U_1;\R^{r+2m})$ and $F_1'''\in\Co^\beta_\loc(U_1;\R^q)$.
\end{enumerate}

Now $\Phi_1^*\Se=(F_1)_*\Se\le\C T\Omega_1$. By \ref{Item::CpxFroPf1::RFroSpan}, $\Phi_1^*\Se\le\Phi_1^*(\Se+\bar\Se)=\Span(\Coorvec{\nu^1},\dots,\Coorvec{\nu^{r+2m}})|_{\Omega_1}=(\C T\Omega_1')\times\Omega_1'''$.

By Lemma \ref{Lem::ODE::S+IsBundle} since $\Se\in\Co^\alpha$, we know $\Se+\bar\Se\in\Co^\alpha$ as well.
Thus by \ref{Item::CpxFroPf1::DF1Reg} and \ref{Item::CpxFroPf1::Phi1Reg},
\begin{enumerate}[parsep=-0.3ex,label=(S1.\arabic*')]\setcounter{enumi}{2}
    \item\label{Item::CpxFroPf1::DF1Reg'} $F_1^*\Coorvec{\nu^1},\dots,F_1^*\Coorvec{\nu^{r+2m}}\in\Co^\alpha_\loc(U_1;T\Mf|_{U_1})$.
    \item\label{Item::CpxFroPf1::Phi1Reg'} $\Phi_1\in\Co^{\alpha+1}_{\nu,\loc} L^\infty_{\lambda,\loc}(\Omega_1',\Omega_1''';\Mf)$ and $\nabla_\nu\Phi_1\in\Co^{\alpha}_{\nu,\loc} L^\infty_{\lambda,\loc}(\Omega_1',\Omega_1''';T\Mf\otimes\R^{r+2m})$.
\end{enumerate}

Applying Lemma \ref{Lem::ODE::PullBackReg} on \ref{Item::CpxFroPf1::DF1Reg'}:
\begin{enumerate}[parsep=-0.3ex,label=(S1.\arabic*)]\setcounter{enumi}{5}
    \item\label{Item::CpxFroPf1::Phi1*SReg} When $\beta=\LogL$: $\Phi_1^*\Se\le(\C T\Omega'_1)\times\Omega'''_1$ is a $\Co^{\alpha,\alpha(1-\eps/2)}$-subbundle.
    
    When $\beta=\logl$: $\Phi_1^*\Se\le(\C T\Omega'_1)\times\Omega'''_1$ is a $\Co^{\alpha,\alpha-}$-subbundle.

When $\beta\ge\Lip$: $\Phi_1^*\Se\le(\C T\Omega'_1)\times\Omega'''_1$ is a $\Co^{\alpha,\alpha}$-subbundle, i.e. $\Phi_1^*\Se\le\C T(\Omega'_1\times\Omega'''_1)$ is a $\Co^{\alpha}$-subbundle.
\end{enumerate}

Now Theorem \ref{KeyThm::EllipticPara} applies to $\Phi_1^*\Se\le (\C T\Omega_1')\times\Omega_1'''$.
We endow $\R^r\times\C^m\times\R^q$ with standard real and complex coordinate system $(t,z,s)=(t^1,\dots,t^r,z^1,\dots,z^m,s^1,\dots,s^q)$. By Theorem \ref{KeyThm::EllipticPara}, we can find a neighborhood $\Omega_2=\Omega_2'\times\Omega_2''\times\Omega_2'''\subseteq\R^r_t\times\C^m_z\times\R^q_s$ of $(0,0,0)$ (which can depend on $\eps$ when $\beta=\LogL$)
and a $\Co^\alpha$-parameterization $\Phi_2=(\Phi_2',\Phi_2'''):\Omega_2\to\Omega_1'\times\Omega_1'''$, such that by taking $U_2:=\Phi_2(\Omega_2)$ and $F_2=(F_2',F_2'',F_2'''):=\Phi_2^\Inv:U_2\to\R^r\times\C^m\times\R^q$,
\begin{enumerate}[parsep=-0.3ex,label=(S1.\arabic*)]\setcounter{enumi}{6}
    \item\label{Item::CpxFroPf1::Phi20} $\Phi_2(0,0,0)=(0,0)$ and $\nabla_{t,z}\Phi_2$ is continuous.
    \item\label{Item::CpxFroPf1::Phi2''} $\Phi_2'''=\Phi_2'''(s):\Omega_2'''\subseteq\R^q_s\to\Omega'''_1$ is a smooth parameterization that does not depend on $(t,z)$.
    \item\label{Item::CpxFroPf1::F2Chart} $F_2'=F_2'(\nu)\in C^\infty$ does not depend on $\lambda$; and $F_2'''=F_2'''(\lambda)\in C^\infty$ does not depend on $\nu$.
    \item\label{Item::CpxFroPf1::KeySpan} $\Phi_1^*\Se|_{U_2}$ has a local basis $(F_2^*\Coorvec{t^1},\dots,F_2^*\Coorvec{t^r},F_2^*\Coorvec{z^1},\dots,F_2^*\Coorvec{z^m})$.
    
    \item\label{Item::CpxFroPf1::Phi2Reg}
    When $\beta=\LogL$: $\Phi_2'\in\Co^{\alpha+1,(2\alpha-1)(1-\eps/2)}(\Omega_2'\times\Omega_2'',\Omega_2''';\Omega_{1,\eps}')$ and $\nabla_{t,z}\Phi_2'\in\Co^{\alpha,(2\alpha-1)(1-\eps/2)}_{(t,z),s}$.
    
    When $\beta=\logl$ or $\alpha=1$: $\Phi_2'\in\Co^{\alpha+1,(2\alpha-1)-}(\Omega_2'\times\Omega_2'',\Omega_2''';\Omega_1')$ and $\nabla_{t,z}\Phi_2'\in\Co^{\alpha,(2\alpha-1)-}_{(t,z),s}$.
    
    When $\beta\ge\Lip$ and $\alpha\neq1$: $\Phi_2'\in\Co^{\alpha+1,\min(\alpha,2\alpha-1)}(\Omega_2'\times\Omega_2'',\Omega_2''';\Omega_1')$ and $\nabla_{t,z}\Phi_2'\in\Co^{\alpha,\min(\alpha-,2\alpha-1)}_{(t,z),s}$.
\end{enumerate}
Moreover for any $U_2'\subseteq\Omega_1'$ and $U_2'''\subseteq\Omega_1'''$ such that $U_2'\times U_2'''\subseteq U_2$, we have regularity estimates:
\begin{enumerate}[parsep=-0.3ex,label=(S1.\arabic*)]\setcounter{enumi}{11}
    \item\label{Item::CpxFroPf1::F2Reg}When $\beta=\LogL$: $F_2''\in\Co^{\alpha+1,(2\alpha-1)(1-\eps/2)}(U_2',U_2''';\C^m)$.
    
    When $\beta=\logl$ or $\alpha=1$: $F_2''\in\Co^{\alpha+1,(2\alpha-1)-}(U_2',U_2''';\C^m)$;
    
    When $\beta\ge\Lip$ and $\alpha\neq1$: $F_2''\in\Co^{\alpha+1,\min(\alpha,2\alpha-1)}(U_2',U_2''';\C^m)$.
    \item\label{Item::CpxFroPf1::DF2Reg}When $\beta=\LogL$: $F_2^*\Coorvec{t^j},F_2^*\Coorvec{z^k},F_2^*\Coorvec{\bar z^k}\in\Co^{\alpha,(2\alpha-1)(1-\eps/2)}(U_2',U_2''';\C TU_2')$ for $1\le j\le r$, $1\le k\le m$.
    
    When $\beta=\logl$: $F_2^*\Coorvec{t^j},F_2^*\Coorvec{z^k},F_2^*\Coorvec{\bar z^k}\in\Co^{\alpha,(2\alpha-1)-}(U_2',U_2''';\C TU_2')$ for $1\le j\le r$ and $1\le k\le m$.
    
    When $\beta\ge\Lip$: $F_2^*\Coorvec{t^j},F_2^*\Coorvec{z^k},F_2^*\Coorvec{\bar z^k}\in\Co^{\alpha,\min(\alpha-,2\alpha-1)}(U_2',U_2''';\C TU_2')$ for $1\le j\le r$ and $1\le k\le m$.
\end{enumerate}

Now we take $\Omega'\subseteq\Omega_2'$, $\Omega''\subseteq\Omega_2''$, $\Omega'''\subseteq\Omega_2'''$ and $\Omega:=\Omega'\times\Omega''\times\Omega'''$ such that $\Phi_2(\Omega)\subseteq U_2'\times U_2'''\subseteq\Phi_2(\Omega_2)$ holds for some $U_2'\times U_2'''\subseteq U_2\cap(\Omega_1'\times\Omega_1''')$. We define $\Phi:=\Phi_1\circ\Phi_2:\Omega\to\Mf$, $U:=\Phi(\Omega)$ and $F=(F',F'',F'''):=\Phi^\Inv$.

Immediately by \ref{Item::CpxFroPf1::Phi10} and \ref{Item::CpxFroPf1::Phi20} we have \ref{Item::CpxFro::0}. 

By \ref{Item::CpxFroPf1::Phi2''} we have $\Phi(t,z,s)=\Phi_1(\Phi_2'(t,z,s),\Phi_2'''(s))$. By chain rule we have
\begin{equation*}
    \nabla_{t,z}\Phi(t,z,s)=((\nabla_\nu\Phi_1)\circ \Phi_2)(t,z,s))\cdot\nabla_{t,z}\Phi_2'(t,z,s)=(\nabla_\nu\Phi_1)\circ (\Phi'_2(t,z,s),\Phi_2'''(s))\cdot\nabla_{t,z}\Phi_2'(t,z,s).
\end{equation*}
Applying Lemma  \ref{Lem::Hold::CompofMixHold} \ref{Item::Hold::CompofMixHold::Comp} with 
%
\ref{Item::CpxFroPf1::Phi1Reg}, \ref{Item::CpxFroPf1::Phi1Reg'} and \ref{Item::CpxFroPf1::Phi2Reg}:
\begin{itemize}[parsep=-0.3ex]
    \item When $\beta=\LogL$ ($\frac12<\alpha\le1$),  $\Phi\in\Co^{\alpha+1,\min((2\alpha-1)(1-\eps/2),1-\eps/2)}_{(t,z),s}=\Co^{\alpha+1,(2\alpha-1)(1-\eps/2)}_{(t,z),s}\subset\Co^{\alpha+1,2\alpha-1-\eps}_{(t,z),s}$,
    
    and $\nabla_{t,z}\Phi\in\Co^{(1-)\circ(\alpha+1),(1-)\circ(2\alpha-1)(1-\eps/2)}_{(t,z),s}\cdot\Co^{\alpha,(2\alpha-1)(1-\eps/2)}_{(t,z),s}\subseteq\Co^{\alpha,(2\alpha-1)(1-\eps/2)-}_{(t,z),s}\subset\Co^{\alpha,2\alpha-1-\eps}_{(t,z),s}$. 
    
    \item When $\beta=\logl$ ($\frac12<\alpha\le1$), or $\beta\ge\Lip$ and $\alpha=1$:  $\Phi\in\Co^{\alpha+1,\min((2\alpha-1)-,1-)}_{(t,z),s}=\Co^{\alpha+1,(2\alpha-1)-}_{(t,z),s}$,
    
    and $\nabla_{t,z}\Phi\in\Co^{(1-)\circ(\alpha+1),(1-)\circ\min((2\alpha-1)-,1-)}_{(t,z),s}\cdot\Co^{\alpha,(2\alpha-1)-}_{(t,z),s}=\Co^{\alpha,(2\alpha-1)-}_{(t,z),s}$. 
    
    \item When $\beta\ge\Lip$ ($\frac12<\alpha\le\beta$) and $\alpha\neq1$, we have $\Phi\in\Co^{\alpha+1,\min((\beta+1)\circ\min(\alpha,2\alpha-1),\beta)}_{(t,z),s}\subseteq\Co^{\alpha+1,\min(\alpha,2\alpha-1)}_{(t,z),s}$,
    
    and $\nabla_{t,z}\Phi\in\Co^{\beta\circ\alpha,\beta\circ\min(\alpha-,2\alpha-1)}_{(t,z),s}\cdot\Co^{\alpha,\min(\alpha-,2\alpha-1)}_{(t,z),s}=\Co^{\alpha,\min(\alpha-,2\alpha-1)}_{(t,z),s}$. 
\end{itemize}
We see that \ref{Item::CpxFro::PhiReg1}, \ref{Item::CpxFro::PhiReg1.5} \ref{Item::CpxFro::PhiReg2} and \ref{Item::CpxFro::PhiReg3} all hold. In particular we obtain \ref{Item::CpxFro::1}.

By \ref{Item::CpxFroPf1::RFroSpan} and \ref{Item::CpxFroPf1::KeySpan} we have $\Se=(\Phi_1)_*\Span(F_2^*\Coorvec t,F_2^*\Coorvec z)=\Span(F_1^*F_2^*\Coorvec t,F_1^*F_2^*\Coorvec z)=\Span(F^*\Coorvec t,F^*\Coorvec z)$, which is \ref{Item::CpxFro::Span}.

Clearly $F'=F_2'\circ F_1'\in\Co^\infty\circ\Co^\kappa=\Co^\kappa$. By construction $F''=F_2''\circ(F_1',F_1''')$ and $F'''=F_2'''\circ F_1'''$. By \ref{Item::CpxFroPf1::F1Reg} and \ref{Item::CpxFroPf1::F2Chart} we have:
\begin{itemize}[parsep=-0.3ex]
    \item When $\beta=\LogL$: $F''\in\Co^{((2\alpha-1)(1-\eps/2)-)\circ(1-\eps/2)}\subseteq\Co^{(2\alpha-1)(1-\eps/2)^2-}\subset\Co^{2\alpha-1-\eps}$, $F'''\in\Co^{1-\eps/2}\subset\Co^{1-\eps}$.
    
    \item When $\beta=\logl$ ($\frac12<\alpha<1$), or $\beta\ge\Lip$ and $\alpha=1$: $F''\in\Co^{((2\alpha-1)-)\circ(1-)}\subseteq\Co^{(2\alpha-1)-}$ and $F'''\in\Co^{1-}$.
    
    \item When $\beta\ge\Lip$ ($\frac12<\alpha\le\beta$) and $\alpha\neq1$, $F''\in\Co^{\min(\alpha,2\alpha-1)\circ\beta}\subseteq\Co^{\min(\alpha,2\alpha-1)}$ and $F'''\in\Co^\beta$.
\end{itemize}
Thus we get the result \ref{Item::CpxFro::FReg1}, \ref{Item::CpxFro::FReg1.5}, \ref{Item::CpxFro::FReg2} and \ref{Item::CpxFro::FReg3}.

By construction for $1\le j\le r$, $F^*\Coorvec{t^j}=F_1^*F_2^*\Coorvec{t^j}=((F_2^*\Coorvec{t^j})\cdot\nabla_\nu\Phi_1)\circ F_1$, since $F_2^*\Coorvec{t^j}$ are sections of $\Phi_1^*(\Se+\bar\Se)=\Span(\Coorvec{\nu^1},\dots,\Coorvec{\nu^{r+2m}})$. Therefore by \ref{Item::CpxFroPf1::Phi1Reg}, \ref{Item::CpxFroPf1::F1Reg} and \ref{Item::CpxFroPf1::DF1Reg} with compositions, we get:
\begin{itemize}[parsep=-0.3ex]
    \item When $\beta=\LogL$ ($\frac12<\alpha\le1$), $F^*\Coorvec{t^j}\in\Co^{(2\alpha-1)(1-\eps/2)}\circ\Co^{1-\eps/2}\subseteq\Co^{(2\alpha-1)(1-\eps/2)^2}\subset\Co^{2\alpha-1-\eps}$.
    \item When $\beta=\logl$ ($\frac12<\alpha\le1$), or $\beta\ge\Lip$ and $\alpha=1$: $F^*\Coorvec{t^j}\in\Co^{(2\alpha-1)-}\circ\Co^{1-}\subseteq\Co^{(2\alpha-1)-}$.
    \item When $\beta\ge\Lip$ ($\frac12<\alpha\le\beta$) and $\alpha\neq1$, $F^*\Coorvec{t^j}\in\Co^{\min(\alpha-,2\alpha-1)}\circ\Co^\beta\subseteq\Co^{\min(\alpha-,2\alpha-1)}$.
\end{itemize}
The same regularity estimates hold for $F^*\Coorvec{z^k},F^*\Coorvec{\bar z^k}$, $1\le k\le m$ as well. We get \ref{Item::CpxFro::FDReg1}, \ref{Item::CpxFro::FDReg1.5},  \ref{Item::CpxFro::FDReg2} and  \ref{Item::CpxFro::FDReg3}.

Finally by Remark \ref{Rmk::CpxFro::RmkofTrueThm} the result \ref{Item::CpxFro::Span} implies \ref{Item::CpxFro::SBotSpan}, \ref{Item::CpxFro::SCapSpan} and \ref{Item::CpxFro::S+Span}. The result \ref{Item::CpxFro::PhiReg3} and \ref{Item::CpxFro::FReg3} imply \ref{Item::CpxFro::ParaCoor}.
\end{proof}

\begin{proof}[Proof of Theorem \ref{KeyThm::CpxFro2}]
We use a construction that is slightly different from the proof of Theorem \ref{KeyThm::CpxFro1}. Instead of applying Theorem \ref{KeyThm::RealFro::ImprovedLogFro}, we use Theorem \ref{KeyThm::RealFro::BLFro} in the first step.

Still we assume $0<\eps<\alpha-\frac12$ in the case $\beta=\LogL$. Also by Lemma \ref{Lem::ODE::S+IsBundle} (since $\Se\in\Co^\alpha$) we know $\Se+\bar\Se\in\Co^\alpha$ as well, which is useful when $\alpha=1$ and $\beta\in\{\LogL,\logl\}$.

We endow $\R^r\times\R^{2m}\times\R^q$ with standard coordinate system $(\mu,\nu,\lambda)=(\mu^1,\dots,\mu^r,\nu^1,\dots,\nu^{2m},\lambda^1,\dots,\lambda^q)$.

Since $\Se+\bar\Se\in\Co^\beta$ has rank $r+2m$ and $\Se\cap\bar\Se\in\Co^{1+}$ has rank $r$, by Theorem \ref{KeyThm::RealFro::BLFro} we can find a neighborhood $\Omega_1=\Omega'_1\times\Omega''_1\times\Omega'''_1\subseteq\R^r_\mu\times\R^{2m}_\nu\times\R^q_\lambda$ of $(0,0,0)$ and a topological parameterization $\Phi_1:\Omega_1\to\Mf$, which does not depend on $\beta$ and $\gamma$, such that by taking $U_1:=\Phi_1(\Omega_1)$ and $F_1=(F_1',F_1'',F_1'''):U_1\to\R^r\times\R^{2m}\times\R^q$ as $F_1:=\Phi_1^\Inv$, we have
\begin{enumerate}[parsep=-0.3ex,label=(S2.\arabic*)]
    \item\label{Item::CpxFroPf2::Phi10} $\Phi_1(0,0,0)=p_0$ and $\frac{\partial\Phi_1}{\partial\mu^1},\dots,\frac{\partial\Phi_1}{\partial\mu^r},\frac{\partial\Phi_1}{\partial\nu^1},\dots,\frac{\partial\Phi_1}{\partial\nu^{2m}}:\Omega_1\to T\Mf$ are continuous. 
    \item\label{Item::CpxFroPf2::Phi1Span} $(\Se\cap\bar\Se)|_{U_1}=\Span(F_1^*\Coorvec{\mu^1},\dots,F_1^*\Coorvec{\mu^r})$, and $(\Se+\bar\Se)|_{U_1}=\Span(F_1^*\Coorvec{\mu^1},\dots,F_1^*\Coorvec{\mu^r},F_1^*\Coorvec{\nu^1},\dots,F_1^*\Coorvec{\nu^{2m}})$.
    
    \item\label{Item::CpxFroPf2::DF1Reg}When $\Se\cap\bar\Se\in\Co^\gamma$ and $\Se+\bar\Se\in\Co^\beta$, automatically $F_1^*\Coorvec{\mu^1},\dots,F_1^*\Coorvec{\mu^j}\in\Co^\gamma$ and $F_1^*\Coorvec{\nu^k}\in\Co^\alpha$ on $U_1$ for $1\le j\le r$ and $1\le k\le m$.
    
    \item\label{Item::CpxFroPf2::Phi1Reg}When $\beta=\LogL$: $\forall\eps>0$ $\exists$ a neighborhood  $\Omega_{1,\eps}''\subseteq\Omega_1''$ of $0$ such that on $\Omega_{1,\eps}:=\Omega_1'\times\Omega_{1,\eps}''\times\Omega_1'''$,
    
    \quad $\Phi_1\in\Co^{\gamma+1,\min(\gamma,2-),1-\eps/2}_{\mu,\nu,\lambda,\loc}$, $\nabla_\mu\Phi_1\in\Co^{\gamma,\min(\gamma,2-),1-\eps/2}_{\mu,\nu,\lambda,\loc}$ and $\nabla_\nu\Phi_1\in\Co^{\gamma,\min(\gamma-1,1-),(1-\eps/2)\min(\gamma-1,1)}_{\mu,\nu,\lambda,\loc}$.
    
    When $\alpha=1$: automatically $\Phi_1\in\Co^{\gamma+1,2,1-\eps/2}_{\mu,\nu,\lambda,\loc}$, $\nabla_\mu\Phi_1\in\Co^{\gamma,2,1-\eps/2}_{\mu,\nu,\lambda,\loc}$ and $\nabla_\nu\Phi_1\in\Co^{\gamma,1,1-\eps/2}_{\mu,\nu,\lambda,\loc}$ on $\Omega_{1,\eps}$.
    
    When $\beta=\logl$: automatically we have $\Phi_1\in\Co^{\gamma+1,\min(\gamma,2-),1-}_{\mu,\nu,\lambda,\loc}$, $\nabla_\mu\Phi_1\in\Co^{\gamma,\min(\gamma,2-),1-}_{\mu,\nu,\lambda,\loc}$,
    
    \quad and $\nabla_\nu\Phi_1\in\Co^{\gamma,\min(\gamma-1,1-),1-}_{\mu,\nu,\lambda,\loc}$ on $\Omega_1$.

    When $\beta\ge\Lip$: automatically we have $\Phi_1\in\Co^{\gamma+1,\min(\gamma,\beta+1),\min(\gamma,\beta)}_{\mu,\nu,\lambda,\loc}$, $\nabla_\mu\Phi_1\in\Co^{\gamma,\min(\gamma,\beta+1),\min(\gamma,\beta)}_{\mu,\nu,\lambda,\loc}$
    
    \quad and $\nabla_\nu\Phi_1\in\Co^{\gamma,\min(\gamma-1,\beta),\min(\gamma,\beta)}_{\mu,\nu,\lambda,\loc}$ on $\Omega_1$.
    
    \item\label{Item::CpxFroPf2::F1Reg}
    When $\beta=\LogL$: $F_1'\in\Co^\kappa_\loc(U_1;\R^r)$, $F_1''\in\Co^\gamma_\loc(U_1;\R^{2m})$ and $F_1'''\in\Co^{1-\eps/2}_\loc(U_1;\R^q)$.
    
    When $\beta=\logl$: automatically  $F_1'\in\Co^\kappa_\loc(U_1;\R^r)$, $F_1''\in\Co^\gamma_\loc(U_1;\R^{2m})$ and $F_1'''\in\Co^{1-}_\loc(U_1;\R^q)$.
    
    When $\beta\ge\Lip$: automatically  $F_1'\in\Co^\kappa_\loc(U_1;\R^r)$, $F_1''\in\Co^\gamma_\loc(U_1;\R^{2m})$ and $F_1'''\in\Co^\beta_\loc(U_1;\R^q)$.
    
    
\end{enumerate}

Here \ref{Item::CpxFroPf2::DF1Reg} holds because $\gamma-1\ge\alpha$ and $\Se+\bar\Se\in\Co^\alpha$. This is simlar to \ref{Item::CpxFroPf1::DF1Reg'}. We avoid using $\Co^{\min(\gamma-1,\beta)}$ which is not defined if $(\beta,\gamma)=(\logl,2)$.

In application, since $\gamma\ge\alpha+1$, the estimates of $\Phi_1$ and $\nabla_\mu\Phi_1$ \ref{Item::CpxFroPf2::Phi1Reg} can be simplified as
\begin{enumerate}[parsep=-0.3ex,label=(S2.4')]
    \item \label{Item::CpxFroPf2::Phi1Reg'}When $\beta=\LogL$: $\Phi_1\in\Co^{\gamma+1,\alpha+1,1-\eps/2}_{\mu,\nu,\lambda,\loc}$ and $\nabla_\mu\Phi_1\in\Co^{\gamma,\alpha+1,1-\eps/2}_{\mu,\nu,\lambda,\loc}$ on $\Omega_{1,\eps}$.
    
    When $\beta=\logl$: $\Phi_1\in\Co^{\gamma+1,\alpha+1,1-}_{\mu,\nu,\lambda,\loc}$ and $\nabla_\mu\Phi_1\in\Co^{\gamma,\alpha+1,1-}_{\mu,\nu,\lambda,\loc}$ on $\Omega_1$.
    
    When $\beta\ge\Lip$: $\Phi_1\in\Co^{\gamma+1,\alpha+1,\beta}_{\mu,\nu,\lambda,\loc}$ and $\nabla_\mu\Phi_1\in\Co^{\gamma,\alpha+1,\beta}_{\mu,\nu,\lambda,\loc}$ on $\Omega_1$.
\end{enumerate}

By \ref{Item::CpxFroPf2::Phi1Span}, $\Se$ is contained in the span of $F_1^*\Coorvec{\mu^1},\dots,F_1^*\Coorvec{\mu^r},F_1^*\Coorvec{\nu^1},\dots,F_1^*\Coorvec{\nu^{2m}}$. By \ref{Item::CpxFroPf2::DF1Reg} and Lemma \ref{Lem::ODE::PullBackReg}, \ref{Item::CpxFroPf1::Phi1*SReg} still holds.

By the standard linear algebra argument we can find a linear complex coordinate system $w=(w^1,\dots,w^m)$ on $\R^{2m}_\nu$ such that
\begin{equation*}
    \textstyle(\Phi_1^*\Se)_{(0,0)}=\Span\big(\Coorvec{\mu^1}|_{(0,0)},\dots,\Coorvec{\mu^r}|_{(0,0)},\Coorvec{w^1}|_{(0,0)},\dots,\Coorvec{w^m}|_{(0,0)}\big)\le \C T_{(0,0)}(\R^r_\mu\times\C^m_w).
\end{equation*}
In particular $(\Phi_1^*\Se)_{(0,0)}\oplus\Span \big(\Coorvec{\bar w^1}|_{(0,0)},\dots,\Coorvec{\bar w^m}|_{(0,0)})=\C T_{(0,0)}(\R^r_\mu\times\C^m_w)$. Here we have identified $\R^{2m}_\nu$ with $\C^m_w$.

Applying Lemma \ref{Lem::ODE::GoodGen} (by taking $\Mf=\Omega_1'\times\Omega_1''$ and $\Nf=\Omega_1'''$ in the lemma), we can find a neighborhood $W=W'\times W''\times W'''\subseteq\Omega_1'\times\Omega_1''\times\Omega_1'''$ of $(0,0,0)\in\Omega$ ($W\subset\Omega_{1,\eps}$ if $\beta=\LogL$) such that $\Phi_1^*\Se$ has the following local basis on $W$,
\begin{equation*}
    X=\begin{pmatrix}X'\\X''
    \end{pmatrix}\begin{pmatrix}I&&A'\\&I&A''\end{pmatrix}\begin{pmatrix}\partial_\mu\\\partial_w\\\partial_{\bar w}\end{pmatrix}.
\end{equation*}where by \ref{Item::CpxFroPf1::Phi1*SReg}, $A'$ and $A''$ are matrix functions defined in $W$ such that:
\begin{equation}\label{Eqn::CpxFroPf2::RegofA}
\begin{aligned}
    &\text{When }\beta=\LogL:&&\!\!\!\!A'\in\Co^{\alpha,\alpha(1-\frac\eps2)}(W'\times W'',W''';\C^{r\times m}),&&\!\!\!\! A''\in\Co^{\alpha,\alpha(1-\frac\eps2)}(W'\times W'',W''';\C^{m\times m}).
    \\
    &\text{When }\beta=\logl:&&\!\!\!\!A'\in\Co^{\alpha,\alpha-}(W'\times W'',W''';\C^{r\times m}),&&\!\!\!\!A''\in\Co^{\alpha,\alpha-}(W'\times W'',W''';\C^{m\times m}).
    \\
    &\text{When }\beta\ge\Lip:&&\!\!\!\!A'\in\Co^{\alpha}(W;\C^{r\times m}),&&\!\!\!\!A''\in\Co^{\alpha}(W;\C^{m\times m}).
\end{aligned}
\end{equation}

Since $\Coorvec{\mu^1},\dots,\Coorvec{\mu^r}$ are sections of $\Phi_1^*(\Se\cap\bar\Se)\subset\Phi_1^*\Se$, by Lemma \ref{Lem::ODE::GoodGen} \ref{Item::ODE::GoodGen::Uniqueness} we have $A'=0$, i.e. $X_j=\Coorvec{\mu^j}$ for $1\le j\le r$.

By Lemma \ref{Lem::ODE::GoodGen} \ref{Item::ODE::GoodGen::InvComm}, since $\Phi_1^*\Se$ is involutive, we see that $[X_j,X_k]=0$ for $1\le j\le r<k\le r+m$. So $\Coorvec{\mu^j}A''(\mu,w,s)=0$ for $1\le j\le r$ in the domain, which means $A''(\mu,w,s)=A''(w,s)$.

Define a complex tangent subbundle $\Tc\le(\C TW'')\times W'''$ as $\Tc:=\Span(\Coorvec w+A''(w,s)\Coorvec{\bar w})$. We see that $\Tc$ is an involutive subbundle of rank $m$,
\begin{equation}\label{Eqn::CpxFroPf2::Tensor}
    \Tc\cap\bar\Tc=0,\quad\Tc+\bar \Tc=(\C TW'')\times W''',\quad(\Phi_1^*\Se)|_W=(\C TW')\otimes_\C\Tc.
\end{equation}
Moreover by \eqref{Eqn::CpxFroPf2::RegofA}, we have
\begin{enumerate}[parsep=-0.3ex,label=(S2.\arabic*)]\setcounter{enumi}{5}
    \item When $\beta=\LogL$: $\Tc\le(\C TW'')\times W'''$ is a $\Co^{\alpha,\alpha(1-\eps/2)}$-subbundle.
    
    When $\beta=\logl$: $\Tc\le(\C TW'')\times W'''$ is a $\Co^{\alpha,\alpha-}$-subbundle.
    
    When $\beta\ge\Lip$: $\Tc\le(\C TW'')\times W'''$ is a $\Co^{\alpha,\alpha}$-subbundle, i.e. $\Tc\le\C (TW''\times W''')$ is a $\Co^{\alpha}$-subbundle.
\end{enumerate}

Endow $\C^m\times\R^q$ with standard coordinate system $(z,s)=(z^1,\dots,z^m,s^1,\dots,s^q)$. By Theorem \ref{KeyThm::EllipticPara} (by taking $r=0$ in this theorem), we can find a neighborhood $\Omega_2=\Omega_2''\times\Omega_2'''\subseteq\C^m_z\times\R^q_s$ of $(0,0)$ and a $\Co^\alpha$-parameterization $\Phi_2=(\Phi_2'',\Phi_2'''):\Omega_2\to W''\times W'''$, such that by taking $U_2:=\Phi_2(\Omega_2)$ and $F_2=(F_2'',F_2'''):U_2\to\C^m\times\R^q$, $F_2:=\Phi_2^\Inv$, we have
\begin{enumerate}[parsep=-0.3ex,label=(S2.\arabic*)]\setcounter{enumi}{6}
    \item\label{Item::CpxFroPf2::Phi20} $\Phi_2(0,0,0)=(0,0)$  and $\nabla_z\Phi_2$ is continuous..
    \item\label{Item::CpxFroPf2::Phi2''} $\Phi_2'''=\Phi_2'''(s):\Omega_2'''\subseteq\R^q_s\to W'''$ is a smooth parameterization that does not depend on $z$.
    \item\label{Item::CpxFroPf2::F2Chart} $F_2'''=F_2'''(\lambda)\in C^\infty$ does not depend on $(\mu,\nu)$.
    \item\label{Item::CpxFroPf2::Phi2Span} $\Tc|_{U_2}$ has a local basis $(F_2^*\Coorvec{z^1},\dots,F_2^*\Coorvec{z^m})$.
    \item\label{Item::CpxFroPf2::Phi2Reg}
    When $\beta=\LogL$: $\Phi_2'\in\Co^{\alpha+1,(2\alpha-1)(1-\eps/2)}(\Omega_2'',\Omega_2''';W'')$ and $\nabla_z\Phi_2'\in\Co^{\alpha,(2\alpha-1)(1-\eps/2)}(\Omega_2'',\Omega_2''';\C^{2m\times 2m})$.

    When $\beta=\logl$, or $\beta\ge\Lip$ and $\alpha=1$: $\Phi_2'\in\Co^{\alpha+1,(2\alpha-1)-}(\Omega_2'',\Omega_2''';W'')$ and $\nabla_z\Phi_2'\in\Co^{\alpha,(2\alpha-1)-}_{z,s}$.
    
    When $\beta\ge\Lip$ and $\alpha\neq1$: $\Phi_2'\in\Co^{\alpha+1,\min(\alpha,2\alpha-1)}(\Omega_2'',\Omega_2''';W'')$ and $\nabla_z\Phi_2'\in\Co^{\alpha,\min(\alpha-,2\alpha-1)}_{z,s}$.
\end{enumerate}

Moreover for any $U_2''\subseteq W''$ and $U_2'''\subseteq W'''$ such that $U_2''\times U_2'''\subseteq U_2$, we have regularity estimates:
\begin{enumerate}[parsep=-0.3ex,label=(S2.\arabic*)]\setcounter{enumi}{11}
    \item\label{Item::CpxFroPf2::F2Reg}When $\beta=\LogL$: $F_2''\in\Co^{\alpha+1,(2\alpha-1)(1-\eps/2)}(U_2'',U_2''';\C^m)$.
    
    When $\beta=\logl$, or $\beta\ge\Lip$ and $\alpha=1$: $F_2''\in\Co^{\alpha+1,(2\alpha-1)-}(U_2'',U_2''';\C^m)$
    
    When $\beta\ge\Lip$ and $\alpha\neq1$: $F_2''\in\Co^{\alpha+1,\min(\alpha,2\alpha-1)}(U_2'',U_2''';\C^m)$.
    \item\label{Item::CpxFroPf2::DF2Reg}When $\beta=\LogL$: $F_2^*\Coorvec{z^k},F_2^*\Coorvec{\bar z^k}\in\Co^{\alpha,(2\alpha-1)(1-\eps/2)}(U_2'',U_2''';\C TU_2'')$ for $1\le k\le m$.
    
    When $\beta=\logl$: $F_2^*\Coorvec{z^k},F_2^*\Coorvec{\bar z^k}\in\Co^{\alpha,(2\alpha-1)-}(U_2'',U_2''';\C TU_2'')$ for $1\le k\le m$.
    
    When $\beta\ge\Lip$: $F_2^*\Coorvec{z^k},F_2^*\Coorvec{\bar z^k}\in\Co^{\alpha,\min(\alpha-,2\alpha-1)}(U_2',U_2''';\C TU_2'')$ for $1\le k\le m$.
\end{enumerate}

Now we take $\Omega''\subseteq\Omega_2''$,  and $\Omega'''\subseteq\Omega_2'''$ such that $\Phi_2(\Omega''\times\Omega''')\subseteq U_2''\times U_2'''\subseteq\Phi_2(\Omega_2)$ holds for some $U_2'\times U_2'''\subseteq U_2\cap(W''\times W''')$.
We define $\Omega':=W'$ and $\Phi:\Omega\to\Mf$ as 
\begin{equation}\label{Eqn::CpxFroPf2::DefPhi}
    \Phi(t,z,s):=\Phi_1(t,\Phi_2(z,s))=\Phi_1(t,\Phi''_2(z,s),\Phi_2'''(s)),\quad t\in\Omega',z\in\Omega'',s\in\Omega'''.
\end{equation}

We define $U:=\Phi(\Omega)$ and $F=(F',F'',F'''):=\Phi^\Inv$.  Immediately by \ref{Item::CpxFroPf2::Phi10} and \ref{Item::CpxFroPf2::Phi20} we have \ref{Item::CpxFro::0}.

Taking the chain rules on \eqref{Eqn::CpxFroPf2::DefPhi} we have
\begin{equation}\label{Eqn::CpxFroPf2::NablaPhiReg}
    \nabla_t\Phi(t,z,s)=\nabla_\mu\Phi_1(t,\Phi_2''(z,s),\Phi_2'''(s)),\quad\nabla_z\Phi(t,z,s)=(\nabla_\nu\Phi_1)(t,\Phi_2''(z,s),\Phi_2'''(s))\cdot\nabla_z\Phi_2(z,s).
\end{equation}

Based on \ref{Item::CpxFroPf2::Phi1Reg}, \ref{Item::CpxFroPf2::Phi1Reg'} and \ref{Item::CpxFroPf2::Phi2Reg}, applying Lemma \ref{Lem::Hold::CompofMixHold} \ref{Item::Hold::CompofMixHold::Comp} to \eqref{Eqn::CpxFroPf2::DefPhi} and \eqref{Eqn::CpxFroPf2::NablaPhiReg}, we get
\begin{itemize}[parsep=-0.3ex]
    \item When $\beta=\LogL$: $\Phi\in\Co^{\gamma+1,(\alpha+1)\circ(\alpha+1),\min((\alpha+1)\circ((2\alpha-1)(1-\eps/2)),1-\eps/2)}_{t,z,s}\subseteq\Co^{\gamma+1,\alpha+1,2\alpha-1-\eps}_{t,z,s}$,
    
    $\nabla_t\Phi\in \Co^{\gamma,(\alpha+1)\circ(\alpha+1),\min((\alpha+1)\circ((2\alpha-1)(1-\eps/2)),1-\eps/2)}_{t,z,s}\subseteq\Co^{\gamma,\alpha+1,2\alpha-1-\eps}_{t,z,s}$,
     
    $\nabla_z\Phi\in \Co^{\gamma,(\alpha+1)\circ\alpha,\min(\gamma-1,1-)\circ((2\alpha-1)(1-\eps/2))}_{(t,z,s)}\cdot\Co^{\alpha,(2\alpha-1)(1-\eps)}_{z,s}\subseteq\Co^{\gamma,\alpha,(2\alpha-1-\eps)\cdot\min(\gamma-1,1)}_{t,z,s}$.
    \item When $\beta=\logl$, or $\beta\ge\Lip$ and $\alpha=1$: $\Phi\in\Co^{\gamma+1,(\alpha+1)\circ(\alpha+1),\min((\alpha+1)\circ((2\alpha-1)-),1-)}_{t,z,s}\subseteq\Co^{\gamma+1,\alpha+1,(2\alpha-1)-}_{t,z,s}$,
    
    $\nabla_t\Phi\in \Co^{\gamma,(\alpha+1)\circ(\alpha+1),\min((\alpha+1)\circ((2\alpha-1)(1-)),1-)}_{t,z,s}\subseteq\Co^{\gamma,\alpha+1,(2\alpha-1)-}_{t,z,s}$,
     
    $\nabla_z\Phi\in \Co^{\gamma,(\alpha+1)\circ\alpha,\min(\gamma-1,1-)\circ((2\alpha-1)(1-))}_{t,z,s}\cdot\Co^{\alpha,(2\alpha-1)-}_{z,s}\subseteq\Co^{\gamma,\alpha,(2\alpha-1)\cdot\min(\gamma-1,1)-}_{t,z,s}$;

    \item When $\beta\ge\Lip$ and $\alpha\neq1$: $\Phi\in\Co^{\gamma+1,(\alpha+1)\circ(\alpha+1),\min(\alpha,2\alpha-1,\beta)}_{t,z,s}\subseteq\Co^{\gamma+1,\alpha+1,\min(\alpha,2\alpha-1)}_{t,z,s}$,
    
    $\nabla_t\Phi\in \Co^{\gamma,(\alpha+1)\circ(\alpha+1),\min(\alpha,2\alpha-1,\beta)}_{t,z,s}\subseteq\Co^{\gamma,\alpha+1,\min(\alpha,2\alpha-1)}_{t,z,s}$,
    
    $\nabla_z\Phi\in \Co^{\gamma,(\alpha+1)\circ\alpha,\min(\gamma-1,\Lip)\circ\min(\alpha,2\alpha-1))}_{t,z,s}\cdot\Co^{\alpha,\min(\alpha-,2\alpha-1)}_{z,s}\subseteq\Co^{\gamma,\alpha,(\gamma-1)\circ\min(\alpha-,2\alpha-1)}_{t,z,s}$
\end{itemize}

Thus we get \ref{Item::CpxFro::ImprovePhi}, \ref{Item::CpxFro::ImproveDPhiDz} and in particular \ref{Item::CpxFro::1}.

\smallskip
By \eqref{Eqn::CpxFroPf2::DefPhi},  we see that $F=(F',F'',F'''):\to \R^r\times\C^m\times\R^q$ satisfies
\begin{equation*}\label{Eqn::CpxFroPf2::EqnofF}
    F'=F_1'=\id_{\R^r}\circ F_1'\quad F''=F_2''\circ F_1''\quad F'''=F_2'''\circ F_1'''.
\end{equation*}

Therefore by \ref{Item::CpxFroPf2::F1Reg}, \ref{Item::CpxFroPf2::F2Chart} and \ref{Item::CpxFroPf2::F2Reg},
\begin{itemize}[parsep=-0.3ex]
    \item When $\beta=\LogL$: $F'\in\Co^\kappa_\loc(U;\R^r)$, $F''\in\Co^{((2\alpha-1)(1-\eps/2))\circ\gamma}_\loc\subset\Co^{2\alpha-1-\eps}_\loc$, $F'''\in\Co^{1-\eps/2}_\loc\subset\Co^{1-\eps}_\loc$.
    
    \item When $\beta=\logl$, or $\beta\ge\Lip$ and $\alpha=1$: $F'\in\Co^\kappa_\loc$, $F''\in\Co^{((2\alpha-1)-)\circ\gamma}_\loc\subset\Co^{(2\alpha-1)-}_\loc$, $F'''\in\Co^{1-}_\loc$.
    
    \item When $\beta\ge\Lip$ and $\alpha=1$: $F'\in\Co^\kappa_\loc$, $F''\in\Co^{\min(\alpha,2\alpha-1)\circ\gamma}_\loc\subset\Co^{\min(\alpha,2\alpha-1)}_\loc$, $F'''\in\Co^\beta_\loc$.
\end{itemize}

This gives \ref{Item::CpxFro::FReg1}, \ref{Item::CpxFro::FReg1.5}, \ref{Item::CpxFro::FReg2} and \ref{Item::CpxFro::FReg3}.

By \eqref{Eqn::CpxFroPf2::NablaPhiReg}, we get $\Coorvec{t^j}\Phi(t,z,s)=\frac{\partial\Phi_1}{\partial\mu^j}(t,\Phi_2''(z,s),\Phi_2'''(s))=(\Phi_1)_*\Coorvec{\mu^j}|_{\Phi(t,z,s)}$ for $1\le j\le r$. So $F^*\Coorvec{t^j}=F_1^*\Coorvec{\mu^j}\in\Co^\gamma$ for $1\le j\le r$, finishing the proof of \ref{Item::CpxFro::ImproveDDt}. 

By \eqref{Eqn::CpxFroPf2::Tensor} and \ref{Item::CpxFroPf2::Phi2Span} we see that $\Phi_1^*\Se=\Span(\Coorvec{\mu},F_2^*\Coorvec{z})$, so by writing $\id_{\R^r}:\R^r_\mu\to\R^r_t$ as the identity we have $\Phi_1^*\Se=\Span(\id_{\R^r}^*\Coorvec{t},F_2^*\Coorvec{z})$. Taking pushforward of $\Phi_1$ we get $\Se|_U=F_1^*\Span(\id_{\R^r}^*\Coorvec{t},F_2^*\Coorvec{z})|_U=\Span(F^*\Coorvec t,F^*\Coorvec z)$, finishing the proof of \ref{Item::CpxFro::Span}.

Using \eqref{Eqn::CpxFroPf2::EqnofF} we have $F=\widehat F_2\circ F_1$ where $\widehat F_2:=(\id_{\R^r},F_2)$. By \ref{Item::CpxFroPf2::DF2Reg} and \ref{Item::CpxFroPf2::DF1Reg} we see that:
\begin{itemize}[parsep=-0.3ex]
    \item When $\beta=\LogL$: $\widehat F_2^*\Coorvec{z^1},\widehat F_2^*\Coorvec{\bar z^1},\dots,\widehat F_2^*\Coorvec{z^m},\widehat F_2^*\Coorvec{\bar z^m}$ are $\Co^{(2\alpha-1)(1-\frac\eps2)}$-linear combinations of $\Coorvec{\nu^1},\dots,\Coorvec{\nu^{2m}}$. 
    
    \quad While $F_1^*\Coorvec{\nu^1},\dots,F_1^*\Coorvec{\nu^{2m}}$ are $\Co^{1-\eps/2}$-vector fields.
    
    \item When $\beta=\logl$: $\widehat F_2^*\Coorvec{z^1},\widehat F_2^*\Coorvec{\bar z^1},\dots,\widehat F_2^*\Coorvec{z^m},\widehat F_2^*\Coorvec{\bar z^m}$ are $\Co^{(2\alpha-1)-}$-linear combinations of $\Coorvec{\nu^1},\dots,\Coorvec{\nu^{2m}}$.
    
    \quad While $F_1^*\Coorvec{\nu^1},\dots,F_1^*\Coorvec{\nu^{2m}}$ are $\Co^{1-}$-vector fields.
    
    \item When $\beta\ge\Lip$: $\widehat F_2^*\Coorvec{z^1},\widehat F_2^*\Coorvec{\bar z^1},\dots,\widehat F_2^*\Coorvec{z^m},\widehat F_2^*\Coorvec{\bar z^m}$ are $\Co^{\min(\alpha-,2\alpha-1)}$-linear combinations of $\Coorvec{\nu^1},\dots,\Coorvec{\nu^{2m}}$.    
    
    \quad While $F_1^*\Coorvec{\nu^1},\dots,F_1^*\Coorvec{\nu^{2m}}$ are $\Co^\beta$-vector fields.
    
\end{itemize}

Say $\widehat F_2^*\Coorvec{z^k}=:\sum_{j=1}^{2m}g_j^k\Coorvec{\nu^j}$. Then $F^*\Coorvec{z^k}=\sum_{j=1}^{2m}(g_j^k\circ F_1)\cdot(F_1^*\Coorvec{\nu^j})$. By \ref{Item::CpxFroPf2::F1Reg}, we see that:
\begin{itemize}[parsep=-0.3ex]
    \item When $\beta=\LogL$: $g_j^k\in\Co^{(2\alpha-1)(1-\eps/2)\circ(1-\eps/2)}_\loc\subset\Co^{2\alpha-1-\eps}_\loc$, thus $F^*\Coorvec{z^k}\in\Co^{2\alpha-1-\eps}_\loc$.
    
    \item When $\beta=\logl$: $g_j^k\in\Co^{((2\alpha-1)-)\circ(1-)}_\loc=\Co^{(2\alpha-1)-}_\loc$, thus $F^*\Coorvec{z^k}\in\Co^{(2\alpha-1)-}_\loc$.
    
    \item When $\beta\ge\Lip$: $g_j^k\in\Co^{\min(\alpha-,2\alpha-1)\circ\beta}_\loc=\Co^{\min(\alpha-,2\alpha-1}_\loc$, thus $F^*\Coorvec{z^k}\in\Co^{\min(\alpha-,2\alpha-1)}_\loc$.
\end{itemize}
The argument for $F^*\Coorvec{\bar z^k}$ are the same. Therefore we get \ref{Item::CpxFro::FDReg1}, \ref{Item::CpxFro::FDReg1.5}, \ref{Item::CpxFro::FDReg2} and \ref{Item::CpxFro::FDReg3}.

The rest of the results follow from Remark \ref{Rmk::CpxFro::RmkofTrueThm} as well. 
\end{proof}

\section{Examples for Some Sharpness}\label{Section::Examples}
In this chapter we discuss certain sharpness of the coordinate charts and coordinate vector fields for the Frobenius-type structures.

\setcounter{subsection}{-1}
\subsection{An overview to estimates of complex Frobenius structures}\label{Section::ExampleOverview}
We only consider the case $\alpha>1$ here.

Fixed four numbers $\alpha,\beta,\gamma,\kappa>1$ such that $\alpha\le\min(\beta,\gamma)$ and $\kappa\ge\max(\beta,\gamma)+1$. Consider all possible complex Frobenius structures $\Se\le \C T\Mf$ among all $\Co^\kappa$-manifolds $\Mf$ such that $\Se\in\Co^\alpha$, $\Se+\bar\Se\in\Co^\beta$ and $\Se\cap\bar\Se\in\Co^\gamma$. 

Fix a point $p\in \Mf$, let $F=(F',F'',F'''):U\subseteq\Mf\to\R^r_t\times\C^m_z\times\R^{n}_s$ be any desired $C^1$-coordinate chart near a fixed point $p$ that represents $\Se$, that is, $F^*\Coorvec{t^1},\dots,F^*\Coorvec{t^r},F^*\Coorvec{z^1},\dots,F^*\Coorvec{z^m}$ span $\Se$ in $U$. 

Recall the notation $\Co^{\mu+}=\bigcup_{\nu>\mu}\Co^\nu$ for $\mu\in\R$.

We list the best H\"older-Zygmund regularities the six objects: $F'$, $F''$, $F'''$, $F^*\Coorvec t$, $F^*\Coorvec z$, $F^*\Coorvec s$, as follows.

\begin{itemize}[parsep=-0.1ex]
    \item $F'\notin\Co^{\kappa+}$: $\Mf$ is $\Co^\kappa$, its coordinate functions are at most $\Co^\kappa$, so $F'\in\Co^\kappa$ is optimal.
    
    \item
 $F'''\notin\Co^{\beta+}$: This is due to the sharpness of real Frobenius theorem that $\Se+\bar\Se\in\Co^\beta$ only guarantees $F'''$ to be $\Co^\beta$. See Proposition \ref{Prop::Exmp::RealFroisSharp} below.

    \item
$F^*\Coorvec t\notin\Co^{\gamma+}$: $F^*\Coorvec t$ spans $\Se\cap\bar\Se$ which is $\Co^\gamma$. If $F^*\Coorvec t\in\Co^{\gamma+}$ then $\Se\cap\bar\Se\in\Co^{\gamma+}$ as well.

    \item
$F''\notin\Co^{\alpha+}$ and $F^*\Coorvec s\notin\Co^{(\alpha-1)+}$: The argument is similar to the sharpness of real Frobenius theorem. See Proposition \ref{Prop::Exmp::Sharpdds} below.

    \item
$F^*\Coorvec z\notin\Co^\alpha$: This is the counter intuitive part in Theorem \ref{MainThm::CpxFro}. See Section \ref{Section::Sharpddz} for the complete deduction and Proposition \ref{Prop::Exmp::SharpddzRed}, Theorem \ref{Thm::Exmp::ProofExampleddz} for the proof. The idea comes from the proof of Theorem \ref{MainThm::CounterNN}.
\end{itemize}

When $\gamma\ge\alpha+1$ (that $\Se\cap\bar\Se$ is at least 1 order more regular than $\Se$), by Theorem \ref{KeyThm::CpxFro1}, all six regularity estimates above are sharp. 

When $\alpha\le\gamma<\alpha+1$, from Theorem \ref{KeyThm::CpxFro1} we only have $F^*\Coorvec t\in\Co^{\alpha-}$, the rest of five estimates are still sharp.

When $\alpha\le\gamma<\alpha+1$ and $\gamma>2$, using Theorem \ref{KeyThm::CpxFro1} we can find a $\Co^\alpha$-coordinate chart $F=(F',F'',F''')$ to $\R^r_t\times\C^m_z\times\R^n_s$ such that $F''\in\Co^\alpha$, $F^*\Coorvec z\in\Co^{\alpha-}$ but $F^*\Coorvec t\in\Co^{\alpha-}$. If we use Theorem \ref{KeyThm::CpxFro2}, then we can get $F^*\Coorvec t\in\Co^{\gamma}$ which is sharp, but neither $F''\in\Co^{\gamma-1}$, $F^*\Coorvec z\in\Co^{(\gamma-1)-}$ nor $F^*\Coorvec s\in\Co^{\gamma-2}$ is sharp. 

We do not know whether we can choose suitable coordinate chart $F$ so that $F^*\Coorvec t\in\Co^{\gamma}$ and $F''\in\Co^\alpha$ holds in the same time.

\subsection{Some elementary examples}
In this subsection we use the conventions $f_u=\partial_uf$, $f_v=\partial_vf$, etc.

For  parameterizations of real tangent subbundles, we can reduce the discussion to the regularity of ODE flows, because of the following:
\begin{lem}\label{Lem::Exmp::ParaToFlow}
Let $X(x,y)=(1,a(x,y))$ be a continuous real vector field in $\R^2_{x,y}$. Let $\V\le T\R^2$ be the rank $1$ real tangent subbundle spanned by $X$. 

Suppose $\Omega\subseteq\R^2_{t,s}$ is a neighborhood of $(0,0)$ and  $\Phi=(\phi,\psi):\Omega\subseteq\R^2_{t,s}\to\R^2_{x,y}$ is a continuous map which is homeomorphism onto its image, such that  $\frac{\partial\Phi}{\partial t}:\Omega\to\R^2$ $(T\R^2)$ is continuous,  and $\frac{\partial\Phi}{\partial t}(t,s)$ spans $\V_{\Phi(t,s)}$ for each $(t,s)\in\Omega$.
Then 
\begin{enumerate}[parsep=-0.3ex,label=(\roman*)]
    \item\label{Item::Exmp::ParaToFlow::GInv} the map $G(t,s):=(\phi(t,s),\psi(0,s))$ is locally homeomorphism near $(0,0)$. And the map $\Psi(u,v):=(u,\psi\circ G^\Inv(u,v))$ satisfies
    \begin{itemize}[nolistsep]
        \item $\Psi$ is homeomorphism onto its image.
        \item $\Psi(0,v)=(0,v)$ and $\frac{\partial\Psi}{\partial u}(u,v)=X(\Psi(u,v))$.
    \end{itemize}
    \item\label{Item::Exmp::ParaToFlow::Calpha}  Suppose in addition there is a  $\alpha\in\{1-\}\cup(1,\infty)$ such that $\Phi$ and $\Phi^\Inv$ are $\Co^\alpha$ near $(0,0)$. Then $\Psi,\Psi^\Inv\in\Co^\alpha$ near $(0,0)$ as well.
\end{enumerate}
\end{lem}
 The case $\alpha=1$ is true as well but we do not use that in applications.
 
When $X$ has ODE uniqueness, we see that $\Psi(t,s)=\exp_X(t,s)$. In general, for example when $X$ is H\"older, $\Psi$ is just an integral flow for $X$.

\begin{proof}
By assumption $\frac{\partial\phi}{\partial t}\in C^0(\Omega)$ is non-vanishing, so for each $v$, the map $\phi(\cdot,v)$ has a $C^1$-inverse. 

Define $G_1(t,s):=(\phi(t,s),s)$ and $h(u,v):=\phi(\cdot,u)^\Inv(v)$. Clearly $h$ is continuous and is $C^1$ in $u$. Therefore $G_1^\Inv(u,v)=(h(u,v),v)$ is continuous on $G(\Omega)$ and is $C^1$ in $u$. We see that $G_1$ is homeomorphism onto its image.

Let $\Psi_1(u,v):=(u,\psi\circ G_1^\Inv(u,v))=\Phi\circ G_1^\Inv(u,v)$. By computation, for $(u,v)$ in the domain $$\textstyle\frac{\partial\Psi_1}{\partial u}(u,v)=(1,\frac{\partial\psi}{\partial t}(h(u,v),v)\cdot \frac{\partial h}{\partial u}(u,v))=\frac{\partial h}{\partial u}(u,v)\cdot\frac{\partial \Phi}{\partial t}(G_1^\Inv(u,v)).$$

Since $\frac{\partial h}{\partial u}$ is non-vanishing, we see that $\frac{\partial\Psi_1}{\partial u}(u,v)$ spans $\V_{\Psi_1(u,v)}$.
On the other hand $\Psi_1(u,v)=(u,\ast)$, we must have $\frac{\partial\Psi_1}{\partial u}(u,v)=(1,\ast)$. Therefore $\frac{\partial\Psi_1}{\partial u}(u,v)=X(\Psi_1(u,v))$.

Since $\Phi$ and $G$ are both homeomorphism, we see that $\Psi_1$ is also homeomorphism, thus $v\mapsto\Psi_1(0,v)$ is a topological embedded curve. On the other hand $\Psi_1(0,v)=(0,\psi(0,v))$, we see that $v\mapsto\psi(0,v)$ is injective continuous, which is homeomorphism between neighborhoods of $0\in\R^1$ and $\psi(0)\in\R^1$. 

Let $G_2(u,v):=(u,\psi(0,v))$, we see that $G_2$ is homeomorphism onto its image, $C^1$ in $u$, and $G=G_2\circ G_1$. Therefore $\Psi(u,v)=\Psi_1(u,\psi(0,\cdot)^\Inv(v))=\Phi\circ G^\Inv(u,v)$ satisfies
\begin{itemize}[nolistsep]
    \item $\Psi$ is homeomorphism  onto its image, and $C^1$ in $u$.
    \item $\Psi(0,v)=\Psi_1(0,\psi(0,\cdot)^\Inv(v))=(0,v)$.
    \item $\frac{\partial\Psi}{\partial u}(u,v)=\frac{\partial\Psi_1}{\partial u}(u,\psi(0,\cdot)^\Inv(v))=X\circ \Psi_1(u,\psi(0,\cdot)^\Inv(v))=X(\Psi(u,v))$.
\end{itemize}

Therefore $\Psi$ is as desired, finishing the proof of \ref{Item::Exmp::ParaToFlow::GInv}.

\medskip\noindent Proof of \ref{Item::Exmp::ParaToFlow::Calpha}: By the assumption $G_1, G_2$ are both $\Co^\alpha$. 

When $\alpha>1$, we see that $\nabla\Phi(0,0)\in\R^{2\times 2}$ has full rank 2, thus $\nabla G_1(0,0),\nabla G_2(0,0)\in\R^{2\times 2}$ are both of full rank as well. By the Inverse Function Theorem, we see that $G_1,G_2\in\Co^\alpha$ has $\Co^\alpha$-inverses near $(0,0)$, thus $\Psi$ is $\Co^\alpha$ and has $\Co^\alpha$-inverse near $(0,0)$.

When $\alpha=1-$, the proof of $G_1\in\Co^{1-}$ is similar to the proof of Lemma \ref{Lem::Hold::QPIFTLow}.

By shrinking $\Omega$ if necessary we can assume $C_0:=(\inf_\Omega|\frac{\partial \phi}{\partial t}|)^{-1}<\infty$. Therefore for $t_1,t_2,s_1,s_2\in\R$ that satisfy $(t_1,s_1),(t_2,s_1),(t_2,s_2)\in\Phi(\Omega)$,
\begin{align*}
    |h(t_1,s_1)-h(t_2,s_2)|\le&|h(t_1,s_1)-h(t_2,s_1)|+|h(t_2,s_1)-h(t_2,s_2)|
    \\\le& C_0|t_1-t_2|+C_0|\phi(h(t_2,s_1),s_1)-\phi(h(t_2,s_2),s_1)|&\textstyle(\sup|\frac{\partial h}{\partial u}|=(\inf|\frac{\partial \phi}{\partial t}|)^{-1}=C_0)
    \\
    =&C_0|t_1-t_2|+C_0|t_2-\phi(h(t_2,s_2),s_1)|&(\phi(h(\cdot,s_1),s_1)=\id)
    \\
    =&C_0|t_1-t_2|+C_0|\phi(h(t_2,s_2),s_2)-\phi(h(t_2,s_2),s_1)|&(\phi(h(\cdot,s_2),s_2)=\id)
    \\\le& C_0|t_1-t_2|+C_0\|\phi\|_{\Co^{1-\eps}}|s_1-s_2|^{1-\eps},&\forall \eps>0.\quad (\text{since }\phi\in \Co^{1-})
\end{align*}
We conclude that $G_1^\Inv\in\Co^{1-}$ in the domain. Therefore by Corollary \ref{Cor::Hold::CompOp} \ref{Item::Hold::CompOp::Log1}, $\Psi_1=\Phi\circ G_1^\Inv$ and $\Psi_1^\Inv=G_1\circ\Phi^\Inv$ are both $\Co^{1-}$.

Now $\Psi_1^\Inv(0,v)=(0,\psi(0,\cdot)^\Inv(v))$ is a $\Co^{1-}$-map, thus $\psi(0,\cdot)^\Inv\in\Co^{1-}$, which means $G_2^\Inv\in\Co^{1-}$. Clearly $G_2\in\Co^{1-}$ (since $\psi\in\Co^{1-}$). By Corollary \ref{Cor::Hold::CompOp} \ref{Item::Hold::CompOp::Log1}, $\Psi=\Psi_1\circ G_2^\Inv$ and $\Psi^\Inv=G_2\circ\Psi_1^\Inv$ are both $\Co^{1-}$ as well, finishing the proof.
\end{proof}

For the non-$C^1$ case we have the following examples for the real Frobenius theorems.

\begin{cor}\label{Cor::Exmp::SharpRealFro}
For a continuous function $f_j:\R\to\R$, we set a vector field $X_j(x,y)=(1,f_j(y))$ in $\R^2_{x,t}$ and define $\V_j\le T\R^2$ be a tangent subbundle spanned by $X_j$.
\begin{enumerate}[parsep=-0.3ex,label=(\roman*)]
    \item\label{Item::Exmp::SharpRealFro::<1} Let $0<\alpha<1$ and let $f_1(y):=-|y|^\alpha$. Then $\V_1$ is not locally integrable via parameterizations (recall Definition \ref{Defn::DisInv::LocInvPara}).
    \item\label{Item::Exmp::SharpRealFro::LogL} Let $f_2(y):=y\log|y|$. Then $\V_2$ is a $\Co^1\subset\Co^\LogL$ subbundle. And if $\Phi$ is a parameterization for $\V_2$ in Definition \ref{Defn::DisInv::LocInvPara} near $(0,0)$, then either $\Phi\notin\Co^{1-}$ or $\Phi^\Inv\notin \Co^{1-}$.
\end{enumerate}

\end{cor}
\begin{proof}
By direct computation, the equation $\Psi_1(0,s)=s$, $\frac{\partial\Psi_1}{\partial t}(t,s)=X_1(\Psi_1(t,s))$ has forward uniqueness solution
\begin{equation*}
    \Psi_1(t,s)=\begin{cases}(t,(|s|^{1-\alpha}-(1-\alpha)t)^\frac1{1-\alpha})&0\le t\le \tfrac{|s|^{1-\alpha}}{1-\alpha},\\0&t\ge \tfrac{|s|^{1-\alpha}}{1-\alpha}.\end{cases}
\end{equation*}

Clearly $\Psi_1$ is not injective in a neighborhood of $(0,0)$, which cannot be locally homeomorphism near $(0,0)$. By Lemma \ref{Lem::Exmp::ParaToFlow}, there is no topological parameterization $\Phi_1(t,s)$ near $(0,0)$ such that $\frac{\partial\Phi_1}{\partial t}$ span $\V_1$. Thus $\V_1$ is not locally integrable via parameterizations, giving \ref{Item::Exmp::SharpRealFro::<1}.

For $X_2$, the ODE flow gives
\begin{equation*}
    e^{tX_2}(0,s)=(t,|s|^{e^u}\sgn s).
\end{equation*}
This map is not $\Co^{1-}$ is any neighborhood of $(0,0)$. Therefore by Lemma \ref{Lem::Exmp::ParaToFlow} \ref{Item::Exmp::ParaToFlow::Calpha}, if $\Phi_2$ is a topological parameterization for $\V_2$ near $(0,0)$, then either $\Phi_2\notin\Co^{1-}$ or $\Phi_2^\Inv\notin\Co^{1-}$.
\end{proof}






For the $C^1$ case, we have the following result for the real Frobenius theorem.

\begin{prop}\label{Prop::Exmp::RealFroisSharp}
    Let $\beta>1$, and let $(u,v)$ be the standard coordinate system for $\R^2$. We define $\V\le T\R^2$ as
    $$\textstyle \V:=\Span_\R T,\quad T:=\Coorvec u+f(v)\Coorvec v,\quad\text{where }f(v):=\max(0,v)^\beta.$$
    
    Suppose $F=(F',F'''):U\subseteq\R^2_{u,v}\to\R^1_t\times\R^1_s$ is a $C^1$-coordinate chart near $(0,0)$ such that $F^*\Coorvec t$ spans $\V|_U$, then $\partial_vF'''\notin\Co^{(\beta-1)+}_\loc(U)$. In particular $F'''\notin\Co^{\beta+}_\loc(U;\R^2)$.
\end{prop}Note that $\V$ has rank 1 so is automatically involutive.
\begin{proof}
By direct computation
\begin{equation*}
    e^{tT}(u,v)=\begin{cases}\big(u+t,\frac v{\left(1-(\beta-1)v^{\beta-1}t\right)^\frac1{\beta-1}}\big)&\text{for }v>0\text{ and }-\infty<t<(\beta-1)v^{1-\beta},
    \\
    (u+t,v)&\text{for }v\le0\text{ and }t\in\R.
    \end{cases}
\end{equation*}

First we show that for every $\delta\neq0$, the function
\begin{equation}\label{Eqn::ODE::DefGDelta}
    g^\delta:I_\delta\to\R,\quad g^\delta(v):=\begin{cases}\frac v{\left(1-(\beta-1)\delta v^{\beta-1}\right)^\frac1{\beta-1}},&v>0,\\v,&v\le0,\end{cases}\quad\text{where }I_\delta:=\begin{cases}
    \left(-\infty,(\frac\delta{\beta-1})^\frac1{1-\beta}\right),&\delta>0,\\\R,&\delta<0,
    \end{cases}
\end{equation}
 is not $\Co^{\beta+}$ near $v=0$. Note that $e^{u T}(0,v)=(u,g^u(v))$ for every $u$.

Indeed by Taylor's expansion,
\begin{equation*}
    \frac v{\left(1-(\beta-1)\delta v^{\beta-1}\right)^\frac1{\beta-1}}=v+\delta v^\beta+\sum_{k=2}^\infty{-\frac1{\beta-1}\choose k}\left(-(\beta-1)\delta\right)^k\cdot v^{k\beta-k+1},\quad\text{converging for } 0\le v<(\tfrac\delta{\beta-1})^\frac1{1-\beta}.
\end{equation*}
Hence $g^\delta(v)=v+\max(0,v)^\beta+\sum_{k=2}^\infty{-\frac1{\beta-1}\choose k}\left(-(\beta-1)\delta\right)^k\cdot\max(0,v)^{k\beta-k+1}$ holds near $v=0$.

Now near $v=0$, $\sum_{k=2}^\infty{-\frac1{\beta-1}\choose k}\left(-(\beta-1)\delta\right)^k\cdot\max(0,v)^{k\beta-k+1}$ is a $\Co^{2\beta-1}_\loc\subsetneq\Co^{\beta+}_\loc$-function, but $\max(0,v)^\beta\notin\Co^{\beta+}$. So $g^\delta\notin\Co^{\beta+}$ near $v=0$.
Therefore the function $(u,v)\mapsto g^u(v)$ is not $\Co^{\beta+}$ near $(u,v)=(0,0)$.

By assumption $F^*\Coorvec t$ spans $\V|_U$. So $\V^\bot|_U=\Span dF''$, which means that we have a transport equation
\begin{equation}\label{Eqn::Exmp::RealFroisSharp::EqnofF'''}
    \textstyle TF'''(u,v)=0\quad\text{i.e.}\quad\partial_uF'''(u,v)=-f(v)\cdot\partial_vF'''(u,v),\quad\text{for all }(u,v)\in U.
\end{equation}

Since $F$ is a $C^1$-chart, we have $(F_u'''(0,0),F_v'''(0,0))\neq(0,0)$. By \eqref{Eqn::Exmp::RealFroisSharp::EqnofF'''} and $f(0)=0$, we have $F_u'''(0,0)=0$. Thus $F_v'''(0,0)\neq0$. 

Therefore $h(v):=F'''(0,v)$ is a $C^1$-function that has non-vanishing derivatives in a neighborhood of $v=0$. By the Inverse Function Theorem, on a smaller neighborhood of $v=0$, $h$ is a $C^1$-diffeomorphism onto its image.

\medskip
Suppose by contrast $F'''_v\in\Co^{(\beta-1)+}$ in a neighborhood of $(u,v)=(0,0)$. By \eqref{Eqn::Exmp::RealFroisSharp::EqnofF'''} since $f\in\Co^\beta_\loc\subset\Co^{(\beta-1)+}$, we have $F'''_u\in\Co^{(\beta-1)+}$ as well, which means $F'''\in\Co^{\beta+}$ near $(u,v)=(0,0)$. Since $h(v)=F'''(0,v)$, we have $h\in\Co^{\beta+}$ in the domain and by the Inverse Function Theorem $h^\Inv\in\Co^{\beta+}$.  

So $h^\Inv\circ F''(u,v)$ is a $\Co^{\beta+}$-function defined in a neighborhood of $(u,v)=(0,0)$. By \eqref{Eqn::Exmp::RealFroisSharp::EqnofF'''} and \eqref{Eqn::ODE::DefGDelta}, $F'''(u,v)=F'''(e^{-uT}(u,v))=F'''(0,g^{-u}(v))$ holds a neighborhood of $0$.
So $h^\Inv\circ F'''(u,v)=g^{-u}(v)$, contradicting to the fact that $(u,v)\mapsto g^u(v)$ is not $\Co^{\beta+}$ near $(0,0)$. This concludes that $F_v'''\notin\Co^{(\beta-1)+}$ near $(u,v)=(0,0)$.
Therefore $F'''\notin\Co^{\beta+}$ and we finish the proof.
\end{proof}

For the complex Frobenius theorem the following example give the sharpness of $F''\notin\Co^{\alpha+}$ and $F^*\Coorvec s\notin\Co^{(\alpha-1)+}$.

\begin{prop}\label{Prop::Exmp::Sharpdds}Endow $\R^3$ with standard coordinate system $(u,v,\theta)$. Let $\alpha>0$. We define $\Se\le\C T\R^3$ as
\begin{equation*}
    \Se:=\Span_\C Z,\quad Z:=\partial_u+ie^{f(\theta)}\partial_v\quad\text{where }f(\theta):=\max(0,\theta)^\alpha.
\end{equation*}

Suppose $F=(F'',F'''):U\subseteq\R^3_{u,v,\theta}\to\C^1_z\times\R^1_s$ is a $C^1$-coordinate chart near $0\in\R^3$ such that that $F(0)=(0,0)$ and $F^*\Coorvec z$ spans $\Se|_U$. Then both $F''\notin\Co^{\alpha+}$ and $F^*\Coorvec s\notin\Co^{(\alpha-1)+}$ near $(0,0,0)$.
\end{prop}
By assumption $Z\in\Co^\alpha$, $\Se$ is involutive since it has rank 1, and $\Se+\bar\Se=\Span(\Coorvec u,\Coorvec v)$ is involutive as well. So $\Se$ is a $\Co^\alpha$-complex Frobenius structure.

\begin{proof}
By assumption $\Se+\bar\Se|_U=\Span(\partial_u,\partial_v)|_U=\Span(F^*\Coorvec z,F^*\Coorvec{\bar z})=(\Span F^*ds)^\bot=(\Span dF''')^\bot$, so $\partial_uF'''=\partial_vF'''=0$ which means $F'''(u,v,\theta)=F'''(\theta)$. So $dF'''=F'''_\theta d\theta$ and $F'''_\theta (\theta)$ is non-vanishing in a neighborhood of $\theta=0$.

Now $\Se^\bot|_U=(\Span Z)^\bot|_U=\Span(F^*d\bar z,F^*ds)=\Span (d\bar F'',dF''')$, so
\begin{equation*}
    Z\bar F''=\bar F''_u+ie^{f(\theta)}\bar F''_v=0,
    \quad\Rightarrow\quad F''_v=-ie^{-f(\theta)} F''_u,\quad\text{in }U.
\end{equation*}
Here we use $F''_u=\partial_u F''$, $F''_v=\partial_v F''$ and $F''_\theta=\partial_\theta F''$. So
\begin{equation}\label{Eqn::Exmp::Sharpdds::ClassicalExp1}
\begin{pmatrix}F^*dz\\F^*d\bar z\\F^*ds\end{pmatrix}=
\begin{pmatrix}
F''_u&-F''_u&F''_\theta\\
\bar F''_u&\bar F''_u&\bar F''_\theta\\
&&F'''_\theta 
\end{pmatrix}
\begin{pmatrix}du\\-ie^{-f(\theta)}dv\\d\theta\end{pmatrix}\ \Rightarrow\
\begin{pmatrix}
\partial_u\\ ie^{f(\theta)}\partial_v\\\partial_\theta
\end{pmatrix}=
\begin{pmatrix}
F''_u&\bar F''_u&\\
-F''_u&\bar F''_u\\
F''_\theta&\bar F''_\theta&F'''_\theta 
\end{pmatrix}
\begin{pmatrix}F^*\partial_z\\F^*\partial_{\bar z}\\F^*\partial_s\end{pmatrix}.
\end{equation}
Applying Cramer's rule for inverse matrices to the right hand equation of \eqref{Eqn::Exmp::Sharpdds::ClassicalExp1}, we have
\begin{equation}\label{Eqn::Exmp::RealFroisSharp::EqnofF*dds}
    F^*\Coorvec s=\frac1{F'''_\theta (\theta)}\bigg(\Coorvec\theta-\frac{\bar F''_\theta}{\bar F''_u}\Big(\Coorvec u+ie^{f(\theta)}\Coorvec v\Big)-\frac{F''_\theta}{F''_u}\Big(\Coorvec u-ie^{f(\theta)}\Coorvec v\Big)\bigg),\quad\text{in }U.
\end{equation}

In particular we see that $F''_u$ and $F''_v=-ie^{-f}F''_u$ are both non-vanishing in the domain $U$.

For fixed $\theta\in\R$,  $\partial_u+ie^{f(\theta)}\partial_v$ is scaled $\overline{\partial}$-operator on $\R^2_{u,v}$ that annihilates $F''(\cdot,\theta_0)$. Therefore by a scaling of the Cauchy Integral Formula, for every $\tilde V'\times V''\Subset U$, $(u,v)\in \tilde V'$, $\theta\in V''$ and $k\ge0$,
\begin{equation*}
    F''(u,v,\theta)=\frac1{2\pi i}\int_{\partial V'}\frac{F''(\xi,\eta,\theta)(d\xi+ie^{-f(\theta)}d\eta)}{(\xi-u)+ie^{-f(\theta)}(\eta-v)},\quad \frac{\partial^kF''(u,v,\theta)}{\partial u^k}=\frac{k!}{2\pi i}\int_{\partial V'}\frac{F''(\xi,\eta,\theta)(d\xi+ie^{-f(\theta)}d\eta)}{((\xi-u)+ie^{-f(\theta)}(\eta-v))^{k+1}}.
\end{equation*}

Since $F''\in C^0( \tilde V'\times V'';\C)$, for every precompact subset $ V'\Subset \tilde V'$ we have
\begin{equation*}
    \|\partial_u^kF''\|_{L^\infty(V'\times V'';\C^{2^k})}\le k!\dist(V',\partial \tilde V')^{-k-1}\|F''\|_{L^\infty( \tilde V'\times V'';\C)}<\infty,\quad\forall k\ge0.
\end{equation*}
Using $F_v''=-ie^{-f(\theta)}F_u''$, we see that $\|\nabla_{u,v}^kF''\|_{L^\infty( V'\times V'';\C^{2^k})}<\infty$ for every $V'\times V''\Subset U$, which means $F''\in \Co^\infty_{u,v}L^\infty_\theta(V',V'';\C)$.

Let $\gamma>0$ be such that $F''\in\Co^\gamma_\loc(U;\C)$. Since $F''\in\Co^\infty_{u,v} L^\infty_\theta(V',V'';\C)$ for every $V'\times V''\Subset U$, by Lemma \ref{Lem::Hold::GradMixHold} we have the following:
\begin{equation}\label{Eqn::Exmp::Sharpdds::Claim}
    \text{For every }\gamma>0,\text{ if }F''\in\Co^\gamma_\loc(U;\C),\text{ then }\nabla_{u,v}F''\in\Co^{\infty,\gamma-}_{(u,v),s}(V',V'';\C^2)\text{ holds for all }V'\times V''\Subset U.
\end{equation} 

Now suppose by contrast that $F''\in\Co^{\alpha+}_\loc(U;\C)$. Say $F''\in\Co^{\alpha+2\eps}_\loc(U;\C)$ for some $\eps>0$. Take $\gamma=\alpha+\eps$ in \eqref{Eqn::Exmp::Sharpdds::Claim}, we have $F''_u,F''_v\in\Co^{\alpha+\eps}(V'\times V'';\C)$ for every precompact neighborhood $V'\times V''\Subset U$ of $0\in\R^2_{u,v}\times\R^1_\theta$.

Since $F''_u,F''_v$ are both non-vanishing, we have $-ie^{f(\theta)}=F''_v/F_u''\in\Co^{\alpha+\eps}$ near $(0,0,0)$. However $-ie^{f}\notin\Co^{\alpha+}$ near $\theta=0$ since $f\notin\Co^{\alpha+}$ near $\theta=0$. Contradiction!

This concludes the proof of $F''\notin\Co^{\alpha+}$ near $(0,0,0)$.

\medskip
To show $F^*\Coorvec s\notin\Co^{(\alpha-1)+}$ near $(0,0,0)$, we can assume $F'''\in\Co^{\alpha+}$, since otherwise $F'''_\theta \notin\Co^{(\alpha-1)+}$ in a neighborhood of $\theta=0$ and by \eqref{Eqn::Exmp::RealFroisSharp::EqnofF*dds} the $\Coorvec\theta$-component of $F^*\Coorvec s$ is not $\Co^{(\alpha-1)+}$ already.


Let $\gamma>0$ be the largest number such that $F''\in\Co^{\gamma-}$ near $(0,0,0)$. By assumption $F''\in C^1$ so $\gamma\ge1$. Since $F''\notin\Co^{\alpha+}$ near $(0,0,0)$, we see that $\gamma\le\alpha$.

By \eqref{Eqn::Exmp::Sharpdds::Claim}, we know $F''_u,F''_v\in\Co^{\gamma-\eps}$ near $(0,0,0)$ for every $\eps>0$, in particular $F''_u,F''_v\in\Co^{\gamma-\frac12}$ near $(0,0,0)$.

We see that $F''_\theta\notin\Co^{(\gamma-1)+}$ near $(0,0,0)$. Suppose otherwise, i.e. $F''_\theta\in\Co^{\gamma-1+\eps}$ for some $0<\eps<\frac12$. Since $F''_u,F''_v\in\Co^{\gamma-\frac12}$, we have $F''\in\Co^{\gamma+\eps}$. This contradicts to the maximality of $\gamma$.

Now we have $F''_\theta/F''_u\notin\Co^{(\gamma-1)+}$ for such $\gamma\ge1$, in particular $F''_\theta/F''_u\notin\Co^{(\alpha-1)+}$. Since $Z=\partial_u+ie^f\partial_v$ and $\bar Z=\partial_u-ie^f\partial_v$ are both $\Co^\alpha$-vector fields in $\R^3$ and we already assumed $F'''_\theta \in\Co^{(\alpha-1)+}$, using \eqref{Eqn::Exmp::RealFroisSharp::EqnofF*dds} again we conclude that $F^*\Coorvec s\notin\Co^{(\alpha-1)+}$ near $(0,0,0)$. This completes the whole proof.
\end{proof}

\subsection{A counterexample to $C^k$-regularity for the Newlander-Nirenberg theorem}\label{Section::CountNN}
In this subsection we prove Theorem \ref{MainThm::CounterNN}. 



For convenience we use the viewpoint of eigenbundle of $J$: Set $\Se_n:=\coprod_p\{v\in \C T_p\R^{2n}:J_pv=iv\}$. So $J\Coorvec{w^j}=i\Coorvec{w^j}$ for all $j$ iff $d\bar w^1,\dots,d\bar w^n$ spans $\V_n^\bot|_U\le\C T^*\R^{2n}|_U$. And $J$ is integrable if and only if $\Se_n$ is involutive. See \cite[Chapter 1]{Involutive} for details.

\medskip
First we can restrict our focus to the 1-dimensional case:
\begin{proof}[Proof of Theorem \ref{MainThm::CounterNN}, 1-dimension $\Rightarrow$ n-dimension] Suppose $J_1$ is a $C^k$-almost complex structure on $\C^1_{z^1}$ (not compatible with the standard complex structure), such that near $0$, there is no $C^{k,1}$-complex coordinate $\varphi$ satisfying $\Se_1^\bot=\Span d\bar \varphi$ in the domain. Here $\Se_1^\bot$ is the dual eigenbundle of $J_1$. 


\medskip
Denote by $\theta=\theta(z^1)$ a $C^k$ 1-form on $\C^1_{z^1}$ that spans $\Se_1^\bot$. 

Consider $\R^{2n}\simeq\C^n_{(z^1,\dots,z^n)}$. We identify $\theta$ as the $C^k$ 1-form on $\R^{2n}$.
Take an $n$-dim almost complex structure on $\R^{2n}$ such that the dual of eigenbundle $\Se_n^\bot$ is spanned by $\theta,d\bar z^2,\dots,d\bar z^n$. 

In other words, $\Se_n$ is the ``tensor'' of $\Se_1$ with the standard complex structure of $\C^{n-1}_{(z^2,\dots,z^n)}$.

If $w=(w^1,\dots,w^n)$ is a corresponding $C^{k,1}$-complex chart for $\Se_n$ near $0$, then there is a $1\le j_0\le n$ such that $d\bar w^{j_0}\not\equiv0\pmod{d\bar z^2,\dots,d\bar z^n}$ near $0$. In other words, we have linear combinations $d\bar w^{j_0}=\lambda_1\theta+\lambda_2d\bar z^2+\dots+\lambda_nd\bar z^n$ for some non-vanishing function $\lambda_1(z^1,\dots,z^n)$ near $z=0$.

Therefore $w^{j_0}(\cdot,0^{n-1})$ is a complex coordinate chart defined near $z^1=0\in\R^2$ whose differential spans $\Se_n^\bot|_{\R^2_{z^1}}\cong\Se_1^\bot$ near $0$. By our assumption on $\Se_1$, we have $w^{j_0}(\cdot,0)\notin C^{k,1}$. So $w^{j_0}\notin C^{k,1}$, which means $w\notin C^{k,1}$.
\end{proof}

Now we focus on the one-dimensional case. Note that a 1-dim structure is automatically integrable.

Fix $k\ge1$. Define an almost complex structure by setting its eigenbundle $\Se_1\le\C T\R^2$ equals to the span of $\Coorvec z+a(z)\Coorvec{\bar z}$, where $a\in C^k(\R^2;\C)$ has compact support that satisfies the following:

\medskip
\begin{enumerate}[nolistsep,label=(\roman*)]
    \item\label{Item::Exmp::CountNN::1} $a\in C^\infty_\loc(\R^2\backslash\{0\};\C)$;
    \item\label{Item::Exmp::CountNN::2} $\partial_z^{-1}\partial_{\bar z}a\notin C^{k-1,1}$ near 0;
    \item\label{Item::Exmp::CountNN::3} $za\in C^{k+1}(\R^2;\C)$ and $z^{-1}a\in C^{k-1}(\R^2;\C)$ (which implies $a(z)=o(|z|)$ as $z\to0$);
    \item\label{Item::Exmp::CountNN::4} $\supp a\subset\B^2$;
    \item\label{Item::Exmp::CountNN::5} $\|a\|_{C^0}<\delta_0$ for some small enough $\delta_0>0$ (take $\delta_0=10^{-1}$ will be okay).
\end{enumerate}

\medskip
Here we take $\partial_z^{-1}$ to be the \textit{conjugated Cauchy-Green operator} on the unit disk $\B^2\subset\C^1$: $$\displaystyle \partial_z^{-1}\phi(z)=\partial_{z,\B^2}^{-1}\phi(z):=\frac1\pi\int_{\B^2}\frac{\phi(\xi+i\eta)d\xi d\eta}{\bar z-\xi+i\eta}.$$

We use notation $\partial_z^{-1}$ because it is an right inverse of $\partial_z$. For $f$ supported in $\B^2\subset\C$, $\partial_z^{-1}f=4\partial_{\bar z}\Ga\ast f$, where $\Ga$ is 2-dimensional Newtonian potential in \eqref{Eqn::Hold::Newtonian}. Therefore by Lemma \ref{Lem::Hold::GreensOp} $\partial_z^{-1}:C^{m,\beta}_c(\B^2;\C)\to C^{m+1,\beta}(\B^2;\C)$ is bounded linear for all $m\in\Z_{\ge0}$, $0<\beta<1$. 

\medskip
We can take  $\supp a\subset\B^2$ such that when $|z|<\frac12$, 
\begin{equation}\label{Item::Exmp::CountNN::a}
    \textstyle a(z):=\frac1{100}\bar z^{k+1}\partial_z\big((-\log|z|)^{\frac12}\big)=\frac1{100}\partial_z\big(\bar z^{k+1}(-\log|z|)^{\frac12}\big).
\end{equation}

Note that for this $a$ we have $a(z)=O\big(|z|^k(-\log|z|)^{-\frac12}\big)=o(|z|^k)$ as $z\to0$.

\begin{remark}
Roughly speaking, Property \ref{Item::Exmp::CountNN::1} $\operatorname{Singsupp}a=\{0\}$ says that the regularity of $a(z)$ corresponds to the vanishing  order of $a$ at $0$. To some degree, by multiplying with $a(z)$, a function gains some regularity at the origin.
\end{remark}

 We check Property \ref{Item::Exmp::CountNN::2} that $\partial_z^{-1}\partial_{\bar z}a\notin C^{k-1,1}$ here.
\begin{lem}\label{Lem::Exmp::CountNN::lem}
Let $a(z)$ be given by \eqref{Item::Exmp::CountNN::a}, and let $\chi\in C_c^\infty(\frac12\B^2)$ satisfies $\chi\equiv1$ in a neighborhood of 0. Then $\partial_z^{-1}(\chi a_{\bar z})\notin C^{k-1,1}$ near $z=0$.
\end{lem}
\begin{proof}Denote $b(z):=\bar z^{k+1}(-\log|z|)^\frac12$, so $b\in C^\infty(\frac12\B^2\backslash\{0\};\C)$, $\chi a=\frac1{100}\chi b_z$ and $\chi a_{\bar z}=\frac1{100}\chi b_{z\bar z}$. 

First we show that $\partial_z^{-1}(\chi b_{z\bar z})-b_{\bar z}\in C^\infty(\frac12\B^2;\C)$. 
We write $$\partial_z^{-1}(\chi b_{z\bar z})- b_{\bar z}=\partial_z^{-1}\partial_z(\chi b_{\bar z})-\chi b_{\bar z}-(1-\chi)b_{\bar z}-\partial_z^{-1}(\chi_zb_{\bar z}).$$

Since $\partial_z\partial_z^{-1}\partial_z(\chi b_{\bar z})=\partial_z(\chi b_{\bar z})$, we know $\partial_z^{-1}\partial_z(\chi b_{\bar z})-\chi b_{\bar z}$ is anti-holomorphic. By Cauchy integral formula we get $\partial_z^{-1}\partial_z(\chi b_{\bar z})-\chi b_{\bar z}\in C^\infty$.

By assumption $0\notin\supp\chi_z$, $0\notin\supp(1-\chi)$ and $b\in C^\infty(\frac12\B^2\backslash\{0\};\C)$, we know $\chi_zb_{\bar z}\in C_c^\infty$ and $(1-\chi)b_{\bar z}\in C^\infty(\frac12\B^2;\C)$. 

Note that when acting on functions supported in the unit disk, $\partial_z^{-1}=\frac1{\pi\bar z}\ast(\cdot)$ is a convolution operator with kernel $\frac1{\pi\bar z}\in L^1$, so $\partial_z^{-1}(\chi_zb_{\bar z})=\frac1{\pi\bar z}\ast(\chi_zb_{\bar z})\in C^\infty$.

It remains to show $b_{\bar z}\notin C^{k-1,1}$ near $0$. Indeed one has $$\partial_{\bar z}^{k+1}(\bar z^{k+1}(-\log|z|)^\frac12)=(k+1)!(-\log|z|)^{\frac12}+O(1),\quad\text{as }z\to0,$$ because by Leibniz rule
\begin{align*}
    &\textstyle\partial_{\bar z}^k\partial_{\bar z}\big(\bar z^{k+1}(-\log|z|)^{\frac12}\big)=\sum_{j=0}^{k+1}{k+1\choose j}\partial_{\bar z}^{k+1-j}(z^{k+1})\cdot\partial_{\bar z}^j(-\log|z|)^{\frac12}\\
    =&\textstyle(k+1)!(-\log|z|)^{\frac12}+\sum_{j=1}^{k+1}O(z^j)O\big(z^{-j}(-\log|z|)^{-\frac12}\big)=(k+1)!(-\log|z|)^{\frac12}+O\big((-\log|z|)^{-\frac12}\big).\qedhere
\end{align*}
\end{proof}
\medskip Now assume $w:\tilde U\subset\R^2\to\C$ is a 1-dim $C^1$-complex coordinate chart defined near $0$ that represents $\Se_1$, then $\Span d\bar w=\Span (d\bar z-adz)|_{\tilde U}=\Se_1^\bot|_{\tilde U}$. So $d\bar w=\bar w_{\bar z}d\bar z+\bar w_zdz=\bar w_{\bar z}(d\bar z-adz)$, that is,
$$\displaystyle\frac{\partial w}{\partial\bar z}(z)+\bar a(z)\frac{\partial w}{\partial z}(z)=0,\qquad z\in {\tilde U}.$$

\begin{remark}
It is worth noticing that $\partial_w\neq \partial_z+a\partial_{\bar z}$. Indeed $\partial_w$ is only a scalar multiple of $\partial_z+a\partial_{\bar z}$.
\end{remark}

Note that $w_z(0)\neq0$ because $(d\bar z-adz)|_0=d\bar z|_0\in\Span d\bar w|_0$. So by multiplying $w_{ z}(0)^{-1}$, we can assume $w_{z}(0)=1$ without loss of generality.
Then $f:=\log\partial_zw$ is a well-defined function in a smaller neighborhood $U\subset\tilde U$ of $0$, which solves 
\begin{equation}\label{Eqn::Exmp::CountNN::MainEQN}
    \frac{\partial f}{\partial\bar z}(z)+\overline{a(z)}\frac{\partial f}{\partial z}(z)=-\frac{\partial\bar a}{\partial z}(z)\quad\Big(=-\overline{\frac{\partial a}{\partial\bar z}(z)} \Big),\qquad z\in U.
\end{equation}

Property \ref{Item::Exmp::CountNN::5} indicates that the operator $\partial_{\bar z}+\bar a\partial_z$ is a first order elliptic operator. Therefore we can consider a second order divergence form elliptic operator $$L:=\partial_z(\partial_{\bar z}+\bar a\partial_z)$$ whose coefficients are $C^k$ globally and are $C^\infty$ outside the origin.

By the classical Schauder's estimate (see \cite[Theorem 4.2]{TaylorPDE3}, or \cite[Chapter 6 \& 8]{GilbargTrudinger}), we have the following:
\begin{lem}[Schauder's interior estimate]\label{Lem::Exmp::CountNN::Schauder} Assume $u,\psi\in C^{0,1}(\B^2;\C)$ satisfy $Lu=\psi_z$. Let $U\subset\R^2$ be a neighborhood of $0$. The following hold:
\begin{enumerate}[nolistsep,label=(\alph*)]
    \item\label{Item::Exmp::CountNN::S1}If $\psi\in C^\infty_\loc(U\backslash\{0\};\C)$, then $u\in C^\infty_\loc(U\backslash\{0\};\C)$.
    \item\label{Item::Exmp::CountNN::S2}If $\psi\in C^{k-1,1}(\R^2;\C)$, then $u\in C^{k,1-\eps}_\loc(\R^2;\C)$ for all $0<\eps<1$.
\end{enumerate}
\end{lem}

Our Theorem \ref{MainThm::CounterNN} for 1-dim case is done by the following proposition:

\begin{prop}\label{Prop::Exmp::CountNN::MainProp}
For any neighborhood $U\subset\R^2_{z}$ of $0$, there is no $f\in C^{k-1,1}(U;\C)$ solving \eqref{Eqn::Exmp::CountNN::MainEQN}.
\end{prop}
\begin{proof}
Suppose there is a neighborhood $U\subset\frac12\B^2$ of the origin, and a solution $f\in C^{k-1,1}(U;\C)$ to \eqref{Eqn::Exmp::CountNN::MainEQN}.

Applying Lemma \ref{Lem::Exmp::CountNN::Schauder} \ref{Item::Exmp::CountNN::S1} on \eqref{Eqn::Exmp::CountNN::MainEQN} with $u=f$ and $\psi=-\bar a_z$, we know that $f\in C^\infty(U\backslash\{0\};\C)$.

Take $\chi\in C_c^\infty(U;[0,1])$ such that $\chi\equiv1$ in a smaller neighborhood of $0$. Denote $$g(z):=\chi(z) f(z),\qquad h(z):=\bar zg(z).$$ So $g,h$ are $C^{k-1,1}$-functions defined in $\R^2$ that are also \textbf{smooth away from 0}, and satisfy the following:
\begin{align}\label{Eqn::Exmp::CountNN::EQNg}
    g_{\bar z}+\bar ag_z&=\chi_{\bar z}f+\chi_z\bar af-\chi\bar a_z,
\\\label{Eqn::Exmp::CountNN::EQNh}
     h_{\bar z}+\bar ah_z&=g+\bar z(\chi_{\bar z}f+\chi_z\bar af)-\chi \bar z\bar a_z.
\end{align}

By construction $\nabla\chi\equiv0$ holds in a neighborhood of $0$, so $\chi_{\bar z}f+\chi_z\bar af\in C_c^\infty(U;\C)$.

\medskip
Under our assumption that $f\in C^{k-1,1}(U;\C)$, then the key is to show that 
\begin{enumerate}[parsep=-0.3ex,label=(\Roman*)]
    \item $h\in C^{k,1-\eps}_c(\R^2;\C)$, $\forall \eps\in(0,1)$. This implies:
    \item\label{Item::Exmp::CountNN::C2} $\bar ag_z\in C^{k-1,1-\eps}_c(\R^2;\C)$, $\forall \eps\in(0,1)$.
\end{enumerate}

\medskip\noindent Proof of (I): By assumption $\bar z(\chi_{\bar z}f+\chi_z\bar af)\in C^\infty_c$, $g\in C^{k-1,1}$, and by Property \ref{Item::Exmp::CountNN::3}, $\chi\bar z\bar a_z\in C^k(\R^2;\C)$. 

Applying Lemma \ref{Lem::Exmp::CountNN::Schauder} \ref{Item::Exmp::CountNN::S2} to \eqref{Eqn::Exmp::CountNN::EQNh}, with $u=h$ and $\psi=g+\bar z(\chi_{\bar z}f+\chi_z\bar a f)-\chi\bar z\bar a_z\in C^k_c(\R^2;\C)$, we get $h\in C^{k,1-\eps}(\B^2;\C)$, for all $0<\eps<1$.

\medskip\noindent Proof of (II): When $k\ge2$,  we know  $z^{-1}a\in C^{k-1}$ for Property \ref{Item::Exmp::CountNN::3}. So for any $\eps\in(0,1)$, one has $\bar ag_z\in C^{k-1,1-\eps}$ because $$\nabla_{z,\bar z}(\bar ag_z)=\bar a\cdot(\partial_zg_z,\partial_{\bar z}g_z)+g_z\nabla_{z,\bar z}\bar a=\bar z^{-1}\bar a\cdot(h_{zz},h_{z\bar z}-g_z)+O(C^{k-2,1})\in C^{k-2,1-\eps}.$$

When $k=1$, for any $z_1,z_2\in\R^2\backslash\{0\}$, note that $g_{\bar z}$ is smooth outside the origin, so $\bar ag_z\in C^{0,1-\eps}$:
    \begin{align*}
    &|\bar ag_z(z_1)-\bar ag_z(z_2)|\le|\bar a(z_1)||g_z(z_1)-g_z(z_2)|+|\bar a(z_1)-\bar a(z_2)||g_z(z_2)|
    \\
    \le&|\bar z_1^{-1}\bar a(z_1)||\bar z_1g_z(z_1)-\bar z_2g_z(z_2)|+|\bar z_1^{-1}\bar a(z_1)||\bar z_1-\bar z_2||g_z(z_2)|+|\bar a(z_1)-\bar a(z_2)||g_z(z_2)|
    \\
    \le&\|z^{-1}a\|_{C^0}\|\nabla h\|_{C^{1-\eps}}|z_1-z_2|^{1-\eps}+\|z^{-1}a\|_{C^0}\|g\|_{C^{0,1}}|z_1-z_2|+\|a\|_{C^1}\|g\|_{C^{0,1}}|z_1-z_2|.
\end{align*}
So for either case of $k$, we have $\bar ag_z\in C^{k-1,1-\eps}(\R^2;\C)$ for all $\eps\in(0,1)$.

\medskip
\textit{Here as a remark, $|\bar ag_z(z_1)-\bar ag_z(z_2)|\lesssim_{a,g,h}|z_1-z_2|^{1-\eps}$ still makes sense when $z_1$ or $z_2=0$, though $\bar g_z(0)$ may not be defined. Indeed $\lim\limits_{z\to0}\bar ag_z(z)$ exists because $\bar ag_z$ itself has bounded $C^{0,1-\eps}$-oscillation on $\B^2\backslash\{0\}$, and then the limit defines the value of $\bar ag_z$ at $z=0$.}

\medskip 
Based on consequence \ref{Item::Exmp::CountNN::C2}, for \eqref{Eqn::Exmp::CountNN::EQNg} we have 
\begin{equation}\label{Eqn::Exmp::CountNN::FinalEqn}
    g=(g-\partial_{\bar z}^{-1}g_{\bar z})+\partial_{\bar z}^{-1}(\chi_{\bar z}f+\chi_z\bar af)-\partial_{\bar z}^{-1}(\bar ag_z)-\partial_{\bar z}^{-1}(\chi\bar a_z),\qquad\text{in }\B^2.
\end{equation}

The right hand side of \eqref{Eqn::Exmp::CountNN::FinalEqn} consists of four terms, the first to the third are all $C^k$, while the last one is not $C^{k-1,1}$. We explain these as follows:
\begin{itemize}[parsep=-0.3ex]
    \item Since $\partial_z(g-\partial_{\bar z}^{-1}g_{\bar z})=0$ in $\B^2$, we know $g-\partial_{\bar z}^{-1}g_{\bar z}$ is anti-holomorphic, which is smooth in $\B^2$.
    \item By assumption $\chi_{\bar z}f+\chi_z\bar af\in C^\infty_c(\B^2;\C)$, so $\partial_{\bar z}^{-1}(\chi_{\bar z}f+\chi_z\bar af)\in C^\infty(\B^2;\C)$ as well.
    \item By consequence \ref{Item::Exmp::CountNN::C2} $\bar ag_z\in C^{k-1,1-\eps}(\R^2;\C)$, for all $\eps\in(0,1)$, so $\partial_{\bar z}^{-1}(\bar ag_z)\in  C^{k,1-\eps}\subset C^k(\B^2;\C)$.
    \item However by Lemma \ref{Lem::Exmp::CountNN::lem}, $\partial_{\bar z}^{-1}(\chi\bar a_z)\notin C^{k-1,1}$ near 0.
\end{itemize}

Combining each term to the right hand side of \eqref{Eqn::Exmp::CountNN::FinalEqn}, we know $g\notin C^{k-1,1}$ near 0. Contradiction!
\end{proof}

\begin{remark}
The key to the proof is the non-surjectivity of $\partial_z:C^k\to C^{k-1}$, which we use to  construct a function $a(z)$ such that $a(0)=0$, $\operatorname{Singsupp}a=\{0\}$, and $\partial_z^{-1}\partial_{\bar z}a\notin C^k$. 
\end{remark}

\subsection{The sharpness of $F^*\Coorvec z$ in Theorem \ref{MainThm::CpxFro}: Proof of Theorem \ref{MainThm::Sharpddz}}\label{Section::Sharpddz}

In this subsection we will construct a rank-1 complex Frobenius structure $\Se$ on $\C^1\times\R^1$, such that $F^*\Coorvec z\notin\Co^\alpha$ whenever $F$ is a coordinate chart representing $\Se$ near the origin. 

We sketch the idea to the proof of Theorem \ref{MainThm::Sharpddz} here. It is adapted from Theorem \ref{MainThm::CounterNN}, where we prove that there is a $C^k$-integrable almost complex structure that does not admit a $C^{k+1}$-coordinate chart.

Endow $\C^1\times\R^1$ with standard coordinate system $(w,\theta)$. We consider a subbundle $\Se\le\C T(\C^1\times\R^1)$ as
\begin{equation*}
    \Se:=\Span(\partial_w+a(w,\theta)\partial_{\bar w}),\quad\text{where }a\in\Co^\alpha(\C^1\times\R^1;\C)\text{ such that }\|a\|_{L^\infty}<1.
\end{equation*}

Clearly $\Se$ is involutive and $\Se+\bar\Se=\Span(\partial_w,\partial_{\bar w})=(\C T\C_w^1)\times\R^1_\theta$ since $\rank \Se=1$ and $\sup_{w,\theta}|a(w,\theta)|<1$. So $\Se$ is a $\Co^\alpha$-complex Frobenius structure.

\smallskip
\noindent\textit{Step 1} (Proposition \ref{Prop::Exmp::SharpddzRed}): We show that, the existence of a coordinate chart $F\in C^1$ near $(0,0)$ such that $F^*\Coorvec z\in\Co^\alpha$ spans $\Se$, is equivalent to the existence of $f\in\Co^\alpha$ near $(w,\theta)=(0,0)$ that solves
\begin{equation}\label{Eqn::Exmp::SharpddzRed::aPDE}
    (\partial_w+a(w,\theta)\partial_{\bar w})f=-\partial_{\bar w}a,\quad\text{or equivalently}\quad f_w+(af)_{\bar w}-a_{\bar w}f=-a_{\bar w}.
\end{equation}

\smallskip
\noindent\textit{Step 2} (Proposition \ref{Prop::Exmp::Expa}): We construct an $a(w,\theta)$ such that $a$ is smooth away from $\{w=0\}=\{0\}\times\R^1_\theta$, and the complex Hilbert transform (also known as the Beurling transform) $\partial_w^{-1}\partial_{\bar w} a\notin\Co^\alpha$ near $w=0$.

\smallskip
\noindent\textit{Step 3} (Theorem \ref{Thm::Exmp::ProofExampleddz}): We show that there is no $f\in\Co^\alpha$ defined near $(w,\theta)=(0,0)$ that solves \eqref{Eqn::Exmp::SharpddzRed::aPDE}. This can be done by showing that:
\begin{itemize}[parsep=-0.3ex]
    \item For every solution $f$ to the homogeneous equation $(\partial_w+a\partial_{\bar w})f=0$ with $f(0,\theta)=0$, $f$ must be $\Co^\alpha$.
    \item There exists a non-$\Co^\alpha$ function $f$ such that $f(0,\theta)=0$ and $f$ solves \eqref{Eqn::Exmp::SharpddzRed::aPDE}.
\end{itemize}
Thus, by the superposition principle, every solution to \eqref{Eqn::Exmp::SharpddzRed::aPDE} cannot be $\Co^\alpha$. By Step 1 we complete the proof.

To make the result more general, we consider the case of mixed regularity, where $a\in\Co^{\alpha,\beta}_{w,\theta}$ for arbitrary $\alpha>\frac12$ and $\beta>0$ (recall again from Lemma \ref{Lem::Hold::CharMixHold} that $\Co^\alpha_{(w,\theta)}=\Co^{\alpha,\alpha}_{w,\theta}$). See Theorem \ref{Thm::Exmp::SharpSum}.

\medskip
We now start the proof. We use the notation $f_w=\partial_wf$ etc.

\begin{prop}\label{Prop::Exmp::SharpddzRed}
Let $\alpha>\frac12,\beta>0$, let $a\in\Co^{\alpha,\beta}(\C^1_w,\R^1_\theta;\C)$ be such that $\sup|a|<1$, and let $\Se=\Span(\partial_w+a(w,\theta)\partial_{\bar w})$ be a $\Co^{\min(\alpha,\beta)}$-complex Frobenius structure with rank $1$. The following \ref{Item::Exmp::SharpddzRed::Chart} and \ref{Item::Exmp::SharpddzRed::PDE} are equivalent:
\begin{enumerate}[parsep=-0.3ex,label=(\Alph*)]
    \item\label{Item::Exmp::SharpddzRed::Chart} There is a continuous map $F:U_1\times V_1\subseteq\C^1_w\times\R^1_\theta\to\C^1_z\times\R^1_s$ near $(0,0)$ such that
    \begin{enumerate}[nolistsep,label=(A.\arabic*)]
        \item\label{Item::Exmp::SharpddzRed::Chart::Reg} $F$ is homeomorphic to its image, $\nabla_zF\in C^0(U_1\times V_1;\C^{2\times 2})$ and $\nabla_w(F^\Inv)\in C^0(U_1\times V_1;\C^{2\times 2})$.
        \item\label{Item::Exmp::SharpddzRed::Chart::Span}  $F^*\Coorvec z\in\Co^{\alpha,\beta}(U_1,V_1;\C^3)$ and $F^*\Coorvec z$ spans $\Se|_{U_1\times V_1}$.
    \end{enumerate}
    \item\label{Item::Exmp::SharpddzRed::PDE} There is a function $f\in\Co^{\alpha,\beta}(U_2, V_2;\C)$ defined in a neighborhood $U_2\times V_2\subseteq\C^1_w\times\R^1_\theta$ of $(0,0)$ that solves \eqref{Eqn::Exmp::SharpddzRed::aPDE}.
\end{enumerate}
\end{prop}
The results are also true for $\alpha$ or $\beta=\infty$, without any change of the proof.

For the application to Theorem \ref{MainThm::Sharpddz}, we only need the contraposition of \ref{Item::Exmp::SharpddzRed::Chart} $\Rightarrow$ \ref{Item::Exmp::SharpddzRed::PDE}.

\begin{proof}
\ref{Item::Exmp::SharpddzRed::Chart} $\Rightarrow$ \ref{Item::Exmp::SharpddzRed::PDE}: Write $F=:(F'',F''')$ and $\Phi=(\Phi'',\Phi'''):=F^\Inv:F(U_1\times V_1)\subseteq\C^1_z\times\R^1_s\to \C^1_w\times\R^1_\theta$. By assumption $\nabla_z\Phi=(\nabla_z\Phi'',\nabla_z\Phi''')\in C^0(F(U_1\times V_1);\C^{2\times 2})$, so $F^*\Coorvec z=\frac{\partial(\re\Phi'',\im\Phi'',\Phi''')}{\partial z}\circ F$ and $F^*\Coorvec {\bar z}=\frac{\partial(\re\Phi'',\im\Phi'',\Phi''')}{\partial\bar z}\circ F$ are well-defined continuous vector fields on $\R^3_{(\re w,\im w,\theta)}$.

By assumption $(\Se+\bar\Se)|_{U_1\times V_1}=\Span(\partial_w,\partial_{\bar w})$ is spanned by $F^*\Coorvec{z}$ and $F^*\Coorvec z$, so $\partial_wF'''=\partial_{\bar w}F'''=0$ in the sense of distributions, which means $F'''(w,\theta)=F'''(\theta)$.

By assumption $d\bar F''=F^*d\bar z$ is continuous in $w$ and is annihilated by $\Span F^*\Coorvec z=\Span(\partial_w+a\partial_{\bar w})$. Therefore
\begin{equation}\label{Eqn::Exmp::SharpddzRed::Exp1}
    (\partial_w+a\partial_{\bar w})\bar F''=\bar F''_w+a\bar F''_{\bar w}=0,\quad\text{i.e.}\quad F''_{\bar w}=-\bar aF''_w.
\end{equation}

So for the transition matrix (cf. \eqref{Eqn::Exmp::Sharpdds::ClassicalExp1}) we have
\begin{equation}\label{Eqn::Exmp::SharpddzRed::Exp2}
\begin{pmatrix}F^*dz\\F^*d\bar z\\F^*ds\end{pmatrix}=
\begin{pmatrix}
F''_w&-\bar aF''_w&F''_\theta\\
-a\bar F''_{\bar w}&\bar F''_{\bar w}&\bar F''_\theta\\
&&F'''_\theta
\end{pmatrix}
\begin{pmatrix}dw\\d\bar w\\d\theta\end{pmatrix}\quad \Rightarrow\quad
\begin{pmatrix}
\partial_w\\ \partial_{\bar w}
\end{pmatrix}=
\begin{pmatrix}
F''_w&-a\bar F''_{\bar w}\\
-\bar aF''_w&\bar F''_{\bar w}
\end{pmatrix}
\begin{pmatrix}F^*\partial_z\\F^*\partial_{\bar z}\end{pmatrix}.
\end{equation}

Applying the Cramer's rule of inverse matrices to the right hand matrix in \eqref{Eqn::Exmp::SharpddzRed::Exp2}, we have 
\begin{equation}\label{Eqn::Exmp::SharpddzRed::EqnF*ddz}
F^*\Coorvec z=\frac{\bar F''_{\bar w}}{(1-|a|^2)|F''_w|^2}(\partial_w+a\partial_{\bar w})=\frac{\partial_w+a\partial_{\bar w}}{(1-|a|^2)F''_w}\quad\text{in}\quad U_1\times V_1.
\end{equation}
Since $F^*\Coorvec z\in\Co^{\alpha,\beta}(U_1, V_1;\C^3)$ and $a\in\Co^{\alpha,\beta}(U_1,V_1;\C)$, we see that $F''_w\in\Co^{\alpha,\beta}_\loc(U_1,V_1;\C)$ and $F''_w\not\equiv0$ in $U_1$.

Take $U_2\times V_2\subseteq U_1\times V_1$ as a neighborhood of $(0,0)$ such that $\sup_{(w,\theta)\in U_2\times V_2}|\frac{F''_w(w,\theta)}{F''_w(0,\theta)}-1|<1$. 
We define 
\begin{equation}\label{Eqn::Exmp::SharpddzRed::DefEqnf}
    f(w,\theta):=\log\frac{\bar F''_{\bar w}(w,\theta)}{\bar F''_{\bar w}(0,\theta)}=\sum_{n=1}^\infty\frac{(-1)^{n-1}}n\Big(\frac{\bar F''_{\bar w}(w,\theta)}{\bar F''_{\bar w}(0,\theta)}-1\Big)^n.
\end{equation}

Since $\big[(w,\theta)\mapsto \frac{F''_w(w,\theta)}{F''_w(0,\theta)}\big]\in\Co^{\alpha,\beta}(U_2,V_2;\C)$, we see that $f\in\Co^{\alpha,\beta}(U_2,V_2;\C)$ as well. So \eqref{Eqn::Exmp::SharpddzRed::aPDE} follows:
\begin{equation}\label{Eqn::Exmp::SharpddzRed::PfTmp1}
    f_w+af_{\bar w}=\frac{\bar F''_{w\bar w}}{\bar F''_{\bar w}}+\frac{a\bar F''_{\bar w\bar w}}{\bar F''_{\bar w}}=\frac{(\bar F''_w+a\bar F''_{\bar w})_{\bar w}}{\bar F''_{\bar w}}-\frac{a_{\bar w}\bar F''_{\bar w}}{\bar F''_{\bar w}}=-a_{\bar w}.
\end{equation}


\noindent\ref{Item::Exmp::SharpddzRed::PDE} $\Rightarrow$ \ref{Item::Exmp::SharpddzRed::Chart}: This is the reverse of the above argument. By shrinking $U_2\times V_2$ we assume $U_2$ and $V_2$ are both convex neighborhoods.

Suppose $f\in\Co^{\alpha,\beta}(U_2,V_2;\C)$ is a solution to \eqref{Eqn::Exmp::SharpddzRed::aPDE}. Take $\mu:=e^{f}$, we know $\mu\in\Co^{\alpha,\beta}(U_2,V_2;\C)$ satisfies
\begin{equation}\label{Eqn::Exmp::SharpddzRed::PfTmp2}
    (a\mu)_{\bar w}=a_{\bar w}e^f+af_{\bar w}e^f=-f_we^f=-\mu_w,\quad\text{in}\quad U_2\times V_2.
\end{equation}
Therefore $\mu(\cdot,\theta)(d\bar w-a(\cdot,\theta)dw)$ is a closed 1-form in $U_2\subseteq\C^1$, for each $\theta\in\R$. We define $F'':U_2\times V_2\to\C$ as 
\begin{equation}\label{Eqn::Exmp::SharpddzRed::DefF''}
    F''(w,\theta):=\int_0^1(w\cdot\bar \mu(tw,\theta)+\bar w\cdot\bar a(tw,\theta)\bar\mu(tw,\theta))dt.
\end{equation}

We see that $F''\in\Co^{\alpha,\beta}_\loc(U_1,V_1;\C)$ and satisfies $d\bar F''=\mu d\bar w-a\mu dw+\bar F''_\theta d\theta$. In other words
\begin{equation*}
    \nabla_w(F'',\bar F'')=\begin{pmatrix}F''_w&\bar F''_w
    \\ F''_{\bar w}&\bar F''_{\bar w}
    \end{pmatrix}=\begin{pmatrix}\bar \mu&-a\mu
    \\ -\bar a\bar\mu&\mu
    \end{pmatrix}\quad\Rightarrow\quad \det\big(\nabla_w(F'',\bar F'')\big)=|\mu|^2(1-|a|^2)>0.
\end{equation*}

Since $a,\mu\in \Co^{\alpha,\beta}$, we see that $F''\in\Co^{\alpha+1,\beta}$ and the matrix function $\nabla_w(F'',\bar F'')\in\Co^{\alpha,\beta}$ is invertible at $(0,0)$. Applying Lemma \ref{Lem::Hold::CompofMixHold} \ref{Item::Hold::CompofMixHold::InvFun} to $(\re F'',\im F'')\in\Co^{\alpha+1,\beta}(U_2,V_2;\R^2)$, we can find a neighborhood $U_1\times V_1\subseteq U_2\times V_2$ such that $F''(\cdot,\theta):U_1\to\C^1_z$ is injective for each $\theta\in V_1$ and its inverse map $\Phi''(z,s)=F''(\cdot,s)^\Inv(z)$ is $\Co^{\alpha+1,\beta}_{z,s}$ near $(F''(0,0),0)=(0,0)$.

Take $F(w,\theta):=(F''(w,\theta),\theta)$ and $\Phi:=F^\Inv$, we see that $\Phi(z,s)=(\Phi''(z,s),s)$. Therefore $F:U_1\times V_1\to\C^1\times\R^1$ is homeomorphic to its image (since $F,\Phi$ are both $\Co^{\alpha+1,\beta}\subset C^0$), and $\nabla_z(F^\Inv)=(\nabla_z\Phi'',0)$ is a continuous map on $F(U_1\times V_1)$. This completes the proof of \ref{Item::Exmp::SharpddzRed::Chart::Reg}.

Now for each $\theta\in V_1$, $F''(\cdot,\theta):U_1\to\C_z$ is a $\Co^{\alpha+1}$-coordinate chart satisfying $d(\bar F''(\cdot,\theta))=\mu(\cdot,\theta)d\bar w-a\mu(\cdot,\theta)dw$ in $U_1$, so $(\partial_w+a\partial_{\bar w})F''=0$, which means $F''(\cdot,\theta)^*d\bar z=d\bar F''(\cdot,\theta)$ spans $\Se|_{U_1\times\{\theta\}}=\Span(\partial_w+a(\cdot,\theta)\partial_{\bar w})$. Taking a dual we get $\Se|_{U_1\times\{\theta\}}=\Span F''(\cdot,\theta)^*\Coorvec z$ for each $\theta\in V_1$.

Now $F(w,\theta)=(F''(w,\theta),\theta)$ and $\Phi(z,s)=(\Phi''(z,s),s)$, so $((\nabla_z\Phi)\circ F)(w,\theta)=((\nabla_z\Phi''(\cdot,\theta))\circ F''(\cdot,\theta))(w)$ hence $F^*\Coorvec z|_{(w,\theta)}=F''(\cdot,\theta)^*\Coorvec z|_w$ spans $\Se_{(w,\theta)}$ for each $w\in U_1$ and $\theta\in V_2$, finishing the proof of \ref{Item::Exmp::SharpddzRed::Chart::Span}.
\end{proof}
\begin{remark}We do not know how to relax the assumption $\alpha>\frac12$ to $\alpha>0$. The difficulty comes from \eqref{Eqn::Exmp::SharpddzRed::PfTmp1} and \eqref{Eqn::Exmp::SharpddzRed::PfTmp2} in the proof, where the intermediate calculations need a product of a $\Co^{\alpha-1}$-function and a $\Co^\alpha$-function. On the other hand, the PDE \eqref{Eqn::Exmp::SharpddzRed::aPDE} can still be defined when $0<\alpha\le\frac12$, by requiring $a$ has some extra Sobolev regularity, for example we can assume $a\in\Co^{\alpha,\beta}(\C^1_w,\R^1_\theta;\C)\cap L^\infty(\R^1_\theta;W^{1,2}(\C^1_w;\C))$ (see Proposition \ref{Prop::Exmp::Expa} \ref{Item::Exmp::Expa::Reg1}).
\end{remark}
For Step 2, the construction of $a(w,\theta)$ requires the Sobolev-valued H\"older space $\Co^\beta_s(V;W^{k,p}_x(U))$. Recall its definition in Definition \ref{Defn::Intro::DefofHold}. 




\begin{lem}\label{Lem::Exmp::HolSobLem}
Let $\beta>0$ and $1<p<\infty$. Let $U\subset\R^n$ be a bounded smooth domain.
\begin{enumerate}[parsep=-0.3ex,label=(\roman*)]
    \item\label{Item::Exmp::HolSobLem::Prod1} The product map $[(f,g)\mapsto fg]:\Co^\beta(\R^q;L^p(U))\times L^\infty\Co^\beta(U,\R^q)\to\Co^\beta(\R^q;L^p(U))$ is bounded bilinear.
    \item\label{Item::Exmp::HolSobLem::Prod2} When $\beta<1$, then the product map $[(f,g)\mapsto fg]:(L^\infty\cap W^{1,2})(U)\times\Co^{\beta-1}(U)\to\Co^{\beta-1}(U)$ is bounded bilinear.
\end{enumerate}
\end{lem}
\begin{proof}
    \ref{Item::Exmp::HolSobLem::Prod1} can follow from the boundedness of the product map $L^p(U)\times L^\infty(U)\to L^p(U)$ and that $\Co^\beta(\R^q)$ is closed under multiplications. Alternatively one can apply Lemma \ref{Lem::Hold::VectPara} by taking $\Pi:L^p(U)\times L^\infty(U)\to L^p(U)$ be the natural product map. We leave the details to the readers.

    
    The proof of \ref{Item::Exmp::HolSobLem::Prod2} follows from \cite[Chapter 2.8.1]{BahouriCheminDanchin}, where in the reference we use Bony's decomposition $uv=T(u,v)+T(v,u)+R(u,v)$ for distributions $u,v$ on $\R^2$ (see \cite[(2.41)]{BahouriCheminDanchin} and Remark \ref{Rmk::Hold::LemParaProd}). By \cite[Theorem 2.82]{BahouriCheminDanchin} (or Lemma \ref{Lem::Hold::ParaBdd}) we have $T:L^\infty(\R^2)\times\Co^{\beta-1}(\R^2)\to\Co^{\beta-1}(\R^2)$ and $T:\Co^{\beta-1}(\R^2)\times L^\infty(\R^2)\to\Co^{\beta-1}(\R^2)$. By \cite[Theorem 2.85]{BahouriCheminDanchin} we have $R:W^{1,2}(\R^2)\times\Co^{\beta-1}(\R^2)\to\Bs_{22}^\beta(\R^2)$. Here $\Bs_{22}^\beta(\R^2)=H^\beta(\R^2)$ is the fractional $L^2$-Sobolev space (see for example \cite[Definition 2.68]{BahouriCheminDanchin} or \cite[Definition 2.3.1/2(i)]{Triebel1}).
    We have Sobolev embedding $\Bs_{22}^{\beta}(\R^2)\hookrightarrow\Co^{\beta-1}(\R^2)$, whose proof can be found in for example \cite[Chapter 2.7.1]{Triebel1}\footnote{We recall the correspondence of Besov space $\Co^\beta(\R^n)=\Bs_{\infty\infty}^\beta(\R^n)$ and Triebel-Lizorkin space $L^p(\R^n)=\mathscr F_{p2}^0(\R^n)$ from \cite[Theorems 2.5.7 and 2.5.12]{Triebel1}.}. Thus $R:W^{1,2}(\R^2)\times\Co^{\beta-1}(\R^2)\to\Co^{\beta-1}(\R^2)$ is bounded.
    
    Therefore $[(u,v)\mapsto uv]:(L^\infty\cap W^{1,2})(\R^2)\times\Co^{\beta-1}(\R^2)\to\Co^{\beta-1}(\R^2)$ is bilinearly bounded. Taking restrictions to $U$ we complete the proof.
\end{proof}

\begin{lem}\label{Lem::Exmp::InvDBar}
Let $\Pc_0$ and $\Pc$ be the zero Dirichlet boundary solution operators in Lemmas \ref{Lem::Hold::DiriSol} and \ref{Lem::Hold::LapInvBdd} with $U=\B^2_w\subset\C^1$ and $V=\R^1_\theta$ in the lemma. Define $\Tc_0:=4\partial_{\bar w}\Pc_0$ and $\Tc:=4\partial_{\bar w}\Pc$. Then 
\begin{enumerate}[parsep=-0.3ex,label=(\roman*)]
    \item\label{Item::Exmp::InvDBar::Bdd} $\partial_w\Tc f=f$ for all $f\in L^\infty(\B^2,\R^1;\C)$ and $\Tc$ has the following boundedness: for $\alpha>0$, $1<p<\infty$ and $\Xs\in\{L^\infty,\Co^\gamma:\gamma>0\}$,
\begin{gather}\label{Eqn::Exmp::InvDBar::Bdd}
    \Tc:\Co^{\alpha-1}\Xs(\B^2,\R^1;\C)\to\Co^\alpha\Xs(\B^2,\R^1;\C)\text{ and } \Tc:\Co^\alpha(\R^1;L^p(\B^2))\to\Co^\alpha(\R^1;W^{1,p}(\B^2)).
\end{gather}
    \item\label{Item::Exmp::InvDBar::Holo} Moreover, if $f\in L^\infty\Co^\beta(\B^2,\R^1;\C)$ satisfies $f|_{\{|w|>\frac12\}}\equiv0$, then
\begin{equation}\label{Eqn::Exmp::InvDBar::Holo}
    f-\Tc\partial_w f\in\Co^\infty\Co^\beta(\B^2,\R^1;\C).
\end{equation}
\end{enumerate}

\end{lem}
\begin{proof}
That $\partial_w\Tc=\id$ and $\Tc:\Co^{\alpha-1}_w\Xs_\theta\to \Co^{\alpha}_w\Xs_\theta$ are the immediate consequences to Lemma \ref{Lem::Hold::LapInvBdd} since $\Delta_{\re w,\im w}=4\frac{\partial^2}{\partial w\partial \bar w}$.

We have boundedness $\Pc_0:L^p(\B^2)\to W^{2,p}(\B^2)$ by \cite[Theorem 9.13]{GilbargTrudinger}, thus $\Tc_0:L^p(\B^2;\C)\to W^{1,p}(\B^2;\C)$. We get $\Tc:\Co^\beta(L^p)\to\Co^\beta(W^{1,p})$ by Lemma \ref{Lem::Hold::OperatorExtension}. Thus we prove \ref{Item::Exmp::InvDBar::Bdd}.

\smallskip
To prove \ref{Item::Exmp::InvDBar::Holo}, we need to show that $\id-\Tc\partial_w:L^\infty(\frac12\B^2;\C)\to\Co^M(\B^2;\C)$ is bounded for all $M>0$. 

Let $g\in L^\infty(\frac12\B^2;\C)$ with zero extension outside $\frac12\B^2$. By formula of Green functions (see for example \cite[Theorem 2.12 and Section 2.2.c]{Evans}) we have for $w\in\B^2$,
\begin{equation*}
    g(w)-\Tc_0\partial_wg(w)=g(w)-\int_{\frac12\B^2}4g(z)\partial^2_{z\bar w}\big(\Ga(z-w)-\Ga\big(|z|w-\tfrac z{|z|}\big)\big)dxdy=\int_{\frac12\B^2}4g(z)\partial^2_{z\bar w}\Ga\big(|z|w-\tfrac z{|z|}\big)dxdy.
\end{equation*}
Here $z=x+iy\in\B^2$ and $\Ga(z)=\frac1{2\pi}\log|z|$ (which is the $n=2$ case in \eqref{Eqn::SecHolLap::NewtonPotent}). The first equality follows from integration by part, the second equality follows from the fact that $\partial_{z\bar w}\Ga(z-w)=\partial_{z\bar z}\Ga(z-w)=\frac14\Delta_{x,y}\Ga(z-w)=\frac{\delta_w(z)}4$.

Clearly $(z,w)\mapsto\Ga(|z|w-z/|z|)$ is a smooth function on $(z,w)\in r\B^2\times\B^2$ for any $r<1$. Thus for any integer $M\ge0$
\begin{equation*}
    \|\nabla^M(g-\Tc_0\partial_wg)\|_{L^\infty(\B^2;\C^{2^k})}\le\|g\|_{L^\infty(\frac12\B^2;\C)}\sup_{|w|<\frac12}\int_{\B^2}\big|\nabla^{M+2}_{z,w}\Ga\big(|z|w-\tfrac z{|z|}\big)\big|dxdy\lesssim_M\|g\|_{L^\infty(\frac12\B^2;\C)}.
\end{equation*}
We conclude that $\id-\Tc_0\partial_w:L^\infty(\frac12\B^2;\C)\to\Co^M(\B^2;\C)$ is bounded for all $M>0$. By Lemma \ref{Lem::Hold::OperatorExtension}, we get the boundedness $\id-\Tc\partial_w:L^\infty\Co^\beta(\frac12\B^2,\R^1;\C)\to\Co^M\Co^\beta(\B^2,\R^1;\C)$. The proof is complete.
\end{proof}

\begin{lem}\label{Lem::Exmp::CalpNonlocal}
Let $\beta>0$ and let $(c_k)_{k=1}^\infty\subset\C$ be a sequence.
Define $g(\theta):=\sum_{k=1}^\infty c_ke^{2\pi i\cdot 2^k\theta}2^{-\beta\theta}$, then
\begin{enumerate}[parsep=-0.3ex,label=(\roman*)]
    \item\label{Item::Exmp::CalpNonlocal::Bdd} $\|g\|_{\Co^\beta(\R;\C)}\approx\sup_{k\ge0}|c_k|$ with the implied constant depending only on $\beta$.
    \item\label{Item::Exmp::CalpNonlocal::Main} If $\varlimsup_{k\to\infty}|c_k|=\infty$, then $g\notin\Co^\beta(I;\C)$ for every interval $I\subseteq\R$.
\end{enumerate}

\end{lem}
\begin{proof}
Let $\{\psi_k\}_{k=0}^\infty$ be a dyadic resolution for $\R$. Recall from \eqref{Eqn::Hold::RmkDyaSupp} $\supp\hat\psi_k\subset(-2^{k-1},-2^{k+1})\cup(2^{k-1},2^{k+1})$ holds for $k\ge1$. Since the Fourier support of $e^{2\pi i2^k\theta}$ is $\{2^k\}$, we have $\psi_k\ast g(\theta)=c_k2^{-\beta\theta}e^{2\pi i\cdot 2^k\theta}$ for all $k\ge1$. Therefore by Lemma \ref{Lem::Hold::HoldChar} \ref{Item::Hold::HoldChar::LPHoldChar}, we get the result \ref{Item::Exmp::CalpNonlocal::Bdd}:
$$\textstyle\|g\|_{\Co^\beta(\R)}\approx_{\psi,\beta}\sup_{k\ge0}2^{k\beta}\|\psi_k\ast g\|_{L^\infty}=\sup_{k\ge1}|c_k|.$$

Now suppose $\varlimsup_{k\to\infty}|c_k|=+\infty$, so $\|g\|_{\Co^\beta(\R)}\approx_{\psi,\beta}=\sup|c_k|=\infty$.

For a given interval $I\subseteq\R$, $I$ contains a dyadic subinterval, say $[ n2^{-k_0},(n+1)2^{-k_0}]$ where $n\in\Z$ and $k_0\in\Z_+$. Take $g^{k_0}(\theta):=\sum_{k=k_0+1}^\infty c_{k}2^{2\pi i2^k\theta}2^{-k\beta}$, so $$\textstyle g(\theta)-g^{k_0}(\theta)=\sum_{k=1}^{k_0}c_ke^{2\pi i2^k\theta}2^{-\beta\theta},$$
must have finite $\Co^\beta(\R)$-norm (though whose bound depends on $k_0$).

Since $g^{k_0}$ is $2^{-k_0}$-periodic, $\|g\|_{\Co^\beta(\R)}=\infty$ and $\|g-g^{k_0}\|_{\Co^\beta(\R)}<\infty$, we know $\|g^{k_0}\|_{\Co^\beta(n2^{-k_0},(n+1)2^{-k_0})}=\infty$. Therefore, we complete the proof of \ref{Item::Exmp::CalpNonlocal::Main} by:
$$\textstyle\|g\|_{\Co^\beta(I)}\ge\|g\|_{\Co^\beta([n2^{-k_0},(n+1)2^{-k_0}])}\ge \|g^{k_0}\|_{\Co^\beta([n2^{-k_0},(n+1)2^{-k_0}])}-\|g-g^{k_0}\|_{\Co^\beta(\R)}=\infty.\qedhere$$
\end{proof}

We now construct our $\Co^{\alpha,\beta}$-function $a(w,\theta)$ and finish the Step 2:
\begin{prop}\label{Prop::Exmp::Expa}
Let $\alpha,\beta>0$. There is a function $ A(w,\theta)$ defined in $\C^1_w\times\R^1_\theta$, that satisfies the following properties:
\begin{enumerate}[parsep=-0.3ex,label=(\roman*)]
    \item\label{Item::Exmp::Expa::Reg1} $ A\in\Co^{\alpha,\beta}(\C^1,\R^1;\C)\cap \Co^\beta(\R^1;W^{1,2}(\C^1;\C))$ and $[(w,\theta)\mapsto w\partial_{\bar w} A(w,\theta)]\in\Co^{\alpha,\beta}(\C^1,\R^1;\C)$. 
    \item\label{Item::Exmp::Expa::Reg2} $[(w,\theta)\mapsto w^{-1} A(w,\theta)]\in \Co^\frac\beta2(\R^1;L^p(\C^1;\C))$ for all  $1<p<\frac4{2-\min(\alpha,2)}$. In particular $w^{-1}A\in\Co^\frac\beta2(\R^1;L^\frac{\alpha+4}2(\C^1;\C))$.
    
    \item\label{Item::Exmp::Expa::Unbdd} $\Tc[\partial_{\bar w} A]\notin L^\infty_w\Co^\beta_\theta$ near $(w,\theta)=(0,0)$, where $\Tc$ be the inverse $\partial_w$-operator in Lemma \ref{Lem::Exmp::InvDBar}.

\end{enumerate}
\end{prop}

The construction is essentially the parameter version of Lemma \ref{Lem::Exmp::CountNN::lem}. 
\begin{remark}
By taking some simple modifications from the construction below (see \eqref{Eqn::Exmp::Expa::Defofa}), one can also prove Proposition \ref{Prop::Exmp::Expa} for the case $\alpha=\infty$, $\beta<\infty$. We leave the details to readers.
\end{remark}

\begin{proof}[Proof of Proposition \ref{Prop::Exmp::Expa}]
First we find a function $b\in C_c^0(\B^2;\C)$ such that,
\begin{enumerate}[nolistsep,label=(b.\arabic*)]
    \item\label{Item::Exmp::Expa::BReg} $b\in  C^\infty_\loc(\B^2\backslash\{0\};\C)$, $b\in W^{1,2}(\B^2;\C)$ and $wb\in C_c^1(\B^2;\C)$. 
    \item\label{Item::Exmp::Expa::BUnbdd} $\Tc_0\partial_{\bar w}b\notin L^\infty$ near $w=0$. Here $\Tc_0$ is the inverse $\partial_w$-operator in Lemma \ref{Lem::Exmp::InvDBar}.
\end{enumerate}

Indeed, we take a function $b$ such that
\begin{equation}\label{Eqn::Exmp::Expa::Defofb}
    b(w):=\partial_w\big(\bar w(-\log|w|)^{\frac13}\big)=-\tfrac16(-\log|w|)^{-\frac23}\tfrac{|w|^2}{w^2},\quad\text{when }|w|<\tfrac12.
\end{equation}
And $b$ is smoothly cutoffed outside $|w|<\frac12$.

Immediately we get $b\in C^\infty_\loc(\B^2\backslash\{0\};\C)$ and $wb(w)\in C^1_c$. And clearly we have $|\nabla b(w)|\lesssim|w|^{-1}(-\log|w|)^{-\frac23}$, so $b\in W^{1,2}$, which gives \ref{Item::Exmp::Expa::BReg}.

To prove \ref{Item::Exmp::Expa::BUnbdd}, since $\partial_w\Tc_0 \partial_{\bar w}b=\partial_{\bar w}b$ it suffices to show that if $U$ is a neighborhood of $0$ and $f:U\subseteq\frac12\B^2\to\C$ is a function satisfies $\partial_wf=\partial_{\bar w}b$, then $f\notin L^\infty$ near $0$.

Indeed by \eqref{Eqn::Exmp::Expa::Defofb} in a neighborhood of $0$ we have $f_w=(\bar w(-\log|w|)^\frac12)_{w\bar w}$, so $f-(\bar w(-\log|w|)^\frac12)_{\bar w}$ is annihilated by $\partial_w$ which is an anti-holmorphic function. Therefore $f-(\bar w(-\log|w|)^\frac12)_{\bar w}$ is a bounded function near $w=0$. On the other hand 
$$(\bar w(-\log|w|)^\frac13)_{\bar w}=(-\log|w|)^\frac13+O(1)\quad \text{as }w\to0,$$
is not bounded near $0$. We conclude that $f\notin L^\infty$ near $0$.

\medskip
Let $\rho_0\in\Sc(\C^1_w)$ be a Schwartz function whose Fourier transform satisfies $\hat\rho_0\in C_c^\infty(\B^2)$ and $\hat\rho_0|_{\frac12\B^2}\equiv1$. Define functions $(\rho_k,c_k)_{k=0}^\infty$ and $ A(w,\theta)$ as 
\begin{equation}\label{Eqn::Exmp::Expa::Defofa}
    \rho_k(w):=2^{n\frac\beta\alpha k}\rho_0(2^{\frac\beta\alpha k}w),\quad c_k(w):=\rho_k\ast b(w)-\rho_k\ast b(0),\quad k\ge0;\quad  A(w,\theta):=\sum_{j=0}^\infty c_j(w)e^{2\pi i2^j\theta}2^{-j\beta}.
\end{equation}

Thus, we have $wc_k=\int_\C(w-z)\rho_k(w-z)b(z)d\Vol(z)=(w\rho_k)\ast b$ for $k\ge0$.
By scaling
\begin{equation}\label{Eqn::Exmp::Expa::EstRho}
    \|\rho_k\|_{L^1}=\|\rho_0\|_{L^1},\quad\|\nabla\rho_k\|_{L^1}=2^{\frac\beta\alpha k}\|\nabla\rho_0\|_{L^1},\quad\|w\nabla\rho_k\|_{L^1}=\|w\nabla\rho_0\|_{L^1}\quad\text{for each }k\ge1.
\end{equation}
Therefore,
\begin{equation}\label{Eqn::Exmp::Expa::EstC}
\begin{gathered}
    \|c_k\|_{L^\infty\cap W^{1,2}}\le2\|\rho_k\ast b\|_{L^\infty\cap W^{1,2}}\lesssim\|\rho_k\|_{L^1}\|b\|_{L^\infty\cap W^{1,2}}\lesssim1;\\
    \|w\partial_{\bar w} c_k\|_{L^\infty}=\|\partial_{\bar w}(w\rho_k)\ast b\|_{L^\infty}\le\|w\partial_{\bar w}\rho_k\|_{L^1}\|b\|_{L^\infty}\lesssim1.
\end{gathered}
\end{equation}

Let $(\phi_j)_{j=0}^\infty$ and $(\psi_k)_{k=0}^\infty$ be two dyadic resolutions for $\C^1_w$ and $\R^1_\theta$ respectively (recall Definition \ref{Defn::Hold::DyadicResolution}). Recall from \eqref{Eqn::Hold::RmkDyaSupp} $\hat\phi_0\in C_c^\infty(2\B^2)$, $\hat\psi_0\in C_c^\infty(2\B^1)$, $\supp\hat\phi_j\subset2^{j+1}\B^2\backslash2^{j-1}\B^2$ and $\supp\hat\psi_j\subset2^{j+1}\B^1\backslash2^{j-1}\B^1$ for $j\ge1$. Since $\supp\hat\rho_0\subset\B^2$, by scaling we have $\supp\hat\rho_k\subset2^{\frac\beta\alpha k}\B^2$. Thus,
\begin{equation}\label{Eqn::Exmp::Expa::Estconv}
    \phi_j\ast_w A(w,\theta)=\sum_{k\ge\max(0,\frac\alpha\beta(j-1))}(\phi_j\ast c_k)(w)e^{2\pi i2^k\theta}2^{-k\beta},\quad \psi_j\ast_\theta A(w,\theta)=c_j(w)e^{2\pi i2^j\theta}2^{-j\beta},\quad j\ge0.
\end{equation}
By \eqref{Eqn::Exmp::Expa::EstC} and \eqref{Eqn::Exmp::Expa::Estconv} we get $ A\in\Co^\alpha L^\infty\cap L^\infty\Co^\beta(\C^1,\R^1;\C)\cap \Co^\beta(\R^1;W^{1,2}(\C^1;\C))$, since for $j\ge0$,
\begin{equation*}
\begin{gathered}
    \|\phi_j\ast_w A\|_{L^\infty(\C^1\times\R^1;\C)}\le\sum_{k\ge\max(0,\frac\alpha\beta(j-1))}\|\phi_j\|_{L^1(\C^1)}\|c_k\|_{L^\infty(\C^1\times\R^1;\C)}2^{-k\beta}\lesssim_{\phi,c,\alpha,\beta}2^{-j\alpha};\\
    \sup_{\theta\in\R^1}\|\psi_j\ast_\theta A(\cdot,\theta)\|_{L^\infty\cap W^{1,2}(\C^1;\C)}\le\|\psi_j\|_{L^1(\R)}\|c_j\|_{L^\infty\cap W^{1,2}(\C^1;\C)}2^{-j\beta}\lesssim_{\psi,c,\beta}2^{-j\beta}.
\end{gathered}
\end{equation*}
Taking supremum over $j\ge0$ we get \ref{Item::Exmp::Expa::Reg1}.

\medskip Using \eqref{Eqn::Exmp::Expa::EstRho} and \eqref{Eqn::Exmp::Expa::EstC} we also have $\|c_k\|_{C^1}\lesssim 2^{\frac\beta\alpha k}$ and $\|w\nabla c_k\|_{C^0}\lesssim1$. Thus 
\begin{equation*}
    |\nabla c_k(w)|\lesssim\min(|w|^{-1},2^{\frac\beta\alpha k})\le|w|^{\min(\frac\alpha2,1)-1}2^{k\frac\beta2\cdot\min(\frac\alpha2,1)}\le|w|^{\min(\frac\alpha2,1)-1}2^{k\frac\beta2} ,\quad\text{for all}\quad k\ge0,\quad w\in\C^1.
\end{equation*}
By \eqref{Eqn::Exmp::Expa::Defofa} we have $c_k(0,\theta)\equiv0 $ for all $\theta\in\R$. Thus by \eqref{Eqn::Exmp::Expa::Estconv},
\begin{equation*}
    |\psi_k\ast_\theta (w^{-1}A)(w,\theta)|=|w^{-1}c_k(w)|2^{-k\beta}\le|\nabla c_k(w)|2^{-k\beta}\le|w|^{\min(\frac\alpha2,1)-1}2^{-k\frac\beta2},\quad k\ge0,\quad w\in\C^2,\quad\theta\in\R^1.
\end{equation*}

On the other hand $|w|^{-\gamma}\in L^{\frac2\gamma-}(\B^2)$ holds for all $0<\gamma<2$. Thus for $1<p<\frac2{1-\min(\alpha/2,1)}=\frac4{2-\min(\alpha,2)}$ we have $\sup_{k\ge0;\theta\in\R^1}2^{k\frac\beta2}\|\psi_k\ast_\theta (w^{-1}A)(\cdot,\theta)\|_{L^p(\C^1;\C)}<\infty$.

Clearly $2<\frac{\alpha+4}2<\frac4{2-\min(\alpha,2)}$, we finish the proof of \ref{Item::Exmp::Expa::Reg2}.

\medskip To prove \ref{Item::Exmp::Expa::Unbdd}, since $\Tc$ only acts on $w$-variable,  by \eqref{Eqn::Exmp::Expa::Defofa} we have 
$$\Tc  A_{\bar w}(w,\theta)=\sum_{k=1}^\infty2^{-k\alpha}e^{2\pi i2^k\theta}\Tc_0\partial_{\bar w} c_k(w)=\sum_{k=1}^\infty2^{-k\alpha}e^{2\pi i2^k\theta}\Tc_0(\rho_k\ast \partial_{\bar w}b)(w),\quad w\in\C^1,\ \theta\in\R.$$

Now $\lim\limits_{k\to\infty}\rho_k\ast(\partial_{\bar w}b)=\partial_{\bar w} b$ holds (as $L^2$-functions), and by \ref{Item::Exmp::Expa::BUnbdd} $\Tc_0 \partial_{\bar w}b\notin L^\infty$ near $w=0$. So for any neighborhood $U\subseteq\C^1$ of $w=0$, we have $\varlimsup_{k\to\infty}\|\Tc_0 \partial_wc_k\|_{L^\infty(U;\C)}=+\infty$. Therefore by Lemma \ref{Lem::Exmp::CalpNonlocal} \ref{Item::Exmp::CalpNonlocal::Main}, for any interval $I\subseteq\R$ containing $\theta=0$, we have $\|\Tc\partial_{\bar w} A\|_{L^\infty\Co^\alpha(U,I;\C)}\gtrsim \varlimsup_{k\to\infty}\|\Tc_0\partial_{\bar w}b\|_{L^\infty(U;\C)}=+\infty$.

So $\Tc \partial_{\bar w} A\notin L^\infty_w\Co^\beta_\theta$ near $(0,0)$, finishing the proof of \ref{Item::Exmp::Expa::Unbdd}.
\end{proof}

Now we prove Step 3, which is the non-existence of $f\in\Co^{\alpha,\beta}$ to \eqref{Eqn::Exmp::SharpddzRed::aPDE}. Cf. Proposition \ref{Prop::Exmp::CountNN::MainProp}. Here we consider all $\alpha,\beta>0$.

\begin{thm}[$\bar\partial$-equation for sharpness of $F^*\Coorvec z\notin\Co^\alpha $]\label{Thm::Exmp::ProofExampleddz}
    Let $\alpha,\beta>0$, let $ A$ be as in Proposition \ref{Prop::Exmp::Expa}. Then there is a constant $\delta=\delta(\alpha,\beta, A)>0$, such that for every neighborhood $U_2\times V_2\subset\C^1_w\times\R^1_\theta$ of $(0,0)$ there is no $f\in\Co^{\alpha,\beta}(U_2,V_2;\C)$ solving \eqref{Eqn::Exmp::SharpddzRed::aPDE} with $a(w,\theta):=\delta\cdot A(w,\theta)$.
\end{thm}
This result contains the largest possible range for $\alpha$ and $\beta$. The proof can be simplified if $\alpha>1$.
\begin{proof}
Let $\Tc:\Co^{\alpha-1} L^\infty(\B^2,\R^1;\C)\to \Co^{\alpha} L^\infty(\B^2,\R^1;\C)$ be the inverse $\partial_w$-operator as in Lemma \ref{Lem::Exmp::InvDBar}. Recall that $\Tc:\Co^\beta(\R^1;L^p(\B^2;\C))\to \Co^\beta(\R^1;W^{1,p}(\B^2;\C))$ is also bounded for $1<p<\infty$. 

We define a linear operator $\Tf=\Tf_{ A}$ on functions in $\B^2\times\C^1$ as
\begin{equation}
    \Tf[\varphi]:=\Tc( A\cdot\partial_{\bar w}\varphi),\quad\varphi\in \Co^{\alpha} L^\infty(\B^2,\R^1;\C).
\end{equation}

By Lemma \ref{Lem::Hold::Product} \ref{Item::Hold::Product::Hold1}, Lemma  \ref{Lem::Exmp::HolSobLem} \ref{Item::Exmp::HolSobLem::Prod1} and \ref{Item::Exmp::HolSobLem::Prod2},
\begin{align*}
    \| A\partial_{\bar w}\varphi\|_{\Co^{\alpha-1} L^\infty(\B^2,\R^1;\C)}&\lesssim_\alpha\|A\|_{\Co^\alpha L^\infty}\|\partial_{\bar w}\varphi\|_{\Co^{\alpha-1} L^\infty}\lesssim_A\|\varphi\|_{\Co^{\alpha} L^\infty(\B^2,\R^1;\C)},&\text{provided }\alpha>\tfrac12;
    \\
    \| A\partial_{\bar w}\varphi\|_{\Co^{\alpha-1} L^\infty(\B^2,\R^1;\C)}&\lesssim_\alpha\|A\|_{L^\infty(\R^1;L^\infty\cap W^{1,2})}\|\partial_{\bar w}\varphi\|_{\Co^{\alpha-1} L^\infty}\lesssim_A\|\varphi\|_{\Co^{\alpha} L^\infty(\B^2,\R^1;\C)},&\text{provided }0<\alpha<1;
    \\
    \| A\partial_{\bar w}\varphi\|_{\Co^\beta(\R^1;L^p(\B^2;\C))}&\lesssim_{\beta}\| A\|_{L^\infty\Co^\beta}\|\partial_{\bar w}\varphi\|_{\Co^\beta(\R^1;L^p(\B^2;\C))}\lesssim_{ A,p}\|\varphi\|_{\Co^\beta(\R^1;W^{1,p}(\B^2;\C))},&\text{for all } 1<p<\infty.
\end{align*}

Taking $p=p_\alpha:=4+16/\alpha$, we conclude that $\Tf$ is a bounded endomorphism on both $\Co^\alpha L^\infty(\B^2,\R^1;\C)$ and $\Co^\alpha L^\infty(\B^2,\R^1;\C)\cap\Co^\beta(\R^1;W^{1,4+\frac{16}\alpha}(\B^2;\C))$.

We now take $\delta=\delta(\alpha,\beta, A)>0$ small such that
\begin{equation}\label{Eqn::Exmp::ProofExampleddz::DefDelta}
    \delta\big(\|\Tf\|_{\Co^{\alpha} L^\infty(\B^2,\R^1;\C)}+\|\Tf\|_{\Co^{\alpha} L^\infty(\B^2,\R^1;\C)\cap\Co^\beta(\R^1;W^{1,4+16/\alpha}(\B^2;\C))}\big)\le\tfrac12.
\end{equation}
Note that $\Tf$ is now contraction map in these two spaces.

We claim that $a(w,\theta)=\delta A(w,\theta)$ is as desired.

\medskip
Suppose by contradiction there is a $U_2\times V_2\ni (0,0)$ and a $f\in\Co^{\alpha,\beta}(U_2,V_2;\C)$ solving \eqref{Eqn::Exmp::SharpddzRed::aPDE}. Take $0<\eps_0<\frac12$ such that $\eps_0\B^2\subseteq U_2 $ and $(-\eps_0,\eps_0)\subseteq V_2$.  We define $\tilde\chi$ and $\chi$ such that
\begin{equation}\label{Eqn::Exmp::ProofExampleddz::Chi}
    \tilde\chi\in C_c^\infty(-\eps_0,\eps_0),\quad\tilde\chi|_{(-\frac12\eps_0,\frac12\eps_0)}\equiv1,\quad\chi\in C_c^\infty(\eps_0\B^2\times(-\eps_0,\eps_0)),\quad\chi(w,\theta):=\tilde\chi(|w|)\tilde\chi(\theta).
\end{equation}
Therefore $\chi\in C_c^\infty(U_2)$ equals to 1 in the neighborhood $\frac{\eps_0}2\B^2\times(-\frac{\eps_0}2,\frac{\eps_0}2)$ of $(0,0)$, and $\nabla_w\chi$ vanishes in the  neighborhood $\frac{\eps_0}2\B^2\times\R^1$ of $\{w=0\}=\{0\}_w\times\R^1_\theta$.

Define $g,h\in\Co^{\alpha,\beta}_c(U_2,V_2;\C)$ as $g:=\chi f$ and $h(w,\theta):=w g(w,\theta)$, so
\begin{equation}\label{Eqn::Exmp::ProofExampleddz::ghPDE}
    g_w+a\cdot g_{\bar w}=(\chi_w+a\chi_{\bar w})f-\chi a_{\bar w},\quad h_w+a\cdot h_{\bar w}=g+w(\chi_w+a\chi_{\bar w})f-\chi wa_{\bar w},\quad\text{in }\B^2_w\times\R^1_\theta.
\end{equation}

\noindent \textit{Claim 1:} $h\in\Co^\beta(\R^1;W^{1,4+\frac{16}\alpha}(\B^2;\C))$ and therefore $h_{\bar w}\in\Co^\beta(\R^1;L^{4+\frac{16}\alpha}(\B^2;\C))$.

Rewriting the second equation of \eqref{Eqn::Exmp::ProofExampleddz::ghPDE} we have
\begin{equation*}
    h=-\Tc(ah_{\bar w})+(h-\Tc h_w)+\Tc(g+w(\chi_w+a\chi_{\bar w})f-\chi wa_{\bar w}).
\end{equation*}

In other words $h=\varphi$ is a solution to the following affine linear equation:
\begin{equation}\label{Eqn::Exmp::ProofExampleddz::EqnVPhi}
    \varphi=-\delta\Tf[\varphi]+(h-\Tc h_w)+\Tc(g+w(\chi_w+a\chi_{\bar w})f-\chi wa_{\bar w}).
\end{equation}

By assumption $g\in\Co^{\alpha,\beta}(\B^2,\R^1;\C)$, $f\nabla_w\chi\in\Co^{\alpha,\beta} (\B^2,\R^1;\C^2) $. By Proposition \ref{Prop::Exmp::Expa} \ref{Item::Exmp::Expa::Reg1} $wa_{\bar w}\in \Co^{\alpha,\beta}(\B^2,\R^1;\C)$. Therefore by \eqref{Eqn::Exmp::InvDBar::Bdd} (using $L^\infty_w\Co^\beta_\theta\subset\Co^\beta_\theta(L^{4+\frac{16}\alpha}_w)$) we get 
\begin{equation}\label{Eqn::Exmp::ProofExampleddz::hTerm1}
    \Tc(g+w(\chi_w+a\chi_{\bar w})f-\chi wa_{\bar w})\in\Tc(\Co^{\alpha,\beta}(\B^2,\R^1;\C))\subset\Co^{\alpha}L^\infty(\B^2,\R^1;\C)\cap \Co^\beta(\R^1;W^{1,4+\frac{16}\alpha}(\B^2;\C)).
\end{equation}

By assumption $h\in\Co^{\alpha,\beta}(\B^2,\R^1;\C)\subset L^\infty\Co^\beta(\B^2,\R^1;\C)$ is supported in $\frac12\B^2_w\times\R^1_\theta$, so by \eqref{Eqn::Exmp::InvDBar::Holo} $h-\Tc h_w\in\Co^\infty\Co^\beta(\B^2,\R^1;\C)$. In particular 
\begin{equation}\label{Eqn::Exmp::ProofExampleddz::hTerm2}
    h-\Tc h_w\in\Co^{\alpha}L^\infty(\B^2,\R^1;\C)\cap \Co^\beta(\R^1;W^{1,4+\frac{16}\alpha}(\B^2;\C)).
\end{equation}
Using \eqref{Eqn::Exmp::ProofExampleddz::hTerm1}, \eqref{Eqn::Exmp::ProofExampleddz::hTerm2} and the assumption \eqref{Eqn::Exmp::ProofExampleddz::DefDelta}, $\big[\varphi\mapsto -\delta\Tf[\varphi]+(h-\Tc h_w)+\Tc(g+w(\chi_w+a\chi_{\bar w})f-\chi wa_{\bar w})\big]$ is a contraction map in both $\Co^{\alpha}L^\infty(\B^2,\R^1;\C)$ and $\Co^{\alpha} L^\infty(\B^2,\R^1;\C)\cap \Co^\beta(\R^1;W^{1,4+\frac{16}\alpha}(\B^2;\C))$. Since \eqref{Eqn::Exmp::ProofExampleddz::EqnVPhi} has a unique fixed point in both spaces, and  $h\in\Co^{\alpha} L^\infty(\B^2,\R^1;\C)$ is a fixed point to \eqref{Eqn::Exmp::ProofExampleddz::EqnVPhi}, we conclude that  $h\in\Co^\beta(\R^1;W^{1,4+\frac{16}\alpha}(\B^2;\C))$. Thus $h_{\bar w}\in\Co^\beta(\R^1;L^{4+\frac{16}\alpha}(\B^2;\C))$, obtaining the Claim 1.

\medskip
\noindent\textit{Claim 2:} $a\cdot g_{\bar w}\in\Co^{\eps-1}_w\Co^\beta_\theta(\B^2,\R^1;\C)$ for some $\eps>0$, and therefore $\Tc(a\cdot g_{\bar w})\in L^\infty\Co^{\beta}(\B^2,\R^1;\C)$.

Let $(\psi_j)_{j=0}^\infty$ be a dyadic resolution for $\R^1_\theta$. By Lemma \ref{Lem::Hold::LemParaProd}, for $l\ge0$, on $\C^1\times\R^1$,
\begin{align*}
    \psi_l\ast_\theta (ag_{\bar w})&=\psi_l\ast_\theta\Big(\sum_{j=l-2}^{l+2}\sum_{j'=0}^{j-3}+\sum_{j=l-2}^{l+2}\sum_{j'=0}^{j-3}+\sum_{\substack{j,j'\ge l-3;|j-j'|\le2}}\Big)(\psi_j\ast_\theta a)\cdot(\psi_{j'}\ast_\theta g_{\bar w})=:P_l^1+P_l^2+P_l^3.
\end{align*}

Our goal is to show that for $\eps=\frac\alpha{8+2\alpha}>0$, $\sup_{l\ge0;\theta\in\R^1}2^{l\beta}\|P_l^\nu(\cdot,\theta)\|_{\Co^{\eps-1}(\C^1;\C)}<\infty$ for $\nu=1,2,3$.

By assumption $a\in \Co^\beta(\R^1;L^\infty\cap W^{1,2}(\C^1;\C))$ and $g_{\bar w}\in\Co^{\alpha-1}L^\infty\cap\Co^{-1}\Co^\beta(\C^1,\R^1;\C)\subset\Co^{\frac\alpha2-1}\Co^\frac\beta2(\C^1,\R^1;\C)$, so we have $\sup_{\theta\in\R^1}\|\psi_j\ast_\theta a(\cdot,\theta)\|_{L^\infty\cap W^{1,2}}\lesssim2^{-j\beta}$ and $\sup_{\theta\in\R^1}\|\psi_{j'}\ast_\theta g_{\bar w}(\cdot,\theta)\|_{\Co^{\frac\alpha2-1}}\lesssim2^{-j'\frac\beta2}$. Thus by Lemma \ref{Lem::Exmp::HolSobLem} \ref{Item::Exmp::HolSobLem::Prod2}, we have for each $\theta\in\R^1$,
\begin{align*}
    \|P^1_l(\cdot,\theta)\|_{\Co^{\frac\alpha2-1}}&\le\|\psi_l\|_{L^1}\sum_{j=l-2}^{l+2}\sum_{j'=0}^{j-3}\|\psi_j\ast_\theta a(\cdot,\theta)\|_{L^\infty\cap W^{1,2}}\|\psi_{j'}\ast_\theta g_{\bar w}(\cdot,\theta)\|_{\Co^{\frac\alpha2-1}}\lesssim_{a,g}\sum_{j=l-2}^{l+2}\sum_{j'=0}^{j-3}2^{-j\beta}2^{-j'\frac\beta2}\lesssim2^{-l\beta};
    \\
    \|P^3_l(\cdot,\theta)\|_{\Co^{\frac\alpha2-1}}&\le\|\psi_l\|_{L^1}\sum_{\substack{j,j'\ge l-3\\|j-j'|\le2}}\|\psi_j\ast_\theta a(\cdot,\theta)\|_{L^\infty\cap W^{1,2}}\|\psi_{j'}\ast_\theta g_{\bar w}(\cdot,\theta)\|_{\Co^{\frac\alpha2-1}}\lesssim_{a,g}\sum_{\substack{j,j'\ge l-3\\|j-j'|\le2}}2^{-j\beta}2^{-j'\frac\beta2}\lesssim2^{-l\frac32\beta}.
\end{align*}
Taking supremum over $\theta\in\R^1$ we get $\sup_{l\ge0;\theta\in\R^1}2^{l\beta}(\|P^1_l(\cdot,\theta)\|_{\Co^{\frac\alpha2-1}}+\|P^3_l(\cdot,\theta)\|_{\Co^{\frac\alpha2-1}})<\infty$, bounding the first and the third term.

For the second sum we use $(\psi_j\ast_\theta a)\cdot(\psi_{j'}\ast_\theta g_{\bar w})=(\psi_j\ast_\theta(w^{-1} a))\cdot(\psi_{j'}\ast_\theta wg_{\bar w})=(\psi_j\ast_\theta(w^{-1} a))\cdot(\psi_{j'}\ast_\theta h_{\bar w})$. By Proposition \ref{Prop::Exmp::Expa} \ref{Item::Exmp::Expa::Reg2} $\sup_{\theta\in\R^1}\|\psi_j\ast_\theta(w^{-1} a)(\cdot,\theta)\|_{L^\frac{\alpha+4}2}\lesssim2^{-j\frac\beta2}$, by Claim 1 $\sup_{\theta\in\R^1}\|\psi_{j'}\ast_\theta h_{\bar w}(\cdot,\theta)\|_{L^{4+\frac{16}\alpha}}\lesssim2^{-j'\beta}$. And since the product map $L^\frac{\alpha+4}2\times L^{4+\frac{16}\alpha}\to L^\frac{4(\alpha+4)}{\alpha+8}$ is bounded, we have
{\small
\begin{align*}
    \|P^2_l(\cdot,\theta)\|_{L^\frac{4(\alpha+4)}{\alpha+8}}&\le\|\psi_l\|_{L^1}\sum_{j'=l-2}^{l+2}\sum_{j=0}^{j'-3}\|\psi_j\ast_\theta (w^{-1}a)(\cdot,\theta)\|_{L^\frac{\alpha+4}2}\|\psi_{j'}\ast_\theta h_{\bar w}(\cdot,\theta)\|_{L^{4+\frac{16}\alpha}}\lesssim\sum_{j'=l-2}^{l+2}\sum_{j=0}^{j'-3}2^{-j\frac\beta2}2^{-j'\beta}\lesssim2^{-l\beta}.
\end{align*}
}

By Sobolev embedding $L^\frac{4(\alpha+4)}{\alpha+8}(\C^1)\hookrightarrow\Co^{\frac\alpha{8+2\alpha}-1}(\C^1)$ (also see \cite[Chapter 2.7.1]{Triebel1}), we conclude that $\sup_{l\ge0;\theta\in\R^1}2^{l\beta}\|P^2_l(\cdot,\theta)\|_{\Co^{\frac\alpha{8+2\alpha}-1}(\C^1;\C)}<\infty$.

Since $\frac\alpha{8+2\alpha}<\frac\alpha2$, we conclude that $ag_{\bar w}\in \Co^{\frac\alpha{8+2\alpha}-1}\Co^\beta(\C^1,\R^1;\C)$. By Lemma \ref{Lem::Exmp::InvDBar} \ref{Item::Exmp::InvDBar::Bdd} we get $\Tc(a\cdot g_{\bar w})\in\Co^\frac\alpha{8+2\alpha}\Co^\beta(\B^2,\R^1;\C)\subset L^\infty\Co^{\beta}(\B^2,\R^1;\C)$, finishing the proof of Claim 2.

\medskip\noindent\textit{Final Step:}
We use $\Tc(a\cdot g_{\bar w})\in L^\infty_w\Co^{\beta}_{\theta}$ to obtain a contradiction.

Now rewriting the first equation in \eqref{Eqn::Exmp::ProofExampleddz::ghPDE} we have
\begin{equation}\label{Eqn::Exmp::ProofExampleddz::PDEofg}
    g=-\Tc(ag_{\bar w})+(g-\Tc g_w)+\Tc((\chi_w+a\chi_{\bar w})f)-\Tc((\chi-1) a_{\bar w})-\Tc( a_{\bar w}).
\end{equation}

By Claim 2, we have $\Tc(ag_{\bar w})\in L^\infty\Co^{\beta}(\B^2,\R^1;\C)$. By assumption $f\in\Co^{\alpha,\beta}_{w,\theta}$, we have $(\chi_w+a\chi_{\bar w})f\in\Co^{\alpha,\beta}(\B^2,\R^1;\C)$, so $\Tc((\chi_w+a\chi_{\bar w})f)\in \Co^1\Co^\beta(\B^2,\R^1;\C)\subset L^\infty\Co^\beta(\B^2,\R^1;\C)$ as well.

Since $g\in\Co^{\alpha,\beta}_{w,\theta}\subset L^\infty_w\Co^\beta_\theta$ is supported in $\{|w|<\frac12\}\times\R^1$, by \eqref{Eqn::Exmp::InvDBar::Holo} $g-\Tc g_w\in L^\infty\Co^\beta(\B^2,\R^1;\C)$  as well.

Let $\mu(\theta):=\tilde\chi(2\theta)$, we have $\mu|_{(-\frac14\eps_0,\frac14\eps_0)}\equiv1$. Thus by \eqref{Eqn::Exmp::ProofExampleddz::Chi}, for every $(w,\theta)\in\frac{\eps_0}2\B^2\times(-\frac{\eps_0}4,\frac{\eps_0}4)$, using $\mu(\chi-1)|_{\frac{\eps_0}2\B^2\times(-\frac{\eps_0}4,\frac{\eps_0}4)}=0$ we have
\begin{equation*}
    \Tc((\chi-1) a_{\bar w})(w,\theta)=\mu(\theta)\Tc((\chi-1) a_{\bar w})(w,\theta)=\Tc(\mu(\chi-1) a_{\bar w})(w,\theta)=\big(\Tc(\mu(\chi-1) a_{\bar w})-\mu(\chi-1) a_{\bar w}\big)(w,\theta).
\end{equation*}

By assumption $\mu(\chi-1) a_{\bar w}$ is supported in $\{\frac{\eps_0}2<|w|<\eps_0\}\times(-\frac{\eps_0}4,\frac{\eps_0}4)$, so by Proposition \ref{Prop::Exmp::Expa} \ref{Item::Exmp::Expa::Reg1} we have $\mu(\chi-1) a_{\bar w}\in L^\infty\Co^\beta(\B^2,\R^1;\C)$. Using Lemma \ref{Lem::Exmp::InvDBar} \ref{Item::Exmp::InvDBar::Holo} we get $\Tc((\chi-1) a_{\bar w})\in L^\infty\Co^\beta(\B^2,\R^1;\C)$ as well.

However, by Proposition \ref{Prop::Exmp::Expa} \ref{Item::Exmp::Expa::Unbdd} we have $\Tc(a_{\bar w})\notin L^\infty_w\Co^\beta_{\theta}$ near $(0,0)$. Therefore near $(w,\theta)=(0,0)$, the left hand side of \eqref{Eqn::Exmp::ProofExampleddz::PDEofg} is in $L^\infty_w\Co^\beta_\theta$ while the right hand side is not, this contradicts to the assumption $f\in\Co^{\alpha,\beta}(U_2,V_2;\C)$ and finishes the proof.
\end{proof}

To summarize this subsection, what we obtained is the following:
\begin{thm}\label{Thm::Exmp::SharpSum}
    Let $\alpha>\frac12$ and $\beta>0$, let $a\in\Co^{\alpha,\beta}(\C^1_w,\R^1_\theta;\C)$ be as in Theorem \ref{Thm::Exmp::ProofExampleddz} (see \eqref{Eqn::Exmp::Expa::Defofa}). Let $\Se:=\Span(\partial_w+a(w,\theta)\partial_{\bar w})$ be a rank 1 complex Frobenius structure in $\C^1\times\R^1$. 
    
    Then for any neighborhood $U\times V\subseteq\C^1_w\times\R^1_\theta$ of $(0,0)$ there is no continuous map $F:U\times V\to\C^1_z\times\R^1_s$ such that $F$ is homeomorphic to its image, $\nabla_wF$ and $\nabla_z(F^\Inv)$ are both continuous, and $F^*\Coorvec z$ spans $\Se|_{U\times V}$.
    
    In particular we get Theorem \ref{MainThm::Sharpddz}: if $\alpha=\beta>1$, then there is no $C^1$-coordinate chart $F:U\times V\subseteq\C^1_z\times\R^1_s$ near $(0,0)$ such that $F^*\Coorvec z$ spans $\Se|_{U\times V}$.
\end{thm}
\begin{proof}
Suppose, by contrast, such continuous map $F$ exists. Then by Proposition \ref{Prop::Exmp::SharpddzRed} \ref{Item::Exmp::SharpddzRed::Chart} $\Rightarrow$ \ref{Item::Exmp::SharpddzRed::PDE} there are a neighborhood $U_2\times V_2\subseteq U\times V$ of $(0,0)$ and a function $f\in\Co^{\alpha,\beta}(U_2,V_2;\C)$ that solves the PDE $f_w+af_{\bar w}=-a_{\bar w}$. However by Theorem \ref{Thm::Exmp::ProofExampleddz} such $f$ would not exist, contradiction.

When $\alpha=\beta>1$, a $C^1$-coordinate chart $z=F(w)$ implies $\nabla_wF\in C^0$ and $\nabla_z(F^\Inv)\in C^0$ automatically, so the assumptions of Proposition \ref{Prop::Exmp::SharpddzRed} \ref{Item::Exmp::SharpddzRed::Chart} are satisfied. By Proposition \ref{Prop::Exmp::SharpddzRed} \ref{Item::Exmp::SharpddzRed::Chart} $\Rightarrow$ \ref{Item::Exmp::SharpddzRed::PDE} and Theorem \ref{Thm::Exmp::ProofExampleddz} the proof follows.
\end{proof}

We can extend Theorem \ref{MainThm::Sharpddz} to general $r,q\ge0$ and $m\ge1$ following from a tensor construction that is similar to the proof of Theorem \ref{MainThm::CounterNN} in Section \ref{Section::CountNN}.
\begin{cor}\label{Cor::Exmp::ExdSharpddz}
Let $r\ge0$ and $m,q\ge1$ be integers. Let $\alpha>1$.  Then there is a $\Co^\alpha$-complex Frobenius structure $\Se$ on $\R^{r+2m+q}$, such that $\rank\Se =r+m$, $\rank(\Se+\bar\Se)=r+2m$, and there does not exist a $C^1$-coordinate chart $F:U\subseteq\R^{r+2m+q}\to\R^r_t\times\C^m_z\times\R^q_s$ near $0\in\R^{r+2m+q}$ such that $F^*\Coorvec{t^1},\dots,F^*\Coorvec{t^r},F^*\Coorvec{z^1},\dots,F^*\Coorvec{z^m}$ spans $\Se|_U$ and $F^*\Coorvec{t^1},\dots,F^*\Coorvec{t^r},F^*\Coorvec {z^1},\dots,F^*\Coorvec{z^m}\in\Co^\alpha(U;\C^{r+2m+q})$.
\end{cor}
\begin{remark}
The result says that if $F$ is a $C^1$-coordinate chart representing $\Se$, then one of the vector field $F^*\Coorvec{t^1},\dots,F^*\Coorvec{t^r},F^*\Coorvec {z^1},\dots,F^*\Coorvec{z^m}$ cannot $\Co^\alpha$. It is not known to the author whether one of $F^*\Coorvec {z^1},\dots,F^*\Coorvec{z^m}$ must be $\Co^{\alpha-}\backslash\Co^\alpha$.
\end{remark}
\begin{proof}
 Note that $\R^{r+2m+q}\simeq\R^r\times\C^m\times\R^q$ as real smooth manifolds. Let $\tau=(\tau^1,\dots,\tau^r)$, $w=(w^1,\dots,w^m)$ and $\theta=(\theta^1,\dots,\theta^q)$ be the standard coordinate systems for $\R^r$, $\C^m$ and $\R^q$, respectively.

Let $a\in\Co^\alpha(\C^1\times\R^1;\C)$ be the function  as in Theorem \ref{Thm::Exmp::ProofExampleddz} (see \eqref{Eqn::Exmp::Expa::Defofa}, with $\beta=\alpha$). Let $Z(\tau,w,\theta):=\Coorvec{w^1}+a(w^1,\theta^1)\Coorvec{\bar w^1}$. We define $\Se\le \C T(\R^r\times\C^m\times\R^q)$ by $$\textstyle\Se:=\Span(\Coorvec{\tau^1},\dots,\Coorvec{\tau^r},Z,\Coorvec{w^2},\dots,\Coorvec{w^m}).$$ Clearly $\rank\Se=r+m$ and $\rank(\Se+\bar\Se)=r+2m$.

Now suppose $U\subseteq\R^r_\tau\times\C^m_w\times\R^q_\theta$ is a neighborhood of $(0,0,0)$ and $F=(F',F'',F'''):U\to\R^r_t\times\C^m_z\times\R^q_s$ is a $C^1$ coordinate chart such that $F(0,0,0)=(0,0,0)$, $\Se|_U=\Span(F^*\Coorvec t,F^*\Coorvec z)$ and $F^*\Coorvec t\in\Co^\alpha(U;\R^{r\times(r+2m+q)})$, $F^*\Coorvec z\in\Co^\alpha(U;\C^{m\times(r+2m+q)})$. 
For convenience we define $A\in\Co^\alpha(\R^r\times\C^m\times\R^q;\C^{m\times m})$ as $A(\tau,w,\theta)=\begin{pmatrix}
a(w^1,\theta^1)&0_{1\times(m-1)}\\0_{(m-1)\times1}&0_{(m-1)\times(m-1)}
\end{pmatrix}$. Since $\Se^\bot|_U=\Span(F^*d\bar z,F^*ds)=\Span(d\bar F'',dF''')$ and $(\Se+\bar\Se)|_U=\Span(\Coorvec\tau,\Coorvec w,\Coorvec{\bar w)}|_U=\Span (dF'',d\bar F'',dF''')^\bot$, similar to \eqref{Eqn::Exmp::SharpddzRed::Exp1}, using matrix notation in Section \ref{Section::Convention} (also see \eqref{Eqn::EllipticPara::GradH}) we have $F''_\tau=0$ and 
\begin{equation*}
    \textstyle F''_{\bar w}+A\cdot F''_w =0\in\C^{m\times m},\quad\text{i.e.}\quad\tfrac{\partial F''^j}{\partial\bar w^k}+\sum_{l=1}^mA_k^l\tfrac{\partial F''^j}{\partial w^l}=\tfrac{\partial F''^j}{\partial\bar w^k}+\delta_k^1\cdot a\cdot\tfrac{\partial F''^j}{\partial w^1}=0\text{ for all }1\le j,k\le m.
\end{equation*}
Therefore similar to \eqref{Eqn::Exmp::SharpddzRed::Exp2} we have
\begin{equation*}
    \begin{pmatrix}
    F^*dt&F^*dz&F^*d\bar z&F^*ds
    \end{pmatrix}
    =
    \begin{pmatrix}
    d\tau&dw&d\bar w&d\theta
    \end{pmatrix}
    \begin{pmatrix}
    F'_\tau
    \\
    F'_w&F''_w&-A\bar F''_{\bar w}
    \\
    F'_{\bar w}&-\bar AF''_w&\bar F''_{\bar w}
    \\
    F'_\theta&F''_\theta&\bar F''_\theta&F'''_\theta
\end{pmatrix}.
\end{equation*}
Taking transposition we get
\begin{equation}\label{Eqn::Exmp::Sharpddz::ExtPf1}
    \begin{pmatrix}
    \partial_\tau\\\partial_w\\\partial_{\bar w}
    \end{pmatrix}
    =
    \begin{pmatrix}
    F'_\tau
    \\
    F'_w&F''_w&-A\bar F''_{\bar w}
    \\
    F'_{\bar w}&-\bar AF''_w&\bar F''_{\bar w}
\end{pmatrix}
\begin{pmatrix}
F^*\partial_t\\F^*\partial_z\\F^*\partial_{\bar z}
\end{pmatrix}.
\end{equation}

By assumption $F^*\Coorvec t,F^*\Coorvec z,F^*\Coorvec{\bar z}$ are all $\Co^\alpha$, so all matrix blocks in \eqref{Eqn::Exmp::Sharpddz::ExtPf1} are $\Co^\alpha$ as well. Since $F^*\Coorvec t=(F'_\tau)^{-1}\cdot\Coorvec\tau$, we see that 
$\begin{pmatrix}F'_w&F''_w&-A\bar F''_{\bar w}
    \\
-\bar AF''_w&\bar F''_{\bar w}
\end{pmatrix}$ is invertible at every point in $U$. Therefore there exists a $1\le j_0\le m$ such that $F''^{j_0}_{w^1}(0,0)=\frac{\partial F''^{j_0}}{\partial w^1}\neq0$.

Let $\Gamma:\C^1_{w^1}\times\R^1_{\theta^1}\to\R^r_\tau\times\C^m_w\times\R^q_\theta$, $\Gamma(w^1,\theta^1):=(0^r,(w^1,0^{m-1}),(\theta^1,0^{q-1}))$ be the natural embedding, and let $U_1:=\{(w^1,\theta^1):\Gamma(w^1,\theta^1)\in U\text{ and }|F''^{j_0}_{ w^1}(\Gamma(w^1,\theta^1))-F''^{j_0}_{w^1}(\Gamma(0,\theta^1))|<\frac12|F''^{j_0}_{w^1}(\Gamma(0,\theta^1))|\}\subseteq\C^1\times\R^1$. We define a function $f:U_1\to\C$ by (cf. \eqref{Eqn::Exmp::SharpddzRed::DefEqnf})
\begin{equation*}
    f(w^1,\theta^1):=\log\frac{\bar F''^{j_0}_{\bar w^1}(\Gamma(w^1,\theta^1))}{\bar F''^{j_0}_{\bar w^1}(\Gamma(0,\theta^1))}=\sum_{n=1}^\infty\frac{(-1)^{n-1}}n\bigg(\frac{\bar F''^{j_0}_{\bar w^1}(\Gamma(w^1,\theta^1))}{\bar F''^{j_0}_{\bar w^1}(\Gamma(0,\theta^1))}-1\bigg)^n.
\end{equation*}

Since $F''_w\in\Co^\alpha$ we see that $f\in\Co^\alpha(U_1;\C)$. Since $Z\bar F''^{j_0}=(\partial_{w^1}+a(w^1,\theta^1)\partial_{\bar w^1})\bar F''^{j_0}=0$, using the same computation to \eqref{Eqn::Exmp::SharpddzRed::PfTmp1} we see that $f$ solves \eqref{Eqn::Exmp::SharpddzRed::aPDE}. Therefore by Proposition \ref{Prop::Exmp::SharpddzRed} and Theorem \ref{Thm::Exmp::ProofExampleddz}, we reach a contradiction, i.e. such $f\in\Co^\alpha$ does not exist. 

To conclude, one of $F^*\Coorvec{t^1},\dots,F^*\Coorvec{t^r},F^*\Coorvec{z^1},\dots,F^*\Coorvec{z^m}$ cannot be $\Co^\alpha$, and we finish the proof. 
\end{proof}

\subsection{Canonical coordinate system lose 1 derivative}\label{Section::CanCoor-1Reg}
As mentioned in Section \ref{Section::Rough1FormRmkCoor}, DeTurck and Kazdan showed that a Riemannian metric may not have optimal regularity with respect to
geodesic normal coordinates \cite[Example 2.3]{DeTurckKazdan}. A natural analog of geodesic normal coordinates for vector fields
are canonical coordinates (of the first kind).
Next, we show that vector fields may not have optimal regularity
with respect to these canonical coordinates.


Given $C^1$-vector fields $X_1,\dots,X_n$ on $\Mf$ that form a basis on the tangent space at every point, recall from Definitions \ref{Defn::ODE::MultiFlow} and \ref{Defn::ODE::CanCoor} that the canonical coordinates at $p\in\Mf$ is the map $\Phi_p(t^1,\dots,t^n):=e^{t^1X_1+\dots+t^nX_n}p$ defined via solving the ordinary differential equation, provided that it is solvable.
When $X_1,\dots,X_n\in\Co^\alpha$ for some $\alpha>1$, classical regularity theorems for ODEs (see Lemma \ref{Lem::ODE::ODEReg}) show that $\Phi_p$ is at least $\Co^\alpha$. 
Therefore, $\Phi_p^{*}X_1,\ldots, \Phi_p^{*}X_n$ are at least $\Co^{\alpha-1}$; which is one derivative less than the original
regularity of $X_1,\ldots, X_n$.  The following result
shows that this loss of one derivative is sometimes inevitable.


\begin{lem}\label{Lem::Exmp::CanonicalCoords}
Endow $\R^2$ with standard coordinate system $(x,y)$. Let $\alpha>1$ and let $X:=\partial_x$ and $Y:=xf(y)\partial_x+\partial_y$ where $f(y):=\alpha \max(0,y)^{\alpha-1}$.

Then we can find a  new $\Co^{\alpha+1}$ atlas $\Af$ on $\R^2$ which is compatible with the standard $\Co^\alpha$-structure on $\R^2$, such that $X,Y$ are $\Co^\alpha_\loc$ with respect to this the new atlas, but for the canonical coordinates $\Phi(t,s):=e^{tX+sY}(0)$ we have $\Phi^*Y\notin\Co^{\alpha-1+\eps}$ near $(0,0)$, in particular the collection $(\Phi^*X,\Phi^*Y)$ is not $\Co^{\alpha-1+\eps}$.
\end{lem}


Note that $X$ and $Y$ form a local basis of the tangent space at every point, and $f\in\Co^{\alpha-1}_\loc(\R)$.
\begin{proof}

First we show the existence of the new atlas $\Af$ with respect to which $X$ and $Y$ are $\Co^\alpha$. In particular $X,Y$ are $C^1$ in $\Af$, so $\Phi$ is well-defined.

Note that $[X,Y]=f(y)\partial_x\in\Co^{\alpha-1}(\R^2;\R^2)$. Specifically, the dual basis of $X,Y$ are 1-forms $\lambda=dx-xf(y)dy,\eta=dy$ which satisfy that $d\lambda=f(y)dx\wedge dy$ and $d\eta=0$ are both $\Co^{\alpha-1}$ 2-forms, so the condition \ref{Item::Rough1Form::BasicCaseforMainThm::b} in Corollary \ref{Cor::Rough1Form::BasicCaseforMainThm}  is satisfied. By Corollary \ref{Cor::Rough1Form::BasicCaseforMainThm} we can find a $\Co^{\alpha+1}$-atlas $\Af$ on $\R^2$ such that $X,Y$ are both $\Co^\alpha_\loc$ on $\Af$.

Since $X,Y$
are $\Co^{\alpha}\subsetneq C^1$ with respect to $\Af$, we see that $\Phi(t,s)$ is well-defined near $(t,s)=(0,0)$. We can compute $\Phi$ in terms of $f$. 

Clearly $\Phi(t,s)=(*,s)$ since $e^{s\partial_y}(x,y)=(x,y+s)$. We can write $\Phi(t,s)=(\phi(t,s),s)$. Define $\Phi(t,s;r)$ for $r\in\R$ as the solution to the ODE $\Coorvec r\Phi(t,s;r)=tX(\Phi(t,s;r))+sY(\Phi(t,s;r))$, $\Phi(t,s;0)=0$. So $\Phi(t,s)=\Phi(t,s;1)$ and we have $\Phi(t,s;r)=(\phi(t,s;r),rs)$ where
\begin{equation*}
    \Coorvec r\phi(t,s;r)=t+sf(rs)\phi(t,s;r),\quad\phi(t,s;0)=0,\quad t,s,r\in\R.
\end{equation*}
Solving this ODE we have
\begin{equation*}
    \phi(t,s;r)=e^{\int_0^r sf(\rho s)d\rho}\int_0^re^{-\int_0^\rho sf(\mu s)d\mu}td\rho,\quad\phi(t,s)=\phi(t,s;1)=t\frac{e^{-\int_0^sf(\rho)d\rho}}s\int_0^se^{\int_0^{\rho }f(\mu)d\mu}d\rho.
\end{equation*}

Now plug in \[f(y)=\begin{cases}\alpha y^{\alpha-1},&y\ge0,\\0&y\le0.\end{cases}\] We have
\begin{equation*}
    \phi(t,s)=t \frac {e^{-s^\alpha}}s\int_0^se^{\rho^\alpha} d\rho\quad\text{when } s>0;\quad \phi(t,s)=t\quad\text{when } s\le0.
\end{equation*}
Thus, $\phi(t,s)=tg(s)$ where 
\begin{equation}\label{Eqn::Rough1Form::RmkCoor::g}
    g(s)=\begin{cases}\frac {e^{-s^\alpha}}s\int_0^se^{\rho^\alpha} d\rho&s>0\\1&s\le0.\end{cases}
\end{equation}

We are going to show $\Phi\in\Co^\alpha_\loc(\R^2;\R^2)$ and $\Phi\notin\Co^{\alpha+\eps}$ near $(0,0)$. To see this, it suffices to show $g\in\Co^\alpha_\loc(\R)$ and $g\notin\Co_{\loc}^{\alpha+\eps}$ near $0$, for every $\eps>0$.

By Taylor's expansion on the exponential function we have, when $s>0$,
\begin{align*}
    g(s)&=e^{-s^\alpha}\frac1s\int_0^se^{\rho^\alpha}d\rho=\sum_{j=0}^\infty\frac{(-1)^js^{j\alpha}}{j!}\frac1s\int_0^s\sum_{k=0}^\infty\frac{\rho^{k\alpha}}{k!}d\rho=\sum_{j,k=0}^\infty\frac{(-1)^j}{j!k!}\frac{s^{j\alpha}s^{k\alpha}}{k\alpha+1}
    \\&=\sum_{l=0}^\infty\Big(\sum_{k=0}^l\frac{(-1)^{l-k}}{k\alpha+1}\Big)s^{l\alpha}=1-\big(1-\tfrac1{\alpha+1}\big)s^\alpha+O(s^{2\alpha}).
\end{align*}
In other words
\begin{equation*}
    g(s)=\sum_{l=0}^\infty\Big(\sum_{k=0}^l\frac{(-1)^{l-k}}{k\alpha+1}\Big)\max(s,0)^{l\alpha}=1-\big(1-\tfrac1{\alpha+1}\big)\max(s,0)^\alpha+\sum_{l=2}^\infty\Big(\sum_{k=0}^l\frac{(-1)^{l-k}}{k\alpha+1}\Big)\max(s,0)^{l\alpha},\quad\forall s\in\R.
\end{equation*}

Note that for $\beta>0$ the function $\max(s,0)^\beta=\begin{cases}s^\beta&s\ge0\\0&s\le0\end{cases}$ is $\Co^\beta_\loc$ but not $\Co^{\beta+\eps}$ near $0$.
Indeed when $0<\beta<1$ we know that $\max(s,0)^\beta\in C^{0,\beta}_\loc=\Co^\beta_\loc$ and is not $C^{0,\beta+\eps}=\Co^{\beta+\eps}$ near $0$ for any $0<\eps<1-\beta$, when $\beta=1$ we know $\max(s,0)\in C^{0,1}_\loc\subsetneq\Co^1_\loc$ and is not $C^1$ near 0 so is not $\Co^{1+\eps}$ for any $\eps>0$. For $\beta>1$ by passing to its derivatives we see that $\max(s,0)^\beta\in\Co^\beta_\loc$ and is not $\Co^{\beta+\eps}$ near $0$.

So the remainder $\sum_{l=2}^\infty\big(\sum_{k=0}^l\frac{(-1)^{l-k}}{k\alpha+1}\big)\max(s,0)^{l\alpha}$ is $\Co^{2\alpha}\subset \Co^{\alpha+\eps}_\loc$ for all $0<\eps\le\alpha$, while the main term $1-\big(1-\tfrac1{\alpha+1}\big)\max(s,0)^\alpha$ is $\Co^\alpha_\loc$ but not $\Co^{\alpha+\eps}$ near 0. Therefore we conclude that $g\in\Co^\alpha(\R)$, but $g\notin\Co^{\alpha+\eps}$ near $s=0$.


Now we know $\Phi\in\Co^{\alpha}_\loc(\R^2;\R^2)$ but $\Phi\notin\Co^{\alpha+\eps}_\loc$ near $(t,s)=(0,0)$. Consider the inverse function of $\Phi$, 
and set $(u(x,y),v(x,y)):=\Phi^{-1}(x,y)$; so that $\Phi^*Y=(Yu,Yv)\circ\Phi$. We have $v(x,y)=y$ and $u(x,y)=\frac1{g(y)}x$. Note that $g(s)>0$ for every $s$, so $y\mapsto\frac1{g(y)}$ is not $\Co^{\alpha+\eps}$ near $y=0$, for any $\eps>0$. Therefore $Yv=x\cdot\partial_y\frac1{g(y)}$ is not $\Co^{\alpha+\eps-1}$ near $(x,y)=(0,0)$. By composing with $\Phi$ which is a $\Co^\alpha_\loc\subset\Co^{\alpha+\eps-1}_\loc$-diffeomorphism (for $0<\eps<1$), we see that $\Phi^*Y$ is not $\Co^{\alpha+\eps-1}$ near $(t,s)=(0,0)$, for every $\eps>0$. 
\end{proof}
\begin{remark}
As a differentiable map $\Phi:(\R^2_{t,s},\mathrm{std})\to(\R^2,\Af)$ between two $\Co^{\alpha+1}$-manifolds, we see that $\Phi$ is not $\Co^{\alpha+\eps}$ near $(0,0)$ for any $\eps>0$. Otherwise since $X$ and $Y$ are $\Co^\alpha_\loc$ on $\Af$, we have $\Phi^*Y\in\Co^{\min(\alpha,\alpha-1+\eps)}=\Co^{\alpha-1+\min(\eps,1)}$ near $(0,0)$, contradicting to Lemma \ref{Lem::Exmp::CanonicalCoords}.
\end{remark}
\begin{remark}
We can generalize Lemma \ref{Lem::Exmp::CanonicalCoords} to higher dimension: on $\R^n$, we can take $X_1:=\Coorvec{x_1}$, $X_2:=x_1\max(0,x_2)^{\alpha-1}\Coorvec{x_1}+\Coorvec{x_2}$ and $X_j:=\Coorvec{x_j}$ for $3\le j\le n$. Thus there exists a coordinate atlas $\Af$ with respect to which $X_1,\dots,X_n$ are $\Co^\alpha$, while for the canonical coordinate map $\Phi$ we still have $\Phi^*X_2\notin\Co^{\alpha-1+\eps}$ near $0\in\R^n$, for all $\eps>0$. We leave the proof to reader.
\end{remark}

\section{A Quantitative Frobenius Theorem}\label{Section::QuantFro}
\subsection{A few remarks on quantitative estimates}

To make the Frobenius theorem quantitative, instead of discussing the (regular) tangent subbundle $\V$, we consider a set of real generators $X_1,\dots,X_m$ of $\V$, that also label the intrinsic scale of $\V$.
It is more practical to focus only on a leaf, unless there are extra intrinsic conditions for $X_1,\dots,X_m$ that allow us to estimate the parameter dependence between different leaves. We refer \cite{StreetQuantFrobenius} to readers for more discussion between qualitative and quantitative Frobenius theorems.

The involutive condition of $\V$ shows that there are functions $c_{ij}^k$ such that $[X_i,X_j]=\sum_{k=1}^mc_{ij}^kX_k$ (see Proposition \ref{Prop::DisInv::CharInv1} \ref{Item::DisInv::CharInv1::Gen2}). Conversely the singular Frobenius theorem (say Theorem \ref{MainThm::SingFro}) shows that if $[X_i,X_j]=\sum_{k=1}^mc_{ij}^kX_k$ and $X_i,c_{ij}^k$ have certain smoothness assumptions, then we can find an (at least $C^1$) submanifold $\Sf$ containing a fixed point, such that $X_1,\dots,X_m$ spans $T\Sf$. In particular for a regular parameterization $\Phi:\Omega\subseteq\R^r\to\Sf$ where $r=\dim\Sf$, we can talk about the pullback vector fields $\Phi^*X_1,\dots,\Phi^*X_m$ on $\Omega$.

Now we can to find a better $\Phi$ such that $\Phi^*X_1,\dots,\Phi^*X_m$ are ``normalized''. More precisely, we want $\Phi$ satisfies the following:
\begin{itemize}[parsep=-0.3ex]
    \item $\Phi^*X_1,\dots,\Phi^*X_m$ are ``at the unit scale'': $\sup\limits_{1\le i\le m}\|\Phi^*X_i\|_{C^0}\approx1$ and $\sup\limits_{i_1,\dots,i_n}\inf\limits_{\Omega}|\Phi^*X_{i_1}\wedge\dots\wedge\Phi^*X_{i_r}|\gtrsim1$.
    \item For a fixed $\alpha>0$, we have $\sup_{1\le i\le m}\|\Phi^*X_i\|_{\Co^\alpha}\approx1$, provided that there exists a $\Phi$ such that $\Phi^*X_1,\dots,\Phi^*X_m\in\Co^\alpha(\Omega;\R^r)$.
\end{itemize}
Here all the implied constant depend only on the diffeomorphic invariant quantities  like the norms $\|c_{ij}^k\|_{\Co^{\alpha-1}_X}$ (see \cite[Section 5.1]{CoordAdaptedII}). 

It is critical that the quantities should not depend on the choice of the original coordinate system.
Otherwise, if the vector fields are given in a coordinate system where the above upper and lower bounds
are bad, the estimates may be bad as well.

The first point is addressed in \cite[Section 3]{StreetQuantFrobenius} where we require $c_{ij}^k\in C^2_X$. The second point is discussed in \cite{CoordAdapted,CoordAdaptedII} where $s$ satisfies $s>2$, and in \cite{StreetYaoVectorFields} we improve the assumption to $s>1$.
We also refer \cite{CoordAdaptedIII} to the case where $\Phi^*X_1,\dots,\Phi^*X_m $ can be real-analytic.

If $X_1,\dots,X_m$ are $C^1$, then Theorem \ref{MainThm::QuantFro} is a corollary to \cite[Theorem 9.1]{StreetYaoVectorFields}. The only thing we need to change in the thesis is to apply the log-Lipschitz version of singular Frobenius theorem (Theorem \ref{MainThm::SingFro}) instead of the $C^1$-version (say \cite[Proposition 3.1]{CoordAdapted}).

\begin{remark}
Alternatively when $X_1,\dots,X_m$ span the tangent spaces at every point, we can relax the log-Lipschitz assumption to merely $X_1,\dots,X_m\in\Co^{\frac12+}$, where $[X_i,X_j]=\sum_kc_{ij}^kX_k$ make sense as distribution. Indeed in this case we do not need the apply singular Frobenius theorem since the total space $\Mf$ itself is the leaf we work on. And by Corollary \ref{Cor::Rough1Form::BasicCaseforMainThm} there is a new atlas $\Af$ for which $X_1,\dots,X_m$ are $\Co^\alpha\subset C^1$. For this case we leave the detail to readers.
\end{remark}

\begin{remark}
There are two ingredients in \cite[Theorem 9.1]{StreetYaoVectorFields} which are not included in Theorem \ref{MainThm::QuantFro}: 
\begin{itemize}[parsep=-0.3ex]
    \item Higher order function spaces along $X$: in Theorem \ref{MainThm::QuantFro} we only discuss the case where $c_{ij}^k\in\Co^s_X$ for $0<s<1$. There is a natural way to generalize the Banach space $\Co^s_X(\Mf)$ to $s\ge1$, see \cite[Sections 2.2 and 8]{CoordAdapted}, also see \cite[Proposition 7.1]{CoordAdaptedII}. 
    \item Interior Schauder's estimate: if $c_{ij}^k\in\Co^{\beta-1}_X$ for some $\beta>\alpha$, then $\Phi^*X_1,\dots,\Phi^*X_m$ are automatically $\Co^\beta$, with their $\Co^\beta$-norm bounded by a constant that depends only on the previous quantities with the upper bound of $\sup_{i,j,k}\|c_{ij}^k\|_{\Co^{\beta-1}_X}$. See \cite[Proposition 9.4 (ii)]{StreetYaoVectorFields}.
\end{itemize}

We do not include these two in the thesis, as they do not add substantial difficulty in the proof.
\end{remark}

\begin{remark}
In Theorem \ref{MainThm::QuantFro}, by considering the coordinate chart $F=(F^1,\dots,F^n):=\Phi^\Inv:\Phi(\B^n)\to\R^n$ in Theorem \ref{MainThm::QuantFro}, we also have $F\in C^{2,\alpha-1}_X(\Phi(\B^n);\R^n)$ (or equivalently $F\in \Co^{\alpha+1}_X(\Phi(\B^n);\R^n)$). Indeed Theorem \ref{MainThm::QuantFro} gives $\|\Phi^*X\|_{\Co^\alpha}\lesssim_{n,m,\alpha,\mu_0,M_0}1$, thus $\|f\|_{\Co^\alpha_X(\Phi(\B^n))}\approx\|f\circ\Phi\|_{\Co^\alpha_{\Phi^*X}(\B^n)}\approx\|f\circ\Phi\|_{\Co^\alpha(\B^n)}$. Using $\Phi^*X=(XF)\circ\Phi$, we have $\|F\|_{\Co^{\alpha+1}_X}=\|F\|_{\Co^\alpha_X}+\|XF\|_{\Co^\alpha_X}\approx\|F\circ\Phi\|_{\Co^{\alpha}}+\|XF\circ\Phi\|_{\Co^\alpha}$. The result $F\in \Co^{\alpha+1}_X(\Phi(\B^n);\R^n)$ then follows.

\end{remark}

Theorem \ref{MainThm::QuantFro}, as well as \cite[Theorem 9.1]{StreetYaoVectorFields}, is still open for the case $0<\alpha\le1$. We remark the difficulty in the formulation of the problem and the proof techniques.

To formulate the correspondent result for $0<\alpha\le1$, if we still use the characterization $[X_i,X_j]=\sum_kc_{ij}^kX_k$ we need to explain the meaning of $c_{ij}^k\in\Co^{\alpha-1}_X(\Mf)$. In \cite{CoordAdapted,CoordAdaptedII} we only define $\Co^s_X(\Mf)$ for $s>0$. If we know $X_1,\dots,X_m\in\Co^\beta_\loc$ for some $\beta>0$, in \cite{StreetYaoVectorFields} we can define $\Co^s_{X,\loc}(\Mf)$ for $s>-\beta$, where we have not define the norm $\|\cdot\|_{\Co^s_X}$ for $-\beta<s\le0$.

One way to define $\|\cdot\|_{\Co^s_X}$ for negative $s$ is to modify the corresponding negative H\"older norm in $\R^n$. For example we can define $\|f\|_{\Co^s_X}:=\inf\{\sum_{j=1}^m\|g_j\|_{\Co^{s+1}_X}:f=\sum_{j=1}^mX_jg_j\text{ as distributions}\}$ for $-\beta<s\le0$. 
But in this setting $c_{ij}^k\in\Co^{\alpha-1}_X(\Mf)$ make sense only when we know $X_1,\dots,X_m\in\Co^{\frac12+}$ under some good coordinate chart, since otherwise $\alpha-1\le-\alpha$.

An alternative way to state the problem is to use Corollary \ref{Cor::Rough1Form::BasicCaseforMainThm} \ref{Item::Rough1Form::BasicCaseforMainThm::b} (also see \cite[Theorem 3.1 (b)]{StreetYaoVectorFields}), where we take $(\lambda^1,\dots,\lambda^n)$ to be the dual basis of $(X_1,\dots,X_n)$, we can use the quantities ``$\|d\lambda^i\|_{\Co^{\alpha-1}_X}$'' instead of the $\|c_{ij}^k\|_{\Co^{\alpha-1}_X}$. In \cite[Lemma 4.8]{StreetYaoVectorFields} we can define ``$\|d\lambda^i\|_{\Co^{\alpha-1}_X}$'' by the infimum of all $\|\eta^i\|_{\Co^\alpha_X}$ such that $d\lambda^i=d\eta^i$, which are well-defined quantities.

In this way we can formulate the problem as the following. For simplicity we only consider the case $m=n$, i.e. $X_1,\dots,X_n$ are $\Co^\alpha$-vector fields that span the tangent space at every point in $\Mf$.
\begin{prob}
Let $0<\alpha\le1$ and $\mu_0,M_0>0$, can we find constants $\hat K(n,\alpha,\mu_0,M_0)$ and $K_0(n,\alpha,\mu_0,M_0)$ such that the following holds?

Let $\Mf$ be a $\Co^{\alpha+1}$-manifold and let $X_1,\dots,X_n$ be $\Co^\alpha$ 1-forms on $\Mf$ that form a basis on the tangent space at every point. Let $g$ be the (unique) $\Co^\alpha$ Riemannian metric of $\Mf$ such that $X_1,\dots,X_n$ are orthonormal at every point, and let $\Co^\alpha$ 1-forms $(\lambda^1,\dots,\lambda^n)$ be the dual basis of $(X_1,\dots,X_n)$. 

Suppose $(\Mf,g)$ is a complete manifold and satisfies the following:
\begin{itemize}[parsep=-0.3ex]
    \item For every $\theta=(\theta^1,\dots,\theta^n)\in\Sp^{n-1}\subset\R^n$ and every $\gamma:[0,\mu_0]\to\Mf$ such that $\dot\gamma(t)=\sum_{j=1}^n\theta^jX_j(\gamma(t))$ for $0\le t\le\mu_0$, we have $\gamma(\mu_0)\neq\gamma(0)$.
    \item For the H\"older norm on $\Co^\alpha_g=\Co^\alpha(\Mf,g)$, there are $\Co^\alpha$ 1-forms $\eta^1,\dots,\eta^n$ such that $d\lambda^i=d\eta^i$ for all $1\le i\le n$ and
    \begin{equation}
        \sum_{j=1}^n\|\lambda^i\|_{\Co^\alpha(\Mf,g;T^*\Mf)}+\|\eta^i\|_{\Co^\alpha(\Mf,g;T^*\Mf)}<M_0.
    \end{equation}
    
    Then for any $p_0\in\Mf$ there is a $\Co^{\alpha+1}$-regular parameterization $\Phi:\B^n\subset\R^n\to\Mf$ such that $\Phi(0)=p_0$ and $\Phi^*X=[\Phi^*X_1,\dots,\Phi^*X_n]^\top$ has the form $\Phi^*X=\hat K\cdot(I+A)\nabla$, where $A(0)=0$ and $\|A\|_{\Co^\alpha(\B^n;\R^{n\times n})}<K_0$.
\end{itemize}
\end{prob}
Note that the qualitative version has done in \cite[Theorem 3.1]{StreetYaoVectorFields}, where the existence of such $\eta^1,\dots,\eta^n$ is done by \cite[Lemma 4.8]{StreetYaoVectorFields}.

In the proof of Theorem \ref{MainThm::QuantFro} we use the canonical coordinates $t\mapsto e^{t\cdot X}(p_0)$ to give the initial normalization. However when $X_1,\dots,X_m\in\Co^\alpha_\loc\backslash\Co^{\alpha+}_\loc$, we may not be able to talk about $e^{t\cdot X}$ (unless $\alpha=1$) due to the non-uniqueness of H\"older ODE flows. To prove the quantitative result we must bypass the method of canonical coordinates. Solving harmonic coordinates on a suitable Carnot-Carath\'eodory ball may be a way.

\subsection{Proof of Theorem \ref{MainThm::QuantFro}}\label{Section::PfQuantFro}

First by Theorem \ref{MainThm::SingFro}, since $X_1,\dots,X_m$ are log-Lipschitz and \eqref{Eqn::Quant::AssumptionM0} is a stronger assumption than \eqref{Eqn::MainThm::SingFro::Inv}. We can find a $n$-dimensional $\Co^{2-}$ submanifold $\Sf\subseteq\Mf$ containing $p_0$ such that $X_1,\dots,X_m$ span the tangent space of $\Sf$ at every point. Thus we can replace $\Mf$ by $\Sf$ in the following discussion.

By assumption $X_1,\dots,X_n$ are linearly independent at every point in $\Sf$, which form a local basis. Let $(\theta^1,\dots,\theta^n)$ be the dual basis of $(X_1,\dots,X_n)$, which are log-Lipschitz 1-forms in $\Sf$. The assumption \eqref{Eqn::Quant::AssumptionM0} implies that 
\begin{equation}\label{Eqn::QuantFro:Pf::dLambda}
    d\theta^k=\sum_{1\le i<j\le n}c_{ij}^k\theta^i\wedge\theta^j,\quad1\le k\le n.
\end{equation}

By Lemma \ref{Lem::ODE::FlowVsMetric} and Remark \ref{Rmk::ODE::FlowVsMetric} we see that $c_{ij}^k\in\Co^{\alpha-1}_\loc(\Mf)$, thus $d\theta^k\in \Co^{\alpha-1}_\loc(\Mf;\wedge^2T^*\Mf)$. By Corollary \ref{Cor::Rough1Form::BasicCaseforMainThm} we can find a $\Co^{\alpha+1}$-atlas $\Af$ for $\Mf$, which is compatible with the standard $\Co^{2-}$ structure of $\Mf$, such that $\theta^1,\dots,\theta^n$ are all $\Co^\alpha$.
Therefore by passing to the atlas $\Af$, we can assume from now on that $\theta^1,\dots,\theta^n,X_1,\dots,X_n$ are all $C^1$.

Let $\mu_1=\mu_1(n,\mu_0,M_0):=\frac1{10}\min(\mu_0,\frac1{n\sqrt n M_0})$.
By Definition \ref{Defn::ODE::AlmostInjRad} and Corollary \ref{Cor::ODE::InjRad}, the canonical coordinate $\Phi_0:B^n(0,\mu_1)\to\Mf$, $\Phi_0(t):=e^{t\cdot X'}(p_0)$ is a $C^1$ regular parameterization. By Proposition \ref{Prop::ODE::InjLemma} \ref{Item::ODE::InjLemma::A} and \ref{Item::ODE::InjLemma::HoldNorm}, we can find a $K_1=K_1(n,m,\mu_0,\alpha,M_0)>0$ such that $\Phi_0^*X'=(I+A_0')\nabla$ where $\|A_0\|_{C^0(B^n(0,\mu_1);\R^{n\times n})}\le \frac12$ and $\|A_0\|_{\Co^{\alpha-1}(B^n(0,\mu_1);\R^{n\times n})}\le K_1$.

Write $\hat\theta^i:=\Phi^*\theta^i$ and $\hat c_{ij}^k:=\Phi^*c_{ij}^k$ for $1\le i,j,k\le n$. By Proposition \ref{Prop::ODE::InjLemma} \ref{Item::ODE::InjLemma::DistNorm} we see that for $1\le i,j,k\le n$, $\|\hat c_{ij}^k\|_{\Co^{\alpha-1}(B^n(0,\mu_1))}\le 2\|\hat c_{ij}^k\|_{\Co^{\alpha-1}_{\Phi_0^*X'}(B^n(0,\mu_1))}\le  2\|c_{ij}^k\|_{\Co^{\alpha-1}_X(\Mf)}\le M_0$. By Proposition \ref{Prop::ODE::InjLemma} \ref{Item::ODE::InjLemma::HoldNorm} we get $\|A_0'\|_{\Co^{\alpha-1}(B^n(0,\mu_1);\R^{n\times n})}\lesssim_{n,\mu_1,M_0}1$. By \eqref{Eqn::QuantFro:Pf::dLambda} we have $d\hat\theta^k=\sum_{i<j}\hat c_{ij}^k\hat\theta^i\wedge\hat\theta^j$ for $1\le k\le n$. Therefore there is a $M_1=M_1(n,\mu_1,\alpha,M_0)>1$ such that $\sum_{1\le k\le n}\|\hat\theta^k\|_{\Co^{\alpha-1}}+\|d\hat\theta^k\|_{\Co^{\alpha-1}}\le M_1$.

Now let $c=c(n,\alpha-1,\alpha)>0$ be the constant in Theorem \ref{KeyThm::Rough1Form}. By Proposition \ref{Prop::Rough1Form::Scaling}, there are a $\kappa_0=\kappa_0(n,\alpha-1,\alpha,\mu_1,c,M_1)\in(0,\mu_0)$ and 1-forms $\lambda^1,\dots,\lambda^n\in\Co^{\alpha-1}(\B^n;T^*\B^n)$ that satisfies the assumptions of Theorem \ref{KeyThm::Rough1Form} and 
\begin{equation}\label{Eqn::QuantFro::ScalingCondition}
    \hat\theta^i|_{\frac13\B^n}=\tfrac1{\kappa_0}(\phi_{\kappa_0}^*\lambda^i)|_{\frac13\B^n}\text{ for }1\le i\le n,\quad \text{ where }\phi_{\kappa_0}:\B^n\to B^n(0,\mu_1),\quad\phi_{\kappa_0}(x)=\kappa_0\cdot x.
\end{equation}

By Theorem \ref{KeyThm::Rough1Form} (also by Remark \ref{Rmk::Rough1From::RmkofThm}), there is a $\Co^\alpha$-diffeomorphism $F:\B^n\to\B^n$, such that $B^n(F(0),\frac16)\subseteq F(\frac13\B^n)$, $F_*\lambda^1,\dots,F_*\lambda^n\in\Co^\alpha(\B^n;T^*\B^n)$ and $\|F_*\lambda-dy\|_{\Co^\alpha(\B^n;\R^{n\times n})}\le\frac1{12}$.

Let $Y=[Y_1,\dots,Y_n]^\top$ be a collection of vector fields which is the dual basis of $(F_*\lambda^1,\dots,F_*\lambda^n)$. We write $Y=:(I+A_1)\nabla$, thus as matrices $F_*\lambda=(I+A_1)^{-1}$. By Lemma \ref{Lem::Hold::Product} \ref{Item::Hold::Product::Hold2} and $\|F_*\lambda-dy\|_{\Co^\alpha(\B^n;\R^{n\times n})}\le\frac1{12}$, we see that
\begin{equation}\label{Eqn::QuantFro::A1Norm}
    \|A_1\|_{\Co^\alpha}=\|(F^*\lambda)^{-1}-I\|_{\Co^\alpha}\le\sum_{j=1}^\infty(2\|F_*\lambda-dy\|_{\Co^\alpha})^j\le\tfrac15.\quad\text{In particular }\|I+A_1(F(0))\|_{C^0}\le\tfrac65.
\end{equation}

Define an affine linear map 
\begin{equation}\label{Eqn::Quant::DefPsi}
    \psi(t):=(I+A_1(F(0)))\cdot\tfrac t8+F(0),\quad t\in\B^n.
\end{equation} 
Note that by \eqref{Eqn::QuantFro::A1Norm} we have $\psi(\B^n)\subseteq\frac65\cdot\frac18\B^n+F(0)\subset B^n(F(0),\frac16)$.

We define  $\Phi_1:\B^n\to B^n(0,\mu_1)$ and $\Phi:\B^n\to\Mf$ by
\begin{equation}\label{Eqn::Quant::DefPhi1}
    \Phi:=\Phi_0\circ\Phi_1,\quad\Phi_1(t):=\phi_{\kappa_0}\circ F^\Inv\circ\psi(t)=\kappa_0\cdot F^\Inv\mleft((I+A_1(F(0)))\cdot\tfrac t8+F(0)\mright),\quad t\in\B^n.
\end{equation}
Here $\Phi_1$ is well-defined because $\psi(\B^n)\subset B^n(F(0),\frac16)$, $B^n(F(0),\frac16)\subseteq F(\frac13\B^n)$ and $\phi_{\kappa_0}(\frac13\B^n)\subset\mu_1\B^n$. Clearly $\Phi_1(0)=\kappa_0\cdot F^\Inv(F(0))=0$ and thus $\Phi(0)=\Phi_0(0)=p_0$.

By \eqref{Eqn::QuantFro::ScalingCondition}, $\lambda^i\big|_{\frac13\B^n}=\frac1{\kappa_0}\cdot(\phi_{\kappa_0}^*\theta^i)\big|_{\frac13\B^n}$, and the fact that $F^\Inv\circ\psi(\B^n)\subseteq\frac13\B^n$, we have
\begin{equation}\label{Eqn::Quant::Phi1*Theta}
\begin{aligned}
    (\Phi_1^*\theta)(t)=[\Phi_1^*\theta^1(t),\dots,\Phi_1^*\theta^n(t)]&=\kappa_0\psi^*(F_*\lambda)=\kappa_0dy\cdot\psi^*((I+A_1(y))^{-1})=\kappa_0\cdot d\psi(t)\cdot(I+A_1(\psi(t)))
    \\
    &=\tfrac{\kappa_0}8dt\cdot(I+A_1(\psi(t)))^{-1}\cdot (I+A_1(F(0))).
\end{aligned}
\end{equation}
Here we use $dy=[dy^1,\dots,dy^n]$ and $d\psi(t)=[d\psi^1(t),\dots,d\psi^n(t)]$ as row vectors in \eqref{Eqn::Quant::Phi1*Theta}.


Since $\Phi_0^*X_1,\dots,\Phi_0^*X_n$ and $\Phi_1^*\theta^1,\dots,\Phi_1^*\theta^n$ are dual basis to each other, we can write
\begin{equation}\label{Eqn::Quant::EqnofPhi*Y}
    \Phi_0^*X=:\tfrac 8{\kappa_0}\cdot(I+ A')\tfrac\partial{\partial t},\quad\text{where } A'(t)=(I+A_1(F(0)))^{-1}\cdot(I+A_1(\psi(t)))-I,\quad t\in\B^n.
\end{equation}

Taking $t=0$ in \eqref{Eqn::Quant::EqnofPhi*Y}, since $\psi(0)=F(0)$, we get $A'(0)=0$.

Let $\widehat K:=\frac8{\kappa_0}$, since $\kappa_0=\kappa_0(n,\alpha,\mu_1,M_1)$ depends only on $n,\alpha,\mu_0,M_0$, we get $\widehat K=\widehat K(n,\alpha,\mu_0,M_0)$ as desired.

Since $\|\frac18(I+A_1)\|_{C^0}\le\frac3{20}<1$, by direct computation using the argument in the proof of  Lemma \ref{Lem::Hold::ScalingLem},  we get $\|A_1\circ\psi\|_{\Co^\alpha}\le\|A_1\|_{\Co^\alpha}\le\frac15$. Using the  power series $A'(t)=A_1(\psi(t))+\sum_{j=1}^\infty(-1)^jA_1(F(0))^j\cdot(1+A_1(\psi(t)))$, we have
\begin{equation*}
    \|A'\|_{\Co^\alpha(\B^n;\R^{n\times n})}\le\|A_1\circ\psi\|_{\Co^\alpha}+\sum_{j=1}^\infty|A_1(F(0))|_{\R^{n\times n}}^j\|I+A_1\circ\psi\|_{\Co^\alpha}\le\tfrac15+\sum_{j=1}^\infty(\tfrac15)^j\cdot\tfrac{6}{5}=\tfrac12.
\end{equation*}
This gives the desire control for $A'$. In particular we see that $\|\Phi^*X'\|_{\Co^\alpha(\B^n;\R^{n\times n})}=\|I+A\|_{\Co^\alpha(\B^n;\R^{n\times n})}\le\frac32$.

Finally the result $\|A''\|_{\Co^\alpha}\lesssim_{n,m,\alpha,\mu_0,M_0}1$ follow from \cite[Lemma 9.33]{CoordAdapted}. Recall that in our notation we only consider use the spaces $C^{0,\alpha-1}_X$ and $C^{1,\alpha-1}_X$.
\hfill\qedsymbol






\cleardoublepage
\phantomsection
\addcontentsline{toc}{chapter}{Bibliography}
\bibliographystyle{amsalpha}
\small\bibliography{references}

\providecommand{\bysame}{\leavevmode\hbox to3em{\hrulefill}\thinspace}
\providecommand{\MR}{\relax\ifhmode\unskip\space\fi MR }
\providecommand{\MRhref}[2]{%
  \href{http://www.ams.org/mathscinet-getitem?mr=#1}{#2}
}
\providecommand{\href}[2]{#2}
\begin{thebibliography}{HvNVW16}

\bibitem[BCD11]{BahouriCheminDanchin}
Hajer Bahouri, Jean-Yves Chemin, and Rapha\"{e}l Danchin, \emph{Fourier
  analysis and nonlinear partial differential equations}, Grundlehren der
  Mathematischen Wissenschaften [Fundamental Principles of Mathematical
  Sciences], vol. 343, Springer, Heidelberg, 2011. \MR{2768550}

\bibitem[BCH08]{Involutive}
Shiferaw Berhanu, Paulo~D. Cordaro, and Jorge Hounie, \emph{An introduction to
  involutive structures}, New Mathematical Monographs, vol.~6, Cambridge
  University Press, Cambridge, 2008. \MR{2397326}

\bibitem[DK81]{DeTurckKazdan}
Dennis~M. DeTurck and Jerry~L. Kazdan, \emph{Some regularity theorems in
  {R}iemannian geometry}, Ann. Sci. \'{E}cole Norm. Sup. (4) \textbf{14}
  (1981), no.~3, 249--260. \MR{644518}

\bibitem[Edw73]{Edwards}
C.~H. Edwards, Jr., \emph{Advanced calculus of several variables}, Academic
  Press (a subsidiary of Harcourt Brace Jovanovich, Publishers), New
  York-London, 1973. \MR{0352341}

\bibitem[Eva10]{Evans}
Lawrence~C. Evans, \emph{Partial differential equations}, second ed., Graduate
  Studies in Mathematics, vol.~19, American Mathematical Society, Providence,
  RI, 2010. \MR{2597943}

\bibitem[FR95]{DirichletBoundedness}
Jens Franke and Thomas Runst, \emph{Regular elliptic boundary value problems in
  {B}esov-{T}riebel-{L}izorkin spaces}, Math. Nachr. \textbf{174} (1995),
  113--149. \MR{1349041}

\bibitem[GKOS01]{GeomDist}
Michael Grosser, Michael Kunzinger, Michael Oberguggenberger, and Roland
  Steinbauer, \emph{Geometric theory of generalized functions with applications
  to general relativity}, Mathematics and its Applications, vol. 537, Kluwer
  Academic Publishers, Dordrecht, 2001. \MR{1883263}

\bibitem[Gon20]{Gong}
Xianghong Gong, \emph{A {F}robenius-{N}irenberg theorem with parameter}, J.
  Reine Angew. Math. \textbf{759} (2020), 101--159. \MR{4058177}

\bibitem[GT01]{GilbargTrudinger}
David Gilbarg and Neil~S. Trudinger, \emph{Elliptic partial differential
  equations of second order}, Classics in Mathematics, Springer-Verlag, Berlin,
  2001, Reprint of the 1998 edition. \MR{1814364}

\bibitem[Her60]{HermannFoliations}
Robert Hermann, \emph{On the differential geometry of foliations}, Ann. of
  Math. (2) \textbf{72} (1960), 445--457. \MR{142130}

\bibitem[H{\"{o}}r03]{Hormander}
Lars H{\"{o}}rmander, \emph{The analysis of linear partial differential
  operators. {I}}, Classics in Mathematics, Springer-Verlag, Berlin, 2003,
  Distribution theory and Fourier analysis, Reprint of the second (1990)
  edition [Springer, Berlin; MR1065993 (91m:35001a)]. \MR{1996773}

\bibitem[HT03]{RoughNN}
C.~Denson Hill and Michael Taylor, \emph{Integrability of rough almost complex
  structures}, J. Geom. Anal. \textbf{13} (2003), no.~1, 163--172. \MR{1967042}

\bibitem[HT07]{RoughComplex}
\bysame, \emph{The complex {F}robenius theorem for rough involutive
  structures}, Trans. Amer. Math. Soc. \textbf{359} (2007), no.~1, 293--322.
  \MR{2247892}

\bibitem[HvNVW16]{AnalBana}
Tuomas Hyt\"{o}nen, Jan van Neerven, Mark Veraar, and Lutz Weis, \emph{Analysis
  in {B}anach spaces. {V}ol. {I}. {M}artingales and {L}ittlewood-{P}aley
  theory}, Ergebnisse der Mathematik und ihrer Grenzgebiete. 3. Folge. A Series
  of Modern Surveys in Mathematics [Results in Mathematics and Related Areas.
  3rd Series. A Series of Modern Surveys in Mathematics], vol.~63, Springer,
  Cham, 2016. \MR{3617205}

\bibitem[Kre16]{KreitThesisGeneralizedHolder}
Damien Kreit, \emph{On generalized {H}\"{o}lder-{Z}ygmund spaces}, [Online]:
  \url{
  https://orbi.uliege.be/bitstream/2268/197514/1/Th%C3%A8se%20Kreit%20Damien.pdf},
  2016, Thesis (Ph.D.)--Universit\'e de Li\`ege, p.~120.

\bibitem[Lav18]{ShortGuideFrobenius}
Sylvain Lavau, \emph{A short guide through integration theorems of generalized
  distributions}, Differential Geom. Appl. \textbf{61} (2018), 42--58.
  \MR{3856749}

\bibitem[Lew14]{Lewis}
Andrew~D. Lewis, \emph{Generalised subbundles and distributions: A
  comprehensive review}, [Online]:
  \url{https://mast.queensu.ca/~andrew/notes/pdf/2011a.pdf}, 2014.

\bibitem[LTW19]{ContFro}
Stefano Luzzatto, Sina Tureli, and Khadim War, \emph{Integrability of
  continuous bundles}, J. Reine Angew. Math. \textbf{752} (2019), 229--264.
  \MR{3975643}

\bibitem[Lun92]{ShortFrobenius}
Albert~T. Lundell, \emph{A short proof of the {F}robenius theorem}, Proc. Amer.
  Math. Soc. \textbf{116} (1992), no.~4, 1131--1133. \MR{1145422}

\bibitem[Mal69]{Malgrange}
B.~Malgrange, \emph{Sur l'int\'{e}grabilit\'{e} des structures
  presque-complexes}, Symposia {M}athematica, {V}ol. {II} ({INDAM}, {R}ome,
  1968), Academic Press, London, 1969, pp.~289--296. \MR{0253383}

\bibitem[MM13]{MontanariMorbidelliSingularFrobenius}
Annamaria Montanari and Daniele Morbidelli, \emph{A {F}robenius-type theorem
  for singular {L}ipschitz distributions}, J. Math. Anal. Appl. \textbf{399}
  (2013), no.~2, 692--700. \MR{2996747}

\bibitem[Mor58]{Analyticity}
Charles~B. Morrey, Jr., \emph{On the analyticity of the solutions of analytic
  non-linear elliptic systems of partial differential equations. {I}.
  {A}nalyticity in the interior}, Amer. J. Math. \textbf{80} (1958), 198--218.
  \MR{106336}

\bibitem[Mun00]{MunkresTopology}
James~R. Munkres, \emph{Topology}, Prentice Hall, Inc., Upper Saddle River, NJ,
  2000, Second edition of [ MR0464128]. \MR{3728284}

\bibitem[Nir57]{Nirenberg}
Louis Nirenberg, \emph{A complex {F}robenius theorem}, Seminars on Analytic
  Functions I, Institute for Advanced Study, Princeton, 1957.

\bibitem[NW63]{NijenhuisWoolf}
Albert Nijenhuis and William~B. Woolf, \emph{Some integration problems in
  almost-complex and complex manifolds}, Ann. of Math. (2) \textbf{77} (1963),
  424--489. \MR{149505}

\bibitem[Osg98]{Osgood}
W.~F. Osgood, \emph{Beweis der {E}xistenz einer {L}\"{o}sung der
  {D}ifferentialgleichung {$\frac{{dy}}{{dx}} = f\left( {x,y} \right)$} ohne
  {H}inzunahme der {C}auchy-{L}ipschitz'schen {B}edingung}, Monatsh. Math.
  Phys. \textbf{9} (1898), no.~1, 331--345. \MR{1546565}

\bibitem[Ost65]{ColorDim}
Phillip~A. Ostrand, \emph{Dimension of metric spaces and {H}ilbert's problem
  {$13$}}, Bull. Amer. Math. Soc. \textbf{71} (1965), 619--622. \MR{177391}

\bibitem[Ram07]{Rampazzo}
Franco Rampazzo, \emph{Frobenius-type theorems for {L}ipschitz distributions},
  J. Differential Equations \textbf{243} (2007), no.~2, 270--300. \MR{2371789}

\bibitem[RS96]{BookProduct}
Thomas Runst and Winfried Sickel, \emph{Sobolev spaces of fractional order,
  {N}emytskij operators, and nonlinear partial differential equations}, De
  Gruyter Series in Nonlinear Analysis and Applications, vol.~3, Walter de
  Gruyter \& Co., Berlin, 1996. \MR{1419319}

\bibitem[RS07]{RampazzoSussmanCommutators}
Franco Rampazzo and H\'{e}ctor~J. Sussmann, \emph{Commutators of flow maps of
  nonsmooth vector fields}, J. Differential Equations \textbf{232} (2007),
  no.~1, 134--175. \MR{2281192}

\bibitem[Sim96]{Simic}
Slobodan Simi\'{c}, \emph{Lipschitz distributions and {A}nosov flows}, Proc.
  Amer. Math. Soc. \textbf{124} (1996), no.~6, 1869--1877. \MR{1328378}

\bibitem[SS18]{CoordAdapted}
Betsy Stovall and Brian Street, \emph{Coordinates adapted to vector fields:
  canonical coordinates}, Geom. Funct. Anal. \textbf{28} (2018), no.~6,
  1780--1862. \MR{3881835}

\bibitem[Str11]{StreetQuantFrobenius}
Brian Street, \emph{Multi-parameter {C}arnot-{C}arath\'{e}odory balls and the
  theorem of {F}robenius}, Rev. Mat. Iberoam. \textbf{27} (2011), no.~2,
  645--732. \MR{2848534}

\bibitem[Str20a]{CoordAdaptedIII}
\bysame, \emph{Coordinates adapted to vector fields {III}: real analyticity},
  Asian J. Math. \textbf{24} (2020), no.~6, 1029--1078. \MR{4312758}

\bibitem[Str20b]{SharpElliptic}
\bysame, \emph{Sharp regularity for the integrability of elliptic structures},
  J. Funct. Anal. \textbf{278} (2020), no.~1, 108290. \MR{4027743}

\bibitem[Str20c]{StreetSubHermitian}
\bysame, \emph{Sub-{H}ermitian geometry and the quantitative
  {N}ewlander-{N}irenberg theorem}, Adv. Math. \textbf{368} (2020), 107137, 87.
  \MR{4085140}

\bibitem[Str21]{CoordAdaptedII}
\bysame, \emph{Coordinates adapted to vector fields {II}: sharp results}, Amer.
  J. Math. \textbf{143} (2021), no.~6, 1791--1840. \MR{4349133}

\bibitem[SY22]{StreetYaoVectorFields}
Brian Street and Liding Yao, \emph{Improving the regularity of vector fields},
  J. Funct. Anal. \textbf{283} (2022), no.~5, Paper No. 109537. \MR{4429573}

\bibitem[Tay11]{TaylorPDE3}
Michael~E. Taylor, \emph{Partial differential equations {III}. {N}onlinear
  equations}, second ed., Applied Mathematical Sciences, vol. 117, Springer,
  New York, 2011. \MR{2744149}

\bibitem[Tr{\`{e}}92]{HypoAnalytic}
Fran\c{c}ois Tr{\`{e}}ves, \emph{Hypo-analytic structures}, Princeton
  Mathematical Series, vol.~40, Princeton University Press, Princeton, NJ,
  1992, Local theory. \MR{1200459}

\bibitem[Tr{\`{e}}06]{TrevesTopologicalSpaces}
\bysame, \emph{Topological vector spaces, distributions and kernels}, Dover
  Publications, Inc., Mineola, NY, 2006, Unabridged republication of the 1967
  original. \MR{2296978}

\bibitem[Tri95]{TriebelInterpolation}
Hans Triebel, \emph{Interpolation theory, function spaces, differential
  operators}, second ed., Johann Ambrosius Barth, Heidelberg, 1995.
  \MR{1328645}

\bibitem[Tri10]{Triebel1}
\bysame, \emph{Theory of function spaces}, Modern Birkh\"{a}user Classics,
  Birkh\"{a}user/Springer Basel AG, Basel, 2010, Reprint of 1983 edition
  [MR0730762], Also published in 1983 by Birkh\"{a}user Verlag [MR0781540].
  \MR{3024598}

\bibitem[Yao20]{YaoLLFro}
Liding Yao, \emph{The {F}robenius theorem for log-{L}ipschitz subbundles},
  arXiv:2004.07288v3, 2020.

\bibitem[Yao21]{LidingCounterNN}
\bysame, \emph{A counterexample to {$C^k$}-regularity for the
  {N}ewlander-{N}irenberg theorem}, Proc. Amer. Math. Soc. \textbf{149} (2021),
  no.~3, 1111--1115. \MR{4211866}

\bibitem[Yao22]{YaoCpxFro}
\bysame, \emph{Sharp {H}\"older regularity for {N}irenberg's complex
  {F}robenius theorem}, arXiv:2202.07729v1, 2022.

\bibitem[Zhu05]{ZhuSpaceBall}
Kehe Zhu, \emph{Spaces of holomorphic functions in the unit ball}, Graduate
  Texts in Mathematics, vol. 226, Springer-Verlag, New York, 2005. \MR{2115155}

\end{thebibliography}



\end{document}